\documentstyle{amsart}

\newtheorem{Theorem}{Theorem}[section]
\newtheorem{Lemma}[Theorem]{Lemma}
\newtheorem{Corollary}[Theorem]{Corollary}
\newtheorem{Proposition}[Theorem]{Proposition}
\newtheorem{Remark}[Theorem]{Remark}

\newtheorem{Example}[Theorem]{Example}

\newtheorem{Definition}[Theorem]{Definition}

\begin{document}

\title{Monomialization  of  morphisms\\
from 3 folds to surfaces}

\author{ Steven Dale Cutkosky}
\thanks{  partially supported by NSF}

\maketitle

\section{Introduction}

Suppose that $X$ is a nonsingular variety, over an algebraically closed
field $k$ of characteristic zero.

If $V\subset X$ is a nonsingular
subvariety, the 
blowup of $V$ is the morphism
$$
\pi:Y=\text{Proj}(\oplus_{n\ge 0}{\cal I}_{V}^n)\rightarrow X.
$$
If $p$ is a closed point of $Y$ and $\pi(p)=q$, 
there exist regular parameters $(x_1,\ldots,x_n)$ at $q$
and regular parameters $(y_1,\ldots,y_n)$ at $p$ such that
$$
x_1=x_2=\cdots x_r=0
$$
(with $r\le n=\text{dim}(X)$) are local equations of $V$ at $q$ and 
$$
x_1=y_1,x_2=y_1y_2,\cdots,x_r=y_1y_r,x_{r+1}=y_{r+1},\ldots,x_n=y_n.
$$
If $V=q$, so that $r=n$, $\pi$ is called the blowup of a point.

Another simple example of a morphism is a monomial morphism,
$\Phi:{\bold A}^n\rightarrow {\bold A}^m$
defined by
$$
\begin{array}{ll}
y_1=&x_1^{a_{11}}\cdots x_n^{a_{1n}}\\
&\vdots\\
y_m=&x_1^{a_{m1}}\cdots x_n^{a_{mn}}
\end{array}
$$
$\Phi$ is dominant if and only if $\text{rank}(a_{ij}) = m$.
This notion of a monomial morphism is a little too restrictive, so we extend it in the following way.

\begin{Definition}(Definition \ref{Def1082}) Suppose that $\Phi:X\rightarrow Y$ is a 
dominant
morphism of nonsingular $k$-varieties (where $k$ is a field of characteristic zero). $\Phi$ is
{\bf monomial}  if for all $p\in X$ there exists an \'etale neighborhood $U$ of $p$,
uniformizing parameters $(x_1,\ldots,x_n)$ on $U$, regular parameters 
$(y_1,\ldots,y_m)$ in ${\cal O}_{Y,\Phi(p)}$, and a matrix $(a_{ij})$ of nonnegative
integers (which necessarily has rank $m$) such that
$$
\begin{array}{ll}
y_1=&x_1^{a_{11}}\cdots x_n^{a_{1n}}\\
&\vdots\\
y_m=&x_1^{a_{m1}}\cdots x_n^{a_{mn}}
\end{array}
$$
\end{Definition}

\vskip .2truein

Suppose that 
\begin{equation}\label{eq*}
 \Phi:X\rightarrow Y
\end{equation}
 is a dominant morphism of   $k$-varieties, where $k$ is a field of characteristic 0. The structure of $\Phi$ is
extremely complicated. However, we can hope to 
construct a commutative diagram 
\begin{equation}\label{eq***}
\begin{array}{lll}
X_1&\stackrel{\Psi}{\rightarrow}&Y_1\\
\downarrow&&\downarrow\\
X&\stackrel{\Phi}{\rightarrow}&Y
\end{array}
\end{equation}
where the vertical maps are products of blowups of nonsingular subvarieties,
to obtain a morphism $\Psi:X_1\rightarrow Y_1$ which has a relatively simple structure.

The most optimistic conclusion we can hope for is to construct a  diagram
(\ref{eq***}) such that $\Psi$ is monomial. 

\begin{Definition}(Definition \ref{Def1082}) Suppose that $\Phi:X\rightarrow Y$ is a 
dominant morphism of  $k$-varieties.  A morphism $\Psi:X_1\rightarrow Y_1$ is a {\bf  monomialization} of $\Phi$
if there are sequences of blowups of nonsingular subvarieties
 $\alpha:X_{1}\rightarrow X$
and $\beta: Y_{1}\rightarrow Y$, and a morphism $\Psi:X_1\rightarrow Y_1$ such that  the diagram
\[
\begin{array}{lll}
X_1&\stackrel{\Psi}{\rightarrow}&Y_1\\
\downarrow&&\downarrow\\
X&\stackrel{\Phi}{\rightarrow}&Y
\end{array}
\]
commutes, and $\Psi$ is  a  monomial morphism. 
\end{Definition} 

In many cases a monomialization or something close to a monomialization exists 
 so it is natural to ask the following question.

\vskip .2truein
\noindent {\bf Question}  {\it Suppose that $\Phi:X\rightarrow Y$ is
a dominant morphism of $k$-varieties (over a field $k$
of characteristic zero). Does there exist a 
monomialization of $\Phi$?}
\vskip .2truein

By resolution of singularities and resolution of indeterminancy, we
easily reduce to the case where $X$ and $Y$ are nonsingular.

The characteristic of $k$ must be zero in the question. If char $k=p>0$, 
a monomialization may not exist even for curves.
$$
t=x^p+x^{p+1}
$$
gives a simple example of a mapping of curves which cannot
be monomialized, since
$\sqrt[p]{1+x}$ is inseparable over $k[x]$.

The obstruction to monomialization in positive characteristic is thus wild ramification.

In  \cite{C2}, we prove that a local analogue of the Question
 has a positive answer  for generically finite morphisms. 
A discussion of these results is given in section \ref{section2}.

In Section \ref{section3}, we outline short proofs of the positive answer to the
question in the previously known cases, a morphism to a curve and a morphism of
surfaces (\cite{AKi}, \cite{CP} in characteristic $p\ge 0$ when no wild ramification is
present).

In this paper we give a positive answer to the question in the case of a dominant morphism from a
3 fold to a surface.

\begin{Theorem}(Theorem \ref{Theorem1079})
Suppose that $\Phi:X\rightarrow S$ is a dominant morphism from a 3 fold $X$ to
a surface $S$ (over an algebraically closed field $k$ of characteristic zero). Then there exist sequences of
blowups of nonsingular subvarieties $X_1\rightarrow X$ and $S_1\rightarrow S$ such
that the induced map $\Phi_1:X_1\rightarrow S_1$ is a monomial morphism.
\end{Theorem}

From this we deduce that it is possible to toroidalize (\cite{KKMSD}, \cite{AK}, 
Definition \ref{Def1093})
a dominant morphism from a 3 fold to a surface. A toroidal morphism $X\rightarrow Y$ is
a morphism which is monomial with respect to fixed SNC divisors on $X$ and $Y$.

\begin{Theorem}(Theorem \ref{Theorem1078}) Suppose that $\Phi:X\rightarrow S$ is
a dominant morphism from a 3 fold $X$ to a surface $S$ (over an algebraically
closed field $k$ of characteristic zero) 
and $D_S$ is a reduced 1 cycle on $S$ such that $E_X=\Phi^{-1}(D_S)_{red}$ contains
$\text{sing}(X)$ and $\text{sing}(\Phi)$. Then there exist  sequences of
blowups of nonsingular subvarieties 
$\pi_1:X_1\rightarrow X$ and $\pi_2:S_1\rightarrow S$  such that the induced morphism $X_1\rightarrow S_1$ is a
toroidal morphism with respect to $\pi_2^{-1}(D_S)_{red}$ and $\pi_1^{-1}(E_X)_{red}$.
\end{Theorem}

Suppose that $\Phi:X\rightarrow Y$ is a dominant morphism of nonsingular $k$-varieties, and
$\text{dim}(Y)>1$.

To begin with, we  point out that monomialization is not a direct consequence of
embedded resolution of singularities and principalization of ideals.

Suppose that $p\in X$ is a point where $\Phi$ is not smooth, and $q=\Phi(p)$. Let 
$(y_1,\ldots,y_m)$ be regular parameters in ${\cal O}_{Y,q}$. By standard theorems
on resolution, we have a sequence of blowups of nonsingular subvarieties $\pi:X_1\rightarrow X$ such that
if $p_1\in\pi^{-1}(p)$, then there exist regular parameters $(x_1,\ldots,x_n)$ in
${\cal O}_{X_1,p_1}$, a matrix $(a_{ij})$ with nonnegative coefficients and units
$\delta_1,\ldots,\delta_m\in{\cal O}_{X_1,p_1}$ such that 
\begin{equation}\label{eq2000}
\begin{array}{ll}
y_1&=x_1^{a_{11}}\cdots x_n^{a_{1n}}\delta_1\\
&\vdots\\
y_m&=x_1^{a_{m1}}\cdots x_n^{a_{mn}}\delta_n
\end{array}
\end{equation}

In general, $p_1$ will lie on a single exceptional component of $\pi$, and $p_1$
will be disjoint from the strict transforms of  codimension 1 subschemes
of $X$ determined by $y_i=0$, $1\le i \le m$, on a neighborhood of $\Phi^{-1}(q)$. In this case
we will have $a_{ij}=0$ if $j>1$ and $(a_{ij})$ will have rank 1.

There thus cannot exist regular parameters $(\overline x_1,\ldots,\overline x_n)$
in $\hat{\cal O}_{X_1,p_1}$ such that 
$$
\begin{array}{ll}
y_1&=\overline x_1^{a_{11}}\cdots\overline x_n^{a_{1n}}\\
&\vdots\\
y_n&=\overline x_1^{a_{m1}}\cdots\overline x_n^{a_{mn}}
\end{array}
$$
since this would imply that $\text{rank}(a_{ij})=m>1$. 

In fact, in general it is necessary to blowup in both $X$ and $Y$ to construct
a monomialization. For instance, if we blowup a point $p$ on a nonsingular surface
$S$, blowup a point on the exceptional curve $E_1$, blowup the intersection point
of the new exceptional curve $E_2$ with the strict transform of $E_1$, then blowup
a general point on the new exceptional curve $E_3$ with exceptional curve $E_4$, we get a birational map
$\pi:S_1\rightarrow S$ such that if $p_1\in E_4$ is a general point we have regular
parameters $(u,v)$ in ${\cal O}_{S,p}$ and regular paramaters $(x,y)$ in $\hat{\cal O}_{S_1,p_1}$
such that
$$
u=x^2, v=\alpha x^3+x^4y.
$$
$\pi$ is not monomial at $p_1$ and further blowups over $S_1$ will produce a morphism
which is further from being monomial.

Suppose that $Y$ is a nonsingular surface. If $\pi_2:Y_1\rightarrow Y$ is a sequence
of blowups of points over $q\in Y$, and $q_1\in\pi_2^{-1}(q)$ is a point which only lies
on a single exceptional component $E$ of $\pi_2$, then there exist regular
parameters $(u,v)$ in ${\cal O}_{Y,q}$ and $(\overline x,\overline y)$ in $\hat{\cal O}_{Y_1,q_1}$
such that
\begin{equation}\label{eq2001}
\begin{array}{ll}
u&=\overline x^a\\
v&=P(\overline x)+\overline x^b\overline y
\end{array}
\end{equation}
where $a,b\in{\bold N}$ and $P(\overline x)$ is a polynomial of degree $\le b$.

If we perform a sequence of blowups of nonsingular subvarieties $\pi_1:X_1\rightarrow X$, and if
$p_1\in (\Phi\circ\pi_1)^{-1}(q)$ is such that $\hat{\cal O}_{X_1,p_1}$ has
regular parameters $(\overline x_1,\overline x_2,\overline x_3,\cdots,\overline x_n)$
such that  
\begin{equation}\label{eq2002}
\begin{array}{ll}
u&=\overline x_1^a\\
v&=P(\overline x_1)+\overline x_1^b\overline x_2
\end{array}
\end{equation}
of the form of (\ref{eq2001}), we will have a factorization $X_1\rightarrow S_1$
which is a morphism in a
neighborhood of $p_1$, and $X_1\rightarrow S_1$ will be monomial at $p_1$.

A strategy for monomializing a dominant morphism from a nonsingular variety $X$ to a
nonsingular surface $S$ is thus to first perform a sequence of blowups of nonsingular
subvarieties 
$\pi_1:X_1\rightarrow X$ so that for all points $p$ of $X_1$, appropriate
regular parameters $(u,v)$ in ${\cal O}_{S_1,q}$ where $q=\Phi\circ\pi_1(p)$ will have simple forms which we will
call prepared (Definition \ref{Def57} if $\text{dim}(X)=3$) which include the form of
(\ref{eq2002}).  This is accomplished if $\text{dim}(X)=3$ in Theorem \ref{Theorem58}.
Almost the entirety of this paper is devoted to proving this Theorem.

An interesting case when the existence of a global monomialization is still open 
is for birational morphisms of nonsingular, characteristic 0 varieties of dimension $\ge 3$.
Such birational maps are known to have a simple structure, since they
 can be factored by alternating
sequences of blowups and blowdowns \cite{AKMW}. A local form of factorization along 
a valuation is proven in Theorem 1.6 \cite{C2}.

\section{Local Monomialization}\label{section2}

A local version of monomialization is proven in \cite{C2}.

Suppose that 
$R\subset S$
is a local homomorphism of local rings essentially of finite type over a field $k$
and that $V$ is a valuation ring of the quotient field $K$ of $S$, such that $V$ dominates $S$. 
Then we can ask if there are sequences of monoidal transforms $R\rightarrow R'$ and $S\rightarrow S'$ such that
$V$ dominates $S'$, $S'$ dominates $R'$, and $R\rightarrow R'$ is a monomial mapping. 
\begin{equation}
\begin{array}{lll}
R'&\rightarrow&S'\subset V\\
\uparrow&&\uparrow\\
R&\rightarrow &S
\end{array}
\end{equation}

A monoidal transform of a local ring $R$ is the local ring $R'$ of a point in the blowup of
a nonsingular subvariety of $\text{spec}(R)$ such that $R'$ dominates $R$. If $R$ is a
regular local ring, then $R'$ is a regular local ring.

 \begin{Theorem}\label{TheoremA}(Monomialization)(Theorem 1.1 \cite{C2})
 Suppose that $R \subset S$ are   regular local rings, essentially of finite type over a 
field $k$ of characteristic zero,
such that  the quotient field $K$ of $S$ is a finite extension of the quotient field $J$
of $R$.

Let $V$ be a valuation ring of $K$ which dominates $S$. 
  Then there
exist sequences of  monoidal transforms 
 $R \rightarrow R'$ and $S \rightarrow S'$
such that $V$ dominates $S'$, $S'$ dominates $R'$ and 
there are regular parameters $(x_1, .... ,x_n)$
in $R'$,  $(y_1, ... ,y_n)$ in $S'$, units $\delta_1,\ldots,\delta_n\in S'$ and a matrix $(a_{ij})$ of
nonnegative integers such that  $\text{Det}(a_{ij}) \ne 0$ and 
\begin{equation}\label{eq1*}
\begin{array}{lll}
x_1 &=& y_1^{a_{11}} ..... y_n^{a_{1n}}\delta_1\\
&&\vdots\\   
x_n &=& y_1^{a_{n1}} ..... y_n^{a_{nn}}\delta_n.
\end{array}
\end{equation}
\end{Theorem}

Thus (since $\text{char}(k)=0$) there exists an etale extension $S'\rightarrow S''$
where $S''$ has regular parameters
$\overline y_1,\ldots,\overline y_n$ such that $x_1,\ldots,x_n$ are pure monomials in $\overline y_1,\ldots,\overline y_n$.

The standard theorems on resolution of singularities allow one to easily find $R'$ and $S'$ such that (\ref{eq1*}) holds,
but, in general, we will not  have the essential condition $\text{Det}(a_{ij}) \ne 0$. The difficulty of the
problem is to achieve this condition.

It is an interesting open problem to prove Theorem \ref{TheoremA} in  positive
characteristic, even in dimension 2. Theorem \ref{TheoremA} implies simultaneous resolution
from above \cite{C3}, which is a key step in a program of Abhyankar's for proving resolution in
positive characteristic. This method is completely worked out by Abhyankar in dimension 2
\cite{Ab1}.

A quasi-complete variety over a field $k$ is an integral finite type $k$-scheme which satisfies the 
existence part of the valuative criterion for properness (c.f. Chapter 0 \cite{H1}
where the notion is called complete).  
Quasi-complete and separated is equivalent to proper.

The construction of a monomialization by  quasi-complete varieties follows from Theorem \ref{TheoremA}.

\begin{Theorem}\label{TheoremB}(Theorem 1.2 \cite{C2}) Let $k$ be a field of characteristic zero, $\Phi:X\rightarrow Y$ a generically finite
morphism of  nonsingular proper  $k$-varieties.   Then there are
  birational morphisms of nonsingular quasi-complete  $k$-varieties $\alpha:X_{1}\rightarrow X$
and $\beta: Y_{1}\rightarrow Y$, and a locally monomial morphism $\Psi:X_1\rightarrow Y_1$ such that  the diagram
\[
\begin{array}{lll}
X_1&{\Psi}\atop{\rightarrow}&Y_1\\
\downarrow&&\downarrow\\
X&{\Phi}\atop{\rightarrow}&Y
\end{array}
\]
commutes and $\alpha$ and $\beta$ are
 locally  products of blowups of nonsingular subvarieties. 
That is, for every $z\in X_{1}$, there
exist affine neighborhoods $V_1$ of $z$, $V$ of $x=\alpha(z)$, such that $\alpha:V_1\rightarrow V$ 
is a  finite product of monoidal transforms, and  
there
exist affine neighborhoods $W_1$ of $\Psi(z)$, $W$ of $y=\alpha(\Psi(z))$, such that $\beta:W_1\rightarrow W$ 
is a  finite product of monoidal transforms.
\end{Theorem}

In this Theorem, a monoidal transform of a nonsingular $k$-scheme $S$ is the map $T\rightarrow S$ induced by an open subset $T$ of
$\text{Proj}(\oplus {\cal I}^n)$, where ${\cal I}$ is the ideal sheaf of a nonsingular subvariety of $S$.

Theorems 1.1 and 1.2 of \cite{C2} are analogues for morphisms of the Theorems on
local uniformization and local resolution of singularities of varieties of Zariski \cite{Z1},
\cite{Z2}.

\section{Monomialization of  Morphisms in Low Dimensions}\label{section3}

We will outline proofs of monomialization in the previously known cases.
Suppose that $k$ is an algebraically closed field of characteristic zero
and $\Phi:X\rightarrow Y$ is a dominant morphism of nonsingular $k$ varieties.

Let $\text{sing}(\Phi)$ be the closed subset of $X$ where $\Phi$ is not
smooth.

If $\Phi$ is a dominant morphism from a  variety to a curve, the existence of a global monomialization 
follows immediately from resolution of singularities.
In fact, it is really a restatement of embedded resolution of hypersurface
singularities.

\begin{Theorem} Suppose that
 $\Phi:X\rightarrow C$ is a dominant morphism from a  $k$-variety to a curve. Then $\Phi$ has
a   monomialization.
\end{Theorem}

\begin{pf}
Suppose that $\Phi:X\rightarrow C$ where $C$ is a nonsingular curve,
$X$ is a nonsingular $n$ fold. 
$\Phi(\text{sing}(\Phi))$  is a finite number of points of $C$, so we 
may fix a regular parameter $t$ at a point in $q\in\Phi(\text{sing}(\Phi))$,
and monomialize the mapping above $q$.

By embedded resolution of hypersurfaces, there exists a sequence
of blowups of nonsingular subvarieties which dominate subvarieties of $\Phi^{-1}(q)$,
$\pi:X_1\rightarrow X$   such that for all $p\in (\Phi\circ\pi)^{-1}(q)$,
there exists regular parameters $(x_1,\ldots,x_n)$ at $p$ such that
$$
t=ux_1^{a_1}\cdots x_n^{a_n}
$$
where $a_1>0$, $u\in {\cal O}_{X_1,p}$ is a unit. If $x_1=\overline x_1 u^{-\frac{1}{a_1}}$, we have
$$
t=\overline x_1^{a_1}\cdots x_n^{a_n}.
$$
\end{pf}

If $\Phi:T\rightarrow S$ is a dominant morphism of surfaces, 
monomialization is not a direct corollary of resolution of singularities.
One proof of monomialization in this case (over $\bold C$) is given by
Akbulut and King in \cite{AKi}.

In our paper \cite{CP} with Oliver Piltant, we show that if $L$ is a perfect field and $\Phi:T\rightarrow S$ is a dominant morphism of $L$-surfaces,
then $\Phi$ can be monomialized if $\Phi$ is unramified. That is, no
wild ramification occurs with respect to any divisorial valuation of 
$L(T)$ over $L(S)$.   
This condition occurs, for instance, if $p\not\,\mid [K:L(S)]$ where $K$ is a 
Galois closure of $L(T)$ over $L(S)$.

 We will now outline
a simple proof of monomialization for morphisms of surfaces (when $k$ is algebraically
closed of characteristic zero).

\begin{Theorem}
Suppose that $\Phi:T\rightarrow S$ is a dominant morphism of surfaces
over $k$. Then $\Phi$ has a  monomialization.
\end{Theorem}

If $\Phi$ is a  monomial mapping, then $\Phi$ comes from an expression
\begin{equation}\label{Ieq5}
\begin{array}{ll}
u&=x^ay^b\\
v&=x^cy^d
\end{array}
\end{equation}
where $ad-bc\ne 0$.

$\text{sing}(\Phi)$ must be contained in $xy=0$. At a point $p$ on $x=0$
we have regular parameters $(\hat x,\hat y)$ in $\hat{\cal O}_{X,p}$
such that 
$$
\hat x=x,\hat y=y-\alpha
$$
for some $\alpha\in k$. If $a>0$ and $c>0$  we have  
\begin{equation}\label{Ieq6}
\begin{array}{ll}
u&=\hat x^a(\hat y+\alpha)^b=\overline x^a\\
v&=\hat x^c(\hat y+\alpha)^d=\beta \overline x^c+\overline x^c\overline y
\end{array}
\end{equation}
where 
$$
\hat x=\overline x(\hat y+\alpha)^{-\frac{b}{a}},
\overline y = (\hat y+\alpha)^{d-\frac{cb}{a}}-\beta
$$
with $\beta = \alpha^{d-\frac{cb}{a}}$.

If $a=0$ or $c=0$ we also obtain a form (\ref{Ieq6}) with respect to regular parameters
$(u_1,v_1)$ in ${\cal O}_{S,\Phi(p)}$.

Thus $\Phi$ is monomial at a point $p$ if and only if there exist
regular parameters in $\hat{\cal O}_{X,p}$ such that one of the forms
(\ref{Ieq5}) or (\ref{Ieq6}) hold.

We will say that $\Phi$ is prepared at $p\in T$
if  there exist regular parameters $(u,v)$ in ${\cal O}_{S,\Phi(p)}$,
 regular parameters $(x,y)$ in $\hat{\cal O}_{T,p}$, and a power series $P$
such that 
one of the following forms  holds at $p$. 
\begin{equation}\label{Ieq7}
\begin{array}{ll}
u&=x^a\\
v&=P(x)+x^cy
\end{array} 
\end{equation}
or 
\begin{equation}\label{Ieq8}
\begin{array}{ll}
u&=(x^ay^b)^m\\
v&=P(x^ay^b)+x^cy^d
\end{array} 
\end{equation}
where $(a,b)=1$ and $ad-bc\ne 0$.

We first observe that by resolution of singularities and indeterminancy,
there exists a commutative diagram
$$
\begin{array}{lll}
T_1&\stackrel{\Phi_1}{\rightarrow}&S_1\\
\downarrow&&\downarrow\\
T&\stackrel{\Phi}{\rightarrow}&S
\end{array}
$$
where the vertical maps are products of blowups of points, 
$\text{sing}(\Phi_1)$ is a simple normal crossings (SNC) divisor,
and for all $p\in \text{sing}(\Phi_1)$, there exist regular
parameters $(u,v)$ at $\Phi_1(p)$ such that $u=0$ is a local equation
of $\text{sing}(\Phi_1)$ at $p$.

The essential observation is that $\Phi_1$ is now prepared.
We give a  simple proof  that appears in \cite{AKi}. 

\begin{Lemma}\label{Lemma2004} $\Phi_1$ is prepared.
\end{Lemma}
\begin{pf} Suppose that $p\in T_1$. With our assumptions, one of the following must hold at $p$. 
\begin{equation}\label{Ieq9}
\begin{array}{ll}
u&=x^a\\
u_xv_y-u_yv_x&=\delta x^e
\end{array}
\end{equation}
where $\delta$ is a unit
or 
\begin{equation}\label{Ieq10}
\begin{array}{ll}
u&=(x^ay^b)^m\\
u_xv_y-u_yv_x&=\delta x^ey^f
\end{array}
\end{equation}
where $a,b,e,f>0$, $(a,b)=1$ and $\delta$ is a unit.

Write $v=\sum a_{ij}x^iy^j$ with $a_{ij}\in k$.
First suppose that (\ref{Ieq9}) holds. Then $ax^{a-1}v_y=\delta x^e$
implies we have the form (\ref{Ieq7})
(after making a change of parameters in $\hat{\cal O}_{T_1,p}$). Now suppose that (\ref{Ieq10}) holds.
$$
u_xv_y-u_yv_x=\sum m(aj-bi)a_{ij}x^{am+i-1}y^{bm+j-1}
=\delta x^ey^f.
$$
Thus 
$$
v=\sum_{aj-bi=0}a_{ij}x^iy^j+\epsilon x^{e+1-am}y^{f+1-bm}
$$
where $\epsilon$ is a unit. After making a change of parameters, 
multiplying $x$ by a unit, and multiplying $y$ by a unit, we get the
form  (\ref{Ieq8}).
\end{pf}

It is now not difficult to construct a  monomialization. 
We must blowup points $q$ on $S_1$ over which some point is not monomial, and 
blowup points on $T_1$ to make $m_q{\cal O}_{T_1}$ principal. If we iterate this
procedure, it can be shown that we construct a commutative diagram
 $$
\begin{array}{lll}
T_2&\stackrel{\Phi_2}{\rightarrow}&S_2\\
\downarrow&&\downarrow\\
T_1&\stackrel{\Phi_1}{\rightarrow}&S_1
\end{array}
$$
such that $\Phi_2$ is monomial. 

\section{An overview of the proof of\\
Monomialization of  morphisms from 3 folds to surfaces} \label{section4}

Suppose that $k$ is an algebraically closed field of characteristic zero, and
$\Phi:X\rightarrow Y$ is a dominant morphism of  nonsingular $k$-varieties.

A natural first step in monomializing a morphism $\Phi:X\rightarrow Y$
is to use resolution of singularities and resolution of indeterminancy to construct a
commutative diagram
$$
\begin{array}{lll}
X_1&\stackrel{\Phi_1}{\rightarrow}&Y_1\\
\downarrow&&\downarrow\\
X&\stackrel{\Phi}{\rightarrow}&Y
\end{array}
$$
where the vertical maps are products of blowups of nonsingular subvarieties, 
$\text{sing}(\Phi_1)$ is a simple normal crossings (SNC) divisor,
and for all $p\in \text{sing}(\Phi_1)$, there exist regular
parameters $(u_1,\ldots,u_n)$ at $\Phi_1(p)$ such that $u_1=0$ is a local equation
of $\text{sing}(\Phi_1)$ at $p$. 

We observed that if $X$ and $Y$ are surfaces, then
$\Phi_1$ is prepared.
Unfortunately, even for morphisms from a 3 fold to a surface,
$\Phi_1$ may be quite complicated (Examples \ref{Example1099}, \ref{Example2002}).

A key step in the local proof of monomialization, 
Theorem \ref{TheoremA}, is to define a new invariant, which measures how far the situation 
is from a specific form which is close to being monomial. In the local valuation theoretic proof 
we make use of special products of monoidal transforms defined by Zariski called Perron transforms \cite{Z2}.
Under appropriate application of Perron transforms our invariant does not increase, and we can 
in fact make the invariant decrease, by an appropriate algorithm.

An essential difficulty globally is that our invariant can increase after a permissible monoidal transform (Example \ref{Example2003}).
This is a significant difference from resolution of singularities, where a foundational result is that
 the multiplicity of an ideal does not
go up under permissible blowups.

We will give a brief overview of the proof of Theorem \ref{Theorem1079} (Monomialization
of morphisms from 3 folds to surfaces). 
\vskip .2truein
{\bf Step 1.} First construct a diagram
\[
\begin{array}{lll}
X'&\stackrel{\Phi'}{\rightarrow}&S'\\
\downarrow&&\downarrow\\
X&\stackrel{\Phi}{\rightarrow}&S
\end{array}
\]
where the vertical maps are products of blowups of nonsingular subvarieties
 such that $X'$, $S'$ are nonsingular, there exist reduced SNCS divisors $D_{S'}$ on $S'$, $E_{X'}=(\Phi')^{-1}(D_{S'})_{red}$ on $X'$ such that 
 $\text{sing}(\Phi')\subset E_{X'}$ and  components of
$E_{X'}$ on $X'$ dominating distinct components of $D_{S'}$ are disjoint. Such a morphism
$\Phi'$ will be called weakly prepared (Definition \ref{Def1085} and Lemma \ref{Lemma1048}).

For all $p\in X'$ there exist regular parameters $(u,v)$ in ${\cal O}_{S',q}$ ($q=\Phi'(p)$)
and regular parameters $(x,y,z)$ in $\hat{\cal O}_{X',p}$ such that $u=0$ is a local
equation of $E_{X'}$,  $u=0$ (or $uv=0$) is a local equation of $D_{S'}$ and 
exactly one of
the following cases hold:
\begin{enumerate}
\item 
$$
u=x^a, v=P(x)+x^bF
$$
 where $x\not \,\mid F$, $F$ has no terms which are monomials in $x$.
\item 
$$
u=(x^ay^b)^m, v=P(x^ay^b)+x^cy^dF
$$
where $(a,b)=1$, 
$x\not\, \mid F$, $y\not\, \mid F$, $x^cy^dF$ has no terms which are monomials in $x^ay^b$.
\item 
$$
u=(x^ay^bz^c)^m, v=P(x^ay^bz^c)+x^dy^ez^fF
$$
where $(a,b,c)=1$, $x \not\, \mid F$, $y \not\, \mid F$, $z \not\,  \mid F$,
$x^dy^ez^fF$ has no terms which are monomials in $x^ay^bz^c$.
\end{enumerate}

The structure of the singularities of $F$ can be very complicated (Examples \ref{Example1099}
and \ref{Example2002}). This is in sharp contrast to the case of a morphism of surfaces
(Lemma \ref{Lemma2004}).

Our main invariant is
$$
\nu(p)=\text{mult}(F).
$$
This invariant is independent of parameters in the forms above.
$$
S_r(X')=\{p\in X'\mid\nu(p)\ge r\}
$$
is a constructible (but not Zariski closed) subset of $X'$ (Proposition \ref{Prop1} and
Example \ref{Ex1}).

\vskip .2truein
{\bf Step 2.} This is the difficult step. We construct a commutative diagram
$$
\begin{array}{lll}
X''&&\\
\downarrow\lambda&\stackrel{\,\,\Phi''}{\searrow}&\\
X'&{\Phi'}\atop{\rightarrow}&S'
\end{array}
$$
so that everywhere we have one of the forms:
\begin{enumerate}
\item $u=x^a, v=P(x)+x^by$,
\item $u=(x^ay^b)^m, v=P(x^ay^b)+x^cy^d$,
\item $u=(x^ay^b)^m, v=P(x^ay^b)+x^cy^dz$,
\item $u=(x^ay^bz^c)^m, v=P(x^ay^bz^c)+x^dy^ez^f$
with
$$
\text{rank}\left\{\begin{array}{lll}
a&b&c\\
d&e&f
\end{array}
\right\}=2.
$$
\end{enumerate}

We impose that further condition that 1. - 4. are compatible with the reduced SNC divisors
$D_{S'}$ and $E_{X''}=(\Phi'')^{-1}(D_{S'})_{red}$.
$u=0$ is a local equation of $E_{X''}$, $u=0$ (or $uv=0$) are local equations of $D_{S'}$
in the above forms.
We will say that $\Phi''$ is prepared
(Definition \ref{Def57}).
This is accomplished in Theorem \ref{Theorem1050}.

We use descending induction on 
$$
r=\text{max}\{t\mid \nu(p)=t\text{ for some }p\in X\}
$$
to achieve the conclusions of the Theorem. A major difficulty is that, unlike in the
case of resolution of singularities, $\nu(p)$ can go up after blowing up a point
or a nonsingular curve
(Example \ref{Example2003}).

However, $\nu(p)$ can go up by at most 1, and some other invariants get better,
or at least no worse. For a local resoution, we reduce to two difficult cases
(Sections \ref{Spec1} and \ref{Spec2}) which we settle by blowing up  generic
curves on $E_{X'}$ through a particular point, and use a generalization of
Abhynakar's Good Point Algorithm (\cite{Ab6}, \cite{Li}) to achieve an improvement.
This depends on arithmetic information which is captured in this algorithm.

\vskip .2truein
{\bf Step 3.} We construct a commutative diagram
\[
\begin{array}{lll}
X'''&\stackrel{\Phi'''}{\rightarrow}&S''\\
\downarrow&&\downarrow\\
X''&\stackrel{\Phi''}{\rightarrow}&S'
\end{array}
\]
such that  $X'''\rightarrow X''$ is a product of blowups of nonsingular curves, 
$S''\rightarrow S'$ is a product of blowups of points and $\Phi'''$ is monomial.
This is accomplished in  Theorem \ref{Theorem1068}.

$\pi:S''\rightarrow S'$ is a sequence of blowups of points. If $q\in S'$ and $q_1\in\pi^{-1}(q)$
then there exist regular parameters $(u,v)$ in ${\cal O}_{S',q}$ and $(u_1,v_1)$ in
$\hat{\cal O}_{S'',q_1}$ such that
$$
\begin{array}{ll}
u&=u_1^a\\
v&=P(u_1)+u_1^bv_1
\end{array}
$$
or
$$
\begin{array}{ll}
u&=(u_1^av_1^b)^m\\
v&=P(u_1^av_1^b)+u_1^cv_1^d
\end{array}
$$
with $ad-bc\ne 0$ and $(a,b)=1$.

If $p\in X''$ is a point of the  form 1. of Step 2, then there exists 
$\overline \pi:S_1\rightarrow S'$
and $q_1\in\overline \pi^{-1}(q)$ with regular parameters $(u_1,v_1)$ in ${\cal O}_{S_1,q_1}$,
$(\overline x,\overline y,\overline z)$ in $\hat{\cal O}_{X'',p}$ such that
$$
\begin{array}{ll}
u_1&=\overline x^{\overline a}\\
v_1&=\overline x^b(\alpha+\overline y).
\end{array}
$$

We have an essentially canonical procedure for achieving Step 3. We blowup on $S'$ the
(finitely many) images of all non monomial points of $X''$, then blowup nonsingular
curves on $X''$ to resolve the indeterminancy of the resulting rational map. An invariant improves.
By induction we eventually construct $\Phi'''$.

\section{Notations}

We will suppose that $k$ is an algebraically closed field of characteristic zero.
By a variety we will mean a separated, integral finite type $k$-scheme.

Suppose that $Z$ is a variety and $p\in Z$. Then $m_p$ will denote the maximal
ideal of ${\cal O}_{Z,p}$. 

\begin{Definition}\label{Def1098}
 A reduced divisor $D$ on a nonsingular variety $Z$ of dimension $n$ is
a simple normal crossing divisor (SNC divisor) if
\begin{enumerate}
\item All components of $D$ are nonsingular.
\item Suppose that $p\in X$. Let $D_1,\ldots, D_s$ be the components of $D$ containing $p$. 
Then $s\le n$ and there exist regular parameters $(x_1,\ldots,x_n)$ in ${\cal O}_{X,p}$ such 
that $x_i=0$ are local equations of $D_i$ at $p$ for $1\le i\le s$.
\end{enumerate}
\end{Definition}

A curve is a 1 dimensional $k$ variety. A surface is a 2 dimensional $k$ variety.
A 3 fold is a 3 dimensional $k$ variety. A point of a variety will mean a closed point.

By a generic point or a generic curve on a variety $Z$, we will mean
a point or a curve which satisfies a good condition which holds
on an open set (in some parametrizing space) of points or curves.

Suppose that $Z$ is a variety and $p\in Z$. the blowup of $p$ or the
quadratic transform of $p$ will denote
$Z_1=\text{Proj}(\oplus_{n\ge 0}m_p^n)$.  If $V\subset Z$ is a
nonsingular subvariety then the blowup of $V$ or the monodial transform
of $Z$ centered at $V$ will denote
$Z_1=\text{Proj}(\oplus_{n\ge 0}{\cal I}_{V}^n)$.

If $R$ is a regular local ring with maximal ideal $m$, then a quadratic
transform of $R$ is $R_1=R[\frac{m}{x}]_{m_1}$ where $0\ne x\in m$ and
$m_1$ is a maximal
ideal of $R_1$.

Suppose that $P(x)=\sum_{i=0}^{\infty}b_ix^i\in k[[x]]$ is a series. Given $t\in \bold N$,
$P_t(x)$ will denote the polynomial
$$
P_t(x)=\sum_{i=0}^tb_ix^i.
$$

Given a series $f(x_1,\ldots,x_n)\in k[[x_1,\ldots,x_n]]$ $\nu(f)$, $\text{mult}(f)$ or $\text{ord}(f)$
will denote the order of $f$.

If $x\in {\bold Q}$, $[x]=n$ if $n\in{\bold N}$, $n\le x<n+1$. $\{x\}=[x]-x$.
The greatest common divisor of $a_1,\ldots, a_n\in\bold N$ will be denoted by
$(a_1,\ldots,a_n)$.

\section{The invariant $\nu$}

\begin{Definition}\label{Def1085}
Suppose that $\Phi_X:X\rightarrow S$ is a  dominant morphism from a 
nonsingular 3 fold
$X$ to a nonsingular surface $S$, with  reduced SNC divisors $D_S$ on $S$ and $E_X$ on $X$ such that
$\Phi_X^{-1}(D_S)_{\text{red}}=E_X$. 
Let $\text{sing}(\Phi_X)$ be the locus of singular points of $\Phi_X$.

We will say that $\Phi_X$ is weakly prepared if 
\begin{enumerate}
\item $\text{sing}(\Phi_X)\subset E_X$
and 
\item If $p\in S$ is a singular point of $D_S$, $C_1$ and $C_2$ are the components of $D_S$ containing $p$,  $T_1$ is a 
component of $E_X$ dominating $C_1$ and $T_2$ is a component of $E_X$ dominating
$C_2$ then $T_1$ and $T_2$ are disjoint.
\end{enumerate}
\end{Definition}

\begin{Lemma}\label{Lemma1048}
Suppose that $\Phi:X\rightarrow S$ is a dominant morphism from
a 3 fold $X$ to a surface $S$, $D_S$ is a reduced Weil divisor on $S$
such that $\text{sing}(\Phi)\subset\Phi^{-1}(D_S)$ and the singular
locus of $X$, $\text{sing}(X)\subset\Phi^{-1}(D_S)$. Then there
exists a commutative diagram
$$
\begin{array}{lll}
X_1&\stackrel{\Phi_1}{\rightarrow}& S_1\\
\downarrow\pi_1&&\downarrow\pi_2\\
X&\stackrel{\Phi}{\rightarrow}&S
\end{array}
$$
such that $\pi_1$ and $\pi_2$ are products of blowups of nonsingular
subvarieties, and if $D_{S_1}=\pi_2^{-1}(D_S)_{\text{red}}$,
$E_{X_1}=(\Phi\circ\pi_1)^{-1}(D_S)_{\text{red}}$, then $\phi_1$
is weakly prepared.
\end{Lemma}

\begin{pf}
By resolution of singularities and resolution of indeterminancy of
mappings \cite{H1}, there exists
$$
\begin{array}{ccc}
\overline X&\stackrel{\overline \Phi}{\rightarrow}&\overline S\\
\downarrow&&\downarrow\\
X&\rightarrow&S
\end{array}
$$
such that $\overline X$ and $\overline S$ are nonsingular, $\pi^{-1}(D_S)_{red}=D_{\overline S}$ and $E_{\overline X}=\overline \Phi^{-1}(D_{\overline S})_{red}$
are SNC divisors.

Suppose that $E_1$ and $E_2$ are components of $E_{\overline X}$ which 
dominate distinct components $C_1$ and $C_2$ of $\overline S$.
If $E_1\cap E_2\ne \emptyset$ then there exists a sequence of blowups 
$\overline X_1\rightarrow \overline X$ with nonsingular centers which
map into $C_1\cap C_2$ with induced map $\overline \Phi_1:\overline X_1\rightarrow \overline S$ such that the strict transform of $E_1$ and
$E_2$ are disjoint on $\overline X_1$, and $E_{\overline X_1}=\overline \Phi_1^*(D_{\overline S})_{\text{red}}$ is a SNC divisor.

One way to construct this is to blow up the conductor of $E_1\cup E_2$
to separate the strict transforms of $E_1$ and $E_2$ (c.f. section 2
of \cite{CS}), and then resolve the singularities of the resulting
variety.

Iterating this procedure, we construct a weakly prepared morphism.
\end{pf}

\begin{Example}\label{Example1099} The structure of weakly prepared morphisms can
be quite complicated.

Consider the germ of maps
$$
u=x^a, v=x^cF
$$
with $a\ge 2$, $c\ge0$
where 
$$
F = x^rz+h(x,y)
$$
where  $h$ is arbitrary.
 The singular locus of this map germ is
the variety definied by the ideal where the jacobian has rank $<2$. That is, the variety with ideal
J= $\sqrt{(x^{a+c-1}F_y,x^{a+c-1}F_z)}$.
Since $F_z=x^r$, we have that $\sqrt{J}=(x)$.
\end{Example}

 Examples of this kind can be used to construct weakly prepared projective morphisms satisfying the assumptions of $\Phi_1$, by resolving the indeterminancy 
of the induced rational map ${\bold P}^3\rightarrow {\bold P}^2$.
A reasonably easy  example to calculate is  
$$
\begin{array}{ll}
u&=x^2\\
v&=y^2+xz.
\end{array}
$$

\begin{Example}\label{Example2002}
Another example of a weakly prepared morphism.
\end{Example}

Consider the (formal) germ of maps
$$
u=x^a, v=x^cF
$$
where 
$$
F = \sum_{i>0, j\ge 0}\frac{i^j}{j!}a_i(x)y^iz^j+x^rz
$$
where $a_i(x)$ are arbitrary series, $a\ge 2$, $c\ge 0$.
The singular locus of this map germ is defined by
$J=x^{a+c-1}(F_y,F_z)$. Since
$F_z-yF_y=x^r$,
$\sqrt{J}=(x)$.

Throughout this section we will suppose that $\Phi_X:X\rightarrow S$ is weakly prepared.

We define permissible parameters $(u,v)$ at points $q\in D_S$ by
the following rules
\vskip .2truein
{\bf 1.} If  $q$ is a nonsingular point of $D_S$, then
regular parameters $(u,v)$ in ${\cal O}_{S,q}$ are permissible parameters at $q$ if 
$u=0$ is a local equation for $D_S$.
Necessarily, $u=0$ is a local equation for $E_X$  in ${\cal O}_{X,p}$ for all $p\in \Phi_X^{-1}(q)$.
\vskip .2truein
{\bf 2.} If  $q$ is a singular point of $D_S$, then regular
parameters $(u,v)$ in ${\cal O}_{S,q}$ are permissible parameters at $q$ if $uv=0$ 
is a local equation for $D_S$ at $q$. Necessarily, $uv=0$ is a local equation for $E_X$ at $p$
for all $p\in \Phi_X^{-1}(q)$ and  either $u=0$ or $v=0$ is a local equation
of $E_X$ at $p$.
\vskip .2truein

\begin{Definition}\label{Def650} Suppose that $(u,v)$ are permissible parameters at $q\in D_S$,
 $p\in \Phi_X^{-1}(q)$ and that $u=0$ is a local equation of $E_X$ at $p$. 
Regular parameters $(x,y,z)$ in $\hat{\cal O}_{X,p}$ are called permissible parameters at $p$ for $(u,v)$ if  $u=x^ay^bz^c$ with $a\ge b\ge c\ge 0$.

If $(x,y,z)$ are permissible parameters at $p$ for $(u,v)$, then one of the following
forms holds for $v$ at $p$.

\begin{enumerate}
\item $p$ is a 1 point:
$$
\begin{array}{ll}
u&=x^a\\
v&=P_p(x) + x^bF_p
\end{array}
$$
where $a>0$, $x\not\,\mid F_p$ and $x^bF_p$ has no terms which are powers of $x$,
\item $p$ is a 2 point:
$$
\begin{array}{ll}
u&=(x^ay^b)^m\\
v&=P_p(x^ay^b) + x^cy^dF_p
\end{array}
$$
where $a,b>0$ $(a,b)=1$, $x,y\not\,\mid F_p$ and $x^cy^dF_p$ has no terms which are powers of $x^ay^b$,
\item $p$ is a 3 point:
$$
\begin{array}{ll}
u&=(x^ay^bz^c)^m\\
v&=P_p(x^ay^bz^c) + x^dy^ez^fF_p
\end{array}
$$
where $a,b,c>0$, $(a,b,c)=1$, $x,y,z\not\,\mid F_p$ and $x^dy^ez^fF_p$ has no terms which are powers of $x^ay^bz^c$.
\end{enumerate}
\end{Definition}

We will say that $(x,y,z)$ are permissible parameters at $p$ and that the above expression
of $V$ is the normalized form of $v$ with respect to these parameters.
We will also say that $F_p$ is normalized with respect to $(x,y,z)$.

The leading form of $F_p$ will be denoted by $L_p$.

With the notation of Definition \ref{Def650}, we see that if $p$ is a 1 point, then
\begin{equation}\label{eq1003}
\hat{\cal I}_{\text{sing}(\Phi_X),p}=\sqrt{x^{a+b-1}\left(\frac{\partial F}{\partial y},
\frac{\partial F}{\partial z}\right)}.
\end{equation}
If $p$ is a 2 point, then 
\begin{equation}\label{eq1004}
\hat{\cal I}_{\text{sing}(\Phi_X),p}=\sqrt{x^{ma+c-1}y^{mb+d-1}\left(
(ad-bc)F+ay\frac{\partial F}{\partial y}-bx\frac{\partial F}{\partial x},
y\frac{\partial F}{\partial z}, x\frac{\partial F}{\partial z}\right)}
\end{equation}
If $p$ is a 3 point then 
\begin{equation}\label{eq1005}
\hat{\cal I}_{\text{sing}(\Phi_X),p}=\sqrt{
x^{ma+d-1}y^{mb+e-1}z^{mc+f-1}\left(
\begin{array}{l}
(ae-bd)zF+ayz\frac{\partial F}{\partial y}-bxz\frac{\partial F}{\partial x},\\
(af-cd)yF+ayz\frac{\partial F}{\partial z}-cxy\frac{\partial F}{\partial x},\\
(bf-ce)xF+bxz\frac{\partial F}{\partial z}-cxy\frac{\partial F}{\partial y}
\end{array}\right)}
\end{equation}

\begin{Definition}\label{Def57}
We will say that permissible parameters $(u,v)$ for $\Phi_X(p)\in D_S$ are prepared at $p\in E_X$ if $u=0$ is a local equation of $E_X$ at $p$ and
there exist permissible parameters $(x,y,z)$ at $p$ such that  one of the following
forms hold: 
\begin{equation}\label{eq1060}
\begin{array}{ll}
u&=x^a\\
v&=P(x)+x^by
\end{array}
\end{equation}
or 
\begin{equation}\label{eq1061}
\begin{array}{ll}
u&=(x^ay^b)^m\\
v&=P(x^ay^b)+x^cy^d
\end{array}
\end{equation}
with $ad-bc \ne 0$ 
\begin{equation}\label{eq1062}
\begin{array}{ll}
u&=(x^ay^b)^m\\
v&=P(x^ay^b)+x^cy^dz
\end{array}
\end{equation}
or 
\begin{equation}\label{eq1063}
\begin{array}{ll}
u&=(x^ay^bz^c)^m\\
v&=P(x^ay^bz^c)+x^dy^ez^f
\end{array}
\end{equation}
with
$$
\text{rank }\left(\begin{array}{lll}
a&b&c\\d&e&f\end{array}\right) =2.
$$
We will say that $\Phi_X$ is prepared with respect to $D_S$ if for every $p\in E_X$
there exist permissible
parameters for $\Phi_X(p)$ which are prepared at $p$.
\end{Definition}

\begin{Lemma}\label{Lemma1}
Suppose that $p\in E_X$, $(u,v)$ are permissible parameters at $q=\Phi_X(p)$ such that 
$u=0$ is a local equation of $E_X$.
Then $r=\nu(F_p)$ is independent of permissible parameters $(x,y,z)$ at $p$ for $(u,v)$.

If $p$ is a 1 point then $\nu(F(0,y,z))$ is independent of permissible parameters $(x,y,z)$ at $p$ for $(u,v)$.
If $p$ is a 1 point and 
$$
F_p=\sum_{i+j+k\ge r}a_{ijk}x^iy^jz^k,
$$
then
$$
\tau(F_p)=\text{max}\{j+k\mid a_{ijk}\ne 0\text{ with }i+j+k=r\}
$$
is independent of permissible parameters $(x,y,z)$ for $(u,v)$ at $p$.

If $p$ is a 2 point, then $\nu(F(0,0,z))$ is independent of permissible parameters $(x,y,z)$ at $p$ for $(u,v)$.
 If $p$ is a 2 point and
$$
F_p=\sum_{i+j+k\ge r}a_{ijk}x^iy^jz^k,
$$
then
$$
\tau(F_p)=\text{max}\{k\mid a_{ijk}\ne 0\text{ with }i+j+k=r\}
$$
is independent of permissible parameters $(x,y,z)$ for $(u,v)$ at $p$.
\end{Lemma}

\begin{pf}
Suppose that $(u,v)$ are permissible parameters at $q$
such that $u=0$ is a local equation of $E_X$ at $p$,
 $(x,y,z)$, $(x_1,y_1,z_1)$
are permissible parameters at $p$ for $(u,v)$.
First suppose that $p$ is a 1 point. We have a (normalized) expression
$$
u=x^a, v=P(x)+x^bF
$$
Thus 
$$
\begin{array}{ll}
x&=\omega x_1\\
y&=y(x_1,y_1,z_1)=b_{21}x_1+b_{22}y_1+b_{23}z_1+\cdots\\
z&=z(x_1,y_1,z_1)=b_{31}x_1+b_{32}y_1+b_{33}z_1+\cdots
\end{array}
$$ 
where $\omega^a=1$ and $b_{22}b_{33}-b_{23}b_{32}\ne 0$, and
$$
u=x_1^a, v=P_1(x_1)+x_1^bF_1
$$
where
$$
\begin{array}{ll}
P_1&=P(\omega x_1)+x_1^b\omega^bF(\omega x_1,y(x_1,0,0),z(x_1,0,0))\\
F_1 &= \omega^b[F(\omega x_1,y(x_1,y_1,z_1),z(x_1,y_1,z_1))-F(\omega x_1,y(x_1,0,0),z(x_1,0,0))]
\end{array}
$$
Substituting into 
$$
F_p=\sum_{i+j+k\ge r}a_{ijk}x^iy^jz^k
$$
we get that $\nu(F)=\nu(F_1)$, $\nu(F(0,y,z))=\nu(F_1(0,y_1,z_1))$ so that $F_1$ is normalized
with respect to $(x_1,y_1,z_1)$, and $\tau(F)=\tau(F_1)$.

Now suppose that $p$ is a 2 point. Then
$$
u=(x^ay^b)^m,
v=P(x^ay^b)+x^cy^dF.
$$
Set $r=\nu(F)$.
We have one of the following two cases.

{\bf case 1} 
$$
\begin{array}{ll}
x&=\alpha x_1\\
y&=\beta y_1\\
z&= z(x_1,y_1,z_1)=a_1x_1+a_2y_1+a_3z_1+\cdots
\end{array}
$$
where $\omega=\alpha^a\beta^b$ satisfies $\omega^m=1$ and $a_3\ne 0$, or

{\bf case 2} 
$$
\begin{array}{ll}
x&=\alpha y_1\\
y&=\beta x_1\\
z&= z(x_1,y_1,z_1)=a_1x_1+a_2y_1+a_3z_1+\cdots
\end{array}
$$
where  $\omega=\alpha^a\beta^b$ satisfies $\omega^m=1$, and $a_3\ne 0$.

In case 1, set $t_0=\text{max}\{\frac{c}{a},\frac{d}{b}\}$.
For $t\ge t_0$, set 
\begin{equation}\label{eq621}
\begin{array}{ll}
b_t &= \left(\frac{\partial^{t(a+b)-c-d}(\alpha^c\beta^dF)}{\partial x_1^{ta-c}\partial y_1^{tb-d}}\right)
\vert_{x_1=y_1=z_1=0}\\
F_1 &= \alpha^c\beta^dF-\sum_{t\ge t_0}b_tx_1^{ta-c}y_1^{tb-d}\\
P_1 &= P(\omega x_1^ay_1^b)+\sum_{t\ge t_0} b_t(x_1^ay_1^b)^t.
\end{array}
\end{equation}
Then 
$$
u=(x_1^ay_1^b)^m,
v= P_1(x_1^ay_1^b)+x_1^cy_1^dF_1.
$$
$F_1$ is normalized with respect to $(x_1,y_1,z_1)$ and $\nu(F(0,0,z))=\nu(F_1(0,0,z_1))$.
Set $\alpha(0,0,0) = \alpha_0$. $\beta(0,0,0)=\beta_0$.
Let $L$, $L_1$ be the respective leading forms of $F$ and $F_1$. Then
$$
L_1=\alpha_0\beta_0L(\alpha_0x_1,\beta_0y_1,a_1x_1+a_2y_2+a_3z_1)
$$
if there does not exist natural numbers $i_0,j_0$ such that
$(c+i_0)b-(d+j_0)a=0$ and $i_0+j_0=r$,
$$
L_1=\alpha_0\beta_0L(\alpha_0x_1,\beta_0y_1,a_1x_1+a_2y_2+a_3z_1)-\overline c x_1^{i_0}y_1^{j_0}
$$
for some $\overline c\in k$,
if there  exist natural numbers $i_0,j_0$ such that
$(c+i_0)b-(d+j_0)a=0$ and $i_0+j_0=r$. Thus $\nu(F)=\nu(F_1)$ and $\tau(F)=\tau(F_1)$.

To verify Case 2, we now need only consider the effect of a substitution
$$
x=y_1, y=x_1.
$$
Finally, suppose that $p$ is a 3 point. We have
$$
u=(x^ay^bz^c)^m,
v=P(x^ay^bz^c)+x^dy^ez^fF
$$
There exists $\sigma\in S_3$,  and unit series $\alpha,\beta,\gamma$
with constant terms $\alpha_0$, $\beta_0$, $\gamma_0$ respectively, such that
$$
x=\alpha w_{\sigma(1)}, y=\beta w_{\sigma(2)}, z=\gamma w_{\sigma(3)}
$$
where 
$$
w_1=x_1, w_2=y_1,w_3=z_1,
$$
and if $\omega = \alpha^a\beta^b\gamma^c$, then $\omega^m=1$. Set 
$t_0=\text{max}\{\frac{d}{a},\frac{e}{b},\frac{f}{c}\}$. For $t\ge t_0$, set
\begin{equation}\label{eq632}
\begin{array}{ll}
b_t&=\frac{\partial^{t(a+b+c)-d-e-f}(\alpha^d\beta^e\gamma^{f}F)}
{\partial w_{\sigma(1)}^{ta-d}\partial w_{\sigma(2)}^{tb-e}\partial w_{\sigma(3)}^{tc-f}}
\mid_{w_{\sigma(1)}=w_{\sigma(2)}=w_{\sigma(3)}=0}\\
F_1&=\alpha^d\beta^e\gamma^fF(\alpha w_{\sigma(1)},\beta w_{\sigma(2)},\gamma w_{\sigma(3)})
-\sum_{t\ge t_0}b_tw_{\sigma(1)}^{ta-d}w_{\sigma(2)}^{tb-e}w_{\sigma(3)}^{tc-f}\\
P_1&=P(\omega w_{\sigma(1)}^aw_{\sigma(2)}^bw_{\sigma(3)}^c)
+\sum_{t\ge t_0}b_tw_{\sigma(1)}^{ta}w_{\sigma(2)}^{tb}w_{\sigma(3)}^{tc}.
\end{array}
\end{equation}
Thus
$$
u = (w_{\sigma(1)}^aw_{\sigma(2)}^bw_{\sigma(3)}^c)^m,
v=P_1(w_{\sigma(1)}^aw_{\sigma(2)}^bw_{\sigma(3)}^c)+w_{\sigma(1)}^dw_{\sigma(2)}^ew_{\sigma(3)}^f
F_1(w_{\sigma(1)},w_{\sigma(2)},w_{\sigma(3)})
$$
where the leading form of $F_1$ is 
$$
L_1=\alpha_0^d\beta_0^e\gamma_0^fF(\alpha_0 w_{\sigma(1)},\beta_0 w_{\sigma(2)},\gamma_0 w_{\sigma(3)})
$$
if there does not exist natural numbers $i_0,j_0,k_0$ such that
$$
(d+i_0)b-(e+j_0)a=0,
(d+i_0)c-(f+k_0)a=0,\text{ and }
i_0+j+0+k_0=r,
$$
$$
L_1=\alpha_0^d\beta_0^e\gamma_0^fF(\alpha_0w_{\sigma(1)},\beta_0w_{\sigma(2)},
\gamma_0w_{\sigma(3)})-\overline cx_1^{i_0}y_1^{j_0}z_1^{k_0}
$$
for some $\overline c\in k$, if there exist natural numbers $i_0,j_0,k_0$
such that 
$$
(d+i_0)b-(e+j_0)a=0,
(d+i_0)c-(f+k_0)a=0\text{ and }
i_0+j_0+k_0=r.
$$
Thus $F_1$ is normalized with respect to $(x_1,y_1,z_1)$ and 
$\nu(L_1)=\nu(L)$.
\end{pf}

\begin{Lemma}\label{Lemma300}
Suppose that $p\in E_X$, $q=\Phi_X(p)$. 
Then $r=\nu(F_p)$ is independent of permissible parameters  $(u,v)$ at $q$  such that 
$u=0$ is a local equation of $E_X$ at $p$.

If $p$ is a 1 point then $\nu(F_p(0,y,z))$ is independent of permissible parameters $(u,v)$ at $q$  such that $u=0$ is a local equation of $E_X$ at $p$.
If $p$ is a 1 point, and 
$$
F_p=\sum_{i+j+k\ge r}a_{ijk}x^iy^jz^k,
$$
then 
$$
\tau(F_p)=\text{max}\{j+k\mid \text{ there exists }a_{ijk}\ne 0\text{ with }
i+j+k=r\}
$$
is independent of permissible parameters $(x,y,z)$ at $p$ for $(u,v)$ such that
$u=0$ is a local equation of $E_X$ at $p$.

If $p$ is a 2 point, then $\nu(F_p(0,0,z))$ is independent of permissible parameters $(u,v)$ at $q$  such that $u=0$ is a local equation of $E_X$ at $p$.
 If $p$ is a 2 point, and 
$$
F_p=\sum_{i+j+k\ge r}a_{ijk}x^iy^jz^k,
$$
then 
$$
\tau(F_p)=\text{max}\{k\mid \text{ there exists }a_{ijk}\ne 0\text{ with }
i+j+k=r\}
$$
is independent of permissible parameters $(x,y,z)$ at $p$ for $(u,v)$
such that $u=0$ is a local equation of $E_X$ at $p$.
\end{Lemma}

\begin{pf} Let $m$ be the maximal ideal of $\hat{\cal O}_{X,p}$.
Suppose that $(u,v)$ and $(u_1,v_1)$ are permissible parameters at $q$ such that
$u=0$ is a local equation of $E_X$ at $p$ and $u_1=0$ is a local equation of $E_X$ at $p$.
We will show that the multiplicities of the Lemma are the same for these two sets of
permissible parameters.

{\bf Case 1} Suppose that $p$ is a 1 point. Then $(u_1,v_1)$ and $(u,v)$ are related by a
composition of changes of parameters of the types of Cases 1.1, 1.2 and 1.3 below.
It thus suffices to prove the Lemma in each of these 3 cases.

{\bf Case 1.1} Suppose that $v_1=u$, $u_1=v$. We have
$$
u=x^a, v=P(x)+x^cF
$$
with $r=\nu(F)>0$. In this case we must have $v=\text{unit }u$ in ${\cal O}_{X,p}$. $v=\text{unit } u$ is equivalent to $p=\overline u(x) x^d$ where $\overline u$ is a unit and $0<d\le c$.
Set
$$
x=\overline x(\overline u(x)+x^{c-d}F)^{\frac{-1}{d}}.
$$
Then $v=\overline x^d$. Set $\tau=\frac{-1}{d}$. Write $\overline u(x) = a_0+a_1x+\cdots$.
$$
\begin{array}{ll}
(\overline u(x)+x^{c-d}F)^{\frac{-1}{d}}&=
\overline u(x)^{\tau}+\tau \overline u(x)^{\tau-1}x^{c-d}F+\frac{\tau(\tau-1)}{2}\overline u(x)^{\tau-2}x^{2(c-d)}F^2+\cdots\\
&\equiv \overline u(x)^{\tau}\text{mod }\overline x^{c-d}m^r
\end{array}
$$
We thus have 
\begin{equation}\label{eq300}
x\equiv \overline x \overline u(x)^{\tau}\text{ mod }\overline x^{c-d+1}m^r.
\end{equation}

Now suppose that $P_0(x,\overline x)$ is a series. By substitution of (\ref{eq300}),
 we see that there exist series $A_1$ and $P_1$ such that
$$
P_0(x,\overline x) \equiv A_1(\overline x)+\overline x P_1(x,\overline x)\text{ mod }\overline x^{c-d+1}m^r
$$
By iteration, we get that there is a polynomial $\overline P(\overline x)$, such that $\overline P(0)=P_0(0,0)$,
\begin{equation}\label{eq301}
P_0(x,\overline x) \equiv  \overline P(\overline x)\text{ mod }\overline x^{c-d+1}m^r.
\end{equation}	
We get from (\ref{eq301}) that 
$$
\overline u(x)\equiv Q(\overline x)\text{ mod }\overline x^{c-d+1}m^r
$$
where $Q(0)=\overline u(0)$. Set $u_0=\overline u(0)$.
We also see that 
$$
x\equiv \overline xQ(\overline x)^{\tau}\text{ mod }\overline x^{c-d+1}m^r
$$
Set $\lambda = \frac{-a}{d}$.
$$
\begin{array}{ll}
u&=\overline x^a(\overline u(x)+x^{c-d}F)^{\lambda}\\
&=\overline x^a[\overline u(x)^{\lambda}+\lambda\overline u(x)^{\lambda-1}x^{c-d}F+\frac{\lambda(\lambda-1)}{2}
\overline u^{\lambda-2}x^{2(c-d)}F^2+\cdots]\\
&\equiv \overline x^a[Q(\overline x)^{\lambda}+\lambda Q(\overline x)^{\lambda -1+\tau(c-d)}
\overline x^{c-d}F
+\frac{\lambda(\lambda-1)}{2}Q(\overline x)^{\lambda-2+2\tau (c-d)}\overline x^{2(c-d)}F^2
+\cdots]\\
&\text{ mod }\overline x^{a+c-d+1}m^r
\end{array}
$$
Thus
$$
\begin{array}{l}
v=\overline x^d\\
u=P_1(\overline x)+\overline x^{a+c-d}F_1(\overline x,y,z)
\end{array}
$$
where $\nu(p_1)=a$, and 
$$
F_1\equiv \lambda u_0^{\lambda-1+\tau(c-d)}F(u_0\overline x,y,z)
+\frac{\lambda(\lambda-1)}{2}u_0^{\lambda-2+2\tau(c-d)}\overline x^{c-d}F(u_0\overline x,y,z)^2
+\cdots\text{ mod }\overline xm^r.
$$
Thus $\nu(F(x,y,z))=\nu(F_1(\overline x,y,z))$,
$\nu(F(0,y,z))=\nu(F_1(0,y,z))$ and $\tau(F)=\tau(F_1)$.
\vskip .2truein

{\bf Case 1.2 } Suppose that $u_1=\alpha u$, $v_1=v$, where $\alpha(u,v)$ is a unit series.
We have
$$
u=x^a, v=p(x)+x^bF
$$
with $r=\nu(F)>0$. Set $\lambda=\frac{-1}{a}$. Define
$x =\overline x\alpha^{\lambda}$, so that $u_1=\overline x^a$. Write
$$
\alpha = \alpha_0(u_1)+\alpha_1(u_1)v+\cdots 
$$
\begin{equation}\label{eq615}
\begin{array}{ll}
\alpha^{\lambda}&=\alpha_0(u_1)^{\lambda}+\lambda\alpha_0(u_1)^{\lambda-1}(\alpha_1(u_1)v+\alpha_2(u_1)v^2+\cdots)\\
&+\frac{\lambda(\lambda-1)}{2}\alpha_0(u_1)^{\lambda-2}(\alpha_1(u_1)v+\alpha_2(u_1)v^2+\cdots)^2+\cdots\\
&\equiv \alpha_0(u_1)^{\lambda}+\lambda\alpha_0(u_1)^{\lambda-1}
(\alpha_1(u_1)P(x)+\alpha_2(u_1)P(x)^2+\cdots)\\
&+\frac{\lambda(\lambda-1)}{2}\alpha_0(u_1)^{\lambda-2}
(\alpha_1(u_1)P(x)+\alpha_2(u_1)P(x)^2+\cdots)^2+\cdots\text{ mod }\overline x^bm^r
\end{array}
\end{equation}
We have an expression 
\begin{equation}\label{eq601}
\alpha^{\lambda}\equiv A_0(\overline x)+\overline xB_0(\alpha^{\lambda},\overline x)\text{ mod }\overline x^bm^r.
\end{equation}
Substitute (\ref{eq615}) into (\ref{eq601}) to get 
$$
\alpha^{\lambda}\equiv A_0(\overline x)+\overline x A_1(\overline x)+\overline x^{2}B_1(\alpha^{\lambda},
\overline x)\text{ mod }\overline x^bm^r
$$
By iteration, we get that there is a series $S(\overline x)$ with $S(0)=\alpha_0(0)^{\lambda}=\overline\alpha$, such that
$$
\alpha^{\lambda}\equiv S(\overline x)\text{ mod }\overline x^bm^r.
$$
and
$$
x=\alpha^{\lambda}\overline x\equiv S(\overline x)\overline x\text{ mod }\overline x^{b+1}m^r.
$$
$$
\begin{array}{ll}
v&= P(x)+x^bF\equiv P(\overline xS(\overline x))+
\overline x^bS(\overline x)^bF(\overline xS(\overline x),y,z)]\text{ mod }\overline x^{b+1}m^r
\end{array}
$$
Thus
$$
\begin{array}{l}
u_1=\overline x^a\\
v=P_1(\overline x)+\overline x^bF_1(\overline x,y,z) 
\end{array}
$$
where
$$
\begin{array}{ll}
F_1&\equiv S(\overline x)^bF(\overline x S(\overline x),y,z)\text{ mod }\overline xm^r\\
&\equiv \overline\alpha^bF(\overline \alpha \overline x,y,z)\text{ mod }\overline xm^r
\end{array}
$$
Thus
$\nu(F)=\nu(F_1)$,  $\nu(F(0,y,z))=\nu(F_1(0,y,z))$ and $\tau(F)=\tau(F_1)$.

\vskip .2truein
{\bf Case 1.3}  Suppose that $u_1=u$, $v_1=\alpha u+\beta v$.
Write
$$
\begin{array}{ll}
\alpha &=\sum \alpha_{ij}u^iv^j\\
\beta &= \sum \beta_{ij}u^iv^j
\end{array}
$$
with $\beta_{00}\ne 0$.

We have
$$
u=x^a, v=P(x)+x^bF
$$
 $r=\nu(F)>0$. 

$$
\begin{array}{ll}
v_1 &= \sum \alpha_{ij}u^{i+1}v^j+\sum \beta_{ij}u^iv^{j+1}\\
&=\sum \alpha_{ij}x^{a(i+1)}(P(x)+x^bF)^j+
\sum \beta_{ij}x^{ai}(P(x)+x^bF)^{j+1}\\
&= \sum \alpha_{ij}x^{a(i+1)}(P(x)^j+jx^bP(x)^{j-1}F+
\frac{j(j-1)}{2}x^{2b}P(x)^{j-2}F^2+\cdots+x^{bj}F^j)\\
&+ \sum \beta_{ij}x^{ai}(P(x)^{j+1}+(j+1)x^bP(x)^{j}F+
\frac{(j+1)j}{2}x^{2b}P(x)^{j-1}F^2+\cdots+x^{b(j+1)}F^{j+1})\\
&=Q(x)+H(x)F+x^{2b}F^2G(x,y,z)
\end{array}
$$
where 
$$
H(x) = \sum \alpha_{ij}x^{a(i+1)}jx^bP(x)^{j-1}+\sum\beta_{ij}x^{ai}(j+1)P(x)^jx^b.
$$
We can further write
$$
H(x)=x^b(\beta_{00}+x\Omega(x)).
$$
$$
\begin{array}{ll}
v_1 &= Q(x)+x^b(\beta_{00}+x\Omega(x))F+x^{2b}F^2G\\
&=P_1(x)+x^bF_1
\end{array}
$$
where $P_1(x)=Q(x)$, and 
$$
F_1=(\beta_{00}+x\Omega(x))F+x^bF^2G.
$$
Thus $\nu(F)=\nu(F_1)$, 
$\nu(F(0,y,z))=\nu(F_1(0,y,z))$ and $\tau(F)=\tau(F_1)$.
\vskip .2truein
{\bf Case 2} Suppose that $p$ is a 2 point. It suffices to prove the Lemma in the
three subcases 2.1, 2.2 and 2.3.
\vskip .2truein

{\bf Case 2.1} Suppose that  $u_1=v$, $v_1=u$.
We have an expression
$$
u=(x^ay^b)^k,
v=P(x^ay^b)+x^cy^dF.
$$
Set $r=\nu(F)$. If 
\begin{equation}\label{eq630}
r=0, c\le \text{ord}(P)a\text{ and }d\le \text{ord}(P)b
\end{equation}
then the multiplicities of the Lemma are the same for the two sets of parameters,
so suppose that (\ref{eq630}) doesn't hold.
$v=x^ey^f\text{ unit}$ (for some $e,f$) implies that there exists $t>0$ such that
$$
v=(x^ay^b)^t(\overline u(x^ay^b)+x^{c-at}y^{d-bt}F)
$$
where $\overline u$ is a unit power series. Set $\tau = \frac{-1}{at}$,
$$
x=\overline x(\overline u(x^ay^b)+x^{c-at}y^{d-bt}F)^{\tau}.
$$
$$
\begin{array}{ll}
(\overline u(x^ay^b)+x^{c-at}y^{d-bt}F)^{\tau}&=
\overline u(x^ay^b)^{\tau}+\tau \overline u(x^ay^b)^{\tau-1}x^{c-at}y^{d-bt}F\\
&+\frac{\tau(\tau-1)}{2}\overline u(x^ay^b)^{\tau-2}x^{2(c-at)}y^{2(d-bt)}F^2+\cdots\\
&\equiv \overline u(x^ay^b)^{\tau} \text{ mod }\overline x^{c-at}y^{d-bt}m^r
\end{array}
$$
\begin{equation}\label{eq302}
\begin{array}{ll}
x^ay^b&=\overline x^ay^b(\overline u(x^ay^b)+x^{c-at}y^{d-bt}F)^{a\tau}\\
&\equiv \overline x^ay^b\overline u(x^ay^b)^{a\tau} \text{ mod }\overline x^{c-at+a}y^{d-bt+b}m^r.
\end{array}
\end{equation}

Now suppose that $P_0(x^ay^b,\overline x^ay^b)$ is a series. By substitution of (\ref{eq302}), we see that
$$
P_0(x^ay^b,\overline x^ay^b) \equiv A_1(\overline x^ay^b)+\overline x^ay^b P_1(x^ay^b,\overline x^ay^b)\text{ mod }
\overline x^{c-at+a}y^{d-bt+b}m^r
$$
By iteration, we get that there is a polynomial $Q(\overline x^ay^b)$, such that $u_0=Q(0)=\overline u(0)$, 
\begin{equation}\label{eq303}
\overline u(x^ay^b) \equiv Q(\overline x^ay^b)\text{ mod } \overline x^{c-at+a}y^{d-bt+b}m^r.
\end{equation}	
we get from (\ref{eq303}) that 
$$
\begin{array}{ll}x&\equiv \overline x\overline u(x^ay^b)\text{ mod }\overline x^{c-at+1}y^{d-bt}m^r\\
&\equiv \overline x Q(\overline x^ay^b)\text{ mod }\overline x^{c-at+1}y^{d-bt}m^r
\end{array}
$$
Set $\lambda = \frac{-k}{t}$.
$$
\begin{array}{ll}
u&=(x^ay^b)^k=(\overline x^ay^b)^k[\overline u(x^ay^b)+x^{c-at}y^{d-bt}F]^{\lambda}\\
&=(\overline x^ay^b)^k[\overline u(x^ay^b)^{\lambda}+\lambda\overline u(x^ay^b)^{\lambda-1}x^{c-at}y^{d-bt}F\\
&+\frac{\lambda(\lambda-1)}{2}\overline u(x^ay^b)^{\lambda-2}x^{2(c-at)}y^{2(d-bt)}F^2+\cdots]\\
&\equiv (\overline x^ay^b)^k[Q(\overline x^ay^b)^{\lambda}+
\lambda Q(\overline x^ay^b)^{\lambda-1+c-at}\overline x^{c-at}y^{d-bt}
F(\overline xQ(\overline x^ay^b),y,z)\\
&+\frac{\lambda(\lambda-1)}{2}Q(\overline x^ay^b)^{\lambda-2+2(c-at)}\overline x^{2(c-at)}y^{2(d-bt)}
F(\overline xQ(\overline x^ay^b),y,z)^2\\
&+\cdots]
\text{ mod }\overline x^{ak+c-at+1}y^{bk+d-bt}m^r
\end{array}
$$
Thus 
$$
\begin{array}{l}
v=(\overline x^ay^b)^t\\
u=P_1(\overline x^ay^b)+\overline x^{ak+c-at}y^{bk+d-bt}F_1(\overline x,y,z)
\end{array}
$$
where
$$
\begin{array}{ll}
F_1(\overline x,y,z)&\equiv  
\lambda Q(\overline x^ay^b)^{\lambda-1+c-at}
F(\overline xQ(\overline x^ay^b),y,z)\\
&+\frac{\lambda(\lambda-1)}{2}Q(\overline x^ay^b)^{\lambda-2+2(c-at)}\overline x^{c-at}y^{d-bt}
F(\overline xQ(\overline x^ay^b),y,z)^2\\
&+\cdots
\text{mod }\overline x m^r\\
&\equiv \lambda u_0^{\lambda-1+c-at}F(\overline x u_0,y,z)
+\frac{\lambda(\lambda-1)}{2}u_0^{\lambda-2+2(c-at)}\overline x^{c-at}y^{d-bt}F(\overline x u_0,y,z)^2\\
&+\cdots
\text{ mod }\overline x m^r.
\end{array}
$$
Thus $\nu(F)=\nu(F_1)$,  $\nu(F_1(0,0,z))=\nu(F(0,0,z))$
 and $\tau(F)=\tau(F_1)$.
\vskip .2truein

{\bf Case 2.2} Suppose that $p$ is a 2 point and that $u_1=\alpha u$, $v_1=v$.
We have an expression
$$
u=(x^ay^b)^k,
v=P(x^ay^b)+x^cy^dF
$$
Set $r=\nu(F)$.
Write 
$$
\alpha=\alpha_0(u_1)+\alpha_1(u_1)v+\cdots
$$
Set $\lambda=\frac{-1}{ak}$,
$$
x=\overline x\alpha^{\lambda}.
$$
We have that
$$
u_1=(\overline x^ay^b)^k.
$$
$$
\begin{array}{ll}
\alpha^{\lambda}&= \alpha_0(u_1)^{\lambda}+\lambda\alpha_0(u_1)^{\lambda-1}(\alpha_1(u_1)v+\alpha_2(u_1)v^2+\cdots)\\
&+\frac{\lambda(\lambda-1)}{2}\alpha_0(u_1)^{\lambda-2}(\alpha_1(u_1)v+\alpha_2(u_1)v^2+\cdots)^2+\cdots\\
&\equiv \alpha_0(u_1)^{\lambda}
+\lambda\alpha_0(u_1)^{\lambda-1}(\alpha_1(u_1)P(x^ay^b)\\
&+\alpha_2(u_1)P(x^ay^b)^{2}+\cdots)
+\frac{\lambda(\lambda-1)}{2}\alpha_0(u_1)^{\lambda-2}\\
&(\alpha_1(u_1)P(x^ay^b)+\alpha_2(u_1)P(x^ay^b)^{2}+\cdots)
+\cdots\text{ mod }\overline x^cy^dm^r
\end{array}
$$

Now suppose that $P_0(x^ay^b,\overline x^ay^b)$ is a series. By substitution of the above equation, we see that
$$
P_0(x^ay^b,\overline x^ay^b) \equiv A_1(\overline x^ay^b)+\overline x^ay^b P_1(x^ay^b,\overline x^ay^b)\text{ mod }
\overline x^{c}y^dm^r
$$
By iteration, we get that there is a polynomial $S(\overline x^ay^b)$, such that $u_0=S(0)=\overline u(0)$, 
\begin{equation}\label{eq304}
\alpha^{\lambda} \equiv S(\overline x^ay^b)\text{ mod } \overline x^{c}y^dm^r.
\end{equation}	
we get from (\ref{eq304}) that 
$$
x=\alpha^{\lambda}\overline x\equiv S(\overline x^ay^b)\overline x\text{ mod }\overline x^{c+1}y^dm^r
$$
$$
\begin{array}{ll}
v&=P(x^ay^b)+x^cy^dF\\
&\equiv P(\overline x^ay^bS(\overline x^ay^b)^a)+\overline x^cy^dS(\overline x^ay^b)^c
F(S(\overline x^ay^b)\overline x,y,z))\\
&\text{ mod }\overline x^{c+1}y^dm^r
\end{array}
$$
so that
$$
\begin{array}{l}
u_1=(\overline x^ay^b)^k\\
v=P_1(\overline x^ay^b)+\overline x^cy^dF_1(\overline x,y,z)
\end{array}
$$
where 
$$
\begin{array}{ll}
F_1&\equiv S(\overline x^ay^b)^cF(\overline xS(\overline x^ay^b),y,z)
\text{ mod }\overline xm^r\\
&\equiv u_0^cF(u_0\overline x,y,z)\text{ mod }\overline x m^r
\end{array}
$$
Thus 
$\nu(F)=\nu(F_1)$, $\nu(F(0,0,z))=\nu(F_1(0,0,z))$,   and $\tau(F)=\tau(F_1)$.
\vskip .2truein

{\bf Case 2.3} Suppose that $p$ is a 2 point and that $u_1=u$, $v_1=\alpha u+\beta v$. 
We have an expression
$$
u=(x^ay^b)^k,
v=P(x^ay^b)+x^cy^dF
$$ 
where  $r=\nu(F)$. 
Write
$$
\begin{array}{ll}
\alpha &=\sum \alpha_{ij}u^iv^j\\
\beta &= \sum \beta_{ij}u^iv^j
\end{array}
$$
with $\beta_{00}\ne 0$.
$$
\begin{array}{ll}
v_1 &= \sum \alpha_{ij}u^{i+1}v^j+\sum \beta_{ij}u^iv^{j+1}\\
&=\sum \alpha_{ij}(x^ay^b)^{(i+1)k}(P(x^ay^b)+x^cy^dF)^j+
\sum \beta_{ij}(x^ay^b)^{ik}(P(x^ay^b)+x^cy^dF)^{j+1}\\
&= \sum \alpha_{ij}(x^ay^b)^{(i+1)k}(P(x^ay^b)^j+j(x^cy^d)P(x^ay^b)^{j-1}F
+\cdots+(x^cy^d)^jF^j)\\
&+ \sum \beta_{ij}(x^ay^b)^{ik}(P(x^ay^b)^{j+1}+(j+1)x^cy^dP(x^ay^b)^{j}F
+\cdots+(x^cy^d)^{j+1}F^{j+1})\\
&=Q(x^ay^b)+H(x,y)F+(x^cy^d)^{2}F^2G(x,y,z)
\end{array}
$$
where 
$$
H(x,y) =x^cy^d(\beta_{0,0}+x^ay^b\Omega(x,y)).
$$
Then
$$
\begin{array}{l}
u=(x^ay^b)^k\\
v_1=Q_1(x^ay^b)+x^cy^dF_1
\end{array}
$$
where
$$
F_1\equiv \beta_{0,0}F\text{ mod }(xym^r+(x^cy^d)m^{2r})
$$

Thus $\nu(F)=\nu(F_1)$, 
$\nu(F(0,y,z))=\nu(F_1(0,y,z))$, 
 and $\tau(F)=\tau(F_1)$.
\vskip .2truein
{\bf Case 3} Suppose that $p$ is a 3 point. It suffices to prove the Lemma in the
three subcases 3.1, 3.2 and 3.3.
\vskip .2truein

{\bf Case 3.1} Suppose that  $u_1=v$, $v_1=u_1$.
$$
u=(x^ay^bz^c)^k,
v=P(x^ay^bz^c)+x^dy^ez^fF
$$
where $r=\nu(F)$, If 
\begin{equation}\label{eq622} 
r=0, d\le\text{ord}(P)a, e\le\text{ord}(P)b\text{ and }f\le\text{ord}(P)c
\end{equation}
then the multiplicities of the Lemma are the same for the two sets of parameters,
so suppose that (\ref{eq622}) doesn't hold.

In this case we must have that $v=x^{\alpha}y^{\beta}z^{\gamma}\text{ unit}$, so that
$P(s)=s^t\overline u(s)$ where $\overline u$ is a unit power series and

$$
v=(x^ay^bz^c)^t[\overline u(x^ay^bz^c)+x^{d-ta}y^{e-bt}z^{f-ct}F].
$$
Set $\tau = \frac{-1}{at}$,
$$
x=\overline x(\overline u(x^ay^bz^c)+x^{d-ta}y^{e-bt}z^{f-ct}F)^{\tau}.
$$
$$
\begin{array}{ll}
(\overline u(x^ay^bz^c)+x^{d-ta}y^{e-bt}z^{f-ct}F)^{\tau}&=
\overline u(x^ay^bz^c)^{\tau}+\tau\overline u(x^ay^bz^c)^{\tau-1}x^{d-ta}y^{e-bt}z^{f-ct}F\\
&+\frac{\tau(\tau-1)}{2}\overline u(x^ay^bz^c)^{\tau-2}x^{2(d-ta)}y^{2(e-bt)}z^{2(f-ct)}F^2+\cdots\\
&\equiv \overline u(x^ay^bz^c)^{\tau}\text{ mod }\overline x^{d-ta}y^{e-bt}z^{f-ct}m^r
\end{array}
$$
$$
\begin{array}{ll}
x^ay^bz^c&=\overline x^ay^bz^c(\overline u(x^ay^bz^c)+x^{d-ta}y^{e-bt}z^{f-ct}F)^{a\tau}\\
&\equiv \overline x^ay^bz^c\overline u(x^ay^bz^c)^{a\tau}
\text{ mod }\overline x^{a+d-ta}y^{e-bt+b}z^{f-ct+c}m^r
\end{array}
$$

Now suppose that $P_0(x^ay^bz^c,\overline x^ay^bz^c)$ is a series. By substitution of the above equation, we see that
$$
P_0(x^ay^bz^c,\overline x^ay^bz^c) \equiv A_1(\overline x^ay^bz^c)+\overline x^ay^bz^c P_1(x^ay^bz^c,\overline x^ay^bz^c)\text{ mod }
\overline x^{a+d-ta}y^{e-bt+b}z^{f-ct+c}m^r
$$
By iteration, we get that there is a polynomial $Q(\overline x^ay^bz^c)$, such that if
$u_0=\overline u(0)$, $Q(0)=\overline u(0)=u_0$,
$$
\overline u(x^ay^bz^c)\equiv Q(\overline x^ay^bz^c)\text{ mod }\overline x^{a+d-ta}y^{e-bt+b}z^{f-ct+c}m^r
$$	
Thus
$$
\begin{array}{ll}
x&\equiv \overline x\overline u(x^ay^bz^c)^{\tau}\text{ mod }\overline x^{d-ta+1}y^{e-bt}z^{f-ct}m^r\\
&\equiv \overline xQ(\overline x^ay^bz^c)^{\tau}\text{ mod }\overline x^{d-ta+1}y^{e-bt}z^{f-ct}m^r
\end{array}
$$
Set $\lambda=\frac{-k}{t}$.
$$
\begin{array}{ll}
u&=(x^ay^bz^c)^k=(\overline x^ay^bz^c)^k(\overline u(x^ay^bz^c)+x^{d-ta}y^{e-bt}z^{f-ct}F)^{\lambda}\\
&=(\overline x^ay^bz^c)^k[\overline u(x^ay^bz^c)^{\lambda}+\lambda\overline u(x^ay^bz^c)^{\lambda-1}
x^{d-ta}y^{e-bt}z^{f-ct}F\\
&+\frac{\lambda(\lambda-1)}{2}\overline u(x^ay^bz^c)^{\lambda-2}x^{2(d-ta)}y^{2(e-bt)}z^{2(f-ct)}F^2+\cdots]\\
&\equiv (\overline x^ay^bz^c)^k[Q(\overline x^ay^bz^c)^a)^{\lambda}\\
&+\lambda Q(\overline x^ay^bz^c)^{\lambda-1+\tau(d-ta)}\overline x^{d-ta}y^{e-bt}z^{f-ct}
F(\overline xQ(\overline x^ay^bz^c)^{\tau},y,z)\\
&+\frac{\lambda(\lambda-1)}{2}Q(\overline x^ay^bz^c)^{\lambda-2+2\tau(d-ta)}
\overline x^{2(d-ta)}y^{2(e-bt)}z^{2(f-ct)}F(\overline xQ(\overline x^ay^bz^c)^{\tau},y,z)^2\\
&+\cdots]\text{ mod }\overline x^{ak+d-ta+1}y^{bk+e-bt}z^{ck+f-ct}m^r
\end{array}
$$
$$
\begin{array}{l}
v=(\overline x^ay^bz^c)^t\\
u=P_1(\overline x^ay^bz^c)+\overline x^{ak+d-ta}y^{bk+e-bt}z^{ck+f-ct}F_1(\overline x,y,z)
\end{array}
$$
where
$$
\begin{array}{ll}
F_1&\equiv \lambda Q(\overline x^ay^bz^c)^{\lambda-1+\tau(d-ta)}
F(\overline xQ(\overline x^ay^bz^c)^{\tau},y,z)\\
&+\frac{\lambda(\lambda-1)}{2}Q(\overline x^ay^bz^c)^a)^{\lambda-2+2\tau(d-ta)}
\overline x^{d-ta}y^{e-bt}z^{f-ct}
F(\overline xQ(\overline x^ay^bz^c)^{\tau},y,z)^2+\cdots
\text{ mod }\overline xm^r\\
&\equiv \lambda u_0^{\lambda-1+\tau(d-ta)}F(u_0^{\tau}\overline x,y,z)\\
&+\frac{\lambda(\lambda-1)}{2}u_0^{\lambda-2+\tau(d-ta)}\overline x^{d-ta}y^{e-bt}z^{f-ct}F(u_0\overline x,y,z)^2+\cdots
\text{ mod }\overline xm^r
\end{array}
$$
Thus $\nu(F_1)=\nu(F)$.
\vskip .2truein
{\bf Case 3.2} 
 Suppose that  $u_1=\alpha u$, $v_1=v$.
$$
u=(x^ay^bz^c)^k,
v=P(x^ay^bz^c)+x^dy^ez^fF
$$
where $r=\nu(F)$. Set $\lambda = \frac{-1}{ak}$, $x=\overline x\alpha^{\lambda}$.
Thus
$$
u_1=(\overline x^ay^bz^c)^k.
$$
Write
$$
\alpha=\alpha_0(u_1)+\alpha_1(u_1)v+\cdots
$$

$$
\begin{array}{ll}
\alpha^{\lambda}&\equiv \alpha_0(u_1)^{\lambda}+\lambda\alpha_0(u_1)^{\lambda-1}(\alpha_1(u_1)v+\alpha_2(u_1)v^2+\cdots)\\
&+\frac{\lambda(\lambda-1)}{2}\alpha_0(u_1)^{\lambda-2}(\alpha_1(u_1)v+\alpha_2(u_1)v^2+\cdots)^2+\cdots\\
&\equiv \alpha_0(u_1)^{\lambda}+\lambda\alpha_0(u_1)^{\lambda-1}(\alpha_1(u_1)P(x^ay^bz^c)+\alpha_2(u_1)P(x^ay^bz^c)^2
+\cdots)\\
&+\frac{\lambda(\lambda-1)}{2}\alpha_0(u_1)^{\lambda-2}(\alpha_1(u_1)P(x^ay^bz^c)+\alpha_2(u_1)P(x^ay^bz^c)^2
+\cdots)^2+\cdots\text{ mod }\overline x^dy^ez^fm^r
\end{array}
$$
Thus
$$
\alpha^{\lambda} \equiv A_0(\overline x^ay^bz^c)+\overline x^ay^bz^c
 B_0(\alpha^{\lambda},\overline x^ay^bz^c)\text{ mod }
\overline x^{d}y^ez^fm^r
$$
Substitute the above equation into itself and iterate to get
$$
\alpha^{\lambda}\equiv S(\overline x^ay^bz^c)\text{ mod }\overline x^dy^ez^fm^r
$$
Set $\overline \alpha=\alpha(0)$. Then $S(0)=\overline \alpha^{\lambda}$.
$$
x =\alpha^{\lambda}\overline x\equiv S(\overline x^ay^bz^c)\overline x\text{ mod }\overline x^{d+1}y^ez^fm^r
$$
$$
\begin{array}{ll}
v&=P(x^ay^bz^c)+x^dy^ez^fF\\
&\equiv P(\overline x^ay^bz^cS(\overline x^ay^bz^c)^a)+\overline x^dy^ez^fS(\overline x^ay^bz^c)^dF(S(\overline x^ay^bz^c)\overline x,y,z)
\text{ mod }\overline x^{d+1}y^ez^fm^r
\end{array}
$$
Thus
$$
\begin{array}{l}
u_1=(\overline x^ay^bz^c)^k\\
v=P_1(\overline x^ay^bz^c)+\overline x^dy^ez^fF_1(\overline x,y,z)
\end{array}
$$
where 
$$
\begin{array}{ll}
F_1&\equiv S(\overline x^ay^bz^c)^dF(S(\overline x^ay^bz^c)\overline x,y,z)\text{ mod }\overline xm^r\\
&\equiv \overline \alpha^{\lambda d}F(\overline\alpha^{\lambda}\overline x,y,z)\text{ mod }\overline xm^r
\end{array}
$$
Thus $\nu(F_1)=\nu(F)$.
\vskip .2truein
{\bf Case 3.3} 
Suppose that $u_1=u$, $v_1=\alpha u+\beta v$. 
We have an expression
$$
u=(x^ay^bz^c)^k,
v=P(x^ay^bz^c)+x^dy^ez^fF
$$
where $r=\nu(F)$.

Write
$$
\begin{array}{ll}
\alpha &=\sum \alpha_{ij}u^iv^j\\
\beta &= \sum \beta_{ij}u^iv^j
\end{array}
$$
with $\beta_{00}\ne 0$.
$$
\begin{array}{ll}
v_1 &= \sum \alpha_{ij}u^{i+1}v^j+\sum \beta_{ij}u^iv^{j+1}\\
&=\sum \alpha_{ij}(x^ay^bz^c)^{(i+1)k}(P(x^ay^bz^c)+x^dy^ez^fF)^j
+\sum \beta_{ij}(x^ay^bz^c)^{ik}(P(x^ay^bz^c)+x^dy^ez^fF)^{j+1}\\
&= \sum \alpha_{ij}(x^ay^bz^c)^{(i+1)k}(P(x^ay^bz^c)^j+j(x^dy^ez^f)P(x^ay^bz^c)^{j-1}F
+\cdots+(x^dy^ez^f)^jF^j)\\
&+ \sum \beta_{ij}(x^ay^bz^c)^{ik}(P(x^ay^bz^c)^{j+1}+(j+1)x^dy^ez^fP(x^ay^bz^c)^{j}F
+\cdots+(x^dy^ez^f)^{j+1}F^{j+1})\\
&=Q(x^ay^bz^c)+HF+(x^dy^ez^f)^2F^2G
\end{array}
$$
where 
$$
H =x^dy^ez^f(\beta_{0,0}+x^ay^bz^c\Omega)
$$
Then
$$
\begin{array}{l}
u=(x^ay^bz^c)^k\\
v_1=Q_1(x^ay^bz^c)+x^dy^ez^fF_1
\end{array}
$$
where
$$
F_1\equiv \beta_{0,0}F\text{ mod }(xyz m^r+x^dy^ez^fm^{2r})
$$

Thus $\nu(F)=\nu(F_1)$.
\end{pf}

By Lemmas \ref{Lemma1} and \ref{Lemma300}, we can make the following definitions, with
the notation of Definition \ref{Def650}.

\begin{Definition}\label{Def1086} Suppose that $p\in E_X$, $(u,v)$ are permissible parameters
at $\Phi_X(p)$, $(x,y,z)$ are permissible parameters
at $p$ for $(u,v)$ such that $u=0$ is a local equation of $E_X$ at $p$. Thus one of the forms of Definition \ref{Def650} holds. 
Define $\nu(p)= \nu(F_p)$. 
If $p$ is a 1 point, define $\gamma(p)=\text{mult}(F_p(0,y,z))$.
If $p$ is a 2 point, define $\gamma(p)=\text{mult}(F_p(0,0,z))$.

Suppose that $p\in X$ is a 1 point such that
$$
\begin{array}{ll}
u&=x^a\\
v&=P(x^a)+x^bF_p\\
F_p&=\sum_{i+j+k\ge r} a_{ijk}x^iy^jz^k
\end{array}
$$
where $\nu(p)=r$. Define
$$
\tau(p)=\text{max}\{j+k\mid \text{ there exits }a_{ijk}\ne 0\text{ with }
i+j+k=r\}.
$$
If $p$ is a 1 point, we have $1\le\tau(p)\le \nu(p)$.
Suppose that $p\in X$ is a 2 point such that
$$
\begin{array}{ll}
u&=(x^ay^b)^m\\
v&=P(x^ay^b)+x^cy^dF_p\\
F_p&=\sum_{i+j+k\ge r} a_{ijk}x^iy^jz^k
\end{array}
$$
where $\nu(p)=r$. Define
$$
\tau(p)=\text{max}\{k\mid \text{ there exits }a_{ijk}\ne 0\text{ with }
i+j+k=r\}.
$$

\end{Definition}

Define 
$$
S_r(X) = \{p\in E_X | \nu(p)\ge r\}.
$$
Let $\overline S_r(X)$ be the Zariski closure of $S_r(X)$ in $X$.
\begin{Definition}\label{Def1064}
A point $p\in E_X$ is resolved if the following condition holds.
\begin{enumerate}
\item If $p$ is a 1 point then $\nu(p)\le 1$.
\item If $p$ is a 2 point then $\gamma(p)\le 1$.
\item If $p$ is a 3 point then $\nu(p)=0$.
\end{enumerate}
\end{Definition}

\begin{Remark}
If $p\in E_X$ is resolved and $(u,v)$ are permissible parameters at $\Phi_X(p)$
such that $u=0$ is a local equation of $E_X$ at $p$, then $(u,v)$ are prepared at $p$.
\end{Remark}

\begin{Lemma}\label{Lemma657}
$S_r(X)\subset \text{sing}(\Phi_X)$ for $r\ge 2$, and all  3 points are contained
in $\text{sing}(\Phi_X)$. If $p\in S_1(X)$ is a 2 point then $p\in \text{sing}(\Phi_X)$.
\end{Lemma}
\begin{pf}
The Lemma is immediate from (\ref{eq1003}), (\ref{eq1004}) and (\ref{eq1005}).
\end{pf}

\begin{Example}\label{Ex1} 
$S_r(X)$ is in general not Zariski closed. Consider the 2 point $p$ with local equations
$$
\begin{array}{ll}
u&=xy\\
v&=x^2y.
\end{array}
$$
$\nu(p)=0$. At 1 points $q$ on the surface $x=0$ there are regular parameters $(x,y_1,z)$
with $y=y_1+\alpha$ for some $0\ne \alpha\in k$. Set $\overline x=x(y_1+\alpha)$.
There are permissible parameters $(\overline x,\overline y,z)$ at $q$ such that
$$
\begin{array}{ll}
u&=\overline x\\
v&=\alpha^{-1}\overline x^2+\overline x^2\overline y.
\end{array}
$$
Thus $\nu(q)=1$.
\end{Example}

\begin{Lemma}\label{Lemma301}
 Suppose that $p\in E_X$ is a 1 point and that $I\subset \hat{\cal O}_{X,p}$ is a reduced ideal
such that if $x=0$ is a local equation of $E_X$ at $p$ then $x\in I$. Then the condition 
$F_p\in I^s$ (with $s\in{\bold N}$) and the condition $F_p\in m_pI^s$ (with $s\in{\bold N}$) are independent of the choice of permissible parameters $(u,v)$ at $\Phi_X(p)$ such that $u=0$ is a local equation
of $E_X$ at $p$, and
permissible parameters $(x,y,z)$ for $(u,v)$ at $p$.
\end{Lemma}

\begin{pf} 
If $I=m_p\hat{\cal O}_{X,p}$, the Lemma follows from Lemmas \ref{Lemma1}
and \ref{Lemma300}. So we assume that $I=(x, f)$ for some series $f(y,z)$.

If $(x,y,z)$ and $(x_1,y_1,z_1)$
 are permissible parameters at $p$ for $(u,v)$ then with the
notation of the proof of Lemma \ref{Lemma1},
$$
F_1=\omega^b[F-F(\omega x_1,y(x_1,0,0),z(x_1,0,0))].
$$
and 
$$
x^{\nu(p)}\mid F(\omega x_1,y(x_1(0,0),z(x_1,0,0)]
$$
implies $F_1\in I^s$ (or $F_1\in mI^s$), $x\in I$ and $s\le\nu(p)$ (or $s\le \nu(p)-1$). 

Now suppose that $(u,v)$, $(u_1,v_1)$ are permissible parameters at $f(p)$.
Suppose that $v_1=u$, $u_1=v$. With the notation of Case 1.1 of the proof of
Lemma \ref{Lemma300}, $F\in I^s$ implies (\ref{eq300}) can be  modified to
$$
x\equiv \overline x \overline u(x)^{\tau}\text{ mod }\overline x^{c-d+1}I^s
$$
and thus
$$
x\equiv \overline x Q(\overline x)^{\tau}\text{ mod }\overline x^{c-d+1}I^s
$$
We thus have
$$
F_1\equiv \lambda u_0^{\lambda-1+\tau(c-d)}F(u_0\overline x,y,z)
+\frac{\lambda(\lambda-1)}{2}u_0^{\lambda-2+2\tau(c-d)}
\overline x^{c-d}F(u_0\overline x,y,z)^2+\cdots
\text{ mod }\overline x I^s
$$
since $I=(x,f(y,z))$ for some $f$, we have $F_1\in I^s$. We have a similar proof when
$F\in mI^s$. We can replace $m^r$ in the formulas of case 1.1 of Lemma
\ref{Lemma300} with $mI^s$.

In the proofs of cases 1.2 and 1.3, we can also replace $m^r$ in all the formulas with $I^s$
(or $mI^s$). Again, since $x\in I$, we get $F_1\in I^s$ (or
$F_1\in mI^s$). 
\end{pf}

\begin{Lemma}\label{Lemma600} Suppose that $C$ is a 2 curve and $p\in C$.
Then the condition $F_p\in \hat{\cal I}_{C,p}^s$, (with $s\in{\bold N}$)
 is independent of permissible parameters
at $\Phi_X(p)$ and $p$.
\end{Lemma}

\begin{pf}
Suppose that $(u,v)$ are permissible parameters at $\Phi_X(p)$. We will first show that the 
condition is independent of permissible parameters for $(u,v)$ at $p$.

If $p$ is a 2 point, this follows from the proof of Lemma \ref{Lemma1}, with the
observation that, in the notation of (\ref{eq621}), $F\in\hat{\cal I}_{C,p}^s$
implies 
$$
\frac{\partial^{t(a+b)-c-d}(\alpha^c\beta^dF)}{\partial x_1^{ta-c}\partial y_1^{tb-d}}
\in\hat{\cal I}_{C,p}^{s-t(a+b)+c+d},
$$
 so that
$\sum b_tx_1^{ta-c}y_1^{tb-d}\in\hat{\cal I}_{C,p}^s$, and thus $F_1\in\hat{\cal I}_{C,p}^s$.

If $p$ is a 3 point, this also follows from the proof of Lemma \ref{Lemma1}. With the
notation of (\ref{eq632}), after possibly permuting the parameters
$(w_{\sigma(1)},w_{\sigma(2)},w_{\sigma(3)})$, we have $\hat{\cal I}_{C,p}=(w_{\sigma(1)},w_{\sigma(2)})$. 

If $G\in \hat{\cal O}_{X,p}$ is a series and $G\in\hat{\cal I}_{C,p}^a$ for some $a$,  we have that
$$
\frac{\partial G}{\partial w_{\sigma(1)}}, 
\frac{\partial G}{\partial w_{\sigma(2)}}
\in \hat{\cal I}_{C,p}^{a-1}
$$
and 
$$
\frac{\partial G}{\partial w_{\sigma(3)}}\in \hat{\cal I}_{C,p}^{a}
$$
Thus 
$$
b_tw_{\sigma(1)}^{ta-d}w_{\sigma(2)}^{tb-e}w_{\sigma(3)}^{tc-f}\in \hat{\cal I}_{C,p}^s
$$
for all $t$, and $F_1\in \hat{\cal I}_{C,p}^s$.

The independence of the conditions from permissible parameters $(u,v)$ at $\Phi_X(p)$
follows from cases 2.1 - 3.3 of Lemma \ref{Lemma300}, with $m^r$ replaced by
$\hat{\cal I}_{C,p}^s=(x,y)^s$ in the formulas of these cases.
\end{pf}

\begin{Example}\label{Example601} If $p$ is a 2 point, the condition
$F_p\in I^s$ where $I\subset\hat{\cal O}_{X,p}$ is a reduced ideal can depend on the
choice of permissible parameters at $p$.
\end{Example}

\begin{pf} Consider
$$
u=xy, v=z^2+xz
$$
the Jacobian is $J=(xz,y(2z+x),x(2z+x))$.
$$
x^2=2xz+x^2-2xz\in J.
$$
 $\sqrt{J}=(x,yz)$. $(x,y,z)$ are permissible parameters for $(u,v)$ at $p$. Let $I=(x,z)$.
$F\in I^2$.

We have other permissible parameters $(x,y,\overline z)$ at $p$,
where $\overline z=z-y$. Then $I=(\overline z+y,x)$. The normalized form
of $v$ with respect to these new parameters is
$$
u=xy,
v=xy+F
$$
where 
$$
F=[(\overline z+y)^2+x\overline z]\not\in I^2.
$$
\end{pf}
\begin{Lemma}\label{Lemma651}
Suppose that $p$ is a 2 point, and $C$ is a curve, making SNCs with the 2 curve through $p$.
Then the condition $F_p\in\hat{\cal I}_{C,p}^s$ with $s\in{\bold N}$ is independent of permissible parameters $(u,v)$ at $\Phi_X(p)$ and permissible parameters
$(x,y,z)$ at $p$ for $(u,v)$ such that $\hat{\cal I}_{C,p}=(x,z)$. 
\end{Lemma}

We will call parameters as in Lemma \ref{Lemma651} permissible parameters for $C$ at $p$.

\begin{pf} Suppose that $(u,v)$ are permissible parameters at $\Phi_X(p)$. We will first
show that this is independent of such permissible parameters at   $p$
for $(u,v)$.
Suppose that $(x,y,z)$ and $(x_1,y_1,z_1)$ are permissible parameters for $(u,v)$ at $p$
such that $\hat{\cal I}_{C,p}=(x,z)=(x_1,z_1)$ and
$$
\begin{array}{ll}
u&=(x^ay^b)^m\\
v&=P(x^ay^b)+x^cy^dF
\end{array}
$$
with $F\in\hat{\cal I}_{C,p}^s=(x,y)^s$. We have
$$
x=\alpha x_1, y=\beta y_1,
z=z(x_1,y_1,z_1)=\omega z_1+\gamma x_1
$$
where $\alpha,\beta,\omega$ are units in $\hat{\cal O}_{X,p}$ and $\gamma\in\hat{\cal O}_{X,p}$.
If $G\in\hat{\cal O}_{X,p}$ is such that $G\in (x_1,z_1)^a$ then
$$
\frac{\partial G}{\partial x_1}\in (x_1,z_1)^{a-1}
$$
and 
$$
\frac{\partial G}{\partial y_1}\in (x_1,z_1)^a.
$$
Thus
$$
\frac{\partial^{t(a+b)-c-d}(\alpha^c\beta^dF)}{\partial x_1^{ta-c}\partial y_1^{tb-d}}
\in (x_1,z_1)^{s-(ta-c)}
$$
In (\ref{eq621}) of Lemma \ref{Lemma1}, we have $b_t=0$ if $s>(ta-c)$, so that
$F_1\in (x_1,z_1)^s$.

The independence of the condition $F\in \hat{\cal I}_{C,p}^s$ from choice of
permissible parameters $(u,v)$ at $\Phi_X(p)$ follows from cases 2.1-2.3 of
Lemma \ref{Lemma300}, with $m^r$ replaced by $\hat{\cal I}_{C,p}^s=(x,z)^s$
is the formulas of these cases.
\end{pf}

Let $B_2(X)$ be the (possibly not closed) curve of 2 points in $X$, $B_3(X)=\{p_1,\ldots,p_r\}$ the set of 3 points in $X$. Let $\overline B_2(X)=B_2(X)\cup B_3(X)$ be the Zariski closure of $B_2(X)$ in $X$.

\begin{Definition}\label{Def1087} Suppose that $Z\subset E_X$ is a reduced closed subscheme of
dimension $\le 1$ and $p\in E_X$. We will  say that $Z$ makes SNCs with $\overline B_2(X)$
at $p$ if
\begin{enumerate}
\item All components of $Z$ are nonsingular at $p$.
\item If $C_1,\ldots, C_s$  are the curves of $Z$ containing $p$ and
$D_1,\ldots, D_t$ are the components of $\overline B_2(X)$ containing $p$, then 
$C_1,\ldots, C_s, D_1,\ldots, D_t$ have independent tangent directions at $p$.
\end{enumerate}
\end{Definition}

We will say that $Z$ makes SNCs with $\overline B_2(X)$ if $Z$ makes SNCs with $\overline B_2(X)$ at $p$ for all $p\in E_X$.

\begin{Definition}\label{Def1089} Suppose that $p\in X$, $U$ is an affine
neighborhood of $p$ in $X$, and $\sigma:V\rightarrow U$ is an \'etale cover. Then 
we will say that $V$ is an \'etale neighborhood of $p$. Suppose that $D\subset X$. We will
write $D\cap V$ to denote $\sigma^{-1}(D\cap U)$.
\end{Definition}

\begin{Definition}\label{Def1088}(c.f. Chapter 3, Section 6 \cite{Mum}.) Suppose that $V$ is an affine $k$-variety. 
$x_1,\ldots, x_n\in\Gamma(V,{\cal O}_V)$ are uniformizing parameters on $V$ if the 
natural morphism $V\rightarrow\text{spec}(k[x_1,\ldots,x_n])$ is \'etale.
\end{Definition}

\begin{Lemma}\label{Lemma1006} Suppose that $(x,y,z)$ are permissible
parameters at $p$ for $(u,v)$ such that $y,z\in{\cal O}_{X,p}$. Then there exists an
affine neighborhood $U$ of $p$ and an \'etale cover $V$ of $U$ such that $(x,y,z)$ are
uniformizing parameters on $V$.
\end{Lemma}

\begin{pf} With the notations of Definition \ref{Def650}, let $\overline a=a$
if $p$ is a 1 point, $\overline a=ma$ if $p$ is a 2 or 3 point. There exists a unit
$\lambda\in{\cal O}_{X,p}$ and $\tilde x\in{\cal O}_{X,p}$ such that $x^{\overline a}
=\lambda\tilde x^{\overline a}$. There exists an affine neighborhood $U_1$ of $p$ such that
$\tilde x,y,z,\lambda\in R=\Gamma(U_1,{\cal O}_{X})$ and $\lambda$ is a unit in $R$. 
Set $S=R[\lambda^{\frac{1}{\overline a}}]$, $V_1=\text{spec}(S)$. $f:V_1\rightarrow U_1$
is an \'etale cover. $k[x,y,z]\rightarrow S$ defines a morphism $g:V_1\rightarrow {\bold A}^3$.
Let $a$ be the origin of ${\bold A}^3$. $q\in g^{-1}(a)$ if and only if $x,y,z\in m_q$ which
holds if and only if $\tilde x,y,z\in m_q$. Thus $g^{-1}(a)=f^{-1}(p)$.
$\hat{\cal O}_{V_1,q}=k[[\tilde x,y,z]]=k[[x,y,z]]$ for all $q\in g^{-1}(a)$.
Thus $g$ is \'etale at all points of $g^{-1}(a)$. Since this is an open condition,
(Proposition 4.5 \cite{SGA}) there
exists a closed set $Z_1$ of $V_1$ which is disjoint from $f^{-1}(p)$ such that
$g\mid (V_1-Z_1)$ is \'etale. Let $U$ be an affine neighborhood of $p$ in $U_1$ which 
is disjoint from the closed set $f(Z_1)$. Let $V=f^{-1}(U)$. Then $V$ is an \'etale cover of $U$
on which $x,y,z$ are uniformizing parameters.
\end{pf}

\begin{Proposition}\label{Prop1} 
$S_r(X)\cap(X-\overline B_2(X))$ is Zariski closed in $X-\overline B_2(X)$ and $S_r(X)\cap B_2(X)$ is Zariski closed in $B_2(X)$. Thus
$S_r(X)$ is a constructible set. 
\end{Proposition}
\begin{pf}
First suppose that $p$ is a 1 point. Then there are 
regular parameters $\tilde x, y, z$ in ${\cal O}_{X,p}$,
permissible parameters $x, y, z$ at $p$, and a unit $\lambda\in{\cal O}_{X,p}$ such that 
$$
\begin{array}{ll}
u&=x^a=\lambda \tilde x^a\\
v&=P(x)+x^bF_p(x,y,z).
\end{array}
$$
 $\tilde x, y, z$ are uniformizing parameters in an affine
neighborhood $U$ of $p$, and there exists an \'etale neighborhood $\sigma:V=\text{spec}
(S)\rightarrow U$ of $p$ such that $(x,y,z)$ are uniformizing parameters on $V$, $x=0$
is a local equation of $E_X\cap V$ in $V$. Let
$$
I=(\frac{\partial^{i+j+k}v}{\partial x^i\partial y^j\partial z^k}\mid j+k>0, i+j+k\le b+r-1)\subset S,
$$
$Z=V(I)\subset V$.

Suppose that $p'\in E_X\cap V$. Then if $\alpha = y(p'), \beta = z(p')$, we have that 
$$
\begin{array}{ll}
u&=x^a\\
v &= \sum  \frac{1}{i!j!k!}\frac{\partial^{i+j+k}v}{\partial x^i\partial y^j\partial z^k}(0,\alpha,\beta)
x^i(y-\alpha)^j(z-\beta)^k
\end{array}
$$
and
$$
v-v(\sigma(p')) = P_{p'}(x)+x^bF_{p'}.
$$
$\nu(p')\ge r$ if and only if $p'\in V(I)$.
Let $Z_1=\sigma(Z)$. $S_r(X)\cap U=Z_1$ is closed in $U$.

Now suppose that $p$ is a 2 point. Then there are regular parameters $\tilde x, y, z$ in ${\cal O}_{X,p}$
and permissible parameters $x, y, z$ at $p$ and a unit $\lambda$ in ${\cal O}_{X,p}$
such that  
$$
\begin{array}{ll}
u&=(x^ay^b)^m=\lambda(\tilde x^ay^b)^m\\
v&=P(x^ay^b)+x^cy^dF_p(x,y,z)
\end{array}
$$

There exists an \'etale neighborhood $\sigma:V=\text{spec}(S)\rightarrow U$ of $p$ such that 
$(x,y,z)$ are uniformizing parameters on $V$, $xy=0$ is a local equation of $E_X\cap V$ in $V$.
Let $C$ be the 2 curve in $X$ containing $p$. Suppose that $p'\in C\cap V$. Then if
$\beta=z(p')$, we have that
$$
\begin{array}{ll}
u&=(x^ay^b)^m\\
v-v(\sigma(p'))&=P_{p'}(x^ay^b)+x^cy^dF_{p'}.
\end{array}
$$
$$
v=\sum \frac{1}{i!j!k!}\frac{\partial^{i+j+k}v}{\partial x^i\partial y^j\partial z^k}(0,0,\beta)
x^iy^j(z-\beta)^k.
$$
Let 
$$
I=\left(\frac{\partial^{i+j+k}v}{\partial x^i\partial y^j\partial z^k}\mid
k>0\text{ or }k=0\text{ and }a(d+j)-b(c+i)=0\text{ and }i+j+k\le c+d+r-1\right)
\subset S,
$$
$p'\in C$ and $\nu(p')\ge r$ if and only if $p'\in V(I)\cap C$.
 $Z=V(I)\subset V$. Let $Z_1=\sigma(Z)$. 
$S_r(X)\cap C\cap U= C\cap Z_1\cap U$ is closed in $C\cap U$.
\end{pf}

\begin{Lemma}\label{Lemma302} Suppose that $p\in E_X$ is a 1 point or a 2 point.
\begin{enumerate}
\item Suppose that $(x,y,z)$ are permissible parameters at $p$, $I\subset
\hat{\cal O}_{x,p}$ is a reduced ideal and 
$F_p\in I^r$ for some $r\ge 2$. Then $\hat{\cal I}_{\overline{S_r(X)},p}\subset I$.
\item Suppose that $(x,y,z)$ are permissible parameters at $p$, $I=(x,f(y,z))\subset
\hat{\cal O}_{x,p}$ is a reduced ideal and $F_p\in(x)+I^2$. Then $\hat{\cal I}_{\overline S_2(X),p}\subset I$.
\end{enumerate}
\end{Lemma}
\begin{pf} Suppose that $F_p\in I^r$ for some $r\ge 2$. First assume that $p$ is a 
1 point. Since $x\in I$ and $r\le\nu(p)$, we can make a permissible change of parameters, and renormalize
to get that $y,z\in {\cal O}_{X,p}$ and $F_p\in I^r$.
\begin{equation}\label{eq968}
\hat{\cal I}_{\overline{S_r(X)},p}=
\sqrt{\left(\frac{\partial^{i+j+k}F_p}{\partial x^i\partial y^j\partial z^k}
\vert i+j+k\le r-1, j+k>0\right)}.
\end{equation}
$F_p\in I^r$ implies
$$
\frac{\partial^{i+j+k}F_p}{\partial x^i\partial y^j\partial z^k}\in I
$$
for all $i+j+k\le r-1$. Thus $\hat{\cal I}_{\overline S_r(X),p}\subset I$.
\vskip .2truein
Suppose that $p$ is a 2 point.
$$
\begin{array}{ll}
u&=(x^ay^b)^m\\
v&=P(x^ay^b)+x^cy^dF_p
\end{array}
$$
with $F_p\in I^r$ and $r\ge 2$,
$F_p\in I^r$ implies 
$$
xy\in \hat{\cal I}_{\text{sing}(\Phi_X),p}
=\sqrt{x^{ma+c-1}y^{mb+d-1}\left((ad-bc)F_p+ay\frac{\partial F_p}{\partial y}
-bx\frac{\partial F_p}{\partial x},y\frac{\partial F_p}{\partial z},
x\frac{\partial F_p}{\partial z}\right)}\subset I,
$$
so that $x\in I$ or $y\in I$. Without loss of generality, $x\in I$.

There exist permissible parameters $(\overline x,\overline y,\overline z)$ at $p$
such that $\overline y,\overline z\in {\cal O}_{X,p}$, $\overline x^a\overline y^b=x^ay^b$
and 
$$
x=\sigma\overline x,
y=\tau \overline y,
z=\overline z+h
$$
for some series $\sigma, \tau, h\in\hat{\cal O}_{X,p}$ with
$$
\sigma\equiv 1\text{ mod }m^{\overline a},
$$
$$
\tau\equiv 1\text{ mod }m^{\overline a}
$$
$$
h\equiv 0\text{ mod }m^{\overline a}
$$
where $m=m_p\hat{\cal O}_{X,p}$,
$$
\overline a\ge \frac{r+c}{a}(a+b)-(c+d).
$$
We have
$$
\begin{array}{ll}
u&=(\overline x^a\overline y^b)^m\\
v&=P(\overline x^a\overline y^b)+\overline x^c\overline y^d
[\sigma^c\tau^dF_p(\sigma\overline x,\tau\overline y,\overline z+h)]
\end{array}
$$
$$
\sigma^c\tau^dF_p(\sigma\overline x,\tau \overline y,\overline z+h)
\equiv F_p(\overline x,\overline y,\overline z)\text{ mod }m^{\overline a}
$$
Let $v=P_1(\overline x^a\overline y^b)+\overline x^c\overline y^dF_1$
be the normalized form of $v$. Since $F_p(\overline x,\overline y,\overline z)$
is normalized, we can only remove terms
$$
(\overline x^a\overline y^b)^t/\overline x^c\overline y^d
$$
with $t(a+b)-(c+d)\ge \overline a$ from
$\sigma^c\tau^dF_p(\sigma\overline x,\tau \overline y,\overline z+h)$
to construct $F_1$. Since this condition implies
$$
at-c\ge  r
$$
we have $F_1\in I^r$. We can thus assume that $y,z\in{\cal O}_{X,p}$.

 Set
$$
w=\frac{v-P_t(x^ay^b)}{x^cy^d}
$$
with $t>c+d+r$. Thus
$$
w=F_p+x^my^mh(x,y)
$$
with $m>r$.   $F_p\in I^r$ implies
$w\in I^r$ which implies that
$$
\frac{\partial^{i+j+k} w}{\partial x^i\partial y^j\partial z^k}\in I
$$
if $i+j+k\le r-1$.

There exists an \'etale neighborhood $\sigma:V\rightarrow U$ of $p$ such that $(x,y,z)$ are uniformizing parameters on $V$, $xy=0$ is a local equation of $E_X\cap V$ in $V$.

Suppose that 
$$
q\in V\left(x,\frac{\partial^{i+j+k}w}{\partial x^i\partial y^j\partial z^k}
\mid i+j+k\le r-1\right)\subset V
$$
is a 1 point in V. $q$ has permissible parameters $(\tilde x,\tilde y,\tilde z)$
defined by 
\begin{equation}\label{eq1047}
\tilde y=y-\alpha,
\tilde z=z-\beta,
\tilde x=xy^{\frac{b}{a}}
\end{equation}
for some $\alpha,\beta\in k$. Thus
$$
\begin{array}{ll}
u&=\tilde x^{am}\\
v&= P_t(\tilde x^a)
+\alpha^{d-\frac{cb}{a}}\tilde x^cw(\tilde x\alpha^{-\frac{b}{a}},\alpha,\beta)\\
&+\tilde x^c\left[
(\tilde y+\alpha)^{d-\frac{cb}{a}}w(x,y,z)-\alpha^{d-\frac{cb}{a}}
w(\tilde x\alpha^{-\frac{b}{a}},\alpha,\beta)\right]
\end{array}
$$
$\nu_q(w)\ge r$ implies $q\in S_r(V)$.

Suppose that 
$$
q\in V\left(x,\frac{\partial^{i+j+k}w}{\partial x^i\partial y^j\partial z^k}
\mid i+j+k\le r-1\right)\subset V
$$
is a 2 point in V.
 $q$ has permissible parameters $(x, y,\tilde z)$
where
$$
\tilde z=z-\beta.
$$
for some $\beta\in k$.
$$
\begin{array}{ll}
u&=(x^ay^b)^m\\
v&=P_t(x^ay^b)+x^cy^d[\sum_{(c+i)b-a(j+d)=0}\frac{\partial^{i+j}w}{\partial x^{i}\partial y^{j}}(q)
x^iy^j]\\
&+x^cy^d\left[w-
\sum_{(c+i)b-a(j+d)=0}\frac{\partial^{i+j}w}{\partial x^i\partial y^j}(q)x^iy^j\right]\\
\end{array}
$$
Again, $\nu_q(w)\ge r$ implies $q\in S_r(V)$. So
$$
V(I)\subset V\left(x,\frac{\partial^{i+j+k}w}{\partial x^i\partial y^j\partial z^k}
\mid i+j+k\le r-1\right)
\subset \overline S_r(V)
$$
implies
$$
\hat{\cal I}_{\overline S_r,p}\subset
\sqrt{\left(x,\frac{\partial^{i+j+k}w}{\partial x^i\partial y^j\partial z^k}
\mid i+j+k\le r-1\right)}
\subset I.
$$

We now prove 2. Suppose that the assumptions of 2. hold. If $p$ is a 1 point, then
(\ref{eq968}) implies $\hat{\cal I}_{\overline S_2(X),p}\subset I$.

If $p$ is a 2 point, then arguing as in the proof of 1., we set
$$
w=\frac{v-P_t(x^ay^b)}{x^cy^d}
$$
and conclude that $\frac{\partial w}{\partial y},\frac{\partial w}{\partial z}\in I$.
Suppose that $q\in V(x,\frac{\partial w}{\partial y},\frac{\partial w}{\partial z})$ is a 1 point.
Then there exists $\overline c\in k$ such that $\text{mult}(w-\overline cx)\ge 2$,
where the multiplicity is computed at $q$. $q$ has permissible parameters as in (\ref{eq1047}).
$$
x\equiv \alpha^{-\frac{b}{a}}\tilde x\text{ mod }m_q^2\hat{\cal O}_{X,q}
$$
implies $\text{mult}(w-\overline c\alpha^{-\frac{b}{a}}\tilde x)\ge 2$, so that $q\in\overline S_2(X)$. We have a simpler argument if
$q\in V(x,\frac{\partial w}{\partial y},\frac{\partial w}{\partial z})$
is a 2 point.
Thus
$$
\hat{\cal I}_{\overline S_2(X),p}\subset \sqrt{(x,\frac{\partial w}{\partial y},\frac{\partial w}{\partial z})}\subset I.
$$ 

\end{pf}

\begin{Lemma}\label{Lemma2} Suppose that $C\subset X$ is a curve and 
there exists $p\in X$
  such that $F_p\in \hat{\cal I}_{C,p}^r$ with $r\ge 2$. Then $C\subset E_X$.
\begin{enumerate}
\item Suppose that $C\subset E_X$ is a curve and there exists $p\in X$ such that 
$F_p\in\hat{\cal I}_{C,p}^r$ with $r\ge 1$. 
Suppose that $q\in C$ is a 1 point. Then 
$F_q\in \hat{\cal I}_{C,q}^r$.
\item Suppose that $C$ is a 2 curve and there exists $p\in C$ such that $F_p\in\hat{\cal I}_{C,p}^r$ with $r\ge 1$. If $q\in C$ is a 2 point, then $F_q\in \hat{\cal I}_{C,q}^r$.
\end{enumerate}
\end{Lemma}

\begin{pf}
We will first show that 1 or 2 hold for all but finitely many $q\in C$. 

First Suppose that $p$ is a 1 point. By Lemma \ref{Lemma301}, we may assume that
$y,z\in {\cal O}_{X,p}$.
$$
u=x^a,
v=P(x)+x^bF_p.
$$
In an \'etale neighborhood $U$ of $p$, $(x,y,z)$ are uniformizing parameters.
Let $I=\hat{\cal I}_{C,p}$.
If $F_p\in I^r$ with $r\ge 2$ then
$$
\frac{\partial F_p}{\partial y}, \frac{\partial F_p}{\partial z}\in I^{r-1}\subset I.
$$

$$
u\in \hat{\cal I}_{\text{sing}(\Phi_X),p}=\sqrt{x^{a-1+b}(\frac{\partial F_p}{\partial y}, \frac{\partial F_p}{\partial z})}\subset I.
$$
implies $x\in I$.

Now assume that $F_p\in\hat{\cal I}_{C,p}^r$ with $r\ge 1$ so that $C \subset E_X$,
either by assumption if $r=1$, or by the above argument if $r\ge 2$. 
Thus $x\in \hat{\cal I}_{C,p}$. Set
$w=\frac{v-P_{b+r}(x)}{x^b}$. 

After possibly replacing $U$ with a smaller \'etale neighborhood of $p$,
there exists a reduced ideal $J=(x,f)\subset \Gamma(U,{\cal O}_U)$ such that
$J\hat{\cal O}_{X,p}= I$. If $q\in V(J)\subset U$, then $q$ has regular parameters $(x,y-\alpha,z-\beta)$
for some $\alpha,\beta\in k$.  $w\in J^r\hat{\cal O}_{X,p}$ implies $w\in J^r$
(since $J$ is a complete intersection implies $J^r$ has no embedded components). Since $\nu(F_q)\ge r$,
we have
$$
F_q=w-\sum_{i\ge r}\frac{1}{i!}\frac{\partial^iw}{\partial x^i}(0,\alpha,\beta)x^i\in J^r\hat{\cal O}_{X,q}.
$$
Thus for all but finitely many $q\in C$, 1. holds.

Now suppose that $p$ is a 2 point. 
$$
u=(x^ay^b)^m,
v=P(x^ay^b)+x^cy^dF
$$
where $F=F_p$.

Suppose that
$F\in \hat{\cal I}_{C,p}^r$ with $r\ge 2$. Then
 $F,\frac{\partial F}{\partial x},\frac{\partial F}{\partial y},\frac{\partial F}{\partial z}\in \hat{\cal I}_{C,p}$. By (\ref{eq1004}),
$$
u\in \hat{\cal I}_{\text{sing}(\Phi_X),p}\subset\sqrt{(F,\frac{\partial F}{\partial x},
\frac{\partial F}{\partial y},\frac{\partial F}{\partial z})}\subset \hat{\cal I}_{C,p}
$$
implies $x\in \hat{\cal I}_{C,p}$ or $y\in \hat{\cal I}_{C,p}$. 

Now assume that $F_p\in\hat{\cal I}_{C,p}^r$ with $r\ge 1$. Then $C\subset E_X$,
either by assumption if $r=1$, or by the above argument if $r\ge 2$.
Thus we have $x$ or $y\in\hat{\cal I}_{C,p}$. Suppose that 
$x\in \hat{\cal I}_{C,p}$. As in the proof of Lemma \ref{Lemma302}, we may assume that
$y,z\in{\cal O}_{X,p}$.
Set $t=c+d+r$,  $w=\frac{v-P_t(x^ay^b)}{x^cy^d}$.
There exists an \'etale neighborhood $U$ of $p$ such that
$u=xy=0$ is a local equation of $E_X$ in $U$, $(x,y,z)$ are uniformizing parameters in $U$,
 and a reduced ideal
$$ 
J=(x,f)\subset \Gamma(U,{\cal O}_U)
$$
such that $J\hat{\cal O}_{X,p}=\hat{\cal I}_{C,p}$. $w\in J^r$ since $J$ is a complete
intersection. 

Suppose that 
$C$ is not a 2 curve, so that $\hat{\cal I}_{C,p}\ne (x,y)$. After possibly replacing $U$
with a smaller \'etale neighborhood of $p$, we can assume that 
$U\cap C\cap \overline B_2(X) = p$.

If $q\ne p$, and $q\in V(J)\subset U$, then $q$ has regular parameters $(x,y-\alpha,z-\beta)$ such that $\alpha\ne 0$.
 $\Phi_X(q)$ has permissible parameters 
$$
u_1=u,
v_1=v-v(\Phi_X(q))
$$
with permissible parameters $(\overline x,\overline y,\overline z)$, defined by
$$
x=\overline x(\overline y+\alpha)^{\frac{-b}{a}},
\overline y=y-\alpha,
\overline z=z-\beta
$$

$(\overline y+\alpha)^{d-\frac{bc}{a}}w\in\hat{\cal I}_{C,q}^r$, $\overline x\in \hat{\cal I}_{C,q}$,
and 
$$
F_q=(\overline y+\alpha)^{d-\frac{bc}{a}}w-\Omega(\overline x)
$$
with $\text{mult}(\Omega)\ge r$, which implies $F_q\in\hat{\cal I}_{C,q}^r$.

Now suppose that $C$ is a 2 curve, so that 
$\hat{\cal  I}_{C,p}=(x,y)$. If $q\in V(J)$, then $(u,v-v(\Phi_X(q)))$ are permissible parameters at $\Phi_X(q)$, and $q$ has
permissible parameters $(x,y,\overline z)$ with $\overline z=z-\alpha$.
$w\in\hat{\cal I}_{C,q}^r$, $x\in\hat{\cal I}_{C,q}$ and
$$
F_q=w-\frac{\Omega(x^ay^b)}{x^cy^d}
$$
for some $\Omega$
with $\text{mult}(\frac{\Omega(x^ay^b)}{x^cy^d})\ge r$. Thus
$F_q\in\hat{\cal I}_{C,q}^r$.

Now suppose that $p$ is a 3 point.
$$
u=(x^ay^bz^c)^m,
v=P(x^ay^bz^c)+x^dy^ez^fF
$$
We can assume that $y,z\in{\cal O}_{X,p}$. 
$F\in \hat{\cal I}_{C,p}^r$ with $r\ge 2$ implies that 
$$
F,\frac{\partial F_p}{\partial x},\frac{\partial F_p}{\partial y}, \frac{\partial F_p}{\partial z}\in \hat{\cal I}_{C,p}
$$
By (\ref{eq1005}), 
$$
u\in\hat{\cal I}_{\text{sing}(f),p}\subset \sqrt{(F,\frac{\partial F_p}{\partial x},\frac{\partial F_p}{\partial y}, \frac{\partial F_p}{\partial z})}
\subset \hat{\cal I}_{C,p}
$$
Thus $x,y$ or $z\in \hat{\cal I}_{C,p}$. 

Now suppose that $F\in\hat{\cal I}_{C,p}^r$ with $r\ge 1$. If $r=1$, then
$x,y$ or $z\in\hat{\cal I}_{C,p}$ by assumption. If $r\ge 2$, then  
$x,y$ or $z\in\hat{\cal I}_{C,p}$ by the above argument. Suppose that $x\in \hat{\cal I}_{C,p}$.
 Set $t=d+e+f+r$,
$$
w=\frac{v-P_t(x^ay^bz^c)}{x^dy^ez^f}.
$$
There exists an \'etale neighborhood $U$ of $p$ such that $(x,y,z)$ are uniformizing parameters in $U$,
$u=xyz=0$ is a local equation of $E_X$ in $U$ and $J=\Gamma(U,{\cal I}_C)=(x,f)$
is a complete intersection. $w\in J^r$
since $J$ is a complete intersection.

Suppose that $C$ is not a 2 curve. Then we can assume that $U\cap \overline{B_2}(X)\cap C=p$.
If $q\in V(J)\subset U$ and $q\ne p$, then $\Phi_X(q)$ has permissible parameters 
$$
u_1=u,
v_1=v-v(\Phi_X(q))
$$
with permissible parameters $(\overline x,\overline y,\overline z)$ at $q$, with
$$
x=\overline x(\overline y+\alpha)^{\frac{-b}{a}}(\overline z+\beta)^{-\frac{c}{a}},
y=\overline y+\alpha,
z=\overline z+\beta
$$
with $\alpha,\beta\ne 0$.
$w\in\hat{\cal I}_{C,q}^r$, $\overline x\in\hat{\cal I}_{C,q}$ and
$$
F_q=(\overline y+\alpha)^{e-\frac{bd}{a}}(\overline z+\beta)^{f-\frac{cd}{a}}w-\Omega(\overline x).
$$
$\text{mult}(\Omega)\ge r$ implies $F_q\in \hat{\cal I}_{C,q}^r$.

Suppose that $C$ is a 2 curve, $\hat{I}_{C,p}=(x,y)$, $q\in V(J)$.
Then 
$$
u_1=u,
v_1=v-v(\Phi_X(q))
$$
are permissible parameters at $\Phi_X(q)$, with permissible parameters $(\overline x,y,
\overline z)$ at $q$,
$$
\overline z=z-\alpha,
x=\overline x(\overline z+\alpha)^{\frac{-c}{a}}
$$
$$
u=(\overline x^ay^b)^m=(\overline x^{\overline a}\overline y^{\overline b})^{\overline m},
$$
with $(\overline a,\overline b)=1$,
 $w\in{\cal I}_{C,q}^r=(\overline x,y)^r$ implies
$$
F_q=(z+\alpha)^{f-\frac{dc}{a}}w-\frac{\Omega(\overline x^{\overline a}y^{\overline b})}{\overline x^dy^e}
$$
with $\text{mult}(\frac{\Omega(\overline x^{\overline a}y^{\overline b})}{\overline x^dy^e})\ge r$. Thus $F_q\in\hat{\cal I}_{C,q}^r$.

We conclude that 1. or 2. hold for all but finitely many $q\in C$. 

Suppose that $q\in C$ is a 1 point. We have at $q$,
$$
u=x^a,
v=P(x)+x^bF
$$
with $x\in\hat{\cal I}_{C,q}$, $y,z\in{\cal O}_{X,q}$. There exists an \'etale neighborhood $U$ of $q$ such that
$(x,y,z)$ are uniformizing parameters on $U$,  $x=0$ is a local equation of $E_X$, 
$J=\Gamma(U,{\cal I}_{C})=(x,f)$ is a complete intersection.
1. holds for all $q\ne q'\in U\cap E_X$ and
$$
w=\frac{v-P_{b+r}(x)}{x^b}\in\Gamma(U,{\cal O}_U).
$$
For $q'\in V(J)\subset U$ with $q'\ne q$, 
$$
u,v_1=v-v(\Phi_X(q'))
$$
are permissible parameters at $\Phi_X(q')$ and $(x,y-\alpha,z-\beta)$ 
are permissible parameters
 for $(u,v_1)$ at $q'$.
$$
F_{q'}=w-P_1(x)\in\hat{\cal I}_{C,q'}^r
$$
where 
$$
P_1(x)=\sum_{i=0}^{\infty}\frac{1}{i!}\frac{\partial^i w}{\partial x^i}(0,\alpha,\beta)
x^i=\sum_{i=0}^{\infty}a_ix^i.
$$
Set
$$\Lambda=w-\sum_{i=0}^{r-1}a_ix^i\in\Gamma(U,{\cal O}_{U}).$$ 
$\Lambda\in\hat{\cal I}_{C,q'}^r$ implies $\Lambda\in J^r$, so that
$\Lambda\in \hat{\cal I}_{C,p}^r$ implies 
$$
\frac{\partial^i\Lambda}{\partial x^i}(0,0,0)=0
$$
 for $i<r$, and
$$
F_p=\Lambda-\sum_{i=r}^{\infty}\frac{1}{r!}\frac{\partial^i\Lambda}{\partial x^i}(0,0,0)
x^i\in\hat{\cal I}_{C,p}^r.
$$
Now suppose that $C$ is a 2 curve. Suppose that $q\in C$ is a 2 point. We have at $q$,
$$
\begin{array}{ll}
u&=(x^ay^b)^m\\
v&=P(x^ay^b)+x^cy^dF
\end{array}
$$
with $\hat{\cal I}_{C,q}=(x,y)$, $y,z\in {\cal O}_{X,q}$.
 There exists an \'etale neighborhood $U$ of $p$ such that
$(x,y,z)$ are uniformizing parameters on $U$, $xy=0$ is a local equation of $E_X\cap U$,
$J=\Gamma(U,{\cal I}_C)=(x,y)$.
2. holds for all 2 points $q\ne q'\in U\cap E_X$.
$$
w=\frac{v-P_{c+d+r}(x^ay^b)}{x^cy^d}\in\Gamma(U,{\cal O}_U).
$$
For $q'\in V(J)\subset U$ with $q'\ne q$,
there exist permissible parameters 
$(x,y,z-\beta)$ at $q'$   for $(u,v-v(\Phi_X(q)))$.
$$
F_{q'}=w-\frac{P_1(x^ay^b)}{x^cy^d}\in \hat{\cal I}_{C,q'}^r
$$
where 
$$
P_1(x^ay^b)=\sum_{i=0}^{\infty}a_i(x^ay^b)^i.
$$
is a series.
Set $\Lambda=w-\frac{\sum_{i=0}^{c+d+r-1}a_i(x^ay^b)^i}{x^cy^d}\in\Gamma(U,{\cal O}_U)$.
$\Lambda\in\hat{\cal I}_{C,q'}^r$ implies
$\Lambda\in J^r$, so that
$\Lambda\in \hat{\cal I}_{C,q}^r$ and 
$$
\frac{\partial^{(a+b)i-c-d}\Lambda}{\partial x^{ai-c}\partial y^{bi-d}}(0,0,0)=0
$$
for $(a+b)i-c-d<r$, so that
$$
F_p=\Lambda-\sum_{(a+b)i-c-d\ge r}\frac{1}{(ai-c)!(bi-d)!}\frac{\partial^{(a+b)i-c-d}\Lambda}
{\partial x^{ai-c}\partial y^{bi-d}}(0,0,0)
x^{ai-c}y^{bi-d}
\in\hat{\cal I}_{C,p}^r.
$$

\end{pf}

\begin{Lemma}\label{Lemma5} Suppose that $r\ge 2$, $C\subset \overline S_r(X)$ is a nonsingular curve, and
$p\in C$ is a 1 point, so that there exist 
 permissible parameters 
$x,y,z$ at $p$  such that 
$$
\begin{array}{ll}
u&=x^a\\
v&=P(x)+x^cF
\end{array}
$$
where $\hat{\cal I}_{C,p} = (x,z)$. Then
$$
F_p = a_r(x,y) + a_{r-1}(x,y)z+\cdots +a_1(x,y)z^{r-1}+g(x,y,z)z^r
$$
where 
$$
x^i \mid a_i\text{ for }1\le i\le r-1,
$$
and
$x^{r-1}\mid a_r$.
\end{Lemma}

\begin{pf}
There exist permissible parameters $(x,\overline y,\overline z)$ at $p$ such that
$\overline y,\overline z\in{\cal O}_{X,p}$ and $(x,\overline z)=\hat{\cal I}_{C,p}$.
Then there exists $a,b\in\hat{\cal O}_{X,p}$ such that $\overline z=ax+bz$
where $b$ is a unit.
Assume that the conclusions of the Lemma are true for the variables
$(x,\overline y,\overline z)$. Then substituting for $x,y,z$ we get the conclusions
of the Lemma for $(x,y,z)$, so we may suppose that $y,z\in {\cal O}_{X,p}$.

 There exists an \'etale neighborhood $U$ of $p$
 such that $x,y,z$ are  uniformizing parameters
in $U$, 
$x=0$ is a local equation of $E_X$ in $U$, $x=z=0$ are equations of $C\cap U$.
 If $p'\in U\cap C$, and $\alpha=y(p')$, then $(x,y_1=y-\alpha, z)$
are permissible parameters at $p'$.
$$
F_{p'} = \sum_{i\ge 0,j+k>0}\frac{1}{(c+i)!j!k!}\frac{\partial^{c+i+j+k} v}{\partial x^{c+i}\partial y^j\partial z^k}(0,\alpha,0)
x^iy_1^jz^k.
$$
$\nu(p')\ge r$ for $p'\in C\cap U$  implies that, if $j+k>0$ and $i+j+k<r$, then
$$
\frac{\partial^{c+i+j+k} v}{\partial x^{c+i}\partial y^j\partial z^k}(0,\alpha,0)=0
$$
for infinitely many $\alpha$, so that
$$
\frac{\partial^{c+i+j+k} v}{\partial x^{c+i}\partial y^j\partial z^k}(0,y,0)=0 
$$
in $U$, if $j+k>0$ and $i+j+k<r$. Thus 
$$
\frac{\partial^{c+i+k} v}{\partial x^{c+i}\partial z^k}(0,y,0)=0
$$
in $U$ if $i+k<r$, $k>0$, so that
$$
\frac{\partial^{c+i+j+k} v}{\partial x^{c+i}\partial y^j\partial z^k}(0,y,0)= 
\frac{\partial^j}{\partial y^j}\left[\frac{\partial^{c+i+k}v}{\partial x^{c+i}\partial z^k}(0,y,0)\right]=0
$$
if $i+k<r$, $k>0$ and $j\ge0$. Thus
$$
\frac{\partial^{c+i+j+k} v}{\partial x^{c+i}\partial y^j\partial z^k}(0,0,0)= 0
$$
if
$k>0$, $i<r-k$, $j\ge 0$,
$$
\frac{\partial^{c+i+1}v}{\partial x^{c+i}\partial y}(0,y,0)=0
$$
if $i<r-1$, so that
$$
\frac{\partial^{c+i+j}v}{\partial x^{c+i}\partial y^j}(0,0,0)=0
$$
if $i<r-1$, $j>0$, and the conclusions of the Lemma follow.
\end{pf}

\begin{Lemma}\label{Lemma6} Suppose that $r\ge 1$, $C\subset X$ is a 2 curve such that $\nu(p)\ge r$ if $p\in C$ is a 2 point ($C\subset\overline S_r(X)$ if $r\ge 2$), and
$p\in C$ is a 2 point, so that there exist 
 permissible parameters 
$x,y,z$ at $p$  such that 
$$
\begin{array}{ll}
u&=(x^ay^b)^m\\
v&=P(x^ay^b)+x^cy^dF
\end{array}
$$
where $\hat{\cal I}_{C,p} = (x, y)$.  Then there exists a series $\tau(z)$ with $\text{mult }
\tau(z)\ge 1$ such that 
$$
F = \left\{\begin{array}{ll} \tau(z)x^{i_0}y^{j_0}+\sum_{i+j\ge r}a_{ij}(z)x^iy^j
&\text{ if there exist nonnegative integers $(i_0,j_0)$ such that}\\
&\text{ $i_0+j_0=r-1$ and $a(d+j_0)-b(c+i_0)=0$}\\
\sum_{i+j\ge r}a_{ij}(z)x^iy^j& \text{ otherwise}
\end{array}\right. 
$$
\end{Lemma}

\begin{pf}
There exist permissible parameters $(\overline x,\overline y,\overline z)$ at $p$ such 
that $\overline y,\overline z\in {\cal O}_{X,p}$. $\sigma x=\overline x$,
$\omega y=\overline y$, $z=\overline z+h$, $\sigma^a\omega^b=1$ with $\sigma,\omega,h\in\hat{\cal O}_{X,p}$,  $\sigma,\omega\equiv 1\text{ mod }m_p^2\hat{\cal O}_{X,p}$,
$h\equiv 0\text{ mod }m_p^2\hat{\cal O}_{X,p}$. Suppose that the
conclusions of the Lemma hold for $(\overline x,\overline y,\overline z)$.
Substituting for $(x,y,z)$ we get the conclusions of the Lemma for $(x,y,z)$.
We may thus assume that $y,z\in{\cal O}_{X,p}$.

There exists an \'etale neighborhood $U$ of $p$ such that $x,y,z$ are  uniformizing parameters
in $U$, $xy=0$ is a local equation of $U\cap E_X$.  Set
$$
w = \frac{v-P_{t}(x^ay^b)}{x^cy^d}
$$
where  $t>c+d+r$. We have 
$w\in\Gamma(U,{\cal O}_{U,p})$ and
$$
\begin{array}{ll}
u&=(x^ay^b)^m\\
v&= P_{t}(x^ay^b)+x^cy^dw.
\end{array}
$$
If $p'\in U\cap C$, and $\alpha=z(p')$, then $(x,y,z_1= z-\alpha)$
are permissible parameters at $p'$.
$$
F_{p'} = \sum_{k>0,\,i,j\ge 0} \frac{1}{i!j!k!}\frac{\partial^{i+j+k} w}{\partial x^i\partial y^j \partial z^k}(0,0,\alpha)
x^iy^jz_1^k +
\sum_{a(i+c)-b(j+d)\ne 0}\frac{1}{i!j!}\frac{\partial^{i+j} w}{\partial x^i \partial y^j}(0,0,\alpha)x^iy^j
$$ 
$\nu(p')\ge r$ for $p'\in C\cap U$ implies that for infinitely many $\alpha$, we have
$$
\frac{\partial^{i+j+k} w}{\partial x^i\partial y^j \partial z^k}(0,0,\alpha)=0\text{ if }
k>0\text{ and }i+j+k<r
$$
and
$$
\frac{\partial^{i+j} w}{\partial x^i \partial y^j}(0,0,\alpha)=0 \text{ if }
a(i+c)-b(j+d)\ne 0\text{ and }i+j<r.
$$
Thus 
$$
\frac{\partial^{i+j+k} w}{\partial x^i\partial y^j \partial z^k}(0,0,z)=0\text{ if }
a(i+c)-b(j+d)\ne 0\text{ and }i+j<r
$$
and
$$
\frac{\partial^{i+j+1} w}{\partial x^i\partial y^j \partial z}(0,0,z)=0\text{ if }
i+j+1<r,
$$
so that
$$
\frac{\partial^{i+j+k}w}{\partial x^i\partial y^j\partial z^k}(0,0,z)=0
$$
if $i+j<r-1$, $k>0$.
Setting $z=0$ in the above equations, we get the statement of the Lemma.
\end{pf}

\begin{Lemma}\label{Lemma7} Suppose that $C\subset X$ is a nonsingular curve
containing a 1 point, 
$p\in C$ is a 2 point such that $C$ makes SNCs with the 2 curve through $p$,
and $(x,y,z)$ are   
permissible parameters 
 at $p$   such that $x=z=0$ are local equations of $C$ at $p$.
\begin{enumerate}
\item Write
$$
\begin{array}{ll}
u&=(x^ay^b)^m\\
v&=P(x^ay^b)+x^cy^dF
\end{array}
$$
  Then if $r\ge 2$ and $C\subset\overline S_r(X)$,
$$
F = x^{r-1}\tau(y)+\sum_{i+k\ge r}a_{ijk}x^iy^jz^k
$$
where $\tau$ is a series with $\text{mult}(\tau(y))\ge 0$.
\item
If there exists a 1 point $q\in C$ such that $F_q\in\hat{\cal I}_{C,q}^r$ with
$r\ge 1$, then $F_p\in \hat{\cal I}_{C,p}^r$.
\end{enumerate}
\end{Lemma}

\begin{pf} 
There exist permissible parameters $(\overline x,\overline y,\overline z)$ at $p$
such that $\overline y,\overline z\in {\cal O}_{X,p}$,
$$
\sigma x=\overline x,
\omega y=\overline y,
\sigma^a\omega^b=1,
$$
 with $\sigma,\omega\in\hat{\cal O}_{X,p}$,
$\sigma,\omega\equiv 1\text{ mod }m_p\hat{\cal O}_{X,p}$,
$\hat{\cal I}_{C,p}=(\overline x,\overline z)$. Then $\overline z=\overline ax+
\overline b z$ for some $\overline a,\overline b\in\hat{\cal O}_{X,p}$.

Suppose that the conclusions of the Lemma hold for $(\overline x,\overline y,\overline z)$.
Substituting back for $(x,y,z)$, we get the conclusions of the Lemma for $(x,y,z)$.
We may thus assume that $y,z\in{\cal O}_{X,p}$

There exists an \'etale neighborhood $U$ of $p$ such that $x,y,z$ are  uniformizing parameters
in $U$, $xy=0$ is a local equation of $E_X\cap U$, $C\cap U=V(x,z)$ in $U$.  Set
$$
w = \frac{v-P_{t}(x^ay^b)}{x^cy^d}.
$$
where  $t>r+c+d$.
We have 
$w\in \Gamma(U,{\cal O}_{U})$,
$$
\begin{array}{ll}
u&=(x^ay^b)^m\\
v&= P_{t}(x^ay^b)+x^cy^dw
\end{array}
$$
If $p'\in U\cap C$, and $\alpha=y(p')$, then $(x,y-\alpha, z)$
are regular parameters in $\hat{\cal O}_{X,p'}$. We  have permissible parameters $\overline x, \overline y, z$ at $p'\ne p$
 defined by 
$$
x = \overline x(\overline y+\alpha)^{-\frac{b}{a}}, y = \overline y+\alpha.
$$
At $p'$, we have 
\begin{equation}\label{eq652}
\begin{array}{ll}
u&=\overline x^{am}\\
v &= P_{t}(\overline x^a) + \overline x^c(\overline y+\alpha)^{d-\frac{cb}{a}}w
\end{array}
\end{equation}

$$
\begin{array}{ll}
x^cy^dw&=\overline x^c(\overline y+\alpha)^{d-\frac{cb}{a}}w \\
&=\overline x^c(\overline y+\alpha)^{d-\frac{cb}{a}}\left[\sum_{i,j,k\ge 0}
\frac{1}{i!j!k!}\frac{\partial^{i+j+k}w}{\partial x^i\partial y^j \partial z^k}(0,\alpha,0)
\overline x^i(\overline y+\alpha)^{-i\frac{b}{a}}\overline y^j z^k\right]\\
&=\sum_{i,k\ge 0}\left[\sum_{j\ge 0}\frac{1}{i!j!k!}\frac{\partial^{i+j+k}w}{\partial x^i\partial y^j \partial z^k}(0,\alpha,0)
\overline y^j\right] (\overline y+\alpha)^{d-\frac{b(c+i)}{a}}\overline x^{c+i}z^k.
\end{array}
$$
Thus

$$
\begin{array}{ll}
F_{p'}&= \sum_{i\ge 0} \left[\sum_{j\ge 0} 
\frac{1}{i!j!}\frac{\partial^{i+j}w}{\partial x^i\partial y^j}(0,\alpha,0)
\overline y^j\right] (\overline y+\alpha)^{d-\frac{b(c+i)}{a}}\overline x^{i}
- \sum_{i\ge 0} 
\frac{1}{i!}\frac{\partial^{i}w}{\partial x^i}(0,\alpha,0)
\alpha^{d-\frac{b(c+i)}{a}}\overline x^{i}\\
&+
\sum_{i\ge 0,k> 0}\left[\sum_{j\ge 0}\frac{1}{i!j!k!}\frac{\partial^{i+j+k}w}{\partial x^i\partial y^j \partial z^k}(0,\alpha,0)
\overline y^j\right] (\overline y+\alpha)^{d-\frac{b(c+i)}{a}}\overline x^{i}z^k.
\end{array}
$$

$\nu(p')=r$ implies that 
\begin{equation}\label{eq10}
\frac{\partial^{i+j+k}w}{\partial x^i\partial y^j \partial z^k}(0,\alpha,0)=0
\end{equation}
if $i+j+k<r$, $k>0$, and for fixed $i<r$ 
\begin{equation}\label{eq11}
\sum_{j<r-i}\frac{1}{i!j!}\frac{\partial^{i+j} w}{\partial x^i\partial y^j}(0,\alpha,0)\overline y^j \equiv
c_{\alpha}^i(\overline y+\alpha)^{\lambda_i}\text{ mod }\overline y^{r-i}
\end{equation}
where $\lambda_i = \frac{b(c+i)}{a}-d$, and the $c_{\alpha}^i\in k$ depend  on $\alpha$ and $i$.
Since (\ref{eq10}) holds for infinitely many $\alpha$, 
$$
\frac{\partial^{i+k}w}{\partial x^i\partial z^k}(0,y,0)=0
$$
if $i+k<r$ and $k>0$.
Thus
$$
\frac{\partial^{i+j+k}w}{\partial x^i\partial y^j\partial z^k}(0,0,0) = \frac{\partial^j}{\partial y^j}
\left[\frac{\partial^{i+k}w}{\partial x^i\partial z^k}(0,y,0)\right](0,0,0)=0
$$
if $i+k<r$, $k>0$.

If $i<r-1$ we have
$$
\frac{1}{i!}\frac{\partial^i w}{\partial x^i}(0,\alpha,0)=c_{\alpha}^i\alpha^{\lambda_i}
$$
and 
$$
\frac{1}{i!}\frac{\partial^{i+1} w}{\partial x^i\partial y}(0,\alpha,0) = c_{\alpha}^i\lambda_i\alpha^{\lambda_i-1}
$$
for infinitely many $\alpha$. Thus
$$
\lambda_i\frac{\partial^i w}{\partial x^i}(0,\alpha,0) = \alpha\frac{\partial}{\partial y}
\frac{\partial^i w}{\partial x^i}(0,\alpha,0)
$$
for infinitely many $\alpha$, and thus for all $\alpha$. Set $\gamma_i(y) = \frac{\partial^i w}{\partial x^i}(0,y,0)$.
We have
$$
\lambda_i\gamma_i(y) = y\frac{d\gamma_i}{dy}.
$$
There is an expansion $\gamma_i(y) = \sum_{j=0}^{\infty}b_jy^j$ with $b_j\in k$.
$\frac{d\gamma_i}{dy}=\sum_{j=1}^{\infty}jb_jy^{j-1}$.
$$
y\frac{d\gamma_i}{dy} = \sum_{j=0}^{\infty}jb_jy^j
$$
$$
\lambda_i\gamma_i-y\frac{d\gamma_i}{dy}=\sum_{j=0}^{\infty}(\lambda_ib_j-jb_j)y^j=0
$$
so that $b_j(\lambda_i-j)=0$ for all $j$, which implies that
$\gamma_i=0$, or $\lambda_i\in \bold N$ and $\gamma_i=b_{\lambda_i}y^{\lambda_i}$.
Suppose that $\lambda_i\in {\bold N}$ and $\gamma_i(y) = \frac{\partial^iw}{\partial x^i}(0,y,0)\ne 0$.
Then
$$
\frac{\partial^{i+\lambda_i}w}{\partial x^i\partial y^{\lambda_i}}(0,0,0)= \lambda_i!b_{\lambda_i}\ne 0.
$$
But 
$$
b(c+i)-a(d+\lambda_i)=b(c+i)-a(\frac{b}{a}(c+i))=0
$$
 implies $i>r$, by our
choice of $t$ in $p_t$ and the assumption that $F$ is normalized, a contradiction.
Thus
$$
\frac{\partial^{i+j}w(0,0,0)}{\partial x^i\partial y^j}=0
$$
if $i<r-1$.

Now suppose that there exists a 1 point $q'\in C$ such that $F_{q'}\in\hat{\cal I}_{C,q'}^r$.
By 1. of Lemma \ref{Lemma2}, $F_q\in\hat{\cal I}_{C,q}^r$ at every 1 point $q\in C$.
With the above notation, (trivially if $r=1$)
$$
F_p=x^{r-1}\tau(y)+\sum_{i+k\ge r}a_{ijk}x^iy^jz^k.
$$
For $p\ne q\in C\cap U$ there exist permissible parameters $(\overline x,\overline y,
z)$ at $q$ such that
$$
x=\overline x(\overline y+\alpha)^{-\frac{b}{a}},
y=\overline y+\alpha
$$
$$
F_q=\overline x^{r-1}\Lambda+\Omega
$$
with
$$
\Lambda=(\overline y+\alpha)^{d-\frac{b}{a}(c+r-1)}\tau(\overline y+\alpha)
-\alpha^{d-(c+r-1)\frac{b}{a}}\tau(\alpha),
$$
$\Omega\in\hat{\cal I}_{C,q}^r$.
$F_q\in\hat{\cal I}_{C,q}^r$ implies $\tau=0$ or
$d-\frac{b}{a}(c+r-1)=0$ and $\tau\in k$. But $ad-b(c+r-1)=0$ and $\tau\in k$ is
not possible since $F$ is normalized. Thus $\tau=0$.

\end{pf}

\begin{Lemma}\label{Lemma8} Suppose that $r\ge 1$ $C\subset X$ is a 2 curve, 
such that $\nu(q)\ge r$ if $q\in C$ is a 2 point ($C\subset\overline S_r(X)$ if $r\ge 2$) and
$p\in C$ is a 3 point, so that there exist  permissible parameters 
$x,y,z$ at $p$  such that 
$$
\begin{array}{ll}
u&=(x^ay^bz^c)^m\\
v&=P(x^ay^bz^c)+x^dy^ez^fF
\end{array}
$$
where $\hat{\cal I}_{C,p} = (x, y)$.  Then
$$
F = \left\{\begin{array}{ll} \tau(z)x^{i_0}y^{j_0}+\sum_{i+j\ge r}a_{ij}(z)x^iy^j
&\text{ if there exist $(i_0,j_0)$ such that}\\
&\text{ $i_0+j_0=r-1$ and $a(e+j_0)-b(d+i_0)=0$}\\
\sum_{i+j\ge r}a_{ij}(z)x^iy^j& \text{ otherwise}
\end{array}\right. 
$$
If there exists a 2 point $q\in C$ such that $F_q\in\hat{\cal I}_{C,q}^r$ and $r\ge 1$,
then $F_p\in\hat{\cal I}_{C,p}^r$
\end{Lemma}

\begin{pf}
There exist permissible parameters $(\overline x,\overline y,\overline z)$ at
$p$ such that $\overline y,\overline z\in{\cal O}_{X,p}$,
$$
\sigma x=\overline x,
\omega y=\overline y,
\mu z=\overline z
$$
for some unit series $\sigma, \omega, \mu\in \hat{\cal O}_{X,p}$. Suppose that the conclusions of 
the Lemma are true for the parameters $(\overline x,\overline y,\overline z)$.
Substituting back for $(x,y,z)$
we get the conclusions of the Lemma for $(x,y,z)$. We may thus assume that 
$y,z\in {\cal O}_{X,p}$.

There exists an \'etale neighborhood $U$ of $p$ such that $(x,y,z)$ are uniformizing parameters in $U$, $xyz=0$ is a local equation of $E_X\cap U$.
Set
$$
w=\frac{v-P_t(x^ay^bz^c)}{x^dy^ez^f}
$$
where $t\ge d+e+f+r$. We have
$w\in \Gamma(U,{\cal O}_U)$ and
$$
u=(x^ay^bz^c)^m,
v=P_t(x^ay^bz^c)+x^dy^ez^fw.
$$
If $p'\in U\cap C$ and $\alpha=z(p')$, then $(x,y,z-\alpha)$
 are regular parameters in $\hat{\cal O}_{X,p'}$. If $\alpha\ne 0$, we have 
permissible parameters $(\overline x,y,\overline z)$ at $p'$ where 
$x=\overline x(\overline z+\alpha)^{\frac{-c}{a}}$, $\overline z=z-\alpha$. At $p'$ we have
$$
u=(\overline x^ay^b)^m,
v=P_t(\overline x^ay^b)+\overline x^dy^e(\overline z+\alpha)^{f-\frac{dc}{a}}w.
$$
$$
\begin{array}{ll}
\overline x^dy^e(\overline z+\alpha)^{f-\frac{cd}{a}}w&=
\overline x^dy^e(\overline z+\alpha)^{f-\frac{cd}{a}}\left[
\sum_{i,j,k\ge 0}\frac{1}{i!j!k!}\frac{\partial^{i+j+k}w}{\partial x^i\partial y^j \partial z^k}(0,0,\alpha)
\overline x^i(\overline z+\alpha)^{\frac{-ic}{a}}y^j\overline z^k\right]\\
&=
\overline x^dy^e\left[
\sum_{i,j\ge 0}\left(\sum_{k\ge 0}\frac{1}{i!j!k!}\frac{\partial^{i+j+k}w}{\partial x^i\partial y^j \partial z^k}(0,0,\alpha)
\overline z^k\right)(\overline z+\alpha)^{f-\frac{c(i+d)}{a}}\overline x^iy^j\right].
\end{array}
$$
Thus
$$
\begin{array}{ll}
F_{p'}&=
\sum_{i,j\ge 0}\left(\sum_{k\ge 0}\frac{1}{i!j!k!}\frac{\partial^{i+j+k}w}{\partial x^i\partial y^j \partial z^k}(0,0,\alpha)
\overline z^k\right)(\overline z+\alpha)^{f-\frac{c(i+d)}{a}}\overline x^iy^j\\
&-\sum_{i,j\text{ such that }b(i+d)-a(j+e)=0}\left(\frac{1}{i!j!}\frac{\partial^{i+j}w}{\partial x^i\partial y^j}(0,0,\alpha)
\alpha^{f-\frac{c(i+d)}{a}}\right)\overline x^iy^j
\end{array}
$$
$\nu(p')=r$ implies 

\begin{equation}\label{eq305}
\frac{\partial^{i+j+k}w}{\partial x^i\partial y^j\partial z^k}(0,0,\alpha)=0\text{ if }i+j+k<r\text{ and }b(i+d)-a(j+e)\ne 0
\end{equation}

and if $b(i+d)-a(j+e)=0$ for fixed $i,j$ with $i+j<r$, 

\begin{equation}\label{eq306}
\sum_{k<r-i-j}\frac{1}{i!j!k!}\frac{\partial^{i+j+k}w}{\partial x^i\partial y^j\partial z^k}(0,0,\alpha)\overline z^k\equiv 
c_{\alpha}^{ij}(\overline z+\alpha)^{\lambda_i}\text{ mod }(\overline z^{r-i-j})
\end{equation}
where $\lambda_i=\frac{c(i+d)}{a}-f$, $c_{\alpha}^{i,j}\in k$ depend  on $\alpha,i$ and $j$.

Since (\ref{eq305}) holds for infinitely many $\alpha$,
$$
\frac{\partial^{i+j+k}w}{\partial x^i\partial y^j\partial z^k}(0,0,z)=0
$$
if $i+j+k<r$ and $b(i+d)-a(j+e)\ne 0$. Thus
$$
\frac{\partial^{i+j+k}w}{\partial x^i\partial y^j\partial z^k}(0,0,0)=\frac{\partial^k}{\partial z^k}\left[
\frac{\partial^{i+j}w}{\partial x^i\partial y^j}(0,0,z)\right](0,0,0)=0
$$
if $i+j<r$, $k\ge 0$ and $b(i+d)-a(j+e)\ne 0$.

If $b(i+d)-a(j+e)=0$ and $i+j<r-1$, we have
$$
\frac{1}{i!j!}\frac{\partial^{i+j}w}{\partial x^i\partial y^j}(0,0,\alpha)=c_{\alpha}^{ij}\alpha^{\lambda_i}
$$
and
$$
\frac{1}{i!j!}\frac{\partial^{i+j+1}w}{\partial x^i\partial y^j\partial z}(0,0,\alpha)=
c_{\alpha}^{ij}\lambda_i\alpha^{\lambda_i-1}
$$
for infinitely many $\alpha$. Thus
$$
\lambda_i\frac{\partial^{i+j}w}{\partial x^i\partial y^j}(0,0,\alpha)
=\alpha\frac{\partial^{i+j+1}w}{\partial x^i\partial y^j\partial z}(0,0,\alpha)
$$
for infinitely many $\alpha$, and thus for all $\alpha$. Set
$$
\gamma_{ij}(z)=\frac{\partial^{i+j}w}{\partial x^i\partial y^j}(0,0,z).
$$
We have an expression
$$
\gamma_{ij}(z)=\sum_{k= 0}^{\infty}b_kz^k
$$
with $b_k\in k$.
$$
\frac{d \gamma_{ij}}{dz}=\sum_{k=1}^{\infty}kb_kz^{k-1}.
$$
$$
\lambda_i\gamma_{ij}-z\frac{d\gamma_{ij}}{dz}=\sum_{k=0}^{\infty}
(\lambda_ib_k-kb_k)z^k=0
$$
implies $b_k(\lambda_i-k)=0$ for all $k$, so that either $\gamma_{ij}=0$, or 
$\lambda_i\in{\bold N}$ and
$\gamma_{ij}=b_{\lambda_i}z^{\lambda_i}$. Suppose that
$$
\gamma_{ij}(z)=\frac{\partial^{i+j}w}{\partial x^i\partial y^j}(0,0,z)\ne 0.
$$
Then
$$
\frac{\partial^{i+j+\lambda_i}w}{\partial x^i\partial y^j\partial z^{\lambda_i}}(0,0,0)
=\lambda_i!b_{\lambda_i}\ne 0
$$
so that we have a nontrivial $x^{d+i}y^{e+j}z^{f+\lambda_i}$ term in $x^dy^ez^fF$. Recall that
$\lambda_i=\frac{c(i+d)}{a}-f$. By assumption, $b(i+d)-a(j+e)=0$. We further have
$a(f+\lambda_i)-c(d+i)=0$,
$$
\begin{array}{ll}
b(f+\lambda_i)-c(e+j)&=b(\frac{c(i+d)}{a})-c(e+j)\\
&=b\frac{c}{a}(d+i)-c\frac{b}{a}(i+d)=0
\end{array}
$$
a contradiction to the assumption that $F$ is normalized. Thus
$$
\frac{\partial^{i+j+k}w}{\partial x^i\partial y^j\partial z^k}(0,0,0)=0
$$
if $b(i+d)-a(j+e)=0$, $i+j<r-1$, $k\ge0$.

Now suppose there exists a 2 point $q'\in C$ such that $F_{q'}\in\hat{\cal I}_{C,q'}^r$.
By 2. of Lemma \ref{Lemma2}, $F_q\in\hat{\cal I}_{C,q}^r$ for all 2 points
$q\in C$. With the above notation, if $F_p\not\in\hat{\cal I}_{C,p}^r$, we have
$$
F_p=\tau(z)x^{i_0}y^{j_0}+\sum_{i+j\ge r}a_{ij}(z)x^iy^j
$$
where $i_0+j_0=r-1$ and $a(e+j_0)-b(d+i_0)=0$, $\tau(z)\ne 0$.

For $p\ne q\in C\cap U$,  there exist permissible parameters $(\overline x,y,\overline z)$ at $q$ such that 
$$
x=\overline x(\overline z+\alpha)^{-\frac{c}{a}},
z=\overline z+\alpha
$$
$$
\begin{array}{ll}
u&=(\overline x^ay^b)^m=(\overline x^{\overline a}y^b)^{\overline m}\\
v&=P_q(\overline x^{\overline a}y^{\overline b})+\overline x^dy^eF_q
\end{array}
$$
with $(\overline a,\overline b)=1$.
$$
F_q=\overline x^{i_0}y^{j_0}\Lambda+\Omega
$$
with $\Omega\in\hat{\cal I}_{C,q}^r$ and
$$
\Lambda=(\overline z+\alpha)^{f-\frac{c(d+i_0)}{a}}\tau(\overline z+\alpha)
-\alpha^{f-\frac{c(d+i_0)}{a}}\tau(\alpha).
$$
$$
F_q\in \hat{\cal I}_{C,q}^r
$$
 implies $\tau=0$, or
$f-\frac{c(d+i_0)}{a}=0$ and $\tau\in k$.
$f-\frac{c(d+i_0)}{a}=0$ and $\tau\in k$ is not possible since $F$ is normalized.

\end{pf}

\begin{Lemma}\label{Lemma961} Suppose that $r\ge 2$, $f_1,\ldots,f_{n-1}$ is
a regular sequence in a $n$ dimensional regular local ring $A$. Let 
$I=(f_1,\ldots,f_{n-1})$. Then
$$
\text{depth } A/(I^r+(f_1)^{r-1})=1
$$
for all $r\ge 2$.
\end{Lemma}
\begin{pf} Let $m$ be the maximal ideal of $A$. There is an exact sequence of $A$ modules
$$
0\rightarrow I^{r-1}/(I^r+(f_1)^{r-1})\rightarrow A/(I^r+(f_1)^{r-1})\rightarrow A/I^{r-1}\rightarrow 0
$$
$I^{r-1}/(I^r+(f_1)^{r-1})$ is a free $A/I$ module, since $f_1,\ldots, f_{n-1}$ is
quasi regular by Theorem 27, \cite{Ma}.
$\text{depth }A/I^t= 1$ for all $t\ge 1$ by Proposition 16.F, \cite{Ma}.
Thus 
$\text{Hom}_A(A/m,A/I^{r-1})=0$ and $\text{Hom}_A(A/m,I^{r-1}/(I^r+(f_1)^{r-1}))=0$, so that
$\text{Hom}_A(A/m,A/(I^r+(f_1)^{r-1}))=0$ and $\text{depth}(A/(I^r+(f_1)^{r-1}))=1$.
\end{pf}

\begin{Lemma}\label{Lemma659} Suppose that $r\ge 2$, $p\in X$ is a 2 point and
$(x,y,z)$ are permissible parameters at $p$ such that $y,z\in {\cal O}_{X,p}$,
$$
\begin{array}{ll}
u&=(x^ay^b)^m\\
v&=P(x^ay^b)+x^cy^dF_p
\end{array}
$$
$p\in C\subset \overline S_r(X)$ is a curve such that $x\in\hat{\cal I}_{C,p}$
and $C$ contains a 1 point $q$. Then there exists a polynomial $g$ such that
one of the following cases hold
\begin{description}
\item[Case 1)] $\overline y^{ad-bc}F_p-g(x\overline y^b)\in
\left( \hat{\cal I}_{C,p}^{r}+(x)^{r-1}\right)
k[[x,\overline y,z]]$ if $ad-bc\ge 0$.
\item[Case 2)] $F_p-g(x\overline y^b)\overline y^{bc-ad}
\in \left(\hat{\cal I}_{C,p}^{r}+(x)^{r-1}\right)k[[x,\overline y,z]]$ if $ad-bc\le 0$
\end{description}
where $y=\overline y^a$.
\end{Lemma}

\begin{pf}
$x,y,z$ are uniformizing parameters in an \'etale neighborhood $U=\text{spec}(R)$
of $p$ and $xy=0$ is a local equation of $E_X\cap U$, $C$ is a complete intersection
in $U$. For $t>c+d+r$,
$$
w=\frac{1}{x^cy^d}[v-P_t(x^ay^b)]\in R.
$$
If $q\in C\cap U$ is a 1 point, such that $C$ is nonsingular at $q$, then $q$ has permissible parameters $(x_1,y_1,z_1)$
where $x=x_1(y_1+\alpha)^{-\frac{b}{a}}$, $y=y_1+\alpha$, $z=z_1+\beta$ for some
$\alpha, \beta\in k$ with $\alpha\ne 0$.
$$
\begin{array}{ll}
u&=x_1^{am}\\
v&=P_q(x_1)+ x_1^cF_q
\end{array}
$$
where 
$$
F_q=(y_1+\alpha)^{\frac{ad-bc}{a}}w-g(x_1)
$$
for some series $g$. $F_q\in\hat{\cal I}_{C,q}^{r}+(x_1)^{r-1}$  (by Lemma \ref{Lemma5})
implies
$$
w-(y_1+\alpha)^{\frac{bc-ad}{a}}g(x_1)\in \hat{\cal I}_{C,q}^{r}+(x_1)^{r-1}.
$$

 Let $y=\overline y^a$, $S=R[\overline y]$,
$\lambda:V= \text{spec}(S)\rightarrow \text{spec}(R)$. 
Suppose that $q'\in \lambda^{-1}(q)$. Let $h(x_1)=g_r(x_1)$.

Suppose that $ad-bc\ge 0$. Then 
$$
\overline y^{ad-bc}w-h(x\overline y^b)\in\hat{\cal I}_{C,q'}^{r}+(x)^{r-1}.
$$
$I=\Gamma(V,{\cal I}_C)$ is a complete intersection in $V$, so that (by Lemma \ref{Lemma961})
$$
\overline y^{ad-bc}w-h(x\overline y^b)\in (I^{r}+(x)^{r-1})S,
$$
and
$$
\overline y^{ad-bc}w-h(x\overline y^b)\in
\left(\hat{\cal I}_{C,p}^{r}+(x)^{r-1}\right)k[[x,\overline y,z]].
$$
$$
\overline y^{ad-bc}w-h(x\overline y^b)\equiv \overline y^{ad-bc}F_p-h(x\overline y^b)
\text{ mod } (x^r)
$$
implies 
$$
\overline y^{ad-bc}F_p-h(x\overline y^b)\in
\left(\hat{\cal I}_{C,p}^{r}+(x)^{r-1}\right)k[[x,\overline y,z]].
$$

The case when $ad-bc\le 0$ is similar.
\end{pf}

\begin{Lemma}\label{Lemma970}
Suppose that $p\in X$ is a 1 or 2 point, $D$ is a generic curve through $p$ on a component of
$E_X$ containing $p$. Then $F_q\not\in \hat{\cal I}_{D,q}$ for $q\in D$ (if $F_q$ is
computed with respect to permissible parameters $(x,y,z)$ at $q$ such that
$x=z=0$ are local equations of $D$ at $q$).
\end{Lemma}

\begin{pf}
By Lemma \ref{Lemma2} and Lemma \ref{Lemma7}, we need only check this at $p$.  When $p$ is a 1 point this
follows from Lemma \ref{Lemma301}. 

Suppose that $p$ is a 2 point, $(x,y,z)$ are
permissible parameters at $p$ such that $x\in\hat{\cal I}_{D,p}$. Then
$$
\begin{array}{ll}
u&=(x^ay^b)^m\\
v&=P(x^ay^b)+x^cy^dF_p
\end{array}
$$
There exists a series 
$$
\overline z = z -\sum_{i=1}^{\infty}\alpha_iy^i
$$
with $\alpha_i \in k$ such that
$\hat{\cal I}_{D,p}=(x,\overline z)$. Let 
$$
v=\overline P(x^ay^b)+x^cy^d\overline F_p
$$
 be the normalized form of $v$
with respect to the permissible parameters $(x,y,\overline z)$. Then
$$
\overline F_p=F_p-\sum b_i x^{\overline a_i}y^{\overline b_i}
$$
with $b_i\in k$ and
$$
a(d+\overline b_i)-b(c+\overline a_i)=0
$$
for all $i$. Suppose that $\overline F_p\in \hat{\cal I}_{D,p}$.
$$
F_p=h(y,z)+x\Omega
$$
with $h\ne 0$. We either have
$$
\overline z\mid h(y,\overline z+\sum\alpha_iy^i)
$$
or there exists $\overline c\in k$, $\overline d\in \bold N$ such that
$a(d+\overline d)-bc=0$ and
$$
\overline z\mid (h(y,\overline z+\sum \alpha_iy^i)-\overline cy^{\overline d}).
$$
Thus either $h(y,\sum\alpha_ky^k)=0$ or
$h(y,\sum \alpha_ky^k)=\overline cy^{\overline d}$.

Since $D$ is generic, we can suppose that $\alpha_1,\alpha_2$ are independent generic points of $k$.

Let $e=\nu(h)$. Write $h=\sum_{i+j\ge e}a_{ij}y^iz^j$.
$$
\begin{array}{ll}
h(y,\sum \alpha_ky^k)&=
\sum_{i+j=e}a_{ij}y^i(\alpha_1^jy^j+j\alpha_1^{j-1}\alpha_2y^{j+1})
+\sum_{i+j=e+1}a_{ij}y^i(\alpha_1 y)^j+ y^{e+2}\Omega\\
&=\left(\sum_{i+j=e}a_{ij}\alpha_1^j\right)y^e
+\left(\sum_{i+j=e}ja_{ij}\alpha_1^{j-1}\alpha_2+\sum_{i+j=e+1}a_{ij}\alpha_1^j\right)y^{e+1}
+y^{e+2}\Omega
\end{array}
$$
We must have 
$h(y,\sum \alpha_ky^k)=\overline c\overline y^{\overline d}$ and $e=\overline d$ since $\sum_{i+j=e}a_{ij}\alpha_1^j\ne 0$ as
$\alpha_1$ is a generic point of $k$. Since $F_p$ is normalized, we must have $a_{e0}=0$,
and $\nu(h)=e$ implies there exists $a_{i_0j_0}\ne 0$ such that $i_0+j_0=e$ and $j_0>0$.

We must have
$$
\left(\sum_{i+j=e} ja_{ij}\alpha_1^{j-1}\right)\alpha_2+\left(\sum_{i+j=e+1}a_{ij}\alpha_1^j
\right)=0.
$$ 
which is a contradiction to the assumption that $\alpha_1,\alpha_2$ are 
independent generic points of $k$.
\end{pf}

\begin{Definition}\label{Def1090}
Suppose that $\Phi_X:X\rightarrow S$ is weakly prepared. 

A monoidal transform $\pi:Y\rightarrow X$ is called weakly permissible if
$\pi$ is the blowup of a point $p$ on $E_X$, or a nonsingular curve $C$ on $E_X$
such that $C$ makes SNCs with $\overline B_2(X)$.

Suppose that $\Phi_X:X\rightarrow S$ is weakly prepared, and $\pi:Y\rightarrow X$ is a
weakly permissible monoidal transform. Define $\Phi_Y=\Phi_X\circ\pi:Y\rightarrow S$,
$E_Y=\pi^{-1}(E_X)_{\text{red}}$. $\Phi_Y$ is weakly prepared.
\end{Definition}

\section{The invariant $\nu$ under quadratic transforms}

Throughout this section we will suppose that $\Phi_X:X\rightarrow S$ is weakly prepared.

\begin{Theorem}\label{Theorem9} Suppose that $\nu(p)=r$, $\pi:X_1\rightarrow X$ is the 
blowup of $p$, $q\in\pi^{-1}(p)$ with $\nu(q)=r_1$.

{\bf Suppose that $p$ is a 1 point. Then}
\begin{enumerate}
\item If $q$ is a 1 point then $r_1\le r$.
\item If $q$ is a 2 point then $r_1\le r$. $r_1=r$ implies $\tau(q)>0$.
\end{enumerate}
{\bf Suppose that  $p$ is a 2 point.  Then}
\begin{enumerate}
\item If $q$ is a 1 point then $r_1\le r+1$. $r_1=r+1$ implies $\gamma(q)=r+1$.
\item If $q$ is a 2 point then $r_1\le r$.
\item If $q$ is a 3 point then $r_1\le r$.
\end{enumerate}
{\bf Suppose that $p$ is a 3 point. Then}
\begin{enumerate}
\item If $q$ is a 1 point then $r_1\le r+1$.
$r_1=r+1$ implies $\gamma(q)=r+1$.
\item If $q$ is a 2 point then $r_1\le r+1$. $r_1=r+1$ implies $\tau(q)>0$.
Furthermore there are permissible parameters $(x_1,y_1,z_1)$ at $q$ such that
$$
\begin{array}{ll}
u&=(x_1^ay_1^b)^m\\
v&=P(x_1^ay_1^b)+x_1^cy_1^dF_1
\end{array}
$$
and the leading form of $F_1$ is
$$
L_1=cy_1^tz_1^{r+1-t}+x_1\Omega
$$
where $0\le t\le r$, $cb-(d+t)a=0$.
In this case, the leading form of $F$ is 
$$
L=y^t(\sum_{i+k=r-t}b_{ik}x^iz^k)
$$
where all $b_{ik}\ne 0$, and there are regular parameters $(\overline x_1,\overline y_1,
\overline z_1)$ in $\hat{\cal O}_{X_1,q}$ such that  
$$
x=\overline x_1,
y=\overline x_1\overline y_1,
z=\overline x_1(\overline z_1+\beta)
$$
for some $0\ne \beta\in k$.
\item If $q$ is a 3 point then $r_1\le r$. If $r_1=r$, and $(x_1,y_1,z_1)$ are permissible parameters at $q$ with
$$
x=x_1, y=x_1y_1, z=x_1z_1,
$$
then the leading form of $F$ is
$$
L=L(y_1,z_1).
$$
\end{enumerate}
\end{Theorem}

\begin{pf}
{\bf Suppose that $p$ is a 1 point} 
$$
\begin{array}{ll}
u&=x^k\\
v&=P(x)+x^cF
\end{array}
$$
 Write $F = \sum_{i+j+k\ge r} a_{ijk}x^iy^jz^k$.

Suppose that $q\in\pi^{-1}(p)$ is a 1 point. Then there are permissible parameters $(x_1,y_1,z_1)$ at $q$  such that
$$
\begin{array}{ll}
x&=x_1\\
y&=x_1(y_1+\alpha)\\
z&=x_1(z_1+\beta)
\end{array}
$$
$$
\begin{array}{ll}
u&=x_1^k\\
v&=P(x_1)+x_1^{c+r}\frac{F}{x_1^r}
\end{array}
$$
\begin{equation}\label{eq17}
\frac{F}{x_1^r} = \sum_{j+k\le r}a_{jk}(y_1+\alpha)^j(z_1+\beta)^k+x_1\Omega
\end{equation}
where $a_{jk}=a_{r-i-j,j,k}$.
Suppose that $\nu(q)>\nu(p)$. Then
$$
\sum_{j+k\le r}a_{jk}(y_1+\alpha)^j(z_1+\beta)^k=\gamma\in k.
$$
and the leading form of $F$ is
$$
L = x_1^r(\sum_{j+k\le r}a_{jk}(y_1+\alpha)^j(z_1+\beta)^k) = \gamma x_1^r = \gamma x^r
$$
a contradiction to the assumption that $F$ is normalized.
Thus $\nu(q)\le\nu(p)$.

Now suppose that $q\in\pi^{-1}(p)$ is a 2 point. Then, after possibly interchanging $y$ and $z$,
 there are permissible parameters $(x_1,y_1,z_1)$ at $q$ such that
$$
\begin{array}{ll}
x&=x_1y_1\\
y&=y_1\\
z&=y_1(z_1+\alpha)
\end{array}
$$
$$
\begin{array}{ll}
u&= x_1^ky_1^k\\
v&=P(x_1y_1) + x_1^cy_1^{c+r}\frac{F}{y_1^r}
\end{array}
$$
$$
\frac{F}{y_1^r} = \sum_{i+k\le r}a_{ik}x_1^i(z_1+\alpha)^k+y_1\Omega
$$
where $a_{ik} = a_{i,r-i-k,k}$.
Suppose that $\nu(q)>\nu(p)$. Then
$$
\sum_{i+k\le r}a_{ik}x_1^i(z_1+\alpha)^k  = \sum \gamma_i x_1^{a_i}
$$
with $a_i+c=c+r$ for all $i$.  $a_i=r$ is the only solution to this equation, so that
if $L$ is the leading form of $F$,
$$
x^cL = \gamma x_1^{r+c}y_1^{r+c} = \gamma x^{r+c}
$$
a contradiction to the assumption that $F$ is normalized. Thus $\nu(q)\le \nu(p)$.

Suppose that $\nu(q)=r$. After making a permissible  change of parameters, we may assume that
$\alpha=0$.  We have 
\begin{equation}\label{eq18}
F_q = \sum_{i+k\le r}a_{ik}x_1^iz_1^k +y_1\Sigma
\end{equation}
Thus $a_{ik} = 0$ if $i+k<r$, and since $L$ is normalized, we must have $a_{ik}\ne 0$ for some $k>0$.
 Thus $\tau(q)>0$.

{\bf Suppose that $p$ is a 2 point} 
$$
\begin{array}{ll}
u&=(x^ay^b)^k\\
v&=P(x^ay^b)+x^cy^dF
\end{array}
$$
 Write $F = \sum_{i+j+k\ge r} a_{ijk}x^iy^jz^k$.

Suppose that $q\in\pi^{-1}(p)$ is a 1 point. Then there are regular parameters $(x_1,y_1,z_1)$ in $\hat{\cal O}_{X_1,q}$ such that
$$
\begin{array}{ll}
x&=x_1\\
y&=x_1(y_1+\alpha)\\
z&=x_1(z_1+\beta)
\end{array}
$$
with $\alpha\ne0$.
$$
u=x_1^{(a+b)k}(y_1+\alpha)^{bk}=\overline x_1^{(a+b)k}
$$
where $x_1=\overline x_1(y_1+\alpha)^{-\frac{b}{a+b}}$.
$$
v=P(\overline x_1^{a+b})+\overline x_1^{c+d+r}(y_1+\alpha)^\lambda\frac{F}{x_1^r}
$$
where $\lambda=d-\frac{b(c+d+r)}{a+b}$. 
\begin{equation}\label{eq19}
\frac{F}{x_1^r}=\sum_{j+k\le r}a_{jk}(y_1+\alpha)^j(z_1+\beta)^k+x_1\Omega
\end{equation}
where $a_{jk}=a_{r-i-j,j,k}$.
Suppose that 
$$
(y_1+\alpha)^{\lambda}\left(\sum_{j+k\le r} a_{jk}(y_1+\alpha)^j(z_1+\beta)^k\right)\equiv\gamma\text{ mod }(y_1,z_1)^{r+2}
$$
for some $\gamma\in k$.
Then 
\begin{equation}\label{eq15}
\sum_{j+k\le r} a_{jk}(y_1+\alpha)^j(z_1+\beta)^k \equiv\gamma(y_1+\alpha)^{-\lambda}\text{ mod }(y_1,z_1)^{r+2}
\end{equation}
 Set 
$$
f(y_1)=(y_1+\alpha)^{-\lambda}=\sum_{i=0}^{\infty} \alpha_iy_1^i
$$
where $\alpha_0=\alpha^{-\lambda}$ and
$$
\alpha_i = \frac{-\lambda(-\lambda-1)\cdots(-\lambda-i+1)}{i!}\alpha^{-\lambda-i}
$$
for $i\ge 1$.
(\ref{eq15}) implies $\alpha_{r+1}=0$, so that $-\lambda\in\{0,1,\ldots,r\}$, and
$$
\sum_{j+k\le r} a_{jk}(y_1+\alpha)^j(z_1+\beta)^k = \gamma(y_1+\alpha)^{-\lambda}
$$
Thus
$$
\sum_{i+j+k=r}a_{ijk}x_1^r(y_1+\alpha)^j(z_1+\beta)^k=\gamma x_1^{r+\lambda}x_1^{-\lambda}(y_1+\alpha)^{-\lambda}
$$
which implies that the leading form of $F$ is 
$$
L=\sum_{i+j+k=r}a_{ijk}x^iy^jz^k = \gamma x^{r+\lambda} y^{-\lambda}
$$
$$
x^cy^dL = \gamma x^{r+c+\lambda}y^{d-\lambda}
$$
$$
\begin{array}{ll}
a(d-\lambda)-b(r+c+\lambda)&=a\left[\frac{b(c+d+r)}{a+b}\right]-b\left[ r+c+d-\frac{b(c+d+r)}{a+b}\right]=0\\
\end{array}
$$
Thus $x^{r+c+\lambda}y^{d-\lambda}$ is a power of $x^ay^b$,
a contradiction to the assumption that $F$ is normalized.  We conclude that $\nu(q)\le \nu(p)+1$.

If $\nu(q) = \nu(p)+1$, we must then have that 
$$
F_1 = (y_1+\alpha)^{\lambda}\left(\sum_{j+k\le r} a_{jk}(y_1+\alpha)^j(z_1+\beta)^k\right)-\gamma + \overline x_1\Sigma
$$
with $\gamma=\alpha^{\lambda}\sum_{j+k\le r}a_{jk}\alpha^j\beta^k$. There is a 
 nonzero degree $r+1$ term in  $F_1(0,y_1,z_1)$, so that $\gamma(q)=r+1$.
 
Now suppose that $q\in\pi^{-1}(p)$ is a 2 point. Then after possibly interchanging $x$ and $y$,
there are permissible parameters $(x_1,y_1,z_1)$ at $q$ such that
$$
\begin{array}{ll}
x&=x_1\\
y&=x_1y_1\\
z&=x_1(z_1+\beta)
\end{array}
$$
with $\beta\ne 0$.
$$
\begin{array}{ll}
u&=(x_1^{a+b}y_1^{b})^k\\
v&=P(x_1^{a+b}y_1^{b})+x_1^{c+d+r}y_1^d\frac{F}{x_1^r}
\end{array}
$$
\begin{equation}\label{eq20}
\frac{F}{x_1^r} = \sum_{j+k\le r} a_{jk}y_1^j(z_1+\beta)^k+x_1\Omega
\end{equation}
with $a_{jk}=a_{r-j-k,j,k}$.
Suppose that 
$$
\sum_{j+k\le r} a_{jk}y_1^j(z_1+\beta)^k = \sum\gamma_i y_1^{t_i}
$$
with $\gamma_i\in k$, $(c+d+r)b-(a+b)(d+t_i) = 0$ for all $i$.
 There is at most one natural number $t=t_i$ which 
is a solution to this equation, which simplifies to 
\begin{equation}\label{eq16}
a(d+t)-b(c+r-t)=0.
\end{equation}
We have
$$
\sum_{j+k\le r} a_{jk}y_1^j(z_1+\beta)^k = \gamma y_1^{t}
$$
so that
$$
\sum_{i+j+k=r}a_{ijk}x_1^ry_1^j(z_1+\beta)^k = \gamma x_1^ry_1^t
$$
Thus $L = \gamma x^{r-t}y^t$. But by (\ref{eq16})
$x^{c+r-t}y^{d+t}$ is a power of $x^ay^b$, a contradiction to the assumption that $F$ is normalized.
Thus $\nu(q)\le r$.

Now suppose that $q\in\pi^{-1}(p)$ is a 3 point. Then there are regular parameters $(x_1,y_1,z_1)$ at $Q$ such that
$$
\begin{array}{ll}
x&=x_1z_1\\
y&=y_1z_1\\
z&=z_1
\end{array}
$$
We have 
\begin{equation}\label{eq21}
\begin{array}{ll}
u&=(x_1^ay_1^bz_1^{a+b})^k\\
v&=P(x_1^ay_1^bz_1^{a+b})+x_1^cy_1^dz_1^{c+d+r}\frac{F}{z_1^r}
\end{array}
\end{equation}
$F_q=\frac{F}{z_1^r}$, so that $\nu(q) = \nu(\frac{F}{z_1^r})\le r$.

{\bf Suppose that $p$ is a 3 point} 
$$
\begin{array}{ll}
u&=(x^ay^bz^c)^k\\
v&=P(x^ay^bz^c)+x^dy^ez^fF
\end{array}
$$
 Write $F = \sum_{i+j+k\ge r} a_{ijk}x^iy^jz^k$.
Suppose that $q\in\pi^{-1}(p)$ is a 1 point. Then there are regular parameters $(x_1,y_1,z_1)$ in $\hat{\cal O}_{X_1,q}$ such that
$$
\begin{array}{ll}
x&=x_1\\
y&=x_1(y_1+\alpha)\\
z&=x_1(z_1+\beta)
\end{array}
$$
with $\alpha,\beta\ne0$.
$$
u=x_1^{(a+b+c)k}(y_1+\alpha)^{bk}(z_1+\beta)^{ck} = \overline x_1^{(a+b+c)k}.
$$
where $\overline x_1$ is defined by
$$
x_1 =\overline x_1(y_1+\alpha)^{-\frac{b}{a+b+c}}(z_1+\beta)^{-\frac{c}{a+b+c}}.
$$
$$
\begin{array}{ll}
v&=P(\overline x_1^{a+b+c})+x_1^{d+e+f+r}(y_1+\alpha)^e(z_1+\beta)^f\frac{F}{x_1^r}\\
&=P(\overline x_1^{a+b+c})+\overline x_1^{d+e+f+r}(y_1+\alpha)^{\lambda_1}(z_1+\beta)^{\lambda_2}\frac{F}{x_1^r}
\end{array}
$$
where 
$$
\begin{array}{ll}
\lambda_1&=e-\frac{b(d+e+f+r)}{a+b+c}\\
\lambda_2&= f-\frac{c(d+e+f+r)}{a+b+c}
\end{array}
$$
\begin{equation}\label{eq22}
\frac{F}{x_1^r} = \sum_{j+k\le r}a_{jk}(y_1+\alpha)^j(z_1+\beta)^k+\overline x_1\Omega
\end{equation}
where $a_{jk}=a_{r-j-k,j,k}$. Suppose that
$$
(y_1+\alpha)^{\lambda_1}(z_1+\beta)^{\lambda_2}\left(\sum_{j+k\le r} a_{jk}(y_1+\alpha)^j(z_1+\beta)^k\right)
\equiv \gamma\text{ mod }(y_1,z_1)^{r+2}
$$
for some $\gamma\in k$. We first observe that we cannot have $\gamma=0$, for $\gamma=0$ implies
$$
 \sum_{j+k\le r} a_{jk}(y_1+\alpha)^j(z_1+\beta)^k
\equiv 0\text{ mod }(y_1,z_1)^{r+2}
$$
which implies
$$
\sum_{j+k\le r} a_{jk}(y_1+\alpha)^j(z_1+\beta)^k=0,
$$
a contradiction. Thus $\gamma\ne 0$. 
\begin{equation}\label{eq12}
\sum_{j+k\le r} a_{jk}(y_1+\alpha)^j(z_1+\beta)^k\equiv \gamma (y_1+\alpha)^{-\lambda_1}(z_1+\beta)^{-\lambda_2}
\text{ mod }(y_1,z_1)^{r+2}
\end{equation}
Set $f(y_1,z_1) = (y_1+\alpha)^{-\lambda_1}(z_1+\beta)^{-\lambda_2}$,
$$
\alpha_{ij} = \frac{1}{i!\,j!}\frac{\partial^{i+j}f}
{\partial y_1^i\partial z_1^j}(0,0)
$$
Then
$$
(y_1+\alpha)^{-\lambda_1}(z_1+\beta)^{-\lambda_2} = \sum \alpha_{ij}y_1^iz_1^j.
$$
$$
\alpha_{ij} = \left\{ \begin{array}{ll}
\left(\frac{-\lambda_1(-\lambda_1-1)\cdots(-\lambda_1-i+1)}{i!}\alpha^{-\lambda_1-i}\right)
\left(\frac{-\lambda_2(-\lambda_2-1)\cdots(-\lambda_2-j+1)}{j!}\beta^{-\lambda_2-j}\right)
&\text{ if }i,j>0\\
\alpha^{-\lambda_1}\left(\frac{-\lambda_2(-\lambda_2-1)\cdots(-\lambda_2-j+1)}{j!}\beta^{-\lambda_2-j}\right)
&\text{ if }i=0,j>0\\
\left(\frac{-\lambda_1(-\lambda_1-1)\cdots(-\lambda_1-i+1)}{i!}\alpha^{-\lambda_1-i}\right)\beta^{-\lambda_2}
&\text{ if }j=0, i>0\\
\alpha^{-\lambda_1}\beta^{-\lambda_2}&\text{ if }i=j=0
\end{array}\right.
$$
Thus $\alpha_{ij} = 0$ for $i+j=r+1$ by (\ref{eq12}), and
$-\lambda_1\in\{0,1,\ldots,r\}$, $-\lambda_2\in\{0,1,\ldots,r\}$ and $-\lambda_1-\lambda_2\le r$.
Thus 
$$
\sum_{j+k\le r} a_{jk}(y_1+\alpha)^j(z_1+\beta)^k = \gamma (y_1+\alpha)^{-\lambda_1}(z_1+\beta)^{-\lambda_2}
$$
so that
$$
\begin{array}{ll}
\sum_{i+j+k= r} a_{ijk}x_1^r(y_1+\alpha)^j(z_1+\beta)^k &= \gamma x_1^r(y_1+\alpha)^{-\lambda_1}(z_1+\beta)^{-\lambda_2}\\
&=\gamma x_1^{r+\lambda_1+\lambda_2}\left[x_1^{-\lambda_1}(y_1+\alpha)^{-\lambda_1}\right]
\left[x_1^{-\lambda_2}(z_1+\beta)^{-\lambda_2}\right]
\end{array}
$$
and the leading form of $F$ is
$$
L = \sum_{i+j+k= r} a_{ijk}x^iy^jz^k = \gamma x^{r+\lambda_1+\lambda_2}y^{-\lambda_1}z^{-\lambda_2}
$$
$$
x^dy^ez^fL=\gamma x^{d+r+\lambda_1+\lambda_2}y^{e-\lambda_1}z^{f-\lambda_2}
$$
Set
$$
\begin{array}{ll}
\underline a &= d+r+e-\frac{b(d+e+f+r)}{a+b+c}+f - \frac{c(d+e+f+r)}{a+b+c}\\
\underline b&= e - \left(e-\frac{b(d+e+f+r)}{a+b+c}\right)\\
\underline c&= f-\left(f-\frac{c(d+e+f+r)}{a+b+c}\right)
\end{array}
$$
Set
$\tau = \frac{d+e+f+r}{a+b+c}$.
$$
\begin{array}{ll}
\underline a &= \frac{(d+e+f+r)(a+b+c)-b(d+e+f+r)-c(d+e+f+r)}{a+b+c}=a\tau\\
\underline b &= b\tau\\
\underline c&= c\tau
\end{array}
$$
$$
\begin{array}{ll}
b\underline a-a\underline b &=(ba-ab)\tau=0\\
a\underline c-c \underline a &= (ac-ca)\tau = 0\\
c\underline b-b\underline c &=(cb-bc)\tau =0
\end{array}
$$
thus $\gamma=0$ since $F$ is normalized. This contradiction shows that 
$\nu(q)\le r+1$.

We have shown that
$$
F_1=(y_1+\alpha)^{\lambda_1}(z_1+\beta)^{\lambda_2}\left(\sum_{j+k\le r} a_{jk}(y_1+\alpha)^j(z_1+\beta)^k\right)-\gamma+\overline x_1\Sigma
$$
with
$$
\gamma=\alpha^{\lambda_1}\beta^{\lambda_2}\sum_{j+k\le r}a_{jk}\alpha^j\beta^k.
$$
Thus $r_1=r+1$ implies there is a nonzero degree $r+1$ term in $F_1(0,y_1,z_1)$ so that
$\gamma(q)=r+1$.

Now suppose that $q\in\pi^{-1}(p)$ is a 2 point. Then after possibly interchanging $x,y,z$,
 there are regular parameters $(x_1,y_1,z_1)$ at $q$ such that
$$
\begin{array}{ll}
x&=x_1\\
y&=x_1y_1\\
z&=x_1(z_1+\beta)
\end{array}
$$
with $\beta\ne0$.
$$
u=x_1^{(a+b+c)k}y_1^{bk}(z_1+\beta)^{ck}
$$
Set
$$
x_1=(z_1+\beta)^{\frac{-c}{a+b+c}}\overline x_1
$$
$$
\begin{array}{ll}
u&=(\overline x_1^{a+b+c}y_1^{b})^k\\
v&=P(\overline x_1^{a+b+c}y_1^b)+x_1^{d+e+f+r}y_1^e(z_1+\beta)^f\frac{F}{x_1^r}\\
&=P(\overline x_1^{a+b+c}y_1^b)+\overline x_1^{d+e+f+r}y_1^e(z_1+\beta)^{\lambda_1}\frac{F}{x_1^r}
\end{array}
$$
where 
\begin{equation}\label{eq14}
\lambda_1 = f - \frac{c(d+e+f+r)}{a+b+c}.
\end{equation}
\begin{equation}\label{eq23}
\frac{F}{x_1^r} = \sum_{j+k\le r}a_{jk}y_1^j(z_1+\beta)^k+\overline x_1\Omega
\end{equation}
where $a_{jk}=a_{r-i-j,j,k}$. 
There is at most one natural number $t$ such that 
\begin{equation}\label{eq13}
(d+e+f+r)b-(e+t)(a+b+c)=0.
\end{equation}

If $\nu(q)>r+1$ 
there exists a $t$ satisfying (\ref{eq13}), and $0\ne \gamma\in k$ such that 
$$
(z_1+\beta)^{\lambda_1}\left(\sum_{j+k\le r} a_{jk}y_1^j(z_1+\beta)^k\right)\equiv  \gamma y_1^t\text{ mod }
(y_1,z_1)^{r+2}.
$$
Thus 
$$
\sum_{j+k\le r} a_{jk}y_1^j(z_1+\beta)^k\equiv \gamma(z_1+\beta)^{-\lambda_1}y_1^t\text{ mod }(y_1,z_1)^{r+2}.
$$
Set 
$$
\tau_j = \frac{-\lambda_1(-\lambda_1-1)\cdots(-\lambda_1-j+1)}{j!}\beta^{-\lambda_1-j}.
$$
$$
\sum_{j+k\le r} a_{jk}y_1^j(z_1+\beta)^k\equiv 
\gamma y_1^t(\sum_{j=0}^{\infty} \tau_jz_1^j)\text{ mod }(y_1,z_1)^{r+2}.
$$
implies
$$
0 = \tau_{r+1-t} = \frac{-\lambda_1(-\lambda_1-1)\cdots(-\lambda_1-r+t)}{(r+1-t)!}\beta^{-\lambda_1-(r+1-t)}
$$
so that $-\lambda_1\in \{0,1,\ldots,r-t\}$ and $t\le r$.
Thus
$$
\sum_{j+k\le r} a_{jk}y_1^j(z_1+\beta)^k=\gamma(z_1+\beta)^{-\lambda_1}y_1^t
$$
$$
\begin{array}{ll}
\sum_{i+j+k=r}a_{ijk}x_1^ry_1^j(z_1+\beta)^k&=\gamma x_1^r(z_1+\beta)^{-\lambda_1}y_1^t\\
&=\gamma x_1^{r-t+\lambda_1}\left[x_1^ty_1^t\right]\left[x_1^{-\lambda_1}(z_1+\beta)^{-\lambda_1}\right]
\end{array}
$$
$$
x^dy^ez^f L = \gamma x^{r-t+\lambda_1+d}y^{t+e}z^{f-\lambda_1}
$$
where $L$ is the leading form of $F$.
Set
$$
\begin{array}{ll}
\underline a &= r-t+\lambda_1+d\\
\underline b &= t+e\\
\underline c &= f-\lambda_1
\end{array}
$$
We have the relations (\ref{eq13}) and (\ref{eq14}).
(\ref{eq13}) implies
$$
t=\frac{(d+e+f+r)b-e(a+b+c)}{a+b+c}.
$$
$$
\begin{array}{ll}
\underline a&=d+r-t+\lambda_1=\frac{a(d+e+f+r)}{a+b+c}\\
\underline b &= t+e=\frac{(d+e+f+r)b}{a+b+c}\\
\underline c &=\frac{c(d+e+f+r)}{a+b+c}
\end{array}
$$ 
$$
0=b\overline a-a\overline b=c\overline a-a\overline c=c\overline b-b\overline c
$$
Thus 
$$
x^{r-t+\lambda_1+d}y^{t+e}z^{f-\lambda_1}=(x^ay^bz^c)^m
$$
for some $m\in \bold N$, a contradiction, since $F$ is normalized. Thus $\nu(q)\le r+1$.

Suppose that $\nu(q)=r+1$.
$$
F_q=\sum_{j\le r}\left[\sum_{k\le r-j}a_{jk}(z_1+\alpha)^{k+\lambda_1}\right]y_1^j-\gamma y_1^t+\overline x_1\Sigma
$$

with $t\le r$, $\gamma\in k$ implies $a_{jk}=0$ if $j\ne t$. Thus
$$
F_q = y_1^t(\sum_{k\le r-t}a_{tk}(z_1+\alpha)^{k+\lambda_1}-\gamma)+\overline x_1\Sigma
$$

implies
$$
L_q = cy_1^tz_1^{r+1-t}+\overline x_1\Omega
$$
where $0\ne c\in k$. 

Now suppose that $q\in\pi^{-1}(p)$ is a 3 point. After possibly permuting $x,y,z$,
 there are permissible parameters $(x_1,y_1,z_1)$ at $q$ such that
$$
\begin{array}{ll}
x&=x_1\\
y&=x_1y_1\\
z&=x_1z_1
\end{array}
$$
$$
\begin{array}{ll}
u&=(x_1^{a+b+c}y_1^{b}z_1^{c})^k\\
v&=P(x_1^{a+b+c}y_1^bz_1^c)+x_1^{d+e+f+r}y_1^ez_1^f\frac{F}{x_1^r}
\end{array}
$$
$\nu(q)=\nu(\frac{F}{x_1^r})\le r$.
\end{pf}

\begin{Example}\label{Example2003}
$\nu(p)$ can go up by 1 after a quadratic transform.
We can construct the example as follows.
\end{Example}

$$
u=xy, v=x^2y
$$
has $F=1$.
blowup $p$ and consider the point $p_1$ above $p$ with regular parameters $(x_1,y_1,z_1)$ defined by
$x=x_1, y=x_1(y_1+\alpha)$, $\alpha\ne 0$, $z=x_1z_1$.
Set $\overline x_1 = x_1(y_1+\alpha)^{\frac{1}{2}}$,
$\overline y_1=(y_1+\alpha)^{-\frac{1}{2}}-\alpha^{-\frac{1}{2}}$.
Then
$$
u=\overline x_1^2, v=\alpha^{-\frac{1}{2}}\overline x_1^3
+\overline x_1^3\overline y_1,
$$
so that $F_1=\overline y_1$.

\begin{Theorem}\label{Theorem13}
Suppose that $\nu(p)=r$, $\pi:X_1\rightarrow X$ is the blowup of $p$, $q\in\pi^{-1}(p)$ with
$r_1=\nu(q)$.

If $p$ is a 1 point  then
\begin{enumerate}
\item If $q$ is a 1 point, then $r_1<r$ if $\tau(p)<r$, and if $r_1=r$ then $\tau(q)=r$.
\item If $q$ is a 1 point, then $\gamma(q)\le r$.
\item If $q$ is a 2 point, and $r_1=r$ then $\tau(p)\le \tau(q)$.
\item If $q$ is a 2 point and $\gamma(p)=r$, then $\gamma(q)\le r$. 
\end{enumerate}

If $p$ is a 2 point and $1\le \tau(p)$ then
\begin{enumerate}
\item If $q$ is a 1 point then $r_1\le r$ and $\gamma(q)\le r$.

\item If $q$ is a 2 point and $r_1 =r$, then $\tau(p)\le\tau(q)$.
\item If $q$ is a 2 point and $\gamma(p)=r$, then $\gamma(q)\le r$.
\item If $q$ is a 3 point then $r_1\le r-\tau(p)$
\end{enumerate}
\end{Theorem}

\begin{pf}
{\bf Suppose that $p$ is a 1 point with $\gamma(p)=r$.}
Suppose that $q\in\pi^{-1}(p)$ is a  1 point, and 
$r_1 =r$. After making a permissible change of parameters we can
assume that
$x = x_1, y=x_1y_1, z=x_1z_1$. We than have, with the notation of (\ref{eq17}).
$$
F_1 = \sum_{j+k\le r}a_{r-j-k,j,k}y_1^jz_1^k+x_1\Omega.
$$

Now suppose that 
$q\in\pi^{-1}(p)$ is a  2 point, and 
$r_1 =r$. After making a permissible change of parameters we can
assume that
$x = x_1y_1, y=y_1, z=y_1z_1$. We than have, with the notation of (\ref{eq18})
$$
F_1 = \sum_{i+k\le r}a_{i,r-i-k,k}x_1^iz_1^k+y_1\Sigma
$$

{\bf Suppose that $p$ is a 2 point and $1\le\tau(p)$.}
Suppose that $q\in\pi^{-1}(p)$ is a  1 point. After making a permissible change of parameters, we have
$x=x_1, y=x_1(y_1+\alpha), z=x_1z_1$
with $\alpha \ne 0$. We then have, with the notation of (\ref{eq19}), 
$$
F_1 = \sum_{j+k\le r}a_{r-j-k,j,k}(y_1+\alpha)^j(y_1+\alpha)^{\lambda}z_1^k - \gamma + \overline x_1\Omega.
$$
There exists $a_{ijk}$ with $i+j+k=r$ and $k=\tau(p)\ge 1$ such that $a_{ijk}\ne0$.
Thus $r_1\le r$ and $\gamma(q)\le r$.

Now suppose that $q\in\pi^{-1}(p)$ is a 2 point. After making a permissible change of parameters, we have
$x=x_1, y=x_1y_1, z=x_1z_1$
We then have, with the notation of (\ref{eq20}), 
$$
F_1 = \sum_{j+k\le r}a_{r-j-k,j,k}y_1^jz_1^k  + x_1\Omega.
$$
and there exist $i,j,k$ such that $i+j+k=r$ and $a_{ijk}\ne0$ with $k=\tau(p)$.
Thus if $r_1=r$, we have $\tau(p)\le\tau(q)$. If $\gamma(p)=r$, we have $\gamma(q)\le r$.

Now suppose that $q\in\pi^{-1}(p)$ is a 3 point. 
Then $x=x_1z_1, y=y_1z_1, z=z_1$
We then have, with the notation of (\ref{eq21}), 
$$
F_1 = \sum_{i+j\le r}a_{i,j,r-i-j}x_1^iy_1^j  + z_1\Omega.
$$
There exists $a_{ijk}$ with $i+j+k=r$ and $k=\tau(p)\ge1$ such that $a_{ijk}\ne0$.
Thus $r_1\le r-\tau(p)$. 
\end{pf}

\begin{Lemma}\label{Lemma10}
Suppose that $r\ge 2$ and $p \in X$ is a 1 point. Suppose that 
$(x,y,z)$ are permissible parameters at $p$ and $C\subset \overline S_r(X)$ is a curve such that
$p\in C$. 
Then $F_p\in \hat{\cal I}_{C,p}^{r}+(x^{r-1})$.  
\end{Lemma}

\begin{pf}
$x\in \hat{\cal I}_{C,p}$ by Lemma \ref{Lemma657}. There exist permissible
parameters $(x,\overline y,\overline z)$ at $p$ such that $\overline y,\overline z\in 
{\cal O}_{X,p}$, $\overline y=y+h_1$, $\overline z=z+h_2$ with $h_1,h_2\in m^r$.

Suppose that the conclusions of the Lemma are true for the parameters $(x,\overline y,
\overline z)$.
$$
\begin{array}{ll}
u&=x^a\\
v&=\overline P(x)+x^b\overline F(x,\overline y,\overline z)
\end{array}
$$
and $\overline F(x,\overline y,\overline z)$ is normalized with respect to the permissible
parameters $(x,\overline y,\overline z)$.
We have an expression
$$
\begin{array}{ll}
u&=x^a\\
v&=P(x)+x^bF(x,y,z)
\end{array}
$$
where $F(x,y,z)$ is normalized with respect to the permissible parameters $(x,y,z)$.
$$
\overline F(x,\overline y,\overline z)=F(x,y,z)+\Omega
$$
with $\Omega$ a series in $x$. Since $F(x,y,z)$ is normalized, 
$\nu(\overline F)=\nu(F)=r$ and only powers of $x$ of
order $\ge r$ can be removed from $\overline F(x,\overline y,\overline z)$ to
normalize to obtain $F(x,y,z)$. Thus the conclusions of the Lemma hold for $(x,y,z)$

We may thus assume that $y,z\in {\cal O}_{X,p}$.

There exists  an \'etale neighborhood $U$ of $p$ such that $(x,y,z)$ are
uniformizing parameters  in $U$, $x=0$ is a local equation of $E_X\cap U$, $C\cap U$
is a complete intersection. Let $R = \Gamma(U,{\cal O}_U)$, $I_C=\Gamma(U,{\cal I}_C)$. Set
$$
w=\frac{v-P_t(x)}{x^b}
$$
where  $t>b+r.$ 
 Thus 
\begin{equation}\label{eq609}
w\in (y,z,x^r){\cal O}_U,p
\end{equation} and 
\begin{equation}\label{29}
w-F\in (x^r)\hat{\cal O}_{X,p}.
\end{equation}
Let $q$ be a smooth point of $C\cap U$. Then there exists $\alpha,\beta\in k$
such that $(x,y-\alpha,z-\beta)$ are permissible parameters at $q$.
Lemma \ref{Lemma5} implies 
$$
F_q\in \hat{\cal I}_{C,q}^{r}+(x^{r-1}).
$$
$$
F_q=w-\sum_{i=0}^{\infty} \frac{\partial^i w}{\partial x^i}(0,\alpha,\beta)x^i
$$
Set
$$
\Lambda = w - \sum_{i<r-1}\frac{\partial^i w}{\partial x^i}(0,\alpha,\beta)x^i.
$$
$\Lambda\in ({\cal I}_{C,q}^r+(x^{r-1}))\hat{\cal O}_{X,q}$ implies by 
Theorem \ref{Theorem9}, Chapter VIII, section 4 \cite{ZS} and Lemma \ref{Lemma961}, 
$$
\Lambda \in ({\cal I}_{C,q}^{r}+(x^{r-1}))\cap R = (I_C^{r}+(x^{r-1})R.
$$
$\nu(p)=r$ and $t>b+r$ implies $\nu(w)\ge r$ and $\frac{\partial^i w}{\partial x^i}(0,\alpha,\beta)=0$ if $i<r-1$, so that $w\in \hat{\cal I}_{C,p}^r+(x^{r-1})$ which 
implies that $F_p\in\hat{\cal I}_{C,p}^r+(x^{r-1})$.
\end{pf}

\begin{Lemma}\label{Lemma11}
Suppose that $p\in X$ is a 1 point and $\nu(p) = \gamma(p)=r\ge 2$.
Then there exists at most one curve $C$ in $\overline S_r(X)$ containing $p$. If $C$ exists then it is nonsingular at $p$.
\end{Lemma}

\begin{pf} Suppose that $(x,y,z)$ are permissible parameters at $p$.
Write $\hat{\cal I}_{C,p}=(x,f(y,z))$. By Lemma \ref{Lemma10},
$F_p\in\hat{\cal I}_{C,p}^{r}+(x^{r-1})$.
$f^r\mid F_p(0,y,z)$ and $\gamma(p)=r$
 implies $\nu(f)=1$ and $C$ is nonsingular at $p$.

Suppose that $D\subset \overline S_r(X)$ is another curve containing $p$. Then $D$
is nonsingular at $p$. Lemma \ref{Lemma5} implies there exist permissible
parameters $(x,y,z)$ at $p$ such that $\hat{\cal I}_{C,p} = (x,z)$, there exist
series $a,b_{ij}$ such that
$$
F_p=x^{r-1}a+\sum_{i+j=r}b_{ij}x^iz^j.
$$
$b_{0r}$ is a unit implies
$$
\frac{\partial^{r-1}F_p}{\partial z^{r-1}}=x\phi+z\psi
$$
where $\psi$ is a unit. Since $F_p\in\hat{\cal I}_{D,p}^{r}+(x^{r-1})$, we have
$$
\frac{\partial^{r-1} F_p}{\partial z^{r-1}}\in \hat{\cal I}_{D,p}
$$
which implies $z\in\hat{\cal I}_{D,p}$, so that $C=D$.
\end{pf}

\begin{Lemma}\label{Lemma3} Suppose that $r\ge 2$, $p$ is a 2 point and 
$C\subset \overline S_r(X)$ is an
irreducible curve containing a 1 point, such that $p\in  C$.  Then $\nu(p)\ge r-1$.
 If $\tau(p)>0$, then $\nu(p)\ge r$.
\end{Lemma}

\begin{pf} First suppose that $C$ is nonsingular at $p$ and is transversal at $p$ to the
2 curve through $p$.
Then the result follows from Lemma \ref{Lemma7}.
Now suppose that $C$ does not make SNCs with the 2 curve through $p$.  Let $s=\nu(p)$. There exists a
sequence of quadratic transforms $\pi:X_1\rightarrow X$ centered at 2 and 3 points
such that the strict transform $C'$ of $C$ makes
SNCs with  $\overline B_2(X_1)$ at a 2 point $p_1 = C'\cap \pi^{-1}(p)$.  We have $s_1=\nu(p_1)\le s+1$,
by Theorems \ref{Theorem9} and \ref{Theorem13}.
$s_1=s+1$ implies $\tau(p_1)>0$ and $\tau(p)=0$. $\tau(p)>0$ 
implies $s_1\le s$ and if we further have $s_1=s$,
then   $\tau(p_1)>0$.

First supppose that $\tau(p)>0$. If $s_1=s$ then $\tau(p_1)>0$ so that $r\le s_1=s$. If  $s_1<s$ then
 $s>s_1\ge r-1$,  which implies $s\ge r$.

Now suppose that $\tau(p)=0$. If $s_1\le s$ then $s\ge s_1\ge r-1$.
If $s_1=s+1$ then $s+1=s_1\ge r$ since $\tau(p_1)>0$, so that $s\ge r-1$.
\end{pf}

\begin{Lemma}\label{Lemma4} Suppose that $r\ge 2$, $p$ is a 3 point and 
$C\subset \overline S_r(X)$ is an
irreducible curve containing a 1 point such that $p\in  C$. Then $\nu(p)\ge r-1$.
\end{Lemma}

\begin{pf} 
 Let $s=\nu(p)$. There exists a
sequence of quadratic transforms $\pi:X_1\rightarrow X$ centered at 2 and 3 points 
 such that the strict transform $C'$ of $C$ makes
SNCs with $\overline B_2(X_1)$ at the 2 point $p_1 = C'\cap \pi^{-1}(p)$.  We have $s_1=\nu(p_1)\le s+1$
by Theorems \ref{Theorem9} and \ref{Theorem13}.

If $s_1=s+1$ then $\tau(p_1)>0$, so that $s_1\ge r$ by Lemma \ref{Lemma3}, so that $s\ge r-1$.
If $s_1\le s$, then $s\ge s_1\ge r-1$ by Lemma \ref{Lemma3}.
\end{pf}

\begin{Theorem}\label{Theorem12}
 Suppose that $p \in X$ has $\nu(p)=r\ge 1$, and $(x,y,z)$ are permissible parameters
at $p$, $\pi:X_1\rightarrow X$ is the blow up of $p$.

{\bf Suppose that $p$ is a 3 point}
\begin{enumerate}
\item Suppose that the leading form $L_p=L(x,y,z)$ depends on $x$, $y$ and $z$.
Then there are no curves $C$ in $\pi^{-1}(p)\cap \overline S_{r+1}(X_1)$.
No 2 curves $C$ of $\pi^{-1}(p)$ satisfy $F_q\in\hat{\cal I}_{C,q}^r$ for $q\in C$.
\item Suppose that $L_p=L(x,y)$ depends on $x$ and $y$. Then the curves in $\pi^{-1}(p)\cap \overline S_{r+1}(X_1)$ are
a finite union of lines passing through a single 3 point of $\pi^{-1}(p)$. No 2 curves $C$ of $\pi^{-1}(p)$ satisfy $F_q\in\hat{\cal I}_{C,q}^r$ for $q\in C$.
\item
Suppose that $L_p=L(x)$ depends on $x$.
Then there are no 1 points in  $\pi^{-1}(p)\cap \overline S_{r+1}(X_1)$ and
 there is at most one curve $C$ in $\pi^{-1}(p)\cap \overline S_{r+1}(X_1)$.
It is the 2 curve $D$ which is the intersection of the strict transform of $x=0$ with $\pi^{-1}(p)$. $D$ is the only 2 curve $C$ in $\pi^{-1}(p)$ such that $F_q\in\hat{\cal I}_{C,q}^r$ for $q\in C$.
\end{enumerate}
{\bf Suppose that $p$ is a 2 point.}
 Then the curves in  $\pi^{-1}(p)\cap \overline S_{r+1}(X_1)$ are a finite union
of lines passing through the 3 point. There are no 2 curves in 
$\pi^{-1}(p)\cap \overline S_{r+1}(X_1)$.
\end{Theorem}

\begin{pf}
First suppose that $p$ is a 3 point and $L_p=L(x,y,z)$ depends
on $x,y$ and $z$. There are no 3 points in $\pi^{-1}(p)\cap \overline{S_{r+1}(X_1)}$
by Lemma \ref{Lemma4} and Lemma \ref{Lemma8} since (by direct calculation) $\nu(q)\le r-1$ at all 3 points in $\pi^{-1}(p)$.
There are no 2 curves in $\pi^{-1}(p)\cap \overline{S_{r+1}(X_1)}$ 
and there are no 2 curves in $\pi^{-1}(p)$ such that $F_q\in\hat{\cal I}_{C,q}^r$ for $q\in C$ by Lemma \ref{Lemma8}.
We will now show that there are no curves in $\overline S_{r+1}(X_1)\cap \pi^{-1}(p)$.
Suppose that there is a curve $C$ in  ${\overline{S_{r+1}(X_1)}}\cap \pi^{-1}(p)$
containing a 1 point. $C$ must contain a 2 point $q$.
$\nu(q) = r$ or $r+1$ by Lemma \ref{Lemma3} and Theorem \ref{Theorem9}.

First suppose that $\nu(q)=r+1$. Then by Theorem \ref{Theorem9}, there exist
permissible parameters $(x,y,z)$ at $p$ such that 
\begin{equation}\label{eq25}
L = y^tf(x,z)
\end{equation}
 for some $t$ with $0< t< r$, (since $L$ depends on $x$, $y$, and $z$). Write
$$
f = \sum_{i+k=r-t}b_{ik}x^iz^k
$$
 At a 1 point of $C$ we have (with the notation of (\ref{eq22})) permissible
parameters $(x_1,y_1,z_1)$ such that 
$x=x_1, y=x_1(y_1+\alpha), z=x_1(z_1+\beta)$ with $\alpha, \beta\ne 0$

\begin{equation}\label{eq24}
(y_1+\alpha)^{\lambda_1}(z_1+\beta)^{\lambda_2}
\left[\sum_{i+j+k=r}a_{ijk}(y_1+\alpha)^j(z_1+\beta)^k\right]
\equiv c_{\alpha,\beta}\text{ mod } (y_1,z_1)^{r+1}
\end{equation}
for some $0\ne c_{\alpha,\beta}\in k$.
Substituting (\ref{eq25}), we have
$$
\begin{array}{ll}
(y_1+\alpha)^{\lambda_1}(z_1+\beta)^{\lambda_2}\left[(y_1+\alpha)^t\sum_{i+k=r-t}b_{ik}(z_1+\beta)^k\right]
&\equiv c_{\alpha,\beta}\text{ mod } (y_1,z_1)^{r+1}\\
(y_1+\alpha)^{t+\lambda_1}\left[\sum_{i+k=r-t}b_{ik}(z_1+\beta)^k(z_1+\beta)^{\lambda_2}\right]
&\equiv c_{\alpha,\beta}\text{ mod } (y_1,z_1)^{r+1}
\end{array}
$$
If $t\ne -\lambda_1$ this is a contradiction, since there is then a nonzero $y_1z_1^s$ term
for some $0\le s\le r-t$. 
Thus $-\lambda_1=t\in\{1, \ldots, r-1\}$. But 
$$
\nu(\sum_{i+k=r-t}b_{ik}(z_1+\beta)^k(z_1+\beta)^{\lambda_2}
-c_{\alpha,\beta})\le r-t+1\le r 
$$
which is contradiction.

Now suppose that $\nu(q)=r$. 
We have from (\ref{eq23}) (in the  proof of Theorem \ref{Theorem9}) that there exist permissible parameters $(x_1,y_1,z_1)$ at $q$ such that
$$
x=x_1,y=x_1y_1,z=x_1(z_1+\beta),
$$
\begin{equation}\label{eq26}
F_q = \left[\sum_{j+k\le r}a_{r-j-k,j,k}y_1^j(z_1+\beta)^k\right](z_1+\beta)^{\lambda_1}-\gamma y_1^t+\overline x_1\Omega 
\end{equation}
if there exists a natural number $t$ such that $(d+e+f+r)b-(e+t)(a+b+c)=0$, and 
\begin{equation}
F_q = \left[\sum_{j+k\le r}a_{r-j-k,j,k}y_1^j(z_1+\beta)^k\right](z_1+\beta)^{\lambda_1}
+\overline x_1\Omega 
\end{equation}
otherwise.
Since $\nu(q)=r$, 
For fixed $j\ne t$, we have
$$
\nu(\sum_{k\le r-j}a_{r-j-k,j,k}(z_1+\beta)^k)\ge r-j
$$
Thus (for fixed $j\ne t$)
$$
\sum_{k\le r-j}a_{r-j-k,j,k}(z_1+\beta)^k=\gamma_jz_1^{r-j}
$$
for some $\gamma_j\in k$ and (for fixed $j\ne t$)
$$
\begin{array}{ll}
\sum_{i+k= r-j}a_{ijk}x^iy^jz^k &= y^j\left[\sum_{k\le r-j}a_{r-j-k,j,k}[x_1^k(z_1+\beta)^k]x_1^{r-j-k}\right]\\
&= y^jx_1^{r-j}\left[\sum_{k\le r-j}a_{r-j-k,j,k}(z_1+\beta)^k\right]\\
&=\gamma_j x_1^{r-j}z_1^{r-j}y^j\\ 
&=\gamma_jx^{r-j}\left(\frac{z}{x}-\beta\right)^{r-j}y^j\\
&=\gamma_j(z-\beta x)^{r-j}y^j
\end{array}
$$
For $j=t$, we have
$$
\nu(\sum_{k\le r-t}a_{r-t-k,t,k}(z_1+\beta)^{k+\lambda_1}-\gamma)\ge r-t
$$
Thus we either have
$$
L_p = \sum_{j\ne t}\gamma_jy^j(z-\beta x)^{r-j}+y^tf(x,z)
$$
where $(d+e+f+r)b-(e+t)(a+b+c)=0$, and $f$ is homogeneous of degree $r-t$,
or
$$
L_p = \sum_{j}\gamma_jy^j(z-\beta x)^{r-j}.
$$
Thus
$$
L_q = \sum\gamma_jy_1^jz_1^{r-j}+x_1\Omega_1.
$$
or
$$
L_q = \sum_{j\ne t}\gamma_jy_1^jz_1^{r-j}+\gamma_ty_1^tz_1^{r-t}+x_1\Omega_1.
$$
for some $\gamma_t\in k$.
If some $\gamma_j\ne 0$ with $j\ne r$, 
$\tau(q)>0$, so that $\nu(q)\ge r+1$ by Lemma \ref{Lemma3}, 
 a contradiction. 

The remaining case is 
\begin{equation}\label{eq603}
L_p = \gamma_r y^r+y^tf(x,z),
\end{equation}
 with $f\ne 0$, $t<r$ and $\gamma_r\ne 0$ if $t=0$, since $L_p$ depends on $x,y$ and $z$.
$$
f = \sum_{i+k=r-t}b_{ik}x^iz^k
$$
At a 1 point of $C$ we have (with the notation of (\ref{eq22})) regular parameters
$(x_1,y_1,z_1)$ such that 
$x=x_1, y=x_1(y_1+\alpha), z=x_1(z_1+\beta)$ with $\alpha, \beta\ne 0$ 
\begin{equation}\label{eq29}
(y_1+\alpha)^{\lambda_1}(z_1+\beta)^{\lambda_2}\left[\sum_{i+j+k=r}a_{ijk}(y_1+\alpha)^j(z_1+\beta)^k\right]
\equiv c_{\alpha,\beta}\text{ mod } (y_1,z_1)^{r+1}
\end{equation}
for some $c_{\alpha,\beta}\in k$.
Substituting (\ref{eq603}), we have
$$
\begin{array}{l}
(y_1+\alpha)^{\lambda_1}(z_1+\beta)^{\lambda_2}
\left[\gamma_r(y_1+\alpha)^r+(y_1+\alpha)^t\sum_{i+k=r-t}b_{ik}(z_1+\beta)^k
\right]
\equiv c_{\alpha,\beta}\text{ mod } (y_1,z_1)^{r+1}\\
\end{array}
$$
$$
\gamma_r(y_1+\alpha)^{r-t}+\sum_{i+k=r-t}b_{ik}(z_1+\beta)^k
\equiv (y_1+\alpha)^{-\lambda_1-t}(z_1+\beta)^{-\lambda_2}c_{\alpha,\beta}
\text{ mod }(y_1,z_1)^{r+1}
$$
The LHS of the last equation has no $y_1z_1$ term which implies $\lambda_1=-t$ or $-\lambda_2=0$.
$\lambda_1=-t$ implies $\gamma_r=0$ and $t>0$,
 $$
\nu(\sum_{i+k=r-t}b_{ik}(z_1+\beta)^k-c_{\alpha,\beta}(z_1+\beta)^{-\lambda_2}
)\le r-t+1\le r 
$$
which is contradiction. $\lambda_2=0$ implies $f=0$, a contradiction.

Now suppose that $p$ is a 3 point and $L_p=L(x,y)$.  Suppose that  $q\in\pi^{-1}(p)$ is a 1 point and $\nu(q)=r+1$.
$\hat{\cal O}_{X_1,q}$ has regular parameters $x_1$, $y_1$, $z_1$ such that 
\begin{equation}\label{eq31}
\begin{array}{ll}
x&=x_1\\
y&=x_1(y_1+\alpha)\\
z&=x_1(z_1+\beta)
\end{array}
\end{equation}
where $\alpha,\beta\ne 0$. Write
$$
L(x,y)=\sum_{i+j=r}a_{ij}x^iy^j.
$$
 (with the notation of (\ref{eq22})) 
\begin{equation}\label{eq604}
(y_1+\alpha)^{\lambda_1}(z_1+\beta)^{\lambda_2}\left[\sum a_{ij}(y_1+\alpha)^j\right] \equiv c_{\alpha,\beta}
\text{ mod }(y_1,z_1)^{r+1}
\end{equation}
which implies 
\begin{equation}\label{eq30}
\sum_{i+j=r} a_{ij}(y_1+\alpha)^j \equiv (y_1+\alpha)^{-\lambda_1}(z_1+\beta)^{-\lambda_2}c_{\alpha,\beta}
\text{ mod }(y_1,z_1)^{r+1}
\end{equation}
so that $\lambda_2 = 0$, and 
\begin{equation}\label{eq953}
\sum_{i+j=r}a_{ij}(y_1+\alpha)^j\equiv (y_1+\alpha)^{-\lambda_1}c_{\alpha}\text{ mod }(y_1)^{r+1}
\end{equation}

We will now show that there exist at most finitely many values of $\alpha$ such that an equation (\ref{eq953}) holds. Set
$$
g(t) = \sum_{i\le r}a_it^i
$$
where $a_i=a_{r-i,i}$. 

Suppose there are infinitely many values of $\alpha$ such that (\ref{eq953}) holds for some $q\in\pi^{-1}(p)$ with value $\beta$
and regular parameters $x_1,y_1,z_1$ in $\hat{\cal O}_{X_1,q}$ 
as in  (\ref{eq31}).  Define $g_{\alpha}$ by
$$
g_{\alpha}(y_1) = g(y_1+\alpha) = g(\frac{y}{x}).
$$
Set $\lambda=\lambda_1$. By assumption, 
\begin{equation}\label{eq32}
g_{\alpha}(y_1) = \sum_{i\le r}a_i(y_1+\alpha)^i\equiv c_{\alpha}(y_1+\alpha)^{-\lambda}\text{ mod } y_1^{r+1}
\end{equation}
We can expand the RHS of (\ref{eq32}) as
$$
\begin{array}{ll}
c_{\alpha}(y_1+\alpha)^{-\lambda}
&=c_{\alpha}\alpha^{-\lambda}+c_{\alpha}(-\lambda)\alpha^{-\lambda-1}y_1\\
&+c_{\alpha}\frac{-\lambda(-\lambda-1)}{2}\alpha^{-\lambda-2}y_1^2+\cdots\\
&+c_{\alpha}\frac{-\lambda(-\lambda-1)\cdots(-\lambda-r+1)}{r!}\alpha^{-\lambda-r}y_1^r+\cdots
\end{array}
$$
We can expand the LHS of (\ref{eq32}) as
$$
\begin{array}{ll}
g_{\alpha}(y_1) &= g_{\alpha}(0)+\frac{dg_{\alpha}}{dy_1}(0)y_1+\frac{1}{2}\frac{d^2g_{\alpha}}{dy_1^2}(0)y_1^2+
\cdots+\frac{1}{r!}\frac{d^rg_{\alpha}}{dy_1^r}(0)y_1^r\\
&= g(\alpha)+\frac{dg}{dt}(\alpha)y_1+\frac{1}{2}\frac{d^2g}{dt^2}(\alpha)y_1^2+
\cdots+\frac{1}{r!}\frac{d^rg}{dt^r}(\alpha)y_1^r
\end{array}
$$
We get that
$$
r!a_r=\frac{d^rg}{dt^r}(\alpha)=c_{\alpha}(-\lambda)(-\lambda-1)\cdots(-\lambda-r+1)\alpha^{-\lambda-r}
$$
which implies that 
$$
c_{\alpha}=\frac{r!a_r}{(-\lambda)(-\lambda-1)\cdots(-\lambda-r+1)\alpha^{-\lambda-r}}
$$
$$
\begin{array}{ll}
g(\alpha) &= c_{\alpha}\alpha^{-\lambda}
=\frac{r!a_r\alpha^{-\lambda}}{(-\lambda)(-\lambda-1)\cdots(-\lambda-r+1)\alpha^{-\lambda-r}}\\
&= \frac{r!a_r\alpha^{r}}{(-\lambda)(-\lambda-1)\cdots(-\lambda-r+1)}
\end{array}
$$
Since this holds for infinitely many $\alpha$, and $g(t)$ is a polynomial, 
$$
g(t) = \frac{r!a_rt^r}{-\lambda(-\lambda-1)\cdots(-\lambda-r+1)}.
$$
Thus 
$$
\sum_{i\le r} a_{r-i,i}t^i =  \frac{r!a_rt^r}{-\lambda(-\lambda-1)\cdots(-\lambda-r+1)}
$$
so that  $a_{r-i,i}=0$ if $i<r$. Thus
$L_p=a_{0r}y^r$, a contradiction to the assumption that $L$ depends on two variables.

Thus the only curves in $\overline{S_{r+1}(X_1)}\cap\pi^{-1}(p)$ which contain a 1 point
are on the strict transforms of $y-\alpha x=0$ for a finite number of nonzero $\alpha$.
These lines contain the  3 point of $X$ which has permissible parameters $(x_1,y_1,z_1)$
defined by $x=x_1z_1, y=y_1z_1$, $z=z_1$.

Since $L_p=L(x,y)$, there is at most one 3 point $q$ in $\pi^{-1}(p)$ with $\nu(q)=r$. Thus
there are no 2 curves in $\overline{S_{r+1}(X_1)}\cap\pi^{-1}(p)$ by Lemma \ref{Lemma8}.

Now suppose that $p$ is a 3 point and $L_p=L(x)$.  Suppose that  $q\in\pi^{-1}(p)$ is a 1 point and $\nu(q)=r+1$.
$\hat{\cal O}_{X_1,q}$ has regular parameters $x_1$, $y_1$, $z_1$ such that 
$$
\begin{array}{ll}
x&=x_1\\
y&=x_1(y_1+\alpha)\\
z&=x_1(z_1+\beta)
\end{array}
$$
where $\alpha,\beta\ne 0$.
$$
L(x) = \overline ax^r
$$
With the notation of (\ref{eq22}), we have
$$
(y_1+\alpha)^{\lambda_1}(z_1+\beta)^{\lambda_2}\overline a \equiv c_{\alpha,\beta}
\text{ mod }(y_1,z_1)^{r+1}
$$
for some $c_{\alpha,\beta}\in k$,
which implies $\lambda_1=\lambda_2=0$.

From equation (\ref{eq22}) we have
$$
\begin{array}{ll}
e&= \frac{b(d+e+f+r)}{a+b+c}\\
f&=\frac{c(d+e+f+r)}{a+b+c}
\end{array}
$$
where
$u=x^ay^bz^c$
and $x^dy^ez^fL_p=\overline ax^{d+r}y^ez^f$.
Thus $ec-fb=0$,  $ae-b(d+r)=0$ and $af-c(d+r)=0$.
It follows that $F_p$ is not normalized, a contradiction. 

The fact that there is at most one curve $C$ in $\pi^{-1}(p)\cap \overline S_{r+1}(X_1)$,
which  is the 2 curve which is the intersection of the strict transform of $x=0$ with $\pi^{-1}(p)$,
follows from Lemma \ref{Lemma8}, since at the 3 point $q$ with permissible parameters
$(x_1,y_1,z_1)$ defined by $x=x_1,y=x_1y_1$, $z=x_1z_1$, $\nu(q)=0$.

Suppose that $p$ is a 2 point. By Theorem \ref{Theorem9}, there are no 2 curves in
$\overline S_{r+1}(X_1)\cap\pi^{-1}(p)$. Suppose that $\overline S_{r+1}(X_1)\cap\pi^{-1}(p)$
contains a 1 point. Then $\tau(p)=0$ by Theorem \ref{Theorem13}. The leading form of
$F_p$ has an expression
$$
L_p=\sum_{i+j=r}a_{ij}x^iy^j.
$$
After possibly interchanging $x$ and $y$, we may assume that $L\ne a_{0r}y^r$.

Suppose that there exist infinitely many distinct values of $\alpha\in k$ such that
there exists a 1 point $q\in \overline S_{r+1}(X_1)\cap\pi^{-1}(p)$ with regular 
parameters $(x_1,y_1,z_1)$ in $\hat{\cal O}_{X_1,q}$ defined by 
$$
x=x_1,
y=x_1(y_1+\alpha),
z=x_1(z_1+\beta)
$$
for some $\alpha,\beta\in k$ with $\alpha\ne 0$, such that $\nu(q)=r+1$.

With the notation of (\ref{eq19}) of Theorem \ref{Theorem9}, there exist $c_{\alpha}\in k$
such that
$$
\sum_{i+j=r}a_{ij}(y_1+\alpha)^j=c_{\alpha}(y_1+\alpha)^{-\lambda}\text{ mod }y_1^{r+1}
$$
Set $g(t)=\sum_{i+j=r}a_{ij}t^j$.
$g(\alpha)=c_{\alpha}\alpha^{-\lambda}$.
$$
r!a_{0r}=\frac{d^rg}{dt^r}(\alpha)
=c_{\alpha}(-\lambda)(-\lambda-1)\cdots(-\lambda-r+1)\alpha^{-\lambda-r}
$$
implies
$$
g(\alpha)=\frac{r!a_{0r}\alpha^r}
{(-\lambda)(-\lambda-1)\cdots(-\lambda-r+1)}
$$
for infinitely many $\alpha$, so that
$$
L_p=\frac{r!a_{0r}}{(-\lambda)(-\lambda-1)\cdots(-\lambda-r+1)}y^r
$$
a contradiction. Thus 1 curves in $\pi^{-1}(p)\cap\overline S_{r+1}(X_1)$ must be
the intersection of the strict transform of $y-\alpha x=0$ and $\pi^{-1}(p)$
for a finite number of  $0\ne\alpha\in k$. These lines intersect in the 3 point of $\pi^{-1}(p)$.
 
\end{pf}

\begin{Lemma}\label{Lemma54}
Suppose that  $r\ge 2$ and $p\in X$ is such that
\begin{enumerate}
\item $\nu(p)\le r$ if $p$ is a 1 point or a 2 point.
\item If $p$ is a 2 point and $\nu(p)= r$, then $\tau(p)>0$.
\item $\nu(p)\le r-1$ if $p$ is a 3 point
\end{enumerate}
and $\pi:Y\rightarrow X$ is the blowup of a point $p\in X$. Then 
$C$ is a line for every curve $C$ in $\overline S_r(Y)\cap
\pi^{-1}(p)$ containing a 1 point.
 Thus $C$ intersects a 2 curve in at most one point, and this intersection must be transversal.

If $p$ is a 1 or 2 point with $\nu(p)=r$ then there is at most one curve $C$ in
$\overline S_r(Y)\cap \pi^{-1}(p)$ containing a 1 point. 
\end{Lemma}

\begin{pf} Suppose that $p$ is a 1 point. Suppose that $q\in \pi^{-1}(p)$ is a
1 point with $\nu(q)=r$. After a permissible change of parameters at $p$,
we have permissible parameters
$x_1,y_1,z_1$   at $q$ defined by 
$$
x=x_1,y=x_1y_1, z=x_1z_1.
$$
Write
$$
F_p=\sum_{i+j+k\ge r}a_{ijk}x^iy^jz^k.
$$
$F_p$ has leading form
$$
L_p=\sum_{i+j+k=r}a_{ijk}x^iy^jk^k.
$$
Thus $L_p=L(y,z)$ depends only on $y$ and $z$.

Suppose that $q'\in\pi^{-1}(p)$ is another 1 point with $\nu(q')=r$, with permissible
parameters $(x_1,y_1,z_1)$ defined by
$$
x=x_1, y=x_1(y_1+\alpha), z=x_1(z_1+\beta)
$$
for some $\alpha,\beta\in k$.
Then there exists a form $L_1$ such that
$$
L_p=L_1(y-\alpha x,z- \beta x)+cx^r
$$
for some $c\in k$. There exist $\alpha_i,\beta_i,\gamma_i,\delta_i\in k$ such that
$$
L_p(y,z)=\prod_{i=1}^r(\alpha_iy-\beta_iz)
$$
$$
L_1(y,z)=\prod_{i=1}^r(\gamma_iy-\delta_iz)
$$
We can also assume that $\alpha_i\beta_j-\alpha_j\beta_i=0$ implies
$\alpha_i=\alpha_j$ and $\beta_i=\beta_j$.
$$
\prod_{i=1}^r(\alpha_iy-\beta_iz)=\prod_{i=1}^r(\gamma_i(y-\alpha x)
-\delta_i(z-\beta x))+cx^r.
$$
Set $x=0$ to get that, after reindexing the $(\gamma_i,\delta_i)$, there exist $0\ne \epsilon_i\in k$ such that 
$$
(\alpha_i,\beta_i)=\epsilon_i(\gamma_i,\delta_i)
$$
for all $i$, and $\prod_{i=1}^r\epsilon_i=1$. Thus
 $$
\prod_{i=1}^r(\alpha_iy-\beta_iz)=\prod_{i=1}^r(\alpha_i(y-\alpha x)
-\beta_i(z-\beta x))+cx^r.
$$
First suppose that there exists $(\alpha_i,\beta_i)$, $(\alpha_j,\beta_j)$ such that
$\alpha_i\beta_j-\alpha_j\beta_i\ne 0$.
Suppose that $\alpha_i\ne 0$. There exist $t<r$ distinct values of
$(\alpha_k, \beta_k)$ such that $\alpha_i\beta_k-\alpha_j\beta_k=0$. Set 
$y=\frac{\beta_iz}{\alpha_i}$ to get
$$
0=(\beta_i\beta-\alpha_i\alpha)^tx^t\prod_{j\mid \alpha_i\beta_j-\alpha_j\beta_i\ne 0}\left(\left(\frac{\alpha_j}{\alpha_i}\beta_i-\beta_j\right)
z+\left(\beta\beta_j-\alpha\alpha_j\right)x\right)+cx^r
$$
We conclude that $\beta\beta_i-\alpha \alpha_i=0$. If $\beta_i\ne 0$, we can set
$z=\frac{\alpha_i}{\beta_i}y$ to again conclude that
$\beta\beta_i-\alpha\alpha_i=0$.
Thus $\beta\beta_i-\alpha\alpha_i=0$ for all $i$, and the 1 points $q'\in\pi^{-1}(p)\cap \overline S_r(X)$
must thus lie on the lines $\gamma_i$ which are the intersection of the strict transform of
$\beta_iz-\alpha_iy=0$ and $\pi^{-1}(p)$.

Thus $q'$ must be in the intersection $\cap \gamma_i\subset\pi^{-1}(p)\cong\bold P^2$,
and there is at most one point $q\in\pi^{-1}(p)$ such that $\nu(q)=r$.

Now suppose that $\alpha_i\beta_j-\alpha_j\beta_i=0$ for all $i,j$. Then after a
permissible change of parameters at $p$, we have $L_p=z^r$.

$L_p=(z-\beta x)^r+cx^r$ and $r\ge 2$ implies $\beta=c=0$, so $q'$ is on the line
$\gamma\subset \pi^{-1}(q)\subset\bold P^2$ which is the intersection of the strict
transform of $z=0$ and $\pi^{-1}(q)$.

Suppose that $p$ is a 2 point such that $\nu(p)=r$ and $\tau(p)>0$. Write
$$
F_p=\sum_{i+j+k\ge r}a_{ijk}x^iy^jz^k.
$$
Suppose there exists a 1 point $q\in \pi^{-1}(p)$ such that $\nu(q)=r$.
After a permissible change of parameters at $p$,
 $q$ has permissible parameters $(\overline x_1,y_1,z_1)$
at $q$ such that, with the notation of (\ref{eq19}) of Theorem \ref{Theorem9}, 
$$
x=x_1, y=x_1(y_1+\alpha), z=x_1z_1
$$
$$
x_1=\overline x_1(y_1+\alpha)^{-\frac{b}{a+b}}
$$
$$
F_q=\sum_{i+j+k=r}a_{ijk}(y_1+\alpha)^{j+\lambda}z_1^k
-\sum_{i+j=r}a_{ij0}\alpha^{j+\lambda}+\overline x_1\Omega.
$$
Let
$$
L_p=\sum_{i+j+k=r}a_{ijk}x^iy^jz^k
$$
be the leading form of $F_p$.
$$
F_q=\sum_{k>0}(\sum_{i+j=r-k}a_{ijk}(y_1+\alpha)^j)(y_1+\alpha)^{\lambda}z_1^k
+(\sum_{i+j=r}a_{ij0}(y_1+\alpha)^{j+\lambda}-\sum_{i+j=r}a_{ij0}\alpha^{j+\lambda})
+\overline x_1\Omega.
$$
$\nu(q)=r$ implies, for fixed $k>0$,
$$
\sum_{i+j=r-k}a_{ijk}x^iy^j=c_k(y-\alpha x)^{r-k}
$$
for some $c_k\in k$, thus
$$
L_p=\sum_{k>0}c_k(y-\alpha x)^{r-k}z^k+G(x,y).
$$
$\tau(p)>0$ implies some $c_k\ne 0$.

Suppose that there exists another 1 point $q'\in\pi^{-1}(p)$ with $\nu(q')=r$.
$\hat{\cal O}_{Y,q'}$ has regular
parameters $(x_1,y_1,z_1)$ such that 
$$
x=x_1, y=x_1(y_1+\overline \alpha), z=x_1(z_1+\overline \beta)
$$
with $\overline\alpha\ne 0$. Then 
$$
L_p=\sum_{k>0} \overline c_k(y-\overline\alpha x)^{r-k}(z-\overline \beta x)^k+\overline G(x,y).
$$
Thus 
\begin{equation}\label{eq513}
\sum_{k>0} c_k(y-\alpha x)^{r-k}z^k=\sum_{k>0}\overline c_k(y-\overline \alpha x)^{r-k}(z-\overline \beta x)^k
+H(x,y).
\end{equation}
Set $x=0$ in (\ref{eq513}) to get $c_k=\overline c_k$ for all $k$. Let 
$$
k_0=\text{max}\{k\mid c_k\ne 0\}=\tau(p).
$$
By assumption, $k_0>0$. 
$$
c_{k_0}(y-\alpha x)^{r-k_0}z^{k_0}=c_{k_0}(y-\overline\alpha x)^{r-k_0}z^{k_0},
$$
and if $k_0>1$,
$$
c_{k_0-1}(y-\alpha x)^{r-k_0+1}z^{k_0-1}
=c_{k_0}(y-\overline \alpha x)^{r-k_0}(-\overline\beta k_0 x)z^{k_0-1}+
c_{k_0-1}(y-\overline\alpha x)^{r-k_0+1}z^{k_0-1}.
$$

If $k_0<r$, then $\alpha=\overline \alpha$ implies
 all 1 points in $\overline S_r(X_1)\cap \pi^{-1}(p)$
are contained in the line which is the intersection of the strict transform of
 $y-\alpha x=0$ and $\pi^{-1}(p)$. This line contains the 3 point of $\pi^{-1}(p)$.

If $k_0=r$ ($\ge 2$),
$$
c_{r-1}(y-\alpha x)=-c_r\overline\beta rx+c_{r-1}(y-\overline \alpha x)
$$
so that
$$
-\alpha c_{r-1}=-c_r\overline \beta r-\overline \alpha c_{r-1}
$$
which implies that all 1 points in $\overline S_r(X_1)\cap \pi^{-1}(p)$ are contained in the
line which is the intersection of the strict transform of
$$
c_rrz+c_{r-1}y-\alpha c_{r-1}x=0
$$
and $\pi^{-1}(p)$. 

Suppose that $p$ is a 2 point or a 3 point,with $\nu(p)=r-1$. Then by Theorem \ref{Theorem12},
the conclusions of the Theorem hold.
\end{pf} 

\section{Permissible Monoidal Transforms Centered at Curves}

Throughout this section we will assume that $\Phi_X:X\rightarrow S$ is weakly prepared.

\begin{Lemma}\label{Lemma655}
Suppose that $C\subset X$ is a 2 curve. Then either 
$F_p\in\hat{\cal I}_{C,p}$ for all $p\in C$ or $F_p\not\in
\hat{\cal I}_{C,p}$ for all $p\in C$.

 Suppose that $r\ge 2$ , $C\subset \overline S_r(X)$
is a 2 curve. Then either $F_p\in\hat{\cal I}_{C,p}^r$ for all $p\in C$ or
$F_p\in \hat{\cal I}_{C,p}^{r-1}$, $F_p\not\in\hat{\cal I}_{C,p}^r$ for all $p\in C$
\end{Lemma}

\begin{pf} This follows from Lemmas \ref{Lemma2}, \ref{Lemma6}, \ref{Lemma8}.
\end{pf}

\begin{Lemma}\label{Lemma653} Suppose that $r\ge 2$ and $C\subset \overline S_r(X)$
is a nonsingular curve such that $C$ contains a 1 point and $C$ makes SNCs with $\overline B_2(X)$.
Then either $F_p\in\hat{\cal I}_{C,p}^r$ with respect to permissible
parameters for $C$ at $p$ for all $p\in C$, or 
$F_p\in\hat{\cal I}_{C,p}^{r-1}$, $F_p\not\in \hat{\cal I}_{C,p}^r$ 
with respect to permissible parameters for $C$ at $p$ for all $p\in C$.
\end{Lemma}
\begin{pf}
This follows from Lemmas \ref{Lemma2}, \ref{Lemma5}, \ref{Lemma7}.
\end{pf}

\begin{Definition}\label{Def1091} Suppose that $r\ge 2$, $p\in X$, $C\subset \overline S_r(X)$ is a  curve which contains $p$ and makes SNCs with $\overline B_2(X)$ at $p$ and $C\not\subset\overline S_{r+1}(X)$.
$C$ is r big at $p$ if $F_p\in\hat{\cal I}_{C,p}^r$ with respect to permissible parameters
for $C$ at $p$. $C$ is r small at $p$ if $C$ is not r big at $p$.

Suppose that $C$ is a 2 curve, $\nu(q)\ge 1$ if $q\in C$ is a 2 point, $C\not\subset\overline S_2(X)$ and $p\in C$. Then $C$ is 1 big at $p$ if $F_p\in\hat{\cal I}_{C,p}$.
$F_p$ is 1 small at $p$ if $C$ is not 1 big at $p$.

Suppose that $r\ge 2$, $C\subset \overline S_r(X)$ is a curve which makes SNCs with $\overline B_2(X)$. We will say that $C$ is r big if $C$ is r big at $p$ for all $p\in C$.
We will say that $C$ is r small if $C$ is r small at $p$ for all $p\in C$.

Suppose that $C$ is a 2 curve, $\nu(q)\ge 1$ if $q\in C$ is a 2 point, $C\not \subset \overline S_2(X)$. We will say that $C$ is 1 big if $C$ is 1 big for all $p\in C$. We will say that $C$ is
1 small if $C$ is 1 small at $p$ for all $p\in C$.
\end{Definition}

\begin{Lemma}\label{Lemma960}
Suppose that $C$ is a 2 curve on $X$, $p\in C$ is a 2 point, $D_1$  and $D_2$ are curves in $E_X$ 
containing $p$ such that $D_1\cup D_2$ makes SNCs with $C$ at $p$. Then there are
regular parameters $(x,y,z)$ in ${\cal O}_{X,p}$ such that 
$$
{\cal I}_{C,p}=(x,y), {\cal I}_{D_1,p}=(x,z), {\cal I}_{D_2,p}=(y,z)
$$
\end{Lemma}

\begin{pf}
There exist regular parameters $(\tilde x,\tilde y,\tilde z)$ in ${\cal O}_{X,p}$, and
$\phi\in{\cal O}_{X,p}$ such that  
$$
{\cal I}_{C,p}=(\tilde x,\tilde y), {\cal I}_{D_1,p}=(\tilde x,\tilde z), 
{\cal I}_{D_2,p}=(\tilde y,\phi)
$$
and $\phi\equiv a\tilde x+c\tilde z\text{ mod }m_p^2$, with
$a,c\in k$, $c\ne 0$. In $\hat{\cal O}_{X,p}$, there exist series $h,g$ such that
$$
\phi=h(\tilde x,\tilde y,\tilde z)\tilde y+g(\tilde x,\tilde z)
$$
$$
g=u(\tilde z-\psi(\tilde x))
$$
where $u$ is a unit, $\psi$ is a series.
$$
\tilde z-\psi(\tilde x)\in\hat{\cal I}_{D_1,p}\cap \hat{\cal I}_{D_2,p}
=({\cal I}_{D_1,p}\cap {\cal I}_{D_2,p})\hat{\cal O}_{X,p},
$$
where the last equality is by Corollary 2 to Theorem 11 of Chapter VIII, section 4 \cite{ZS}).
Suppose that 
$$
{\cal I}_{D_1,p}\cap {\cal I}_{D_2,p}=(f_1,\ldots,f_n).
$$
$$
\tilde z-\psi(\tilde x)=\sum \lambda_i f_i
$$
implies there exists $f\in {\cal I}_{D_1,p}\cap {\cal I}_{D_2,p}$ such that
$$
f\equiv \overline a\tilde x+\overline c \tilde z\text{ mod }m_p^2
$$
where $\overline a,\overline c\in k$, $\overline c\ne 0$. Since $f\in(\tilde x,\tilde z)$, we have
$f=\lambda\tilde x+\tau\tilde z$ where $\tau$ is a unit. Thus
$(\tilde x,\tilde z)=(\tilde x,f)$.
 Since $f\in(\tilde y,\phi)$, we have
$f=\alpha\tilde y+\beta\phi$ where $\beta$ is a unit. Thus
$(\tilde y,\phi)=(\tilde y,f)$. $(\tilde x,\tilde y, f)$ are the desired
regular parameters.
\end{pf}

\begin{Lemma}\label{Lemma17}
Suppose that $p\in X$ is a 1 point or a 2 point with $\gamma(p)=r\ge 2$,
 and  $(u,v)$ are permissible parameters at $\Phi_X(p)$, such that $u=0$
is a local equation of $E_X$ at $p$. Then there exist regular parameters
$(\tilde x, y, \tilde z)$ in $R ={\cal O}_{X,p}$ and permissible parameters $(x,y,z)$
at $p$  with $x=\gamma\tilde x$, $z=\sigma\tilde z$ for some series $\gamma, \sigma\in
\hat{\cal O}_{X,p}$
such that if $p$ is a 1 point, 
\begin{equation}\label{eq37}
\begin{array}{ll}
u&=x^a\\
v&=P(x)+x^cF
\end{array}
\end{equation}
with $F = \tau z^r+\sum_{i=2}^ra_i(x,y)z^{r-i}$, $\tau$ a unit and some $a_i\ne 0$.

Further suppose that $\overline S_r(X)\cup\overline B_2(X)$ makes SNCs at $p$.
Then there is at most one curve $D$ in $\overline S_r(X)$ through $p$. If $D$ exists, 
we can choose $(x,y,z)$ so that 
$x=0, z=0$ are local equations of $D$ at $p$.

If $p$ is a 2 point, 
\begin{equation}\label{eq38}
\begin{array}{ll}
u&=(x^ay^b)^m\\
v&=P(x^ay^b)+x^cy^dF
\end{array}
\end{equation}
with $F = \tau z^r+\sum_{i=2}^ra_i(x,y)z^{r-i}$, $\tau$ a unit, and some $a_i\ne 0$.

Further suppose that $\overline S_r(X)\cup\overline B_2(X)$ makes SNCs at $p$.
Then there are at most 2 curves $D_1$ and $D_2$ in $\overline S_r(X)$ through $p$.
If $D_1$ exists (or if $D_1$ and $D_2$ exist) then we can choose $(x,y,z)$ so that
$x=0, z=0$ are local equations of
$D_1$ at $p$ ($ x=0,  z=0$ are local equations of $D_1$ at $p$ and
$y=0, z=0$ are local equations of $D_2$ at $p$).
\end{Lemma}

\begin{pf}
There exist regular parameters $\tilde x,y,\tilde z$ in $R$, and permissible parameters
$(x=\gamma\tilde x,y,\tilde z)$ at $p$  such that
$u=x^a$ or $u=(x^ay^b)^m$ in $\hat R$, and $\nu(F(0,0,\tilde z))=r$.

If $p$ is a 1 point, then there exists at most one curve $D$ in $\overline S_r(X)$
containing $p$ by Lemma \ref{Lemma11}. If $D$ exists, we may assume that
$ x=0, \tilde z=0$ are local equations of $D$ at $p$. If $p$ is a 2 point,
then there exist at most 2 curves $D_1$ and $D_2$ in $\overline S_r(X)$.
If $D_1$ (or $D_1$ and $D_2$ exist) we may assume that $x=0, \tilde z=0$ are local
equations of $D_1$ at $p$ (or $x=0, \tilde z=0$ are local equations of $D_1$ at $p$
and $y=0, \tilde z=0$ are local equations of $D_2$ at $p$ by Lemma \ref{Lemma960}).

 Set
$$
\overline z = \frac{\partial^{r-1} F}{\partial \tilde z^{r-1}}= \omega(\tilde z-\phi(x,y))
$$
where $\omega$ is a unit by the formal implicit function theorem. Set 
$z_1=\tilde z-\phi(x,y)$,
$G(x,y,z_1) = F(x,y,\tilde z)$.

Suppose that $p$ is a 1 point and there exists a curve $D\subset \overline S_r(X)$
containing $p$, so that $D$ has local equations $x=0, \tilde z=0$. Then
$F_p\in\hat{\cal I}_{D,p}^r+(x^{r-1})$ by Lemma \ref{Lemma5}, so that
$$
\frac{\partial^{r-1}F}{\partial \tilde z^{r-1}}\in\hat{\cal I}_{D,p}
$$
and $x\mid\phi(x,y)$.
Thus $x=0, z_1=0$ are local equations of $D$ at $p$.

Suppose that $p$ is a 2 point and there exist curves $D_1,D_2\subset\overline S_r(X)$
containing $p$, so that $D_1$ has local equations $x=0,\tilde z=0$ and $D_2$
has local equations $y=0, \tilde z=0$. $F_p\in\hat{\cal I}_{D_1,p}^r+(x^{r-1})$
and $F_p\in\hat{\cal I}_{D_2,p}^r+(y^{r-1})$ by Lemma \ref{Lemma7}. Thus
$$
\frac{\partial^{r-1} F}{\partial\tilde z^{r-1}}\in\hat{\cal I}_{D_1,p}
$$
and $\frac{\partial^{r-1}F}{\partial\tilde z^{r-1}}\in\hat{\cal I}_{D_2,p}$,
so that $xy\mid\phi(x,y)$, and $x=0, z_1=0$ are local equations of $D_1$ at $p$
and $y=0,z_1=0$ are local equations of $D_2$ at $p$.

$$
G=G(x,y,0)+\frac{\partial G}{\partial z_1}(x,y,0)z_1+\cdots + \frac{1}{(r-1)!}\frac{\partial^{r-1}G}{\partial z_1^{r-1}}(x,y,0)z_1^{r-1}+\frac{1}{r!}
\frac{\partial^r G}{\partial z_1^r}(x,y,0)z_1^r+\cdots 
$$
$$
\frac{\partial^{r-1}G}{\partial z_1^{r-1}}(x,y,0)=
\frac{\partial^{r-1}F}{\partial \tilde z^{r-1}}(x,y,\phi(x,y))=0
$$
$$
\frac{\partial^r G}{\partial z_1^r}(x,y,0)=\frac{\partial^r F}{\partial \tilde z^r}(x,y,\phi(x,y))
$$
is a unit.
Thus with the regular parameters $(\tilde x,y,\overline z)$ in $R$ and permissible parameters $(x,y,z_1)$ at $p$, $F$ has the desired form.

We cannot have $a_i=0$ for all $i$, since $r\ge 2$ and  $x$ or 
$xy\in\sqrt{\hat{\cal I}_{\text{sing}(\Phi_X),p}}$.
\end{pf}

\begin{Lemma}\label{Lemma500}
Suppose that $r\ge 2$, $C\subset X$ is a 2 curve such that
$C$ is $r-1$ big or $r$ small, 
 $\pi:X_1\rightarrow X$
is the blowup of $C$.
\begin{enumerate}
\item 
\begin{enumerate}
\item If $q\in C$ is a 2 point with $\nu(q)=r-1$ and $q_1\in\pi^{-1}(q)$, then
\begin{enumerate}
\item If $q_1$ is a 1 point then $\nu(q_1)\le r$ and  $\gamma(q_1)\le r$.
\item If $q_1$ is a 2 point then $\nu(q_1)\le r-1$.
\end{enumerate}
\item If $q\in C$ is a 2 point with $\nu(q)=r$, $\tau(q)>0$ and $q_1\in\pi^{-1}(q)$, then
\begin{enumerate}
\item If $q_1$ is a 1 point then $\nu(q_1)\le r$. $\nu(q_1)=r$ implies $\gamma(q_1)=r$.
\item If $q_1$ is a 2 point then $\nu(q_1)\le r$. $\nu(q_1)=r$ implies $\tau(q_1)>0$.
\end{enumerate}
\item If $q\in C$ is a 3 point with $\nu(q)=r-1$ and $q_1\in\pi^{-1}(q)$ then
\begin{enumerate}
\item $q_1$ a 2 point implies $\nu(q_1)\le r$ and  $\gamma(q_1)\le r$.
\item $q_1$ a 3 point implies $\nu(q_1)\le r-1$.
\end{enumerate}
\end{enumerate}
\item Suppose that $C\subset \overline S_r(X)$ ( so that $C$ is r small). 
If $q\in C$ is a 2 point with $\nu(q)=r$, $\tau(q)>0$ and $q_1\in \pi^{-1}(q)$, then
\begin{enumerate}
\item If $q_1$ is a 1 point then $q_1$ is resolved.
\item If $q_1$ is a 2 point then $\nu(q_1)\le r$. $\nu(q_1)=r$ implies $\tau(q)>0$.
\end{enumerate}
\end{enumerate}
\end{Lemma}

\begin{pf}
Suppose that $q\in C$ is a 2 point with $\nu(q)=r-1$, and $q$ has permissible parameters
$(x,y,z)$ with
$$
\begin{array}{ll}
u&=(x^ay^b)^m\\
v&=P(x^ay^b)+x^cy^dF_q\\
L_q&=\sum_{i+j=r-1}a_{ij}x^iy^j
\end{array}
$$

Suppose that $q_1\in\pi^{-1}(q)$ and $\hat{\cal O}_{Y_1,q_1}$ has regular
parameters $(x_1,y_1,z)$ such that
$$
x=x_1, y=x_1(y_1+\alpha)
$$
with $\alpha\ne 0$. Set
$$
x_1=\overline x_1(y_1+\alpha)^{-\frac{b}{a+b}}
$$
\begin{equation}\label{eq623}
\begin{array}{ll}
u&=\overline x_1^{(a+b)m}\\
v&=P_{q_1}(\overline x_1)+\overline x_1^{c+d+r-1}F_{q_1}\\
F_{q_1}&=\sum_{i+j=r-1}a_{ij}(y_1+\alpha)^{\lambda+j}-\sum_{i+j=r-1}a_{ij}\alpha^{\lambda+j}
+\overline x_1\Omega
+zG,
\end{array}
\end{equation}
$$
\lambda=d-\frac{b(c+d+r-1)}{a+b}.
$$
Thus $\nu(q_1)\le r$ and $\gamma(q_1)\le r$.

Suppose that $q_1\in \pi^{-1}(q)$ and $q_1$ has permissible 
parameters $(x_1,y_1,z)$ such that 
$$
x=x_1, y=x_1y_1.
$$
Then
$$
\begin{array}{ll}
u&=(x_1^{a+b}y_1^b)^m\\
v&=P(x_1^{a+b}y_1^b)+x_1^{c+d+r-1}F_{q_1}\\
F_{q_1}&=\sum_{i+j=r-1}a_{ij}y_1^j+x_1\Omega+zG
\end{array}
$$

implies that $\nu(q_1)\le r-1$.

A similar argument holds at the point $q_1\in\pi^{-1}(q)$ with permissible parameters
$(x_1,y_1,z)$ such that $x=x_1y_1, y=y_1$.

Suppose that $q\in C$ is a 2 point with $\nu(q)=r$ and $\tau(q)>0$. Then $q$ has permissible
parameters $(x,y,z)$ with
$$
\begin{array}{ll}
u&=(x^ay^b)^m\\
v&=P(x^ay^b)+x^cy^dF_q\\
L_q&=z(\sum_{i+j=r-1}a_{ij1}x^iy^j)+\sum_{i+j=r}a_{ij0}x^iy^j
\end{array}
$$
with some $a_{ij1}\ne 0$.

Suppose that $q_1\in\pi^{-1}(q)$ and $\hat{\cal O}_{Y_1,q_1}$ has regular parameters
$(x_1,y_1,z)$ such that $x=x_1,y=x_1(y_1+\alpha)$ with $\alpha\ne 0$.
Set 
$$
x_1=\overline x_1(y_1+\alpha)^{-\frac{b}{a+b}}.
$$
$$
\begin{array}{ll}
u&=\overline x_1^{(a+b)m}\\
v&=P_{q_1}(\overline x_1)+\overline x_1^{c+d+r-1}F_{q_1}
\end{array}
$$
with 
\begin{equation}\label{eq624}
\begin{array}{ll}
F_{q_1}&=z(\sum_{i+j=r-1}a_{ij1}(y_1+\alpha)^j)(y_1+\alpha)^{\lambda}
+\overline x_1\Omega+z^2G,
\end{array}
\end{equation}
$$
\lambda=d-\frac{b(c+d+r-1)}{a+b}.
$$
Thus $\nu(q_1)\le r$ and $\nu(q_1)=r$ implies $\gamma(q_1)=r$.

Suppose that $q_1\in\pi^{-1}(q)$ and $q_1$ has permissible
parameters $(x_1,y_1,z)$ such that
$$
x=x_1,
y=x_1y_1.
$$
Then
$$
\begin{array}{ll}
u&=(x_1^{a+b}y_1^b)^m\\
v&=P(x_1^{a+b}y_1^b)+x_1^{c+d+r-1}y_1^dF_{q_1}\\
F_{q_1}&=\sum_{i+j=r-1}za_{ij1}y_1^j+x_1\Omega+z^2G
\end{array}
$$
implies $\nu(q_1)\le r$ and $\nu(q_1)=r$ implies $\tau(q_1)>0$.

A similar analysis holds at the point $q_1\in \pi^{-1}(q)$ with permissible
parameters $(x_1,y_1,z)$ such that $x=x_1y_1,y=y_1$.

Suppose that $q\in C$ is a 3 point with $\nu(q)=r-1$,
$$
\begin{array}{ll}
u&=(x^ay^bz^c)^m\\
v&=P(x^ay^bz^c)+x^dy^ez^fF_q\\
F_q&=\sum_{i+j\ge r-1, k\ge 0}a_{ijk}x^iy^jz^k
\end{array}
$$
some $a_{ij0}\ne 0$ with $i+j=r-1$.

Suppose that $q_1\in\pi^{-1}(q)$ is a 2 point.
$$
x=x_1,
y=x_1(y_1+\alpha)
$$
with $\alpha\ne 0$.
$$
x_1=\overline x_1(y_1+\alpha)^{-\frac{b}{a+b}}
$$
$$
\begin{array}{ll}
u&=(\overline x_1^{a+b}z^c)^m=(\overline x_1^{\overline a}z^{\overline c})^{\overline m}\\
v&=P_{q_1}(\overline x_1^{\overline a}z^{\overline c})+\overline x_1^{d+r-1+e}z^fF_{q_1}
\end{array}
$$
with
$$
\lambda=e-\frac{b(d+r-1+e)}{a+b}, (\overline a,\overline c)=1
$$
\begin{equation}\label{eq625}
F_{q_1}=(y_1+\alpha)^{\lambda}\frac{F_q}{x_1^{r-1}}
-\frac{g(\overline x_1^{\overline a}z^{\overline c})}{\overline x_1^{d+r-1+e}z^f}
\end{equation}
Thus
$$
F_{q_1}=\sum_{i+j=r-1}a_{ij0}(y_1+\alpha)^{j+\lambda}+zG+\overline x_1\Omega,
$$
or
$$
F_{q_1}=\sum_{i+j=r-1}a_{ij0}(y_1+\alpha)^{j+\lambda}
-\sum a_{ij0}\alpha^{j+\lambda}+zG+\overline x_1\Omega,
$$
implies $\nu(q_1)\le r$, and $\gamma(q_1)\le r$.

Suppose that $q\in C$ is a 2 point with $\nu(q)=r$ and $\tau(q)>0$
and $C\subset \overline S_r(X)$.
 By Lemma \ref{Lemma6}, $q$ has  permissible 
parameters $(x,y,z)$ with 
\begin{equation}\label{eq957}
\begin{array}{ll}
u&= (x^ay^b)^m\\
v&= P(x^ay^b)+x^cy^dF_q\\
F_q &= \sum_{i+j\ge r, k\ge 0}c_{ijk}x^iy^jz^k + \overline c z x^{i_0}y^{j_0}
\end{array}
\end{equation}
where  $i_0+j_0=r-1$, $(c+i_0)b-a(d+j_0)=0$, $\overline c\ne0$. 

Suppose that $q_1\in \pi^{-1}(q)$, and $\hat{\cal O}_{Y_1,q_1}$ has regular parameters $(x_1,y_1,z)$
such that 
$$
x=x_1, y=x_1(y_1+\alpha)
$$
with $\alpha\ne 0$. Set $ x_1 = \overline x_1(y_1+\alpha)^{-\frac{b}{a+b}}$. Then
$$
\begin{array}{ll}
u&=\overline x_1^{(a+b)m}\\
v&= P_{q_1}(\overline x_1)+\overline x_1^{c+d+r-1}F_{q_1}\\
F_{q_1} &= \overline c z (y_1+\alpha)^{\lambda+j_0} 
+\overline x_1\Omega\
\end{array}
$$
where $\lambda = d-\frac{b(c+d+r-1)}{a+b}$. Thus $q_1$ is resolved.

Suppose that $q_1\in \pi^{-1}(q)$, and $\hat{\cal O}_{Y_1,q_1}$ has regular parameters $(x_1,y_1,z)$
such that 
$$
x=x_1, y=x_1y_1.
$$
 Then
$$
\begin{array}{ll}
u&=(x_1^{a+b}y_1^b)^m\\
v&= P(x_1^{a+b}y_1^b)+x_1^{c+d+r-1}y_1^dF_{q_1}\\
F_{q_1} &=\frac{F_q}{x_1^{r-1}}= \overline c z y_1^{j_0} 
+x_1\Omega
\end{array}
$$
$q_1$ satisfies the conclusions of the Theorem since $j_0\le r-1$.

Suppose that $q_1\in \pi^{-1}(q)$, and $\hat{\cal O}_{Y_1,q_1}$ has regular parameters $(x_1,y_1,z)$
such that 
$$
x=x_1y_1, y=y_1.
$$
 Then an argument similar to the above case shows that $q_1$ satisfies the conclusions of the Theorem (since $i_0\le r-1$).

\end{pf}

\begin{Lemma}\label{Lemma501}
Suppose that $r\ge 2$, $C\subset X$ is a 2 curve such that
$C$ is r-1 big and 
\begin{enumerate}
\item $p\in C$ a 2 point implies $\nu(p)\le r$, and if $\nu(p)=r$ then $\tau(p)>0$.
\item $p\in C$ a 3 point implies $\nu(p)\le r-1$.
\end{enumerate}
Suppose that $\pi:X_1\rightarrow X$ is the blowup of $C$. Then 
$$
\pi^{-1}(C)\cap \overline{S_r(X_1)}
$$
contains at most one curve. If $D\subset \pi^{-1}(C)\cap \overline S_r(X_1)$ is
a curve, then $D$ is a section over $C$, and $D$ contains a 1 point.

Suppose that $D\subset\pi^{-1}(C)\cap \overline S_r(X_1)$ is a curve (which is necessarily
a section over $C$). Supppose that $q\in C$ is a 2 point such that $\nu(q)=r-1$. Then
$\pi^{-1}(q)\cap D$ is a 1 point.
\end{Lemma}

\begin{pf}
Suppose that $q\in C$ is a 2 point with $\nu(q)=r-1$. Suppose, with the notation of
(\ref{eq623}) of Lemma \ref{Lemma500}, that there exists $q_1\in\pi^{-1}(q)$ with
$\nu(q_1)=r$. Then there exist regular parameters $(x_1,y_1,z)$ in $\hat{\cal O}_{X_1,q_1}$ such that 
$$
x=x_1, y=x_1(y_1+\alpha)
$$
with $\alpha\ne 0$, and $\gamma_{\alpha}\in k$ such that
$$
\sum_{i+j=r-1}a_{ij}(y_1+\alpha)^j\equiv \gamma_{\alpha}(y_1+\alpha)^{-\lambda}\text{ mod }
y_1^r.
$$
$-\lambda\not\in\{0,\ldots,r-1\}$ since $F_q$ is normalized.

Set $g(t)=\sum_{i+j=r-1}a_{ij}t^j$. We have
$$
\frac{1}{i!}\frac{d^ig}{dt^i}(\alpha)=\gamma_{\alpha}\left(
\frac{-\lambda(-\lambda-1)\cdots(-\lambda-i+1)}{i!}\right)\alpha^{-\lambda-i}
$$
for $i\le r-1$. Thus
$$
a_{0,r-1}=\frac{1}{(r-1)!}\frac{d^{r-1}g}{dt^{r-1}}(\alpha)
=\gamma_{\alpha}\left(\frac{-\lambda(-\lambda-1)\cdots(-\lambda-r+2)}{(r-1)!}\right)
\alpha^{-\lambda-r+1}
$$
$$
a_{1r-2}+(r-1)a_{0r-1}\alpha=\frac{1}{(r-2)!}\frac{d^{r-2}g}{dt^{r-2}}(\alpha)
=\gamma_{\alpha}\left(\frac{-\lambda(-\lambda-1)\cdots(-\lambda-r+3)}{(r-2)!}
\right)\alpha^{-\lambda-r+2}
$$
$$
\begin{array}{ll}
a_{1r-2}&=\gamma_{\alpha}\left[
\frac{-\lambda(-\lambda-1)\cdots(-\lambda-r+3)-(-\lambda)(-\lambda-1)\cdots
(-\lambda-r+2)}{(r-2)!}\right]\alpha^{-\lambda-r+2}\\
&=\gamma_{\alpha}[\frac{\lambda(-\lambda-1)\cdots(-\lambda-r+3)(-\lambda-r+1)}
{(r-2)!}]\alpha^{-\lambda-r+2}
\end{array}
$$
$$
\frac{-\lambda(-\lambda-1)\cdots(-\lambda-r+2)}{(r-1)!}\ne 0
$$
and
$$
\frac{\lambda(-\lambda-1)\cdots(-\lambda-r+3)(-\lambda-r+1)}
{(r-2)!}\ne 0
$$
since $-\lambda\not\in\{0,\ldots, r-1\}$.

If $q_2\in\pi^{-1}(q)$ has $\nu(q_2)=r$, and $q_2\ne q_1$, then there exist
$\alpha\ne\beta\in k$ such that 
$$
a_{0,r-1}
=\gamma_{\beta}\left(\frac{-\lambda(-\lambda-1)\cdots(-\lambda-r+2)}{(r-1)!}\right)
\beta^{-\lambda-r+1}
$$
$$
a_{1r-2}=\gamma_{\beta}\left[
\frac{-\lambda(-\lambda-1)\cdots(-\lambda-r+3)-(-\lambda)(-\lambda-1)\cdots
(-\lambda-r+2)}{(r-2)!}\right]\beta^{-\lambda-r+2}
$$
which implies that 
$$
\gamma_{\beta}=\gamma_{\alpha}(\frac{\alpha}{\beta})^{-\lambda-r+1}
$$
and
$$
\gamma_{\alpha}\alpha^{-\lambda-r+2}=\gamma_{\beta}\beta^{-\lambda-r+2}
=\gamma_{\alpha}\alpha^{-\lambda-r+1}\beta.
$$
so that $\alpha=\beta$.

Thus there is at most one point $q_1\in\pi^{-1}(q)$ with $\nu(q_1)=r$.
$q_1$, if it exists, is a 1 point.

Suppose that $q\in C$ is 2 point with $\nu(q)=r$ and $\tau(q)>0$. Suppose, with the notation of
(\ref{eq624}) of Lemma \ref{Lemma500}, that there exists a 1 point
$q_1\in\pi^{-1}(q)$ with $\nu(q_1)=r$. Then there exist regular parameters
$(x_1,y_1,z)$ in $\hat{\cal O}_{X_1,q_1}$ such that $x=x_1, y=x_1(y_1+\alpha)$.

Set $g(t)=\sum_{i+j=r-1}a_{ij1}t^j$. By (\ref{eq624}),
$$
\nu(\sum_{i+j=r-1}a_{ij1}(y_1+\alpha)^j)=r-1.
$$
which implies
$g(t+\alpha)=a_{0,r-1,1}t^{r-1}$ which implies 
$g(t)=a_{0,r-1,1}(t-\alpha)^{r-1}$.

Thus there is at most one 1 point $q_1\in\pi^{-1}(q)$ with $\nu(q_1)=r$.

Suppose that $q\in C$ is a 3 point. Suppose that, with the notation of (\ref{eq625}),
of Lemma \ref{Lemma500}, that there exists a 2 point $q_1\in\pi^{-1}(q)$ with $\nu(q_1)=r$.
Then there exist regular parameters $(x_1,y_1,z)$ in $\hat{\cal O}_{X_1,q_1}$
such that 
$$
x=x_1, y=x_1(y_1+\alpha)
$$
with $\alpha\ne 0$, and $\gamma_{\alpha}\in k$ such that
$$
\sum_{i+j=r-1}a_{ij0}(y_1+\alpha)^j\equiv\gamma_{\alpha}(y_1+\alpha)^{-\lambda}
\text{ mod }y_1^r.
$$
As in the argument for the case when $q$ is a 2 point with $\nu(q)=r-1$, we can conclude
that there is at most one point $q_1\in\pi^{-1}(q)$ with $\nu(q_1)=r$. $q_1$, if it exists,
is a 2 point.

Suppose that $D\subset \pi^{-1}(C)\cap \overline S_r(X_1)$ is a curve, which is
necessarily a section over $C$, and $q\in C$ is a 2 point such that $\nu(q)=r-1$.
Suppose there exists a 2 point $q'\in\pi^{-1}(q)$ such that $q'\in D$. Then $\nu(q')=r-1$
by Lemma \ref{Lemma3}, so that (by the proof of Lemma \ref{Lemma500}), there exist permissible
parameters $(x,y,z)$ at $q'$ such that
$$
\begin{array}{ll}
u&=(x^ay^b)^m\\
v&=P(x^ay^b)+x^cy^dF_{q'}
\end{array}
$$
$$
F_{q'}=y^{r-1}+x\Omega+zG
$$
and there exists an irreducible series $f(y,z)$ such that $\hat{\cal I}_{D,q'}=(x,f(y,z))$.

Case 1 or Case 2 of Lemma \ref{Lemma659} must hold. Suppose that Case 1 holds.
Set $x=0$ in the formula of Case 1 to get that there exists a series $h(\overline y,z)$
such that 
$$
\overline y^{ad-bc}(\overline y^{(r-1)a}+zG(0,\overline y^a,z))=hf(\overline y^a,z)^r
$$
$\overline y\not\,\mid f(\overline y^a,z)$ implies 
$$
a(r-1)\ge \nu(f(\overline y^a,0)^r)\ge ar
$$
a contradiction.

Now suppose that Case 2 of Lemma \ref{Lemma659} holds. $a(r-1)\ne bc-ad$ since
$F_{q'}$ is normalized. Set $x=0$ in the formula of Case 2 to get that there exists a series $h(\overline y,z)$
such that 
$$
\overline y^{a(r-1)}+zG(0,\overline y^a,z)-g(0)\overline y^{bc-ad}=hf(\overline y^a,z)^r
$$
$$
0\ne \overline y^{a(r-1)}-g(0)\overline y^{bc-ad}=h(\overline y,0)f(\overline y^a,0)^r
$$
Thus 
$$
a(r-1)\ge \nu(\overline y^{a(r-1)}-g(0)\overline y^{bc-ad})\ge r \nu(f(\overline y^a,0))\ge ra
$$
which is a contradiction.

\end{pf}

\begin{Lemma}\label{Lemma654}
Suppose that  $r\ge 2$ and $C\subset \overline S_r(X)$ is a
curve containing a 1 point such that $C$ is r big. let $\pi:X_1\rightarrow X$ be the blowup 
of $C$.
\begin{enumerate}
\item Suppose that $p\in C$ is a 1 point with $\nu(p)=\gamma(p)=r$, and
$q\in\pi^{-1}(p)$. Then
\begin{enumerate}
\item If $q$ is a 1 point then $\gamma(q)\le r$. There is at most one 1 point $q\in\pi^{-1}(p)$ such that $\gamma(q)>r-1$.
\item If $q$ is a 2 point then $\nu(q)=0$.
\end{enumerate}
\item Suppose that $p\in C$ is a 1 point with $\nu(p)=r$, $\gamma(p)\ne r$,
and $q\in\pi^{-1}(p)$. Then
\begin{enumerate}
\item If $q$ is a  1 point then $\gamma(q)<r$,
\item If $q\in\pi^{-1}(p)$ is a 2 point then $\nu(q)\le r-1$.
\end{enumerate}
\item Suppose that $p\in C$ is a 2 point such that $\gamma(p)=\nu(p)=r$, and
$q\in \pi^{-1}(p)$. Then
\begin{enumerate}
\item If $q$ is a 2 point then $\nu(q)\le r$ and $\gamma(q)\le r$.
\item There is at most one 2 point $q\in \pi^{-1}(p)$ such that $\gamma(q)>r-1$.
\item If $q$ is a 3 point then $\nu(q)=0$.
\end{enumerate}
\item Suppose that $p\in C$ is a 2 point with $\nu(p)=r$ and $\tau(p)>0$, and
$q\in\pi^{-1}(p)$. Then
\begin{enumerate}
\item If $q$ is a 2 point then $\gamma(q)\le r$.
\item If $q$ is the 3 point then $\nu(q)\le r-\tau(p)$.
\end{enumerate}
\end{enumerate}
\end{Lemma}

\begin{pf}
 First suppose that $p\in C$ is a 1 point such that $\nu(p)=\gamma(p)=r$.
 We have permissible
parameters $(x,y,z)$ at $p$ such that $\hat{\cal I}_{C,p}=(x,z)$, 
\begin{equation}\label{eq973}
\begin{array}{ll}
u&=x^a\\
v&=P(x)+x^bF_p\\
F_p&=\tau z^r+\sum_{i=2}^r\overline a_i(x,y)x^iz^{r-i}
\end{array}
\end{equation}
where $\tau$ is a unit by  Lemma \ref{Lemma17}.

Suppose that $q\in\pi^{-1}(p)$ and $q$ is a 1 point. Then $q$ has permissible
parameters $(x_1,y,z_1)$ such that
$x=x_1$, $z=x_1(z_1+\alpha)$. Then $\nu(F_q(0,0,z_1))\le r$
and $\nu(F_q(0,0,z_1))<r$ if $\alpha\ne 0$.

If $q\in\pi^{-1}(p)$ is the 2 point then $q$ has permissible parameters
$(x_1,y,z_1)$ such that $x=x_1z_1, z=z_1$. Then $F_q=\frac{F_p}{z_1^r}$ is a unit.

Now suppose that $p\in C$ is a 1 point with $\nu(p)=r$ and $\gamma(p)\ne r$.
Suppose that $q\in\pi^{-1}(p)$ is a 1 point. Then there exist permissible parameters
$(x,y,z)$ at $p$ such that $x=z=0$ are local equations of $C$ at $p$, and permissible
parameters $(x_1,y,z_1)$ at $q$ such that $x=x_1,z=x_1z_1$.
$$
F_p=\sum_{i+j\ge r}a_{ij}(y)x^iz^j
$$
where $a_{r0}(0)=0$, $a_{0r}(0)=0$, and $a_{ij}(0)\ne 0$ for some $i,j$ with
$i+j=r$.
$$
F_q=\frac{F_p}{x_1^r}=\left(\sum_{i+j=r}a_{ij}(0)z_1^j\right)+x_1\Omega+yG
$$
implies $\nu(F_q(0,0,z_1))\le r-1$.

At the 2 point $q\in\pi^{-1}(p)$, there exist permissible parameters
$(x,y,z)$ as above, and permissible parameters $(x_1,y,z_1)$ at $q$ such that $x=x_1z_1$,
$z=z_1$,
$$
F_q=\frac{F_p}{z_1^r}=\sum_{i+j=r}a_{ij}(0)x_1^i+z_1\Omega+yG
$$
where $a_{ij}(0)\ne 0$ for some $i\le r-1$.

Now suppose that $p\in C$ is a 2 point such that
 $\nu(p)=r$ and $\gamma(p)=r$.
$$
\begin{array}{ll}
u&=(x^ay^b)^m\\
v&=P(x^ay^b)+x^cy^dF
\end{array}
$$
 $\hat{\cal I}_{C,p}=(x,z)$.
After a permissible change of parameters, we have by Lemma \ref{Lemma17} 
\begin{equation}\label{eq516}
F=\tau z^r+a_2(x,y)z^{r-2}+\cdots+a_r(x,y)
\end{equation}
where $\tau$ is a unit and $x^i\mid a_i$ for all $i$.

If $p_1\in \pi^{-1}(p)$ has permissible parameters $(x_1,y,z_1)$ with
$$
x=x_1z_1, z=z_1
$$
then
$$
F_1=\tau+x_1\Omega
$$
so that $p_1$ is resolved.
 Suppose that $p_1\in\pi^{-1}(p)$
has regular parameters 
$$
x=x_1,
z=x_1(z_1+\alpha)
$$
$$
\frac{F}{x_1^r}=\tau(z_1+\alpha)^r+\frac{a_2(x,y)}{x^2}(z_1+\alpha)^{r-2}+\cdots
+\frac{a_r(x,y)}{x^r}
$$
Thus $\nu(F_1(0,0,z_1))\le r$ and $\nu(F_1(0,0,z_1))\le r-1$ if $\alpha\ne 0$.

Suppose that $p\in C$ is a 2 point such that $\nu(p)=r$ and $\tau(p)>0$.
$$
\begin{array}{ll}
u&=(x^ay^b)^m\\
v&=P(x^ay^b)+x^cy^dF\\
\end{array}
$$
where $F\in\hat{\cal I}_{C,p}^r=(x,z)^r$.
$$
F = \sum_{i+k\ge r}a_{ijk}x^iy^jz^k.
$$
 Suppose that $q\in\pi^{-1}(p)$
is a 2 point. After a permissible change of parameters, replacing $z$ with
$z-\alpha x$, $q_1$ has permissible parameters $(x_1,y_1,z_1)$ such that
$$
x=x_1, z=x_1z_1
$$
$$
\begin{array}{ll}
u&=(x_1^ay^b)^m\\
v&=P(x_1^ay_1^b)+x_1^{c+r}y^dF_q
\end{array}
$$
$$
F_q=\frac{F}{x_1^r}=\sum_{i+k\ge r} a_{ijk}x_1^{i+k-r}y^jz_1^k
$$
$$
F_q=\sum_{i+k=r}a_{i0k}z_1^k+x_1\Omega_1+yG
$$
Thus  $\gamma(q)\le r$.

Now suppose that $q\in\pi^{-1}(p)$ has permissible parameters
$$
x=x_1z_1, z=z_1
$$
so that $q$ is a 3 point.
$$
\begin{array}{ll}
u&=(x_1^ay^bz_1^a)^m\\
v&=P(x_1^ay^bz_1^a)+x_1^{c}y^dz_1^{c+r}F_q
\end{array}
$$

$$
F_q=\frac{F}{z_1^r}=\sum_{i+k\ge r}a_{ijk}x_1^iy^jz_1^{i+k-r}
$$
$$
F_q=\sum_{i+k=r}a_{i0k}x_1^i+yG+z_1\Omega
$$
$a_{r-k,0,k}\ne0$ if $k=\tau(p)$
 which implies that 
$\nu(q)\le r-\tau(p)$.
\end{pf}

\begin{Lemma}\label{Lemma975}
Suppose that $r\ge 2$, $C\subset\overline S_r(X)$ is a  curve  containing a 1 point such that $C$ is r small.
\begin{enumerate}
\item Let $\pi:Y\rightarrow X$ be the monodial transform centered at $C$.
\begin{enumerate}
\item Suppose that $p\in C$ is a generic point. If $q\in \pi^{-1}(p)$ is a 1 point
then $\nu(q)=1$. If $q\in \pi^{-1}(p)$ is the 2 point then $\nu(q)\le r$
and $\nu(q)=r$ implies $\tau(q)>0$.
\item Suppose that $p\in C$ is a 2 point such that $\nu(q)=r-1$. If $q\in\pi^{-1}(p)$
is a 2 point then $\nu(q)=0$. If $q\in\pi^{-1}(p)$ is a 3 point then $\nu(q)\le r-1$.
\end{enumerate}
\item Suppose that $p\in C$ is a 1 point such that $\nu(p)=r$ or a 2 point such that
$\nu(p)=r$ and $\tau(p)>0$. Then there exists a finite sequence of quadratic
transforms $\sigma:Z\rightarrow X$ centered at points over $p$ such that 
if $q\in\sigma^{-1}(p)$ is a 1 point then $\nu(q)\le r$. $\nu(q)=r$ implies $\gamma(q)=r$.
If $q\in \sigma^{-1}(p)$ is a 2 point then $\nu(q)\le r$. $\nu(q)=r$ implies $\tau(q)>0$.
If $q\in\sigma^{-1}(p)$ is a 3 point then $\nu(q)\le r-1$. The strict transform of $C$
intersects $\sigma^{-1}(p)$ in a 2 point $p'$ such that $\nu(p')=r-1$
\end{enumerate}
\end{Lemma}

\begin{pf}
Suppose that $p\in C$ is a 2 point. By Lemma \ref{Lemma7}, there are permissible parameters $(x,y,z)$ at $p$ 
with $\hat{\cal I}_{C,p}=(x,z)$ such that 
\begin{equation}\label{eq156}
\begin{array}{ll}
u&=(x^ay^b)^m\\
v&=P(x^ay^b)+x^cy^dF_p\\
F_p&= x^{r-1}y^n+\sum_{i+k\ge r}a_{ijk}x^iy^jz^k
\end{array}
\end{equation}
with $n\ge 0$.  Suppose that $\nu(p)=r$ and $\tau(p)>0$. Then $n>0$ 
and $a_{i0k}\ne 0$ for some $i,k$ with $i+k=r$ and $k>0$.
Let $\pi':X'\rightarrow X$ be the blowup of $p$.  
Perform $n$ quadratic transforms,
$\pi_1:X_1\rightarrow X$, centered at the 2 point which is the intersection of the strict transform
of $C$ and the exceptional divisor. Then by Theorem \ref{Theorem13}
\begin{enumerate}
\item All 1 points $q$ in $\pi_1^{-1}(p)$ with $\nu(q)=r$ have $\gamma(q)=r$.
\item All 2 points $q\in\pi^{-1}(p)$ with $\nu(q)=r$ have $\tau(q)>0$.
\item All 3 points $q\in\pi^{-1}(p)$ have $\nu(q)\le r-1$.
\end{enumerate}

If $C_1$ is the strict transform of $C$, and $q$ is the exceptional point on $C_1$, then
there are permissible parameters $(x_1,y_1,z_1)$ at $q$ such that
$$
x=x_1y_1^n,
y=y_1,
z=z_1y_1^n
$$
\begin{equation}\label{eq81}
\begin{array}{ll}
u&=(x_1^ay_1^{na+b})^m\\
F_{q}&= x_1^{r-1}+\sum_{i+k\ge r}a_{ijk}x_1^iy_1^{n(i+k-r)+j}z_1^k
\end{array}
\end{equation}
where $\hat{\cal I}_{C_1,q}=(x_1,z_1)$.
Thus $\nu(q)=r-1$.

Suppose that $p\in C$ is a 1 point with $\nu(p)=r$. Then by Lemma \ref{Lemma5} there are regular parameters $(x,y,z)$ in $\hat{\cal O}_{X,p}$
such that $\hat{\cal I}_{C,p}=(x,z)$, 
\begin{equation}\label{eq83}
\begin{array}{ll}
u&=x^a\\
F_p&=x^{r-1}y^n+\sum_{i+k\ge r}a_{ijk}x^iy^jz^k
\end{array}
\end{equation}
with $n\ge 1$. There are only finitely many 1 points in $C$
 such that $n>1$.

Suppose that $n>1$. Then $a_{i0k}\ne 0$ for some $a_{i0k}$ with 
$i+k=r$ and $k>0$, so that
$\tau(p)>0$.
Let $\lambda:Z\rightarrow X_2$ be the sequence of $n$ quadratic transforms centered 
first at $p$, and then
at the intersection of the strict transform of $C$ and the exceptional fiber. 

Let $C'$ be the strict transform of $C$ on $Z$. Let $q'$ be the exceptional point of $\lambda$ on $C'$. By Theorems  \ref{Theorem9} and \ref{Theorem13}, the conclusions of 2. of the Theorem
hold at all points above $p$, except possibly at $q'$.
$q'$ has permissible parameters $(x_1,y_1,z_1)$ such that
$$
x=x_1y_1^n,
y=y_1,
z=z_1y_1^n.
$$
$$
\begin{array}{ll}
u&=(x_1y_1^n)^{a}\\
F_{q'}&=\frac{F_{q}}{y_1^{nr}}=x_1^{r-1}+\sum_{i+k\ge r}a_{ijk}x_1^iy_1^{(i+k-r)n+j}z_1^k
\end{array}
$$
Thus  $\nu(q')=r-1$ and $C'$ has the form (\ref{eq156}) with $n=0$ at $q'$.

Let $\pi:Y\rightarrow X$ be the blowup of $C$. Suppose that $p\in C$, and $p$ is a 2 point such that $\nu(p)=r-1$ so that (\ref{eq156}) with $n=0$ holds at $p$. 
Suppose that $q\in \pi^{-1}(p)$, and $q$ has permissible parameters $(x_1,y_1,z_1)$
such that
$$
x=x_1,
z=x_1(z_1+\alpha)
$$
After making a permissible change of variables, replacing $z$ with $z-\alpha x$, we may assume that $\alpha=0$.
Then $F_{q}=\frac{F_{p}}{x_1^{r-1}}$, so that $\nu(q)=0$.

Suppose that $q\in \pi^{-1}(p)$, and $q$ has permissible parameters $(x_1,y,z_1)$
such that
$$
x=x_1z_1,
z=z_1
$$
$$
\begin{array}{ll}
u&=(x_1^ay_1^bz_1^a)^m\\
F_{q}&= \frac{F_{p}}{z_1^{r-1}} =x_1^{r-1}+\sum_{i+k\ge r}a_{ijk}x_1^iy^{j}z_1^{i+k-(r-1)}
\end{array}
$$
so that $\nu(q)\le r-1$.

Now suppose that $p\in C$ is a generic point, so that  (\ref{eq83}) holds with $n=1$
at $p$. 
Suppose that $q\in \pi^{-1}(p)$, and $q$ has permissible parameters $(x_1,y,z_1)$
such that
$$
x=x_1,
z=x_1(z_1+\alpha)
$$
After making a permissible change of variables, replacing $z$ with $z-\alpha x$, we may assume that $\alpha=0$.
Then $F_{q}=\frac{F_{p}}{x_1^{r-1}}$, so that $\nu(q)=1$.

Suppose that $q\in \pi^{-1}(p)$, and $q$ has permissible parameters $(x_1,y,z_1)$
such that
$$
x=x_1z_1,
z=z_1
$$
$$
\begin{array}{ll}
u&=(x_1z_1)^a\\
F_{q}&= \frac{F_{p}}{z_1^{r-1}} =x_1^{r-1}y+\sum_{i+k\ge r}a_{ijk}x_1^iy^{j}z_1^{i+k-(r-1)}
\end{array}
$$
so that $\nu(q)\le r$, $\nu(q)=r$ implies $\tau(q)>0$.

\end{pf}

\begin{Lemma}\label{Lemma41}  Suppose that $r\ge 2$, 
$C\subset\overline S_r(X)$ is a curve containing a 1 point such that $C$ is r small and
$\gamma(q)=r$ for $q\in C$.
\begin{enumerate}
\item Let $\pi:Y\rightarrow X$ be the monoidal transform centered at $C$.
\begin{enumerate}
\item Suppose that $p\in C$ is a generic point. Then $\nu(q)\le 1$ if $q\in\pi^{-1}(p)$.
\item Suppose that $p\in C$ is a 2 point such that $\nu(p)=r-1$. Suppose that $q\in\pi^{-1}(p)$.
Then $\nu(q)=0$ if $q$ is a 2 point, and $\nu(q)\le 1$ if $q$ is a 3 point.
\end{enumerate}
\item Suppose that $p\in C$. Then there exists a finite sequence of quadratic
transforms $\sigma:Z\rightarrow X$ centered at points over $p$ such that $\nu(q)\le r$,
and $\gamma(q)\le r$ if $q\in\sigma^{-1}(p)$ is a 1 or 2 point. $\nu(q)=0$ if $q$ is a 
3 point, and the strict transform of $C$ intersects $\sigma^{-1}(p)$ in a 2 point
$p'$ such that $\nu(p')=r-1$ and $\gamma(p')=r$.
\end{enumerate}
\end{Lemma}

\begin{pf} Suppose that $p\in C$ is a 2 point. By Lemma \ref{Lemma7}, there exist
permissible parameters $(x,y,z)$ at $p$ such that $x=z=0$ are local equations of
$C$ at $p$.  
\begin{equation}\label{eq950}
\begin{array}{ll}
u&=(x^ay^b)^m\\
v&=P(x^ay^b)+x^cy^dF_p
\end{array}
\end{equation}
$$
F_{p}=\tau'(y)x^{r-1}y^n+\sum_{i+k\ge r}a_{ijk}x^iy^jz^k
$$
with $\tau'$ a unit, $n\ge 0$, and $a_{00r}\ne 0$.

Suppose that $\nu(p)=r-1$. Then $n=0$ and $\tau(p)=0$.
Let $\pi:Y\rightarrow X$ be the monoidal transform centered at $C$.

If $q\in\pi^{-1}(p)$ is a 2 point, then after a permissible change of parameters at $p$,
we have that $q$ has permissible parameters $(x_1,y,z_1)$ such that $x=x_1$, $z=x_1z_1$.
$$
F_q=\frac{F_p}{x_1^{r-1}}=\tau'(y)+x_1\Omega
$$
so that $\nu(q)=0$.

If $q\in\pi^{-1}(p)$ is the 3 point, there exist permissible parameters $(x_1,y,z_1)$
at $q$ such that $x=x_1z_1$, $z=z_1$.
$$
F_q=\frac{F_p}{z_1^{r-1}}=\tau'(y)x_1^{r-1}+\sum_{i+k\ge r}a_{ijk}x_1^iy^{j}z_1^{i+k-(r-1)}
$$
$a_{00r}\ne 0$ implies $\nu(q)\le 1$.

Suppose that $p$ is a 2 point and $\nu(p)=r$. Let $\sigma:Y\rightarrow X$ be the quadratic transform with
center $p$. Suppose that $q\in\sigma^{-1}(p)$ is a 1 point or a 2 point.
Then by Theorem \ref{Theorem13}, $\nu(q)\le r$ and $\gamma(q)\le r$.
If $q$ is a 3 point then $\nu(q)=0$. At the 2 point $q$ on the strict transform of
$C$, we have permissible parameters $(x_1,y_1,z_1)$ such that $x=x_1y_1$,
$y=y_1$, $z=x_1y_1$. $x_1=z_1=0$ are local equations of the strict transform of $C$
at $q$. 
$$
F_q=\tau'(y)x_1^{r-1}y_1^{n-1}+\sum_{i+k\ge r}a_{ijk}x_1^iy_1^{i+j+k-r}z_1^k
$$
of the form of (\ref{eq950}) with $n$ decreased by 1.

By induction on $n$, we achieve the conclusion of 2. after a finite sequence of
quadratic transforms.

Suppose that $p\in C$ is a 1 point. By Lemma \ref{Lemma5}, there exist
permissible parameters $(x,y,z)$ at $p$ such that $x=z=0$ are local equations
of $C$ at $p$, 
\begin{equation}\label{eq647}
\begin{array}{ll}
u&=x^a\\
v&=P(x)+x^cF_p\\
F_p&=x^{r-1}\tau'(y)y^n+\sum_{i+k\ge r}a_{ijk}x^iy^jz^k
\end{array}
\end{equation}
with $\tau'$ a unit, $n\ge 1$, $a_{00r}\ne 0$.

Suppose that $p\in C$ is a generic point so that $n=1$.
Let $\pi:Y\rightarrow X$ be the monoidal transform centered at $C$.
 If $q\in\pi^{-1}(p)$ is a 
1 point, then after making a permissible change of parameters at $p$, there
are permissible parameters $(x_1,y,z_1)$ at $q$ such that 
$x=x_1,z=x_1z_1$.
$$
F_q=\frac{F_p}{x_1^{r-1}}=y\tau'(y)+x_1\Omega
$$
implies $\nu(q)=1$.
If $q\in\pi^{-1}(p)$ is the 2 point, then there are permissible parameters
$(x_1,y,z_1)$ at $q$ such that $x=x_1z_1$, $z=z_1$.
$$
F_q=\frac{F_p}{z_1^{r-1}}=x_1^{r-1}\tau'(y)y+\sum_{i+k\ge r}a_{ijk}
x_1^iy_1^jz_1^{i+k-(r-1)}
$$
which implies that $\nu(q)\le 1$ since $a_{00r}\ne 0$.

Suppose that $n\ge 2$ in (\ref{eq647}). Let $\sigma:Y\rightarrow X$ be the
quadratic transform with center $p$.

Suppose that $q\in\sigma^{-1}(p)$. Then $q$ is a 1 or 2 point, and $\nu(q)\le r$,
$\gamma(q)\le r$
by Theorem \ref{Theorem13}.

At the 2 point $q\in\pi^{-1}(p)$ which is contained in the strict transform of $C$,
there are permissible parameters $(x_1,y_1,z_1)$ at $q$ such that
$x=x_1y_1, y=y_1, z=y_1z_1$.
$$
\begin{array}{ll}
u&=(x_1y_1)^a\\
F_q&=\frac{F_p}{y_1^r}=x_1^{r-1}y_1^{n-1}\tau'(y_1)+\sum_{i+k\ge r}a_{ijk}x_1^i
y_1^{i+j+k-r}z_1^k
\end{array}
$$
The strict transform of $C$ has local equations $x_1=z_1=0$ at $q$. We are thus at  a
point of the form of (\ref{eq950}) with $n$ decreased by 1.

We thus achieve the conclusions of 2. after a finite number of quadratic transforms.
\end{pf}

\begin{Lemma}\label{Lemma97} Suppose that $r=2$ in Lemma \ref{Lemma41},
$C\subset\overline S_2(X)$ is a curve containing a 1 point such that $C$ is 2 small,
$\gamma(p)=2$ if $p\in C$, $\nu(p)=1$ if $p\in C$ is a 2 point and $p$ is a generic
point of $C$ ($n=1$ in (\ref{eq647})) if $p\in C$ is a 1 point. Let
$\pi:Y\rightarrow X$ be the monodial transform centered at $C$.
Suppose that there exists a 2 point $p\in C$ such that $\nu(p)=r-1=1$, and
$q\in\pi^{-1}(p)$ is a 3 point such that $\nu(q)=1$, or $p\in C$ is a generic point
of $C$ ($n=1$ in (\ref{eq647})) and $q\in\pi^{-1}(p)$ is a 2 point such that $\nu(q)=1$. Let $\overline C$ be the
2 curve through $q$ which is a section over $C$. Then $F_{q'}\in\hat{\cal I}_{\overline C,
q'}$ for all $q'\in\overline C$.

Suppose that there does exist a 2 curve $\overline C$ which is a section over $C$ such that
$F_{q'}\in\hat{\cal I}_{\overline C,q'}$ for $q'\in\overline C$. Let
$\pi_1:Z\rightarrow Y$ be the blowup of $\overline C$. Then 
\begin{enumerate}
\item Suppose that $q\in\overline C$
is a 2 point such that $q\in\pi^{-1}(p)$ where $p$ is a generic point of $C$
($n=1$ in (\ref{eq647})), and $q'\in\pi_1^{-1}(q)$.
\begin{enumerate}
\item If $q'$ is a 1 point then $\nu(q')=1$.
\item If $q'\in\pi_1^{-1}(q)$ is a 2 point then $\gamma(q')\le 1$.
\end{enumerate}
\item Suppose that $q\in\overline C$ is a 
3 point such that $q\in\pi^{-1}(p)$ where $p\in C$ is a 2 point such that $\nu(p)=1$ and $q'\in\pi^{-1}(q)$.
\begin{enumerate}
\item If $q'$ is a 2 point then $\gamma(q')\le 1$.
\item If $q'$ is a 3 point then $\nu(q')=0$.
\end{enumerate}
\end{enumerate}
\end{Lemma}
\begin{pf}
If $p\in C$  is a 2 point with $\nu(p)=1$ (and $\gamma(p)=2$), then there exist permissible parameters $(x,y,z)$ at $p$ such that
$$
\begin{array}{ll}
u&=(x^ay^b)^m\\
v&=P(x^ay^b)+x^cy^dF_p\\
F_p&=x+z^2
\end{array}
$$
where $x=z=0$ are local equations of $C$ at $p$. There exist permissible parameters
$(x_1,y,z_1)$ at the 3 point $q\in\pi^{-1}(p)$ such that $x=x_1z_1, z=z_1$. 
\begin{equation}\label{eq1007}
\begin{array}{ll}
u&=(x_1^ay^bz_1^a)^m\\
v&=P(x_1^ay^bz_1^a)+x_1^cy^dz_1^{c+1}F_q\\
F_q&=x_1+z_1
\end{array}
\end{equation}
and $x_1=z_1=0$ are local equations of $\overline C$ at $q$. We have
$F_q\in\hat{\cal I}_{\overline C,q}$ which implies $F_{q'}\in\hat{\cal I}_{\overline C,q'}$
if $q'\in\overline C$ by Lemma \ref{Lemma655}.

If $p\in C$ is a generic point, then $p$ is a 1 point and there exist permissible parameters $(x,y,z)$ at $p$
such that 
$$
F_p=xy+\sum_{i+k\ge2} a_{ijk}x^iy^jz^k
$$
with $a_{002}\ne 0$
and $x=z=0$ are local equations of $C$ at $p$. There exist permissible parameters
$(x_1,z_1,y)$ at the 2 point $q\in\pi^{-1}(p)$ such that $x=x_1z_1, z=z_1$.
$$
F_q=x_1y_1+z_1(\sum_{i+k=2}a_{ijk}x_1^iy_1^j)+z_1^2\Omega
$$
and $x_1=z_1=0$ are local equations of $\overline C$ at $q$. Since $a_{002}\ne 0$, there exist
permissible parameters $(x_1,\overline z_1,y_1)$ at $q$ such that 
\begin{equation}\label{eq1008}
F_q=x_1y_1+\overline z_1
\end{equation}
and $x_1=\overline z_1=0$ are local equations of $\overline C$ at $q$.

 We have $F_q\in \hat{\cal I}_{
\overline C,q}$ implies $F_{q'}\in\hat{\cal I}_{\overline C,q'}$ for all
$q'\in\overline C$ by Lemma \ref{Lemma655}.

1. follows from (\ref{eq1008}).

Suppose that $q\in\overline C$ is a 3 point, with permissible parameters $(x_1,y,z_1)$
such that (\ref{eq1007}) holds at $q$. Suppose that $q'\in\pi_1^{-1}(q)$. If $q'$ is a
3 point, then $\nu(q')=0$. Suppose that $q'$ is a 2 point. Then there exist regular
parameters $(x_2,y,z_2)$ in $\hat{\cal O}_{Z,q'}$ such that
$$
x_1=x_2, z_1=x_2(z_2+\alpha)
$$
with $\alpha\ne 0$.
$$
u=(x_2^{2a}y^b(z_2+\alpha)^a)^m=(\overline x_2^{2a}y^b)^m=(\overline x_2^{\overline a}
y^{\overline b})^{\overline m}
$$
where $x_2=\overline x_2(z_2+\alpha)^{-\frac{1}{2}}$, $(\overline a,\overline b)=1$.
$$
v=P_{q'}(\overline x_2^{\overline a}y^{\overline b})+\overline x_2^{2c+2}y^d(1+\alpha+z_2)
$$
Thus $\gamma(q')\le 1$ and 2. follows.

\end{pf}

\section{Power series in 2 variables}

\begin{Lemma}\label{Lemma1026}
Suppose that $R=k[[x,y]]$ is a power series ring in  two variables and $u(x,y),v(x,y)\in R$ 
are series. Suppose that $R\rightarrow R'$ is a quadratic transform. Set $R_1=\hat{R'}$.
Then $(u,v)$ are analytically independent in $R_1$ if and only if $u$ and $v$ are
analytically independent in $R$.
\end{Lemma}

\begin{pf} By Zariski's Subspace Theorem (Theorem 10.6 \cite{Ab5}), $R\rightarrow R_1$ is an inclusion, and the
Lemma follows.
\end{pf}

\begin{Lemma}\label{Lemma1027}
Suppose that $R=k[[x,y]]$ is a power series ring in two variables over an algebraically
closed field $k$ of characteristic 0 and
$u(x,y), v(x,y)\in R$ are series such that either
$$
u=x^a
$$
or
$$
u=(x^ay^b)^m
$$
with $(a,b)=1$. Then $u$ and $v$ are analytically dependent if and only if
there exists a series $p(t)$ such that $v=p(x)$ in the first case and
$v=p(x^ay^b)$ in the second case.
\end{Lemma}

\begin{pf} First suppose that $u=x^a$ and $v=p(x)$ is a series. Let $\omega$ be a primitive 
$a$-th root of unity. 
$$
0=\prod_{i=0}^{a-1}(v-p(\omega^ix))\in k[[u,v]]
$$
implies $u$ and $v$ are analytically dependent.

Now suppose that $u=x^a$ and $u,v$ are analytically dependent.  Suppose that $v$ is not
a series in $x$. Write
$$
v=q(x)+x^bF
$$
where $q(x)$ is a polynomial, $x\not\,\mid F$ and $F(0,y)$ is a nonzero series with no constant 
term.
$$
F(0,y)=y^r\mu(y)
$$
for some $r>0$ where $\mu(y)$ is a unit series. $x$ and $x^bF$ are thus analytically
dependent, and $x$ and $F$ are analytically dependent. There  exists an 
irreducible series
$$
p(s,t)=\sum a_{ij}s^it^j
$$
such that 
$$
0=\sum a_{ij}x^iF^j
$$
which implies that
$$
0=\sum a_{0j}F(0,y)^j=\sum a_{0j}y^{rj}\mu(y)^j,
$$
a contradiction, since $p(s,t)$ irreducible implies some $a_{0j}\ne 0$. Thus
$v$ is a series in $x$.

Now suppose that $u=(x^ay^b)^m$ and $v=p(x^ay^b)$ is a series in $x^ay^b$. Let
$\omega$ be a primitive $m$-th root of unity.
$$
0=\prod_{i=0}^{m-1}(v-p(\omega^ix^ay^b))\in k[[u,v]]
$$
implies that $u,v$ are analytically dependent.

Suppose that 
$$
u=(x^ay^b)^m
$$
and $u,v$ are analytically dependent.
Consider the quadratic transform 
$$
R\rightarrow R_1=R[x_1,y_1]_{(x_1,y_1)}
$$
where $x=x_1, y=x_1(y_1+1)$. $\hat{R_1}$ has regular parameters $(\overline x_1,y_1)$
where
$$
x_1=\overline x_1(y_1+1)^{-\frac{b}{a+b}}.
$$
Thus in $\hat R_1$, $u=\overline x_1^{(a+b)m}$. Since $u,v$ must be analytically
dependent in $\hat{R_1}$, there exists a series $q(\overline x_1)$ such that
$v=q(\overline x_1)$, by the first part of the proof.

Suppose that $v$ is not a series in $x^ay^b$. Write $v=\sum a_{ij}x^iy^j$.
There exists a smallest $r$ such that there exists $i_0,j_0$ such that $i_0+j_0=r$,
$bi_0-aj_0\ne 0$ and $a_{i_0j_0}\ne 0$.
$$
v=\sum_{i+j<r, aj-bi=0}a_{ij}\overline x_1^{i+j}
+\overline x_1^r(\sum_{i+j=r}a_{ij}(y_1+\alpha)^{j-\frac{br}{a+b}})+\overline x_1^{r+1}\Omega.
$$
$$
\sum_{i+j=r}a_{ij}(y_1+\alpha)^{j-\frac{br}{a+b}}\in k
$$
implies
$$
(y_1+\alpha)^{\frac{-br}{a+b}}(\sum_{i+j=r} a_{ij}(y_1+\alpha)^j)=c\in k
$$
so that
$$
\sum_{i+j=r} a_{ij}(y_1+\alpha)^j=c(y_1+\alpha)^{\frac{br}{a+b}}.
$$
Thus
$$
\frac{br}{a+b}\in\{0,1,\cdots,r\}
$$
and  $j_0=\frac{br}{a+b}$. This implies that $aj_0-bi_0=0$, a contradiction. Thus
$v$ is a series in $x^ay^b$.
\end{pf}

\begin{Lemma}\label{Lemma1028} Suppose that $R=k[[x,y]]$ is a power series in 
two variables  over an algebraically closed field $k$ of characteristic 0, $u=x^a$ or $u=(x^ay^b)^m$, and $(u,v)$ are analytically independent.
let $\pi:X\rightarrow \text{spec}(R)$ be the blowup of $m=(x,y)$. Then for all but
finitely many points $q\in\pi^{-1}(m)$ there exist regular parameters
$(\overline x,\overline y)$ in $\hat{\cal O}_{X,q}$ such that there is an expansion
\begin{equation}\label{eq1029}
\begin{array}{ll}
u&=\overline x^a\\
v&=P(\overline x)+\overline x^b\overline y
\end{array}
\end{equation}
\end{Lemma}

\begin{pf} First suppose that $u=x^a$. Write $v=P(x)+x^bF$ where $x\not\,\mid F$
and $F$ has no terms which are powers of $x$. Write
$$
F=\sum_{i+j\ge r}a_{ij}x^iy^j
$$ 
where $r=\nu(F)$. There exists $j_0>0$ such that $i_0+j_0=r$ and $a_{i_0j_0}\ne 0$.
For all but one point $q\in\pi^{-1}(m)$ there are regular parameters $(x_1,y_1)$ in
$\hat{\cal O}_{X,q}$ such that 
$$
x=x_1, y=x_1(y_1+\alpha)
$$
with $\alpha\in k$.
\begin{equation}\label{eq1030}
\begin{array}{ll}
u&=x_1^a\\
v&=P(x_1)=x_1^{b+r}(\sum_{i+j=r}a_{ij}(y_1+\alpha)^j+x_1\Omega)
\end{array}
\end{equation}

$v$ has an expansion (\ref{eq1029}) if and only if
\begin{equation}
\begin{array}{ll}
\frac{d}{dy_1}(\sum_{i+j=r}a_{ij}(y_1+\alpha)^j)\mid_{y_1=0}=\sum_{j\le r}ja_{r-j,j}(-\alpha)^{j-1}\ne 0.
\end{array}
\end{equation}
Since
$$
\sum_{i+j=r}ja_{r-j,j}(-\alpha)^{j-1}
$$
has at most finitely many roots, all but finitely many $q\in\pi^{-1}(m)$ have an
expansion (\ref{eq1029}).

Now suppose that $u=(x^ay^b)^m$. Write 
$$
v=P(x^ay^b)+ x^cy^dF
$$
where $x,y\not\,\mid F$ and $x^cy^dF$ has no terms which are powers of $x^ay^b$. Write
$$
F=\sum_{i+j\ge r}a_{ij}x^iy^j
$$
where $r=\nu(F)$. 

For all but two points $q\in\pi^{-1}(m)$ there are regular parameters $(x_1,y_1)$
in $\hat{\cal O}_{X,q}$ such that
$$
x=x_1, y=x_1(y_1+\alpha)
$$
with $\alpha\ne 0$. There are regular parameters $(\overline x_1,y_1)$ in
$\hat{\cal O}_{X,q}$ such that
$$
x_1=\overline x_1(y_1+\alpha)^{-\frac{b}{a+b}}.
$$
$$
\begin{array}{ll}
u&=\overline x_1^{(a+b)m}\\
v&=P(\overline x_1^{(a+b)m})+\overline x_1^{c+d+r}(y_1+\alpha)^{\lambda}
\frac{F}{x_1^r}
\end{array}
$$
where 
$$
\lambda=d-\frac{b(c+d+r)}{a+b}
$$
$$
(y_1+\alpha)^{\lambda}\frac{F}{x_1^r}=\sum_{i+j=r}a_{ij}(y_1+\alpha)^{j+\lambda}
+\overline x_1\Omega.
$$
$v$ does not have an expression (\ref{eq1029}) at $q$ if and only if there exists
$c_\alpha\in k$ such that
$$
\sum_{j=0}^ra_{r-j,j}(y_1+\alpha)^j
\equiv c_{\alpha}(y_1+\alpha)^{-\lambda}\text{ mod }(y_1)^2.
$$
Set $a_j=a_{r-j,j}$. Suppose that $q$ does not have a form (\ref{eq1029}). Then
$$
\sum_{j=0}^ra_j\alpha^j=c_{\alpha}\alpha^{-\lambda}
$$
and
$$
\sum_{j=0}^rja_j\alpha^{j-1}=-c_{\alpha}\lambda \alpha^{-\lambda-1}
$$
implies 
\begin{equation}\label{eq1032}
(-\lambda)\sum_{j=0}^ra_j\alpha^j=\sum_{j=0}^rja_j\alpha^j.
\end{equation}
If there are infinitely many values of $\alpha$ satisfying (\ref{eq1032}), then
$(-\lambda-j)a_j=0$ for $0\le j\le r$, which implies that $-\lambda\in\{0,\ldots,r\}$
and the leading form of $F$ is
$$
L=\sum_{i+j=r}a_{ij}x^iy^j=a_{r+\lambda,-\lambda}x^{r+\lambda}y^{-\lambda}.
$$
Thus $x^cy^dF$ has a nonzero $x^{c+r+\lambda}y^{d-\lambda}$ term.
$$
\begin{array}{ll}
a(d-\lambda)-b(c+r+\lambda)&=ad-b(c+r)-(a+b)\lambda\\
&=ad-b(c+r)-(a+b)(d-\frac{b(c+d+r)}{a+b})\\
&=ad-b(c+r)-(a+b)d+b(c+d+r)=0
\end{array}
$$
which is impossible since $F$ is normalized (contains no terms which are powers of
$x^ay^b$). Thus there are at most a finite number of points $q\in\pi^{-1}(m)$ where
the form (\ref{eq1029}) does not hold.
\end{pf}

\begin{Theorem}\label{Theorem965}
Suppose that $k$ is an algebraically closed field of characteristic zero,
$B$ is a powerseries ring in 2 variables over $k$.
Suppose that $u,v\in B$ are analytically independent, and there exist regular
parameters $(x,y)$ in $B$ such that
$u=x^a$ or $u=x^ay^b$.
 Let $A=\text{spec}(B)$. Then there exists a sequence of quadratic transforms
$\pi:X\rightarrow A$ such that for all points $q\in X$,
there exist regular parameters $(\overline x,\overline y)$ in 
$\hat{\cal O}_{X,q}$ such that either 
\begin{equation}\label{eq976}
\begin{array}{ll}
u&=\overline x^a\\
v&=P(\overline x)+\overline x^b\overline y^c
\end{array}
\end{equation}
or 
\begin{equation}\label{eq977}
\begin{array}{ll}
u&=(\overline x^a\overline y^b)^m\\
v&=P(\overline x^a\overline y^b)+\overline x^c\overline y^d
\end{array}
\end{equation}
where $(a,b)=1$ and $ad-bc\ne 0$.
\end{Theorem}

Theorem \ref{Theorem965} will follow from Theorem \ref{T3}.
Throughout this section, we will use the notations of the statement of Theorem 
\ref{Theorem965}.

If $A\rightarrow  \text{spec}(k[[u,v]])$ is weakly prepared, then a stronger result
than the conclusions of Theorem \ref{Theorem965} are true in $B$.

\begin{Remark} With the assumptions of Theorem \ref{Theorem965}, further suppose that
$$
\sqrt{(\frac{\partial u}{\partial x}\frac{\partial v}{\partial y}-\frac{\partial u}{\partial y}
\frac{\partial v}{\partial x})}=\sqrt{(u)}.
$$
Then there exist regular parameters $(\overline x,\overline y)$ in $B$, and a power
series $P$ in $B$  such that
one of the following forms  holds. 
\begin{equation}\label{eq1069}
\begin{array}{ll}
u&=\overline x^a\\
v&=P(\overline x)+\overline x^c\overline y
\end{array} 
\end{equation}
\begin{equation}\label{eq1070} 
\begin{array}{ll}
u&=(\overline x^a\overline y^b)^m\\
v&=P(\overline x)+\overline x^c\overline y^d
\end{array} 
\end{equation}
where $(a,b)=1$ and $ad-bc\ne 0$.
\end{Remark}

\begin{pf} (7.4 \cite{AKi}) With our assumptions, one of the following must hold. 
\begin{equation}\label{eq1065}
\begin{array}{ll}
u&=x^a\\
u_xv_y-u_yv_x&=\delta x^e
\end{array}
\end{equation}
where $\delta$ is a unit
or 
\begin{equation}\label{eq1068}
\begin{array}{ll}
u&=(x^ay^b)^m\\
u_xv_y-u_yv_x&=\delta x^ey^f
\end{array}
\end{equation}
where $a,b,e,f>0$, $(a,b)=1$ and $\delta$ is a unit.

Write $v=\sum a_{ij}x^iy^j$.
First suppose that (\ref{eq1065}) holds. Then $ax^{a-1}v_y=\delta x^e$
implies we  have the form (\ref{eq1069}). Now suppose that (\ref{eq1068}) holds.
$$
u_xv_y-u_yv_x=\sum m(aj-bi)a_{ij}x^{am+i-1}y^{bm+j-1}
=\delta x^ey^f.
$$
Thus 
$$
v=\sum_{aj-bi=0}a_{ij}x^iy^j+\epsilon x^cy^d
$$
where $\epsilon$ is a unit. After making a change of variables, 
multiplying $x$ by a unit, and multiplying $y$ by a unit, we get the
form  (\ref{eq1070}).
\end{pf}

\begin{Definition} \label{Def1092} Suppose that $\Phi:X\rightarrow A$ is a product of quadratic
transforms, $p$ is a  point of $X$. We will say that $(u,v)$ are 1-resolved at $p$
if there exist regular parameters $(x,y)$ in $\hat{\cal O}_{X,p}$ such that one of the
forms (\ref{eq976}) or (\ref{eq977}) hold at $p$.
\end{Definition}

For the rest of this section, we will assume that
$$
\Phi:X\rightarrow A
$$
is a sequence of quadratic transforms.

 Suppose that $p\in X$ is a
point. 
 Then there are regular parameters $(x,y)$ of  $\hat{\cal O}_{X,p}$
 such that 
$u=x^{\overline a}y^{\overline b}$, and $\overline a>0$, $\overline b\ge 0$. 

Suppose that $\overline b>0$. 
 Let $m = (\overline a,\overline b)$, let $a = \frac{\overline a}{m}$, $b = \frac{\overline b}{m}$.
There are power series $P(t)$ and $F(x,y)$ such that $x$ does not divide $F$, $y$ does not divide $F$, $x^cy^dF$ has no
nonzero terms which are powers of $x^ay^b$ and
(in  $\hat{\cal O}_{X,p}$) 
\begin{equation}\label{I}
\begin{array}{l}
u=(x^{a}y^{b})^m\\
v=P(x^ay^b)+x^cy^dF(x,y)
\end{array}
\end{equation}
In this case, we will say that $p$ is a 2 point.

If $\overline b=0$, there are power series $P(t)$ and $F(x,y)$ such that $x$ does not divide $F$,
 $F$ has no nonzero terms which are
powers of $x$ and
(in  $\hat{\cal O}_{X,p}$) 
\begin{equation}\label{II}
\begin{array}{l}
u=x^{a}\\
v=P(x)+x^cF(x,y)
\end{array}
\end{equation}
In this case we will say that $p$ is a 1 point.

Suppose that $p\in X$, and $(x,y)$ are regular parameters in $\hat{\cal O}_{X,p}$
such that $(u,v)$ have one of the forms (\ref{I}) or (\ref{II}). Set 
\[
\overline \nu(p) = \left\{ \begin{array}{ll}
\text{mult}(F)-1 & \text{if $p$ is a 1 point}\\
\text{mult}(F) & \text{if $p$ is a 2 point}
\end{array}
\right.
\]
\begin{Lemma}\label{L3} 
$\overline \nu(p)$ is independent of the choice of regular parameters $(x,y)$ in (\ref{I}) or (\ref{II}).
\end{Lemma}

\begin{pf}
First suppose that $p$ is a 2 point. To express $u$ and $v$ in the form (\ref{I}) we can only  make a permissible change of variables in $x$ and $y$,
where a permissible change of variables is one of the following two forms: 
\begin{equation}\label{eqS1}
x=\omega_x\overline x, y=\omega_y\overline y\text{ where }\omega_x^{ma}\omega_y^{mb}=1
\end{equation}
or 
\begin{equation}\label{eqS2}
y=\omega_y\overline x, x=\omega_x\overline y\text{ where }\omega_x^{ma}\omega_y^{mb}=1
\end{equation}
where $\omega_x$, $\omega_y$ are unit series.
$\overline \nu(p)$ does not change after a change of variables of one of these forms.

Now suppose $p$ is a 1 point. To preserve the form (\ref{II}) we can only  make a permissible change of variables,
where a permissible change of variables is  of the form: 
\begin{equation}\label{eqS3}
x=\omega_x\overline x, y=\phi(\overline x,\overline y) \text{ where }\text{mult}(\phi(0,\overline y))=1.
\end{equation}
and  $\omega_x\in k$ is an $a$-th root of unity. Then
$$
\phi(\overline x,\overline y) = \overline \phi(\overline x,\overline y)(\overline y+\psi(\overline x))
$$
where $\overline \phi$ is a unit. Write 
$$
\psi(\overline x) = \sum b_i\overline x^i.
$$
 \begin{equation}\label{eqS4}
\begin{array}{l}
u=\overline x^{a}\\
v=\overline P(\overline x)+\overline x^c\overline F(\overline x,\overline y)
\end{array}
\end{equation}
where
$$
\overline F = \omega_x^c(F(\omega_x\overline x,\phi(\overline x,\overline y)) - F(\omega_x\overline x,\phi(\overline x,0)))
$$
$$
\overline P(\overline x) = P(\omega_x\overline x)+\overline x^c\omega_x^cF(\omega_x\overline x,\phi(\overline x,0))
$$
Suppose that the leading form of $F$ is
$$
L = \sum_{i+j= r}a_{ij}x^iy^j.
$$
The leading form $\overline L$ of $\overline F$ is then
$$
\omega_x^c(\sum_{i+j=r}a_{ij}\omega_x\overline x^i(e(\overline y-b_1\overline x))^j
-\sum_{i+j=r}a_{ij}\omega_x\overline x^i(-eb_1\overline x)^j)
$$
where $e=\overline\phi(0,0)$. $\overline L$  is nonzero since $a_{ij}\ne 0$ for some $j>0$.
\end{pf}

$(u,v)$ are 1-resolved at a 2 point $p$ if and only if $F$ is a unit.
$(u,v)$ are 1-resolved at a 1 point $p$ if and only if 
$F(x,y)=g(x,y)^d+h(x)$ for some series $g(x,y)$ with $\text{mult}(g(0,y))=1$, and
positive integer $d$.

\begin{Theorem}\label{T1} 
Suppose that  $g:X_1\rightarrow X$ is a quadratic transform, centered at a  point $p$ of $X$, and 
$p_1\in X_1$ is a  point such that $g(p_1)=p$.  Then
\[
\overline \nu(p_1)\le \overline \nu(p).
\]
If $(u,v)$ are 1-resolved at $p$  then $(u,v)$ are 1-resolved at $p_1$.
\end{Theorem}

\begin{pf}
First suppose that $p$ is a 2 point. Write 
$$
F = \sum_{i+j\ge r}a_{ij}x^iy^j.
$$
in $\hat{\cal O}_{X,p}$, where $r=\text{mult}(F)=\overline \nu(p)$.
Suppose that $\hat{\cal O}_{X_1,p_1}$ has regular parameters $(x_1,y_1)$ such that
$x=x_1, y=x_1(y_1+\alpha)$ with $\alpha\ne 0$.
Define $\overline x_1$ by
$$
x_1 = \overline x_1(y_1+\alpha)^{\frac{-b}{a+b}}.
$$
Then $(\overline x_1, y_1)$ are regular parameters in $\hat{\cal O}_{X_1,p_1}$.
$$
u=x_1^{m(a+b)}(y_1+\alpha)^{mb} = \overline x_1^{m(a+b)}.
$$
$$
v= P(\overline x_1^{a+b})+\overline x_1^{c+d+r}(y_1+\alpha)^{\lambda}(\frac{F}{x_1^r}) 
$$
where
$\lambda = d-\frac{b(c+d+r)}{a+b}$.
$$
F = \sum_{i+j=r}a_{ij}x_1^r(y_1+\alpha)^j+x_1^{r+1}\Omega.
$$
$$
\frac{F}{x_1^r} = \sum_{j=0}^ra_j(y_1+\alpha)^j+ x_1\Omega.
$$
where $a_j=a_{r-j,j}$.
We have
\begin{equation}
\begin{array}{l}
u=\overline x_1^{m(a+b)}\\
v=\overline P(\overline x_1)+\overline x_1^{c+d+r}\overline F(\overline x_1,y_1)
\end{array}
\end{equation}
where 
$$
\overline P = P(\overline x_1^{a+b})+\overline x_1^{c+d+r}(\alpha)^{\lambda}
(\frac{F(\alpha^{\frac{-b}{a+b}}\overline x_1,\alpha^{\frac{a}{a+b}}\overline x_1)}{\alpha^{\frac{-rb}{a+b}}\overline x_1^r})
$$
$$
\overline F = 
(y_1+\alpha)^{\lambda}
(\frac{F((y_1+\alpha)^{\frac{-b}{a+b}}\overline x_1,(y_1+\alpha)^{\frac{a}{a+b}}\overline x_1)}
{(y_1+\alpha)^{\frac{-rb}{a+b}}\overline x_1^r})
-
(\alpha)^{\lambda}
(\frac{F(\alpha^{\frac{-b}{a+b}}\overline x_1,\alpha^{\frac{a}{a+b}}\overline x_1)}{\alpha^{\frac{-rb}{a+b}}\overline x_1^r})
$$
 Set
$$
\beta = (\sum_{j=0}^ra_j\alpha^j)\alpha^{\lambda}.
$$
Suppose that $\overline \nu(p_1)> \overline \nu(p)=r$, so that 
$$
\text{mult}(\overline F)\ge \text{mult}(F)+2=r+2.
$$
 Then
$$
(y_1+\alpha)^{\lambda}(\sum_{j=0}^ra_j(y_1+\alpha)^j)-\beta\equiv 0\text{ mod }(y_1)^{r+2}.
$$
$\beta\ne 0$ since $\sum_{j=0}^ra_j(y_1+\alpha)^j\ne 0$.
We have 
\begin{equation}\label{eqS5}
\sum_{j=0}^ra_j(y_1+\alpha)^j\equiv \beta(y_1+\alpha)^{-\lambda}\text{ mod }(y_1)^{r+2}
\end{equation}

First suppose that $-\lambda\in\{0,1,\ldots,r\}$. Then
$$
\sum_{j=0}^ra_j(y_1+\alpha)^j=\beta(y_1+\alpha)^{-\lambda}
$$
where $t = -\lambda\le r$. Thus the leading form of $F$ is
$$
\begin{array}{lll}
L &=& \sum_{i+j=r}a_{ij}x_1^r(y_1+\alpha)^j\\
&=& \beta x_1^r(y_1+\alpha)^{-\lambda}\\
&=& \beta x^{r+\lambda}y^{-\lambda}\\
&=&\beta x^{r-t}y^t
\end{array}
$$
So the leading form of $F$ is $\beta x^{r-t}y^t$. Thus $\beta x^{c+r-t}y^{d+t}$ is a nonzero term of $x^cy^dF$.
Since
$$
t=\frac{b(c+d+r)}{a+b}-d,
$$
we have
$$
b(c+r-t)-a(d+t) = 0
$$
so that $x^{c+r-t}y^{d+t}$ is a power of $x^ay^b$, a contradiction.

We must then have $-\lambda \not\in \{0,1\ldots,r\}$. But then the $y_1^{r+1}$ coefficient of $\beta(y_1+\alpha)^{-\lambda}$
is non zero, a contradiction to (\ref{eqS5}).

Now suppose that $p$ is a 2 point and $\hat{\cal O}_{X_1,p_1}$ has regular parameters $(x_1,y_1)$ such that
$x=x_1, y=x_1y_1$. Write 
$$
F = \sum_{i+j\ge r}a_{ij}x^iy^j.
$$
Then 
\begin{equation}
\begin{array}{l}
u= x_1^{m(a+b)}y_1^{mb}\\
v=\overline P( x_1^{a+b}y_1^b)+ x_1^{c+d+r}y_1^d\overline F(x_1,y_1)
\end{array}
\end{equation}
where
$\overline P = P$, $\overline F = \frac{F}{x_1^r}$.
We need only check that
$ \frac{F}{x_1^r}$ has no nonzero $x_1^{\alpha}y_1^{\beta}$ terms with 
$b(c+d+r+\alpha) = (a+b)(d+\beta)$.
We have that $a_{ij} = 0$ if $b(c+i)-a(d+j)=0$.
$$
\frac{F}{x_1^r} = \sum a_{ij}x_1^{i+j-r}y_1^j.
$$
Suppose that $b(c+d+r+\alpha) = (a+b)(d+\beta)$. Set $i=\alpha-\beta+r$, $j=\beta$.
Then $b(c+i)-a(d+j)=0$, and $a_{ij}=0$. But this is the coefficient of $x_1^{\alpha}y_1^{\beta}$ in $\frac{F}{x_1^r}$. We have 
$$
\text{mult}(\overline F)\le\text{mult}(F).
$$

The above argument also works, by interchanging the variables $x$ and $y$,
 in the case where  $p$ is a 2 point and $\hat{\cal O}_{X_1,p_1}$ has regular parameters $(x_1,y_1)$ such that
$x=x_1y_1, y=y_1$.

Now suppose that  $p$ is a 1 point and $\hat{\cal O}_{X_1,p_1}$ has regular parameters $(x_1,y_1)$ such that
$x=x_1y_1, y=y_1$. Write 
$$
F = \sum_{i+j\ge r}a_{ij}x^iy^j.
$$
Then 
\begin{equation}
\begin{array}{l}
u= x_1^{a}y_1^{a}\\
v=\overline P(x_1y_1)+ x_1^{c}y_1^{c+r}\overline F(x_1,y_1)
\end{array}
\end{equation} 
where $\overline P=P$, $\overline F=\frac{F}{y_1^r}$. We must show that $\overline F$ has no nonzero terms $x_1^{\alpha}y_1^{\beta}$ terms
with $\alpha=r+\beta$. But this is impossible since $F$ has no nonzero $x^i$ terms, with $i\ge 0$.

The leading form of $\overline F$ is 
$$
\overline F= \sum_{i=0}^{r-1}a_{i,r-i}x_1^i+y_1\Omega
$$
since $a_{r0}=0$, where some $a_{ij}\ne 0$ with $i+j=r$, $j>0$.
Thus $\text{mult}(\overline F)\le \text{mult}(F) -1$.

Now suppose that $p$ is a 1 point and  $\hat{\cal O}_{X_1,p_1}$ has regular parameters $(x_1,y_1)$ such that
$x=x_1, y=x_1(y_1+\alpha)$.  By making if necessary a permissible change of variables at $p$,
replacing $y$ with $y-\alpha x$, we may assume that $x=x_1, y=x_1y_1$. Write 
$$
F = \sum_{i+j\ge r}a_{ij}x^iy^j.
$$
where $a_{i0}=0$ for all $i$.
$$
\begin{array}{l}
u=x_1^a\\
v=\overline P(x_1)+x_1^{c+r}\overline F(x_1,y_1)
\end{array}
$$
where $\overline P = P$, $\overline F = \frac{F}{x_1^r}$.
$\overline F$ has no nonzero terms which are powers of $x_1$.
Thus 
$$
\text{mult}(\overline F)\le\text{mult}(F).
$$
\end{pf}

Suppose that $p\in X$. Set 
\[
\sigma(p) = \left\{ \begin{array}{ll}
 0& \text{if $p$ is a 1 point and }\text{mult}(F)=\text{mult}(F(0,y))\\
\frac{1}{2} & \text{if $p$ is a 2 point}\\
1& \text{ if $p$ is a 1 point and }\text{mult}(F)<\text{mult}(F(0,y))
\end{array}
\right.
\]

\begin{Lemma}\label{L6}
$\sigma(p)$ is independent of the choice of permissible parameters $(x,y)$ at $p$.
\end{Lemma}

\begin{pf}
The proof of Lemma \ref{L3} shows that $\text{mult}(F(0,y))$ is independent of the choice of 
permissible parameters at a 1 point.
\end{pf}

\begin{Lemma}\label{L4} 
Suppose that $g:X_1\rightarrow X$ is a quadratic transform, centered at a point $p$ of $X$, and 
$p_1\in X_1$ is a  point such that $g(p_1)=p$.  Further suppose that $p$ is a 2 point, $p_1$ is 
a 1 point
and $\overline \nu(p_1)=\overline \nu(p)$. Then $\sigma(p_1) = 0$.
\end{Lemma}

\begin{pf}
$\hat{\cal O}_{X_1,p_1}$ has regular parameters $(x_1,y_1)$ such that
$x=x_1, y=x_1(y_1+\alpha)$ with $\alpha\ne 0$.
Let $r=\text{mult}(F)$.
 $\text{mult}(F_1)=\text{mult}(F)+1=r+1$. Let
$$
F = \sum_{i+j\ge r}a_{ij}x^iy^j
$$
As in the analysis leading to (\ref{eqS5}),
\begin{equation}
F_1\equiv (y_1+\alpha)^{\lambda}\left(\sum_{j=0}^ra_{r-j,j}(y_1+\alpha)^j\right)- \beta
\text{ mod }(x_1,y_1^{r+1})
\end{equation}
for some $\beta\in k$. 
If $F_1\equiv 0\text{ mod } (x_1,y_1^{r+1})$, then $\beta\ne 0$ and $-\lambda\not\in 
\{0,1,\ldots,r\}$, as  in the
proof of Theorem \ref{T1}. Then
\begin{equation}
\begin{array}{ll}
F_1&\equiv (y_1+\alpha)^{\lambda}\left(-\beta \frac{-\lambda(-\lambda-1)\cdots(-\lambda-r)}{(r+1)!}\alpha^{-\lambda-r-1}\right)y_1^{r+1}\text{ mod }(x_1,y_1^{r+2})\\
&\equiv -\beta \frac{-\lambda(-\lambda-1)\cdots(-\lambda-r)}{(r+1)!}\alpha^{-r-1}y_1^{r+1}\text{ mod }(x_1,y_1^{r+2})
\end{array}
\end{equation}
Thus $\text{mult}(F_1(0,y_1,z_1))=r+1$.
\end{pf}
 
\begin{Lemma}\label{L5} 
Suppose that $g:X_1\rightarrow X$ is a quadratic transform, centered at a 1 point $p$ of $X$  and 
$p_1$ is a  point above $p$ such that $g(p_1)=p$. 

If $p_1$ is a 1 point and $\overline \nu(p_1)=\overline \nu(p)$, then $\sigma(p_1) = 0$.
If $\sigma(p)=0$ and $p_1$ is a 2 point then $\overline \nu(p_1) = 0$. \end{Lemma}

\begin{pf}
First suppose that $\sigma(p)=0$ and $p_1$ is a 2 point. Then $\hat{\cal O}_{X_1,p_1}$ has regular parameters $(x_1,y_1)$ such that
$x=x_1y_1, y=y_1$.
$$
F = \sum_{i+j\ge r}a_{ij}x^iy^j
$$
with $a_{0r}\ne 0$.
$$
F_1= \sum_{i=0}^{r-1}a_{i,r-i}x_1^i+y_1\Omega.
$$
is then a unit.

Now suppose that $p_1$ is a 1 point and $\nu(p_1)=\nu(p)$. 
After appropriate choice of permissible variables $(x,y)$ at $p$, 
$\hat{\cal O}_{X_1,p_1}$ has regular parameters $(x_1,y_1)$ such that
$x=x_1, y=x_1y_1$.
Set $r=\text{mult}(F)=\text{mult}(F_1)$.
Then $F_1 = \frac{F(x_1,x_1y_1)}{x_1^r}$ and
$\text{mult}(F_1(0,y_1)) =r$.
\end{pf}

\begin{Theorem}\label{T2} 
Suppose that $g:X_1\rightarrow X$ is a quadratic transform, centered at a point $p$ of $X$, and 
$p_1$ is a closed point such that $g(p_1)=p$.  If  $\overline\nu(p_1)=\overline\nu(p)$, then $\sigma(p_1)\le \sigma(p)$.
\end{Theorem}

\begin{pf}
This is immediate from Lemmas \ref{L4} and \ref{L5}.
\end{pf}

Suppose that 
$$
F = \sum_{i+j\ge r}a_{ij}x^iy^j
$$
has multiplicity $r$.
Define
$$
\delta(F;x,y) = \text{min}(\frac{i}{r-j}\mid j<r, a_{ij}\ne 0).
$$
$\delta(F;x,y)=\infty$ if and only if $F = y^r\omega$, where $\omega$ is a unit.
If $\delta(F;x,y)<\infty$, then $\delta(F;x,y)\in \frac{1}{r!}{\bold N}$.

Suppose that $p\in X$. If $(x,y)$ are permissible parameters at $p$ with one of the forms (\ref{I}) or (\ref{II}), set
$$
\delta(p;x,y) = \delta(F;x,y).
$$ 
Then set 
\[
\delta(p) = 
\text{sup}(\delta(p;x,y)) 
\]
where the sup is over all permissible parameters at $p$.
Note that if $p$ is a 2 point, then
$$
\delta(p)=\text{max}(\delta(p;x,y),\delta(p;y,x)) 
$$ 
if  $(x,y)$ are a particular choice of permissible parameters at $p$.

If $p$ is a 2 point and $\nu(p)>0$, then $\delta(p)<\infty$.  If $p$ is a 1 point and $\sigma(p)=1$, then $\delta(p)=1$, since 
$\delta(p;x,y)=1$ 
for all permissible parameters $(x,y)$.

\begin{Lemma} \label{L8} 
Suppose that $p$ is a 1 point, $\sigma(p) = 0$ and $(x,y)$ are fixed permissible parameters at $p$. Then
there exists a  power series $t(x)$ such that 
$$
\delta(p) = \delta(p;x, y-t(x)).
$$
If $\delta(p)<\infty$, then $t(x)$ is a polynomial.

$\delta(p) > \delta = \delta(p;x,y)$
if and only if  $\delta\in \bold N$ and
$$
\sum_{i+\delta j=r\delta}a_{ij}x^iy^j= \tau(y-c x^{\delta})^r+\lambda x^{r\delta}
$$
for some $\tau,c,\lambda\in k$ with $c\ne 0$ (so that $\lambda=-\tau(-c)^r$).
\end{Lemma}

\begin{pf}
Suppose that $(\overline x, \overline y)$ are also permisible parameters at $p$. Then
$\overline x = \lambda x$, with $\lambda^{a}=1$ and $\overline y = \overline \phi (y-t(x))$
for some unit series $\overline \phi$ and series $t(x)$. 
$$
\delta(p;\overline x,\overline y) = \delta(p;x,y-t(x))
$$
Thus 
\begin{equation}\label{eqS7}
\delta(p) = \text{sup}(\delta(p;x,y-t(x)\mid t(x) \text{ is a polynomial of positive order}). 
\end{equation}

Let $\delta=\delta(p;x,y)$,
$$
\overline L = \sum_{i+\delta j=r\delta}a_{ij}x^iy^j.
$$
so that 
$$
F = \overline L+ \sum_{i+\delta j>r\delta}a_{ij}x^iy^j.
$$
Suppose that 
$$
\overline L = \tau(y-c x^{\delta})^r+\lambda x^{r\delta}
$$
for some $\tau,c,\lambda\in k$ with $0\ne c$.
Set $y_1 = y-cx^{\delta}$. Then $\delta \in \bold N$ and $\delta(p;x,y_1)> \delta(p;x,y)$ since
$$
F_1 = \tau y_1^r + \sum_{i+\delta j>r\delta}\overline a_{ij}x^iy_1^j.
$$
where 
$$
v=P_1(x)+x_1^cF_1
$$
is the normalized form of $v$ with respect to $(x,y_1)$.
We can repeat this process, with $y$ replaced by $y-cx^{\delta}$. The process will either
 produce a polynomial $t(x)$
such that if $y_1 = y-t(x)$, and $\delta_1 = \delta(p;x,y_1)$, then $\delta_1\not \in \bold N$, or $\delta_1\in \bold N$ and
\begin{equation}\label{eqS8}
\sum_{i+\delta_1 j=r\delta_1}\overline a_{ij}x^iy_1^j\neq \tau(y_1-c x^{\delta_1})^r+\lambda x^{r\delta_1}
\end{equation}
for any $\tau,c,\lambda\in k$ with $0\ne c$, or we will produce a series $t(x)$ such that if $y_1 = y-t(x)$, then $\delta(p;x,y_1)=\delta(p)=\infty$, so that
$F_1 = y_1^r\phi$, where $\phi$ is a unit series. 

Suppose that we have produced $y_1$ such that 
$\delta(p;x,y_1)\not\in{\bold N}$ or $\delta(p;x,y_1)\in{\bold N}$
and (\ref{eqS8}) holds. We will show that $\delta(p)=\delta(p;x,y_1)$.
Suppose that $\delta_1 = \delta(p;x,y_1)<\delta(p)$. By (\ref{eqS7}), there is a polynomial 
$$
t(x)=\sum e_ix^i
$$
 such that if $y_2 = y_1-t(x)$, then $\delta(p;x,y_2)>\delta(p;x,y_1)$. Substitute $y_1 = y_2+t(x)$ into
$$
F_1 =\sum_{i+\delta_1 j=r\delta_1}\overline a_{ij}x^iy_1^j + \sum_{i+\delta_1 j>r\delta_1}\overline a_{ij}x^iy_1^j,
$$
and normalize with respect to the permissible parameters to get
$$
v=P_2(x)+x^cF_2(x,y_2).
$$
Let $\overline d=\text{ord }(t(x))$. $x^iy_1^j=x^i(y_2+t(x))^j$ has nonzero 
$x^{i+m\overline d}y_2^{j-m}$ terms with $0\le m\le j$, and may have other nonzero
$x^{i+m\overline d+\gamma}y_2^{j-m}$ terms with $0\le m\le j$, $\gamma\ge 0$.

Suppose that $\overline d<\delta_1=\delta(p;x,y_1)$. The expansion of $y_1^r$ has a 
nontrivial $x^{\overline d}y_1^{r-1}$ term. Suppose that $x^iy_1^j$ is such that
its expansion has a nontrivial $x^{\overline d}y_1^{r-1}$ term. Then
$\overline d=i+m\overline d+\gamma$, $r-1=j-m$ with $0\le m\le j$, $i,\gamma\ge 0$.
$\overline d(1-m)=i+\gamma\ge 0$ implies $m=0$ or 1. $m=0$ implies $j=r-1$, $i\le\overline d$.
$\overline a_{ij}=0$ in this case since 
$$
i+\delta_1j\le \overline d+\delta_1(r-1)<\delta_1r.
$$
$m=1$ implies $i=0$, $j=r$. Thus there exists a nontrivial $x^{\overline d}y_1^{r-1}$ term
in $F_2(x,y_1)$ which implies that $\delta_2<\delta_1$, a contradiction. Thus $\overline d\ge
\delta_1$.

We then see that if $i+\delta_1j>r\delta_1$, then all terms $x^{\alpha}y^{\beta}$
in the expansion of $x^iy_1^j=x^i(y_2+t(x))^j$ satisfy $\alpha+\delta_1\beta>r\delta_1$.
Since $\delta_2<\delta_1$,
we see that 
$$
\sum_{i+\delta_1 j=r\delta_1}\overline a_{ij}x^i(y_2+t(x))^j =
\left\{\begin{array}{l}
cy_2^r+
\text{terms with }i+\delta_1 j>r\delta_1\text{ if }\delta_1\not\in\bold N,\\
cy_2^r+dx^{r\delta_1} +
\text{terms with }i+\delta_1 j>r\delta_1\text{ if } \delta_1\in \bold N.
\end{array}
\right.
$$
Thus $\text{mult}(t)=\delta_1$ and
$$
\sum_{i+\delta_1 j=r\delta_1}\overline a_{ij}x^iy_1^j = c(y_1-e_{\delta_1}x^{\delta_1})^r+dx^{r\delta_1}
$$
a contradiction.
\end{pf}

\begin{Lemma}\label{L7} 
Suppose that $g:X_1\rightarrow X$ is a quadratic transform, centered at a point $p$ of $X$, and 
$p_1\in X_1$ is a closed point above $p$ such that $g(p_1)=p$ and $\overline\nu(p_1)=\overline\nu(p)$. 

Suppose that $p$ and $p_1$ are both 2 points. Then $\delta(p_1)=\delta(p)-1$.

Suppose that $p$ and $p_1$ are both 1 points, $\sigma(p)=0$ and $\delta(p)<\infty$. Then $\delta(p_1)=\delta(p)-1$.
\end{Lemma}

\begin{pf} Suppose that $r=\text{mult}(F)$.

First suppose that $p$ and $p_1$ are both 2 points. Then $p$ has permissible  parameters $(x,y)$ and
$\hat{\cal O}_{X_1,p_1}$ has permissible parameters $(x_1,y_1)$ such that
$x=x_1, y=x_1y_1$.
Since $F_1=\frac{F}{x_1^r}$, $\delta(p_1;x_1,y_1) = \delta(p;x,y)-1$.
Since $\overline\nu(p_1)=\overline\nu(p)$, we have $F=\sum a_{ij}x^iy^j$ with $a_{ij}=0$ if $i+j\le r$ and $j<r$.
Thus $a_{0r}\ne 0$, so that $\delta(p;y,x)=1$ and $\delta(p;x,y)>1$. Thus
$\delta(p)=\delta(p;x,y)$. Since $\text{mult}(F_1)=r$ and $\text{mult}(F_1(0,y_1))=r$,
$\delta(p_1;y_1,x_1)=1$ and $\delta(p_1;x_1,y_1)\ge 1$. Then $\delta(p_1)=
\delta(p_1;x_1,y_1)=\delta(p)-1$.

Now suppose that $p$ and $p_1$ are both 1 points, $\sigma(p)=0$  and $\delta(p)<\infty$.
We can suppose that we have permissible coordinates $(x,y)$ at $p$ such that
$\delta=\delta(p)=\delta(F;x,y)$ and
$\text{mult}(F(0,y))=\text{mult}(F)$.
$p_1$ has permissible parameters $(x_1,y_1)$ such that
$x=x_1,y=x_1(y_1+\gamma)$ for some $\gamma\in k$. 

First suppose that $\gamma\ne 0$. 
$$
F_1 = \sum_{i+j=r}a_{ij}(y_1+\gamma)^j-\overline a +x_1\Omega
$$
where 
$$
\overline a = \sum_{i+j=r}a_{ij}\gamma^j.
$$
$\text{mult}(F_1)=\text{mult}(F)$ implies
$$
\sum_{i+j=r}a_{ij}(y_1+\gamma)^j-\overline a=a_{0r} y_1^r.
$$
Thus
$$
 \sum_{i+j= r}a_{ij}x^iy^j = a_{0r}(y-\gamma x)^r+\overline a x^r.
$$
 This is a contradiction to the assumption that $\delta(p;x,y)=\delta(p)$ 
by Lemma \ref{L8}. 

Now suppose that $\gamma=0$. Then $F_1= \frac{F}{x_1^r}$ and  $\delta(p_1;x_1,y_1)=\delta(p;x,y)-1$.
If $\delta(p_1;x_1,y_1)<\delta(p_1)$, then we must also have $\delta(p;x,y)<\delta(p)$ By Lemma \ref{L8}.
Thus $\delta(p_1)=\delta(p)-1$. 
\end{pf}

If $p\in X$ is a 2 point, then $(u,v)$ are 1-resolved at $p$  precisely when $\overline \nu(p)=0$. 
If $p\in X$ is a 1 point then $(u,v)$ are 1-resolved at $p$ precisely when $\delta(p)=\infty$. Thus $(u,v)$ are not 1-resolved at $p\in X$ if and only if $\overline\nu(p)>0$ and
$\delta(p)<\infty$.

We can define an invariant 
$$
\text{Inv}(p) = (\overline\nu(p),\sigma(p),\delta(p))
$$
for $p\in X$.

\begin{Theorem}\label{T3} 
Suppose that $g:X_1\rightarrow X$ is a quadratic transform, centered at a point $p$ of $X$, and 
$p_1\in X_1$ is  such that $g(p_1)=p$. Suppose that $\overline\nu(p)>0$ and $\delta(p)<\infty$. Then
$$
\text{Inv}(p_1) < \text{Inv}(p)
$$
in the lexicographic ordering.
\end{Theorem} 

\begin{pf} The Theorem follows from Theorem \ref{T1}, Lemmas \ref{L4}, \ref{L5}, \ref{L7}.
\end{pf}

The proof of Theorem \ref{Theorem965} is immediate from Theorem \ref{T3}.

 \begin{Lemma}\label{Lemma1024} Suppose that $f(x,y)\in T_0=k[[x,y]]$
is a series. Suppose that we have an infinite sequence of quadratic transforms
$$
T_0\rightarrow T_1\rightarrow\cdots\rightarrow T_n\rightarrow\cdots
$$
Then there exists $n_0$ such that $n\ge n_0$ implies there exist regular parameters $(x_n,y_n)$ in $T_n$, $\alpha_n,\beta_n\in
\bold N$ and a unit  $u_n\in T_n$ such that $f=x_n^{\alpha_n}y_n^{\beta_n}u_n$.
\end{Lemma}

\begin{pf} This follows directly from Zariski's proof of resolution of surface singularities
along a valuation (\cite{Z1}), or can be deduced easily after blowing up enough to make $f=0$ a
SNC divisor.
\end{pf}

\begin{Lemma}\label{Lemma23}
Suppose that $\alpha_j+\beta_j\ge j$, $x^{\alpha_j}y^{\beta_j}\in T_0=k[x,y]_{(x,y)}$  for 
$2\le j\le r$ (or $1\le j\le r$).
Suppose that we have a sequence of quadratic transforms
$$
T_0\rightarrow T_1\rightarrow \cdots \rightarrow T_n\rightarrow\cdots
$$
where each  $T_n$ has regular parameters $(x_n,y_n)$ such that
either $x_{n-1} = x_n$, $y_{n-1}=x_ny_n$, or $x_{n-1}=x_ny_n$, $y_{n-1}=y_n$.
There are natural numbers $\alpha_{n,i}$, $\beta_{n,i}$ such that
$$
 x^{\alpha_j}y^{\beta_j}= x_n^{\alpha_{n,i}}y_n^{\beta_{n,i}}.
$$
Define 
$$
\delta_{n,i,j} = \left(\frac{\alpha_{n,i}}{i}-\frac{\alpha_{n,j}}{j}\right)
\left(\frac{\beta_{n,i}}{i}-\frac{\beta_{n,j}}{j}\right)
$$
Then
\begin{enumerate}
\item $\delta_{n+1,i,j}\ge \delta_{n,i,j}$
\item  $\delta_{n,i,j}<0$ implies $\delta_{n+1,i,j}-\delta_{n,i,j}\ge \frac{1}{r^4}$.
\end{enumerate}
\end{Lemma}

\begin{pf}
We will first verify 1. Suppose that  $x_{n} = x_{n+1}y_{n+1}$, $y_{n}=y_{n+1}$. 
The proof when $x_{n}=x_{n+1}$, $y_{n}=x_{n+1}y_{n+1}$ is the same. 1. is immediate from
$$
\delta_{n+1,i,j} = \delta_{n,i,j} + \left(\frac{\alpha_{n,i}}{i}-\frac{\alpha_{n,j}}{j}\right)^2
$$
Now suppose that $\delta_{n,i,j}<0$. 
Then $\left(\frac{\alpha_{n,i}}{i}-\frac{\alpha_{n,j}}{j}\right)$
and 
$\left(\frac{\beta_{n,i}}{i}-\frac{\beta_{n,j}}{j}\right)$ are nonzero.
We can suppose that $x_{n} = x_{n+1}y_{n+1}$, $y_{n}=y_{n+1}$.  
$$
\begin{array}{ll}
\delta_{n+1,i,j} - \delta_{n,i,j} &= \left(\frac{\alpha_{n,i}}{i}-\frac{\alpha_{n,j}}{j}\right)^2\\
&=\left(\frac{j\alpha_{n,i}-i\alpha_{n,j}}{ij}\right)^2\ge \frac{1}{r^4}
\end{array}
$$
since $i,j \le r$ implies $(ij)^2\le r^4$.
\end{pf}

\begin{Corollary}\label{Corollary24}
Suppose that $\alpha_j+\beta_j\ge j$ and $x^{\alpha_j}y^{\beta_j}\in T_0=k[x,y]_{(x,y)}$  for $2\le j\le r$ (or $1\le j\le r$)
and
$$
T_0\rightarrow T_1\rightarrow \cdots \rightarrow T_n\rightarrow\cdots
$$
is a sequence of quadratic transformations as in the statement of Lemma \ref{Lemma23}.
Then
\begin{enumerate}
\item There exists $n_0$ and $i$ such that $n\ge n_0$ implies 
$$
\frac{\alpha_{n,i}}{i}\le \frac{\alpha_{n,j}}{j} \text{ and }\frac{\beta_{n,i}}{i}\le \frac{\beta_{n,j}}{j}
$$
for $2\le j\le r$ (or $1\le j\le r$).
\item There exists an $n_1\ge n_0$ such that 
$$
\left\{\frac{\alpha_{n_1,i}}{i}\right\} +  \left\{\frac{\beta_{n_1,i}}{i}\right\}<1
$$
\end{enumerate}
\end{Corollary}

\begin{pf} By Lemma \ref{Lemma23}, there exists $n_0$ such that $n\ge n_0$ implies
$\delta_{n,i,j}\ge 0$ for all $i,j$. 
Let $\lambda_1 = \text{min}\left(\frac{\alpha_{n,j}}{j}\right)$.
Let $\lambda_2 = \text{min}\left(\frac{\beta_{n,j}}{j}\text{ such that } \frac{\alpha_{n,j}}{j}=\lambda_1\right)$.
Choose $i$ such that $\frac{\alpha_{n,i}}{i}= \lambda_1$, $\frac{\beta_{n,i}}{i}=\lambda_2$.
Then 
$$\frac{\alpha_{n,i}}{i}\le \frac{\alpha_{n,j}}{j},\,\,\frac{\beta_{n,i}}{i}\le \frac{\beta_{n,j}}{j}
$$
for $2\le j\le r$ (or $1\le j\le r$).

Now we will prove 2. Suppose that $n\ge n_0$. 
Then 
$$
\left\{\frac{\alpha_{n,i}}{i}\right\} +  \left\{\frac{\beta_{n,i}}{i}\right\}\in \frac{1}{i}\bold N.
$$
Suppose that
$$
\left\{\frac{\alpha_{n,i}}{i}\right\} +  \left\{\frac{\beta_{n,i}}{i}\right\}\ge 1.
$$
Without loss of generality,
$$
x_{n}=x_{n+1}y_{n+1}, y_{n}=y_{n+1}.
$$
Then
$$
\begin{array}{ll}
\left\{\frac{\alpha_{n+1,i}}{i}\right\} +  \left\{\frac{\beta_{n+1,i}}{i}\right\}
 &=
\left\{\frac{\alpha_{n,i}}{i}\right\} +\left\{\frac{\alpha_{n,i}}{i}\right\} +  \left\{\frac{\beta_{n,i}}{i}\right\}- 1
\\
&< \left\{\frac{\alpha_{n,i}}{i}\right\} +  \left\{\frac{\beta_{n,i}}{i}\right\}
\end{array}
$$
Thus there exists $n_1\ge n_0$ such that 2. holds.
\end{pf}

\begin{Remark}\label{Remark1049} The conditions $\alpha_i+\beta_i\ge i$ and
$\left\{\frac{\alpha_i}{i}\right\}+\left\{\frac{\beta_i}{i}\right\}<1$ imply either
$\alpha_i\ge i$ or $\beta_i\ge i$.
\end{Remark}

Lemmas \ref{Lemma966} and \ref{Lemma967} are used in Abhyankar's Good Point proof of
resolution of singularities \cite{Ab6}, \cite{Li}.

\begin{Lemma}\label{Lemma966}
Suppose that $\alpha_{j_0}+\beta_{j_0}\ge j$, 
$(\alpha_{j_0},\beta_{j_0})$ are nonnegative integers  for $1\le j\le r$.
Suppose that we have pairs of nonnegative integers $(\alpha_{n,j},\beta_{n,j})$
for all positive $n$ and $1\le j\le r$ such that either
$$
(\alpha_{n+1,j},\beta_{n+1,j})=(\alpha_{n,j}+\beta_{n,j}-j,\beta_{n,j})
$$
or
$$
(\alpha_{n+1,j},\beta_{n+1,j})=(\alpha_{n,j},\alpha_{n,j}+\beta_{n,j}-j).
$$
Define 
$$
\delta_{n,i,j} = \left(\frac{\alpha_{n,i}}{i}-\frac{\alpha_{n,j}}{j}\right)
\left(\frac{\beta_{n,i}}{i}-\frac{\beta_{n,j}}{j}\right)
$$
Then
\begin{enumerate}
\item $\delta_{n+1,i,j}\ge \delta_{n,i,j}$
\item  $\delta_{n,i,j}<0$ implies $\delta_{n+1,i,j}-\delta_{n,i,j}\ge \frac{1}{r^4}$.
\end{enumerate}
\end{Lemma}

\begin{Lemma}\label{Lemma967}
Suppose that 
the assumptions are as in Lemma \ref{Lemma966}. 

Suppose that $\alpha_j+\beta_j\ge j$, $x^{\alpha_j}y^{\beta_j}\in T_0=k[x,y]_{(x,y)}$
for $1\le j\le r$. Suppose that we have a possibly infinite sequence of quadratic
transforms
$$
T_0\rightarrow T_1\rightarrow\cdots\rightarrow T_n\rightarrow \cdots
$$
where each $T_n$ has regular parameters $(x_n,y_n)$ such that either
$x_{n-1}=x_n,y_{n-1}=x_ny_n$ or
$x_{n-1}=x_ny_n, y_{n-1}=y_n$ and $(\alpha_n,\beta_n)$ are defined by the
respective rules of Lemma \ref{Lemma966}. Then
\begin{enumerate} 
\item There exists $n_0$ and $i$ such that $n\ge n_0$ implies
$$
\frac{\alpha_{n,i}}{i}\le\frac{\alpha_{n,j}}{j}\text{ and }
\frac{\beta_{n,i}}{i}\le\frac{\beta_{n,j}}{j}
$$
for $1\le j\le r$.
\item
There exists $n_1\ge n_0$
such that 
$$
\left\{\frac{\alpha_{n_1,i}}{i}\right\} +  \left\{\frac{\beta_{n_1,i}}{i}\right\}<1
$$
\end{enumerate}
\end{Lemma}

\section{$\bf{A_r(X)}$}

Throughout this section we will assume that $\Phi_X:X\rightarrow S$ is weakly prepared.

\begin{Definition}\label{Def1093} Suppose that $r\ge 2$.
 $\overline A_r(X)$ holds  if 
\begin{enumerate}
\item $\nu(p)\le r$ if $p\in X$ is a 1 point or a 2 point.
\item If $p\in X$ is a 1 point and $\nu(p)= r$, then $\gamma(p)=r$.
\item If $p\in X$ is a 2 point and $\nu(p)= r$, then $\tau(p)>0$.
\item $\nu(p)\le r-1$ if $p\in X$ is a 3 point
\end{enumerate}
\end{Definition}

\begin{Definition}\label{Def1094} Suppose that $r\ge 2$.
 $A_r(X)$ holds if
\begin{enumerate}
\item $\overline A_r(X)$ holds.
\item $\overline S_r(X)$ is a union of nonsingular curves and isolated points.
\item $\overline S_r(X)\cap (X-\overline B_2(X))$ is smooth.
\item $\overline S_r(X)$ makes SNCs with $\overline B_2(X)$ on the open set
$X-B_3(X)$.
\item The curves in $\overline S_r(X)$ passing through a 3 point $q\in X$
 have distinct tangent directions at $q$. (They are however, allowed to be tangent to
a 2 curve).
\end{enumerate}
\end{Definition}

\begin{Definition}\label{Def2005} Suppose that $A_r(X)$ holds. A weakly permissible monoidal transform
$\pi:X_1\rightarrow X$ is called permissible if $\pi$ is the blowup of a point,
a 2 curve or a curve $C$ containing a 1 point such that $C\cup\overline S_r(X)$
makes SNCs with $\overline B_2(X)$ at all points of $C$.
\end{Definition}

\begin{Remark}

\begin{enumerate}
\item If $A_r(X)$ holds and $\pi:X_1\rightarrow X$ is a permissible monodial transform, then the
strict transform of $\overline S_r(X)$ on $X_1$ makes SNCs with $\overline B_2(X_1)$ at
1 and 2 points, and has distinct tangent directions at 3 points.
\item If $\pi:X_1\rightarrow X$ is a quadratic transform centered at a point $p\in X$
with $\nu(p)=r$ and $A_r(X)$ holds, then $A_r(X_1)$ holds.
\item If $A_r(X)$ holds and all 3 points $q$ of $X$ satisfy $\nu(q)\le r-2$, then 
$\overline S_r(X)$ makes SNCs with $\overline B_2(X)$.
\end{enumerate}
\end{Remark}

The Remark follows from  Lemmas \ref{Lemma54} and \ref{Lemma4},
and the observation that the strict transforms of nonsingular curves with distinct tangent
directions at a point $p$ intersect the exceptional fiber of the blowup of $p$
transversally in distinct points.

\section{Reduction of $\nu$ in a special case}\label{Spec1}

Throughout this section we will assume that $\Phi_X:X\rightarrow S$ is weakly prepared.

\begin{Lemma}\label{Lemma18}
Suppose that $r\ge 2$ and $A_r(X)$ holds, $p\in X$ is a 1 point or a 2 point with $\nu(p)=\gamma(p)=r$.
Let $R ={\cal O}_{X,p}$. Suppose that $(x,y,z)$ are permissible parameters at $p$ as in Lemma \ref{Lemma17}.
Then
there exists a finite sequence of permissible monodial transforms $\pi:Y\rightarrow \text{Spec}(\hat{R})$
centered at sections over $C = V(x,y)$, such that for $q\in \pi^{-1}(p)$,
there exist permissible parameters $(\overline x,\overline y,z)$ at $q$
such that  $F_q$ has one of the following forms.
\begin{equation}\label{eq42}
\begin{array}{ll}
u&=\overline x^a\\
v&=P(\overline x)+\overline x^cF_q
\end{array}
\end{equation}
or
\begin{equation}\label{eq43}
\begin{array}{ll}
u&=(\overline x^a\overline y^b)^m\\
v&=P(\overline x^a\overline y^b)+\overline x^c\overline y^dF_q
\end{array}
\end{equation}
with
$$
F_q = 
\tau z^r+\sum_{i=2}^{r-1}\overline a_i(\overline x,\overline y)\overline x^{\alpha_i}\overline y^{\beta_i}z^{r-i}+\epsilon \overline x^{\alpha_r}
\overline y^{\beta_r}
$$
with $\tau$ a unit, $\overline a_i$ a unit (or zero), 
$\alpha_i+\beta_i\ge i$ for all $i$, and $\epsilon = 0$ or 1.
\end{Lemma}

\begin{pf} We have one of the forms (\ref{eq37}) or (\ref{eq38}) of Lemma \ref{Lemma17}
at $p$.
By Lemma \ref{Lemma1027} and Theorem \ref{Theorem965} applied to
$$
\overline u=x^a, \overline v=P(x)+x^cF_p(x,y,0)
$$
or
$$
\overline u=(x^ay^b)^m, \overline v = P(x^ay^b)+x^cy^dF_p(x,y,0)
$$
 there exists a sequence of permissible blowups of sections over $C$ such that 
for all $q$ over $p$, there are permissible parameters $(\overline x,\overline y,z)$ at $q$ such that
\begin{equation}\label{eq44}
\begin{array}{ll}
u&=\overline x^a\\
v&=P(\overline x)+\overline x^c(\tau z^r+\sum_{i=2}^{r-1} a_i(\overline x,
\overline y)z^{r-i}+\epsilon \overline x^{e_0}\overline y^{f_0})
\end{array}
\end{equation}
with $\tau$ a unit, $\epsilon = 0$ or 1 and $f_0>0$, or 
\begin{equation}\label{eq45}
\begin{array}{ll}
u&=(\overline x^a\overline y^b)^m\\
v&=P(\overline x^a\overline y^b)+\overline x^c\overline y^d(\tau z^r+\sum_{i=2}^{r-1} a_i(\overline x,\overline y)z^{r-i}+\epsilon \overline x^{e_0}\overline y^{f_0})
\end{array}
\end{equation}
with $\tau$ a unit, $\epsilon = 0$ or 1 and $a(d+f_0)-b(c+e_0)\ne 0$.

By further permissible blowing up (of sections over $C$) we can make 
$$
u\prod_{2\le i\le r-1, a_i\ne0}a_i=0
$$
 a SNC divisor,
while preserving the forms (\ref{eq44}) and (\ref{eq45}). At points $q$ over $p$
satisfying (\ref{eq45}) we have then achieved the
conclusions of the Lemma. 

Suppose that $q$ is a point over $p$ satisfying (\ref{eq44}) such that the conclusions of the Lemma do not hold.
We then have $\epsilon=1$, $f_0>0$ and  
$$
u\prod_{2\le i\le r, a_i\ne0}a_i=0
$$
 is not a SNC divisor.   Since 
$$
u\prod_{2\le i\le r-1, a_i\ne0}a_i=0
$$
is a SNC divisor, there exists a nonzero, nonunit series $g(\overline x)$ such that 
\begin{equation}\label{eq500}
a_i = \overline a_i(\overline x,\overline y) \overline x^{\alpha_i}(\overline y
-g(\overline x))^{\beta_i}
\end{equation}
for $2\le i\le r-1$, where the $\overline a_i$ are units (or 0), and some $\beta_i>0$
with $\overline a_i\ne 0$.
If $f_0=1$, we can set $\tilde y = \overline y-g(\overline x)$ and renormalize 
with respect to $(\overline x,\tilde y,z)$ to get in the form of the
conclusions of the Lemma. 

Otherwise $f_0>1$. Let $t=\nu(g(\overline x))$, so that 
$$
g(\overline x) = \alpha 
\overline x^t+\text{ higher order terms }
$$
 for some $0\ne \alpha$. Now blow up 
$V(\overline x,\overline y)$.
under $\overline x=x_1y_1$, $\overline y=y_1$, we have
$$
\begin{array}{ll}
u&=x_1^ay_1^a\\
v&=P(x_1y_1)+x_1^cy_1^c(\tau z^r+\sum_{i=2}^{r-1}
\tilde a_i(x_1,y_1)x_1^{\alpha_i}y_1^{\beta_i+\alpha_i}z^{r-i}+x_1^{e_0}y_1^{e_0+f_0})
\end{array}
$$
in the form of the conclusions of the Lemma. Under $\overline x=x_1$, $\overline y=x_1(y_1+\beta)$, with $\beta\ne 0$, we have
$$
\begin{array}{ll}
u&=x_1^a\\
v&=P(x_1)+x_1^c(\tau z^r+\sum_{i=2}^{r-1}\overline a_i x_1^{\alpha_i+\beta_i}
(y_1+\beta-\frac{g(x_1)}{x_1})^{\beta_i}z^{r-i}
+ x_1^{e_0+f_0}(y_1+\beta)^{f_0})\\
&=P(x_1)+\beta^{f_0}x_1^{c+e_0+f_0}
+x_1^c(\tau z^r+\sum_{i=2}^{r-1}\overline a_i x_1^{\alpha_i+\beta_i}
((\overline y_1+\beta^{f_0})^{\frac{1}{f_0}}-\frac{g(x_1)}{x_1})^{\beta_i}z^{r-i}
+ x_1^{e_0+f_0}\overline y_1)
\end{array}
$$
where $\overline y_1 = (y_1+\beta)^{f_0}-\beta^{f_0}$. If we are not in the
form of the conclusions of the Lemma,
then
$$
(\overline y_1+\beta^{f_0})^{\frac{1}{f_0}}-\frac{g(x_1)}{x_1}
=a(x_1,\overline y_1)(\overline y_1-\phi(x_1))
$$
where $\nu(\phi)\ge 1$. We can
make a change of variable in $\overline y_1$,
replacing $\overline y_1$ with $\overline y_1-\phi(x_1)$,
 and renormalize, to get in the form of the conclusions of the Lemma.

Under $\overline x=x_1$, $\overline y=x_1y_1$, we have
$$
\begin{array}{ll}
u&=x_1^a\\
v&= P(x_1)+x_1^c(\tau z^r + \sum_{i=2}^{r-1}\overline a_i(x,y) x_1^{\alpha_i+\beta_i}
(y_1-\frac{g(x_1)}{x_1})^{\beta_i}z^{r-i}
+x_1^{e_0+f_0}y_1^{f_0})
\end{array}
$$
the coefficients of $z^i$ are in the form of (\ref{eq500}), but we have a reduction 
$\nu\left(\frac{g(x_1)}{x_1}\right)=t-1$. If $\nu\left(\frac{g(x_1)}{x_1}\right)=0$ we are in 
the form of the conclusions of the Lemma.
Thus after $t$ blowups, centered at the intersection of the strict transform of the surface $y=0$
with the exceptional divisor,  we achieve the conclusions of the Lemma.
\end{pf}

\begin{Theorem}\label{Theorem19}
Suppose that $r\ge 2$ and $A_r(X)$ holds, $p\in X$ is a 1 point or a 2 point with $\nu(p)=\gamma(p)=r$.
Let $R ={\cal O}_{X,p}$. Suppose that $(x,y,z)$ are permissible parameters at $p$
 as in Lemma \ref{Lemma17}, where $z=\sigma\tilde z$ for some $\tilde z\in R$
and unit $\sigma\in \hat R$.
Then
there exists a finite sequence of permissible monodial transforms $\pi:Y\rightarrow \text{Spec}(\hat{R})$
centered at sections over $C = V(x,y)$, such that for $q\in \pi^{-1}(p)$,
$q$ has permissible parameters $(\overline x,\overline y, z)$ such that 
 $F_q$ has one of the following forms:
 \begin{enumerate}
\item
\begin{equation}\label{eq39}
\begin{array}{ll}
u&=\overline x^a\\
v&=P(\overline x)+\overline x^cF_q\text{ with }\\
F_q &= 
\tau z^r+\sum_{i=2}^{r-1}\overline a_i(\overline x,\overline y)
\overline x^{\alpha_i}z^{r-i}+\epsilon \overline x^{\alpha_r}\overline y
\end{array}
\end{equation}
where $\tau$ is a unit, $\alpha_i\ge i$ for $2\le i\le r-1$, $\alpha_r\ge r-1$, $\epsilon =0$ or $1$, and $\overline a_i$ are units (or 0),
or
\item
\begin{equation}\label{eq40}
\begin{array}{ll}
u&=(\overline x^a\overline y^b)^m\\
v&=P(\overline x^a\overline y^b)+\overline x^c\overline y^dF_q\text{ with }\\
F_q &= 
\tau z^r+\sum_{j=2}^{r-1}\overline a_j(\overline x,\overline y)\overline x^{\alpha_j}\overline y^{\beta_j}z^{r-j}+\epsilon \overline x^{\alpha_r}\overline y^{\beta_r}
\end{array}
\end{equation}
where $\tau$ is a unit,
$\alpha_j+\beta_j\ge j$  and $\overline a_j$ are units or 0 for all $j$, $\epsilon=0$ or 1,
there exists an $i$ such that $\overline a_i\ne 0$, $2\le i\le r$ and 
$$
\frac{\alpha_i}{i}\le \frac{\alpha_j}{j},\,\,\frac{\beta_i}{i}\le \frac{\beta_j}{j}
$$
for $2\le j\le r$. We further have
$$
\left\{\frac{\alpha_i}{i}\right\} +  \left\{\frac{\beta_i}{i}\right\}<1
$$
or
\item
\begin{equation}\label{eq41}
\begin{array}{ll}
u&=\overline x^a\\
v&=P(\overline x)+\overline x^bF_q\text{ with }\\
F_q &= 
\tau z^r+\sum_{j=2}^{r-1}\overline a_j(\overline x,\overline y)\overline x^{\alpha_j}\overline y^{\beta_j}z^{r-j}+\epsilon \overline x^{\alpha_r}\overline y^{\beta_r}
\end{array}
\end{equation}
where $\tau$ is a unit,
$\alpha_j+\beta_j\ge j$ and $\overline a_j$ are units or 0 for all $j$, $\epsilon=0$ or 1,
there exists an $i$ such that $\overline a_i\ne 0$, $2\le i\le r$ and 
$$
\frac{\alpha_i}{i}\le \frac{\alpha_j}{j},\,\,\frac{\beta_i}{i}\le \frac{\beta_j}{j}
$$
for $2\le j\le r$. We further have
$$
\left\{\frac{\alpha_i}{i}\right\} +  \left\{\frac{\beta_i}{i}\right\}<1.
$$
 \end{enumerate}
\end{Theorem}

\begin{pf}
We can first construct a sequence of monoidal transforms $\pi:Y\rightarrow\text{Spec}(\hat R)$
satisfying the conclusions of Lemma \ref{Lemma18}.
(\ref{eq39}) holds at all but finitely many points $q\in\pi^{-1}(p)$. 

Suppose that $q\in \pi^{-1}(p)$ and (\ref{eq42}) holds at $q$, with $\epsilon = 1$, but $F_q$ is not in the form of 
(\ref{eq39}) or (\ref{eq41}).
Perform a monoidal transform $\pi':Y'\rightarrow Y$ centered at the section over $C$ through $q$ with local equations $\overline x=\overline y=0$ in (\ref{eq42}).
Suppose that $q'\in (\pi')^{-1}(q)$. Suppose that there are permissible parameters  $(x_1,y_1,z)$  at $q'$ 
such that $\overline x=x_1$, $\overline y=x_1(y_1+\alpha)$ where $\alpha\ne 0$. Then
$$
F_q = \tau z^r + \sum_{i=2}^{r-1} \overline a_i x_1^{\alpha_i+\beta_i}(y_1+\alpha)^{\beta_i}z^{r-i}
+ x_1^{\alpha_r+\beta_r}(y_1+\alpha)^{\beta_r}
$$
Set $\tilde y_1 = (y_1+\alpha)^{\beta_r}-\alpha^{\beta_r}$.
Then
$$
F_{q'} =  \tau z^r + \sum_{i=2}^{r-1} \tilde a_i(x_1,\tilde y_1) x_1^{\tilde \alpha_i}z^{r-i}
+ x_1^{\tilde \alpha_r}\tilde y_1
$$
in the form of (\ref{eq39}).
Thus the only points $q'\in (\pi')^{-1}(q)$ which might not satisfy the conclusions of Theorem \ref{Theorem19}
are the points $q'$ which have regular parameters $(x_1,y_1,z)$ such that
$\overline x=x_1$, $\overline y=x_1y_1$ or $\overline x=x_1y_1$, $\overline y=y_1$.

The analysis of the case when (\ref{eq42}) holds at $q$, with $\epsilon = 0$, is simpler. We again conclude that
the only points in the blow up of the curve with local equations $\overline x=
\overline y=0$
 above $q$ which may not satisfy the conclusions of Theorem \ref{Theorem19}
are the points which have regular parameters $(x_1,y_1,z)$ such that
$\overline x=x_1$, $\overline y=x_1y_1$ or $\overline x=x_1y_1$, $\overline y=y_1$.

Suppose that $q\in \pi^{-1}(p)$ and (\ref{eq43}) holds at $q$, with $\epsilon = 1$, but $F_q$ is not in the form of 
(\ref{eq40}).
Perform a monoidal transform $\pi':Y'\rightarrow Y$ centered at the section over $C$ through $q$ with local equations $\overline x=\overline y=0$.
Suppose that $q'\in (\pi')^{-1}(q)$. Suppose that there are regular parameters  $(x_1,y_1,z)$  at $q'$ 
such that $\overline x=x_1$, $\overline y=x_1(y_1+\alpha)$ where $\alpha\ne 0$. Then

$$
\begin{array}{ll}
u&= \overline x_1^{(a+b)m}\\
(y_1+\alpha)^{\lambda}F_{q}&= \tau (y_1+\alpha)^{\lambda} z^r + \sum_{i=2}^{r-1}\overline a_i (y_1+\alpha)^{\beta_i+\lambda
-(\alpha_i+\beta_i)\frac{b}{a+b}}\overline x_1^{\alpha_i+\beta_i}z^{r-i}\\
&+ \overline x_1^{\alpha_r+\beta_r}(y_1+\alpha)^{\lambda+\beta_r
-(\alpha_r+\beta_r)\frac{b}{a+b}}.
\end{array}
$$
Thus 
$$
F_{q'}= \tilde \tau z^r + \sum_{i=2}^{r-1} \tilde a_i(\overline x_1,\overline y_1)\overline x_1^{\tilde \alpha_i}z^{r-i}
+ \overline x_1^{\tilde \alpha_r}\overline y_1
$$
 is in the form of (\ref{eq39}), where
$x_1=\overline x_1(y_1+\alpha)^{-\frac{b}{a+b}}$,
$\overline y_1= (y_1+\alpha)^{\lambda_1} - \alpha^{\lambda_1}$,
where 
$\lambda = d-\frac{b(c+d)}{a+b}$,
$\lambda_1=\lambda+\beta_r-(\alpha_r+\beta_r)\frac{b}{a+b}\ne 0$ since $F_q$ is normalized
implies 
$$
a(d+\beta_r)-b(c+\alpha_r)\ne 0.
$$

Thus the only points $q'\in (\pi')^{-1}(q)$ which might not satisfy the conclusions of Theorem \ref{Theorem19}
are the points which have permissible parameters $(x_1,y_1,z)$ such that
$\overline x=x_1$, $\overline y=x_1y_1$ or $\overline x=x_1y_1$, $\overline y=y_1$.

The analysis of the case when (\ref{eq43}) holds at $q$, with $\epsilon = 0$, is simpler. We again conclude that
the only points in the blow up of the curve with local equations $\overline x=\overline y=0$ above $q$ which may not satisfy the conclusions of Theorem \ref{Theorem19}
are the points which have regular parameters $(x_1,y_1,z)$ such that
$\overline x=x_1$, $\overline y=x_1y_1$ or $\overline x=x_1y_1$, $\overline y=y_1$.

We can  construct a sequence of monoidal transforms 
$$
Y_n\rightarrow\cdots\rightarrow Y_1\rightarrow Y
$$
with maps $\pi_i:Y_i\rightarrow Y$ such that $Y_i\rightarrow Y_{i-1}$ are
 centered at sections $C_i$ over $C$,
such that $\pi_i^{-1}(p)\cap C_i$ does not satisfy (\ref{eq39}),  (\ref{eq40}) or (\ref{eq41}).
By the above analysis, Lemma \ref{Lemma23} and Corollary \ref{Corollary24}, we reach the conclusions of the theorem
after a finite number of blowups.
\end{pf}

\begin{Remark}\label{Remark656}
In (\ref{eq41}) of Theorem \ref{Theorem19}, we must have $\beta_j<j$ for some $j$.
\end{Remark}

\begin{pf}
$$
u\in\hat{\cal I}_{\text{sing}(\Phi_X),q}\subset \hat{\cal I}_{\overline{S_r(X)},q}
$$
by Lemma \ref{Lemma657}. Thus $\overline x\in\hat{\cal I}_{\overline S_r(X),q}$.
$\beta_j\ge j$ for all $j$ implies $F_q\in (y,z)^r$, so that
$\hat{\cal I}_{\overline S_r(X),q}\subset (y,z)$ by Lemma \ref{Lemma302}, a contradiction.
\end{pf}

\begin{Theorem}\label{Theorem21}
Suppose that $r\ge2$ and $A_r(X)$ holds. 
Suppose that $p\in X$ is a 1 point or a 2 point with $\nu(p)=\gamma(p)=r$.
Let $R ={\cal O}_{X,p}$. Suppose that $\pi:Y\rightarrow \text{Spec}(\hat{R})$
is the sequence of monoidal transforms of sections over 
the curve $C$ with local equations $x=y=0$ at $p$
 of Theorem \ref{Theorem19}.
Suppose that $t>r$ is a positive integer.
Then there exists a sequence of permissible monoidal transforms $\overline \pi:\overline Y\rightarrow \text{spec}(R)$
of sections over $C$ such that for all $q\in \overline \pi^{-1}(p)$, $F_q$ is equivalent
mod $(\overline x,z)^t$ to a form (\ref{eq39}) or (\ref{eq41}) or $F_q$ is equivalent mod  
$(\overline x\overline y,z)^t$ to a form (\ref{eq40}), where $(\overline x,\overline y,z)$ are permissible parameters
for $u,v$ at $q$, and $z=\sigma\tilde z$ for some $\tilde z\in R$
and unit $\sigma\in\hat R$.

$\overline\pi$ extends to a sequence of permissible monoidal transforms $\overline U\rightarrow U$ over an affine neighborhood $U$ of $p$. $\overline S_r(\overline U)$ is the union of the
curves in $\overline \pi^{-1}(p)$ and the strict transforms of the curves $D$ or
$D_1,D_2$ (if they exist) in the notation of Lemma \ref{Lemma17}. $\overline S_r(\overline U)$
makes SNCs with $\overline B_2(\overline U)$.

\end{Theorem}

\begin{pf}
Let $m_0$ be the maximal ideal of $\hat R$. We can after possibly
replacing 
\begin{equation}\label{eq646}
\tilde x\text{ with }\tilde x\omega,
\end{equation}
 where $\omega$ is a unit in $R$,
assume that $x=\gamma\tilde x$ with $\gamma\equiv 1\text{ mod }m_0^{t}$.
We can factor
$$
Y=Y_{n'}\rightarrow \cdots \rightarrow Y_2\rightarrow Y_1\rightarrow \text{spec}(\hat R)=Y_0
$$
so that each map is a permissible monoidal transform. In fact, if $S_0=\text{spec}(k[x,y])$, there exists a sequence of
quadratic transforms
$$
S_{n'}\rightarrow \cdots\rightarrow S_2\rightarrow S_1\rightarrow \text{spec}(k[x,y])=S_0
$$
centered over $(x,y)$
such that $Y_i=S_i\times_{S_0}Y_0$ for all $i$.  Set $\overline S_0=\text{spec}(k[\tilde x,y])$.
$\tilde x\rightarrow x$ induces an isomorphism $\overline S_0\cong S_0$. We have then a sequence of quadratic
transforms
$$
\overline S_{n'}\rightarrow \cdots\rightarrow \overline S_1\rightarrow \overline S_0
$$
where $\overline S_i=S_i\times _{S_0}\overline S_0$ and isomorphisms $S_i\cong \overline S_i$.
Set $\overline Y_0=\text{spec}(R)$. We have a natural map $\overline Y_0\rightarrow \overline S_0$.
Define a sequence of permissible monoidal transforms
$$
\overline Y_{n'}\rightarrow \cdots\rightarrow \overline Y_1\rightarrow \overline Y_0
$$
by $\overline Y_i=\overline S_i\times_{\overline S_0}\overline Y_0$.
We have a commutative diagram 
\begin{equation}\label{eq501}
\begin{array}{lllllllll}
Y_{n'}&\rightarrow &Y_{n'-1}&\rightarrow&\cdots&\rightarrow&Y_1&\rightarrow&Y_0\\
\downarrow&&\downarrow&&&&\downarrow&&\downarrow\\
S_{n'}&\rightarrow&S_{n'-1}&\rightarrow&\cdots&\rightarrow&S_1&\rightarrow &S_0\\
\uparrow&&\uparrow&&&&\uparrow&&\uparrow\\
\overline S_{n'}&\rightarrow&\overline S_{n'-1}&\rightarrow&\cdots&\rightarrow&\overline S_1&\rightarrow &\overline S_0\\
\uparrow&&\uparrow&&&&\uparrow&&\uparrow\\
\overline Y_{n'}&\rightarrow &\overline Y_{n'-1}&\rightarrow&\cdots&\rightarrow&\overline Y_1&\rightarrow&\overline Y_0\\
\end{array}
\end{equation}
The maps  $\overline S_i\rightarrow S_i$ are isomorphisms, and we have maps 
$S_i\times_{S_0}\hat S_0\rightarrow Y_i$,
 $\overline S_i\times_{\overline S_0}\hat{\overline S_0}\rightarrow \overline Y_i$
induced by the natural projections
$$
k[[x,y,z]]\rightarrow k[[x,y]]
\text{ and }
k[[\tilde x,y,z]]\rightarrow k[[\tilde x,y]]
$$
 so that the diagrams 
\begin{equation}\label{eq633}
\begin{array}{lll}
&&Y_i\\
&\swarrow&\\
S_i&&\uparrow\\
&\nwarrow&\\
&&S_i\times_{S_0}\hat S_0
\end{array}
\end{equation}
and
$$
\begin{array}{lll}
&&\overline Y_i\\
&\swarrow&\\
\overline S_i&&\uparrow\\
&\nwarrow&\\
&&\overline S_i\times_{\overline S_0}\hat S_0
\end{array}
$$
commute.

Suppose that $\tilde q\in \overline Y_{n'}$ is a closed point. (\ref{eq501})
and (\ref{eq633}) identifies $\tilde q$ with a closed point $\tilde p\in Y_{n'}$,
and closed points $\overline q\in \overline S_{n'}$, $\overline p\in S_{n'}$.
 We have commutative diagrams: 
\begin{equation}\label{eq502}
\begin{array}{lll}
\hat{\cal O}_{S_{n'},\overline p}&\cong &\hat{\cal O}_{\overline S_{n'},\overline q}\\
\uparrow&&\uparrow\\
k[[x,y]]&\cong&k[[\tilde x,y]]
\end{array}
\end{equation}
induced by $x\rightarrow \tilde x$. This induces  commutative diagrams:
$$
\begin{array}{llll}
\hat{\cal O}_{Y_{n'},\tilde p} = \hat{\cal O}_{S_{n'},\overline p}[[z]] &\cong &\hat{\cal O}_{\overline S_{n'},\overline q}[[z]]
=&\hat{\cal O}_{\overline Y_{n'},\tilde q}\\
\uparrow&&&\uparrow\\
\hat R=k[[x,y,z]]&\cong&&k[[\tilde x,y,z]]
\end{array}.
$$

Suppose that $F(x,y,z)\in \hat R$.
If $(\overline x,\overline y)$ are regular parameters in $\hat{\cal O}_{ S_{n'},\overline p}$, which are
identified with regular parameters $\hat x,\hat y$ in 
$\hat{\cal O}_{\overline S_{n'},\overline q}$ by (\ref{eq502}), and
$F(x,y,z)=G(\overline x,\overline y,z)\in \hat{\cal O}_{Y_{n'},\tilde p}$, then
$F(\tilde x,y,z)=G(\hat x,\hat y,z)\in \hat{\cal O}_{\overline Y_{n'},\tilde q}$.

Since $x=\gamma\tilde x$, $\gamma\equiv 1\text{ mod }m_0^{t}$, and 
$$
F(x,y,z)=F(\gamma\tilde x,y,z)\equiv F(\tilde x,y,z)\text{ mod }m_0^{t}
$$
implies
$$
F(x,y,z)\equiv G(\hat x,\hat y,z)\text{ mod }m_0^t\hat{\cal O}_{\overline Y_{n'},\tilde q}.
$$
Let $m$ be the maximal ideal of $R$.

Suppose that $\tilde p\in Y_{n'}$ is a 1 point, $u=\overline x^a$ in $\hat{\cal O}_{Y_{n'},\tilde p}$. Then (since we can assume that $Y\ne \text{spec}(\hat R)$)
$\overline x\mid x$ and $\overline x\mid y$
 in $\hat{\cal O}_{Y_{n'},\tilde p}$ implies $m_0^{t}\subset (\overline x,z)^{t}$, so that $m^{t}\subset (\hat x,z)^{t}{\cal O}_{\overline Y_{n'},\tilde q}$.

Suppose that $\tilde p\in Y_{n'}$ is a 2 point, $u=(\overline x^a\overline y^b)^m$ in $\hat{\cal O}_{Y_{n'},\tilde p}$. If $\overline x=0,\overline y=0$ are both
local equations of components of the exceptional locus of $Y_{n'}\rightarrow Y_0$, we have 
$\overline x\overline y\mid  x$,
$\overline x\overline y\mid y$ which implies that $m_0^{t}\subset (\overline x\overline y,z)^{t}$, so that
$m^{t}\subset (\hat x\hat y,z)^{t}$.

If one of $\overline x=0, \overline y=0$ is not a local equation of the exceptional locus, then we have regular parameters
$(x',y')$ in $\hat{\cal O}_{ S_{n'},\overline p}$ such that
$x=x'(y')^{\overline b}, y=y'$ (or $x=x', y=y'(x')^{\overline b}$). In the first case we have
$$
\begin{array}{ll}
\hat{\cal O}_{Y_{n'},\tilde p} &=k[[\frac{ x}{y^{\overline b}},y,z]]\\
&=k[[\frac{\gamma\tilde x}{y^{\overline b}},y,z]]\\
&=k[[\frac{\tilde x}{y^{\overline b}},y,z]]=\hat{\cal O}_{\overline Y_{n'},\tilde q}.
\end{array}
$$
Thus we have $\hat{\cal O}_{Y_{n'},\tilde p}=\hat{\cal O}_{\overline Y_{n'},\tilde q}$.
In the second case, we also have $\hat{\cal O}_{Y_{n'},\tilde p}=\hat{\cal O}_{\overline Y_{n'},\tilde q}$.

Let $\overline\pi:\overline Y=\overline Y_{n'}\rightarrow\overline Y_0$ be the morphism of 
the bottom row of (\ref{eq501}).
Suppose that $p_0\in\overline Y_0$ is a 1 point, so that in $\hat{\cal O}_{\overline Y_0,p_0}$,
$$
u=x^a,v=P_{p_0}(x)+x^{c_0}F_{p_0}.
$$
Suppose that $p'\in \overline \pi^{-1}(p_0)$. Let $q$ be the corresponding closed point of $Y=Y_{n'}$.
 Suppose that we have permissible parameters $(\overline x,\overline y,z)$ in $\hat{\cal O}_{Y,q}$ such that 
\begin{equation}\label{eq503}
u=\overline x^{\overline a}, v=P_{q}(\overline x)+\overline x^cF_{q}(\overline x,\overline y,z)
\end{equation}
of the form  (\ref{eq39}) or (\ref{eq41})
with $\frac{\overline a}{a}$, $\frac{c}{c_0}\in\bold N$, $x=\overline x^{\frac{\overline a}{a}}=\overline x^{\frac{c}{c_0}}$. Let $(\tilde x_*,\tilde y_*,z)$ be the corresponding
 regular parameters at $p'$ (by the identification (\ref{eq502})).
$$
u=x^a=\gamma^a\tilde x^a=\gamma^a\tilde x_*^{\overline a}=(\tilde x_*')^{\overline a}
$$
where we define
$$
\tilde x_*=\tilde x_*'\gamma^{-\frac{a}{\overline a}}\equiv \tilde x'_*\text{ mod }m_0^t\hat{\cal O}_{\overline Y,p'}.
$$
Thus $(\tilde x'_*,\tilde y_*,z)$ are permissible parameters for $(u,v)$ in $\hat{\cal O}_{\overline Y,p'}$

There exists a series $\tilde P_q(\overline x)$ such that 
$$
F_{p_0}(x,y,z)=F_q(\overline x,\overline y,z)+\tilde P_q(\overline x).
$$
Thus
$$
\begin{array}{ll}
F_{p_0}(x,y,z)&\equiv F_q(\tilde x_*,\tilde y_*,z)+\tilde P_q(\tilde x_*)
\text{ mod }m_0^{t}\hat{\cal O}_{\overline Y,p'}\\
&\equiv F_q(\tilde x_*',\tilde y_*,z)+\tilde P_q(\tilde x_*')\text{ mod }m_0^{t}
\hat{\cal O}_{\overline Y,p'}
\end{array}
$$
$$
P_q(\overline x)=P_{p_0}(\overline x^{\frac{c}{c_0}})+\overline x^c\tilde P_q(\overline x)
$$
implies
$$
\begin{array}{ll}
v&=P_{p_0}(x)+x^{c_0}F_{p_0}(x,y,z)\\
&=P_{p_0}((\tilde x_*')^{\frac{c}{c_0}})+(\tilde x_*')^cF_{p_0}(x,y,z)\\
&=P_q(\tilde x_*')+(\tilde x_*')^c(F_q(\tilde x_*',\tilde y_*,z)+h)
\end{array}
$$
with $h\in m_0^{t}{\cal O}_{\overline Y,p'}$.

The  case when $p_0\in\overline Y_0$ is a 1 point and (\ref{eq40}) holds 
in $\hat{\cal O}_{Y,q}$
is 
a combination of the case when $p_0$ is a 1 point and the form (\ref{eq39}) or (\ref{eq41})
holds in $\hat{\cal O}_{Y,q}$, and the following case.

Suppose that $p_0\in\overline Y_0$ is a 2 point so that in $\hat{\cal O}_{\overline Y_0,p_0}$,
$$
\begin{array}{ll}
u&=(x^ay^b)^{m_0}\\
v&=P_{p_0}(x^ay^b)+x^{c_0}y^{d_0}F_{p_0}
\end{array}
$$

Suppose that $p'\in \overline \pi^{-1}(p_0)$. Let $q$ be the corresponding closed point in 
$Y=Y_{n'}$.

Suppose that we have permissible parameters $(\overline x,\overline y,z)$ at $q$ such that
$$
\begin{array}{ll}
u&=(\overline x^{\overline a}\overline y^{\overline b})^{m_1}\\
v&=P_{q}(\overline x^{\overline a}\overline y^{\overline b})
+\overline x^{c}\overline y^{d}F_q(\overline x,\overline y,z)
\end{array}
$$
of the form of (\ref{eq40}). We have
$$
\begin{array}{ll}
x&=\overline x^{a_1}\overline y^{b_1}\gamma_1(\overline x,\overline y)\\
y&=\overline x^{a_2}\overline y^{b_2}\gamma_2(\overline x,\overline y)
\end{array}
$$
where $\gamma_1,\gamma_2$ are units in $\hat{\cal O}_{Y,q}$, such that
$$
(aa_1+ba_2)m_0=\overline am_1, (ab_1+bb_2)m_0=\overline bm_1
$$
where $m_0\mid m_1$, $c_0a_1+d_0a_2=c, c_0b_1+d_0b_2=d$.

We have
$$
x^ay^b=(\overline x^{\overline a}\overline y^{\overline b})^{\frac{m_1}{m_0}}
$$
and
$$
x^{c_0}y^{d_0}=\overline x^c\overline y^d\phi(\overline x,\overline y)
$$
where $\phi=\gamma_1^{c_0}\gamma_2^{d_0}$.
There exists a series $\tilde P_q(\overline x^{\overline a}\overline y^{\overline b})$
such that
$$
F_q(\overline x,\overline y,z)=\phi(\overline x,\overline y)F_{p_0}(x,y,z)
-\frac{\tilde P_q(\overline x^{\overline a}\overline y^{\overline b})}
{\overline x^c\overline y^d}.
$$
$$
v=P_{p_0}((\overline x^{\overline a}\overline y^{\overline b})^{\frac{m_1}{m_0}})
+\tilde P_q(\overline x^{\overline a}\overline y^{\overline b})
+\overline x^c\overline y^dF_q(\overline x,\overline y,z)
$$
implies
$$
P_q(\overline x^a\overline y^b)=P_{p_0}((\overline x^{\overline a}\overline y^{\overline b})
^{\frac{m_1}{m_0}})+\tilde P_q(\overline x^{\overline a}\overline y^{\overline b}).
$$
Let $(\tilde x_*,\tilde y_*,z)$ be the corresponding regular parameters at $p'$
to $(\overline x,\overline y,\overline z)$, by the
identification of (\ref{eq502}). Define $\tilde x_*'$ by
$$
\tilde x_*=\tilde x_*'\gamma^{-\frac{m_0a}{m_1\overline a}}
\equiv \tilde x_*'\text{ mod }m_0^t{\cal O}_{\overline Y,p'}.
$$
$$
x^ay^b=(\overline x^{\overline a}\overline y^{\overline b})^{\frac{m_1}{m_0}}
$$
implies
$$
x^ay^b=\gamma^a\tilde x^ay^b=\gamma^a(\tilde x_*^{\overline a}\tilde y_*^{\overline b})^{\frac{m_1}{m_0}}=((\tilde x_*')^{\overline a}\tilde y_*^{\overline b})^{\frac{m_1}{m_0}}
$$
Thus $(\tilde x_*',\tilde y_*,z)$ are permissible parameters for $(u,v)$ in $\hat{\cal O}_{\overline Y,p'}$. $x^{c_0}y^{d_0}=\overline x^c\overline y^d\phi$ implies
$$
x^{c_0}y^{d_0}=\gamma^{c_0}\tilde x^{c_0}y^{d_0}=\gamma^{c_0}\tilde x_*^c\tilde y_*^d
\phi(\tilde x_*,\tilde y_*)
=(\tilde x_*')^c\tilde y_*^d\tilde\phi
$$
with 
$$
\tilde\phi\equiv \phi(\tilde x_*',\tilde y_*)\text{ mod }m_0^{t}{\cal O}_{\overline Y,p'}.
$$
$$
\begin{array}{ll}
F_{p_0}(x,y,z)&\equiv \phi(\tilde x_*,\tilde y_*)^{-1}(F_q(\tilde x_*,\tilde y_*,z)+
\frac{\tilde P_q(\tilde x_*^{\overline a}\tilde y_*^{\overline b})}{\tilde x_*^c\tilde y_*^d})
\text{ mod }m_0^{t}\hat{\cal O}_{\overline Y,p'}\\
&\equiv \tilde \phi(\tilde x_*',\tilde y_*)^{-1}(F_q(\tilde x_*',\tilde y_*,z)+
\frac{\tilde P_q((\tilde x_*')^{\overline a}\tilde y_*^{\overline b})}{(\tilde x_*')^c\tilde y_*^d})
\text{ mod }m_0^{t}\hat{\cal O}_{\overline Y,p'}
\end{array}
$$
$$
\begin{array}{ll}
u&=((\tilde x_*')^{\overline a}\tilde y_*^{\overline b})^{m_1}\\
v&=P_{p_0}(((\tilde x_*')^{\overline a}\tilde y_*^{\overline b})^{\frac{m_1}{m_0}})
+(\tilde x_*')^c(\tilde y_*^d)\tilde\phi F_{p_0}(x,y,z)\\
&=P_{p_0}(((\tilde x_*')^{\overline a}\tilde y_*^{\overline b})^{\frac{m_1}{m_0}})
+\tilde P_q((\tilde x_*')^{\overline a}\tilde y_*^{\overline b})
+(\tilde x_*')^c(\tilde y_*^d)[ F_q(\tilde x_*',\tilde y_*,z)+h]\\
&=P_q((\tilde x_*')^{\overline a}\tilde y_*^{\overline b})+
(\tilde x_*')^c(\tilde y_*^d)[ F_q(\tilde x_*',\tilde y_*,z)+h]
\end{array}
$$
with $h\in m_0^{t}\hat{\cal O}_{\overline Y,p'}$.

The case when $p_0\in\overline Y_0$ is a 2 point and (\ref{eq39}) or (\ref{eq41})
holds in $\hat{\cal O}_{Y,q}$ is similar to the case when (\ref{eq40}) holds in $\hat{\cal O}_{Y,q}$.

Suppose that $p'$ is a generic point of $C$.
If $(x',y',z')$ are permissible parameters at $p'$ such that $x'=y'=0$ are local
equations of $C$ at $p'$, then $\nu(F_{p'}(0,0,z'))\le 1$,
 so that after extending $\overline \pi$ to a
sequence of permissible blowups $\overline U\rightarrow U$ over a
small affine neighborhood $U$ of $p$, $\gamma(p*)\le 1$ at all points $p*$ of
$\overline \pi^{-1}(p')$. Thus the curves in $\overline S_r(\overline U)$ must be components of
$\overline \pi^{-1}(p)$, and the strict transforms of the curves $D$ or $D_1$, $D_2$
in $\overline S_r(\overline U)$,
(if they exist), with the notation of Lemma \ref{Lemma17}.

Thus a curve $E$ in $\overline S_r(\overline U)$ must have local equations as asserted by the 
Theorem.
\end{pf}

\begin{Theorem}\label{Theorem22}
Suppose that $r\ge 3$ and $A_r(X)$ holds.
Suppose that $p\in X$ is a 1 point or a 2 point with $\nu(p)=\gamma(p)=r$.
Let $R ={\cal O}_{X,p}$. Suppose that $\pi:Y_p\rightarrow \text{Spec}(\hat R)$
is the sequence of monoidal transforms of sections over the curve $\overline C$
with local equations $\tilde x=y=0$ of Theorem \ref{Theorem19}. For $q\in\pi^{-1}(p)$, define
$$
l_q=\left\{\begin{array}{ll}
\left(\text{min}_{2\le i\le r}\{[\frac{\alpha_i}{i}]\}+3\right)r
& \text{if $F_q$ is a form $(\ref{eq39})$ or $(\ref{eq41})$.}\\
\left(\text{min}_{2\le i\le r}\{[\frac{\alpha_i}{i}]\}+
\text{min}_{2\le i\le r}\{[\frac{\beta_i}{i}]\}+3\right)r
& \text{if $F_q$  is a form $(\ref{eq40})$.}
\end{array}\right.
$$
let $l=\text{max}\{l_q\mid q\in  \pi^{-1}(p)\}$.

Suppose that $t\ge l$. Let $\overline\pi:\overline Y_p\rightarrow \text{spec}(R)$ be
the sequence of monodial transforms of Theorem \ref{Theorem21}. Let
$$
\cdots\rightarrow Y_n\rightarrow\cdots\rightarrow Y_1\rightarrow \overline Y_p
$$
be a sequence of permissible monoidal transforms centered at curves $C$ in $\overline S_r$
such that $C$ is r big. Then there exists
$n_0<\infty$ such that
$$
V_p=Y_{n_0}\stackrel{\pi_1}{\rightarrow}\overline Y_p\rightarrow\text{spec}(R)
$$
extends to a permissible sequence of monoidal transforms 
$$
\overline U_1\rightarrow \overline U\rightarrow U
$$
over an affine neighborhood $U$ of $p$, in the notation of Theorem \ref{Theorem21},
such that $\overline S_r(\overline U_1)$ contains no curves $C$ such that 
$C$ is r big. Let
$$
\cdots\rightarrow Z_n\rightarrow\cdots\rightarrow V_p
$$
be a permissible sequence of monoidal transforms centered at  curves $C$ in $\overline S_r$
such that $C$ is r small. Then there exists $n_1<\infty$
such that $\pi_2:Z_p=Z_{n_1}\rightarrow V_p$ extends to a permissible sequence of
monoidal transforms
$$
\overline U_2\rightarrow \overline U_1\rightarrow \overline U\rightarrow U
$$
over an affine neighborhood $U$ of $p$ such that $\overline S_r(\overline U_2)=\emptyset$.

Finally, there exists a sequence of quadratic transforms $\pi_3:W_p\rightarrow Z_{p}$
which extends to a permissible sequence of monoidal transforms
$$
\overline U_3\rightarrow\overline U_2\rightarrow \overline U_1\rightarrow \overline U\rightarrow U
$$
over an affine neighborhood $U$ of $p$ such that $\overline S_r(\overline U_3)=\emptyset$,
and if 
$$
q\in(\overline\pi\circ\pi_1\circ\pi_2\circ\pi_3)^{-1}(p),
$$
\begin{enumerate}
\item $\nu(q)\le r-1$ if $q$ is a 1 or 2 point.
\item If $q$ is a 2 point and $\nu(q)=r-1$, then $\tau(q)>0$.
\item $\nu(q)\le r-2$ if $q$ is a 3 point.
\end{enumerate}
\end{Theorem}

\begin{pf} Let $\overline Y=\overline Y_p$.
If $q\in \overline \pi^{-1}(p)$, we have permissible parameters $(x,y,z)$ in $\hat{\cal O}_{\overline Y,q}$ for
$(u,v)$ with forms obtained from those of (\ref{eq39}), (\ref{eq40}), (\ref{eq41}) by modifying $F_q$ by adding an 
appropriate series $h$ to $F_q$. 

(\ref{eq39}) is modified by changing $F_q$ to 

\begin{equation}\label{eq634}
F_q=\tau z^r+\sum_{j=2}^{r-1}\overline a_j(x,y)x^{\alpha_j}z^{r-j}+\epsilon x^{\alpha_r}y+h
\end{equation}
with $h\in (x,z)^t$. By assumption, there exists $\tilde z\in{\cal O}_{\overline Y,q}$ and
a unit $\sigma\in\hat{\cal O}_{\overline Y,q}$ such that $z=\sigma\tilde z$. Then
 $x=z=0$ defines a germ of an algebraic curve $D$ at $q$. We also assume that $\nu(q)=r$.

Given a form (\ref{eq634}) at $q$, suppose that $D\subset\overline S_r(\overline Y)$ is a curve
such that $q\in D$. By assumption, $\overline S_r(\overline Y)$ makes SNCs with
$\overline B_2(\overline Y)$. Since  $D$ is nonsingular at $q$, $x\in\hat{\cal I}_{D,q}$ and $F_q\in\hat{\cal I}_{D,q}^r+(x)^{r-1}$ by Lemma \ref{Lemma5} implies
$$
\frac{\partial^{r-1} F_q}{\partial z^{r-1}}\in\hat{\cal I}_{D,q},
$$
 so that
$z\in \hat{\cal I}_{D,q}$ and $x=z=0$ are local equations of $D$ at $q$.

 (\ref{eq40}) is modified by changing $F_q$ to 

\begin{equation}\label{eq505}
F_q=\tau z^r+\sum_{j=2}^{r-1}\overline a_j(x,y)x^{\alpha_j}y^{\beta_j}z^{r-j}+\epsilon x^{\alpha_r}y^{\beta_r}+h
\end{equation}
with $h\in(xy,z)^t$. By assumption, there exists $\tilde z\in {\cal O}_{\overline Y,q}$
and a unit $\sigma\in\hat{\cal O}_{\overline Y,q}$ such that $z=\sigma\tilde z$. Then
$x=z=0$ and $x=y=0$ define germs of algebraic curves $D_1$ and $D_2$ at $q$.
We also assume $\nu(q)=r$.

Given a form (\ref{eq505}) at $q$, suppose that $D\subset \overline S_r(\overline Y)$ is
a curve such that $q\in D$. Since (by assumption) $\overline S_r(\overline Y)$ makes SNCs with $\overline B_2(\overline Y)$, either $x$ or $y\in\hat{\cal I}_{D,q}$, and by Lemma \ref{Lemma7}, there exist
$d_i\in k$ such that 
$$
F_q-\frac{1}{x^cy^d}\sum d_i(x^ay^b)^i\in\hat{\cal I}_{D,q}^r+(x)^{r-1}
$$
or
$$
F_q-\frac{1}{x^cy^d}\sum d_i(x^ay^b)^i\in\hat{\cal I}_{D,q}^r+(y)^{r-1}.
$$
$$
\frac{\partial^{r-1}F_q}{\partial z^{r-1}}\in\hat{\cal I}_{D,q}
$$
implies $z\in \hat{\cal I}_{D,q}$ so that either $x=z=0$ or $y=z=0$ are local equations
of $D$ at $q$.

 (\ref{eq41}) is modified by changing $F_q$ to 
\begin{equation}\label{eq506}
F_q=\tau z^r+\sum_{j=2}^{r-1}\overline a_j(x,y)x^{\alpha_j}y^{\beta_j}z^{r-j}+\epsilon x^{\alpha_r}y^{\beta_r}+h
\end{equation}
with $h\in(x,z)^t$. By assumption, there exists $\tilde z\in{\cal O}_{\overline Y,q}$
and a unit $\sigma\in\hat{\cal O}_{\overline Y,q}$ such that $z=\sigma\tilde z$.
Then $x=z=0$ defines a germ of an algebraic curve $D$ at $q$. We also assume $\nu(q)=r$.

Given a form (\ref{eq506}) at $q$, suppose that $D\subset \overline S_r(\overline Y)$
is a curve such that $q\in D$.
By assumption, $\overline S_r(\overline Y)$ makes SNCs with
$\overline B_2(\overline Y)$.  As in the analysis of the case when (\ref{eq634}) holds,
we conclude that $x=z=0$ are local equations of $D$ at $q$.

Suppose that $D\subset\overline S_r(\overline Y)$ is a curve such that $D$ is r big. Let $\pi':Y'\rightarrow \overline Y$ be the blowup of $D$.
By assumption, $\pi'$ is a permissible monodial transform. $\overline S_r(Y')$ makes
SNCs with $\overline B_2(Y')$ by Lemma \ref{Lemma654}.

First suppose that $q\in \overline \pi^{-1}(p)\cap D$, and that $u,v$ have the form of (\ref{eq634}). Then $x=z=0$ are local equations of $D$ at $q$.
  $F_q\in (x,z)^r$ implies $\alpha_r\ge r$
if $\epsilon=1$. 
 Suppose that $q'\in (\pi')^{-1}(q)$. First suppose that
$q'$ has permissible parameters $(x_1,y,z_1)$ where $x=x_1$, $z=x_1(z_1+\alpha)$ for some $\alpha\ne 0$.
Substituting into $F_q$, we get $\nu(F_{q'}(0,0,z_1))\le r-1$.
Suppose that $q'$ has regular parameters $(x_1,y,z_1)$ where $x=x_1z_1$, $z=z_1$. Then $F_{q'}$ is a unit.
The remaining case is when $q'$ has permissible parameters $(x_1,y,z_1)$ where $x=x_1$, $z=x_1z_1$.
Then 
\begin{equation}\label{eq47}
\begin{array}{ll}
u&=x_1^a\\
F_{q'} &= \tau z_1^r + \sum_{j=2}^{r-1}\overline a_j(x_1,y)x_1^{\alpha_j'}z_1^{r-j}+\epsilon x_1^{\alpha_r'}y+h_1
\end{array}
\end{equation}
where $\alpha_j' = \alpha_j-j$ for $2\le j\le r$ and $h_1\in (x_1,z_1)^{t-r}$.
We either have a reduction in multiplicity $\nu(q')<r$, or $\nu(q')=r$ and we are back in the form of (\ref{eq634}) with
a reduction in  the $\alpha_j$ by $j$, and a decrease of $t$ by $r$. 
By Lemma \ref{Lemma654}, $\overline S_r(Y')\cup \overline B_2(Y')$ makes SNCs in a 
neighborhood of $(\pi')^{-1}(q)$. Since $S_r(Y')$ is closed in the open set of
1 points of $E_{Y'}$, and $r\ge 2$, by Lemma \ref{Lemma4}, $\overline S_r(Y')\cap (\pi')^{-1}(q)=\emptyset$
or it is the point $q'$ of (\ref{eq47}), if $\nu(q')=r$.
Suppose that $\nu(q')=r$ in (\ref{eq47}). Since there exists a unit series $\sigma$ such that
$\sigma z\in {\cal O}_{\overline Y,q}$, 
there exists a unit series $\sigma'$ such that $\sigma'z_1\in{\cal O}_{Y',q'}$.

Now suppose that $q\in \overline \pi^{-1}(p)\cap D$, and that $u,v$ have the form of (\ref{eq506}).
  Then we  have that $x=z=0$  are local equations of $D$ at $q$. 
Since $F_q\in\hat{\cal I}_{D,q}^r$, we have $\alpha_j\ge j$ for all $j$.

Suppose that $q'\in (\pi')^{-1}(q)$. First suppose that $q'$ 
 has permissible parameters $(x_1,y,z_1)$ where $x=x_1$, $z=x_1(z_1+\alpha)$ for some $\alpha\ne 0$.
substituting into $F_q$, we get $\nu(F_{q'}(0,0,z_1))\le r-1$.
Suppose that $q'$ has regular parameters $(x_1,y,z_1)$ where $x=x_1z_1$, $z=z_1$. Then $F_{q'}$ is a unit.
The remaining case is when $q'$ has regular parameters $(x_1,y,z_1)$ where $x=x_1$, $z=x_1z_1$.
Then 
\begin{equation}\label{eq48}
\begin{array}{ll}
u&=x_1^a\\
F_{q'} &= \tau z_1^r + \sum_{j=2}^{r-1}\overline a_j(x_1,y)x_1^{\alpha_j'}y^{\beta_j}z_1^{r-j}
+\epsilon x_1^{\alpha_r'}y^{\beta_r}+h_1
\end{array}
\end{equation}
where $\alpha_j' = \alpha_j-j$ for $2\le j\le r$ and $h_1\in (x_1,z_1)^{t-r}$.
We either have a reduction in multiplicity $\nu(q')<r$, or we are back in the form of (\ref{eq506}) with
a reduction in  $\alpha_i$ by $i$, and a decrease of $t$ by $r$.
As in the analysis of (\ref{eq634}), we either have 
$\overline S_r(Y')\cap (\pi')^{-1}(q)=\emptyset$, or $\overline S_r(Y')\cap(\pi')^{-1}(q)$
is the single point $q'$ of (\ref{eq48}). In this case there exists a unit series $\sigma'$ such 
that $\sigma' z_1\in {\cal O}_{Y',q'}$.

Now suppose that $q\in \overline\pi^{-1}(p)\cap D$, and that $u,v$ have the form of (\ref{eq505}). 
Then either $x=z=0$ or $y=z=0$ are local equations of $D$ at $q$. We may suppose that
$x=z=0$ are local equations of $D$ at $q$, so that
 $F_q\in (x,z)^r$.

 Suppose that $q'\in (\pi')^{-1}(q)$. First suppose that $q'$ 
 has permissible parameters $(x_1,y,z_1)$ where $x=x_1$, $z=x_1(z_1+\alpha)$ for some $\alpha\ne 0$.
substituting into $F_q$, we get $u=(x_1^ay^b)^m$ and $\nu(F_{q'}(0,0,z_1))\le r-1$.
Suppose that $q'$ has permissible parameters $(x_1,y,z_1)$ where $x=x_1z_1$, $z=z_1$. Then $F_{q'}$ is a unit.
The remaining case is when $q'$ has permissible parameters $(x_1,y,z_1)$ where $x=x_1$, $z=x_1z_1$.
Then 
\begin{equation}\label{eq49}
\begin{array}{ll}
u&=(x_1^ay^b)^m\\
v&=P(x_1^ay^b)+x_1^{c+r}y^dF_{q'}\\
F_{q'} &= \tau z_1^r + \sum_{j=2}^{r-1}\overline a_j(x_1,y)x_1^{\alpha_j'}y^{\beta_j}z_1^{r-j}
+\epsilon x_1^{\alpha_r'}y^{\beta_r}+h_1
\end{array}
\end{equation}
where $\alpha_j' = \alpha_j-j$ for $2\le j\le r$ and $h_1\in (x_1y_1,z_1)^{t-r}$.
We either have a reduction in multiplicity $\nu(q')<r$, or we are back in the form of (\ref{eq505}) with
a reduction of $\alpha_i$ by $i$ and a decrease of $t$ by $r$. 
In this case, there exists a unit series $\sigma'$ such that $\sigma' z_1\in{\cal O}_{Y',q'}$.

Suppose that $q'\in \overline S_r(Y')\cap (\pi')^{-1}(q)$ is a 2 point with $\nu(q')=r-1$,
and $q'$ lies on a curve $E$ in $\overline S_r(Y'')$. $E$  is transversal to the
2 curve at $q'$ by Lemma \ref{Lemma654}. By Lemma \ref{Lemma7}, there exist $b_t\in k$ such that
$$
F_{q'}+\frac{1}{x_1^{c+r}y^d}\sum b_t(x_1^ay^b)^t\in\hat{\cal I}_{E,q'}^r+(x_1)^{r-1},
$$
if $x_1\in\hat{\cal I}_{E,q'}$ or the series is in $\hat{\cal I}_{E,q'}^r+(y)^{r-1}$ if 
$y\in\hat{\cal I}_{E,q'}$. Then $\frac{\partial^{r-2}F_{q'}}{\partial z^{r-2}}\in\hat{\cal I}_{E,q'}^2$.
Suppose that $q'$ has permissible parameters $(x_1,y,z_1)$ such that
$x=x_1,z=x_1(z_1+\alpha)$ with $\alpha\ne 0$. 
$$
F_{q'}=\Lambda z_1^{r-1}+\sum_{i=1}^{r-1}\tilde a_i(x_1,y)z_1^{r-1-i}
$$
where $\Lambda$ is a unit. $x_1^i\mid \tilde a_i$ (or $y^i\mid\tilde a_i$) for
$1\le i\le r-1$ since $F_{q'}\in \hat{\cal I}_{E,q'}^{r-1}$. 
Then $z_1\in\hat{\cal I}_{E,q'}^2+(x_1)$
which is impossible.

 Thus, by Lemma \ref{Lemma4}, the only possible point in $\overline S_r(Y')\cap (\pi')^{-1}(q)$
is the point $q'$ of (\ref{eq49}). 
If there is a curve $E\subset \overline S_r(Y')$ containing $q'$, then we must have
$z_1\in \hat{\cal I}_{E,q'}$
 since 
$$
\frac{\partial^{r-1} F_{q'}}{\partial z^{r-1}}\in\hat{\cal I}_{E,q'}.
$$
Thus $E$ has local equations $x_1=z_1=0$ or $y=z_1=0$.

After any sequence of permissible monoidal transforms, centered at r big curves $C\subset\overline S_r$,  we eventually obtain $\pi_1:V_p\rightarrow \overline Y$
where there are no r big curves $C$ in $\overline S_r(V_p)$
and $\overline S_r(V_p)$ makes SNCs with $\overline B_2(V_p)$.

Further, if $q\in(\overline \pi\circ \pi_1)^{-1}(p)$, and either $q\in\overline S_r(V_p)$ or one of
1. - 3. of the conclusions of $W_p$ fail at $q$ then $q$ must satisfy one of 
(\ref{eq634}), (\ref{eq505}) or (\ref{eq506}) (with $\nu(q)=r$ or $\nu(q)=r-1$).

Suppose that (\ref{eq505}) holds at $q\in(\overline \pi\circ\pi_1)^{-1}(p)$,
and $\nu(q)=r$. Then we either have $\alpha_j\ge j$
for all $j$ or $\beta_j\ge j$ for all $j$ by Remark \ref{Remark1049}. If $\alpha_j\ge j$ for all $j$, then $F_q\in(x,z)^r$. If 
$\beta_j\ge j$ for all $j$, then $F_q\in(y,z)^r$. Since $x=z=0$ and $y=z=0$ are local
equations of curves on $V_p$, in either case we have a curve $D\subset\overline S_r(V_p)$
such that $D$ is r big by Lemma \ref{Lemma653}. Thus (\ref{eq505})
cannot hold on $V_p$ with $\nu(q)=r$.

Suppose that (\ref{eq506}) holds at $q\in(\overline \pi\circ\pi_1)^{-1}(p)$,
and $\nu(q)=r$. Then we either have $\alpha_j\ge j$
for all $j$ or $\beta_j\ge j$ for all $j$ by Remark \ref{Remark1049}. By Remark \ref{Remark656}, $\beta_j<j$ for some
$j$. Thus $F_q\in(x,z)^r$.
Since $x=z=0$ are local equations of a curve in $V_p$, By Lemma \ref{Lemma653},
we must have a curve $D\subset\overline S_r(V_p)$ such that $D$ is r big.
Thus (\ref{eq506}) cannot hold on $V_p$, with $\nu(q)=r$.

  The only points on $(\overline\pi\circ\pi_1)^{-1}(p)$ where the conclusions of the Theorem do not
hold are at points $q'$ over $p$ where  one of (\ref{eq972}) or (\ref{eq50}) following hold.

\begin{equation}\label{eq972}
\begin{array}{ll}
u&=x^a\\
v&=P(x)+x^cF_q\\
F_q&=\tau z^r+\sum_{j=2}^{r-1}\overline a_j(x,y)x^{\alpha_j}z^{r-j}
+\epsilon x^{\alpha_r}y+h
\end{array}
\end{equation}
where $\nu(q')=r$, some $\alpha_j<j$,  and $h\in (x,z)^{3r}$. Further, there exists a series
$\sigma$ such that $\sigma z\in{\cal O}_{V_p,q'}$.  The other possiblity is that 
$q'$ has permissible parameters $(x_1,y_1,z_1)$ of the form 
of (\ref{eq49}), with 
\begin{equation}\label{eq50}
\begin{array}{ll}
u&=(x_1^ay_1^b)^m\\
v&=P(x_1^ay_1^b)+x_1^cy_1^dF_{q'}\\
F_{q'} &= \tau z_1^r + \sum_{j=2}^{r-1}\overline a_j(x_1,y_1)x_1^{\alpha_j'}y_1^{\beta_j'}z_1^{r-j}
+ \epsilon x_1^{\alpha_r'}y_1^{\beta_r'}+h_1
\end{array}
\end{equation} 
with $\nu(q')=r-1$, 
$$
h_1\in (x_1y_1,z_1)^{3r}
$$
and there exists $i$ such that
$$
\frac{\alpha_i'}{i}\le\frac{\alpha_j'}{j},\,\,\frac{\beta_i'}{i}\le \frac{\beta_j'}{j}
$$
for $2\le j\le r$ and
$$
\left\{\frac{\alpha_i'}{i}\right\}+\left\{\frac{\beta_i'}{i}\right\}<1.
$$
Further, there exists a series $\sigma$ such that $\sigma z_1\in {\cal O}_{V_p,q'}$.

Suppose that $D\subset \overline S_r(V_p)$ is a curve (which is necessarily r small). Let $\pi':Z_1\rightarrow V_p$ be the blowup of
$D$.

Suppose that $q'\in D\cap(\overline\pi\circ\pi_1)^{-1}(p)$. $q'$ can only be a point of
the form  of (\ref{eq972}) or (\ref{eq50}).

Suppose that $q'$ satisfies (\ref{eq972}). Then $x=z=0$ are local equations of $D$ at $q'$,
and by Lemma \ref{Lemma5}, $\epsilon =1$,
$\alpha_r=r-1$ and $\alpha_j\ge j$ if $j\ne r-1$.

Suppose that $q''\in (\pi')^{-1}(q')$. Suppose that $\hat{\cal O}_{Z_1,q''}$
 has regular parameters $(x_1,y,z_1)$ where $x=x_1$, $z=x_1(z_1+\alpha)$ for some $\alpha\in k$.
$$
\begin{array}{ll}
u&=x_1^a\\
F_{q''} &= \tau x_1(z_1+\alpha)^r +\sum_{j=2}^{r-1}\overline a_j(x_1,y)x_1^{\alpha_j-j+1}
(z_1+\alpha)^{r-j}+y -g(x_1)+h_1
\end{array}
$$
with $h_1\in (x_1,z_1)^{2r}$
for some series $g(x_1)$. $q''$ 
 is resolved, since $\nu(F_{q'}(0,y,0))=1$.
If $q''$  has permissible parameters $(x_1,y,z_1)$ where $x=x_1z_1$, $z=z_1$, 
\begin{equation}\label{eq1009}
\begin{array}{ll}
u&=x_1^az_1^a\\
F_{q''} &= \tau z_1 +\sum_{j=2}^{r-1}\overline a_j(x_1z_1,y)x_1^{\alpha_j}
z_1^{\alpha_j-j+1}+x_1^{r-1}y+h_1
\end{array}
\end{equation}
with $h_1\in(x_1,z_1)^{2r}$. $\nu(q'')\le r-2$
since $\nu(F_{q'}(0,0,z_1))=1$ and $r\ge 3$. 
By Lemma \ref{Lemma3}, $(\pi')^{-1}(q')\cap \overline S_r(Z_1)=\emptyset$, and the
conclusions of 1. - 3. of the Theorem hold on $(\pi')^{-1}(q')$.

Suppose that $q'\in(\overline\pi\circ\pi_1)^{-1}(p)$ satisfies (\ref{eq50}). Then either $x_1=z_1=0$ or $y_1=z_1=0$ are local equations
of $D$ at $q'$. Without loss of generality, assume that $y_1=z_1=0$ are local
equations of $D$ at $q'$. Then we have $\beta_j'\ge j$ if $2\le j\le r-1$ and 
$\epsilon =1$, $\beta_r'=r-1$ by 
Lemma \ref{Lemma7}. $\nu(q')=r-1$ implies $\alpha_r'=0$.

We have 
\begin{equation}\label{eq51}
\begin{array}{ll}
u&=(x_1^ay_1^b)^m\\
v&=P(x_1^ay_1^b)+x_1^{c}y_1^dF_{q'}\text{ where}\\
F_{q'}&= \tau z_1^r+\sum_{j=2}^{r-1}\overline a_j(x_1,y_1)x_1^{\alpha_j'}y_1^{\beta_j'}z_1^{r-j}
+y_1^{r-1}+h_1
\end{array}
\end{equation}
with $\beta'_j\ge j$ for all $j$,
$$
h_1\in (x_1y_1,z_1)^{3r}
$$

 Suppose that $q''\in (\pi')^{-1}(q')$, and $q''$ has regular parameters $(x_1, y_2, z_2)$ defined by
$$
y_1=y_2, z_1=y_2(z_2+\alpha)
$$
Then
$$
\begin{array}{ll}
u&=(x_1^ay_2^b)^m\\
v&=P_{q''}(x_1^ay_2^b)+x_1^{c}y_2^{d+r-1}F_{q''}\\
\frac{F_{q'}}{y_2^{r-1}} &=
\tau(z_2+\alpha)^ry_2+\sum_{j=2}^{r-1}\overline a_j x_1^{\alpha_j'}y_2^{\beta_j'-j+1}(z_2+\alpha)^{r-j}+1+h_2\\
&= 1 + y_2\Omega 
\end{array}
$$
since
$$
h_2\in (y_2)^{2r}.
$$
$a(d+r-1)-bc\ne 0$ since $F_{q'}$ is normalized.
Thus $F_{q''}= 1+y_2\Omega'$ is a unit.

Suppose that $q'' \in (\pi')^{-1}(q')$, and $q''$ has regular parameters $(x_1, y_2, z_2)$ defined by
$$
y_1=y_2z_2, z_1=z_2
$$
Then 
\begin{equation}\label{eq52}
\begin{array}{ll}
F_{q''}=\frac{F_{q'}}{z_2^{r-1}} &=
\tau z_2+\sum_{j=2}^{r-1}\overline a_j x_1^{\alpha_j'}y_2^{\beta_j'}z_2^{\beta'_j-j+1}+y_2^{r-1}+h_2\\
&= \tau z_2+ + y_2+z_2^2\Omega' 
\end{array}
\end{equation}
since
$$
h_2\in (z_2)^{2r}
$$
Thus $\nu(q'')=1$. 
$\nu(q'')\le r-2$ since $r\ge 3$.
Thus $(\pi')^{-1}(q')\cap \overline S_r(Z_1)=\emptyset$
by Lemma \ref{Lemma4}, and the conclusions 1. - 3. of the
Theorem hold on $(\pi')^{-1}(q')$.

We  thus construct a permissible sequence of monodial transforms $\pi_2:Z_p\rightarrow V_p$
centered at the strict transforms of curves $C\subset \overline S_r(V_p)$
which are r small
 so that $\overline S_r(Z_p)$
contains no curves.
$Z_p\rightarrow \overline Y_p$ extends to $\overline U_2\rightarrow U$ in the notation
of the Theorem.

Suppose that $q'\in(\overline \pi\circ\pi_1\circ\pi_2)^{-1}(p)$ does not satisfy the
conclusions of the Theorem. Then $q$ must either satisfy (\ref{eq972}) or (\ref{eq50}).

We cannot have that (\ref{eq972}) holds at $q'$, since then $\nu(q')=r$, which implies
that $\alpha_j\ge j$ for $j\ge 2$ and $\alpha_{r-1}\ge r-1$, so that $x=z=0$ are local equations of a curve $D$ in $\overline S_r(\overline U_2)$.
We further see that $\overline S_r(\overline U_2)=\emptyset$.

Suppose that (\ref{eq50}) holds at $q'$. Then $\nu(q')=r-1$ and 2. of the conclusions
of the Theorem does not hold, so that $\tau(q')=0$. Thus $\alpha_j'+\beta_j'\ge j$
for $j\ne r$ and $\epsilon =1$, $\alpha_r'+\beta_r'=r-1$.

First suppose that (\ref{eq50}) holds, with $\tau(q')=0$, $\alpha_r'$ and $\beta_r'\ne 0$. 
Let $\pi'':W_1\rightarrow Z_p$ be the quadratic transform with center $q'$. Suppose that $q''\in (\pi'')^{-1}(q')$ and $\hat{\cal O}_{W_1,q''}$
 has regular parameters $(x_2,y_2,z_2)$ such that
$$
x_1=x_2, y_1=x_2(y_2+\alpha), z_1=x_2(z_2+\beta).
$$
with  $\alpha\ne 0$.
Set $x_2 = \overline x_2(y_2+\alpha)^{\frac{-b}{a+b}}$. 
there exists
$$
h_2\in (x_2)^{2r}
$$
such that
$$
\begin{array}{ll}
\frac{F_{q'}}{x_2^{r-1}}&= \tau x_2(z_2+\beta)^r+\sum_{j=2}^{r-1}\overline a_j x_2^{\alpha_j'+\beta_j'-j+1}
(y_2+\alpha)^{\beta_j'}(z_2+\beta)^{r-j}+ (y_2+\alpha)^{\beta'_r}+h_2\\
&=(y_2+\alpha)^{\beta'_r}+\overline x_2\Omega
\end{array}
$$
Thus 
$$
\begin{array}{ll}
u&=(\overline x_2^{a+b})^m\\
v&= P_{q''}(\overline x_2) + \overline x_2^{c+d+r-1}F_{q''}\\
F_{q''}&=(y_2+\alpha)^{\lambda+\beta_r'}
-\alpha^{\lambda+\beta_r'}+\overline x_2\Omega'
\end{array}
$$
where $\lambda= d-\frac{b(c+d+r-1)}{a+b}$. Since $F_{q'}$ is normalized,
$$
\begin{array}{ll}
(a+b)(\beta_r'+d)-b(c+d+r-1)
&=a(d+\beta_r')-b(c+r-1-\beta_r')\\
&=a(d+\beta_r')-b(c+\alpha_r')\ne 0
\end{array}
$$
Thus $\lambda+\beta_r'\ne 0$.
$q''$ is thus a resolved point.

Suppose that $q''\in (\pi'')^{-1}(q')$
 has permissible parameters $(x_2,y_2,z_2)$ such that
$$
x_1=x_2, y_1=x_2y_2, z_1=x_2(z_2+\beta).
$$
There exists
$$
h_2\in (x_2)^{2r}
$$
such that
$$
\begin{array}{ll}
\frac{F_{q'}}{x_2^{r-1}}
&= \tau x_2(z_2+\beta)^r+\sum_{j=2}^{r-1}\overline a_j x_2^{\alpha_j'+\beta_j'-j+1}
y_2^{\beta_j'}(z_2+\beta)^{r-j}+ y_2^{\beta'_r}+h_2\\
&=y_2^{\beta_r'}+ x_2\Omega
\end{array}
$$
$$
b(c+d+r-1)-(a+b)(d+\beta_r')\ne 0
$$
since $F_{q'}$ is normalized.

Thus
$$
\begin{array}{ll}
u&= (x_2^{a+b}y_2^b)^m\\
v&= P_{q''}(x_2^{a+b}y_2^b) + x_2^{c+d+r-1}y_2^dF_{q''}\\
F_{q''}&=y_2^{\beta'_r}
+ x_2\Omega'
\end{array}
$$
$\nu(F_{q''})\le \beta_r'<r-1$.

Suppose that $q''\in (\pi'')^{-1}(q')$
 has regular parameters $(x_2,y_2,z_2)$ such that
$$
x_1=x_2y_2, y_1=y_2, z_1=y_2(z_2+\beta).
$$
Then
$$
\begin{array}{ll}
u&= (x_2^{a}y_2^{a+b})^m\\
v&= P_{q''}(x_2^{a}y_2^{a+b}) + x_2^cy_2^{c+d+r-1}F_{q''}\\
F_{q''}&=x_2^{\alpha'_r}
+y_2\Omega'
\end{array}
$$
since $F_{q'}$ is normalized. $\nu(F_{q''})\le\alpha_r'<r-1$.

The remaining point in  $q''\in (\pi'')^{-1}(q')$
 has regular parameters $(x_2,y_2,z_2)$ such that
$$
x_1=x_2z_2, y_1=y_2z_2, z_1=z_2.
$$
There exists 
$$
h_2\in (z_2)^{2r}
$$
such that 
\begin{equation}\label{eq1010}
\begin{array}{ll}
F_{q''}&= \frac{F_{q'}}{z_1^{r-1}}=
\tau z_2 + \sum_{j=2}^{r-1}\overline a_j x_2^{\alpha'_j}y_2^{\beta'_j}z_2^{\alpha'_j+\beta'_j-j+1}+x_1^{\alpha_r'}y_2^{\beta'_r}+h_2\\
&\equiv \tau z_2\text{ mod }(x_2, y_2, z_2^2)
\end{array}
\end{equation}
Thus $q''$ is a 3 point with $\nu(F_{q''})=1\le r-2$, since $r\ge 3$.

$\nu(q'')<r-1$ for $q''\in (\pi'')^{-1}(q')$, so that 
$(\pi'')^{-1}(q')\cap\overline S_r(Y'')=\emptyset$, and the conclusions of 1.-3. of Theorem \ref{Theorem22} hold on $(\pi'')^{-1}(q')$.

Now suppose that 
(\ref{eq50}) holds, with $\tau(q')=0$ and $\alpha_r'=0$ or $\beta_r'= 0$. Since the 2 cases are symmetric, we may assume that $\alpha_r'=0$.

We thus have $\beta_r'=r-1$ (and $\alpha_j'+\beta_j'\ge j$ for $j\ne m$). Suppose that $i\ne r$. Then $\frac{\alpha_i'}{i}\le
\frac{\alpha_r'}{r}$ implies $\alpha_i'=0$ and $\frac{\beta_i'}{i}\le\frac{r-1}{r}<1$
implies $\beta_i'<i$, so that $\tau(q')>0$, a contradiction. We thus have $i=r$.
$$
\frac{r-1}{r}=\frac{\beta_r'}{r}\le\frac{\beta_j'}{j}
$$
for all $j$ implies $\beta_j'\ge j-\frac{j}{r}$ for $2\le j<r$. Since $\beta_j'\in\bold N$,
we have $\beta_j'\ge j$, and the curve $D$ with local equations $y_1=z_1=0$ at $q'$ is such
that $D\subset \overline S_r(Z_p)$, a contradiction.

\end{pf}

\begin{Theorem}\label{Theorem1011}
Suppose that $r=2$ and $A_2(X)$ holds.
Suppose that $p\in X$ is a 1 point or a 2 point with $\nu(p)=\gamma(p)=2$.
Let $R ={\cal O}_{X,p}$. Suppose that $\pi:Y_p\rightarrow \text{Spec}(\hat R)$
is the sequence of monoidal transforms of sections over the curve $\overline C$
with local equations $\tilde x=y=0$ of Theorem \ref{Theorem19}. For $q\in\pi^{-1}(p)$, define
$$
l_q=\left\{\begin{array}{ll}
([\frac{\alpha_2}{2}]+3)2
& \text{if $F_q$ is a form $(\ref{eq39})$ or $(\ref{eq41})$.}\\
([\frac{\alpha_2}{2}]+
[\frac{\beta_2}{2}]+3)2
& \text{if $F_q$  is a form $(\ref{eq40})$.}
\end{array}\right.
$$
let $l=\text{max}\{l_q\mid q\in \pi^{-1}(p)\}$.

Suppose that $t\ge l$. Let $\overline\pi:\overline Y_p\rightarrow \text{spec}(R)$ be
the sequence of monodial transforms of Theorem \ref{Theorem21}. Let
$$
\cdots\rightarrow Y_n\rightarrow\cdots\rightarrow Y_1\rightarrow \overline Y_p
$$
be a sequence of permissible monoidal transforms centered at curves $C$ in $\overline S_2$
such that $C$ is 2 big. Then there exists
$n_0<\infty$ such that
$$
V_p=Y_{n_0}\stackrel{\pi_1}{\rightarrow}\overline Y_p\rightarrow\text{spec}(R)
$$
extends to a permissible sequence of monoidal transforms 
$$
\overline U_1\rightarrow \overline U\rightarrow U
$$
over an affine neighborhood $U$ of $p$, in the notation of Theorem \ref{Theorem21},
such that $\overline S_2(\overline U_1)$ contains no curves $C$ such that 
$C$ is 2 big. Let
$$
\cdots\rightarrow Z_n\rightarrow\cdots\rightarrow V_p
$$
be a permissible sequence of monoidal transforms centered at  curves $C$ in $\overline S_2$
such that $C$ is 2 small. Then there exists $n_1<\infty$
such that $\pi_2:Z_p=Z_{n_1}\rightarrow V_p$ extends to a permissible sequence of
monoidal transforms
$$
\overline U_2\rightarrow \overline U_1\rightarrow \overline U\rightarrow U
$$
over an affine neighborhood $U$ of $p$ such that $\overline S_2(\overline U_2)=\emptyset$.

Finally, there exists a sequence of quadratic transforms and monodial transforms 
centered at strict transforms of 2 curves $C$ on $Z_p$ 
such that $C$ is 1 big and $C$ is a section over a 2 small curve blown up in $Z_p\rightarrow V_p$,   $\pi_3:W_p\rightarrow Z_{p}$
which extends to a permissible sequence of monoidal transforms
$$
\overline U_3\rightarrow\overline U_2\rightarrow \overline U_1\rightarrow \overline U\rightarrow U
$$
over an affine neighborhood $U$ of $p$ such that $\overline S_2(\overline U_3)=\emptyset$,
and if 
$$
q\in(\overline\pi\circ\pi_1\circ\pi_2\circ\pi_3)^{-1}(p)
$$
 then $q$ is resolved.
\end{Theorem}

\begin{pf} The analysis of Theorem \ref{Theorem22} is valid for $r=2$, except in 
(\ref{eq1009}), (\ref{eq52}) and (\ref{eq1010}).

The situation of (\ref{eq1010}) cannot occur when $r=2$, since this comes from the case when (\ref{eq50}) holds, with $\tau(q')=0$, $\alpha_r'$ and $\beta_r'\ne 0$. Since
$\alpha_r'+\beta_r'=r-1=1$, this case cannot occur.

Suppose that a case (\ref{eq1009}) occurs in 
$$
(\overline \pi\circ \pi_1\circ \pi_2):Z_p\rightarrow \text{spec}(R).
$$
 Then we have a 2 point $q''\in (\overline\pi\circ\pi_1\circ\pi_2)^{-1}(p)$ such that 
$$
\begin{array}{ll}
u&=(x_1^az_1^b)^m\\
v&=P_{q''}(x_1^az_1^b)+x_1^cz_1^dF_{q''}\\
F_{q''}&=z_1+x_1y_1+h_1
\end{array}
$$
with $h_1\in(x_1,z_1)^4$.

Let $C$ be the 2 curve on $Z_p$ with local equations $x_1=z_1=0$. $C$ is 1 big by Lemma \ref{Lemma655}.

Let $\pi':W_1\rightarrow Z_p$ be the blowup of $C$. Suppose that $\overline q\in(\pi')^{-1}(q'')$ is a 
1 point. Then there exist regular parameters $(x_2,y_1,z_2)$ in $\hat{\cal O}_{W_1,\overline q}$  such that
$$
x_1=x_2, z_1=x_2(z_2+\alpha)
$$
with $\alpha\ne 0$. Set
$$
x_2=\overline x_2(z_1+\alpha)^{-\frac{b}{a+b}}.
$$
$$
\begin{array}{ll}
u&=\overline x_2^{(a+b)m}\\
v&=P_{q''}(\overline x_2^{a+b})+\overline x_2^{c+d+1}(z_2+\alpha)^{d-\frac{b(c+d+1)}{a+b}}
(z_2+\alpha+y_1+x_1\Omega)
\end{array}
$$
Thus $\nu(F_{\overline q}(0,y_1,z_2))=1$ and $\overline q$ is resolved.

Suppose that $\overline q\in(\pi')^{-1}(q'')$ is the 2 point with permissible parameters
$$
x_1=x_2, z_1=x_2z_2.
$$
Then $F_{\overline q}=\frac{F_{q''}}{x_2}=z_2+y_1+x_1\Omega$ and $\overline q$ is resolved.

If $\overline q\in (\pi')^{-1}(q'')$ is the 2 point with permissible parameters
$$
x_1=x_2z_2, z_1=z_2
$$
then
$$
F_{\overline q}=\frac{F_{q''}}{z_2}=1+x_2y_1+z_2\Omega
$$
and is resolved.

Suppose that a case (\ref{eq52}) occurs in 
$$
(\overline \pi\circ\pi_1\circ\pi_2):Z_p\rightarrow\text{spec}(R).
$$
Then we have a 3 point $q''\in(\overline\pi\circ\pi_1\circ\pi_2)^{-1}(p)$ such that
$$
\begin{array}{ll}
u&=(x_2^ay_2^bz_2^c)^m\\
v&=P_{q''}(x_2^ay_2^bz_2^c)+x_2^dy_2^ez_2^fF_{q''}\\
F_{q''}&=z_2+y_2+h_1
\end{array}
$$
with $h_1\in (z_2)^4$.

Let $C$ be the 2 curve on $Z_p$ with local equations $y_2=z_2=0$. $C$ is 1 big by Lemma
\ref{Lemma655}. Let $\pi':W_1\rightarrow Z_p$ be the blowup of $C$. Suppose that
$\overline q\in(\pi')^{-1}(q'')$ is a 2 point. Then there exist regular parameters
$(x_2,y_3,z_3)$ in $\hat{\cal O}_{W_1,\overline q}$ such that
$$
y_2=y_3, z_2=y_3(z_3+\alpha)
$$
with $\alpha\ne 0$. Set $y_3=\overline y_3(z_3+\alpha)^{-\frac{c}{b+c}}$.
$$
\begin{array}{ll}
u&=(x_2^a\overline y_3^{b+c})^m=(x_2^{\overline a}\overline y_3^{\overline b})^{\overline m}\\
v&=P_{q''}(x_2^a\overline y_3^{b+c})+x_2^d\overline y_3^{e+f+1}(z_3+\alpha)^{f-\frac{c(e+f+1)}{b+c}}(z_3+\alpha+1+y_3^3\Omega)
\end{array}
$$
with $(\overline a,\overline b)=1$. Thus $\nu(F_{\overline q}(0,0,z_3))=1$ and $\overline q$ is resolved.

Suppose that $\overline q\in (\pi')^{-1}(q'')$ is the 3 point with permissible
parameters $y_2=y_3$, $z_2=y_3z_3$. Then
$$
F_{\overline q}=\frac{F_{q''}}{y_3}=z_3+1+y_3^3\Omega
$$
and $\overline q$ is resolved.

Suppose that $\overline q\in(\pi')^{-1}(q'')$ is the 3 point with permissible
parameters $y_2=y_3z_3$, $z_2=z_3$. Then
$$
F_{\overline q}=\frac{F_{q''}}{z_3}=1+y_3+z_3^3\Omega
$$
and $\overline q$ is resolved.

\end{pf}

\begin{Theorem}\label{Theorem51}
Suppose that $r\ge 3$ and $A_r(X)$ holds.
Suppose that  $p\in X$ is a 2 point such that
$\nu(p)=r-1$, $\tau(p)=0$  and  $\gamma(p)=r$.
Let $\overline C$ be the 2 curve containing $p$. Let $R = {\cal O}_{X,p}$.

There exists a sequence of permissible monoidal transforms centered at sections over $\overline C$,
$\overline Y_p\rightarrow \text{spec}(R)$, which extends to a sequence of permissible
monoidal transforms $\overline U\rightarrow U$ where $U$ is an affine neighborhood of $p$,
with the following property.

 Let
$$
\cdots\rightarrow Y_n\rightarrow\cdots\rightarrow Y_1\rightarrow \overline Y_p
$$
be a sequence of permissible monoidal transforms centered at curves $C$ in $\overline S_r$
such that $C$ is r big. Then there exists
$n_0<\infty$ such that
$$
V_p=Y_{n_0}\stackrel{\pi_1}{\rightarrow}\overline Y\rightarrow\text{spec}(R)
$$
extends to a permissible sequence of monoidal transforms 
$$
\overline U_1\rightarrow \overline U\rightarrow U
$$
over an affine neighborhood $U$ of $p$, in the notation of Theorem \ref{Theorem21},
such that $\overline S_r(\overline U_1)$ contains no curves $C$ such that $C$ is r big. Let
$$
\cdots\rightarrow Z_n\rightarrow\cdots\rightarrow V_p
$$
be a permissible sequence of monoidal transforms centered at curves $C$ in $\overline S_r$
such that $C$ is r small. Then there exists $n_1<\infty$
such that $\pi_2:Z_p=Z_{n_0}\rightarrow V_p$ extends to a permissible sequence of
monoidal transforms
$$
\overline U_2\rightarrow \overline U_1\rightarrow \overline U\rightarrow U
$$
over an affine neighborhood $U$ of $p$ such that $\overline S_r(\overline U_2)=\emptyset$.

Finally, there exists a sequence of quadratic transforms $\pi_3:W_p\rightarrow Z_{p}$
which extends to a permissible sequence of monoidal transforms
$$
\overline U_3\rightarrow\overline U_2\rightarrow \overline U_1\rightarrow \overline U\rightarrow U
$$
over an affine neighborhood $U$ of $p$ such that $\overline S_r(\overline U_3)=\emptyset$,
and if 
$$
q\in(\overline\pi\circ\pi_1\circ\pi_2\circ\pi_3)^{-1}(p),
$$
\begin{enumerate}
\item $\nu(q)\le r-1$ if $q$ is a 1 or 2 point.
\item If $q$ is a 2 point and $\nu(q)=r-1$, then $\tau(q)>0$.
\item $\nu(q)\le r-2$ if $q$ is a 3 point.
\end{enumerate}
\end{Theorem}

\begin{Theorem}\label{Theorem1012}
Suppose that $r=2$ and $A_2(X)$ holds
Suppose that $p\in X$ is a  2 point such that $\nu(p)=1$, $\tau(p)=0$ and $\gamma(p)=2$.
Let $R ={\cal O}_{X,p}$. 

There exists a sequence of permissible monoidal transforms centered at sections over $\overline C$, $\overline Y_p\rightarrow \text{spec}(R)$ which extends to a sequence of permissible
monoidal transforms $\overline U\rightarrow U$, where $U$ is an affine neighborhood of
$p$, with the following property.

 Let
$$
\cdots\rightarrow Y_n\rightarrow\cdots\rightarrow Y_1\rightarrow \overline Y_p
$$
be a sequence of permissible monoidal transforms centered at curves $C$ in $\overline S_2$
such that $C$ is 2 big. Then there exists
$n_0<\infty$ such that
$$
V_p=Y_{n_0}\stackrel{\pi_1}{\rightarrow}\overline Y_p\rightarrow\text{spec}(R)
$$
extends to a permissible sequence of monoidal transforms 
$$
\overline U_1\rightarrow \overline U\rightarrow U
$$
over an affine neighborhood $U$ of $p$, in the notation of Theorem \ref{Theorem21},
such that $\overline S_2(\overline U_1)$ contains no curves $C$ such that $C$ is 2 big. Let
$$
\cdots\rightarrow Z_n\rightarrow\cdots\rightarrow V_p
$$
be a permissible sequence of monoidal transforms centered at  curves $C$ in $\overline S_2$
such that $C$ is 2 small. Then there exists $n_1<\infty$
such that $\pi_2:Z_p=Z_{n_1}\rightarrow Y_{n_0}$ extends to a permissible sequence of
monoidal transforms
$$
\overline U_2\rightarrow \overline U_1\rightarrow \overline U\rightarrow U
$$
over an affine neighborhood $U$ of $p$ such that $\overline S_2(\overline U_2)=\emptyset$.

Finally, there exists a sequence of quadratic transforms and monodial transforms 
centered at strict transforms of 2 curves $C$ on $Z_p$  such that 
$C$ is 1 big and $C$ is a section over a 2 small curve blown up in $Z_p\rightarrow V_p$,
$\pi_3:W_p\rightarrow V_p$
which extends to a permissible sequence of monoidal transforms
$$
\overline U_3\rightarrow\overline U_2\rightarrow \overline U_1\rightarrow \overline U\rightarrow U
$$
over an affine neighborhood $U$ of $p$ such that $\overline S_2(\overline U_3)=\emptyset$,
and if 
$$
q\in(\overline\pi\circ\pi_1\circ\pi_2\circ\pi_3)^{-1}(p)
$$
 then $q$ is resolved.
\end{Theorem}

\begin{pf}(of Theorems \ref{Theorem51} and \ref{Theorem1012}) The conclusions of Lemma \ref{Lemma17} hold at $p$.

The conclusions of Lemma \ref{Lemma18} must be modified to:
$\alpha_j+\beta_j\ge j$ for $2\le j\le r-2$, $\alpha_r+\beta_r\ge r-1$.

In the conclusions of Theorem \ref{Theorem19}, we must add a fourth case: 
\begin{equation}\label{eq638}
\begin{array}{ll}
u&=(x^ay^b)^m\\
v&=P(x^ay^b)+x^cy^dF_q\\
F_q&=\tau z^r+\sum_{i=2}^{r-1}\overline a_i(x,y)x^{\alpha_i}y^{\beta_i}z^{r-i}
+y^{r-1}
\end{array}
\end{equation}
with $\beta_i\ge i$ for all $i$.

Theorem \ref{Theorem21} must be modified by adding the case
$F_q$ equivalent $\text{ mod }(\overline x\overline y,z)^t$ to a form (\ref{eq638}).
 The proof of Theorem \ref{Theorem22}
must be modified by adding an analysis of (\ref{eq638}).
Such $q$ are not effected by blowing up r big curves, so the construction of
$\pi_1:V_p\rightarrow \overline Y_p$ is as in the proof of Theorem \ref{Theorem22}.
Suppose that $q\in (\overline\pi\circ\pi_1)^{-1}(p)$ satisfies (\ref{eq638}).
 There is a unique r small curve $D \subset \overline S_r(V_p)$
containing $q$, which has local equations $y=z=0$. Let $\pi':Z_1\rightarrow \overline Y_p$ be the blowup of $D$. Then if $r\ge 3$, all points of $(\pi')^{-1}(q)$ satisfy 1. - 3. of the 
conclusions of the Theorem, and $\overline S_r(Z_1)\cap (\pi')^{-1}(q)=\emptyset$.

If $r=2$, and $q'\in(\pi')^{-1}(q)$ is the 3 point, then there exist permissible
parameters $(x,y_1,z_1)$ at $q'$ such that 
$$
F_{q'}=z_1+y_1+z_1^2\Omega.
$$
If $C$ is the 2 curve with local
equations $y_1=z_1=0$, then $C$ is 1 big, and if $\pi':W_1\rightarrow Z_p$ is the blowup of $C$, then all points
of $(\pi')^{-1}(q')$ are resolved.

\end{pf}

\section{reduction of $\nu$ in a second special case}\label{Spec2}

Throughout this section, we will assume that $\Phi_X:X\rightarrow S$ is weakly prepared.

\begin{Theorem}\label{Theorem27}(Theorem27)
Suppose that $r\ge 2$, $A_r(X)$ holds,
 $p\in X$ is a 2 point with $\nu(p)=r-1$, 
$C$ is a generic curve through $p$,
and there are permissible parameters $(x,y,z)$ at $p$ for $(u,v)$ (with $y,z\in{\cal O}_{X,p}$)
such that $L_p(x,0,0)\ne 0$,  and $C$ has local equations $y=z=0$ at $p$.
Let $R ={\cal O}_{X,p}$. 
We have an expression at $p$
\begin{equation}\label{eq658}
\begin{array}{ll}
u&=(x^ay^b)^m\\
v&=P(x^ay^b)+x^cy^dF_p\\
F_p&=\tau x^{r-1}+\sum_{i=1}^{r-1}\tilde a_i(y,z)x^{r-i-1}
\end{array}
\end{equation}
where $\tau$ is a unit and  $\nu(\tilde a_i)\ge i$ for all $i$.
Then
there exists a finite sequence of permissible monodial transforms $\pi:Y\rightarrow \text{Spec}(R)$
centered at sections over $C$, such that for $q\in \pi^{-1}(p)$, 
there exist permissible parameters $(\overline x,\overline y,\overline z)$ at $q$ such that
$F_q$ has one of the following forms.

\begin{enumerate}
\item
\begin{equation}\label{eq62}
\begin{array}{ll}
u&=(\overline x^a\overline y^b)^m\\
v&=P(\overline x^a \overline y^b)+\overline x^c\overline y^dF_q\text{ with }\\
F_q &= 
\tau \overline x^{r-1}+\sum_{j=1}^{r-1}
\overline y^{d_j}\Lambda_j(\overline y,\overline z)\overline x^{r-1-j}
\end{array}
\end{equation}
where $\tau$ is a unit, $\Lambda_j(\overline y,\overline z)=0$ or $e_j=\nu(\Lambda_j(0,\overline z))=0$ or $1$
and  $d_j+e_j\ge j$ for all $j$.
\item
\begin{equation}\label{eq63}
\begin{array}{ll}
u&=(\overline x^a\overline y^b\overline z^c)^m\\
v&=P(\overline x^a\overline y^b\overline z^c)+\overline x^d\overline y^e\overline z^fF_q\text{ with }\\
F_q &= 
\tau \overline x^{r-1}+\sum_{j=1}^{r-1}\overline a_j(\overline y,\overline z)\overline y^{d_j}\overline z^{e_j}\overline x^{r-1-j}
\end{array}
\end{equation}
where $\tau$ is a unit, $d_j+e_j\ge j$, $\overline a_j$ are units (or zero) for all $j$ and
there exists an $i$ such that $1\le i\le r-1$, $\overline a_i\ne 0$ and 
$$
\frac{d_i}{i}\le  \frac{d_j}{j},    \frac{e_i}{i} \le \frac{e_j}{j}
$$
for $1\le j\le r-1$. We further have
$$
\left\{\frac{d_i}{i}\right\} +  \left\{\frac{e_i}{i}\right\}<1.
$$
\item
\begin{equation}\label{eq64}
\begin{array}{ll}
u&=(\overline x^a\overline y^b)^m\\
v&=P(\overline x^a\overline y^b)+\overline x^c\overline y^dF_q\text{ with }\\
F_q &= 
\tau \overline x^{r-1}+\sum_{j=1}^{r-1}\overline a_j(\overline y,\overline z)
\overline y^{d_j}\overline z^{e_j}\overline x^{r-1-j}
\end{array}
\end{equation}
where $\tau$ is a unit, $d_j+e_j\ge j$, $\overline a_j$ are units (or zero)  for all $j$ and
there exists an $i$ such that $1\le i\le r-1$, $\overline a_i\ne 0$ and 
$$
\frac{d_i}{i}\le  \frac{d_j}{j},\frac{e_i}{i} \le \frac{e_j}{j}
$$
for $1\le j\le r-1$. We further have
$$
\left\{\frac{d_i}{i}\right\} +  \left\{\frac{e_i}{i}\right\}<1.
$$
\end{enumerate}

In all these cases $\overline x=x$ and $\overline x=0$ is a local equation at $q$ of the strict transform
of the component of $E_X$ with local equation $x=0$ at $p$.

There exists an affine neighborhood $U$ of $p$ such that $Y\rightarrow\text{spec}(R)$
extends to a sequence of permissible monoidal transforms $\overline U\rightarrow U$ such that
$A_r(\overline U)$ holds.

Suppose that $r\ge 3$. Let $D_i$ be the curves in $\overline S_{r-1}(X)$ which contain $p$, and such that $x\in \hat{\cal I}_{D_i,p}$. We further have that the  strict transforms $\overline D_i$ of the $D_i$ on $\overline U$ are nonsingular, disjoint and make SNCs with $\overline B_2(\overline U)$.

\end{Theorem}

\begin{pf} Set $S_0=\text{Spec}(k[y,z])$, $Y_0=\text{Spec}(R)$.
 Consider the sequence of monoidal transforms centered at sections over $C$ 
\begin{equation}\label{eq507}
\cdots\rightarrow Y_n\rightarrow Y_{n-1}\rightarrow Y_{n-2}\rightarrow\cdots\rightarrow Y_1\rightarrow Y_0
\end{equation}
where the sequence is obtained from a sequence of quadratic transforms 
\begin{equation}\label{eq508}
\cdots\rightarrow S_n\rightarrow S_{n-1}\rightarrow\cdots\rightarrow S_1\rightarrow S_0
\end{equation}
over the closed point $p_0$ with local equations $y=z=0$ in $S_0$, and (\ref{eq507}) is obtained from (\ref{eq508}) by base change with $Y_0\rightarrow S_0$, so that $Y_i\cong S_i\times_{S_0}Y_0$. 

The  map $\hat{S_0}\rightarrow \hat{Y_0}$ obtained from the natural projection
$$
k[[x,y,z]]\rightarrow k[[y,z]]
$$
induces maps 
$S_i\times_{S_0}\hat S_0\rightarrow Y_i\times_{Y_0}\hat{Y_0}$
such that the composed map
$$
S_i\times_{S_0}\hat S_0\rightarrow Y_i\times_{Y_0}\hat{Y_0}
\rightarrow S_i\times_{S_0}\hat S_0
$$
is an isomorphism for all $i$.

We can thus identify the center of the quadratic transform
 $S_{i+1}\rightarrow S_{i}$
 with a point $p_i\in Y_i$ over $p$. A section over $C$ through $p_i$ is blown up 
in (\ref{eq507}) only if none of the forms (\ref{eq62}), (\ref{eq63})
 or (\ref{eq64}) hold at $p_i$.

We will show that (\ref{eq507}) is finite, so that there exists $n$ such that $Y_n$ satisfies the conclusions of the theorem.

Suppose that (\ref{eq507}) is not finite. Then we may assume that there exists an infinite sequence of points $p_0, p_1,p_2,\ldots,p_n,\ldots$
such that  $Y_{i+1}\rightarrow Y_i$
is a permissible monoidal transform, centered at a section $C_i$ over $C$,
containing $p_i$, such that $p_i$ maps to $p_{i-1}$ for all $i$,
and $F_{p_i}$ does not satisfy (\ref{eq62}), (\ref{eq63}) or (\ref{eq64}) for any $i$.

Each point $p_i$ has permissible parameters $(x,y_i,z_i)$ for $(u,v)$, such that
one of the following cases hold.
\vskip .2truein
{\bf Case 1} $p_i$ is a 2 point
$$
u=(x^{a_i}y_i^{b_i})^{m_i}, v=P_i(x^{a_i}y_i^{b_i})+x^cy_i^{d_i}F_i
$$
with $a_im_i=am$, and permissible parameters at $p_{i+1}$ are as in one of the following
cases.

\noindent {\bf Case 1a}
$$
y_i=y_{i+1}, z_i=y_{i+1}(z_{i+1}+\alpha_{i+1})
$$
{\bf Case 1b} 
$$
y_i=y_{i+1}z_{i+1},
z_i=z_{i+1}
$$

{\bf Case 2} $p_i$ is a 3 point
$$
\begin{array}{ll}
u&=(x^{a_i}\omega_i^{k_i})^{m_i}, \omega_i=y_i^{\overline b_i}z_i^{\overline c_i}\\
v&=P_i(x^{a_i}\omega_i^{k_i})+x^c\omega_i^{ d_i}F_i
\end{array}
$$
with $a_im_i=am$, $(\overline b_i,\overline c_i)=1$, $(a_i,k_i)=1$, and permissible
parameters at $p_{i+1}$ are as in one of the following cases.

\noindent {\bf Case 2a}
$$
y_i=y_{i+1}, z_i=y_{i+1}z_{i+1}
$$
{\bf Case 2b} 
$$
y_i=y_{i+1}z_{i+1}, z_i=z_{i+1}
$$
{\bf Case 2c}
$$
y_i=y_{i+1}(z_{i+1}+\alpha_{i+1})^{-\frac{\overline c_i}
{\overline b_i+\overline c_i}},
z_i=y_{i+1}(z_{i+1}+\alpha_{i+1})^{\frac{\overline b_i}{\overline b_i+\overline c_i}}
$$
with $\alpha_{i+1}\ne 0$. In Case 2c, $y_{i+1}, z_{i+1}$ are constructed from the monoidal transform
$$
y_{i}=\overline y_{i+1}, z_i=\overline y_{i+1}(\overline z_{i+1}+\alpha_{i+1}).
$$
Then define
$$
\overline y_{i+1}=y_{i+1}(z_{i+1}+\alpha_{i+1})^{-\frac{\overline c_i}
{\overline b_i+\overline c_i}},
\overline z_{i+1}=z_{i+1}.
$$

If $p_i$ is a 2 point, then $y$ is a power of $y_i$, and if
$q_i$ is a 3 point, then $y$ is a monomial in $y_i$ and $z_i$. 
If $p_i$ is a 2 point, then there is a series $g_i$ such that
$$
F_i=F_p-\frac{g_i(x^{a_i}y_i^{b_i})}{x^cy_i^{d_i}}.
$$
If $p_i$ is a 3 point, then there is a series $g_i$ such that 
$$
F_i=F_p-\frac{g_i(x^{a_i}\omega_i^{k_i})}{x^c\omega_i^{d_i}}.
$$
In either case, we have an expression
 $$
F_i=F_{p_i}=\tau'x^{r-1}+\sum_{j=1}^{r-1}a_j'(y_j,z_j)x^{r-j-1}.
$$

We will show that $\tau'$ is a unit. Suppose not. First suppose that $p_i$ is a 
2 point. Then $a_id_i-b_i(c+r-1)=0$.
$$
(x^ay^b)^m=(x^{a_i}y^{b_i})^{m_i}
$$
 implies $y=y_i^{\frac{b_im_i}{bm}}$.
$x^{c+r-1}y^d=x^{c+r-1}y_i^{d_i}$ implies 
$$
d_i=\frac{b_im_id}{mb}.
$$
Thus $ad-b(c+r-1)=0$, a contradiction to the assumption that $F_p$ is normalized.

Now suppose that $p_i$ is a 3 point and $\tau'$ is not a unit. Then $a_id_i-(c+r-1)k_i=0$.
$$
(x^ay^b)^m=(x^{a_i}\omega_i^{k_i})^{m_i}
$$
 implies $y=\omega_i^{\frac{k_im_i}{mb}}$.
$x^{c+r-1}y^d=x^{c+r-1}\omega_i^{d_i}$ implies 
$$
d_i=\frac{dk_im_i}{mb}.
$$
 Thus
$ad-(c+r-1)b=0$, a contradiction to the assumption that $F_p$ is normalized.

If $p_i$ is a 2 point,
$$
a_j'=\left\{\begin{array}{ll}
\tilde a_j-\tilde cy_i^{t_j(i)}&\text{ if }a_i( d_i+t_j(i))-b_i(c+r-j-1)=0
\text{ with }t_j(i)\in\bold N\\
&\text{ and }\tilde c\text{ is the coefficient of }y_i^{t_j(i)}\text{ in the
expansion of $\tilde a_j$ in terms of $y_i, z_i$}\\
\tilde a_j&\text{ if }a_i(d_i+t)-b_i(c+r-j-1)\ne 0\text{ for any }t\in \bold N.
\end{array}\right.
$$
If $p_i$ is a 3 point,
$$
a_j'=\left\{\begin{array}{ll}
\tilde a_j-\tilde c\omega_i^{t_j(i)}&\text{ if }a_i( d_i+t_j(i))-k_i(c+r-j-1)=0
\text{ with }t_j(i)\in\bold N\\
&\text{ and }\tilde c\text{ is the coefficient of }\omega_i^{t_j(i)}\text{ in the
expansion of } \tilde a_j \text{ in terms of }y_i, z_i\\
\tilde a_j&\text{ if }a_i( d_i+t)-k_i(c+r-j-1)\ne 0\text{ for any }t\in \bold N.
\end{array}\right.
$$

In particular, there exists at most one value of $t_j(i)$ such that a term can be removed from
 any $\tilde a_j$. 

Set $\overline u=y^b$, $\overline v_j=y^d\tilde a_j(y,z)$. We have
$$
\begin{array}{ll}
u&=(x^a\overline u)^m\\
v&=P(x^a\overline u)+\tau x^{c+r-1}y^d+
\sum_{j=1}^{r-1}\overline v_jx^{c+r-1-j}.
\end{array}
$$
By Theorem \ref{Theorem965}, Lemma \ref{Lemma1027} and Theorem \ref{T1}, there exists $i_0$ such that for $i\ge i_0$ in (\ref{eq508}),
one of the following forms holds at $p_i$, for $1\le j\le r-1$.

If $p_i$ is a 2 point,  
\begin{equation}\label{eq644}
\begin{array}{ll}
\overline u&=y_i^{b_im_i}\\
\overline v_j&=P_{ji}(y_i)+y_i^{d_j(i)}\psi_{ji}(y_i,z_i)^{e_j(i)}
\end{array}
\end{equation}
where $e_j(i)\ge 1$, $\nu(\psi_{ji}(0,z_i))=1$ or $\psi_{ji}=0$ for $1\le j\le r-1$.

If $p_i$ is a 3 point, 
\begin{equation}\label{eq645}
\begin{array}{ll}
\overline u&=(\omega_i^{k_i})^{m_i}=(y_i^{\overline b_i}z_i^{\overline c_i})^{k_im_i}\\
\overline v_j&=P_{ji}(\omega_i)
+y_i^{d_j(i)}z_i^{e_j(i)}\phi_{ji}(y_i,z_i)
\end{array}
\end{equation}
where $\phi_{ji}(y_i,z_i)$ is a unit and
$d_j(i)\overline c_i-\overline b_ie_j(i)\ne 0$,
or $\phi_{ji}=0$ for $1\le j\le r-1$. For $i$ sufficiently large,  we have that $\overline u\,\overline v_j=0$ are
SNC divisors for $1\le j\le r-1$ (by Lemma \ref{Lemma1024}). 

Suppose that (\ref{eq644}) holds at $p_i$ with $e_j(i)=1$ or $\phi_{ji}=0$
for all $j$. If there exists $t_j(i)\in\bold N$ such that 
$$
a_i(d_i+t_j(i))
-b_i(c+r-j-1)=0,
$$
 then the normalized form of
$x^{c+r-1-j}\overline v_j$ at $p_i$ is
$$
x^{c+r-1-j}\overline v_j-\tilde c x^{c+r-1-j}y_i^{t_j(i)+d_i}=x^{c+r-1-j}y_i^{\lambda_j(i)}
\Lambda_{ji}(y_i,z_i)
$$
where $\Lambda_{ji}$ is a unit, zero, or $\nu(\Lambda_{ji}(0,z_i))=1$. 

If there does not exist $t\in\bold N$ such that
$$
a_i(d_i+t)-b_i(c+r-j-1)=0,
$$
 then
$$
x^{c+r-1-j}\overline v_j=x^{c+r-1-j}y_i^{\lambda_j(i)}\Lambda_{ji}(y_i,z_i)
$$
where $\Lambda_{ji}$ is a unit, 0 or $\nu(\Lambda_{ji}(0,z_i))=1$.

Thus if (\ref{eq644}) holds at $p_i$, with $e_j(i)=1$ or $\phi_{ji}=0$ for all $j$,
(\ref{eq62}) holds at $p_i$.

If $p_i$ is a 3 point, so that (\ref{eq645}) holds, and $p_{i+1}$ is a 2 point,
then (\ref{eq644}) holds at $p_{i+1}$, with $e_j(i)=1$ or $\psi_{ji}=0$ for
$1\le j\le r-1$, so that (\ref{eq62}) holds at $p_{i+1}$.

We are reduced to the 2 cases where either for all $i\ge i_0$ in (\ref{eq508}) all
monoidal transforms are of the forms 2a or 2b, or for all $i\ge i_0$ in
(\ref{eq508}), all monoidal transforms are of the form 1a with some $e_j(i)>1$.

If all monoidal transforms are of the form 2a or 2b for $i\ge i_0$, then $F_i=F_{i_0}$ for 
all $i\ge i_0$. If
$$
F_{i_0}=\tau_0x^{r-1}+\sum_{j=1}^{r-1}\tilde a_j(y_{i_0},z_{i_0})x^{r-1-j}
$$
then for $i>>i_0$, $\tilde a_j(y_{i_0},z_{i_0})$ is a monomial in $y_i$ and $z_i$
times a unit for all $j$ by Lemma \ref{Lemma1024}. By Lemmas
\ref{Lemma23} and Corollary \ref{Corollary24}, there exists $i_1>i_0$ such that
(\ref{eq63}) holds at $p_{i_1}$.

Suppose that all monoidal transforms are of the form 1a for $i\ge i_0$ and some
$e_j(i)>1$. There exists a permissible change of parameters $(x,y_{i_0},\overline z_{i_0})$
such that
$$
\overline z_{i_0}=z_{i_0}-\tilde p(y_{i_0})
$$
for some series $\tilde p$, such that $(x, y_i,\overline z_i)$ are permissible parameters
at $p_i$ for all $i\ge i_0$ with $y_i=y_{i_0}$, $\overline z_{i_0}=y_i^{i-i_0}\overline z_i$.
Let 
$$
F=\overline P(y_i)+\overline F_i
$$
be the normalized form of $F$ with respect to the parameters $(x,y_i,\overline z_i)$.
Then $\overline F_i=\overline F_{i_0}$ for $i\ge i_0$. If
$$
\overline F_{i_0}=\tau_0x^{r-1}+\sum_{j=1}^{r-1}\tilde a_j(y_{i_0},\overline z_{i_0})x^{r-1-j}
$$
then for $i>>i_0$, $\tilde a_j(y_{i_0},\overline z_{i_0})$ is a monomial in
$y_i$ and $\overline z_i$ times a unit for all $j$ by Lemma \ref{Lemma1024}.
By Lemma \ref{Lemma23} and Corollary \ref{Corollary24} there exists $i_1>i_0$ such
that (\ref{eq64}) holds at $p_{i_1}$.

$C$ generic implies $F_q$ is resolved for $q\in C$ a generic point. There exists an
affine neighborhood $U$ of $p$ such that $Y\rightarrow \text{spec}(R)$ extends to a
permissible sequence of monoidal transforms of sections over $C$, $\overline U\rightarrow
U$ such that $\overline S_r(\overline U)\cup \overline B_2(\overline U)$ is 
contained in the
union of $\overline B_2(\overline U)$ and the strict transform of $\overline S_r(U)$.
Thus $A_r(\overline U)$ holds.

If $r\ge 3$, we can choose $i_0$ sufficiently large in obtaining the forms of
(\ref{eq644}) and (\ref{eq645}) so that the strict transforms of the curves $D_i$ in 
$\overline S_{r-1}(U)$ such that $x\in \hat{\cal I}_{D_i,p}$ are disjoint and
make SNCs with $\overline B_2(\overline U)$. 

\end{pf}

\begin{Theorem}\label{Theorem28}
Suppose that $r\ge 3$, $A_r(X)$ holds, $p$ is a 2 point with $\nu(p)=r-1$,  
and $L(x,0,0)\ne 0$, as in the assumptions of Theorem \ref{Theorem27}.
Let $R ={\cal O}_{X,p}$. Suppose that  $\pi:Y_p\rightarrow \text{Spec}(R)$ is the morphism of
Theorem \ref{Theorem27}. 

Let 
$$
\cdots\rightarrow Y_n\rightarrow\cdots\rightarrow Y_1\rightarrow Y_p
$$
be a sequence of permissible monodial transforms centered at 2 curves $D$
such that $D$ is r-1 big.
Then there exists $n_0<\infty$ such that 
$$
V_p=Y_{n_0}\stackrel{\pi_1}{\rightarrow} Y_p\rightarrow\text{spec}(R)
$$
extends to a permissible sequence of monodial transforms 
$$
\overline U_1\rightarrow\overline U\rightarrow U
$$
over an affine neighborhood $U$ of $p$ (with the notation of Theorem \ref{Theorem27})
such that $\overline U_1$ contains no 2 curves $D$ such that
 $D$ is r-1 big or r small, and for $q\in\overline U_1$,
\begin{enumerate}
\item If $q$ is a 1 or a 2 point then $\nu(q)\le r$. $\nu(q)=r$ implies $\gamma(q)=r$.
\item If $q$ is a 3 point then $\nu(q)\le r-2$.
\item $\overline S_r(\overline U_1)$ makes SNCs with $\overline B_2(\overline U_1)$.
\end{enumerate}
There exists a sequence of quadratic transforms $W_p\rightarrow V_p$ such that if
$Z_p\rightarrow W_p$ is the sequence of monodial transforms (in any order) centered
 at the strict transforms of
curves $C$ in $\overline S_r(X)$
then 
$$
Z_p\rightarrow W_p\rightarrow V_p\rightarrow Y_p\rightarrow\text{spec}(R)
$$
extends to a permissible sequence of monodial transforms
$$
\overline\pi:\overline U_2\rightarrow \overline U_1\rightarrow \overline U\rightarrow U
$$
over an affine neighborhood of $p$ such that 
$\overline U_2$ contains no 2 curves $D$ such that 
$D$ is r-1 big or r small. $\overline S_r(\overline U_2)$
makes SNCs with $\overline B_2(\overline U_2)$, and 
if $q\in \overline \pi^{-1}(p)$,
\begin{description}
\item[1'.] $\nu(q)\le r$ if $q$ is a 1 or 2 point. $\nu(q)=r$ implies $\gamma(q)=r$.
\item[2'.] If $q$ is a 2 point and $\nu(q)=r-1$, then either $\tau(q)>0$ or 
$\gamma(q)=r$ or 
$\tau(q)=0$ and (\ref{eq64}) holds at $q$ with $0<d_i<i$, $e_i=i$ and
$\overline S_{r-1}(Y_1)$ contains a single curve $D$ containing $q$, and containing a 1 point, which has local equations
$x=z=0$ at $q$.
\item[3'.] $\nu(q)\le r-2$ if $q$ is a 3 point.
\end{description}
\end{Theorem}

\begin{pf}
Suppose that there exists a 2 curve $D\subset Y=Y_p$ such that $D$ is r-1 big.
Let $\pi_1:Y_1\rightarrow Y$ be the blowup of $D$. Then $A_r(Y_1)$ holds by Lemmas
\ref{Lemma500} and
\ref{Lemma501}, since $\nu(q)=r-1$ for all $q\in D$. 

Suppose that $q\in D$ and (\ref{eq62}) holds at $q$. Then $d_j\ge j$ for all $j$.
Suppose that $q'\in \pi_1^{-1}(q)$ and 
$\hat{\cal O}_{Y_1,q'}$ has regular parameters $(x_1,y_1,\overline z)$ such that
$$
\overline x=x_1, \overline y=x_1(y_1+\alpha)
$$
with $\alpha \ne 0$. Then $q'$ is a 1 point, so that $\nu(q')\le r$ and $\nu(q')=r$
implies $\gamma(q')=r$ by Lemma \ref{Lemma500}.

Suppose that $q'\in \pi_1^{-1}(q)$ and 
$q'$ has permissible parameters $(x_1,y_1,\overline z)$ such that
$$
\overline x=x_1, \overline y=x_1y_1
$$ 
Then
$$
\begin{array}{ll}
u&=(x_1^{a+b}y_1^b)^m\\
v&=P(x_1^{a+b}y_1^b) + x_1^{c+d+r-1}y_1^dF_{q'}\\
F_{q'}&=\frac{F_q}{x_1^{r-1}}= \tau + \sum_{j=1}^{r-1}y_1^{d_j}\Lambda_j(x_1y_1,\overline z)x_1^{d_j-j}
\end{array}
$$
so that $\nu(F_{q'})=0$.

Suppose that $q'\in \pi_1^{-1}(q)$ and 
$q'$ has permissible parameters $(x_1,y_1,\overline z)$ such that
$$
\overline x=x_1y_1, \overline y=y_1
$$ 
Then
$$
\begin{array}{ll}
u&=(x_1^{a}y_1^{a+b})^m\\
v&=P(x_1^{a}y_1^{a+b}) + x_1^{c}y_1^{c+d+r-1}F_{q'}\\
F_{q'}&=\frac{F_q}{y_1^{r-1}}= \tau x_1^{r-1} + \sum_{j=1}^{r-1}y_1^{d_j-j}\Lambda_j(y_1,\overline z)x_1^{r-1-j}
\end{array}
$$
so that either $\nu(F_{q'})<r-1$, or we are back in the form (\ref{eq62}) but the $d_j$ have decreased by $j$.

Suppose that $q\in D$ and (\ref{eq63}) holds at $q$.  Without loss of generality,
we may assume that $D$ has local equations $\overline x=\overline z=0$ at $q$.
Then $e_j\ge j$ for all $j$.
 Suppose that $q'\in \pi_1^{-1}(q)$ and 
$\hat{\cal O}_{Y_1,q'}$ has regular parameters $(x_1,\overline y,z_1)$ such that
$$
\overline x=x_1, \overline z=x_1(z_1+\alpha)
$$
with $\alpha \ne 0$.
Set 
$$
x_1 = \overline x_1(z_1+\alpha)^{-\frac{c}{a+c}}
$$
Set $\lambda = f+ (d+f+r-1)(\frac{-c}{a+c})$, 
$$
G(x_1,\overline y,z_1) = \frac{(z_1+\alpha)^{\lambda}F_q}{x_1^{r-1}}=(z_1+\alpha)^{\lambda}\tau+\sum_{j=1}^{r-1} (z_1+\alpha)^{\lambda+e_j}
\overline a_j(\overline y_,x_1(z_1+\alpha))\overline y^{d_j} x_1^{e_j-j}
$$
Then
$$
\begin{array}{ll}
u&=(\overline x_1^{a'}\overline y^{b'})^{m'}\\
v&=P((\overline x_1^{a'}\overline y^{b'})^{\frac{m'}{m}})+\overline x_1^{d+f+r-1}\overline y^eG
\end{array}
$$
$$
G(0,0,z_1) = (z_1+\alpha)^{\lambda}\left[\tau_0 + \sum_{d_j=0,e_j=j}\overline  a_j(0,0)(z_1+\alpha)^j\right]
$$
where $\tau_0=\tau(0,0,0)$,
$(a+c)m=a'm'$, $bm=b'm'$, $(a',b')=1$. 
$$
F_{q'}(0,0,z_1)=\left\{\begin{array}{ll}
G(0,0,z_1)&\text{ if }a'e-b'(d+f+r-1)\ne 0\\
G(0,0,z_1)-G(0,0,0)&\text{ if }a'e-b'(d+f+r-1)=0
\end{array}\right.
$$
Thus $\nu(F_{q'}(0,0,z_1))\le r$, except possibly if $\lambda=0$ and
$a'e-b'(d+f+r-1)=0$. Then we have 
\begin{equation}\label{eq1013}
af-c(d+r-1)=0
\end{equation} and 
\begin{equation}\label{eq1014}
[ae-b(d+r-1)]+[ce-fb]=0
\end{equation}
with $a,b,c>0$. Substituting $d+r-1=\frac{af}{c}$ into (\ref{eq1014}), we get
$$
(\frac{a}{c}+1)(ce-bf)=0
$$
so that 
\begin{equation}\label{eq1015}
ce-bf=0
\end{equation} 
and 
\begin{equation}\label{eq1016}
ae-b(d+r-1)=0
\end{equation}
(\ref{eq1013}), (\ref{eq1015}) and (\ref{eq1016}) cannot all hold since $F_q$ is
normalized. Thus $\nu(F_{q'}(0,0,z_1))\le r$ and
 $\nu(q')\le r$, $\gamma(q')\le r$.

Suppose that $q'\in \pi_1^{-1}(q)$ and 
$q'$ has permissible parameters $(x_1,\overline y,z_1)$ such that
$$
\overline x=x_1, \overline z=x_1z_1
$$ 
Then $\nu(q') = 0$.

Suppose that $q'\in \pi_1^{-1}(q)$ and 
$q'$ has permissible parameters $(x_1,y_1,\overline z)$ such that
$$
\overline x=x_1z_1, \overline z=z_1
$$ 
Then
we either have a 3 point with $\nu(q')<r-1$, or we are back in the form of (\ref{eq63}) with 
$e_i$ decreased by $i$.

Suppose that $q\in D$ and (\ref{eq64}) holds at $q$. Then $d_j\ge j$ for all $j$.
 Suppose that $q'\in \pi_1^{-1}(q)$ and 
$\hat{\cal O}_{Y_1,q'}$ has regular parameters $(x_1,y_1,\overline z)$ such that
$$
\overline x=x_1, \overline y=x_1(y_1+\alpha)
$$
with $\alpha \ne 0$.
Then $q'$ is a 1 point so that $\nu(q')\le r$ and $\gamma(q')\le r$
by Lemma \ref{Lemma500}.

Suppose that $q'\in \pi_1^{-1}(q)$ and 
$q'$ has permissible parameters $(x_1,y_1,\overline z)$ such that
$$
\overline x=x_1, \overline y=x_1y_1
$$ 
Then $\nu(q') = 0$.

Suppose that $q'\in \pi_1^{-1}(q)$ and 
$q'$ has permissible parameters $(x_1,y_1,\overline z)$ such that
$$
\overline x=x_1y_1, \overline y=y_1
$$ 
Then
we either have a 2 point with $\nu(q')<r-1$, or we are back in the form of (\ref{eq64}) with 
$d_i$ decreased by $i$.

After a finite number of blowups of 2 curves $\pi_1:Y_{n_0}\rightarrow Y$ we have that there
are no 2 curves $D$ on $Y_{n_0}$ such that $D$ is r-1 big. 
By Lemmas \ref{Lemma500}, \ref{Lemma501}, and since all 3 points $q\in(\pi\circ\pi_1)^{-1}(p)$
have $\nu(q)\le r-2$ (if $q$ is a 3 point and $\nu(q)=r-1$ so that $q$  satisfies (\ref{eq63}), then 
either $d_j\ge j$ for all $j$ or $e_j\ge j$ for all $j$), there exists a neighborhood $\overline U_1$ with the properties
asserted by the statement of the Theorem.

The only 2 points $q\in(\pi\circ\pi_1)^{-1}(p)$ where $\nu(q)=r-1$ and $\gamma(q)>r$ either
satisfy (\ref{eq62}) with some $d_j=j-1$ and $e_j=1$ (so that $\tau(q)\ge 1$) or
satisfy (\ref{eq64}) with $\nu(q)=r-1$ and $d_i<i$ so that $e_i\ge i$. 
Furthermore, $\overline x=0$ is 
a local equation at $q$ of the strict transform of the surface with local equation $x=0$ at $p$.

Suppose that $D\subset \overline S_r(X)$ is a curve containing $p$, and $\overline D$ is
the strict transform of $D$ on $\overline U_1$. $\overline D$ can only  intersect 
$(\pi\circ\pi_1)^{-1}(p)$ at  points $q$ such that $q$ is a 2 point and $\nu(q)=r$ or $\nu(q)=r-1$ by Lemma
\ref{Lemma4} and Lemma \ref{Lemma3}. We must either have $\gamma(q)\le r$ or $q$ satisfies (\ref{eq62}) or (\ref{eq64}).

Suppose that $q\in \overline D\cap (\pi\circ\pi_1)^{-1}(p)$ satisfies (\ref{eq62}) or
(\ref{eq64}) and $\overline y\in \hat{\cal I}_{\overline D,q}$. Then there exist
$a_j\in k$ such that 
$$
F_q-\sum a_j\frac{(\overline x^a\overline y^b)^j}{\overline x^c\overline y^d}\in (\overline y)^{r-1}+(\overline y,
f(\overline x,\overline z))^r
$$
where $\overline y=f(\overline x,\overline z)=0$ are local equations of $\overline D$
at $q$ (by Lemma \ref{Lemma7}). This is impossible by the form of $F_q$. Thus
$\overline x\in\hat{\cal I}_{\overline D,q}$.

Suppose that $q$ satisfies (\ref{eq62}). We have 
$$
\hat{\cal I}_{\overline D,q}=(\overline x,\overline z-\phi(\overline y))
$$
for some series $\phi$. There exist $a_j\in k$ such that,
when renormalizing with respect to these new parameters, 
$$
F_q-\sum a_j\frac{(\overline x^a\overline y^b)^j}{\overline x^c \overline y^d}\in (\overline x)^{r-1}+(\overline x,\overline z-\phi(\overline y))^r
$$
by Lemma \ref{Lemma7}.
Setting $\overline x=0$ in 
$$
F_q-\sum a_j\frac{(\overline x^a\overline y^b)^j}{\overline x^c\overline y^d},
$$
we get $\overline y^{d_{r-1}}\Lambda_{r-1}(\overline y,\overline z)$ or
$\overline y^{d_{r-1}}\Lambda_{r-1}(\overline y,\overline z)+\tilde c\overline y^{\overline n}
$ for some $\tilde c\in k$, $\overline n\in \bold N$.
Thus 
$$
(\overline z-\phi(\overline y))^r\mid \overline y^{d_{r-1}}\Lambda_{r-1}(\overline y,\overline z)
$$
or
$$
(\overline z-\phi(\overline y))^r\mid \overline y^{d_{r-1}}\Lambda_{r-1}(\overline y,\overline z)+\tilde c\overline y^{\overline n}
$$

which is nonzero since $\overline x\not\,\mid F_q$. As $\nu(\Lambda_{r-1}(\overline y,\overline z))\le 1$, this is a contradiction. Thus $q$ cannot have the form of (\ref{eq62}). 

Suppose that $q$ satisfies (\ref{eq64}). We have $\hat{\cal I}_{\overline D,q}
=(\overline x,\overline z-\phi(\overline y))$ for some series $\phi$. There exists
$a_j \in k$ such that
$$
F_q-\sum a_j\frac{(\overline x^a\overline y^b)^j}{\overline x^c\overline y^d}\in (\overline x)^{r-1}+(\overline x,\overline z
-\phi(\overline y))^r
$$
by Lemma \ref{Lemma7}.
Setting $\overline x=0$ in 
$$
F_q-\sum a_j\frac{(\overline x^a\overline y^b)^j}{\overline x^c\overline y^d},
$$
we get $\overline a_{r-1}(\overline y,\overline z)\overline y^{d_{r-1}}\overline z^{e_{r-1}}$
or $\overline a_{r-1}(\overline y,\overline z)\overline y^{d_{r-1}}\overline z^{e_{r-1}}
+\tilde c\overline y^{\overline n}$ for some $\tilde c\in k$, $\overline n\in \bold N$.
Thus
$$
(\overline z-\phi(\overline y))^r\mid \overline y^{d_{r-1}}\overline z^{e_{r-1}}
$$
or
$$
(\overline z-\phi(\overline y))^r\mid \overline a_{r-1}(\overline y,\overline z)\overline y^{d_{r-1}}\overline z^{e_{r-1}}+\tilde c\overline y^{\overline n}
$$
which is nonzero since $\overline x\not\,\mid F_q$.  In either case,
we have
$$
(\overline z-\phi(\overline y))^{r-1}\mid\overline y^{d_{r-1}}\overline z^{e_{r-1}-1}
$$
since
$$
\frac{\partial}{\partial \overline z}(\overline a_{r-1}\overline y^{d_{r-1}}\overline z^{e_{r-1}}
+\tilde c\overline y^{\overline n})
=\overline z^{e_{r-1}-1}\overline y^{d_{r-1}}(e_{r-1}\overline a_{r-1}+\frac{\partial \overline a_{r-1}}{\partial \overline z}\overline z),
$$
which implies that  $e_{r-1}\ge r$ and $\phi(\overline y)=0$. Thus $\overline x=\overline z=0$
are local equations of $\overline D$ at $q$.

Suppose that $\overline D$ is such that $\overline D$ is r small and $q\in (\pi\circ\pi_1)^{-1}(p)\cap\overline D$ satisfies
$\nu(q)=r$ and $\gamma(q)=r$. By Lemma \ref{Lemma41}, there exists a sequence of quadratic
transforms $\sigma_1:W_1\rightarrow Y_{n_0}$ such that the strict transform $\tilde D$
of $\overline D$ intersects $\sigma_1^{-1}(q)$ in a 2 point $q'$ such that $\nu(q')=r-1$ and
$\gamma(q')=r$. Furthermore, there are no 2 curves $C$ in $\sigma_1^{-1}(q)$ such that
$C$ is r-1 big, and 1'. - 3'. of the conclusions of the
Theorem hold at all points of $\sigma_1^{-1}(q)$. Thus there exists a sequence of
quadratic transforms $\sigma:W\rightarrow Y_{n_0}$, centered at 2 points
$\{q_1,\ldots,q_m\}$ such that $\nu(q_i)=r$ and $\gamma(q_i)=r$ on the strict
transform $\overline D$ of curves $D$ in $\overline S_r(X)$ containing $p$ such that if
$\overline D\subset \overline S_r(W)$ is the strict transform of a curve $D\subset \overline
S_r(X)$ containing $p$, and $\overline D$ is r small then $\overline D$ intersects $(\pi\circ\pi_1\circ\sigma)^{-1}(p)$
in  2 points $q$ of the form of (\ref{eq64}), and in 2 points $q$ such that $\nu(q)=r-1$
and $\gamma(q)=r$. $W$ contains no 2 curves $C$ such that $C$ is r-1 big, $\overline S_r(W)$ makes SNCs with $\overline B_2(W)$ and $\gamma(q)\le r$ for all exceptional 1 and 2 points of $\sigma$, $\nu(q)=0$ for all exceptional 3 points
of $\sigma$.

Suppose that $\overline D\subset W$ is the strict transform of a curve $D$ in $\overline
S_r(X)$ containing $p$. First suppose that $\overline D$ is r big. If $q\in \overline D\cap (\pi\circ\pi_1\circ\sigma)^{-1}(p)$,
then $q$ must be a 2 point with $\nu(q)=\gamma(q)=r$. Suppose that $\lambda_1:Z_1\rightarrow W$
is the blowup of $\overline D$. By Lemma \ref{Lemma654}, 1'. - 3'.
of the conclusions of the Theorem hold on $Z_1$, and
the conclusions of $\overline U_2$ hold in a neighborhood of $\lambda_1^{-1}(q)$.

If $\overline D$ is r small, then if
$q\in\overline D\cap (\pi\circ\pi_1\circ\sigma)^{-1}(p)$, $q$ must be either a 
2 point where $\nu(q)=r-1$ and $\gamma(q)=r$ or $q$ satisfies (\ref{eq64}) and 
$\overline x=\overline z=0$ are local equations of $\overline D$ at $q$.

Let $\lambda_1:Z_1\rightarrow W$ be the blowup of $\overline D$. If 
$q\in\overline D\cap (\pi\circ\pi_1\circ\sigma)^{-1}(p)$ is a 2 point such that
$\nu(p)=r-1$ and $\gamma(p)=r$, then  1.' - 3.' of the
conclusions of the Theorem hold and the conclusions of $\overline U_2$ hold
in a neighborhood of $\lambda_1^{-1}(q)$ (since $r\ge 3$)  by Lemma \ref{Lemma41}.

Suppose that $q\in\overline D\cap(\pi\circ\pi_1\circ\sigma)^{-1}(p)$ satisfies
(\ref{eq64}). $\overline x=\overline z=0$ are local equations of $\overline D$
at $q$ and $d_i<i$. Since $\overline D\subset \overline S_r(W)$, we have $e_i>i$.

Since 
$A_r(X)$ holds, $\gamma(q')=r$ if $q'\in \overline D$ is a 1 point. Then $e_{r-1}=r$
in (\ref{eq64}).
 Suppose that $q'\in\lambda_1^{-1}(q)$, $q'$ has
permissible parameters $(x_1, \overline y, z_1)$ such that 
$$
\overline x=x_1, \overline z=x_1(z_1+\alpha)
$$
$q'$ is a 2 point. $e_i>i$ implies
$$
\frac{F_p}{x_1^{r-1}}=\tau+x_1\Omega
$$
so that $\nu(q)=0$.

Suppose that $q'\in \lambda_1^{-1}(q)$ has permissible parameters $(x_1,\overline y,z_1)$
such that 
$$
\overline x=x_1z_1, \overline z=z_1
$$
Then we either have a 3 point with $\nu(q')<r-1$, or we are in the form of
(\ref{eq63}) with $e_i$ decreased by $i$, $d_i<i$, and $\overline x=\overline z=0$ is a local equation of a 2 curve which is a section over $D$. 
Since we have $e_{r-1}=1$ in (\ref{eq63}), we must have $e_i<i$ (since $r\ge 2$), $d_i<i$ so that
$\nu(q')<r-1$ by Remark \ref{Remark1049}.

Thus the conclusions of $\overline U_2$ hold in a neighborhood of $\lambda_1^{-1}(q)$ by
Lemma \ref{Lemma41}.

Then if $\overline\lambda:Z\rightarrow W$ is the sequence of monodial transforms (in any order)
centered at the strict transforms of curves $C$ in $\overline S_r(X)$, the conclusions
of 1'. - 3'. of the Theorem hold, except possibly at a finite number of points $q$ of the
form of (\ref{eq64}) with $d_i<i$ and $e_i\ge i$. There are no 2 curves $C\subset Z$
such that $C$ is r-1 big.

If $q\in(\pi\circ\pi_1\circ\sigma\circ\overline\lambda)^{-1}(p)$
does not satisfy one of 1' - 3' of the conclusions of the Theorem then $q$ satisfies (\ref{eq64}), $\tau(q)=0$ and $d_i<0$ so that $e_i\ge i$.
 $\overline x=0$ is then a local equation of the surface with local equation
$x=0$ at $p$.

 $e_i\ge i$ implies $F_q\in(\overline x,\overline z)^{r-1}$.
Since $r\ge 3$, this implies (by Lemma \ref{Lemma302}) that there exists an algebraic curve 
$\overline D\subset \overline S_{r-1}(Z)$ such that $\overline x=\overline z=0$ are
local equations of a formal branch of $\overline D$. $\overline D$ is necessarily
the strict transform of a curve $D\subset \overline S_{r-1}(U)$, since
$\overline x\in\hat{\cal I}_{\overline D,q}$ and $\overline y\not\in
\hat{\cal I}_{\overline D,q}$. $\overline D$ is thus nonsingular at $q$, by the
conclusions of Theorem \ref{Theorem27}. Thus $\overline x=\overline z=0$ are
local equations of an algebraic curve $\overline D$ at $q$. If $e_i>i$, then
$\overline D\subset\overline S_r(Z)$ and if $e_i=i$, then $\overline D\subset
\overline S_{r-1}(Z)$.

Since the strict transforms of all curves in $\overline S_r(X)$ have been blown up
in the map $Z\rightarrow W$, we must have $e_i=i$.

Thus at $q$, (\ref{eq64}) holds, $e_i=i$, $d_i<i$ and $\tau(q)=0$.
Let $T$ be the component of $E_X$ with local equation $x=0$ at $p$.
$\tau(q)=0$ implies $d_i>0$, which implies $d_j>0$ for all $j$, so that
$F_q=\tau \overline x^{r-1}+\overline y\Omega$. Since $\overline x=0$ is a local equation of the strict transform
$T'$ of $T$, the only curve in $\overline S_{r-1}(Z)\cap T'$ containing $q$ is the curve $\overline D$ with
local equations $\overline x=\overline z=0$, since a curve in $\overline S_{r-1}(Z)\cap T'$ containing $q$ must be the strict transform of a curve in $\overline S_{r-1}(U)\cap T$.

$$
\begin{array}{ll}
u&=(\overline x^a\overline y^b)^m\\
v&=P(\overline x^a\overline y^b)+\overline x^c\overline y^dF
\end{array}
$$
and $F_q=\tau \overline x^{r-1}+\overline y\Omega$, where $\tau$ is a unit. 

We will show that there does not exist a curve $C\subset\overline S_{r-1}(Z)$
containing $q$ (and a 1 point) such that $\overline y\in\hat{\cal I}_{C,q}$.

After a permissible change of parameters, we may assume that $\overline y,\overline z\in{\cal O}_{X,q}$ with 
$$
F_q=\tau\overline x^{r-1}+\overline y\Omega.
$$
${\cal I}_{C,q}=(\overline y,g(\overline x,\overline z))$. Set $\overline x=\tilde x^b$. By Lemma \ref{Lemma659},
either there exists a series $f$ such that 
\begin{equation}\label{eq642}
\tilde x^{bc-ad}F_q-f(\tilde x^a \overline y)\in \left((\overline y,g(\tilde x^b,\overline z))^{r-1}+(\overline y)^{r-2}\right)k[[\tilde x,
\overline y,\overline z]]
\end{equation}
if $bc-ad\ge 0$ or 
\begin{equation}\label{eq643}
F_q-f(\tilde x^a\overline y)\tilde x^{ad-bc}\in \left((\overline y,g(\tilde x^b,\overline z))^{r-1}+(\overline y)^{r-2}\right)k[[\tilde x,
\overline y,\overline z]]
\end{equation}
if $ad-bc>0$.

If (\ref{eq642}) holds, since $\nu(q)>0$, we have $\nu(f)>0$, which implies
$$
g(\tilde x^b,\overline z)^{r-1}\mid \tilde x^{b(r-1)+bc-ad},
$$
 and $g=\overline x$, a contradiction since $C$ is then a 2 curve.

Suppose that (\ref{eq643}) holds. Let $\overline c=f(0)$.
$$
\tau(\overline x,0,\overline z)\tilde x^{b(r-1)}-\overline c\tilde x^{ad-bc}=h(\tilde x,\overline z)
g(\tilde x^b,\overline z)^{r-1}
$$
for some series $h$. If $ad-bc=b(r-1)$ then $ad-b(c+r-1)=0$, a contradiction to
the assumption that $F_q$ is normalized.

Let 
$$
\overline d= \left\{\begin{array}{ll}
\text{min}\{b(r-1), ad-bc\}&\text{if }\overline c\ne 0\\
b(r-1)&\text{if }\overline c=0.
\end{array}\right.
$$
$$
\tau(\overline x,0,\overline z)\tilde x^{b(r-1)}-\overline c\tilde x^{ad-bc}=\Lambda \tilde x^{\overline d}
$$
where $\Lambda$ is a unit. Then $g(\tilde x^b,\overline z)^{r-1}$ is a power of $\tilde x$,
and $g=\tilde x$, a contradiction since $C$ is then a 2 curve.

Thus the conclusions of 1'. - 3'. of the Theorem hold and the conclusions of $\overline U_2$
hold on $Z$.
\end{pf}

\begin{Theorem}\label{Theorem1018}
Suppose that $r=2$, $A_2(X)$ holds, $p$ is a 2 point with $\nu(p)=r-1=1$,  
and $L(x,0,0)\ne 0$, as in the assumptions of Theorem \ref{Theorem27}.
Let $R ={\cal O}_{X,p}$. Suppose that  $\pi:Y_p\rightarrow \text{Spec}(R)$ is the morphism of
Theorem \ref{Theorem27}. 

Let 
$$
\cdots\rightarrow Y_n\rightarrow\cdots\rightarrow Y_1\rightarrow Y_p
$$
be a sequence of permissible monodial transforms centered at 2 curves $D$
such that $D$ is 1 big.
Then there exists $n_0<\infty$ such that 
$$
V_p=Y_{n_0}\stackrel{\pi_1}{\rightarrow} Y_p\rightarrow\text{spec}(R)
$$
extends to a permissible sequence of monodial transforms 
$$
\overline U_1\rightarrow\overline U\rightarrow U
$$
over an affine neighborhood $U$ of $p$ (with the notation of Theorem \ref{Theorem27})
such that $\overline U_1$ contains no 2 curves $D$ such that
 $D$ is 1 big or 2 small, and for $q\in\overline U_1$,
\begin{enumerate}
\item If $q$ is a 1 or a 2 point then $\nu(q)\le 2$. $\nu(q)=2$ implies $\gamma(q)=2$.
\item If $q$ is a 3 point then $\nu(q)=0$.
\item $\overline S_2(\overline U_1)$ makes SNCs with $\overline B_2(\overline U_1)$.
\end{enumerate}
There exists a sequence of quadratic transforms $W_p\rightarrow V_p$ such that if
$Z_p\rightarrow W_p$ is the sequence of monodial transforms (in any order) centered
 at the strict transforms $C'$ of
curves $C$ in $\overline S_2(X)$, followed by monoidial transforms centered at any
2 curves $\overline C$ which are sections over $C'$ such that $\overline C$ is 1 big,  
then 
$$
Z_p\rightarrow W_p\rightarrow V_p\rightarrow Y_p\rightarrow\text{spec}(R)
$$
extends to a permissible sequence of monodial transforms
$$
\overline\pi:\overline U_2\rightarrow \overline U_1\rightarrow \overline U\rightarrow U
$$
over an affine neighborhood of $p$ such that 
$\overline U_2$ contains no 2 curves $D$ such that 
$D$ is 1 big or 2 small.
$\overline S_2(\overline U_2)$
makes SNCs with $\overline B_2(\overline U_2)$, and 
if $q\in \overline\pi^{-1}(p)$,
\begin{description}
\item[1'.] $\nu(q)\le 2$ if $q$ is a 1 or 2 point. $\nu(q)=2$ implies $\gamma(q)=2$.
\item[2'.] If $q$ is a 2 point and $\nu(q)=1$, then either $q$ is resolved or 
$\gamma(q)=2$. 
\item[3'.] $\nu(q)=0$ if $q$ is a 3 point.
\end{description}
\end{Theorem}

\begin{pf}
We can construct 
$$
W=W_p\stackrel{\sigma}{\rightarrow}V_p=Y_{n_0}\stackrel{\pi_1}{\rightarrow}Y=Y_p
\stackrel{\pi}{\rightarrow}\text{spec}(R)
$$
exactly as in the proof of Theorem \ref{Theorem28}.

If $q\in(\pi\circ\pi_1\circ\sigma)^{-1}(p)$ satisfies \ref{eq64}), then 
\begin{equation}\label{eq1017}
\begin{array}{ll}
u&=(\overline x^a\overline y^b)^m\\
F_q&=\overline x+\overline z^{e_1}.
\end{array}
\end{equation}

If $D\subset\overline S_2(X)$ contains $p$, and $\overline D$ is the strict 
transform of $D$ on $W$, and $q\in\overline D\cap (\pi\circ\pi_1\circ\sigma)^{-1}(p)$,
then either $q$ satisfies (\ref{eq1017}),
$$
F_q=\overline x+\overline z^{e_1}
$$
with $e_1\ge 2$ and $\overline x=\overline z=0$ are local equations of $\overline D$
at $q$, or $q$ is a 2 point with $\nu(q)=1$, $\gamma(q)=2$, or $\overline D$ is 2 big
 and $q$ is a 2 point with $\nu(q)=\gamma(q)=2$.

Suppose that $\overline D\subset\overline S_2(W)$ is the strict transform of
$D\subset\overline S_2(X)$ such that $p\in D$. Let $\lambda_1:Z_1\rightarrow W$
be the blowup of $\overline D$.

First suppose that $\overline D$ is 2 big.
Suppose that $q\in\overline D\cap (\pi\circ\pi_1\circ\sigma)^{-1}(p)$. Then
$q$ is a 2 point with $\nu(q)=\gamma(q)=2$. By Lemma \ref{Lemma654}, 1' - 3' of the
conclusions of the Theorem hold on $Z_1$, and the conclusions of $\overline U_2$
hold in a neighborhood of $\lambda_1^{-1}(q)$.

Suppose that $\overline D$ is 2 small. If
$q\in\overline D\cap(\pi\circ\pi_1\circ\sigma)^{-1}(p)$, then either $q$ is a 2 point
with $\nu(q)=1$, $\gamma(q)=2$, or $q$ satisfies (\ref{eq1017}) with $e_1\ge 2$, $\overline
x=\overline z=0$ are local equations of $\overline D$ at $q$.

If $q\in\overline D\cap (\pi\circ\pi_1\circ\sigma)^{-1}(p)$ is a 2 point with
$\nu(q)=1$, $\gamma(q)=2$, then by Lemma \ref{Lemma41} and Lemma \ref{Lemma97},
either 1' - 3' of the conclusions of the Theorem and the conclusions of $\overline U_2$
hold in a neighborhood of $\lambda_1^{-1}(q)$ or there exists a 2 curve $\overline C$
which is a section over $\overline D$, such that if $\lambda_2:Z_2\rightarrow Z_1$ is the
blowup of $\overline C$, then 1' - 3' of the conclusions of the Theorem hold and the
conclusions of $\overline U_2$ hold in a neighborhood of $(\lambda_1\circ\lambda_2)^{-1}(q)$.

Suppose that $q\in\overline D\cap (\pi\circ\pi_1\circ\sigma)^{-1}(p)$ satisfies (\ref{eq1017})
with $e_1\ge 2$, $\overline x=\overline z=0$ are local equations of $\overline D$ at $q$.
Since $A_2(X)$ holds, $\gamma(q')=2$ if $q'\in\overline D$ is a 1 point, so that $e_1=2$.
If $q'\in\lambda_1^{-1}(q)$ we have $\nu(q')=0$ except if $q'$ is the 3 point with
permissible parameters $(x_1,\overline y,z_1)$ such that $\overline x=x_1z_1$,
$\overline z=z_1$. Then $F_{q'}=x_1+z_1$.

Let $\overline C$ be the 2 curve through $q'$ with local equations $x_1=z_1=0$ at $q'$.
$\overline C$ is a section over $\overline D$. By Lemma \ref{Lemma655}, $F_a\in\hat{\cal I}_{\overline C,a}$ for all $a\in \overline C$, so that $\overline C$ is 1 big.

Let $\lambda_2:Z_2\rightarrow Z_1$ be the blowup of $\overline C$. Then 1' - 3' of the
conclusions of the Theorem, and the conclusions of $\overline U_2$ hold in a
neighborhood of $(\lambda_1\circ\lambda_2)^{-1}(q)$.

Then if $\overline \lambda:Z\rightarrow W$ is the sequence of monodial transforms
(in any order) centered at the strict transform of curves $C$ in $\overline S_2(X)$,
followed by the monodial transforms centered at 2 curves $\overline C$ which are
sections over $C'$ such that $\overline C$ is 1 big,
we have that there are no 2 curves $C\subset Z$ such that $\overline C$ is 1 big.

The only points of $Z$ which may not satisfy the conclusions of 1' - 3' of the Theorem
are the 2 points $q\in (\pi\circ\pi_1\circ\sigma\circ\overline\lambda)^{-1}(p)$
which satisfy (\ref{eq1017}) with $\tau(q)=0$.
Then (after a permissible change of parameters)
$$
\begin{array}{ll}
u&=(\overline x^a\overline y^b)^m\\
F_q&=\overline x+\overline z^{e_r}
\end{array}
$$
with $e_r\ge 2$.

By Lemma \ref{Lemma302}, there exists an algebraic curve $\overline D\subset \overline S_2(Z)$
such that $\overline x=\overline z=0$ are local equations of $D$ at $q$. Since
$\overline x=0$ is a local equation of the strict transform of the component of
$E_X$ with local equation $x=0$, and $\overline D$ is not contained in the component of
$E_X$ with local equation $\overline y=0$, $\overline D$ is the strict transform of a curve
in $\overline S_r(X)$ containing $p$, a contradiction to the construction of $Z$.
Thus $q$ is a resolved point.

\end{pf}

\begin{Theorem}\label{Theorem30}(Theorem30)
Suppose that $r\ge 2$, $A_r(X)$ holds, $p\in X$ is a 3 point with $\nu(p)=r-1$, 
and
we have permissible parameters $(x,y,z)$ at $p$ for $u,v$ (with $y,z\in {\cal O}_{X,p}$)
 such that 
\begin{equation}\label{eq509}
\begin{array}{ll}
u&= (x^by^{a+nb}z^{a+(n+1)b})^m\\
v&=P(x^by^{a+nb}z^{a+(n+1)b})+x^dy^{c+n(d+r)}z^{c+(n+1)(d+r)}F_p\\
F_p&=\tau x^{r-1}+\sum_{i=1}^{r-1}a_i(y,z)x^{r-i-1}
\end{array}
\end{equation}
with $n\ge 0$, $a,b>0$,
 $L(x,0,0)\ne 0$, so that $\tau$ is a unit, and (eq510)
\begin{equation}\label{eq510}
a(d+r-1)-bc=0
\end{equation}
or 
\begin{equation}\label{eq511}
\begin{array}{ll}
u&=(x^by^az^b)^m\\
v&=P(x^by^az^b)+x^dy^cz^{d+r}F_p\\
F_p&=\tau x^{r-1}+\sum_{i=1}^{r-1} a_i(y,z)x^{r-i-1}
\end{array}
\end{equation}
with $n\ge 0$, $a,b>0$,
 $L(x,0,0)\ne 0$, so that $\tau$ is a unit, and 
\begin{equation}\label{eq512}
a(d+r-1)-bc=0
\end{equation}

Let $R ={\cal O}_{X,p}$. Let $C$ be the 2 curve  with local equations
$y=z=0$ at $p$.
Then
there exists a finite sequence of permissible monodial transforms $\pi:Y\rightarrow \text{Spec}(R)$
centered at sections over $C = V(y,z)$, such that for $q\in \pi^{-1}(p)$, $F_q$ has one of the  forms
(\ref{eq62}), (\ref{eq63}) or (\ref{eq64}) of Theorem \ref{Theorem27}.

In all these cases $\overline x=x$ and $\overline x=0$ is a local equation at $q$ of the strict
transform of the component of $E_X$ with local equation $x=0$ at $p$.

There exists an affine neighborhood $U$ of $p$ such that $Y\rightarrow\text{spec}(R)$
extends to a sequence of permissible monoidal transforms $\overline U\rightarrow U$
such that $A_r(\overline U)$ holds.

$\gamma(q)\le 1$ at a generic point $q$ of $C$, so
that all points $q'$ on the fiber of the blowup of $C$ over $q$ are resolved.

Suppose that $r\ge 3$. Let $D_i$ be the curves in $\overline S_{r-1}(X)$ which contain $p$, and $x\in\hat{\cal
I}_{D_i,p}$. We further have  that the strict transforms $\overline D_i$ of the $D_i$
on $\overline U$ are nonsingular, disjoint, and make SNCs with $\overline B_2(\overline U)$.

\end{Theorem}

\begin{pf}
We modify the proof of Theorem \ref{Theorem27} to prove this Theorem. In the sequence of
(\ref{eq507}) we must add a new case,
\begin{description}
\item[Case 0]
$$
\begin{array}{ll}
u&=(x^{a_i}(y_i^{\overline b_i}z_i^{\overline c_i})^{k_i})^{m_i}\\
v&=P_i(x^{a_i}(y_i^{\overline b_i}z_i^{\overline c_i})^{k_i})+x^cy_i^{\overline d_i}
z_i^{\overline e_i}F_i
\end{array}
$$
\item[Case 0a]
$$
y_i=y_{i+1}, z_i=y_{i+1}z_{i+1}
$$
\item[Case 0b] 
$$
y_i=y_{i+1}z_{i+1}, z_i=z_{i+1}
$$
\item[Case 0c]
$$
y_i=y_{i+1}(z_{i+1}+\alpha_{i+1})^{-\frac{\overline c_i}{\overline b_i+\overline c_i}},
$$
$$
z_i=y_{i+1}(z_{i+1}+\alpha_{i+1})^{\frac{\overline b_i}{\overline b_i+\overline c_i}}
$$
\end{description}

In the sequence (\ref{eq507}), $Y_0$ has the form Case 0 (and not Case 2). The
transformations of type 0a and type 0b produce a $p_{i+1}$ of the type of Case 0,
and $F_i=F_{i+1}$.

If all transformations in (\ref{eq507}) are of types 0a or 0b, then we
eventually get a $p_i$ of type (\ref{eq63}) by Lemmas \ref{Lemma1024}, \ref{Lemma23} and
Corollary \ref{Corollary24}.

Otherwise, we eventually reach a first $p_i$ where $p_{i+1}$ is obtained by a
transformation of type 0c. We have $F_i=F_p$.
$$
\begin{array}{ll}
u&=(x^{a_i}y_{i+1}^{(\overline b_i+\overline c_i)k_i})^{m_i}=(x^{a_{i+1}}y_{i+1}^{b_{i+1}})^{m_{i+1}}\\
v&=P_i(x^{a_{i}}y_{i+1}^{(b_{i}+\overline c_i)k_i})
+x^cy_{i+1}^{d_{i+1}}(z_{i+1}+\alpha_{i+1})^{\lambda}F_p
\end{array}
$$
$$
\lambda=\frac{\overline e_i\overline b_i-\overline d_i\overline c_i}{\overline b_i+\overline c_i},
$$
 $d_{i+1}=\overline d_i+\overline e_i$.
Thus $p_{i+1}$ has the form of Case 1, with
$$
F_{i+1}=(z_{i+1}+\alpha_{i+1})^{\lambda} F_p-\frac{g_{i+1}(x^{a_{i+1}}y_{i+1}^{b_{i+1}})}
{x^cy_{i+1}^{d_{i+1}}}
$$
Thus $p_{i+1}$ satisfies the assumptions of Theorem \ref{Theorem27}, provided
$\nu(F_{i+1}(x,0,0))=r-1$.

$p_i$ has permissible parameters $(x,y_1,z_1)$ with $y_1=y_i$, $z_1=z_i$ such that

$$
y=y_1^{\alpha}z_1^{\beta}, z=y_1^{\gamma}z_1^{\delta}
$$
with $\alpha\delta-\beta\gamma=\pm 1$,

$p_{i+1}$ has permissible parameters $(x,y_2,z_2)$ with $y_2=y_{i+1}$, $z_2=z_{i+1}$
such that 

$$
y_1=y_2, z_1= y_2(z_2+\overline \alpha)
$$
with $\overline\alpha\ne 0$.

First suppose that we are in the situation of (\ref{eq509}) and (\ref{eq510}).
Set
$$
\lambda_1=\alpha(a+nb)+\gamma(a+(n+1)b)+\beta(a+nb)+\delta(a+(n+1)b),
$$ 
$$
\lambda_2= \beta(\alpha+nb)+\delta(a+(n+1)b).
$$
$$
y_2=\overline y_2(z_2+\overline \alpha)^{-\frac{\lambda_2}{\lambda_1}}
$$
$$
\begin{array}{ll}
u&=(x^by_1^{\alpha(a+nb)+\gamma(a+(n+1)b)}z_1^{\beta(a+nb)+\delta(a+(n+1)b)})^m\\
&=(x^by_2^{\lambda_1}
(z_2+\overline \alpha)^{\lambda_2})^m\\
&=(x^b\overline y_2^{\lambda_1})^m
\end{array}
$$
$$
\begin{array}{ll}
x^{d+r-1}y^{c+n(d+r)}z^{c+(n+1)(d+r)}&=
x^{d+r-1}y_1^{\alpha(c+n(d+r))+\gamma(c+(n+1)(d+r))}z_1^{\beta(c+n(d+r))
+\delta(c+(n+1)(d+r))}\\
&=x^{d+r-1}y_2^{\lambda_3}(z_2+\overline \alpha)^{\lambda_4}\\
&=x^{d+r-1}\overline y_2^{\lambda_3}(z_2+\overline\alpha)^{\lambda_4-
\frac{\lambda_2\lambda_3}{\lambda_1}}
\end{array}
$$
where
$$
\lambda_3=\alpha(c+n(d+r))+\gamma(c+(n+1)(d+r))+\beta(c+n(d+r))+
\delta(c+(n+1)(d+r)),
$$
$$
\lambda_4=\beta(c+n(d+r))
+\delta(c+(n+1)(d+r)).
$$
$$
\begin{array}{ll}
&b\lambda_3-(d+r-1)\lambda_1\\
&=b[(\alpha+\beta)(c+n(d+r))+(\gamma+\delta)(c+(n+1)(d+r))]\\
&-(d+r-1)[(\alpha+\beta)(a+nb)+(\gamma+\delta)(a+(n+1)b)]\\
&=(\alpha+\beta+\gamma+\delta)[bc-(d+r-1)a]
+nb(\alpha+\beta+\gamma+\delta)+b(\gamma+\delta)\\
&=n(\alpha+\beta+\gamma+\delta)b+(\gamma+\delta)b>0.
\end{array}
$$
Since we cannot remove
$$
\tau_0\overline \alpha^{\lambda_4-\frac{\lambda_2\lambda_3}{\lambda_1}}
x^{d+r-1}\overline y_2^{\lambda_3}
$$
from $x^dy^{c+n(d+r)}z^{c+(n+1)(d+r)}F_p$, where $\tau_0=\tau(0,0,0)$,
 when normalizing to obtain $F_{p_{i+1}}$,
we must have $\nu(F_{p_{i+1}}(x,0,0))=r-1$.

Now suppose that we are in the situation of (\ref{eq511}) and (\ref{eq512}).
Set
$$
\lambda_1=\alpha a+\gamma b+\beta a+\delta b,
$$
$$
\lambda_2=\beta a+\delta b.
$$
$$
y_2=\overline y_2(z_2+\overline \alpha)^{-\frac{\lambda_2}{\lambda_1}}
$$
$$
\begin{array}{ll}
u&=(x^by_1^{\alpha a+\gamma b}z_1^{\beta a+\delta b})^m\\
&=(x^by_2^{\lambda_1}(z_2+\overline \alpha)^{\lambda_2})^m\\
&=(x^b\overline y_2^{\lambda_1})^m
\end{array}
$$
Set
$$
\lambda_3=\alpha c+\gamma(d+r)+\beta c+\delta(d+r),
$$
$$
\lambda_4=\beta c+\delta(d+r).
$$
$$
\begin{array}{ll}
x^{d+r-1}y^cz^{d+r}&= x^{d+r-1}y_1^{\alpha c+\gamma(d+r)}z_1^{\beta c+\delta(d+r)}\\
&=x^{d+r-1}y_2^{\lambda_3}(z_2+\overline\alpha)^{\lambda_4}\\
&=x^{d+r-1}\overline y_2^{\lambda_3}(z_2+\overline \alpha)^{\lambda_4
-\frac{\lambda_2\lambda_3}{\lambda_1}}
\end{array}
$$
$$
\begin{array}{ll}
b\lambda_3-\lambda_1(d+r-1)
&=b(\alpha c+\gamma(d+r)+\beta c+\delta(d+r))
-(\alpha a+\gamma b+\beta a+\delta b)(d+r-1)\\
&=(\alpha+\beta)[bc-a(d+r-1)]
+(\gamma+\delta)[b(d+r)-b(d+r-1)]\\
&=(\gamma+\delta)b\ne 0
\end{array}
$$
Thus $\nu(F_{p_{i+1}}(x,0,0))=r-1$ in this case also. The proof now preceeds as in Theorem
\ref{Theorem27}.
\end{pf}

\begin{Theorem}\label{Theorem29} 
Suppose that $r\ge 3$, $A_r(X)$ holds, $p$ is a 3 point with $\nu(p)=r-1$, and
we have permissible parameters $(x,y,z)$ at $p$ for $u,v$ (with $y,z\in {\cal O}_{X,p}$)
 such that 
$$
\begin{array}{ll}
u&= (x^by^{a+nb}z^{a+(n+1)b})^m\\
v&=P(x^by^{a+nb}z^{a+(n+1)b})+x^dy^{c+n(d+r)}z^{c+(n+1)(d+r)}F_p\\
F_p&=\tau x^{r-1}+\sum_{i=1}^{r-1}a_i(y,z)x^{r-i-1}
\end{array}
$$
with $n\ge 0$, $a,b>0$,
$\nu(p)=r-1$, $L(x,0,0)\ne 0$, so that $\tau$ is a unit, and 
$$
a(d+r-1)-bc=0
$$
or 
$$
\begin{array}{ll}
u&=(x^by^az^b)^m\\
v&=P(x^by^az^b)+x^dy^cz^{d+r}F_q\\
F_p&=\tau x^{r-1}+\sum a_i(y,z)x^{r-i-1}
\end{array}
$$
with $n\ge 0$, $a,b>0$,
$\nu(p)=r-1$, $L(x,0,0)\ne 0$, so that $\tau$ is a unit, and 
$$
a(d+r-1)-bc=0
$$

Let $R ={\cal O}_{X,p}$.

 Suppose that  $\pi:Y_p\rightarrow \text{Spec}(R)$ is the morphism of
Theorem \ref{Theorem30}.

Let 
$$
\cdots\rightarrow Y_n\rightarrow\cdots\rightarrow Y_1\rightarrow Y_p
$$
be a sequence of permissible monodial transforms centered at 2 curves $D$
such that  $D$ is r-1 big. Then there exists $n_0<\infty$ such that 
$$
V_p=Y_{n_0}\stackrel{\pi_1}{\rightarrow} Y\rightarrow\text{spec}(R)
$$
extends to a permissible sequence of monodial transforms 
$$
\overline U_1\rightarrow\overline U\rightarrow U
$$
over an affine neighborhood $U$ of $p$ (with the notation of Theorem \ref{Theorem30})
such that $\overline U_1$ contains no 2 curves $D$ such that 
$D$ is r-1 big or r small, and for $q\in\overline U_1$,
\begin{enumerate}
\item If $q$ is a 1 or a 2 point then $\nu(q)\le r$. $\nu(q)=r$ implies $\gamma(q)=r$.
\item If $q$ is a 3 point then $\nu(q)\le r-2$.
\item $\overline S_r(\overline U_1)$ makes SNCs with $\overline B_2(\overline U_1)$.
\end{enumerate}
There exists a sequence of quadratic transforms $W_p\rightarrow V_p$ such that if
$Z_p\rightarrow W_p$ is the sequence of monodial transforms (in any order) centered
 at the strict transforms of
curves $C$ in $\overline S_r(X)$
then 
$$
Z_p\rightarrow W_p\rightarrow V_p\rightarrow Y_p\rightarrow\text{spec}(R)
$$
extends to a permissible sequence of monodial transforms
$$
\overline\pi:\overline U_2\rightarrow \overline U_1\rightarrow \overline U\rightarrow U
$$
over an affine neighborhood of $p$ such that 
$\overline U_2$ contains no 2 curves $D$ such that 
$D$ is r-1 big or r small.
$\overline S_r(\overline U_2)$
makes SNCs with $\overline B_2(\overline U_2)$, and 
if $q\in \overline \pi^{-1}(p)$,
\begin{description}
\item[1'.] $\nu(q)\le r$ if $q$ is a 1 or 2 point. $\nu(q)=r$ implies $\gamma(q)=r$.
\item[2'.] If $q$ is a 2 point and $\nu(q)=r-1$, then either $\tau(q)>0$ or $\gamma(q)=r$ or
 $\tau(q)=0$ and (\ref{eq64}) holds at $q$ with $0<d_i<i$, $e_i=i$ and
$\overline S_{r-1}(Y_1)$ contains a single curve $D$ containing $q$, and containing a 1 point, which has local equations
$x=z=0$ at $q$.
\item[3'.] $\nu(q)\le r-2$ if $q$ is a 3 point.
\end{description}

If $p\not\in\overline S_r(X)$, then $Z=Y_{n_0}$.
\end{Theorem}

\begin{pf}
The proof of Theorem \ref{Theorem28} applied to the conclusions of Theorem \ref{Theorem30} proves this theorem.
\end{pf}

\begin{Theorem}\label{Theorem1019} 
Suppose that $r=2$, $A_2(X)$ holds, $p$ is a 3 point with $\nu(p)=r-1=1$, and
we have permissible parameters $(x,y,z)$ at $p$ for $u,v$ (with $y,z\in {\cal O}_{X,p}$)
 such that 
$$
\begin{array}{ll}
u&= (x^by^{a+nb}z^{a+(n+1)b})^m\\
v&=P(x^by^{a+nb}z^{a+(n+1)b})+x^dy^{c+n(d+2)}z^{c+(n+1)(d+2)}F_p\\
F_p&=\tau x+a_1(y,z)
\end{array}
$$
with $n\ge 0$, $a,b>0$,
 $L(x,0,0)\ne 0$ so that $\tau$ is a unit and 
$$
a(d+r-1)-bc=a(d+1)-bc=0
$$
or 
$$
\begin{array}{ll}
u&=(x^by^az^b)^m\\
v&=P(x^by^az^b)+x^dy^cz^{d+2}F_q\\
F_p&=\tau x+a_1(y,z)
\end{array}
$$
with $n\ge 0$, $a,b>0$,
 $L(x,0,0)\ne 0$ so that $\tau$ is a unit and 
$$
a(d+r-1)-bc=a(d+1)-bc=0
$$

Let $R ={\cal O}_{X,p}$.

 Suppose that  $\pi:Y_p\rightarrow \text{Spec}(R)$ is the morphism of
Theorem \ref{Theorem30}.

Let 
$$
\cdots\rightarrow Y_n\rightarrow\cdots\rightarrow Y_1\rightarrow Y_p
$$
be a sequence of permissible monodial transforms centered at 2 curves $D$
such that  $D$ is 1 big.  Then there exists $n_0<\infty$ such that 
$$
V_p=Y_{n_0}\stackrel{\pi_1}{\rightarrow} Y_p\rightarrow\text{spec}(R)
$$
extends to a permissible sequence of monodial transforms 
$$
\overline U_1\rightarrow\overline U\rightarrow U
$$
over an affine neighborhood $U$ of $p$ (with the notation of Theorem \ref{Theorem30})
such that $\overline U_1$ contains no 2 curves $D$ such that 
$D$ is 1 big or 2 small, and for $q\in\overline U_1$,
\begin{enumerate}
\item If $q$ is a 1 or a 2 point then $\nu(q)\le 2$. $\nu(q)=2$ implies $\gamma(q)=2$.
\item If $q$ is a 3 point then $\nu(q)=0$.
\item $\overline S_2(\overline U_1)$ makes SNCs with $\overline B_2(\overline U_1)$.
\end{enumerate}
There exists a sequence of quadratic transforms $W_p\rightarrow V_p$ such that if
$Z_p\rightarrow W_p$ is the sequence of monodial transforms (in any order) centered
 at the strict transforms $C'$ of
curves $C$ in $\overline S_2(X)$, followed by  monodial transforms centered at any 2 curves
$\overline C$ which are sections over $C'$ such that $\overline C$ is 1 big,
then 
$$
Z_p\rightarrow W_p\rightarrow V_p\rightarrow Y_p\rightarrow\text{spec}(R)
$$
extends to a permissible sequence of monodial transforms
$$
\overline\pi:\overline U_2\rightarrow \overline U_1\rightarrow \overline U\rightarrow U
$$
over an affine neighborhood of $p$ such that 
$\overline U_2$ contains no 2 curves $D$ such that 
$D$ is 1 big or 2 small.
$\overline S_2(\overline U_2)$
makes SNCs with $\overline B_2(\overline U_2)$, and 
if $q\in \overline\pi^{-1}(p)$,
\begin{description}
\item[1'.] $\nu(q)\le 2$ if $q$ is a 1 or 2 point. $\nu(q)=2$ implies $\gamma(q)=2$.
\item[2'.] If $q$ is a 2 point and $\nu(q)=1$, then either $q$ is resolved or $\gamma(q)=2$.
 \item[3'.] $\nu(q)=0$ if $q$ is a 3 point.
\end{description}

If $p\not\in\overline S_2(X)$, then $Z_p=Y_{n_0}$.
\end{Theorem}

\begin{pf} The proof of Theorem \ref{Theorem1018} applied to the conclusions of Theorem
\ref{Theorem30} proves the Theorem.
\end{pf}

\section{Resolution 1}

Throughout this section we will assume that $\Phi_X:X \rightarrow S$ is weakly prepared.

In this chapter we will need to consider the following condition on a 2 point 
 $p\in X$  such that $\nu(p)=r$ and $\tau(p)=1$. The condition is that  
$\Phi_X(p)$ has permissible parameters $(u,v)$ such that $u=0$ is a local
equation of $E_X$ at $p$ and $p$ has permissible parameters $(x,y,z)$ for $(u,v)$
such that(eq998)
\begin{equation}\label{eq998}
\begin{array}{ll}
u&=(x^ay^b)^m\\
v &= P(x^ay^b)+x^cy^dF_p\\
\end{array}
\end{equation}
and $L_p$ contains a nonzero $y^{r-1}z$ term with  $a(d+r-1)-bc=0$.
Up to interchanging $x$ and $y$, this condition is independent of permissible
parameters at $p$ for $(u,v)$.

\begin{Lemma}\label{Lemma962}
Suppose that $X$ satisfies $A_r(X)$, with $r\ge 2$, and $C$ is a 2 curve on $X$ such that $C\subset
\overline S_r(X)$. Suppose that
$$
\cdots\rightarrow X_n\rightarrow \cdots \rightarrow X_1\rightarrow X
$$
is a sequence of permissible monodial transforms centered at 2 curves $C_i$ such that
$C_i\subset\overline S_r(X_i)$  are sections over $C$. Then this sequence is finite.
That is, there exists $n<\infty$
such that $X_n$ contains no 2 curve $C_n$ with this property.
\end{Lemma}

\begin{pf} Since $A_r(X)$ holds, $C$ must be r small.
Suppose that $q\in C$ is a 2 point and the sequence has infinite length. Let $q_n$ be the point on $X_n$ which is the intersection 
of the fiber over $q$ and $C_n$.  With the notations of (\ref{eq957})
in the proof of Lemma \ref{Lemma500}, there are permissible
parameters $(x,y,z)$ at $q$ such that
$$
\begin{array}{ll}
u&=(x^ay^b)^m\\
F_q&=\overline czy^{r-1}+\Sigma_{i+j\ge r,k\ge 0} c_{ijk}x^iy^jz^k.
\end{array}
$$
For all $n$ there are permissible parameters $(x_n,y_n,z)$ at $q_n$ such that
$$
x=x_n, y=x_n^ny_n
$$
$$
F_{q_n}=\frac{F_q}{x_1^{n(r-1)}}
$$
If the sequence has infinite length, then $c_{ijk}=0$ if $j<r-1$, so that $y\mid F_q$,
a contradiction to the assumption that $F_q$ is normalized.
\end{pf}

\begin{Lemma}\label{Lemma1021} Suppose that $X$ satisfies $A_r(X)$ with 
$r\ge 2$ and $C$ is a 2 curve on $X$ such that $C$ is r-1 big. Suppose that
$$
\cdots\rightarrow X_n\rightarrow\cdots\rightarrow X_1\rightarrow X
$$
is a sequence of permissible monodial transforms, centered at 2 curves $C_i$ such that
$C_i$ is a section over $C$ and $C_i$ is r-1 big.
Then the sequence is finite. That is, there exists $n<\infty$ such that $X_n$
contains no 2 curves $C_n$ with this property.
\end{Lemma}

\begin{pf}
Suppose that $p\in C$ is a 2 point such that $\nu(p)=r-1$. $p$ has permissible
parameters $(x,y,z)$ such that
$$
\begin{array}{ll}
u&=(x^ay^b)^m\\
v&=P(x^ay^b)+x^cy^dF_p\\
F_p&=\sum_{i+j\ge r-1}a_{ij}(z)x^iy^j
\end{array}
$$
Let $\pi:Y\rightarrow\text{spec}({\cal O}_{X,p})$ be the blowup of $C$.

Suppose that $q\in\pi^{-1}(p)$ is the 2 point with permissible parameters
$(x_1,y_1,z_1)$
$$
x=x_1y_1, y=y_1
$$
\begin{equation}\label{eq616}
\begin{array}{ll}
u&=(x_1^ay_1^{a+b})^m\\
v&=P(x_1^ay_1^{a+b})+x_1^cy_1^{c+d+r-1}(\sum_{i+j=r-1}a_{ij}(z)x_1^i+y_1\Omega)
\end{array}
\end{equation}
$\nu(q)\le r-2$ unless $L_p=x^{r-1}$. Then $\nu(q)=r-1$.

Suppose that $L_p=x^{r-1}$. 
$$
F_p=\tau x^{r-1}+\sum_{i=1}^{r-1}y^{b_i}a_i(y,z)x^{r-1-i}
$$
where $\tau$ is a unit and $x\not\,\mid a_i$, $y\not\,\mid a_i$, $b_i\ge i$ for all $i$.

Suppose that the 2 curve $C_1\subset Y$ containing $q$ is such that $C_1$ is r-1 big. $C_1 = V(x_1,y_1)$ in (\ref{eq616}).
$$
F_q=\tau x_1^{r-1}+\sum_{i=1}^{r-1}y_1^{b_i-i}a_i(y_1,z)x_1^{r-1-i}.
$$

By induction on $b_i$, after a finite number of blowups of 2 curves, we reach
$\lambda: Z\rightarrow\text{spec}({\cal O}_{X,p})$ such that if $D$ is a 2
curve in $Z$ which is a section over $C$, then $D$ is not r-1 big.
\end{pf}

\begin{Definition}
Suppose that $X$ satisfies $A_r(X)$ with $r\ge 2$. A 2 point $p\in X$ contained in a 2 curve $C$ is called
bad if $\nu(p)=r$, $\tau(p)=1$ and one of the following holds.
\begin{enumerate}
\item $C\not\subset \overline S_r(X)$.
\item $C\subset\overline S_r(X)$ is r small and there exists a sequence of monodial transforms 
$$
X_n\rightarrow X_{n-1}\rightarrow\cdots\rightarrow X_1\rightarrow\text{spec}({\cal O}_{X,p})
$$
and 2 curves $C_i\subset X_i$
which are sections over $C$ such that 
$C_i\subset\overline S_r(X_i)$ is r small for $i<n$,  $C_n\not\subset \overline S_r(X_i)$,
$X_{i+1}\rightarrow X_i$ is centered at $C_i$ if $i<n$, and if
$p_n$ is the point on $C_n$ over $p$ then $\nu(p_n)=r$.
\item There exists a curve $D\subset\overline S_r(X)$ such that $D$ contains a 1 point
$p\in D$, and $D$ is r big at $p$.
\end{enumerate}
\end{Definition}

Suppose that $r\ge 2$ and $A_r(X)$ holds. Then there are only finitely many bad 2 points on $X$.

\begin{Lemma}\label{Lemma663} Suppose that $X$ satisfies $\overline A_r(X)$
with $r\ge 2$. Then
there exists a sequence of quadratic transforms $X_1\rightarrow X$ such that
$A_r(X_1)$ holds.
\end{Lemma}

\begin{pf} Let $\pi:X_1\rightarrow X$ be a sequence of quadratic transforms so that
the strict transform of $\overline S_r(X)$ makes SNCs with $\overline B_2(X)$. Then
$\overline A_r(X_1)$ holds  by Theorems \ref{Theorem9} and \ref{Theorem13},
and $A_r(X_1)$ holds by Lemma \ref{Lemma54} and Theorem \ref{Theorem12}.
\end{pf}

\begin{Lemma}\label{Lemma34}
Suppose that $A_r(X)$ holds with $r\ge 2$ and $p\in X$ is a bad 2 point. Suppose that there does not exist
a curve $D\subset\overline S_r(X)$ such that $p\in D$, $D$ contains a 1 point and $D$ is r big at $p$.  Then there exists
a sequence of quadratic transforms $\pi:X_1\rightarrow X$ centered at 2 points over $p$ such that
$A_r(X_1)$ holds and all 2 points $q\in\pi^{-1}(p)$ are good. 
\end{Lemma}

\begin{pf}
 There exist permissible parameters $(x,y,z)$ at $p$ such that
$$
\begin{array}{ll}
u&=(x^ay^b)^m\\
v&=P(x^ay^b)+x^cy^dF_p\\
F_p&=\sum_{i+j+k\ge r}a_{ijk}x^iy^jz^k
\end{array}
$$
Let $\pi_1:X_1\rightarrow X$ be the blowup of $p$. By Theorems \ref{Theorem9}, \ref{Theorem13}
and Lemma \ref{Lemma54}, $A_r(X_1)$ holds.
Suppose that $q\in \pi^{-1}(p)$ is a 2 point such that  $\nu(q)=r$ and $\tau(q)=1$. After a permissible
change of parameters at $p$, we may assume that $q$ has permissible parameters $(\overline x_1,
\overline y_1,\overline z_1)$ such that 
$x=\overline x_1\overline y_1, y=\overline y_1, z=\overline y_1\overline z_1$.
$\tau(q)=1$ and $\nu(q)=r$  implies that, after replacing $z$ by a constant times
$z$, that
$$
L_p=L(x,z)=\overline d x^r+x^{r-1}z
$$
for some $\overline d\in k$.

Suppose that there exists a 2 point $q'\in\pi^{-1}(p)$ such that $\nu(q')=r$ and $q'$
has permissible parameters $x',y',z'$ such that $x=x', y=x'y', z=x'(z'+\alpha)$
for some $\alpha\in k$. Then there exists a form $L_1$ and $\overline c\in k$ such that
$$
L_p=\left\{\begin{array}{ll} L_1+\overline c x^{\overline a}y^{\overline b}&\text{ if there exists }
\overline a,\overline b\in{\bold N}\text{ such that }\overline a+\overline b=r\\
&\text{ and }a(d+\overline b)-b(c+\overline a)=0\\
L_1&\text{otherwise.}
\end{array}\right.
$$
where
$$
L_1=L_1(y,z-\alpha x)=ey^r+f(z-\alpha x)y^{r-1}
$$
for some $e,f\in k$, with $f\ne 0$. This is not possible, since $r\ge 2$.
Thus all 2 points $q_1\in\pi^{-1}(p)$ with $\nu(q_1)=r$ have permissible
parameters $(x_1,y_1,z_1)$ such that 
$$
x=x_1y_1, y=y_1, z=y_1(z_1+\alpha)
$$
for some $\alpha\in k$. There exist at most finitely many bad 2 points $q_1\in\pi_1^{-1}(p)$.

Consider the following sequence of quadratic transforms
$$
\cdots\rightarrow X_n\rightarrow X_{n-1}\rightarrow\cdots\rightarrow X_1\rightarrow X
$$
with maps $\lambda_n:X_n\rightarrow X$,
where $\pi_i:X_i\rightarrow X_{i-1}$ is the blowup of all bad 2 points in $\lambda_{i-1}^{-1}(p)$.
We will show that there exists $n<\infty$ such that $\lambda_n^{-1}(p)$ contains no
bad 2 points. Suppose not. Then there exist bad 2 points $q_i\in X_i$ such that
$\pi_i(q_i)=q_{i-1}$ for all $i$.

$q_1$ has permissible parameters $(x_1,y_1,z_1)$ such that 
$$
x=x_1y_1, y=y_1, z=y_1(z_1+\alpha_1)
$$
for some $\alpha_1\in k$. $\nu(q_1)=r$ implies
$$
L_{q_1}=\overline d_1 x_1^r+x_1^{r-1}z_1+y_1\Omega_1
$$
for some $\overline d_1\in k$ and series $\Omega_1$. Since $\nu(q_2)=r$, and because of the existence of the $x_1^{r-1}z_1$ term in $L_{q_1}$, $q_2$ must have permissible parameters $(x_2,y_2,z_2)$ such that
$$
x_1=x_2y_2, y_1=y_2, z_1=y_2(z_2+\alpha_2)
$$
for some $\alpha_2\in k$. $\nu(q_2)=r$ implies
$$
L_{q_2}=\overline d_2x_2^r+x_2^{r-1}z_2+y_2\Omega_2.
$$
We see that there exists a series $\sigma(y)=\sum_{i=1}^{\infty}\alpha_iy^i$
such that if we replace $z$ with $\tilde z=z-\sigma(y)$, we have permissible
parameters $(x_n,y_n,\tilde z_n)$ at $q_n$ such that
$$
x=x_ny_n^n, y=y_n, \tilde z=\tilde z_ny_n^n.
$$
Then
$$
F_{q_n} = \frac{F_p(x_ny_n^n,y_n,y_n^n\tilde z_n)}{y_n^{rn}}
$$ 
for all $n$,
so that  $a_{ijk}=0$ if $i+k<r$ and $F_p\in (x,\tilde z)^r$. By Lemma \ref{Lemma302}, 
$\hat{\cal I}_{\overline S_r(X),p}\subset (x,\tilde z)$, so that
$x=\tilde z=0$ are local equations of a curve $D\subset \overline S_r(X)$, since
$A_r(X)$ holds. $D$ is r big at $p$ by Lemma \ref{Lemma653}.
\end{pf}

\begin{Theorem}\label{Theorem665} Suppose that $A_r(X)$ holds
with $r\ge 2$.
 Then there exists a sequence of permissible monodial
transforms $X_1\rightarrow X$ such that 
the following properties hold:
\begin{enumerate}
\item $A_r(X_1)$ holds. 
\item All bad 2 points  $p\in X_1$ satisfy (\ref{eq998}).
\item  Suppose that $D\subset\overline S_r(X)$ is a curve which is r big at a 1 point. Then there
exists at most one 3 point $q\in D$ and $D$ has a tangent direction at $q$ distinct from
those of $\overline B_2(X)$ at $q$. Furthermore, if $D$ is not r big, there exists only one 2 point
$q\in D$. If $C$ is the 2 curve containing $q$, then $C$ is not r-1 big or r small.
\item If $C$ is a r small or r-1 big 2 curve containing a 2 point $p$ such that $p\in D$
where $D$ is a curve containing a 1 point and $D$ is r big at $p$, then $D$ is r big.
\end{enumerate}
\end{Theorem}

\begin{pf}
By Theorems \ref{Theorem9} and \ref{Theorem12}, there exists a sequence of quadratic transforms
$\pi_1:X_1\rightarrow X$ centered at 3 points so that $A_r(X_1)$ holds and if $q\in X_1$ is a bad 2 point such that
(\ref{eq998}) doesn't hold, and there exists a curve $D\subset\overline S_r(X_1)$ such that
$q\in D$ and $D$ is r big at $q$, then $D$ is r big. All exceptional 2 points for $\pi_1$ which
are bad must satisfy (\ref{eq998}).

Let $\pi_2:X_2\rightarrow X_1$ be the blowup of such a $D$. By Lemma \ref{Lemma654}, $A_r(X_2)$ holds and all 2 points in $\pi_2^{-1}(D)$ are good. We  have that if $q\in X_2$ is a bad 2 point such that (\ref{eq998}) doesn't hold, and there exists a curve
$D\subset \overline S_r(X_2)$ such that $q\in D$ and $D$ is r big at $q$, then $D$ is r big.
If such a $D$ exists, let $\pi_3:X_3\rightarrow X_2$ be the blowup of $D$.

After a finite sequence of blowups, we then obtain $\lambda_1:Z_1\rightarrow X$ such that
$A_r(Z_1)$ holds, and if $q\in Z_1$ is a bad 2 point which doesn't satisfy (\ref{eq998}),
then there doesn't exist a curve $D\subset \overline S_r(Z_1)$ such that $D$ is r big at $q$.
By Lemma \ref{Lemma34}, there exists a sequence of quadratic transforms $\lambda_2:Z_2\rightarrow Z_1$ such that $A_r(Z_2)$ holds, and if $q\in Z_2$ is a bad
2 point, then (\ref{eq998}) holds at $q$.

By Theorems \ref{Theorem9} and \ref{Theorem12}, there exists a sequence of
quadratic transforms $\lambda_3:Z_2\rightarrow Z_2$ centered at 3 points such that the conclusions of the
Theorem hold on $Z_3$. 
\end{pf}

\begin{Theorem}\label{Theorem666} Suppose that the conclusions of Theorem
\ref{Theorem665} hold on $X$. Then there exists a finite sequence of quadratic transforms
centered at 3 points $X_1\rightarrow X$ such that 
\begin{enumerate}
\item $A_r(X_1)$ holds.
\item All bad 2 points $p\in X_1$  satisfy (\ref{eq998}).
\item Suppose that $D\subset\overline S_r(X)$ is a curve which is r big at a 1 point. Then there
exists at most one 3 point $q\in D$ and $D$ has a tangent direction at $q$ distinct from
those of $\overline B_2(X)$ at $q$. If $D$ is not r big, there exists only one 2 point
$q\in D$. If $C$ is the 2 curve containing $q$, then $C$ is not r-1 big or r small.
\item If $C$ is a r small or r-1 big 2 curve containing a 2 point $p$ such that $p\in D$
where $D$ is a curve r big at $p$, then $D$ is r big.
\item If $q\in X_1$ is a 3 point with $\nu(q)=r-1$, then 
 either there are permissible parameters $(x,y,z)$ at $q$ such that 
\begin{equation}\label{eq954}
L_q\text{ depends on both $y$ and $z$ and }F_q\in (y,z)^{r-1}
\end{equation}
or
there are permissible parameters $(x,y,z)$
 at $q$
such that 
\begin{equation}\label{eq668}
F_q = \tau y^{r-1}+\sum_{j=1}^{r-1} a_j(x,z)x^{\alpha_j}z^{\beta_j}y^{r-1-j}
\end{equation}
where $\tau$ is a unit, $a_j$ are units (or zero), $\alpha_j+\beta_j\ge j$ for all $j$, and there exists $i$ such that
$$
\frac{\alpha_i}{i}\le \frac{\alpha_j}{j},\frac{\beta_i}{i}\le\frac{\beta_j}{j}
$$
for all $j$, and
$$
\left\{\frac{\alpha_i}{i}\right\}+\left\{\frac{\beta_i}{i}\right\}<1.
$$
\end{enumerate}
\end{Theorem}

\begin{pf}
X satisfies 1. - 4. of the conclusions of the Theorem. By Theorems \ref{Theorem9} and
\ref{Theorem12}, 1. - 4. are stable under quadratic transforms centered at 3 points.

Suppose that $\pi_1:X_1\rightarrow X$ is the blow up of a 3 point $p$. 

If $L_p$ depends on all three variables $x,y,z$ then
$\nu(q)\le r-2$ for all 3 points $q\in \pi^{-1}(p)$.

Suppose that  $L_p$ depends  on both $y$ and $z$.  Then
$\nu(q)\le r-2$ for all 3 points $q\in \pi^{-1}(p)$, except possibly under the quadratic transform
$$
x=x_1, y=x_1y_1, z=x_1z_1.
$$
At this 3 point $q$, 
$$
F_{q} = \frac{F_p}{x_1^{r-1}} = L_p(y_1,z_1)+x_1\Omega.
$$
By a sequence of quadratic transforms centered at 3 points, we can get the Theorem to hold above $p$, except
 possibly along an infinite
sequence 
$$
R = {\cal O}_{X,p} \rightarrow R_1 \rightarrow \cdots \rightarrow R_n\rightarrow \cdots
$$
where for all $n$ $R_n$ has permissible parameters $(x_n,y_n,z_n)$ with
$$
x=x_n, y=x_n^ny_n, z=x_n^nz_n
$$
and $\nu(\frac{F_{p}(x_n,x_n^ny_n,x_n^nz_n)}{x^{n(r-1)}})=r-1$. Thus $F_p\in (y,z)^{r-1}$.

Now suppose that $L_p$ depends only on $y$. Then 
\begin{equation}\label{eq955}
F_p= \tau y^{r-1} + \sum_{i=1}^{r-1} a_i(x,z)y^{r-1-i}
\end{equation}
where $\tau$ is a unit.

If $p_1\in\pi_1^{-1}(p)$ is a 3 point with $\nu(q)=r-1$, then $p_1$ has permissible parameters
$(x_1,y_1,z_1)$ of one of the following 2 forms:
$$
x=x_1,y=x_1y_1,z=x_1z_1
$$
or
$$
x=x_1z_1, y=y_1z_1, z=z_1
$$
and 
\begin{equation}\label{eq956}
F_{p_1}=\tau y_1^{r-1}+\sum_{i=1}^{r-1}\frac{a_i(x_1,x_1z_1)}{x_1^i}y_1^{r-1-i}
\end{equation}
or
$$
F_{p_1}=\tau y_1^{r-1}+\sum_{i=1}^{r-1}\frac{a_i(x_1z_1,z_1)}{z_1^i}y_1^{r-1-i}.
$$
Suppose that
$$
\cdots\rightarrow X_n\rightarrow X_{n-1}\rightarrow \cdots\rightarrow X_1\rightarrow X
$$
is a sequence of quadratic transforms, $\pi_i:X_i\rightarrow X_{i-1}$ with induced
maps $\lambda_i:X_i\rightarrow X$ such that  for all $i$, $\pi_i$ is the blowup of
a 3 point $p_{i-1}$, with $\nu(p_{i-1})=r-1$ and $\pi_{i-1}(p_{i-1})=p_{i-2}$.

We will show that there exists $n$ such that $p_n$ satisfies (\ref{eq668}).
Each $p_i$ has permissible parameters $(x_i,y_i.z_i)$ such that either
$$
x_{i-1}=x_i, y_{i-1}=x_iy_i, z_{i-1}=x_iz_i
$$
or 
$$
x_{i-1}=x_iz_i, y_{i-1}=y_iz_i, z_{i-1}=z_i.
$$
By (\ref{eq956}), and resolution of plane curve singularities, there exists $n_0$ such that
$n\ge n_0$ implies $F_{p_n}$ has the form
$$
F_{p_n}=\tau y_n^{r-1}+\sum_{i=1}^{r-1}a_{n,i}(x_n,z_n)x_n^{\alpha_{ni}}z_n^{\beta_{ni}}
y_n^{r-1-i}
$$
where $a_{n,i}(x_n,z_n)$ are units, and either
$$
(\alpha_{n+1,i},\beta_{n+1,i})=(\alpha_{n,i}+\beta_{n,i}-i,\beta_{n,i})
$$
or
$$
(\alpha_{n+1,i},\beta_{n+1,i})=(\alpha_{n,i},\alpha_{n,i}+\beta_{n,i}-i)
$$
 for all $i$.
The proof now follows from Lemmas \ref{Lemma966} and \ref{Lemma967}.

\end{pf}

\begin{Theorem}\label{TheoremE2} Suppose that $A_r(X)$ holds 
with $r\ge2$ and $p\in X$ is a 2 point such that
(\ref{eq998}) holds.
 Then either 
\begin{enumerate}
\item There exists a sequence of quadratic transforms $\pi:Y\rightarrow X$ centered at
points over $p$ such that
\begin{enumerate}
\item If $q\in\pi^{-1}(p)$ is a 1 point then $\nu(q)\le r$. $\nu(q)=r$ implies 
$\gamma(q)=r$.
\item If $q\in\pi^{-1}(p)$ is a 2 point then $\nu(p)\le r$. $\nu(p)=r$ implies
$\tau(p)\ge 2$.
\item $q\in\pi^{-1}(p)$ a 3 point implies $\nu(q)\le r-1$. $\nu(q)=r-1$
implies $q$ satisfies the assumptions of (\ref{eq509}) and (\ref{eq510})
of Theorem \ref{Theorem30}. If $D_q$ is the 2 curve containing $q$, with local
equations $y=z=0$ at $q$, in the notation of Theorem \ref{Theorem30}, then
$F_{q'}$ is resolved for all $q\ne q'\in D_q$.
\item $A_r(Y)$ holds.
\end{enumerate}
or
\item There exists a curve $C\subset\overline S_r(X)$ which is r big at $p$.
Then there exists an affine neighborhood $U$ of $p$ such
that the blowup of $C\cap U$,
$\pi:Z\rightarrow U$ is a permissible monodial transform
such that 
\begin{enumerate}
\item If $q\in\pi^{-1}(p)$ is a 2 point then $\nu(q)=0$.
\item If $q\in\pi^{-1}(p)$ is the 3 point then either $\nu(q)\le r-2$ or
$q$ satisfies the assumptions of
(\ref{eq511})and (\ref{eq512}) of Theorem \ref{Theorem30}.
$D_q=\pi^{-1}(p)$ is the 2 curve with local equations $y=z=0$ at $q$ in the notation of
 Theorem \ref{Theorem30}.
\item $A_r(Z)$ holds.
\end{enumerate}
\end{enumerate}
\end{Theorem}

\begin{pf} We first assume that the assumption of 2. doesn't hold.
There are permissible parameters $(x,y,z)$ at $p$ such that
\begin{equation}\label{eq626}
\begin{array}{ll}
u&=(x^ay^b)^m\\
v&=P(x^ay^b)+x^cy^dF_p.
\end{array}
\end{equation}
By assumption, $L_p$ has the form 
\begin{equation}\label{eq627}
L_p=f(x,y)+zg(x,y)
\end{equation}
and $L_p$ 
contains a $y^{r-1}z$ term with 
\begin{equation}\label{eq628} 
a(d+r-1)-bc=0.
\end{equation}

Let $\pi_1:X_1\rightarrow X$ be the blowup of $p$.

Suppose that there exists a 2 point $q\in\pi_1^{-1}(p)$ such that $\nu(q)=r$ and
$q$ has permissible parameters $(x_1,y_1,z_1)$ such that 
$$
x=x_1y_1, y=y_1, z=y_1(z_1+\alpha)
$$
After a permissible change of parameters, we may assume that $\alpha=0$. Then $F_q=\frac{F_p}{y_1^r}$ and $\nu(F_q(0,0,z_1))\le 1$, so that $q$ is resolved.

Suppose that $p_1\in \pi_1^{-1}(p)$ is a 2 point such that $\nu(p_1)=r$ and $p_1$ has
permissible parameters $(x_1,y_1,z_1)$ such that
$$
x=x_1, y=x_1y_1, z=x_1(z_1+\alpha).
$$
After making a permissible change of parameters, we may assume that $\alpha=0$.
Then $L_p$ depends only on $y$ and $z$, and 
$$
L_p=\overline ey^r+\overline bzy^{r-1}
$$
for some $\overline b,\overline e\in k$ with $\overline b\ne 0$.

Suppose that $q\in\pi_1^{-1}(p)$ is another 2 point such that $\nu(q)=r$, and $q$
has permissible paramters $(x_1,y_1,z_1)$ such that
$$
x=x_1,y=x_1y_1,z=x_1(z_1+\alpha)
$$
with $\alpha\ne 0$.

Then there exists a form $L$ such that 
$$
L_p=\left\{\begin{array}{l}
L(y,z-\alpha x)\text{ or}\\
L(y,z-\alpha x)+\overline cx^{\overline \alpha}y^{\overline \beta}\\
\text{ for some $\overline c\in k$, such that 
$a(d+\overline \beta)-b(c+\overline \alpha)=0, \overline \alpha+\overline \beta=r$.}
\end{array}\right.
$$
$$
L_p=\overline ey^{r}+\overline b\alpha xy^{r-1}+\overline b(z-\alpha x)y^{r-1}
$$
implies $L(y,z-\alpha x)=L_p-\overline b\alpha x y^{r-1}$, but

 $\overline \alpha=1,\overline \beta=r-1$ is not possible, since 
\begin{equation}\label{eq629}
a(d+r-1)-b(c+1)=-b\ne 0
\end{equation}

 Thus $p_1$ is the unique 2 point $q\in \pi_1^{-1}(p)$
with $\nu(q)=r$. There are permissible parameters $(x_1,y_1,z_1)$ at $p_1$ such that
$$
x=x_1,
y=x_1y_1,
z=x_1z_1
$$
and $L_p=L(y,z)=\overline e y^r+\overline by^{r-1}z$, $L_{p_1}=L(y_1,z_1)+x_1\Omega$. 

If $p'\in \pi^{-1}(p)$ is the 3 point, then $\nu(p')\le r-1$.
If $p'\in\pi^{-1}(p)$ is a 1 point then $\nu(p')\le r$ and $\nu(p')=r$ implies $\gamma(p')=r$
by Theorem \ref{Theorem13}.

By Theorem \ref{Theorem13}, $\tau(p_1)\ge 1$.
Suppose that $\nu(p_1)=r$ and $\tau(p_1)=1$. Let $\pi_2:X_2\rightarrow X_1$
be the blowup of $p_1$.

We can make an analysis of $\pi_2^{-1}(p_1)$ which is similar to that of $\pi_1^{-1}(p)$.
(\ref{eq626}) is replaced with
$$
\begin{array}{ll}
u&=(x^{a+b}y^b)^m\\
v&=P(x^{a+b}y^b)+x^{c+d+r}y^dF
\end{array}
$$
$\tau(p_1)=1$ implies $L_{p_1}$ has the form of (\ref{eq627}). (\ref{eq628}) holds
at $p_1$. (\ref{eq629})
is then modified to:
 $\overline \alpha=1,\overline \beta=r-1$ is not possible since
$$
\begin{array}{l}
(a+b)(d+r-1)-b(c+d+r+1)\\
=a(d+r-1)+b(d+r)-b-bc-b(d+r)-b=-2b\ne 0
\end{array}
$$

We conclude that there is at most one 2 point $p_2\in \pi_2^{-1}(p_1)$ with $\nu(p_2)=r$, and after replacing $z$ with $z-\alpha x^2$ for some $\alpha\in k$, we have that $p_2$
has permissible parameters $(x_2,y_2,z_2)$ such that
$$
x=x_2, y=x_2^2y_2, z=x_2^2z_2.
$$

Suppose that we can construct an infinite sequence of quadratic transforms 
$$
\cdots\rightarrow X_n\rightarrow \cdots \rightarrow X_1\rightarrow X
$$
Where $\pi_n:X_{n+1}\rightarrow X_n$ is the blowup of a 2 point $p_n \in X_n$ over $p$ 
such that $\nu(p_n)=r$ and $\tau(p_n)=1$. 
Then there exists a series $\sigma(x)$ such that if we make a formal change of variables, replacing
$z$ with $z-\sigma(x)$, we get that there are permissible parameters   $(x_n,y_n,z_n)$
at $p_n$
such that
$$
x=x_1,
y=x_1^ny_1,
z=x_1^nz_1
$$

$p_n$ must be the only 2 point of $\pi_n^{-1}(p_{n-1})$ such that $\nu(p_n)=r$
and $\tau(p_n)=1$. To show this, the argument following (\ref{eq626}) is
modified by replacing (\ref{eq626}) with
$$
\begin{array}{ll}
u&=(x^{a+nb}y^b)^m\\
v&=P(x^{a+nb}y^b)+x^{c+n(d+r)}y^dF
\end{array}
$$
where (\ref{eq628}) holds at $p_n$.
(\ref{eq629}) is modified to:

$\overline \alpha=1,\overline \beta=r-1$ is not possible since
$$
\begin{array}{l}
(a+nb)(d+r-1)-b(c+n(d+r)+1)\\
=a(d+r-1)+nb(d+r-1)-bc-nb(d+r)-b\\
=-nb-b=-(n+1)b\ne0
\end{array}
$$

$\nu(p_n)=\nu(\frac{F_p}{x_1^{nr}})=r$ for all $n$ implies that $F_p\in (y,z)^r$.
$\hat{\cal I}_{\overline S_r(X),p}\subset (y,z)$ by Lemma \ref{Lemma302}, so that
since $A_r(X)$ holds, $y=z=0$ are local equations at $p$ of a curve $C\subset\overline S_r(X)$,
which is r big at $p$ by Lemma \ref{Lemma653},  a contradiction to the assumption that the assumption of 2. doesn't hold.
Thus there exists a sequence of quadratic transforms 
$\pi:Y\rightarrow X$ such that if $q\in \pi^{-1}(p)$ is a 2 point, then $\nu(q)\le r$,
and $\nu(q)=r$ implies $\tau(q)>1$. By Theorem \ref{Theorem13} and Lemma \ref{Lemma54}, $A_r(Y)$ holds.

Suppose that $q\in \pi^{-1}(p)$ is a 3 point such that $\nu(q)=r-1$. 
There 
exist permissible parameters $(x_2,y_2,z_2)$ at $q$, and there exists a 2 point $p_n\in X_n$ such that $\nu(p_n)=r$, $\tau(p_n)=1$ and 
$p_n$ has permissible parameters $(x_1,y_1,z_1)$ such that 
\begin{equation}\label{eq1020}
x=x_1,
y=x_1^ny_1,
z=x_1^nz_1 
\end{equation}
$$
x_1=x_2z_2,
y_1=y_2z_2,
z_1=z_2
$$
\begin{equation}\label{eq649}
\begin{array}{ll}
u&=(x_1^{a+nb}y_1^b)^m\\
v&=P(x_1^{a+nb}y_1^b)+x_1^{c+n(d+r)}y_1^d(\overline by_1^{r-1}z_1+\overline e y_1^r+x_1\Omega)
\end{array}
\end{equation}
with $\overline b\ne 0$.
\begin{equation}\label{eq98}
\begin{array}{ll}
u&=(x_2^{a+nb}y_2^bz_2^{a+(n+1)b})^m\\
v&=P(x_2^{a+nb}y_2^bz_2^{a+(n+1)b})+ x_2^{c+n(d+r)}y_2^dz_2^{c+(n+1)(d+r)}(\overline by_2^{r-1}+\cdots)
\end{array}
\end{equation}
with $a(d+r-1)-bc=0$. Thus $q$ satisfies the assumptions of (\ref{eq509})
and (\ref{eq510}) of 
Theorem \ref{Theorem30}, with $x=y_2,y=x_2,z=z_2$.

On $Y$ we have:
\begin{enumerate}
\item If $q\in\pi^{-1}(p)$ is a 1 point then $\nu(p)\le r$.
$\nu(p)=r$ implies $\gamma(p)=r$.
\item If $q\in \pi^{-1}(p)$ is a 2 point then $\nu(p)\le r$.
$\nu(p)=r$ implies $\tau(p)\ge 2$.
\item $q\in \pi^{-1}(p)$ a 3 point implies $\nu(q)\le r-1$. $\nu(q)=r-1$ implies
$q$ satisfies the assumptions of (\ref{eq509}) and (\ref{eq510}) of Theorem \ref{Theorem30}
\item $A_r(Y)$ holds.
\end{enumerate}

Let $T$ be the set of 3 points $q\in\pi^{-1}(p)$ such that $\nu(q)=r-1$, so that
$q$ satisfies (\ref{eq509}) and (\ref{eq510}) of Theorem \ref{Theorem30}.

Suppose that $q\in T$.
In the factorization $Y\rightarrow X$ by quadratic transforms there exists a factorization
$Y\rightarrow X_n\rightarrow X$ such that $q$ is an exceptional point on the blowup of a
2 point $p_n$ of the form of (\ref{eq1020}) on $X_n$.
Let $\tau:Y\rightarrow X_n$ be this map, $D_q$ be the nonsingular  curve 
$D_q\subset\tau^{-1}(p_n)$
such that $\hat{\cal I}_{D_q,q}=(x_2,z_2)$. The other points $q'\in D_q$ have regular parameters
$(x',y',z')$ with the notation of (\ref{eq1020}) such that 
$$
x_1=x'y',
y_1=y',
z_1=y'(z'+\alpha)
$$
with $\alpha\in k$. 
At such a 2 point $q'$, we have
$\nu(F_{q'}(0,0,z'))\le 1$ by (\ref{eq649}). Thus $q'$ is resolved.
\vskip .2truein

We now prove 2. There are permissible parameters $(x,y,z)$ at $p$ such that $y=z=0$ are local equations of $C$ at $p$ and
$$
\begin{array}{ll}
u&=(x^ay^b)^m\\
v&=P(x^ay^b)+x^cy^dF_p.
\end{array}
$$
There exists $\overline a\in k$ such that 
$$
L_p=\overline a y^r+zy^{r-1}
$$
Let $\pi_1:Y\rightarrow\text{spec}({\cal O}_{X,p})$ be the
blowup of $C$. Suppose that $q\in\pi_1^{-1}(p)$ is a 2 point.
$q$ has permissible parameters $(x,y_1,z_1)$ such that
$$
y=y_1,z=y_1(z_1+\alpha)
$$
for some $\alpha\in k$.
$$
\begin{array}{ll}
u&=(x^ay_1^b)^m\\
v&=p_1(x^ay_1^b)+x^cy_1^{d+r}(z_1+\alpha'+x\Omega'+y_1\Omega'')
\end{array}
$$
with $\alpha'\in k$, so that $q$ is resolved. At the 3 point $q\in\pi^{-1}(p)$,
there are permissible parameters $(x,y_1,z_1)$ such that 
$$
y=y_1z_1, z=z_1.
$$
$$
\begin{array}{ll}
u&=(x^ay_1^bz_1^b)^m\\
v&=P(x^ay_1^bz_1^b)+x^cy_1^dz_1^{d+r}F_q\\
F_q&=y_1^{r-1}+\overline a y_1^r+x\Omega'+z\Omega''
\end{array}
$$
with $a(d+r-1)-bc=0$. Either $\nu(q)\le r-2$ or $\nu(q)=r-1$ and $q$ satisfies the
assumptions of (\ref{eq511}) and (\ref{eq512}) of Theorem \ref{Theorem30}  (with $x=y_1, y=x_1, z=z_1$).

 The curve $C$ blown up in
Theorem \ref{Theorem30} is the fiber  $\pi^{-1}(p)$, which is resolved away from $q$. 
There exists an affine neighborhood $U$ of $q$ such that if $Z\rightarrow U$ is the
blowup of $C\cap U$, then
$A_r(Z)$ holds by Lemma \ref{Lemma654}.
\end{pf}

\begin{Theorem}\label{Theorem667} Suppose that $A_r(X)$ holds with $r\ge 2$.
Then there exists a finite sequence of permissible monodial transforms
$X_1\rightarrow X$ such that 
\begin{enumerate}
\item $A_r(X_1)$ holds.
\item  If $p\in X_1$ is a 3 point, then $\nu(p)\le r-2$.
\item If $p\in X_1$ is a 2 point such that $\nu(p)=r$ and $\tau(p)=1$ then 
$p$ has permissible parameters $(x,y,z)$
such that
\begin{equation}\label{eq87}
\begin{array}{ll}
u&=(x^ay^b)^m\\
v &= P(x^ay^b)+x^cy^dF_q\\
\end{array}
\end{equation}
and $L_p$ contains a nonzero $y^{r-1}z$ term with  $a(d+r-1)-bc=0$.
\item $\overline S_r(X_1)$ makes SNCs with $\overline B_2(X_1)$.
\item If $C$ is a 2 curve on $X_1$, then $C$ is not r small or r-1 big.
\end{enumerate}
\end{Theorem}

\begin{pf} We may assume that the conclusions of Theorem \ref{Theorem666} hold on $X$.
Let $\{D_1,\ldots, D_n\}$ be the curves $D$ in $X$ which intersect a r-1 big or r small
2 curve at a 2 point such that $D$ is r big there. By assumption, $D_1,\ldots, D_n$ are r big.

Let $\sigma_1:W_1\rightarrow X$ be the blowup of $D_1$. By Theorem \ref{TheoremE2}
and Lemma \ref{Lemma654} $A_r(Z_1)$ holds. $\sigma_1^{-1}(D_1)$ contains no bad 2 points.
If $q_1\in\sigma_1^{-1}(D_1)$ is a 3 point with $\nu(q_1)=r-1$, then $q_1\in\sigma_1^{-1}(q)$
where $q$ is a bad 2 point. In this case, 2. (b) of Theorem \ref{TheoremE2} holds at $q_1$.

Let $\{\overline D_2,\ldots, \overline D_n\}$ be the strict transforms of $\{D_2,\ldots,
D_n\}$ on $W_1$. These curves are all r big, and are the curves $D$ in $Z_1$ which
intersect a r-1 big or r small curve at a 2 point such that $D$ is r big there.

We can blowup successively the strict transforms of $\overline D_2,\ldots, \overline D_n$
by a map $\lambda:W\rightarrow X$ to get a $W$ such that $A_r(W)$ holds, the exceptional
locus of $\lambda$ contains no bad 2 points, and if $q$ is an exceptional 3 point with
$\nu(q)=r-1$ then $q$ must satisfy 2. (b) of Theorem \ref{TheoremE2}.

Furthermore, if $C$ is an r-1 big or r small 2 curve on $W$, and $p\in C$ is a bad 2 point, then
there does not exist a curve $D\subset\overline S_r(W)$ such that $D$ is r big at $p$.

By Theorem \ref{TheoremE2}, after performing a sequence of quadratic transforms $X_1\rightarrow
Z$ over bad 2 points $p\in X$ such that  if $C$ is the 2 curve
containing $p\in X$ then $C$ is r-1 big or r small, we have
\begin{description}
\item[1'.] The conclusions of 1. - 2. of Theorem \ref{Theorem666} hold.
\item[2.'] Suppose that $C$ is a 2 curve such that $C$ is r-1 big or r small.
If $p\in C$ is a 2 point, then $p$ is good.
\item[3'.] If $p$ is a 3 point such that $\nu(p)=r-1$, then either
\begin{description}
\item[(a)] $p$ satisfies the assumptions of (\ref{eq509}) and (\ref{eq510}) of Theorem
\ref{Theorem30}. If $D_p$ is the 2 curve containing $p$, with local equations
$y=z=0$ at $p$, in the notation of Theorem \ref{Theorem30}, then $F_q$ is resolved for
all $p\ne q\in D_p$. or
\item[(b)] $p$ satisfies the assumptions of (\ref{eq511}) and (\ref{eq512}) of Theorem
\ref{Theorem30}. If $D_p$ is the 2 curve containing $p$, with local equations
$y=z=0$ at $p$, in the notation of Theorem \ref{Theorem30}, then $F_q$ is resolved for
all $p\ne q\in D_p$. or
\item[(c)] $p$ satisfies (\ref{eq954}) or (\ref{eq668}) of Theorem \ref{Theorem666}.
\end{description}
\end{description}

Let $\{p_1,\ldots,p_n\}$ be the 3 points of $X_1$ which satisfy 3'. (a) or 3' (b). By Theorem 
\ref{Theorem30}, there exist sequences of permissible monodial transforms, over
sections of $D_{p_i}$ for $1\le i\le n$,
$$
\lambda_{p_i}:Y_{p_i}\rightarrow \text{spec}({\cal O}_{X_1,p_i})
$$
such that the conclusions of Theorem \ref{Theorem30} hold.

Since $D_{p_i}$ is resolved at points $p_i\ne q\in D_{p_i}$, the only obstruction to 
extending $\lambda_{p_i}$ to a permissible sequence of monodial transforms of sections
over $D_{p_i}$ in $X_1$ is if the corresponding sections over $D_{p_i}$
in $X_1$ do not make SNCs with 2 curves. This difficulty can be resolved by performing
quadratic transforms at the points where the section does not make SNCs with 2 curves,
since these points are necessarily resolved.

We can thus extend the 
$$
Y_{p_i}\rightarrow \text{spec}({\cal O}_{X_1,p_i})
$$
to a sequence of permissible monoidal transforms
$$
\lambda:Y\rightarrow X_1
$$
such that $Y\times_{X_1}\text{spec}({\cal O}_{X_1,p_i})\cong Y_{p_i}$
for $1\le i\le n$,  $A_r(Y)$ holds, 1'. and 2'. hold on $Y-\lambda^{-1}(\{p_1,\ldots,p_n\})$, and if
$q\in Y-\lambda^{-1}(\{p_1,\ldots,p_n\})$ is a 3 point such that $\nu(q)=r-1$, then 
$q$ satisfies (\ref{eq954}) or (\ref{eq668}) of Theorem \ref{Theorem666}.
Let
$$
\cdots\rightarrow Z_n\rightarrow \cdots \rightarrow Z_1\rightarrow Y
$$
be a sequence of permissible monodial transforms such that $Z_n\rightarrow Z_{n-1}$
is the blowup of a 2 curve $C$ such that $C$ is r-1 big or r small. 

We will show that there exists $n<\infty$ such that $Z_n$ does not contain a 2 curve
$C$ such that $C$ is r-1 big or r small, and that
$Z_n$ satisfies the conclusions of the Theorem.

By Theorem \ref{Theorem29}, this holds above a neighborhood of
$\lambda^{-1}(\{p_1,\ldots, p_n\})$. We must verify this condition over
$\overline Y=Y-\lambda^{-1}(\{p_1,\ldots,p_n\})$.

Suppose that $C$ is a 2 curve on $\overline Y$, such that $C$ is r-1 big or r small  and $p\in C$ is a 3 point.
Then all 3 points $q$ on $C\subset \overline Y$ have $\nu(q)=r-1$, 
and satisfy (\ref{eq954}) or (\ref{eq668}). 

 Let 
$\pi:Y_1\rightarrow \overline Y$ be the blowup of $C$.
The assumption that all 2 points of $C$ are good and Lemma \ref{Lemma500} imply that
$q\in\pi^{-1}(C)$ a 2 point implies $\nu(q)\le r$ and if $\nu(q)=r$ then either
$\gamma(q)=r$ or  $\tau(q)=1$,  $q$ is a good 2 point, $C\subset \overline S_r(X)$ and if
$\overline C$ is the 2 curve containing $q$, then $\overline C$ is a section over $C$ such
that $\overline C\subset\overline S_r(X)$.  Lemma \ref{Lemma500}, the assumption that all
2 points are good, and Lemma \ref{Lemma501}
imply $A_r(Y_1)$ holds. Further, by Lemma \ref{Lemma500}, if $\overline C\subset\pi^{-1}(C)$
is a 2 curve such that $\overline C$ is r-1 big or r small, then
$\overline C$ is a section over $C$. All 2 points of $\overline C$ are good points.

Suppose that $q\in C$ is a 3 point with permissible regular parameters $(x,y,z)$
such that $y=z=0$ are local equations of $C$ and
$$
\begin{array}{ll}
u&=(x^ay^bz^c)^m\\
v&=P(x^ay^bz^c)+x^dy^ez^fF_q.
\end{array}
$$
$F_q\in (y,z)^{r-1}$.

Suppose that $q_1\in \pi^{-1}(q)$, and $q_1$ has permissible parameters $(x,y_1,z_1)$
such that 
$$
y=y_1, z=y_1z_1
$$
If (\ref{eq954}) holds at $q$,
$$
F_{q_1}=\frac{F_q}{y_1^{r-1}}=L_q(1,z_1)+y_1\Omega+x\Lambda
$$
implies $\nu(q_1)\le r-2$.  If (\ref{eq668}) holds at $q$ we must have  $\beta_j\ge j$ for all $j$, and
$$
F_{q'}=\tau+z\Omega'
$$
so that $\nu(q_1)=0$.

Now suppose that $q_1\in\pi^{-1}(q)$ and $q_1$ has permissible parameters $(x,y_1,z_1)$
such that 
$$
y=y_1z_1, z=z_1
$$
If (\ref{eq954}) holds at $q$,
$$
F_{q_1}=\frac{F_q}{z_1^{r-1}}=L_q(y_1,1)+z_1\Omega+x\Lambda
$$
implies $\nu(q_1)\le r-2$.  If (\ref{eq668}) holds at $q$,
 we must have  $\beta_j\ge j$ for all $j$, and
$$
F_{q_1}=\tau y_1^{r-1}+\sum_{j=1}^{r-1} a_j(x,z_1)x^{\alpha_j}z_1^{\beta_j-j}y_1^{r-1-j}
$$
so that either $\nu(q_1)\le r-2$, or $q_1$ has the form of (\ref{eq668}) also, but
with $\frac{\beta_i}{i}$ decreased by 1.

By Lemma \ref{Lemma962} and Lemma \ref{Lemma1021}, after a finite number of blowups of 2 curves $C$ such that $C$ is r small or r-1 big, we reach $\tilde Y\rightarrow
\overline Y$ such that $\tilde Y$ contains no 2 curves $C$ such that $C$ is r small or r-1 big. Since all 3 points $q$ of $\tilde Y$ with $\nu(q)=r-1$ must satisfy (\ref{eq954}) or (\ref{eq668}), which implies that $\alpha_j\ge j$ for all $j$ or $\beta_j\ge j$ for all $j$, so that
there exists a 2 curve through $q$ which is r small or r-1 big, we must have $\nu(q)\le r-2$ if $q\in \tilde Y$ is a 3 point.

\end{pf}

\section{Resolution 2}

Throughout this section we will assume that $\Phi_X:X\rightarrow S$ is weakly prepared.

We define a new condition on $X$
\begin{Definition}\label{Def1095}  Suppose that $r\ge 2$.
We will say that $C_r(X)$  holds  if:
\begin{enumerate}
\item If $p\in X$ is a 1 point then $\nu(p)\le r$. If $\nu(p)=r$ then $\gamma(p)=r$.
\item If $p$ is a 2 point then $\nu(p)\le r$. 
If $\nu(p)=r$ then $\gamma(p)=r$. If $\nu(p)=r-1$ then one of the following three cases
must hold: 
\begin{enumerate}
\item $\tau(p)>0$ or
\item  $\gamma(p)=r$ or
\item $r\ge 3$, $\nu(p)=r-1$, $\tau(p)=0$, $p\not\in\overline S_r(X)$, there exists a unique curve 
$D\subset\overline S_{r-1}(X)$  containing a 1 point such that $p\in D$, and permissible parameters
$(x,y,z)$ at $p$  such that 
$x=z=0$ are local equations of $D$,
\begin{equation}\label{eq641}
\begin{array}{ll}
u&=(x^ay^b)^m\\
v&=P(x^ay^b)+x^cy^dF_p\\
F_p&=\tau x^{r-1}+\sum_{j=1}^{r-1}\overline a_j(y,z)y^{d_j}z^{e_j}x^{r-1-j}
\end{array}
\end{equation}
where $\tau$ is a unit, $\overline a_j$ are units (or 0). There exists $i$
such that $\overline a_i\ne 0$, $e_i=i$, $0<d_i<i$,
$$
\frac{d_i}{i}\le \frac{d_j}{j},\frac{e_i}{i}\le \frac{e_j}{j}
$$
for all $j$, and
$$
\left\{\frac{d_i}{i}\right\}+\left\{\frac{e_i}{i}\right\}<1.
$$
\end{enumerate}
\item If $p$ is a 3 point then $\nu(p)\le r-2$.
\item $\overline S_r(X)$ makes SNCs with $\overline B_2(X)$.
 \end{enumerate}
\end{Definition}

\begin{Remark} If $C_r(X)$ holds then there does not exist a 2 curve $C$ on $X$ such
that $C$ is r small or r-1 big.
\end{Remark}

\begin{Theorem}\label{TheoremE1} Suppose that $r\ge 2$, $A_r(X)$ holds,
$p\in X$ is a 2 point such that
$\nu(p)=r$ and $2\le\tau(p)<r$,  then either
\begin{enumerate}
\item There exists a sequence of quadratic transforms $\pi:Y\rightarrow X$ over $p$ such that 
\begin{enumerate}
\item $A_r(Y)$ holds.
\item If $q\in \pi^{-1}(p)$ is a 1 point then $\nu(q)\le r$.
$\nu(q)=r$  implies $\gamma(q)=r$.
\item If $q\in\pi^{-1}(p)$ is a 2 point then $\nu(q)\le r-1$.
\item If $q\in\pi^{-1}(p)$ is a 3 point, then $\nu(q)\le r-2$.
\item If $D\subset\pi^{-1}(p)$ is a 2 curve,  then $D$ is not r small or r-1 big.
\end{enumerate} or
\item There exists a curve $C\subset \overline S_r(X)$ such that $p\in C$
 and $C$ is r big at $p$.
 There exists an affine neighborhood
$U$ of $p$ such that the blowup of $C\cap U$,
 $\pi:Y\rightarrow U$ is a permissible monodial transform such that
\begin{enumerate} 
\item $A_r(Y)$ holds.
\item If $q\in\pi^{-1}(p)$ is a 2 point, then $\nu(q)\le r-1$.
\item If $q\in \pi^{-1}(p)$ is the 3 point, then $\nu(q)\le r-2$.
\item The 2 curve $D=\pi^{-1}(p)$  is not r small or r-1 big.
\end{enumerate}
\end{enumerate}
In either case, if $X$ satisfies the conclusions of Theorem \ref{Theorem667}, then $Y$ satisfies the
conclusions of Theorem \ref{Theorem667}.
\end{Theorem}

\begin{pf} $p$ has permissible parameters $(x,y,z)$ such that 
 $$
\begin{array}{ll}
u&=(x^ay^b)^m\\
v&=P(x^ay^b)+x^cy^dF_p\\
F_p&=\sum_{i+j+k\ge r} a_{ijk}x^iy^jz^k
\end{array}
$$

Suppose that there does not exist a curve $C\subset\overline S_r(X)$ such that 
$C$ is r big at $p$.

Let $\pi:X_1\rightarrow X$ be the blowup of $p$. 
We will first show that (a), (b) and (d)  of 1. hold on $X_1$ and if $q\in\pi^{-1}(p)$ is a 2 point with $\nu(q)=r$ then $\tau(q)\ge \tau(p)$.
This follows from Theorem \ref{Theorem9}, Theorem \ref{Theorem13} and Lemma \ref{Lemma54}.
All exceptional 2 curves $D$ of $\pi$ contain a 3 point $q$ such that $\nu(q)\le r-2$.
(e) thus holds by Lemmas \ref{Lemma655} and \ref{Lemma4}.

By Lemma \ref{Lemma655} there are at most finitely many 2 points $q\in\pi^{-1}(p)$ such that $\nu(q)=r$.  Suppose that there exists a 2 point $q\in\pi^{-1}(p)$ and $\nu(q)=r$.
After a permissible change of parameters at $p$, we have permissible
parameters $(x_1,y_1,z_1)$ at $q$ such that $x=x_1,y=x_1y_1,z=x_1z_1$. 
$L_p=L_p(y,z)$ depends on both $y$ and $z$.

Supppose there also exists a 2 point $q'\in\pi^{-1}(p)$ such that 
$\nu(q')=r$ and $q'$ has permissible parameters $(x',y',z')$ such that
$$
x=x'y', y=y', z=y'(z'+\alpha)
$$
for some $\alpha\in k$. Then there exists a form $L(x,z-\alpha y)$
such that 
$$
L_p(y,z)=\left\{\begin{array}{ll}
L(x,z-\alpha y)+\overline cx^{\overline a}y^{ \overline b}
&\text{ if there exists }\overline a,\overline b\in {\bold N} \text{ such that }\\
&\overline a+\overline b=r, a(d+\overline b)-b(c+\overline a)=0\\
L(x,z-\alpha y)&\text{ otherwise}
\end{array}
\right.
$$
Thus 
$$
L_p=\overline d(z-\alpha y)^r+\overline cy^r
$$
for some $\overline d,\overline c\in k$ with $\overline d\ne 0$, a contradiction to the assumption that
$\tau(p)<r$. Let
$$
\cdots \rightarrow Y_n\rightarrow Y_{n-1}\rightarrow \cdots\rightarrow Y_1\rightarrow X
$$
be the sequence of quadratic transforms $\pi_n:Y_n\rightarrow Y_{n-1}$
constructed by blowing up all 2 points $q'$ on $Y_n$ which lie  over $p$ and have $\nu(q')=r$.

Suppose that this sequence has infinite length. Then there exists $q_n\in Y_n$ such that
$\pi_n(q_n)=q_{n-1}$ and $\nu(q_n)=r$ for all $n$. There exists a series
$\phi(x)=\sum \alpha_ix^i$ such that after replacing $z$ with $z-\phi(x)$, $q_n$ has permissible
parameters $(x_n,y_n,z_n)$ such that 
$$
x=x_n, y=x_n^ny_n, z=x_n^nz_n
$$
and 
$$
F_{q_n}=L_q(y_n,z_n)+x_n\Omega_n.
$$
$F_{q_n}=\frac{F_q}{x_n^{nr}}$ for all $n>0$ implies $F_q\in (y,z)^r$.

$\hat{\cal I}_{\overline S_r(X),p}\subset(y,z)$ by Lemma \ref{Lemma302}.
Since $\overline S_r(X)$ makes SNCs with $\overline B_2(X)$ at $p$, $y=z=0$ are local equations
at $p$ of a curve $C\subset\overline S_r(X)$.

Now suppose that there exists a curve $C\subset\overline S_r(X)$ such that $p\in C$
and $C$ is r big at $p$. There exists
an affine neighborhood $U$ of $p$ such that $C\cap U$ makes SNCs with $\overline B_2(U)$.
 Let 
$\pi:Y\rightarrow U$ be the blowup of $C\cap U$. There exist permissible parameters
$(x,y,z)$ at $p$ such that $y=z=0$ are local equations of $C$ at $p$,
$$
\begin{array}{ll}
u&=(x^ay^b)^m\\
F_p&=\sum_{i+j\ge r}a_{ij}(x)y^iz^j.
\end{array}
$$
At the 3 point  $q\in\pi^{-1}(p)$, there are permissible parameters $(x,y_1,z_1)$ 
such that
$$
y=y_1z_1,z=z_1
$$
$$
F_q=\frac{F_p}{z_1^r}=\sum_{i+j=r}a_{ij}(0)y_1^i+z_1G+x\Omega
$$
we have $\nu(q)\le r-2$, since $2\le\tau(p)$.

At a 2 point $q\in\pi^{-1}(p)$, after a permissible change of variables at $p$,
there exist permissible parameters $(x,y_1,z_1)$ at $q$ such that 
$$
y=y_1,z=y_1z_1
$$
$$
F_q=\frac{F_p}{y_1^r}=\sum_{i+j=r}a_{ij}(0)z_1^j+y_1G+x\Omega.
$$
$\nu(q)<r$ and $\gamma(q)<r$ since $\tau(p)<r$.
Furthermore, if $D=\pi^{-1}(p)$, then
$F_q\not\in\hat{\cal I}_{D,q}^{r-1}$. 
There exists a possibly smaller affine neighborhood $U$ of $p$ such that $A_r(Y)$ holds by
 Lemma \ref{Lemma654}.

\end{pf}

\begin{Theorem}\label{TheoremE3} Suppose that the conclusions of 
Theorem \ref{Theorem667} hold on $X$ with $r\ge 2$, $p\in X$ is a 2 point
with permissible parameters $(x,y,z)$ such that 
$$
\begin{array}{ll}
u&=(x^ay^b)^m\\
v&=P(x^ay^b)+x^cy^dF_p
\end{array}
$$
and $\nu(p)=r-1$,  $\tau(p)=0$, $L_p=f(x,y)$ depends on both $x$ and $y$.
Then 
there exists a sequence of quadratic transforms $\pi:Z\rightarrow X$ over $p$
such that
\begin{enumerate}
\item $q\in \pi^{-1}(p)$ a 1 point or a 2 point implies that $\nu(q)\le r$.
$\nu(q)=r$ implies $\gamma(q)=r$.
\item $q\in\pi^{-1}(p)$ a 2 point with $\nu(q)=r-1$ implies that $\tau(q)>0$  or 
$\gamma(q)=r$.
\item $q\in \pi^{-1}(p)$ a 3 point implies $\nu(q)\le r-2$.
\item The conclusions of Theorem \ref{Theorem667} hold on $Z$
\end{enumerate} 
\end{Theorem}

\begin{pf}
Let 
\begin{equation}\label{eq89}
\pi:X_1\rightarrow X
\end{equation}
be the blowup of $p$. 

If $p_1\in\pi^{-1}(p)$ is a 1 point then $\nu(p_1)\le r$ and $\nu(p_1)=r$ implies
$\gamma(p_1)=r$ by Theorem \ref{Theorem9}.
If $p_1\in \pi^{-1}(p)$ is a 2 point then we must have $\nu(p_1)\le r-2$,
by our assumption on $f$. Suppose that $p_1\in \pi^{-1}(p)$ is the 3 point. Then
$\nu(p_1)\le r-1$ and $p_1$ has permissible parameters $(x_1,y_1,z_1)$ such that
$$
x=x_1z_1,
y=y_1z_1,
z=z_1
$$
Suppose that $\nu(p_1)=r-1$.
Then 
\begin{equation}\label{eq90}
L_{p_1}=f(x_1,y_1)+z_1\Omega.
\end{equation}

Let 
$$
F_p = \sum_{i+j+k\ge r-1} a_{ijk}x^iy^jz^k.
$$
Suppose that we can construct an infinite sequence of quadratic transforms
$$
\cdots\rightarrow X_n\rightarrow \cdots\rightarrow X_1\rightarrow X
$$
where $X_{n+1}\rightarrow X_{n}$ is the blowup of a 3 point $p_n$ lying over $p_{n-1}$
with $\nu(p_n)=r-1$. Then $p_n$ has permissible parameters $(x_n,y_n,z_n)$
such that
$$
x=x_nz_n^n,
y=y_nz_n^n,
z=z_n
$$
and
$$
F_{p_n} = \frac{F_p}{z_n^{n(r-1)}} = 
\sum a_{ijk} x_n^iy_n^jz_n^{(k+ n(i+j-r+1))}.
$$
Thus $a_{ijk}=0$ if $i+j<r-1$, which implies that  $F_p\in (x,y)^{r-1}$,
a contradiction since the conclusions of Theorem \ref{Theorem667} hold.

Thus by Theorem \ref{Theorem9} and Lemma \ref{Lemma54} there exists a finite sequence of quadratic transforms
$$
\pi: X_m\rightarrow \cdots\rightarrow X_1\rightarrow X
$$
where $X_{n+1}\rightarrow X_{n}$ is the blowup of a 3 point $p_n$ lying over $p_{n-1}$
with $\nu(p_n)=r-1$, such that $A_r(X_m)$ holds,  $\nu(q)\le r-2$ if 
$q\in\pi^{-1}(p)$ is a 3 point, and if $q\in\pi^{-1}(p)$ is a 1 point then $\nu(q)\le r$,
$\nu(q)=r$ implies $\gamma(q)=r$.  Suppose that $C$ is a 2 curve which is exceptional
for $\pi$. Then $C$ is not r-1 big or r small since $C$ must 
contain a 3 point $q'$ with $\nu(q')\le r-2$. Suppose that $q\in\pi^{-1}(p)$ is a 2 point and $\nu(q)\ge r-1$.
Then there exists a largest $n$ such that $q$ maps to a 3 point
$p_n\in X_n$. The point $q$ is then a 2 point on  $X_{n+1}$. $p_n$ has permissible parameters
$(x_1,y_1,z_1)$ such that
$$
x=x_1z_1^n,
y=y_1z_1^n,
z=z_1
$$
By assumption, $\nu(p_n)=r-1$. Write
$$
f = \sum_{i+j=r-1}a_{ij}x^iy^j.
$$
We then have 

\begin{equation}\label{eq91}
\begin{array}{ll}
u&=(x_1^ay_1^bz_1^{n(a+b)})^m\\
v&=P(x_1^ay_1^bz_1^{n(a+b)})+x_1^cy_1^dz_1^{n(c+d+r-1)}F_{p_n}\\
F_{p_n}&=\frac{F_p}{z_1^{n(r-1)}}=\sum_{i+j=r-1}a_{ij}x_1^iy_1^j+
\sum_{i+j+k=r-1, k>0}b_{ijk}x_1^iy_1^jz_1^k+\Omega
\end{array}
\end{equation}
with $\nu(\Omega)\ge r$.

Since $q$ is a 2 point, $\hat{\cal O}_{X_{n+1},q}$ has regular parameters $(x_2,y_2,z_2)$ of one of the
following forms: 
\begin{equation}\label{eq92}
x_1=x_2,
y_1=x_2(y_2+\alpha),
z_1=x_2z_2
\end{equation}
with $\alpha\ne 0$, or 
\begin{equation}\label{eq93}
x_1=x_2,
y_1=x_2y_2,
z_1=x_2(z_2+\beta)
\end{equation}
with $\beta\ne 0$, or 
\begin{equation}\label{eq94}
x_1=x_2y_2,
y_1=y_2,
z_1=y_2(z_2+\beta)
\end{equation}
with $\beta\ne 0$.

First suppose that (\ref{eq93}) holds. (\ref{eq94}) is symmetrical, and the analysis of that case is the same.
Set 
$$
x_2=\overline x_2(z_2+\beta)^{-\frac{n}{n+1}}.
$$
$$
\begin{array}{ll}
u&=(\overline x_2^{(n+1)(a+b)}y_2^b)^m=(\overline x_2^{\overline a}y_2^{\overline b})^{\overline m}\\
v&= P_q(\overline x_2^{\overline a}y_2^{\overline b})+\overline x_2^{(n+1)(c+d+r-1)}y_2^dF_q\\
F_q&=
[\sum_{i+j=r-1}a_{ij}y_2^j\\
&+\sum_{k>0, i+j+k=r-1}b_{ijk}y_2^j(z_2+\beta)^k+\overline x_2\Omega^1]
-\sum c_i\overline x_2^{\overline a_i}y_2^{\overline b_i}
\end{array}
$$
where $(\overline a,\overline b)=1$,
$$
(n+1)(a+b)(\overline b_i+d)-b((n+1)(c+d+r-1)+\overline a_i)=0.
$$
If some $b_{ijk}\ne 0$ in (\ref{eq91}), we have $\nu(q)\le r-1$ and $\nu(q)=r-1$
implies that $\tau(q)>0$. So suppose that all $b_{ijk}=0$ in (\ref{eq91}).
If  $\nu(q)\ge r-1$, then we must have 
$$
\sum a_{ij}y_2^j = a_{i_0j_0}y_2^{j_0}+a_{0,r-1}y_2^{r-1}
$$
where $a_{0,r-1}$ could be zero, $0<i_0$, $a_{i_0,j_0}\ne0$ (since $f(x,y)$ depends on
$x$ and $y$) and 
$$
(n+1)(c+d+r-1)b-(n+1)(a+b)(j_0+d)=0.
$$
Thus 
\begin{equation}\label{eq96}
(c+r-1-j_0)b-a(d+j_0)=0
\end{equation}

We then have
$$
x^cy^df = a_{i_0j_0}x^{c+r-1-j_0}y^{d+j_0}+a_{0,r-1}x^cy^{d+r-1}
$$
which  is normalized, so that $a_{i_0,j_0}=0$ by (\ref{eq96}).
 This contradiction shows that we must have that $\nu(q)\le r-2$ in this case.

Suppose that (\ref{eq92}) holds. Substitute (\ref{eq92}) into (\ref{eq91}).
Set 
$$
x_2 = \overline x_2(y_2+\alpha)^{-\frac{b}{(n+1)(a+b)}}
$$
$$
\begin{array}{ll}
u&= (\overline x_2^{(n+1)}z_2^{n})^{m(a+b)}\\
v&=P(\overline x_2^{(n+1)(a+b)}z_2^{n(a+b)})+\overline x_2^{(n+1)(c+d+r-1)}z_2^{n(c+d+r-1)}G
\end{array}
$$
where 
$$
G = (y_2+\alpha)^{\lambda}\left(\sum_{i+j=r-1} a_{ij}(y_2+\alpha)^j + 
\sum_{k>0, i+j+k=r-1}b_{ijk}(y_2+\alpha)^jz_2^k
+\overline x_2\Omega_2\right),
$$
with 
$$
\lambda = -\frac{b}{a+b}(c+d+r-1)+d.
$$

The only term which can be removed from the first sum 
$$
(y_2+\alpha)^{\lambda}\left(\sum_{i+j=r-1} a_{ij}(y_2+\alpha)^j)\right)
$$ 
of $G$ in obtaining  $F_q$ 
is the constant term. Thus $\gamma(q)\le r$.
\end{pf}

\begin{Theorem}\label{TheoremE6} Suppose that $r\ge 2$ and the conclusions of Theorem \ref{Theorem667} hold on $X$,
so that if $C$ is a 2 curve, then $C$ is not r small or r-1 big.
 Suppose that $p\in X$ is a 1 or a 2 point and $D$ is a generic
curve through $p$ on a component of $E_X$. Then there exists a sequence of quadratic transforms  
centered over  a 
finite number of points on the strict transform of $D$,
but not in the fiber over $p$, $\pi:X_1\rightarrow X$, such that the following
conditions hold.
\begin{enumerate}
\item There exists a neighborhood $U$ of $D-p$ such that $C_r(\pi^{-1}(U))$ holds.
The case 2. (c) of $C_r$ does not occur in $\pi^{-1}(U)$.
\item Let $D'$ be the strict transform of $D$ on $X_1$. Suppose that $q\in D'-p$, and 
$(x,y,z)$ are permissible parameters at $q$ such that $x=z=0$ are local equations
of $D'$ at $q$. If $q$ is a 1 point then $\nu(F_q(0,y,0))=1$. If $q$ is a 2 point then
$\nu(q)=0$.
\item The conclusions of Theorem \ref{Theorem667} hold on $X_1$.
\end{enumerate}
\end{Theorem}

\begin{pf}
Suppose that $q\in D$ is a 1 point. Then we can find permissible parameters $(x,y,z)$ at $q$
such that $x=y=0$ are local equations of $D$ at $q$. The multiplicity
$$
\phi(q)=\nu(F_q(0,0,z))
$$
is independent of such permissible parameters at $q$. Furthermore, the set
$$
\{q\in D\cap (X-\overline B_2(X))\mid \phi(q)\ge2\}
$$
is Zariski closed in $D\cap(X-\overline B_2(X))$.
By Lemma \ref{Lemma970}, $F_q\not\in\hat{\cal I}_{D,q}$ if $q\in D$.
  At most 1 points $q$ on $D$,   $\phi(q)=1$. Thus there are at most a finite number of
points $q\in D-p$ such that the conclusions of the Theorem do not hold at $q$.
\vskip .2truein
{\bf 1)}
Suppose that $q\in D-p$ and $\nu(q)=r$. Then $q$ is a generic point on a curve
$C$ of $\overline S_r(X)$.   $q$ is a 1 point.
\vskip .2truein 
{\bf 1a)} Suppose that $C$ is r big.
Then there exist permissible parameters $(x,y,z)$ at $q$ such that   
$$
\begin{array}{ll}
u&=x^a\\
v&=P(x)+x^bF_q\\
F_q&=\sum_{i+k\ge r, j\ge 0}a_{ijk}x^iy^jz^k
\end{array}
$$
where $x=z=0$ are local equations of $C$ at $q$, $x=y=0$ are local equations of $D$ at $q$. 
 $\gamma(q)=r$ implies $a_{00r}\ne 0$.

Let $\pi:Y\rightarrow X$ be the blowup of $q$.
Then $\nu(q')=0$ if $q'$ is the point on the intersection of the strict transform
of $D$ and $\pi^{-1}(q)$. Points of $\pi^{-1}(q)$ satisfy the condition of $C_r$ by Theorem
\ref{Theorem13} and Lemma \ref{Lemma54}.

\vskip .2truein
{\bf 1b)} 
Suppose that $C$ is r small.  By Lemma
\ref{Lemma5},
$$
F_q=\sum_{i+j\ge r}a_{ij}(y)x^iz^j+\tau(y)x^{r-1}
$$
where $x=z=0$ are local equations of $C$ at $q$, $x=y=0$ are local equations of $D$ at $q$,
(with $\nu(\tau)\ge 1$). Since $q$ is a generic point of $C$, $\nu(\tau)=1$, and after a permissible change of parameters, we have $\tau=y$.
$\gamma(q)=r$ implies $a_{0r}(y)$ is a unit. Let $\pi:Y\rightarrow 
X$ be the blowup of $q$. Then $\nu(q')=0$ if $q'$ is the point
on the intersection of the strict transform of $D$ and $\pi^{-1}(q)$. Points of $\pi^{-1}(q)$ satisfy the condition of $C_r$ by Theorem \ref{Theorem13} and Lemma \ref{Lemma54}.
\vskip .2truein
{\bf 2)} Suppose that $q\in D-p$, $\nu(q)=r-1$ and the conclusions of the Theorem
do not hold at $q$.
\vskip .2truein
{\bf 2a)} Suppose that  $q$ is a 1 point and 
$r\ge 3$. Then $q$ is  a general point on a
curve $C$  in $\overline S_{r-1}(X)$. There are permissible parameters $(x,y,z)$ at $q$ such that 
$x=z=0$ are local equations of $D$ at $q$. 

Let 
$\pi_1:X_1\rightarrow X$ be the blowup of $p$.
Theorem \ref{Theorem9} implies $\nu(q')\le r-1$ for all $q'\in\pi_1^{-1}(q)$ and $q'\in\pi_1^{-1}(q)$ a 2 point with $\nu(q')=r-1$ implies $\tau(q')>0$.

At the 2 point $q_1\in\pi_1^{-1}(q)$ on the strict transform of $D$, there are permissible
parameters $(x_1,y_1,z_1)$ such that 
$$
x=x_1y_1,
y=y_1,
z=y_1z_1
$$

Suppose that $\nu(q_1)=r-1$. We must have $\tau(q_1)>0$. Let $\pi_1:X_2\rightarrow X_1$ be the blowup
of $q_1$. By Theorem \ref{Theorem13}, if $q'\in\pi^{-1}(q_1)$, then if $q'$ is a 
1 point $\nu(q')\le r-1$. If $q'$ is a 2 point, $\nu(q')\le r-1$, $\nu(q')=r-1$ implies
$\tau(q')>0$. $q'$ a 3 point implies $\nu(q')\le r-2$. Let $q_2\in\pi_1^{-1}(q)$ be the 2 
point on the
strict transform of $D$. There are permissible parameters $(x_2,y_2,z_2)$ at $q_2$ such that
$$
x_1=x_2y_2,
y_1=y_2,
z_1=y_2z_2.
$$
If $\nu(q_2)=r-1$, then $\tau(q_2)>0$.

Suppose that we can construct an infinite sequence of quadratic transforms
$$
\cdots\rightarrow X_n\rightarrow X_{n-1}\rightarrow \cdots\rightarrow X_1\rightarrow X
$$
centered at the point $q_n$ on the strict transform of $D$ over $q$ on $X_n$, where $q_n$ are blown up as long as $\nu(q_n)=r-1$.

By Theorem \ref{Theorem13}, all points $q'$ on $X_n$ lieing over $p$ satisfy
$\nu(q')\le r-1$, $\nu(q')\le r-2$ If $q'$ is a 3 point and if $q'$ is a 2 point with $\nu(q')=r-1$. Then $\tau(q')>0$.

Suppose that $\nu(q_n)=r-1$ for all $n$. Then $q_n$ has permissible parameters $(x_1,y_1,z_1)$
such that 
$$
x=x_1y_1^n,y=y_1,z=z_1y_1^n
$$

$$
F_q=\sum_{i+j+k\ge r-1}a_{ijk}x^iy^jz^k
$$
$$
F_{q_n}=\frac{F_q}{y_1^{n(r-1)}}=\sum_{i+j+k\ge r-1} a_{ijk}x_1^iy_1^{j+n(i+k-(r-1))}z_1^k
$$
implies $a_{ijk}=0$ if $i+k<r-1$, so that $F_q\in(x,z)^{r-1}$.  This is a contradiction since 
$F_q\not\in\hat{\cal I}_{D,q}$.

Thus after a finite sequence of quadratic transforms,
 $\pi:Z\rightarrow X$
the strict transform of $D$ intersects $\pi^{-1}(q)$ at a 2 point $q_1$ with
 $\nu(q_1)<r-1$, so we are in case 3) below.
\vskip .2truein
{\bf 2b)} 
Suppose that $q$ is a 2 point. Suppose that $C$ is the 2 curve through $q$.
 By Lemma \ref{Lemma6},
 our assumption that $C$ is not r-1 big, 
and since $q$ is a generic point of $C$, we have  $\tau(q)>0$.  There exist permissible parameters $(x,y,z)$ at $q$
such that 
\begin{equation}\label{eq1023}
\begin{array}{ll}
u&=(x^ay^b)^m\\
v&=P(x^ay^b)+x^cy^dF_q\\
F_q&=\sum_{i+j\ge r-1,k\ge 0}a_{ijk}x^iy^jz^k+zx^{i_0}y^{j_0}
\end{array}
\end{equation}
with $i_0+j_0=r-2$.
$x=z=0$  are local equations of $D$ at $q$ and $\tau(q)>0$.

Let $\pi:X_1\rightarrow X$ be the blowup
of $q$. Then $\nu(q_1)\le r-1$  at all points  $q_1\in\pi^{-1}(q)$,
$\nu(q_1)\le r-2$ if $q_1\in\pi^{-1}(q)$ is a 3 point
and all 2 points $q_1$ of 
$\pi^{-1}(q)$ with $\nu(q_1)=r-1$ satisfy $\tau(q_1)>0$ by Theorem \ref{Theorem13}.

The strict transform of $D$ intersects $\pi^{-1}(q)$ at a 2 point $q'$ such that 
$$
x=x_1y_1, y=y_1, z=y_1z_1
$$
$$
F_{q_1}=\frac{F_q}{y_1^{r-1}}=\sum_{i+j\ge r-1,k\ge 0}a_{ijk}x_1^iy_1^{(i+j+k)-(r-1)}z_1^k+z_1x_1^{i_0}
$$
$x_1=z_1=0$ are local equations of the strict transform of $D$ at $q'$. If $\nu(q')\le r-2$
we are in case 3). Otherwise, $q_1$ is a 2 point with $\nu(q_1)=r-1$ and $\tau(q_1)>0$
(so that $i_0=r-2$).
Let
$$
\cdots\rightarrow X_n\rightarrow X_{n-1}\rightarrow \cdots\rightarrow X_1\rightarrow X
$$
be the sequence of quadratic transforms centered at the point $q_n$ on the strict
transform of $D$ over $q$ on $X_n$ where points $q_n$ are blownup as long as
$\nu(q_n)=r-1$. By Theorem \ref{Theorem13}, all points $q'$ on $X_n$ lieing over
$p$ satisfy $\nu(q')\le r-1$, $\nu(q')\le r-2$ if $q'$ is a 3 point, and if $q'$
is a 2 point with $\nu(q')=r-1$, then $\tau(q')>0$.

Suppose that $\nu(q_n)=r-1$ for all $n$. $q_n$ has permissible parameters $(x_1,y_1,z_1)$
such that 
$$
\begin{array}{ll}
x=x_1y_1^n,
y=y_1,
z=z_1y_1^n
\end{array}
$$
$$
F_{q_n}=\frac{F_q}{y_1^{n(r-1)}}=\sum_{i+j\ge r-1,k\ge 0}a_{ijk}x_1^iy_1^{j+n(i+k-(r-1))}
z_i^k+z_1x_1^{r-2}
$$
Thus $i+k-(r-1)\ge 0$ whenever $a_{ijk}\ne 0$ and $F_p\in(x,z)^{r-1}$, so that $F_p\in\hat{\cal I}_{D,p}$, which is a contradiction.

 After a finite sequence
of quadratic transforms, $\pi':X'\rightarrow X$, the strict transform of $D$ thus intersects
$(\pi')^{-1}(q)$ at a 2 point $q'$ with $\nu(q')\le  r-2$, so the result follows from
Case 3).
\vskip .2truein
{\bf 3)}
Suppose that $q\in D-p$, $\nu(q)\le r-2$, and the conclusions of the Theorem do not
hold at $q$. $q$ is a 1 point or a 2 point and $D$ makes SNCs with
the 2 curve through $C$. The result then follows from a similar but slightly simpler
argument to that of Case 2, by Theorems \ref{Theorem9} and \ref{Theorem13}.
\end{pf}

\begin{Theorem}\label{TheoremE7}
 Suppose that $X$ satisfies the conclusions of Theorem \ref{Theorem667} with $r\ge 2$.
 Then there
exists a sequence of permissible monoidal transforms
$X_1\rightarrow X$ such that $C_r(X_1)$ holds. 
\end{Theorem}

\begin{pf}
Let $T$ be the finite set of 2 points $p$ on $X$ such that (\ref{eq87}) holds at $p$,  and
$p\not\in D$ for any r big curve $D$ which contains a 1 point.

By 1. of Theorem \ref{TheoremE2}, there exists a sequence of quadratic transforms
$\pi_0:X_0\rightarrow X$ centered over points $p\in T$ such that

\begin{enumerate}
\item If $p\in X_0$ is a 1 point then $\nu(p)\le r$, $\nu(p)=r$ implies $\gamma(p)=r$.
\item If $p\in X_0$ is a 2 point then $\nu(p)\le r$. $\nu(p)=r$ implies
$\tau(p)\ge 2$, or (\ref{eq87}) holds at $p$ and there exists an r big
curve $D\subset\overline S_r(X_0)$ containing $p$.
\item If $p\in X_0$ is a 3 point, then $\nu(p)\le r-1$. $\nu(p)=r-1$ implies $p$
satisfies the assumptions of (\ref{eq509}) and (\ref{eq510}) of Theorem \ref{Theorem30}.
If $D_p$ is the 2 curve containing $p$ with local equations $y=z=0$ at $p$ in the
notation of Theorem \ref{Theorem30}, then $F_q$ is resolved for all $p\ne q\in D_p$.
\item $A_r(X_0)$ holds.
\item If $C$ is a 2 curve on $X_0$, then $C$ is not r small. If $C$ is r-1 big,
  then $\nu(p)=r-1$ for all $p\in C$.
\end{enumerate}

 Let $T_1$ be the
3 points of $X_0$ which satisfy (\ref{eq509}) and (\ref{eq510}) of Theorem
\ref{Theorem30}.

For $p\in T_1$, let $\lambda_p:Y_p\rightarrow \text{spec}({\cal O}_{X_0,p})$
be the sequence of monodial transforms centered over sections of $D_p$
such that the conclusions of Theorem \ref{Theorem30} hold. By Theorem \ref{Theorem29}
(or Theorem \ref{Theorem1019} if $r=2$),
there exists a sequence of monodial transforms $V_p\rightarrow Y_p$
centered at 2 curves $C$ such that
$C$ is r-1 big so that $V_p$ satisfies the conclusions
of 1. - 3. of Theorem \ref{Theorem29} (or Theorem \ref{Theorem1019} if $r=2$). 

Since $D_p$ is resolved at all points $p\ne q\in D_p$, the only obstruction to
extending $\lambda_p$ to a permissible sequence of monodial transforms of sections over $D_p$
in $X_0$ is if the corresponding sections over $D_p$ do not make SNCs with 2 curves.
This difficulty can be removed by performing quadratic transforms at the (resolved) points
where the section does not make SNCs with the 2 curves.

By  5. above and Lemmas
 \ref{Lemma500},  \ref{Lemma501} and \ref{Lemma1021}, we can thus construct a sequence
of permissible monodial transforms $\pi_0':X_0'\rightarrow X_0$ such that 
$$
X_0'\times_{X_0}\text{spec}({\cal O}_{X_0,p})\cong V_p
$$
 for $p\in T_1$, 

and $X_0'$ satisfies:
\begin{enumerate}
\item $p\in X_0'$ a 1 point implies $\nu(p)\le r$. $\nu(p)=r$ implies $\gamma(p)=r$.
\item $p\in X_0'$ a 2 point implies $\nu(p)\le r$. $\nu(p)=r$ implies $\tau(p)\ge 2$
or (\ref{eq87}) holds at $p$, and there exists a r big curve $D\subset\overline S_r(X_0)$
containing $p$.
\item $p\in X_0'$ a 3 point implies $\nu(p)\le r-2$.
\item $\overline S_r(X_0')$ makes SNCSs with $\overline B_2(X_0')$.
\item If $C$ is a 2 curve on $X_0$, then $C$ is not r small or r-1 big.
\end{enumerate}

Let $\gamma_1,\ldots,\gamma_n$ be the r big curves in $\overline S_r(X_0')$.
Each $\gamma_i$ necessarily contains 
a 1 point.

Let $\pi:X_1\rightarrow X_0'$ be the sequence of monodial transforms (in any order)
centered at the (strict transforms of) $\gamma_1,\ldots,\gamma_n$.  

By Lemma \ref{Lemma654}, 2. of Theorem \ref{TheoremE1}, and 2. of Theorem
\ref{TheoremE2},
\begin{enumerate}
\item If $p\in X_1$ is a 1 point then $\nu(p)\le r$. $\nu(p)=r$ implies $\gamma(p)=r$.
\item If $p\in X_1$ is a 2 point then $\nu(p)\le r$.
$\nu(p)=r$ implies $\tau(p)\ge 2$. If $\nu(p)=r$ and $\tau(p)<r$, then $p$ does not lie on a r big curve $E$ in
$\overline S_r(X_1)$. 
\item $p\in X_1$ a 3 point implies $\nu(p)\le r-1$.
$\nu(p)=r-1$ implies  $p$ satisfies the assumptions of (\ref{eq511}) and (\ref{eq512}) of Theorem
\ref{Theorem30}. If $D_p$ is the 2 curve containing $p$ with local equations
$y=z=0$ at $p$ in the notation of Theorem \ref{Theorem30}, then $F_q$ is resolved for
all $p\ne q\in D_p$.
\item $A_r(X_1)$ holds.
\item If $C$ is a 2 curve on $X_0$, then $C$ is not r small or r-1 big. 
\item There are only finitely many 2 points $p\in X_1$ such that $\nu(p)=r$.
\end{enumerate}
\vskip .2truein
6. is a consequence of 5. Let $T_2$ be the set of 3 points on $X_1$ satisfying (\ref{eq511}) and
(\ref{eq512}) of Theorem \ref{Theorem30}.

 By Theorem \ref{TheoremE1}, there exists a 
sequence of quadratic transforms $\pi_2:X_2\rightarrow X_1$ centered over the 2 points $p$
of $X_1$ with $\nu(p)=r$ and $2\le \tau(p)<r$ such that
\begin{enumerate}
\item If $p\in X_2$ is a 1 or 2 point then $\nu(p)\le r$. If $\nu(p)=r$ then $\gamma(p)=r$.
\item If $p\in X_2$ is a 3 point then $\nu(p)\le r-1$. If $\nu(p)=r-1$, then
 $p\in T_2$.
\item $A_r(X_2)$ holds.
\item If $D\subset X_2$ is a 2 curve, then $D$ is not r small or r-1 big.
\end{enumerate}

For $p\in T_2$, let 
\begin{equation}\label{eq650}
\lambda^p:Y_p\rightarrow \text{spec}({\cal O}_{X_2,p})
\end{equation}
be a sequence of permissible monodial transforms over sections of $D_p$ such that the
conclusions of Theorem \ref{Theorem30} hold.

Since $D_p$ is resolved at all points $p\ne q\in D_p$, the only obstruction to 
extending $\lambda^p$ to a permissible sequence of monodial transforms of sections over $D_p$ in 
$X_2$ is if the corresponding sections over $D_p$ do not make SNCs with 2 curves. This
difficulty can be removed by performing quadratic transforms at the points where the section
does not make SNCs with the 2 curves.

By Theorems \ref{Theorem9} and \ref{Theorem13}, we can thus construct a permissible
sequence of monodial transforms $\pi_3:X_3\rightarrow X_2$ such that
\begin{enumerate}
\item If $p\in T_2$, then $X_3\times_{X_2}\text{spec}({\cal O}_{X_2,p})\cong Y_p$.
\item If $q\in X_3-\pi_3^{-1}(T_2)$, and $q$ is a 1 or 2 point then $\nu(q)\le r$.
If $\nu(q)=r$, then $\gamma(q)=r$.
\item If $q\in X_3-\pi_3^{-1}(T_2)$ and $q$ is a 3 point then $\nu(q)\le r-2$.
\item $A_r(X_3)$ holds. 
\item If $D\subset X_3-\pi_3^{-1}(T_2)$ is a 2 curve, then $D$ is not r small or r-1 big.
\end{enumerate}

By 5. and Theorem  \ref{Theorem29} (or Theorem \ref{Theorem1019} if $r=2$), we can perform a sequence of
permissible monodial transforms $\sigma:\overline Z_2\rightarrow X_3$ centered at
r-1 big 2 curves $C$  to get that
\begin{enumerate}
\item If $p\in\overline Z_2$ is a 1 or 2 point, then $\nu(p)\le r$.
$\nu(p)=r$ implies $\gamma(p)=r$.
\item If $p\in\overline Z_2$ is a 3 point, then $\nu(p)\le r-2$.
\item There are no 2 curves $C$ in $\overline Z_2$ which are r small or r-1 big.
\item $\overline S_r(\overline Z_2)$ makes SNCs with $\overline B_2(Z_2)$.
\end{enumerate}

Since 3. holds, there are only finitely many 2 points $\{q_1,\ldots, q_m\}$ on $\overline Z_2$ such that
$\nu(q_i)=r-1$ and $\tau(q_i)=0$.

By Theorem \ref{TheoremE3},
 we can perform a sequence of
quadratic transforms $\sigma_1:W_1\rightarrow \overline Z_2$
over  the finitely many 2 points $q_i$ 
in $\overline Z_2$ such that $\nu(q_i)=r-1$, $\tau(q_i)=0$ and $L_q$ depends on both
$x$ and $y$ (where $(x,y,z)$ are permissible parameters at $q_i$)  so that 
\begin{enumerate}
\item $\nu(q)\le r$ and $\nu(q)=r$ implies $\gamma(q)=r$ at 1 and 2 points of $W_1$.
\item $\nu(q)\le r-2$ at 3 points of $W_1$
\item If $q\in W_1$ is a 2 point with $\nu(q)=r-1$
and $\tau(q)=0$, then either $\gamma(q)=r$ or there exist permissible parameters $(x,y,z)$ at $q$ such that
$L_q$ depends only on $x$.
\item There are no 2 curves $C$ on $W_1$ which are r small or r-1 big.
\item $\overline S_r(W_1)$ makes SNCs with $\overline B_2(W_1)$.
\end{enumerate}

Over the (finitely many) points $\{a_1,\ldots,a_n\}$ of $W_1$ which are 2 points with $\nu(p)=r-1$, $\gamma(p)>r$,
 $\tau(p)=0$, and there exist permissible parameters $(x,y,z)$ at $a_i$ such that $L_{a_i}$ depends only on $x$, by Theorem \ref{Theorem27}, there exist 
sequences of permissible monodial transforms
$$
 Y_{a_i}\rightarrow \text{spec}({\cal O}_{W_1,a_i})
$$
where $Y_{a_i}\rightarrow\text{spec}({\cal O}_{W_1,a_i})$ is a sequence of
blowups of sections over a general curve through $a_i$, and 
satisfies the conclusions of Theorem \ref{Theorem27}.

By Theorem \ref{TheoremE6} there exists a sequence of permissible monodial transforms $\tilde \pi:W_2\rightarrow W_1$
such that 
$$
W_2\times_{W_1}\text{spec}({\cal O}_{W_1,a_i})\cong Y_{a_i}
$$
for all $i$, and  for $q\in\tilde\pi^{-1}(W_2-\{a_1,\ldots,a_n\})$,
\begin{enumerate}
\item $\nu(q)\le r$, $\nu(q)=r$ implies  $\gamma(q)\le r$ if $q$ is a 1 or 2 point.
\item $q$ a 2 point and $\nu(q)=r-1$ implies $\tau(q)>0$ or $\gamma(q)=r$.
\item $\nu(q)\le r-2$ if $q$ is a 3 point.
\item There are no 2 curves $C$ in $\tilde\pi^{-1}(W_2-\{a_1,\ldots,a_n\})$ 
which are r small or r-1 big.
\item $A_r(W_2)$ holds.
\end{enumerate}

Since all 2 curves $C\subset W_2$ which are  r-1 big must map to some $a_i$ by 
4., there exists a sequence of
permissible monodial transforms $\pi_3:W_3\rightarrow W_2$ centered at 2 curves $C$
which are   r-1 big  such that
$$
W_3-(\tilde\pi\circ\pi_3)^{-1}(\{a_,\ldots,a_n\})\cong W_2-\tilde\pi^{-1}(\{a_1,\ldots,a_n\})
$$
and $W_3\times_{W_1}\text{spec}({\cal O}_{W_1,a_i})$ satisfies the conclusions of $V_{a_i}$
of Theorem \ref{Theorem28} (or of $V_{a_i}$ of Theorem \ref{Theorem1018}).

Then the sequence of quadratic transforms $W_{a_i}\rightarrow V_{a_i}$ of
Theorem \ref{Theorem28} (or of $W_{a_i}\rightarrow V_{a_i}$ of Theorem \ref{Theorem1018}) extend to $\pi_4:W_4\rightarrow W_3$ such that
$$
W_4-(\tilde\pi\circ\pi_3\circ\pi_4)^{-1}(\{a_1,\ldots,a_n\})\cong W_2
-\tilde\pi^{-1}(\{a_1,\ldots,a_n\})
$$
and $W_4\times_{W_1}\text{spec}({\cal O}_{W_1,a_i})\cong W_{a_i}$ for $1\le i\le n$
in the notation of Theorem \ref{Theorem28} (or of  Theorem \ref{Theorem1018}).

Now assume that $r\ge 3$.
Let $\{D\}$ be the strict transform on $W_4$ of the curves $\{\overline D\}$ in $\overline S_r(W_1)$
which contain some $a_i$.   Each $\overline D$ contains a 1 point
and $\overline D$ is  r small
since $\overline D$ makes SNCs with $\overline B_2(W_1)$ and $\nu(a_i)=r-1$.
By Theorem \ref{Theorem28} and Lemma \ref{Lemma41},
and since by Lemma \ref{Lemma7} there does not exist a 2 point $q\in D$ such that
$\nu(q)=r-1$ and $\tau(q)>0$, there exists a finite sequence of quadratic transforms
$W_5\rightarrow W_4$ centered at points disjoint from any fiber over some $a_i$
such that if $W_6\rightarrow W_5$ is a sequence of monodial transforms centered at the
strict transforms of the $D$ then $C_r(W_6)$ holds.

 Now suppose that $r=2$. Let $\{D\}$ be the strict transforms on $W_4$ of the curves
$\{\overline D\}$ in 
$\overline S_r(W_1)$ which contain some $a_i$. Each $\overline D$ contains a 1 point and $\overline D$ is not r big. By Theorem \ref{Theorem1018}, Lemmas \ref{Lemma41} and \ref{Lemma97}
there exists a sequence of quadratic transforms $W_5\rightarrow W_4$ centered at points 
disjoint from any fiber over some $a_i$ such that if $W_6\rightarrow W_5$ is 
a sequence of monodial transforms centered at the strict transforms of the $D$, and then followed by a sequence of monodial transforms $W_7\rightarrow W_6$ centered at  the strict transforms
of 2 curves $C$ on $W_6$ which are sections over one of the $D$ blownup in $W_6\rightarrow
W_5$ and such that $C$ is 1 big, then $C_r(W_6)$ holds. 
\end{pf}

\begin{Theorem}\label{TheoremE8} Suppose that  $C_r(X)$ holds  with $r\ge 2$.
Then there exists a sequence of quadratic transforms
$\pi:X_1\rightarrow X$ such that  $C_r(X_1)$ holds and
if $C$ is a 2 curve on $X_1$ such that $C$ contains a 2 point $p$ with $\nu(p)=r$ and 
$p$ lies on a curve $D$ in $\overline S_r(X_1)$, then for $p\ne q\in C-B_3(X)$, $\nu(q)\le r-1$ and
$\nu(q)=r-1$ implies $\gamma(q)=r-1$.
\end{Theorem}

\begin{pf} 
Let $\pi:X_1\rightarrow X$ be the product of quadratic transforms centered at all 2 points
$q\in X$ such that $\nu(q)=r$ and $q$ is on a curve $D\subset\overline S_r(X)$.

 Suppose that $q\in X$ is a 2 point on a 2 curve $C$ such that $q\in D$,
for a curve $D\subset
\overline{S}_r(X)$. There exist permissible parameters $(x,y,z)$ at $q$ such that
$$
\begin{array}{ll}
u&=(x^ay^b)^m\\
v&=P(x^ay^b)+x^cy^dF_q
\end{array}
$$
where $x=z=0$ are local equations of $D$ at $q$.
Suppose that
$\nu(q)=r$ so that $\gamma(q)=r$. By Lemma \ref{Lemma17}, after a permissible change of variables,
$$
F_q=\tau z^r+\sum_{i=2}^ra_i(x,y)z^{r-i}
$$
with $\tau$ a unit, $\nu(a_i)\ge i$.
At a 1 point $q'\in\pi^{-1}(q)$
we have $\nu(q')\le r$, $\nu(q')=r$ implies $\gamma(q')=r$.

The 2 points $q'\in\pi^{-1}(q)$ have permissible parameters $(x_1,y_1,z_1)$ such that
$$
x=x_1y_1, y=y_1, z=y_1(z_1+\alpha)
$$
or 
$$
x=x_1, y=x_1y_1, z=x_1(z_1+\alpha).
$$
In either case
$$
F_{q'}=\tau(z_1+\alpha)^r+\text{ terms of order }\le r-2\text{ in }z_1
$$
implies $\gamma(q')\le r$, $\gamma(q')\le r-1$ if $\alpha\ne 0$.
Thus each exceptional curve $\overline C$ of $\pi$ contains at most one 2 point $q'$
such that $\nu(q')=r$.
At the 3 point $q'\in \pi^{-1}(q)$,
$$
x=x_1z_1,
y=y_1z_1, z=z_1
$$
and $\nu(q')=0$.

Thus by Lemma \ref{Lemma54}, $C_r(X_1)$ holds. 
Let $C'$ be the strict transform of $C$, $D'$ the strict transform of $D$ on $X_1$.
$C'$ and $D'$ are disjoint. $C'$ intersects $\pi_1^{-1}(q)$ at the 3 point $q'$ with 
$\nu(q')=0$.
 By Lemma \ref{Lemma54}
there is at most one curve $E$ in ${\overline S_r(X_1)}$ such that $E\subset \pi^{-1}(q)$,
and $E$ intersects each 2 curve in at most one point.
If $E$ intersects an exceptional 2 curve $\overline C$ in a point $q'$ such that $\nu(q')=r$,
then for $q'\ne q''\in\overline C$, $\gamma(q'')\le r-1$ (by the above analysis). 
Thus each exceptional 2 curve $\overline C$ for $\pi_1$ satisfies the conditions
of the conclusions of the Theorem.

The strict transform $C'$ of a 2 curve $C$ on $X_1$ contains no 2 points $q$ with $\nu(q)=r$ which are contained in a curve $D$ in $\overline S_r(X_1)$.

\end{pf}

\section{resolution 3 }

Throughout this section we will assume that $\Phi_X:X\rightarrow S$ is weakly prepared.
 
\begin{Lemma}\label{Lemma42}
 Suppose that
$C\subset X$ is a 2 curve. Suppose that $t$ is a natural number or $\infty$.
Then the set
$$
\{q\in C \vert\, q \text{ is a 2 point and }\gamma(q)\ge t\}
$$
is Zariski closed in $C-B_3(X)$. 
\end{Lemma}

\begin{pf} Suppose that $p\in C$ is a 2 point.
There exist permissible parameters $(x,y,z)$ at $p$ such that $(x,y,z)$
are uniformizing parameters in an \'etale neighborhood  $U$ of $p$ in $X$.
At $p$,
$$
\begin{array}{ll}
u&=(x^ay^b)^m\\
v&=P(x^ay^b)+x^cy^dF(x,y,z).
\end{array}
$$
Set 
$$
w=\frac{v-P_{\lambda}(x^ay^b)}{x^cy^d}
$$
with $\lambda>c+d$.  $w\in \Gamma(U,{\cal O}_X)$. If $q\in C\cap U$, there are permissible parameters $(x,y,z_q=z-\alpha)$ 
at $q$ for some $\alpha\in k$.
There exist $a_i(q)\in k$ such that
$$
F_q=w-\sum_{i=0}^{\infty}a_i(q)\frac{(x^ay^b)^i}{x^cy^d}.
$$
$$
\{q\in C\cap U\mid \nu(F_q(0,0,z_q)\ge t\}
=\left\{\begin{array}{ll}
\{q\in C\cap U\mid \frac{\partial^iw}{\partial z^i}(0,0,\alpha)=0, 0\le i<t\}
&\text{ if }ad-bc\ne 0\\
\{q\in C\cap U\mid \frac{\partial^iw}{\partial z^i}(0,0,\alpha)=0, 0<i<t\}
&\text{ if }ad-bc= 0
\end{array}
\right.
$$
is Zariski closed.

Since $U$ is an \'etale cover of an affine neighborhood $V$ of $p$,
$$
\{q\in C\mid q\text{ is a 2 point and }\gamma(q)\ge t\}\cap V
$$
is Zariski closed in $V\cap C$.
\end{pf}

\begin{Lemma}\label{Lemma43}
Suppose that $C$ is a 2 curve  and there exists $p\in C$ with
 permissible parameters $(x_p,y_p,z_p)$ at $p$ such that
$x_p=y_p=0$ are local equations of $C$ at $p$ and 
$\nu(F_p(0,0,z_p))<\infty$. If $q\in  C$  then 
$\nu(F_q(0,0,z_q))<\infty$, where 
$(x_q,y_q,z_q)$ are permissible parameters at $q$ and
$x_q=0, y_q=0$ are local equations for $C$ at $q$.
\end{Lemma}

\begin{pf}
If $\nu(F_q(0,0,z_q))=\infty$, then $F_q\in \hat{\cal I}_{\overline C,q}$ so that
$F_p\in\hat{\cal I}_{\overline C,p}$ for all $p\in C$ by Lemma \ref{Lemma655}.
 Thus $\nu(F_p(0,0,z_p))=\infty$
for all $p\in C$, a contradiction.
\end{pf}

\begin{Theorem}\label{Theorem45}
Suppose that $C_r(X)$ holds with $r\ge 2$  and
the conclusions of Theorem \ref{TheoremE8} hold on $X$.
Then there exists a sequence of permissible monoidal transforms  $\pi:Y\rightarrow X$
 centered at r big curves $C$ in
 $\overline S_r$ such that $C_r(Y)$ and the conclusions of Theorem \ref{TheoremE8} hold on $Y$ and if $D$ is a curve in 
$\overline{S_r(Y)}$, then $D$ is not r big. 
\end{Theorem}

\begin{pf}
Suppose that the $C\subset \overline S_r(X)$ is r big.
$C$ must contain a 1 point.
Let
$\pi:Y\rightarrow X$ be the blowup of $C$. 

By Lemma \ref{Lemma654}, $C_r(Y)$ holds and the conclusions of Theorem \ref{TheoremE8} hold
on $Y$. There is at most one curve $D\subset \overline S_r(Y)\cap \pi^{-1}(C)$. If this 
curve exists it must be a section over $C$.

 Let $p\in C$ be a 1 point.  As in (\ref{eq973}) of the proof of Lemma \ref{Lemma654},
there exist permissible parameters $(x,y,z)$ at $p$
such that 
$\hat{\cal I}_{C,p} =(x,z)$, 
\begin{equation}\label{eq651}
\begin{array}{ll}
u&=x^a\\
F_p&=\tau z^r+\sum_{i=2}^ra_i(x,y)z^{r-i}
\end{array}
\end{equation}
where $\tau$ is a unit, $x^i\mid a_i$ for $2\le i\le r$.

As shown in the proof of Lemma \ref{Lemma654},  the only point $q\in\pi^{-1}(p)$ which could be in $\overline S_r(Y)$ is the
1 point with permissible parameters $x=x_1, z=x_1z_1$. 
\begin{equation}\label{eq951}
\begin{array}{ll}
u&=x_1^a\\
F_q&=\tau z_1^r+\sum_{i=2}^r\frac{a_i(x_1,y)}{x_1^i}z_1^{r-i}
\end{array}
\end{equation}

In this case, (\ref{eq951}) has the form of (\ref{eq651}) with 
$$
\text{min}\{\frac{j}{i}\text{ such that } x^j\mid a_i, x^{j+1}\not\,\mid a_i\text{ for }2\le i\le r\}
$$
decreased by 1.

By induction on 
$$
\text{min}\{\frac{j}{i}\text{ such that } x^j\mid a_i, x^{j+1}\not\,\mid a_i\text{ for }2\le i\le r\}
$$
 we can construct a sequence of permissible blowups of r big curves in $\overline S_r$ 
such that the conclusions of the Theorem hold.
\end{pf}

\begin{Theorem}\label{Theorem46}
Suppose that $C_r(X)$ holds with $r\ge 2$,  
the conclusions of Theorem \ref{TheoremE8} hold on $X$
 and if $C$ is a curve in $\overline S_r(X)$,
then $C$ is not r big. 
Suppose that $p\in \overline S_r(X)$ is a 1 point,
 $D$ is a general curve through $p$. 
For a 1 point $q\in D$, define 
$$
\epsilon(D,q) = \nu(F_q(0,0,z))
$$
where 
$(x,y,z)$ are permissible parameters at $q$ so that $\hat{\cal I}_{D,q}=(x,y)$ and
$$
\begin{array}{ll}
u&=x^a\\
v&=P(x)+x^cF_q.
\end{array}
$$

Then there exists a sequence of blowups of points on the strict
transform of $D$, but not
at $p$, $\lambda:Z\rightarrow X$ such that
\begin{enumerate}
\item $C_r(Z)$  and the conclusions of Theorem \ref{TheoremE8} hold on $Z$.
\item Let $\tilde D$ be the strict transform of $D$ on $Z$. Then $\epsilon(\tilde D,q)= 1$ for all  1 points
$q\ne p$ on $\tilde D$, $\nu(q)=0$ if $q\in \tilde D$ is a 2 point, and there are no 3 points on $\tilde D$.
\item Suppose that  $\tilde p$ is a fundamental point of $\lambda$.
\begin{enumerate}
\item If $q\in \lambda^{-1}(\tilde p)$ is a 1 point then $\nu(q)\le r-1$.
\item If $q\in \lambda^{-1}(\tilde p)$ is a 2 point, then $\nu(q)\le r$.
If $\nu(q)=r$ then $\gamma(q)=r$. If   $\nu(q) = r-1$, but  $\tau(q)=0$,
 then $\gamma(q)=r$ and $q$ is on the strict transform of a curve in $\overline S_r(X)$.
\item If $q\in\lambda^{-1}(\tilde p)$ is a 3 point then $\nu(q)\le r-2$.
\end{enumerate}
\item There does not exist a curve $C\subset\overline S_r(Z)$ such that
$C$ is r big.
\end{enumerate}
\end{Theorem}

\begin{pf} The existence of $Z$ and the validity of $C_r(Z)$ and 2. follow from Theorem \ref{TheoremE6}.

Suppose that $D$ contains a 1 point $q\ne p$ such that $\nu(q)=r$. Then $D$ intersects a curve $C$ in $\overline S_r(X)$
transversally at $q$, and  $q$ is  a generic point of $C$. By Lemma
\ref{Lemma5} and Lemma \ref{Lemma17}, since $C$ is not r big, $\gamma(q)=r$ and $q$
is a generic point of $C$, there are thus
 permissible parameters $(x,y,z)$ at $q$ such that
$$
\begin{array}{ll}
u&=x^a\\
F_q&= \tau z^r+\sum_{i=2}^{r-1}a_i(x,y)x^{\alpha_i}z^{r-i}+x^{r-1}y
\end{array} 
$$
where $\alpha_i\ge i$, $\tau$ is a unit and $x\not\,\mid a_i$ for $2\le i\le r-1$,
 $\hat{\cal I}_{C,q}=(x,z)$,
$\hat{\cal I}_{D,q}=(x,y)$.

Let $\pi_1:X_1\rightarrow X$ be the blowup of $q$.

Suppose that $q_1\in\pi_1^{-1}(q)$ and there are permissible parameters $(x_1,y_1,z_1)$ at $q_1$
such that 
$$
x=x_1,
y=x_1(y_1+\alpha), z=x_1(z_1+\beta).
$$
Then
$$
\begin{array}{ll}
u&=x_1^a\\
\frac{F_q}{x_1^r} & = \tau(z_1+\beta)^r+\sum_{i=2}^{r-1}a_ix_1^{\alpha_i-i}(z_1+\beta)^{r-i}+(y_1+\alpha)
\end{array}
$$
Thus, after normalizing to get $F_{q'}$, we have that
  $\gamma(q')\le r-1$ if $\beta\ne 0$, and $\gamma(q')=1$, if $\beta = 0$.

 Suppose that $q_1\in \pi_1^{-1}(q)$ and
there are permissible parameters $(x_1,y_1,z_1)$ at $q_1$ so that
$$
x=x_1y_1,
y=y_1,
z=y_1(z_1+\alpha)
$$
If $\alpha\ne 0$,
then $q'$ is a 2 point with $\gamma(q')\le r-1$.
If $\alpha=0$, then $q'$ is a 2 point on the strict transform of $C$,
$\nu(q')\le r-1$  and $\gamma(q')\le r$.

Suppose that $q_1\in \pi_1^{-1}(q)$ and
there are regular parameters $(x_1,y_1,z_1)$ in $\hat{\cal O}_{X_1,q_1}$ so that
$$
x=x_1z_1,
y=y_1z_1,
z=z_1
$$
Then  $q'$ is a 2 point with $\nu(q')=0$. $q'$ is the point in $\pi_1^{-1}(q)$ on the 
strict transform of $D$. Thus 3. holds for $\lambda^{-1}(q)$.

If $q\in D$ is a 1 point with $\nu(q)\le r-1$ or a 2 point with $\nu(q)\le r-2$, 3. for $\lambda^{-1}(q)$ follows from Theorems \ref{Theorem9} and \ref{Theorem13}
(or the proof of Theorem \ref{TheoremE6}).

If $q\in D$ is a 2 point with $\nu(q)=r-1$, then $q$ is a generic point of a 2 curve $C
\subset \overline S_{r-1}(X)$ (or such that $F_q\in\hat{\cal I}_{C,q}$ if $r=2$). Since $C_r(X)$ holds, $C$ is r-1 small. Since $q$ is a generic point of $C$, we must have $\tau(q)>0$
by Lemma \ref{Lemma6}. 1. - 3. for $\lambda^{-1}(q)$ then follow
from Theorem \ref{Theorem9} and Theorem \ref{Theorem13}.

The conclusions of  Theorem \ref{TheoremE8} hold since this condition is stable under
quadratic transforms.
\end{pf}

\begin{Definition}\label{Def1096} Suppose that $C$ is a 2 curve of $X$. Then $C$ satisfies (E) if
 For $q\in C$,
\begin{enumerate}
\item $\nu(q)=0$ if $q$ is a 3 point.
\item $\gamma(q)\le 1$ at all but finitely many 2 points $q\in C$, where either
\begin{enumerate}
\item $\nu(q)=\gamma(q)=r$ or
\item $\nu(q)=r-1$, $\gamma(q)=r$ and $\tau(q)=0$.
\end{enumerate}
\end{enumerate}
\end{Definition}

\begin{Theorem}\label{Theorem48}
 Suppose that $C_r(X)$ holds with $r\ge 2$, 
the conclusions of Theorem \ref{TheoremE8} hold on $X$
 and $C$ is a 2 curve of $X$ containing a 2 point $p$ such that either
$\nu(p)=\gamma(p)=r$, or $\nu(p)=r-1$, $\gamma(p)=r$ and $\tau(q)=0$. Then there exists a sequence of 
quadratic
transforms $\pi:Y\rightarrow X$ (over points in $C$) such that the following properties hold.
Let $\tilde C$ be the strict transform
of $C$.  Suppose that $q$ is an exceptional point of $\pi$. Then
\begin{enumerate}
\item If $q$ is a 1 point, then $\nu(q)\le r-1$.
\item If $q$ is a 2 point, then $\nu(q)\le r-1$. If $\nu(q)=r-1$ then $\tau(q)>0$.
\item If $q$ is a 3 point then $\nu(q)\le r-2$
\end{enumerate}
Furthermore, $\tilde C$ satisfies (E), $C_r(Y)$ holds and the conclusions of
Theorem \ref{TheoremE8} hold on $Y$. 
\end{Theorem}

\begin{pf}
By our assumption on $p$, and Lemma \ref{Lemma655}, $F_{q'}\not\in\hat{\cal I}_{C,q'}$ at
all  points $q'\in C$. There thus cannot exist $q'\in C$ which satisfies (\ref{eq641}),
since then $F_{q'}\in\hat{\cal I}_{C,q'}$.

Suppose that $q\in C$. Then  there exist permissible parameters
$(x,y,z)$ at $q$ such that $\hat{\cal I}_{C,q} = (x,y)$ and
$F_q\not\in\hat{\cal I}_{C,q}$.

First suppose that $q$ is a 3 point and $\nu(q)>0$. Suppose that $\pi_1:X_1\rightarrow X$
is the blowup of $q$. Suppose that  $q'\in \pi_1^{-1}(q)$. Since $\nu(q)\le r-2$, we have that
$\nu(q')\le r-1$ if $q'$ is a 1 point, $\nu(q')\le r-1$ if $q'$ is a 2 point, $\nu(q')=r-1$ implies $\tau(q')>0$,
and $\nu(q')\le r-2$ if $q'$ is a 3 point by Theorem \ref{Theorem9}. The strict transform of $C$ intersects $\pi_1^{-1}(q)$
in a 3 point.

Consider the infinite sequence of blowups of points
$X_{n+1}\rightarrow X_n$, centered at the points $q_n$ on
the strict transform of $C$ on $X_n$ over $q$, 
$$
\cdots\rightarrow X_n\rightarrow \cdots\rightarrow X_1\rightarrow X
$$
$q_n$ has permissible parameters $(x_n,y_n,z_n)$ defined by
$$
x = x_nz_n^n,
y=y_nz_n^n,
z=z_n
$$
$$
F_{q_n} = \frac{F_q}{z_n^{\sum_{i=1}^ns_i}}
$$
where $s_i=\nu(q_i)$.
Since  $F_q\not\in (x,y)$ 
we have that $\nu(q_n)=0$ for all sufficiently large $n$. 

Let 
$$
W = \{q\in C \vert q\text{ is a 2 point with }\gamma(q)>1\}.
$$
 W is a finite set by Lemma \ref{Lemma42}, and since $\gamma(q)\le 1$ at a generic point of $C$. 

Suppose that $q\in W$ is a 2 point. Suppose that  either $\nu(q)\le r-2$ or $q$ is such that 
$\nu(q)=r-1$ and $\tau(q)>0$. Then arguing as in the case when $q$ is a 3 point,
and using Theorems \ref{Theorem9} and \ref{Theorem13} we can produce a sequence 
of blowups of points
$\pi:X_{m}\rightarrow X$, centered at the points $q_n$ on
the strict transform of $C$ on $X_n$ over $q$, such that $\nu(q_m)=0$, and the conclusions of the Theorem hold 
in a neighborhood of  $\pi^{-1}(q)$.  
\end{pf}

\begin{Theorem}\label{Theorem52}
Suppose that $C_r(X)$ holds with $r\ge 2$ and the conclusions of Theorem \ref{TheoremE8} hold on $X$.
 Then there exists a permissible sequence of blowups $\pi:Y\rightarrow X$ such that
for $p\in Y$
\begin{enumerate}
\item $\nu(p)\le r-1$ if $p$ is a 1 point or a 2 point.
\item If $p$ is a 2 point and $\nu(p)= r-1$, then $\tau(p)>0$ or $r\ge 3$ and (\ref{eq641})
holds at $p$.
\item $\nu(p)\le r-2$ if $p$ is a 3 point
\end{enumerate}
\end{Theorem}

\begin{pf} By Theorem \ref{Theorem45},
we can assume that
if $C$ is a curve in $\overline S_r(X)$ then $C$ is not r big.
Furthermore, since $C_r(X)$ holds, each curve in $\overline S_r(X)$ contains a 1 point.
There are finitely many 1 points $\{p_1,\ldots, p_m\}$ in $X$ such that each $p_i$ is in 
$\overline S_r(X)$, and $p_i$ is
either an isolated point in $\overline S_r(X)$, or is a special point of a curve in $\overline S_r(X)$ (A special 1 point on a curve in $\overline S_r(X)$ is a point which is not
generic in the sense that the conclusions of 1. (a) of Lemma \ref{Lemma41} do not hold). Let $D_{p_i}$
be a general curve through $p_i$ for $1\le i\le m$. By Theorem \ref{Theorem46},
after possibly performing a finite sequence of quadratic transforms 
 at points $\ne p_i$ on the $D_{p_i}$, we may assume that
 $\epsilon(D_{p_i},q)= 1$ for all  1 points
$q\ne p_i$ on $D_{p_i}$, $\nu(q)=0$ if $q\in D_{p_i}$ is a 2 point, and there are no 3 points on any $D_{p_i}$.

There are no exceptional 1 points in $\overline S_r(X)$  created by the sequence of blowups in
Theorem \ref{Theorem46}.

For each $p_i$, let $t_{p_i}$ be the number $l$ computed in Theorem \ref{Theorem22}
(or Theorem \ref{Theorem1011} if $r=2$) for $p_i$.
By Theorem \ref{Theorem21}, for $1\le i\le m$, there exist sequences of monoidal transforms
$$
\lambda^{p_i}:X_1(i)\rightarrow \text{spec}({\cal O}_{X,p_i}) 
$$
where $X_1(i)\rightarrow \text{spec}({\cal O}_{X,p_i})$ is a sequence of 
permissible monodial transforms centered at sections over $D_{p_i}$, 
such that the conclusions of Theorem \ref{Theorem21} hold on $X_1(i)$ (with $t\ge t_{p_i}$).
Since $D_{p_i}$ is resolved at all points $p_i\ne q\in D_{p_i}$, the only obstruction
to extending $\lambda^{p_i}$ to a permissible sequence of monodial transforms of sections
over $D_{p_i}$ in $X$ is if the corresponding sections over $D_{p_i}$ in $X$ do
not make SNCs with the 2 curves. This difficulty can be resolved by performing quadratic
transforms at the points where the section does not make SNCs with the 2 curves.

By Theorem \ref{Theorem48}, we can then perform a sequence of quadratic transforms
centered at 2 points, so that if $C$ is a 2 curve containing a point $p$ such that
either 
$$
\nu(p)=\gamma(p)=r
$$
 or
$$
\nu(p)=r-1, \gamma(p)=r\text{ and }\tau(p)=0
$$
then $C$ satisfies (E). There are no exceptional 1 points in $\overline S_r(X)$ in this
sequence of blowups.

 We can thus extend the maps
$X_1(i)\rightarrow \text{spec}({\cal O}_{X,p_i})$ to a sequence of permissible 
monodial transforms centered at   sections over $D_{p_i}$ and points,
$\lambda:Y\rightarrow X$.  $C_r(Y-\lambda^{-1}(\{p_1,\ldots,p_m\}))$ holds and
the conclusions of Theorem \ref{TheoremE8} hold on $Y-\lambda^{-1}(\{p_1,\ldots,p_m\})$, there are no special or isolated 1
points in $\overline S_r(Y)-\lambda^{-1}(\{p_1,\ldots,p_m\})$,
$$
Y\times_X\text{spec}({\cal O}_{X,p_i})\cong X_1(i)
$$
for $1\le i\le m$
and if $q\in \overline S_r(Y)-\lambda^{-1}(\{p_1,\ldots,p_m\})$
is a 2 point with 
$$
\nu(q)=r, \gamma(q)=r,
$$
or
$$
\nu(q)=r-1, \gamma(q)=r\text{ and }\tau(q)=0
$$
 then the 2 curve containing $q$ satisfies
(E).

Let $\{q_1,\ldots,q_n\}$ be the 2 points of $Y$ such that $\nu(q_i)=r$ and $q_i$ is
contained in a curve of $\overline S_r(Y)$. Let $C_{q_i}$ be the 2 curve containing
$q_i$ for $1\le i\le m$. Then $\gamma(q)\le 1$ if $q_i\ne q\in C_{q_i}$ is a 2 point,
and $\gamma(q)=0$ if $q\in C_{q_i}$ is a 3 point since
$C_{q_i}$ satisfies (E) and the conclusions
of Theorem \ref{TheoremE8} hold on $Y$. 
For each $q_i$, let $t_{q_i}$ be the number $l$ computed in Theorem \ref{Theorem22}
(or Theorem \ref{Theorem1011} if $r=2$) for $p_i$.
For $1\le i\le n$, let
$$
\lambda^{q_i}:Y_1(i)\rightarrow \text{spec}({\cal O}_{Y,q_i})
$$
be a permissible sequence of monodial transforms centered at sections over $C_{q_i}$
such that the conclusions of Theorem \ref{Theorem21} hold on $Y_1(i)$ (with $t\ge t_{q_i}$).

Since $C_{q_i}$ is resolved at all points $q_i\ne q\in C_{q_i}$, the only obstruction
to extending $\lambda^{q_i}$ to a permissible sequence of monodial transforms of
sections over $C_{p_i}$ in $Y$ is if the corresponding sections over $C_{p_i}$ in $Y$
do not make SNCs with the 2 curves. This difficulty can be resolved by performing
quadratic transforms at the points where the section does not make SNCs with the 2 curves.

We can thus extend the $Y_1(i)\rightarrow \text{spec}({\cal O}_{Y,q_i})$ to a
sequence of permissible monodial transforms centered at points and sections over $C_{q_i}$,
$\phi:Z\rightarrow Y$, so that
$$
Z\times_Y\text{spec}({\cal O}_{Y,q_i})\cong Y_1(i)
$$
for $1\le i\le n$, and if
$$
Z_0=Z-\phi^{-1}(\{q_1,\ldots,q_n\})-(\lambda\circ\phi)^{-1}(\{p_1,\ldots,p_m\})
$$
then  $C_r(Z_0)$ holds and $Z_0$ satisfies the conclusions of Theorem \ref{TheoremE8}.
$Z_0$ contains no special or isolated 1 points. 
If $D\subset\overline S_r(Z)$ and $D$ is r big then $D\cap Z_0=\emptyset$. If $q\in Z_0$ is a 2 point which does
not satisfy the conclusions of the Theorem, then $\gamma(q)=r$ and $\nu(q)=r-1$ or $\nu(q)=r$. If 
$C\subset Z$ is the 2 curve containing $q$, then $C$ satisfies (E). If $\nu(q)=r$, then $q$ is not contained in a curve $D\subset \overline 
S_r(Z)$.

By Lemma \ref{Lemma41} and Theorem
\ref{Theorem22} if $r\ge 3$ (or Lemma \ref{Lemma41}, Lemma \ref{Lemma97} and Theorem \ref{Theorem1011} if $r=2$)  there exists a sequence of permissible monodial transforms
$\psi:W\rightarrow Z$ consisting of a sequence of blowups of r big curves 
$D\subset\overline S_r$,
followed by a sequence of blowups of r small curves  $D\subset\overline S_r$, and finally followed by a 
sequence of quadratic transforms if $r\ge 3$ (or quadratic transforms and monodial transforms centered at 2 curves
$C$ such that $C$ is 1 big and $C$ is a section over a 2 small curve blown up in
constructing $\phi$ if $r=2$) such that $C_r(W)$ holds, $\overline S_r(W)$ is a
finite union of 2 points, and the conclusions of the Theorem hold everywhere in $W$,
except possibly at a finite  number of 2 points $p$. If $C$ is a 2 curve on $W$ containing
a 2 point $p$ where the theorem fails to hold, then (E) holds on $C$.
In particular, there are no 2 curves in $\overline S_r(W)$.

Suppose that $C\subset W$ is a 2 curve containing a 2 point such that 
 the conclusions of
the theorem do not hold. Then $C$ satisfies (E).

Let $\{q_1,\ldots,q_s\}$ be the two points on  $C$ such that 
\begin{enumerate}
\item $\nu(q_i)=r$, $\gamma(q_i)=r$ or
\item $\nu(q_i)=r-1$, $\gamma(q_i)=r$ and $\tau(q_i)=0$.
\end{enumerate}

$\gamma(q)\le 1$ if $q\in C-\{q_1,\ldots,q_s\}$ is a 2 point, and $\nu(q)=0$ if
$q\in C$ is a 3 point. No $q_i$ is contained in a curve in $\overline S_r(W)$, since
$\overline S_r(W)$ is finite.

We will now show that
there exists an affine neighborhood $U$ of $\{q_1,\ldots,q_s\}$ and uniformizing
parameters $\tilde x,y,\tilde z$ on $U$ such that $\tilde x=y=0$ are local equations of $C$ on $U$,
and all points of $C\cap U$ are 2 points.

Let $A_1$ and $A_2$ be the components of $E_W$ such that $C$ is a connected component of
$A_1\cap A_2$. There exist very ample
divisors $H_1,H_2,H_3,H_4$ on $W$ such that $A_1\sim H_1-H_2$, $A_2\sim H_3-H_4$ and
$q_i\not\in H_j$ for $1\le i\le s$, $1\le j\le 4$.

Let $U=W-(H_1\cup H_2\cup H_3\cup H_4)$. $U=\text{spec}(A)$ is affine and there
exist $\tilde x, y\in A$ such that $\tilde x=0$ is an equation for $A_1\cap U$ in $U$,
$y=0$ is an equation for $A_2\cap U$ in $U$. After possibly replacing $U$ with a 
smaller affine neighborhood of $\{q_1,\ldots,q_s\}$, we may assume that $U\cap C=U\cap A_1\cap 
A_2$ and $E_W\cap U=(A_1\cup A_2)\cap U$.

There exists a morphism $\pi:C\rightarrow \bold P^1$ such that $\pi$ is \'etale
over $\pi(q_i)$, $1\le i\le s$, and $\pi(q_i)\ne\infty$ for any $i$ (We can take $\pi$
to be a generic projection).
Let $z$ be a coordinate on $\bold P^1-\{\infty \}$.
After replacing $U$ with a possibly smaller affine neighborhood of $\{q_1,\ldots, q_s\}$
we have an inclusion $\pi^*:k[z]\rightarrow A$, so that $U\rightarrow \text{spec}(k[\tilde x,y,z])$
is \'etale. 

 There exists a component $E$ of $D_S$ such that $\Phi_W(A_1)\subset E$ and
$\Phi_W(A_2)\subset E$ (since $\Phi_W:W\rightarrow S$ is weakly prepared). There exists an affine neighborhood $\overline V$ of
$\{\Phi_W(q_1),\ldots,\Phi_W(q_s)\}$ in $S$ and $u\in\Gamma(\overline V,{\cal O_S})$ such that
$u=0$ is a local equation of $D_S$. Then $u=0$ is a local equation of $E_W$ in $\Phi_W^{-1}(\overline V)\cap U$. Thus if we replace $U$ with $\Phi_W^{-1}(\overline V)\cap U$, $u$ extends to a system of permissible parameters at $\Phi_W(p)$ for all $p\in C\cap\Phi_W^{-1}(\overline V)\cap U$.

There exist $a,b\in\bold N$ such that $u=\tilde x^ay^b\overline \gamma$ where
$\overline \gamma\in A$  is a unit in $A$.
Let $\gamma=\overline\gamma^{\frac{1}{a}}$, $B=A[\gamma]$, $V=\text{spec}(B)$.
Then $h:V\rightarrow U$ is \'etale. 

Let $x=\gamma \tilde x$. $k[x,y,z]\rightarrow B$ defines a morphism $g:V\rightarrow \bold A^3$.  $q\in g^{-1}(x=0)$ if and only if $x\in m_q$ which holds if and only if $\tilde x\in m_q$. Thus $g^{-1}(x=0)=h^{-1}(\tilde x=0)$.
$g$ is \'etale at all points of $g^{-1}(x=0)$. Since this is an open condition (c.f. Prop 4.5 SGA1)
there exists a Zariski closed subset $Z_1$ of $V$ which is disjoint from $h^{-1}(\tilde x=0)$
such that $g\mid V-Z_1$ is \'etale. Let $U_1$ be an affine neighborhood of $\{q_1,\ldots,q_s\}$ in $U$ which is disjoint from $h(Z_1)$. Let $V_1=h^{-1}(U_1)$.

After replacing $U$ with $U_1$ and $V$ with $V_1$, we have that $V\rightarrow U$ is an
\'etale cover 
 and   $(x,y,z)$ are uniformizing parameters on $V$.

There exist $v_i\in{\cal O}_{S,\Phi_W(q_i)}$ such that $(u,v_i)$ are permissible parameters at $\Phi_W(q_i)$ and $u=0$ is a local equation of $E_W$ at $q_i$
for $1\le i\le s$. For each $q_i$ there exist  $z_i$ such that $(x,y,z_i)$ are permissible parameters at $q_i$ for $(u,v_i)$ for $1\le i\le s$
which satisfy the conclusions of Lemma \ref{Lemma17}.

The morphism $\pi_1:Y_{q_1}\rightarrow\text{spec}(\hat{\cal O}_{W,q_1})$ of
Theorem \ref{Theorem19} (or Theorem \ref{Theorem51} if $\nu(q_i)=r-1$) extends to a sequence of
permissible monodial transforms $\pi_1:\tilde Y_1\rightarrow V$ centered at sections over
$C$.
$$
\tilde Y_1\times_{\text{spec}({\cal O}_{W,q_2})}\text{spec}(\hat{\cal O}_{W,q_2})
\rightarrow \text{spec}(\hat{\cal O}_{W,q_2})
$$
extends to a sequence of permissible  monodial transforms $Y_{q_2}\rightarrow
\text{spec}(\hat{\cal O}_{W,q_2})$ of the form of the conclusions of Theorem \ref{Theorem19}
(or Theorem \ref{Theorem51}).

$Y_{q_2}\rightarrow \text{spec}(\hat{\cal O}_{W,q_2})$ extends to a sequence of
permissible monodial transforms
$$
\tilde Y_2\stackrel{\pi_2}{\rightarrow}\tilde Y_1-\pi_1^{-1}(h^{-1}(q_1))\rightarrow V-\{h^{-1}(q_1)\}.
$$
Preceeding inductively, we extend
$$
\tilde Y_{s-1}\times_{\text{spec}({\cal O}_{W,q_s})}\text{spec}(\hat{\cal O}_{W,q_s})
\rightarrow \text{spec}(\hat{\cal O}_{W,q_s})
$$
to a sequence of permissible monodial transforms $Y_{q_s}\rightarrow
\text{spec}(\hat{\cal O}_{W,q_s})$ of the form of the conclusions of Theorem \ref{Theorem19}
(or Theorem \ref{Theorem51}).

$Y_{q_s}\rightarrow \text{spec}(\hat{\cal O}_{W,q_s})$ extends to a sequence of
permissible monodial transforms
$$
\tilde Y_s\stackrel{\pi_s}{\rightarrow}\tilde Y_{s-1}
-(\pi_1\circ\cdots\circ\pi_{s-1})^{-1}(h^{-1}(q_{s-1}))
\rightarrow V-\{q_1,\ldots,q_{s-1}\}.
$$

For $1\le i\le s$, let $t_{q_i}$ be the value of $l$ in the statement of Theorem
\ref{Theorem22} (or Theorem \ref{Theorem51}) for the point $q_i$.

Let $\omega_i\in A$, $1\le i\le s$ be such that
$$
\omega_i\equiv \gamma\text{ mod }m_{q_i}^{t_{q_i}}\hat{\cal O}_{W,q_i}.
$$
By   the Chinese Remainder Theorem, there exists
$\omega\in  A$ such that after possibly replacing $U$ with
a smaller affine neighborhood of $\{q_1,\ldots,q_s\}$, we have that 
$(\omega \tilde x,y,\tilde z)$ are uniformizing parameters on $U$ and 
$$
\omega\tilde x\equiv x\text{ mod }m_{q_i}^{t_{q_i}}\hat{\cal O}_{W,q_i}
$$
for $1\le i\le s$.

We can thus replace $\tilde x$ with $\omega \tilde x$ in (\ref{eq646}) of Theorem
\ref{Theorem21} for $1\le i\le s$.
With this choice of $\tilde x$, The map
$\overline Y_{q_1}\rightarrow \text{spec}({\cal O}_{W,q_1})$ of Theorem \ref{Theorem21}
(or Theorem \ref{Theorem51}) 
then satisfies the assumptions of Theorem \ref{Theorem22}, and extends to a
permissible sequence of monodial transforms centered at sections over $C$ 
$$
\lambda_1:\hat Y_1\rightarrow U.
$$

The map $\overline Y_{q_2}\rightarrow \text{spec}({\cal O}_{W,q_2})$ of Theorem
\ref{Theorem21} (or Theorem \ref{Theorem51}) satisfies the assumptions of Theorem
\ref{Theorem22} (or Theorem \ref{Theorem51}) and extends to a permissible sequence of monodial transforms
centered at sections over $C$, 
$$
\hat Y_2\stackrel{\lambda_2}{\rightarrow}\hat Y_1-\lambda_1^{-1}(q_1)\rightarrow U-\{q_1\}.
$$ 

Preceeding inductively, the map $\overline Y_{q_s}\rightarrow \text{spec}({\cal O}_{W,q_s})$
of Theorem \ref{Theorem21} (or Theorem \ref{Theorem51}) satisfies the assumptions of
Theorem \ref{Theorem22}, and extends to a permissible sequence of monodial transforms
centered at sections over $C$
$$
\hat Y_s\stackrel{\lambda_s}{\rightarrow}\hat Y_{s-1}
-(\lambda_1\circ\cdots\circ\lambda_{s-1})^{-1}(q_{s-1})
\rightarrow U-\{q_1,\ldots,q_{s-1}\}.
$$

Since all points of $C-U$ are resolved, there exists a sequence of permissible
monodial transforms $\tilde \lambda_1:Y_1'\rightarrow W$ consisting of
quadratic transforms centered at points over $C-U$ and permissible monodial
transforms centered at sections of $C$ such that all points of $\tilde\lambda_1^{-1}
(C-U)$ are resolved, and $Y_1'\times_WU\cong\hat Y_1$.

By Theorems \ref{Theorem22} and \ref{Theorem51} (or Theorems \ref{Theorem1011} and
Theorem \ref{Theorem1012} if $r=2$), there exists a sequence of 
permissible monodial transforms $\tilde Z_1\rightarrow Y_1'$
(with induced maps $\psi_1:\tilde Z_1\rightarrow W$) centered over points and
curves which map to $q_1$ such that all points of $\psi_1^{-1}(q_1)$ satisfy the
conclusions of the Theorem.

By Theorem \ref{Theorem9} and \ref{Theorem13},
$$
\hat Y_2\rightarrow \hat Y_1-\lambda_1^{-1}(q_1)
$$
then extends to a sequence of permissible monodial transforms 
$\tilde\lambda_2:Y_2'\rightarrow\tilde Z_1$ consisting of quadratic transforms
centered at points over $C-(U-\{q_1\})$ and permissible monodial transforms centered
at sections over $C$ such that all points of  
$(\tilde \lambda_1\circ\tilde\lambda_2)^{-1}(C-U)$ are resolved, all points of
$(\tilde\lambda_1\circ\tilde\lambda_2)^{-1}(q_1)$ satisfy the conclusions of the Theorem
and
$$
Y_2'\times_W(U-\{q_1\})\cong\hat Y_2.
$$
By Theorems \ref{Theorem22} and \ref{Theorem51} (or Theorems \ref{Theorem1011} and
\ref{Theorem1012} if $r=2$), there exists a sequence of permissible
monodial transforms $\tilde Z_2\rightarrow Y_2'$ with induced maps
$\psi_2:\tilde Z_2\rightarrow W$ centered at points and curves that map to $q_2$ 
such that all points of $\psi_2^{-1}(q_2)$ satisfy the conclusions of the Theorem.

By induction on $s$, we can then construct a sequence of permissible monodial transforms
$\psi_s:\tilde Z_s\rightarrow W$ centered at points and curves supported over $C$
such that all points of $\psi_s^{-1}(C)$ satisfy the conclusions of the Theorem, and
all points of $\psi_s^{-1}(C-\{q_1,\ldots,q_s\})$ are resolved.

By induction on the number of 2 curves $C\subset W$ which contain a 2 point which
does not satisfy the conclusions of the Theorem, we can construct a sequence of permissible
monodial transforms $\overline W\rightarrow W$ such that $\overline W$ satisfies the conclusions
of the Theorem.

\end{pf}

\section{Resolution 4}

Throughout this section we will assume that $\Phi_X:X\rightarrow S$ is weakly prepared.

\begin{Theorem}\label{Theorem33}
Suppose that $r\ge 1$ and 
for $p\in X$,
\begin{enumerate}
\item $\nu(p)\le r$ if $p$ is a 1 point or a 2 point.
\item If $p$ is a 2 point and $\nu(p)= r$, then $\tau(p)>0$   or $r\ge 2$,
$\tau(p)=0$ and there exists a unique curve $D\subset \overline S_r(X)$ 
(containing a 1 point) such that
$p\in D$, and permissible parameters $(x,y,z)$ at $p$ such that 
$x=z=0$ are local equations of $D$.
\begin{equation}\label{eq640}
\begin{array}{ll}
u&=(x^ay^b)^m\\
v&=P(x^ay^b)+x^cy^dF_p\\
F_p&=\tau x^r+\sum_{j=1}^r\overline a_j(y,z)y^{d_j}z^{e_j}x^{r-j}
\end{array}
\end{equation}

where $\tau$ is a unit, $\overline a_j$ are units (or 0), there exists $i$ such that
$\overline a_i\ne 0$, $e_i=i$, $0<d_i<i$,
$$
\frac{d_i}{i}\le \frac{d_j}{j}, \frac{e_i}{i}\le \frac{e_j}{j}
$$
for all $j$ and
$$
\left\{\frac{d_i}{i}\right\}+\left\{\frac{e_i}{i}\right\}<1.
$$

\item $\nu(p)\le r-1$ if $p$ is a 3 point
\end{enumerate}
 Then there exists a sequence of permissible monoidal transforms $\pi:X_1\rightarrow X$  such that 
$\overline A_r(X_1)$ holds. That is,
\begin{enumerate}
\item $\nu(p)\le r$ if $p\in X$ is a 1 point or a 2 point.
\item If $p\in X$ is a 1 point and $\nu(p)= r$, then $\gamma(p)=r$.
\item If $p\in X$ is a 2 point and $\nu(p)= r$, then $\tau(p)>0$.
\item $\nu(p)\le r-1$ if $p\in X$ is a 3 point
\end{enumerate}
\end{Theorem}

\begin{pf} If $p$ is a 1 point such that $\nu(p)=1$, then $\gamma(p)=1$. Thus $\overline A_r(X)$ holds if $r=1$. For the rest of the proof we will assume that $r\ge 2$.

Let
$$
W(X)=\{p\in \text{1 points of }X\vert \nu(p)=r\text{ and }\gamma(p)>r\}.
$$
$W(X)$ is Zariski closed in the open subset of
1 points of $X$. Let $\overline W(X)$ be the Zariski closure of $W(X)$ in $X$.

Suppose that $p\in\overline W(X)$ is a point where $\overline W(X)$ does not make
SNCs with $\overline B_2(X)$. Then  $p$ can not satisfy (\ref{eq640}).
Let $\pi:X_1\rightarrow X$ be the quadratic transform with center $p$. By Theorems
\ref{Theorem9} and \ref{Theorem13}, all points of $\pi^{-1}(p)$ satisfy the
assumptions of Theorem \ref{Theorem33}, and there are no points of $\pi^{-1}(p)$
which satisfy (\ref{eq640}). If $p$ is a 1 point, $x\mid L_p$ implies $\nu(q)\le r-1$ if
$q\in\pi^{-1}(p)$ is a 1 point, so $\pi^{-1}(p)$ contains no curves of $\overline W(X_1)$. By Theorems \ref{Theorem9} and \ref{Theorem13},
$\pi^{-1}(p)$ contains no curves of $\overline W(X_1)$
if $p$ is a 2 or 3 point.

Thus there exists a sequence of quadratic transforms $\pi:X_1\rightarrow X$
such that $\overline W(X_1)$ is a disjoint union of nonsingular curves and
isolated points, $X_1$ satisfies the assumptions of Theorem \ref{Theorem33}, and
$\overline W(X_1)$ makes SNCs with $\overline B_2(X_1)$. By Theorems
\ref{Theorem9} and \ref{Theorem13} and Lemma \ref{Lemma54}, we can further assume that
$\overline S_r(X_1)$ makes SNCs with $\overline B_2(X_1)$, except possibly at some
3 points of $X_1$, and if $C\subset\overline S_r(X_1)$ is a curve which contains
a 2 point satisfying (\ref{eq640}), then $C$ contains no 3 points. We can then
without loss of generality assume that $X=X_1$.

Suppose that $C\subset \overline W(X)$ is a curve. $C$ makes SNCs with the locus of 2 curves.
We either have that $C$ is r big or r small.

For a curve $C$, or isolated point $p$ in $\overline W(X)$, 
We will show that we can construct a sequence of monoidal transforms
$\pi:Y\rightarrow X$, centered at points and curves over $C$ (or over $p$), such that 
the assumptions of the theorem hold  on  $Y$, and 2. of the conclusions of the theorem hold at points over $C$
(over $p$). 

We can then iterate this process to obtain $Z\rightarrow X$ such that the assumptions of
the theorem hold on $Z$, and if $p\in Z$ is a 1 point with $\nu(p)=r$, then $\gamma(p)=r$.

{\bf Suppose that $C$ is r small}
Since $C$ is r small, (\ref{eq640}) cannot hold at any $p\in C$. By Lemma \ref{Lemma975}, we
can  construct a sequence of monoidal transforms
$\pi:Y\rightarrow X$, centered at points on $C$ and the strict transform of $C$, such that 
the assumptions of the theorem hold  on  $Y$, and the conclusions of the theorem hold at points of $\pi^{-1}(C)$.

{\bf Suppose that $C$ is r big}

Let $\pi:X_1\rightarrow X$ be the blowup of $C$. We will show that the 
assumptions of the theorem and 2. of then conclusions
of the theorem hold at points above $C$.

Suppose that $p\in C$ is a 2 point with $\tau(p)>0$ or a 1 point. Then all points
of $\pi^{-1}(p)$ satisfy the conclusions of the Theorem by Lemma \ref{Lemma654}.

Suppose that $p\in C$ is a 2 point such that (\ref{eq640}) holds. Then $x=z=0$ are local
equations of $C$ at $p$.

Suppose that $q\in\pi^{-1}(p)$ is a 2 point. $q$ has permissible parameters $(x_1,y,z_1)$
such that 
$x=x_1, z=x_1(z_1+\alpha)$.
$$
\begin{array}{ll}
u&=(x_1^ay^b)^m\\
v&=P(x_1^ay^b)+x_1^{c+r}y^d\frac{F_p}{x_1^r}.
\end{array}
$$
$$
F_p=\tau x^r+y\Omega
$$
implies
$$
\frac{F_p}{x_1^r}=\tau+y\frac{\Omega}{x_1^r}
$$
$ad-b(c+r)\ne 0$, since $F_p$ is normalized, which implies that
$\nu(q)=0$.

Suppose that $q\in\pi^{-1}(p)$ is the 3 point. $q$ has permissible parameters $(x_1,y,z_1)$
such that
$$
x=x_1z_1, z=z_1.
$$
$$
\begin{array}{ll}
u&=(x_1^ay^bz_1^a)^m\\
v&=P(x^ay^bz_1^a)+x_1^cy^dz_1^{c+r}F_q
\end{array}
$$
where
$$
F_q=\frac{F_p}{z_1^{r}}=\tau x_1^r+\sum_{j=1}^r\overline a_j(y,z_1)y^{d_j}x_1^{r-j}z_1^{e_j-j}
$$
By assumption
$d_i+r-i+e_i-i<r$. Thus $\nu(q)\le r-1$.

{\bf Suppose that $p$ is an isolated point in $\overline W(X)$}

 There are permissible parameters $(x,y,z)$ at $p$ such that
$$
\begin{array}{ll}
u&=x^a\\
v&= P(x)+x^cF_p\\
L_p &=x^t\Omega(x,y,z)
\end{array}
$$
with $0<t<r$, $x\not\,\mid \Omega$.

Let $\pi_1:X_1\rightarrow X$ be the blowup of $p$.
If $q\in\pi_1^{-1}(p)$ is a 2 point then $\nu(q)\le r$ and $\nu(q)=r$ implies $\tau(q)>0$
by Theorem \ref{Theorem9}.
If $q\in\pi_1^{-1}(p)$ is a 1 point then $\nu(q)\le r-t<r$

We are now reduced to assuming that $\overline W(X)=\emptyset$, so that $\gamma(p)=r$ if
$p\in X$ is a 1 point with $\nu(p)=r$.

Now suppose that $p\in X$ satisfies (\ref{eq640}) so that the curve $D$ in $\overline S_r(X)$
that $p$ lies on satisfies $\gamma(q)=r$ if $q\in D$ is a 1 point.

By our initial reduction, we may assume that $D$ is nonsingular, and makes SNCs with
$\overline B_2(X)$. Since $F_p\in\hat{\cal I}_{C,p}^r$, $C$ is r big.

Let $\pi:X_1\rightarrow X$ be the blowup of $D$.  If $p\in D$ is a 2 point
with $\tau(p)>0$ or a 1 point, then all points of $\pi^{-1}(p)$ satisfy the conclusions of the Theorem by Lemma \ref{Lemma654}.
 The case when $p$ satisfies (\ref{eq640}) is exactly as in
the case when $C\subset \overline W(X)$ is r big.
\end{pf}

\section{Proof of the Main Theorem}

\begin{Theorem}\label{Theorem1000} Suppose that 
$\Phi_X:X\rightarrow S$ is weakly prepared, $r\ge 2$ and $\overline A_r(X)$ holds. Then there
exists a permissible sequence of monodial transforms $Y\rightarrow X$ such that
$\overline A_{r-1}(Y)$ holds.
\end{Theorem}

\begin{pf} The Theorem follows from successive application of Lemma \ref{Lemma663}
and Theorems \ref{Theorem667}, \ref{TheoremE7}, \ref{TheoremE8}, \ref{Theorem52} and 
\ref{Theorem33}
\end{pf}

\begin{Theorem}\label{Theorem1050} Suppose that $\Phi_X:X\rightarrow S$ is weakly prepared. Then there exists
a sequence of permissible monoidal transforms $Y\rightarrow X$ such that
$\Phi_Y:Y\rightarrow S$ is prepared.
\end{Theorem}

\begin{pf}
For $r>>0$ $\overline A_r(X)$ holds by Zariski's Subspace Theorem
(Theorem 10.6 \cite{Ab5}). The theorem then follows from successive application
of Theorem \ref{Theorem1000}, and the fact that $\overline A_1(X)$ holds
if and only if $\Phi_X:X\rightarrow S$ is prepared.
\end{pf}

\begin{Theorem}\label{Theorem58}
 Suppose that $\Phi:X\rightarrow S$ is a  dominant morphism from a  
3 fold to a  surface and $D_S\subset S$ is a reduced 1 cycle  such that
$E_X=\Phi^{-1}(D_S)_{red}$ contains $\text{sing}(X)$ and $\text{sing}(\Phi)$. Then there exist sequences of monoidal 
transforms with nonsingular centers $\pi_1:S_1\rightarrow S$ and $\pi_2:X_1\rightarrow X$ such that 
$\Phi_{X_1}:X_1\rightarrow S_1$ is prepared with respect to
$D_{S_1}=\pi_2^{-1}(D_S)_{\text{red}}$.
\end{Theorem}

\begin{pf} This follows from Lemma \ref{Lemma1048} and Theorem \ref{Theorem1050}.
\end{pf}

\section{Monomialization}

Throughout this section we will suppose that $\Phi:X\rightarrow S$ is a  dominant morphism from a
nonsingular 3 fold to a nonsingular surface, $D_S$ is a reduced SNC divisor on $S$,
$E_X=\Phi^{-1}(D_S)_{\text{red}}$ is a SNC divisor on $X$.

If $p\in E_X$ we will say that $p$ is a 1, 2 or 3 point depending on if
$p$ is contained in 1, 2 or 3 components of $E_X$. $q\in D_S$ will be called
a 1 or 2 point depending on if $q$ is contained in 1 or 2 components of $D_S$.

Regular parameters $(u,v)$ in ${\cal O}_{X,p}$ with $q\in D_S$ are permissible if:
\begin{enumerate}
\item  $u=0$ is a local equation of $D_S$ if $q$ is a 1 point or
\item $uv=0$ is a local equation of $D_S$ if $q$ is a 2 point.
\end{enumerate}

\begin{Definition}\label{Def1071}
We will say that $\Phi$ is Strongly Prepared at $p\in X$ (with respect to $D_S$) if one on the following
forms hold.
\begin{enumerate}
\item $\Phi$ is prepared at $p$ (as defined in Definition \ref{Def57}) or
\item There exist permissible parameters $(u,v)$ at $q$ and regular parameters $(x,y,z)$ in $\hat{\cal O}_{X,p}$
such that one of the following hold:
\begin{enumerate}
\item $p$ is a 2 point and 
$$
u=x^a, v=y^b.
$$
\item $p$ is a 3 point and 
$$
u=x^a, v=y^bz^c
$$
 (with $a,b,c>0$).
\item $p$ is a 3 point and 
$$
u=x^ay^b, v=y^cz^d
$$
(with $a,b,c,d>0$).
\end{enumerate}
\end{enumerate}
\end{Definition}

Suppose that $p\in X$ is strongly prepared and $(u,v)$ are permissible parameters at $\Phi(p)$.
Regular parameters $(x,y,z)$ in $\hat{\cal O}_{X,p}$ are called $*$-permissible parameters at $p$ 
for $(u,v)$ if one of the forms of Definition \ref{Def1071} holds in $\hat{\cal O}_{X,p}$.
We will also say that $(u,v)$ are strongly prepared at $p$. If a form 1. holds at $p$,
$*$-permissible parameters are permissible as defined in Definition \ref{Def650}.

Throughout this section we will assume that $\Phi:X\rightarrow S$ is strongly prepared.

\begin{Lemma}\label{Lemma1035} Suppose that ${\cal O}_{X,p}\rightarrow R$ is
finite \'etale, and there exists $\overline x,\overline y,\overline z\in R$ such that
$(\overline x,\overline y,\overline z)$ are regular parameters in $R_q$ for all
primes $q\subset R$ such that $q\cap {\cal O}_{X,p}=m_p$. Then there exists an \'etale
neighborhood $U$ of $p$ such that $(\overline x,\overline y,\overline z)$ are
uniformizing parameters on $U$.
\end{Lemma}

\begin{pf} There exists an affine neighborhood $V_1=\text{spec}(A)$ of $p\in X$ and
a finite \'etale extension $B$ of $A$ such that $B\otimes_AA_{m_p}\cong R$. Set
$U_1=\text{spec}(B)$. Let $\pi:U_1\rightarrow V_1$ be the natural map. There exists
an open neighborhood $U_2$ of $\pi^{-1}(p)$ such that $(\overline x,\overline y,\overline z)$
are uniformizing parameters on $U_2$. Let $Z=U_1-U_2$. Set $U_3=U_1-\pi^{-1}(W)$.
$U_3\rightarrow V_2=V_1-W$ is finite \'etale.
Thus there exists an \'etale neighborhood $U$ of $p$ where $(\overline x,\overline y,\overline z)$
are uniformizing parameters.
\end{pf}

\begin{Lemma}\label{Lemma1036}
Suppose that permissible parameters $(u,v)$ for $\Phi(p)\in D_S$ are 
strongly prepared at $p\in E_X$.
Then there exist $*$-permissible parameters $(x,y,z)$ at $p$ such that $(x,y,z)$ are uniformizing
parameters on an \'etale neighborhood of $p$, and one of the following forms hold:
\begin{enumerate}
\item $p$ is a 1 point, $u=0$ is a local equation of $E_X$ and
$$
\begin{array}{ll}
u&=x^a\\
v&=P(x)+x^by
\end{array}
$$
where $P(x)$ is a polynomial of degree $\le b$.
\item $p$ is a 2 point, $u=0$ is a local equation of $E_X$ and
$$
\begin{array}{ll}
u&=(x^ay^b)^m\\
v&=P(x^ay^b)+x^cy^d
\end{array}
$$
where $(a,b)=1$, $ad-bc\ne 0$, $P(t)$ is a polynomial of degree $\le\left[\text{ max }\{\frac{c}{a},\frac{d}{b}
\}\right]$.
\item $p$ is a 2 point, $u=0$ is a local equation of $E_X$ and
$$
\begin{array}{ll}
u&=(x^ay^b)^m\\
v&=P(x^ay^b)+x^cy^dz
\end{array}
$$
where $(a,b)=1$, $P(t)$ is a polynomial of degree $\le \left[\text{ max }\{\frac{c}{a},\frac{d}{b}\}\right]$.
\item $p$ is a 3 point, $u=0$ is a local equation of $E_X$ and
$$
\begin{array}{ll}
u&=(x^ay^bz^c)^m\\
v&=P(x^ay^bz^c)+x^dy^ez^f
\end{array}
$$
where $(a,b,c)=1$, $P(t)$ is a polynomial of degree $\le \left[\text{ max }\{\frac{d}{a},\frac{e}{b},
\frac{f}{c}\}\right]$.
\item $p$ is a 2 point, $uv=0$ is a local equation of $E_X$ and 
$$
u=x^a, v=y^b.
$$
\item $p$ is a 3 point, $uv=0$ is a local equation of $E_X$ and 
$$
u=x^ay^b, v=z^c
$$
 (with $a,b,c>0$).
\item $p$ is a 3 point, $uv=0$ is a local equation of $E_X$ and 
$$
u=x^ay^b, v=y^cz^d
$$
(with $a,b,c,d>0$).
\end{enumerate}
\end{Lemma}

\begin{pf}
Suppose there exist regular parameters $(x,y,z)$ in $\hat{\cal O}_{X,p}$ such that
$$
\begin{array}{ll}
u&=x^a\\
v&=P(x)+x^by.
\end{array}
$$
There exist $\alpha\in \hat{\cal O}_{X,p}$ and $\overline x\in{\cal O}_{X,p}$
such that $x=\alpha\overline x$, and $\alpha^a\in{\cal O}_{X,p}$. 
 Set $R={\cal O}_{X,p}[\alpha]$.
Let $L$ be the quotient field of $R$.
$R$ is finite \'etale over ${\cal O}_{X,p}$.
$$
v-P_b(\alpha \overline x)=\alpha^by+
\frac{(P(\alpha \overline x)-P_b(\alpha \overline x)}{\overline x^b})\overline x^b
$$
implies
$$
\overline y=\frac{v-P_b(\alpha\overline x)}{\overline x^b}\in \hat{(R_q)}\cap L=R_q
$$
(by Lemma 2.1 \cite{C2}) for all maximal ideals $q\subset R$. Thus
$\overline y\in \cap R_q=R$. Choose $\overline z\in{\cal O}_{X,p}$ such that 
$$
z\equiv \overline z\text{ mod }m_p^2\hat{\cal O}_{X,p}.
$$
Then $m_pR=(\overline x,\overline y,\overline z)$. 
By Lemma \ref{Lemma1035} there exists an \'etale neighborhood $U$ of $p$ such that
$(\overline x,\overline y,\overline z)$ are uniformizing parameters on $U$.

Suppose there exist regular parameters $(x,y,z)$ in $\hat{\cal O}_{X,p}$ such that
$$
\begin{array}{ll}
u&=(x^ay^b)^m\\
v&=P(x^ay^b)+x^cy^d
\end{array}
$$
There exists $\alpha_1,\alpha_2\in\hat{\cal O}_{X,p}$ and $\overline x,\overline y\in
{\cal O}_{X,p}$ such that $x=\alpha_1\overline x, y=\alpha_2\overline y$.
Set 
$$
e=\left[\text{ max }\{\frac{c}{a},\frac{d}{b}\}\right].
$$
$$
u=(\alpha_1^a\alpha_2^b)^m(\overline x^a\overline y^b)^m.
$$
Set $\gamma=\alpha_1^a\alpha_2^b$. Let $K$ be the quotient field of ${\cal O}_{X,p}$.
$$
\gamma^m=\frac{u}{(\overline x^a\overline y^b)^m}\in\hat{\cal O}_{X,p}\cap K={\cal O}_{X,p}.
$$
Set $R={\cal O}_{X,p}[\gamma]$. $R$ is finite \'etale over ${\cal O}_{X,p}$.
Let $L$ be the quotient field of $R$. Set
$$
\omega=\alpha_1^c\alpha_2^d+\frac{P(\alpha_1^a\alpha_2^b\overline x^a\overline y^b)-P_{e}(\alpha_1^a\alpha_2^b\overline x^a\overline y^b)}{\overline x^c\overline y^d}
=\frac{v-P_e(\alpha_1^a\alpha_2^b\overline x^a\overline y^b)}{\overline x^c\overline y^d}
\in \hat{(R_q)}\cap L=R_q
$$
for all maximal ideals $q$ of $R$. Thus $\omega\in \cap R_q=R$. Set $f=ad-bc$.
Set
$$
\tilde x=(\gamma^d\omega^{-b})^{\frac{1}{f}}\overline x,
$$
$$
\tilde y =(\gamma^{-c}\omega^a)^{\frac{1}{f}}\overline y.
$$
$$
\begin{array}{ll}
u&=(\tilde x^a\tilde y^b)^m\\
v&=P_e(\tilde x^a\tilde y^b)+\tilde x^c\tilde y^d.
\end{array}
$$
Choose $\tilde z\in {\cal O}_{X,p}$ such that $z\equiv \tilde z\text{ mod }m_p^2\hat{\cal O}_{X,p}$. Then  
$$
\tilde x,\tilde y,\tilde z\in R_1=R[(\gamma^d\omega^{-b})^{\frac{1}{f}},
(\gamma^{-c}\omega^a)^{\frac{1}{f}}]
$$
are regular parameters at all maximal ideals of $R_1$. By Lemma \ref{Lemma1035}, there
exists an \'etale neighborhood $U$ of $p$ such that $\tilde x,\tilde y,\tilde z$ are uniformizing
parameters on $U$.

Suppose there exist regular parameters $(x,y,z)$ in $\hat{\cal O}_{X,p}$ such that
$$
\begin{array}{ll}
u&=(x^ay^b)^m\\
v&=P(x^ay^b)+x^cy^dz.
\end{array}
$$
There exist $\alpha_1,\alpha_2\in\hat{\cal O}_{X,p}$ and $\overline x,\overline y\in{\cal O}_{X,p}$ such that $x=\alpha_1\overline x, y=\alpha_2\overline y$. Set $e=\left[\text{ max }\{
\frac{c}{a},\frac{d}{b}\}\right]$. Let $K$ be the quotient field of ${\cal O}_{X,p}$.
$$
u=(\alpha_1^a\alpha_2^b)^m(\overline x^a\overline y^b)^m.
$$
Set $\gamma=\alpha_1^a\alpha_2^b$.
$$
\gamma^m=\frac{u}{(\overline x^a\overline y^b)^m}\in \hat{\cal O}_{X,p}\cap K={\cal O}_{X,p}.
$$
Set $R={\cal O}_{X,p}[\gamma]$. $R$ is finite \'etale over ${\cal O}_{X,p}$.
Let $L$ be the quotient field of $R$. Set
$$
\overline z=\frac{v-P_e(\alpha_1^a\alpha_2^b\overline x^a\overline y^b)}{\overline x^c\overline y^d}\in \hat{(R_q)}\cap L=R_q
$$
for all maximal ideals $q$ of $R$. Thus $\overline z\in\cap R_q=R$. 
Set $\tilde x=\alpha_1\alpha_2^{\frac{b}{a}}\overline x$. 
$\tilde x, \overline y,\overline z\in R_1=R[\alpha_1\alpha_2^{\frac{b}{a}}]$
and $(\tilde x,\overline y,\overline z)=m_pR$. By Lemma \ref{Lemma1035} there
exists an \'etale neighborhood $U$ of $p$ such that $(\tilde x,\overline y,\overline z)$
are uniformizing parameters on $U$.

Suppose there exist regular parameters $(x,y,z)$ in $\hat{\cal O}_{X,p}$ such that
$$
\begin{array}{ll}
u&=(x^ay^bz^c)^m\\
v&=P(x^ay^bz^c)+x^dy^ez^f.
\end{array}
$$
There exist $\alpha_1,\alpha_2,\alpha_3\in\hat{\cal O}_{X,p}$ and $\overline x,\overline y,
\overline z\in {\cal O}_{X,p}$ such that $x=\alpha_1\overline x$, $y=\alpha_2\overline y$,
$z=\alpha_3\overline z$. Set $g=\left[\text{ max }\{\frac{d}{a},\frac{e}{b},\frac{f}{c}\}\right]$.
$$
u=(\alpha_1^a\alpha_2^b\alpha_3^c)^m(\overline x^a\overline y^b\overline z^c)^m.
$$
Set $\gamma=\alpha_1^a\alpha_2^b\alpha_3^c$. Let $K$ be the quotient field of ${\cal O}_{X,p}$.
$$
\gamma^m=\frac{u}{(\overline x^a\overline y^b\overline z^c)^m}\in \hat {\cal O}_{X,p}\cap K
={\cal O}_{X,p}
$$
Set $R={\cal O}_{X,p}[\gamma]$. $R$ is finite \'etale over ${\cal O}_{X,p}$.
Let $L$ be the quotient field of $R$. Set
$$
\omega=\frac{v-P_g(\alpha_1^a\alpha_2^b\alpha_3^c\overline x^a\overline y^b\overline z^c)}
{\overline x^d\overline y^d\overline z^f}\in \hat{(R_q)}\cap L= R_q
$$
for all maximal ideals $q$ of $R$. After possibly permuting $x,y,z$, we can assume that
$h=ae-bd\ne 0$. Set
$$
\tilde x=(\gamma^e\omega^{-b})^{\frac{1}{h}}\overline x,
\tilde y = (\gamma^{-d}\omega^a)^{\frac{1}{h}}\overline y.
$$
Set $R_1=R[(\gamma^e\omega^{-b})^{\frac{1}{h}},(\gamma^{-d}\omega^a)^{\frac{1}{h}}]$.
$\tilde x,\tilde y,\overline z\in R_1$.
$$
\begin{array}{ll}
u&=(\tilde x^a\tilde y^b \overline z^c)^m\\
v&=P_g(\tilde x^a\tilde y^b\overline z^c)+\tilde x^d\tilde y^e\overline z^f
\end{array}
$$
$(\tilde x,\tilde y,\overline z)=m_pR_1$. By Lemma \ref{Lemma1035} there exists an
\'etale neighborhood $U$ of $p$ such that $(\tilde x,\tilde y,\overline z)$ are uniformizing
parameters on $U$.

The arguments for the remaining cases 5., 6. and 7. are easier.
\end{pf}

\begin{Remark}\label{Remark1037}
Suppose that $p\in X$ is a prepared 3 point, so that
$$
\begin{array}{ll}
u&=(x^ay^bz^c)^m\\
v&=P(x^ay^bz^c)+x^dy^ez^f,
\end{array}
$$
$u=0$ is a local equation of $E_X$ and 
$$
\text{rank}\left(\begin{array}{lll}
a&b&c\\
d&e&f
\end{array}\right)=2.
$$
Then at most one of $ae-bd, af-cd,bf-ce$ is zero.
\end{Remark}

\begin{pf} By assumption, $a,b,c$ are all nonnegative.
Suppose that two of these forms are zero. After permuting $x,y,z$, we may
assume that $ae-bd=0$ and $af-cd=0$. Then
$e=\frac{bd}{a}$, $f=\frac{cd}{a}$ and 
$bf-ce=\frac{bcd}{a}-\frac{cbd}{a}=0$, a contradiction.
\end{pf}

\begin{Definition}\label{Def59}
Suppose that $\Phi:X\rightarrow S$ is strongly prepared with respect to $D_S$.
 Suppose that $p\in E_X$  We will say that $p$ is a good point for $\Phi$ if there exist permissible
parameters $(u,v)$ at $\Phi(p)$ and $*$-permissible parameters $(x,y,z)$ at $p$ for $(u,v)$ such that one of the following
forms hold:

$p$ is a 3 point, $u=0$ is a local equation of $E_X$ at $p$ and 
\begin{equation}\label{eq105}
\begin{array}{ll}
u& =x^ay^bz^c\\
v&=x^dy^ez^f
\end{array}
\end{equation}
with 
$$
\text{rank}\left(\begin{array}{lll}a&b&c\\d&e&f\end{array}\right)=2
$$

$p$ is a 3 point, $uv=0$ is a local equation of $E_X$ at $p$, 
\begin{equation}\label{eq1053}
\begin{array}{ll}
u&= x^ay^b\\
v&=z^c
\end{array}
\end{equation}

$p$ is a 3 point, $uv=0$ is a local equation of $E_X$ at $p$ 
\begin{equation}\label{eq1058}
\begin{array}{ll}
u&= x^ay^b\\
v&=y^cz^d
\end{array}
\end{equation}
with $a,b,c,d>0$.

$p$ is a 2 point, $u=0$ is a local equation of $E_X$ at $p$,  
\begin{equation}\label{eq106}
\begin{array}{ll}
u&=x^ay^b\\
v&=x^cy^d
\end{array}
\end{equation}
with $ad-bc\ne 0$

$p$ is a 2 point, $u=0$ is a local equation of $E_X$ at $p$ and there exists $\alpha\in k$ such that 
\begin{equation}\label{eq610}
\begin{array}{ll}
u&=(x^ay^b)^m\\
v&=\alpha(x^ay^b)^t+(x^ay^b)^tz
\end{array}
\end{equation}
with $(a,b)=1$.

$p$ is a 2 point, $u=0$ is a local equation of $E_X$ at $p$ 
\begin{equation}\label{eq1051}
\begin{array}{ll}
u&=x^ay^b\\
v&=x^cy^dz
\end{array}
\end{equation}
with $ad-bc\ne 0$

$p$ is a 2 point, $uv=0$ is a local equation of $E_X$ at $p$ 
\begin{equation}\label{eq1052}
\begin{array}{ll}
u&=x^a\\
v&=y^b
\end{array}
\end{equation}

$p$ is a 1 point, $u=0$ is a local equation of $E_X$ at $p$  and there exists $\alpha\in k$ such that 
\begin{equation}\label{eq107}
\begin{array}{ll}
u&=x^a\\
v&=\alpha x^c+x^cy
\end{array}
\end{equation}
$p\in X$ will be called a bad point if $p$ is not a good point.
\end{Definition}

\begin{Remark}\label{Remark1072}(Remark1072)
Suppose that $p\in X$ is a good point of one of the forms (\ref{eq105}), (\ref{eq1053}), (\ref{eq1058}), (\ref{eq106}),
(\ref{eq610}), (\ref{eq1051}), (\ref{eq1052}) or (\ref{eq107}). Then (as in Lemma \ref{Lemma1036}) there exist $*$-permissible parameters $(x,y,z)$
at $p$ such that $(x,y,z)$ are uniformizing parameters on an \'etale neighborhood
of $p$, and one of the forms (\ref{eq105}), (\ref{eq1053}), (\ref{eq1058}), (\ref{eq106}),
(\ref{eq610}), (\ref{eq1051}), (\ref{eq1052}) or (\ref{eq107}) hold.
\end{Remark}

Suppose that $p\in X$ is a 1 point and $(u,v)$ are permissible parameters at $\Phi(p)$,
$(x,y,z)$ are $*$-permissible parameters at $p$ for $(u,v)$ such that 
$$
\begin{array}{ll}
u&=x^a\\
v&=P(x)+x^cy.
\end{array}
$$
with $\text{deg }(P)\le c$.
Set $d=\text{ord }(P)\in {\bold N}\cup\{\infty\}$. 

Suppose that $(u_1,v_1)$ are also $*$-permissible parameters at $\Phi(p)$ and $(x_1,y_1,z_1)$
are permissible parameters at $p$ for $(u_1,v_1)$ such that
$$
\begin{array}{ll}
u_1&=x_1^{a_1}\\
v_1&=P_1(x_1)+x_1^{c_1}y_1.
\end{array}
$$
with $\text{deg }(P_1)\le c_1$.
Set $d_1=\text{ord }(P_1)\in {\bold N}\cup\{\infty\}$. 

We will compare $a,c,d$ and $a_1,c_1,d_1$.

If $(u,v)=(u_1,v_1)$ then there exists an $a$-th root of unity $\omega\in k$ such that
$x=\omega x_1$, so that $c=c_1$, and 
$$
P_1(x_1)=P(\omega x_1).
$$
Thus $a=a_1$, $c=c_1$, $d=d_1$.

Suppose that $(u,v)$ and $(u_1,v_1)$ are related by a change of parameters of the
type of Case 1.1 of the proof of Lemma \ref{Lemma300}. This case can only occur
if $\Phi(p)$ is a 2 point.
We have $v_1=u$ and $u_1=v$. Then $d=\text{ord}(P)\le c$. The analysis of Case 1.1 in Lemma
\ref{Lemma300} shows that there are $*$-permissible parameters $(\overline x,\overline y,\overline z)$ for $(u_1,v_1)$ such that 
$$
\begin{array}{ll}
u_1&=v=\overline x^d\\
v_1&=u=\overline P(\overline x)+\overline x^{a+c-d}\overline y
\end{array}
$$
where $\text{ord}(\overline P)=a$.
Thus $a_1=d$, $c_1=a+c-d$ and $d_1=a$.

Suppose that $(u,v)$ and $(u_1,v_1)$ are related by a change of parameters of the
type of Case 1.2 of the proof of Lemma \ref{Lemma300}. 
We have $u_1=\alpha u$ and $v_1=v$ where $\alpha(u,v)$ is a unit series. 
 The analysis of Case 1.2 in Lemma
\ref{Lemma300} shows that there are $*$-permissible parameters $(\overline x,\overline y,\overline z)$ for $(u_1,v)$ such that 
$$
\begin{array}{ll}
u_1&=\overline x^a\\
v&=\overline P(\overline x)+\overline x^{c}\overline y
\end{array}
$$
where $\text{ord}(\overline P)=d$.
Thus $a_1=a$, $c_1=c$ and $d_1=d$.

Suppose that $(u,v)$ and $(u_1,v_1)$ are related by a change of parameters of the
type of Case 1.3 of the proof of Lemma \ref{Lemma300}. 
We have $u_1=u$ and $v_1=\alpha u+\beta v$ where $\alpha(u,v), \beta(u,v)$ are series, $\beta$ is a unit series. If $\Phi(p)$ is a 2 point then $\alpha=0$. The analysis of Case 1.3 in Lemma
\ref{Lemma300} shows that there are $*$-permissible parameters $(\overline x,\overline y,\overline z)$ for $(u_1,v_1)$ such that 
$$
\begin{array}{ll}
u_1&=x^a\\
v_1&=\overline P(x)+x^{c}\overline y
\end{array}
$$
such that
$$
\begin{array}{ll}
\overline P(x)&=\sum \alpha_{ij}x^{a(i+1)}P(x)^j+\sum \beta_{ij}x^{a_i}P(x)^{j+1}\\
&\equiv \beta_{00}P(x)+\sum\alpha_{i0}x^{a(i+1)}\text{ mod }x^{d+1}.
\end{array}
$$
$a_1=a$, $c_1=c$, $d_1\le d$ if $a\not\,\mid d$.
If $\alpha=0$, we have $a_1=a$, $c_1=c$, $d_1=d$.

Suppose that $E$ is a component of $E_X$, $p\in E$, $f\in\hat{\cal O}_{X,p}$, $x=0$ is a local
equation of $E$ at $p$. Then define 
$$
\nu_E(f)=\text{max }\{n\text{ such that }x^n\mid f\}.
$$

\begin{Definition}\label{Def61}(Def61)
Suppose that $p\in X$ is a 1 point, and $E$ is the component of $E_X$ containing $p$.
Suppose that $(u,v)$ are permissible parameters at $\Phi(p)$ such that $u=0$ is a local equation of $E$ at $p$. If $(x,y,z)$ are $*$-permissible parameters at $p$ for $(u,v)$, then there is an expression
$$
\begin{array}{ll}
u&=x^a\\
v&=P(x)+x^cy.
\end{array}
$$
For fixed $(u,v)$, $a,c$ and $\nu_E(v)$ are independent of the choice of permissible
parameters $(x,y,z)$ for $(u,v)$.
Define
$$
A(\Phi,p)=\text{ min }(c-\nu_E(v))
$$
where the minimum is over permissible parameters $(u,v)$ at $\Phi(p)$ such that
$u=0$ is a local equation of $E$ at $p$.

If $A(\Phi,p)>0$, define
$$
C(\Phi,p)=\text{ min }(c-\nu_E(v),\nu_E(v)+a)
$$
where the minimum (in the lexicographic order) is over permissible parameters
$(u,v)$ at $\Phi(p)$ such that $u=0$ is a local equation of $E$ at $p$.
\end{Definition}

Suppose that $E$ is a component of $E_X$, $p\in E$ is a 1 point. Suppose that
$(u,v)$ are permissible parameters for $\Phi(p)=q$ such that $u=0$ is a local
equation of $E$ at $p$, $(x,y,z)$ are $*$-permissible parameters for $(u,v)$ at $p$. There is an 
expression 
\begin{equation}\label{eq1066}
\begin{array}{ll}
u&=x^a\\
v&=P(x)+x^cy.
\end{array}
\end{equation}
$c>0$ is equivalent to $\Phi(E)=q$. $c=0$ is equivalent to $\Phi(E)$ is a component of
$D_S$ with local equation $u=0$ at $q$.

Suppose that $\Phi(E)=q$ is a 1 point on $S$. By the discussion before Definition \ref{Def61},
$A(\Phi,p)=c-\nu_E(v)$ if and only if $a\not\,\mid \text{ ord}(P)$ or $c=\text{ord}(P)$.
If $P(x)=\sum a_ix^i$, we can make a permissible change of parameters at $q$, replacing
$v$ with $v-\sum a_{ia}u^i$ to achieve
$A(\Phi,p)=c-\nu_E(v)$.

Suppose that $\Phi(E)=q$ is a 2 point on $S$. By the discussion before Definition
\ref{Def61},
$A(\Phi,p)=c-\nu_E(v)$. 

Suppose that $\Phi(E)$ is a component $D$ of $D_S$. This is equivalent to  $c=0$ in (\ref{eq1066}).
Then
$$
0=A(\Phi,p)=c-\nu_E(v).
$$
In all these cases, if $A(\Phi,p)=c-\nu_E(v)>0$, then we have
$$
C(\Phi,p)=(c-\nu_E(v),a+\nu_E(v)).
$$
and there exists an open neighborhood $U$ of $p$ such that $A(\Phi,p')=A(\Phi,p)$ for
all $p'\in E\cap U$ and $C(\Phi,p')=C(\Phi,p)$ if $A(\Phi,p)>0$. Then $A(\Phi,p')=A(\Phi,p)$
and $C(\Phi,p')=C(\Phi,p)$  at all 1 points $p'\in E$. We can
then define
$$
A(\Phi,E)=A(\Phi,p)
$$
and
$$
C(\Phi,E)=C(\Phi,p)
$$
for $p\in E$ a 1 point.
\begin{Lemma}\label{Lemma1067}
Suppose that $p\in X$ is a 2 point and $E_1$, $E_2$ are the components of $E_X$
containing $p$. Then there exist permissible parameters $(u,v)$ at $q=\Phi(p)$ and
permissible parameters $(x,y,z)$ for $(u,v)$ at $p$ such that, if $p$
satisfies (\ref{eq1061}) of Definition \ref{Def57},
$$
\begin{array}{ll}
u&=(x^ay^b)^k\\
v&=P(x^ay^b)+x^cy^d
\end{array}
$$
or if $p$ satisfies (\ref{eq1062}) of Definition \ref{Def57},
$$
\begin{array}{ll}
u&=(x^ay^b)^k\\
v&=P(x^ay^b)+x^cy^dz
\end{array}
$$
where $x=0$ is a local equation of $E_1$, $y=0$ is a local equation of $E_2$, then
$$
A(\Phi,E_1 )=c-\nu_{E_1}(v), A(\Phi,E_2)=d-\nu_{E_2}(v).
$$
If $A(\Phi,E_1)>0$ then
$$
C(\Phi,E_1)=(c-\nu_{E_1}(v),\nu_{E_1}(v)+ak).
$$
If $A(\Phi,E_2)>0$ then
$$
C(\Phi,E_2)=(d-\nu_{E_2}(v),\nu_{E_2}(v)+bk).
$$

Suppose that $p\in X$ is a 3 point,  $p$ satisifes (\ref{eq1063}) of Definiton
\ref{Def57}, and
$E_1$, $E_2$, $E_3$ are the components of $E_X$ containing $p$. Then there exist
permissible parameters $(u,v)$ at $q=\Phi(p)$ and permissible parameters
$(x,y,z)$ for $(u,v)$ at $p$ such that 
$$
\begin{array}{ll}
u&=(x^ay^bz^c)^m\\
v&=P(x^ay^bz^c)+x^dy^ez^f
\end{array}
$$
where $x=0$ is a local equation of $E_1$, $y=0$ is a local equation of $E_2$,
$z=0$ is a local equation of $E_3$, and
$$
A(\phi,E_1)=d-\nu_{E_1}(v),
A(\Phi,E_2)=e-\nu_{E_2}(v),
A(\Phi,E_3)=f-\nu_{E_3}(v).
$$
If $A(\Phi,E_1)>0$, then 
$$
C(\Phi,E_1)=(d-\nu_{E_1}(v),\nu_{E_1}(v)+am),
$$
If $A(\Phi,E_2)>0$, then 
$$
C(\Phi,E_2)=(e-\nu_{E_2}(v),\nu_{E_2}(v)+bm),
$$
If $A(\Phi_{E_3})>0$, then 
$$
C(\Phi,E_3)=(f-\nu_{E_3}(v),\nu_{E_3}(v)+cm),
$$
\end{Lemma}

\begin{pf} Suppose that $p\in X$ is a 2 point satisfying (\ref{eq1061}),
$(u,v)$ are permissible parameters at $q$ and $(x,y,z)$ are uniformizing parameters
for $(u,v)$ at $p$ such that
$$
\begin{array}{ll}
u&=(x^ay^b)^k\\
v&=P(x^ay^b)+x^cy^d
\end{array}
$$
and $(x,y,z)$ are uniformizing parameters 
on an \'etale neighborhood of $p$. Let
$P(t)=\sum a_it^i$. If $q$ is a 1 point, then we can replace $v$ with 
$v-\sum_i a_{ik}u^i$, so that $k\not\,\mid \text{ord }(P)$.

If $c=0$, then $0=A(\Phi,E_1)=c-\nu_{E_1}(v)$, and if $d=0$, then
$0=A(\Phi,E_2)=d-\nu_{E_2}(v)$.

Suppose that $c>0$. Then $\Phi(E_1)=q$. If $p'$ is a 1 point on $E_1$ near $p$ then
there exist $*$-permissible parameters $(\overline x,\overline y, \overline z)$ at $p'$
such that
$$
\begin{array}{ll}
u&=\overline x^{ak}\\
v&=P_{p'}(\overline x)+\overline x^c\overline y
\end{array}
$$
where $P_{p'}(\overline x)=P(\overline x^a)+\alpha\overline x^c$
for some nonzero $\alpha\in k$.

If $q$ is a 1 point we have $ak\not\,\mid \text{ord }(P_{p'})$ or
$c=\text{ord }(P_{p'})$.
By the discussion before Definition \ref{Def61}, we have that $A(\Phi,E_1)=c-\nu_{E_1}(v)$,
and if $A(\Phi,E_1)>0$, then
$$
C(\Phi,E_1)=(c-\nu_{E_1}(v),\nu_{E_1}(v)+ak).
$$

A similar argument shows that $A(\Phi,E_2)=d-\nu_{E_2}(v)$ if $d>0$, 
and if $A(\Phi,E_2)>0$, then
$$
C(\Phi,E_2)=(d-\nu_{E_2}(v),\nu_{E_2}(v)+bk).
$$

If $p$ satisfies (\ref{eq1062}) or (\ref{eq1063}) then the proof is similar.
\end{pf}

\begin{Remark}\label{Remark1081} If $p$ is a 1 point then
$A(\Phi,p)=0$ if and only if $p$ is a good point.
\end{Remark}

Set 
$$
A(\Phi)=\text{ max }\{A(\Phi,E)\mid E \text{ is a component of }E_X\}.
$$
If $A(\Phi)>0$, define
$$
C(\Phi)=\text{ max }\{C(\Phi,E)\mid E \text{ is a component of }E_X\}.
$$

\begin{Lemma}\label{Lemma1038}
Suppose that $p\in X$ is a 1 point, $(u,v)$ are permissible parameters at $\Phi(p)$
such that $u=0$ is a local equation of $E_X$ at $p$, $(x,y,z)$ are $*$-permissible
parameters at $p$ for $(u,v)$ such that
$$
\begin{array}{ll}
u&=x^a\\
v&=P(x)+x^cy
\end{array}
$$
with $\text{deg}(P)\le c$.
Set $d=\text{ord}(P)\in{\bold N}\cup\{\infty\}$.
\begin{enumerate}
\item Suppose that $\Phi(p)$ is a 1 point. Then $p$ is a bad point if $d<c$ and 
$a\not\,\mid d$.
\item Suppose that $\Phi(p)$ is a 2 point. Then $p$ is a bad point if $d<c$.
\end{enumerate}
\end{Lemma}

\begin{pf}
Suppose that $(u_1,v_1)$ are permissible parameters at $q=\Phi(p)$ such that
$u_1=0$ is a local equation of $E_X$ at $p$, and $(u_1,v_1)$ realize $p$ as a bad point.

If $q$ is a 1 point then there exist  series $\overline \alpha$, $\overline \beta$,
$\overline \gamma$ in $u,v$ such that
$$
\begin{array}{ll}
u_1&=\overline\alpha u\\
v_1&=\overline \beta u+\overline \gamma v.
\end{array}
$$
Thus $(u_1,v_1)$ is obtained by transformations of the form of Case 1.2 and Case 1.3
of Lemma \ref{Lemma300}. The conclusions of the Lemma now follow from the
analysis preceeding Definition \ref{Def61}.

Suppose that $q$ is a 2 point. Then there exist unit series $\overline \alpha$,
$\overline \beta$ in $u,v$ such that
$$
\begin{array}{ll}
u_1&=\overline \alpha u\\
v_1&=\overline \beta v
\end{array}
$$
or
$$
\begin{array}{ll}
u_1&=\overline\alpha v\\
v_1&=\overline\beta u.
\end{array}
$$
In the first case we have, with the notation preceeding Definition \ref{Def61}, that
$a_1=a$, $c_1=c$ and $d_1=d$ so that $d_1<c_1$. In the second case we have
$a_1=d$, $c_1=a+c-d$ and $d_1=a$ so that $d_1<c_1$. $p$ is thus  a bad point.
\end{pf}

\begin{Theorem}\label{Theorem60}
Suppose that $\Phi:X\rightarrow S$ is strongly prepared. Then
the locus of bad points in $X$ is a Zariski closed set of pure codimension 1, consisting of 
a union of components of $E_X$.
\end{Theorem}

\begin{pf}
We will first show that 
the good points of $X$ are a Zariski open set in $E_X$.

Suppose that $p\in X$ is a good 3 point. Then
there exists an open neighborhood $U$ of $p$, uniformizing parameters
$(x,y,z)$ in an \'etale cover of $U$ such that $u=0$ is a local equation of $E_X$ in $U$
and
$$
\begin{array}{ll}
u&=x^ay^bz^c\\
v&=x^dy^ez^f.
\end{array}
$$
If $q\in U$
is a 2 point, then we have (after possibly permuting $x,y,z$) that $(x,y,z_1)$ are regular
parameters at $q$ where $z=z_1+\alpha$ (with $\alpha\ne 0$). Set $x=x_1(z_1+\alpha)^{-\frac{c}{a}}$.
Then $(x_1,y,z_1)$ are permissible parameters at $q$, and
$$
\begin{array}{ll}
u&=x_1^ay^b\\
v&=x_1^dy^e(z_1+\alpha)^{f-\frac{cd}{a}}.
\end{array}
$$
If $ae-bd\ne 0$, we can make a permissible change of variables $(\tilde x,\tilde y,\tilde z)$ at $q$ to get
$$
\begin{array}{ll}
u&=\tilde x^a\tilde y^b\\
v&=\tilde x^d\tilde y^e.
\end{array}
$$
If $ae-bd=0$ then $f-\frac{cd}{a}\ne0$, so that we can make a permissible change of parameters  to get
$$
\begin{array}{ll}
u&=(x_1^{a_1}y^{b_1})^k\\
v&=\beta(x_1^{a_1}y^{b_1})^t+(x_1^{a_1}y^{b_1})^tz_1.
\end{array}
$$
If $q\in U$
is a 1 point, then we have (after possibly permuting $x,y,z$) that $(x,y_1,z_1)$ are regular
paramaters at $q$ where $y=y_1+\alpha$, $z=z_1+\beta$ (with $\alpha,\beta\ne 0$). Set $x=x_1(y_1+\alpha)^{-\frac{b}{a}}(z_1+\beta)^{-\frac{c}{a}}$.
Then $(x_1,y_1,z_1)$ are permissible parameters at $q$, and
$$
\begin{array}{ll}
u&=x_1^a\\
v&=\gamma x_1^d + x_1^d(\gamma_1y_1+\gamma_2z_1+\cdots)
\end{array}
$$
where $\gamma,\gamma_1,\gamma_2\in k$, $\gamma\ne 0$ and either $\gamma_1\ne0$ or $\gamma_2\ne0$ since we cannot have both 
$e-\frac{db}{a}=0$ and $f-\frac{dc}{a}=0$.
 Thus all points in $U$ are good points.

Suppose that $p\in X$ is a good 2 point and (\ref{eq610}) holds at $p$.
Then
there exists an open neighborhood $U$ of $p$, uniformizing parameters
$(x,y,z)$ in an \'etale cover of $U$ such that $u=0$ is a local equation of $E_X$ in $U$
and
$$
\begin{array}{ll}
u&=(x^ay^b)^k\\
v&=\overline \beta(x^ay^b)^t+(x^ay^b)^tz
\end{array}
$$
If $q\in U$ is a 2 point, then we have that $(x,y,z_1)$ are permissible parameters at
$q$ where $z=z_1+\alpha$ and $q$ is a good point.

If $q\in U$
is a 1 point, then we have (after possibly permuting $x,y$) that $(x,y_1,z_1)$ are regular
parameters at $q$ where $y=y_1+\alpha$, $z=z_1+\beta$ (with $\alpha\ne 0$). Set $x=x_1(y_1+\alpha)^{-\frac{b}{a}}$.
Then $(x_1,y_1,z_1)$ are permissible parameters at $q$, and
$$
\begin{array}{ll}
u&=x_1^{ak}\\
v&=(\overline\beta +\beta)x_1^{at}+x_1^{at}z_1
\end{array}
$$
Thus all points in $U$ are good points.

Suppose that $p\in X$ is a good 2 point, and (\ref{eq106}) holds at $p$. 
Then
there exists an open neighborhood $U$ of $p$, uniformizing parameters
$(x,y,z)$ in an \'etale cover of $U$ such that $u=0$ is a local equation of $E_X$ in $U$
and
$$
\begin{array}{ll}
u&=x^ay^b\\
v&=x^cy^d
\end{array}
$$
where $ad-bc\ne0$.
If $q\in U$
is a 2 point, then we have  that $(x,y,z_1)$ are permissible parameters
 at $q$ where $z=z_1+\alpha$, and $q$ is a good point.

If $q\in U$
is a 1 point, then we have (after possibly permuting $x,y$) that $(x,y_1,z_1)$ are regular
paramaters at $q$ where $y=y_1+\alpha$, $z=z_1+\beta$ (with $\alpha\ne 0$). Set $x=x_1(y_1+\alpha)^{-\frac{b}{a}}$.
Set $\gamma=\alpha^{d-\frac{bc}{a}}$, $\tilde y_1=(y_1+\alpha)^{d-\frac{bc}{a}}-\gamma$.
Then $(x_1,\tilde y_1,z_1)$ are permissible parameters at $q$, and
$$
\begin{array}{ll}
u&=x_1^{a}\\
v&=\gamma x_1^c+x_1^c\tilde y_1
\end{array}
$$
Thus all points in $U$ are good points.

If $p$ is a good  point satisfying (\ref{eq1053}), (\ref{eq1058}), (\ref{eq106}), 
(\ref{eq1051}), (\ref{eq1052}) or (\ref{eq107}),
 a similar argument shows that there is a Zariski open neighborhood $U$ of $p$ of good points.
\vskip .2truein
We will now show that the bad points of $X$ have pure codimension 1 in $X$.
It suffices to show that any bad point lies on a surface of bad points.

First suppose that $p$ is a bad 3 point. 
Then
there exists an open neighborhood $U$ of $p$, uniformizing parameters
$(x,y,z)$ in an \'etale cover of $U$ such that $u=0$ is a local equation of $E_X$ in $U$
and
$$
\begin{array}{ll}
u&=(x^ay^bz^c)^k\\
v&=P(x^ay^bz^c)+x^dy^ez^f
\end{array}
$$
where (after possibly permuting $x,y,z$) we have
$$
\text{max}\{\frac{d}{a},\frac{e}{b},\frac{f}{c}\} = \frac{f}{c}
$$
Thus $\text{ord}(P)<\frac{f}{c}$, since $\text{ord }(P)\ge \frac{f}{c}$ implies that $x^dy^ez^f \vert P(x^ay^bz^c)$,
and $p$ is thus a good point.

If $\Phi(p)$ is a 1 point, we can make a permissible change of parameters
so that we have that $k\not\,\mid \text{ord }(P)$.

Let $q\in U$
be a 1 point on the surface $z=0$. $c,f>0$ imply $z=0$ is a local equation of a
component of $E_X$ which maps to $\Phi(p)$. There are regular parameters $(x_1,y_1,z)$ at $q$ where
$x=x_1+\alpha$, $y=y_1+\beta$ with $\alpha,\beta\ne 0$. There are permissible parameters $(x_1,y_1,z_1)$ at $q$ where
$$
z=(x_1+\alpha)^{-\frac{a}{c}}(y_1+\beta)^{-\frac{b}{c}}z_1
$$
$$
\begin{array}{ll}
u&=z_1^{ck}\\
v&=P(z_1^c)+ (x_1+\alpha)^{d-\frac{fa}{c}}(y_1+\beta)^{e-\frac{fb}{c}}z_1^f\\
&=P(z_1^c)+\alpha^{d-\frac{fa}{c}}\beta^{e-\frac{fb}{c}}z_1^f+z_1^f(\gamma_1 x_1+\gamma_2y_1+\cdots)
\end{array}
$$
where $\gamma_1,\gamma_2\in k$ and $\gamma_1$ or $\gamma_2\ne0$. 

Suppose that $\Phi(p)$ is a 1 point. Then $\Phi(q)=\Phi(p)$ is a 1 point.
$q$ is a bad point by Lemma \ref{Lemma1038}, since $ck\not\,\mid c\text{ ord}(P)$ and $c\text{ ord}(P)<f$

Suppose that $\Phi(p)$ is a 2 point. Then $\Phi(q)=\Phi(p)$ is a 2 point. 
$q$ is a bad point by Lemma \ref{Lemma1038} since $c\text{ ord }(P)<f$. 

Suppose that $p$ is a bad 2 point satisfying (\ref{eq1062}). 
There exists an open neighborhood $U$ of $p$ and uniformizing parameters
$(x,y,z)$ on an \'etale cover of $U$ such that
$$
\begin{array}{ll}
u&= (x^ay^b)^k\\
v&=P(x^ay^b)+x^cy^dz.
\end{array}
$$
We can (after possibly permuting $x,y$) assume that $ad-bc\ge 0$. Since $p$ is a bad point, $\text{ord }(P)<\frac{d}{b}$.
If $\Phi(p)$ is a 1 point we can make a permissible change of parameters so that $k\not\,\mid \text{ord }(P)$.
 Let $q\in U$
be a 1 point on the surface $y=0$. $b,d>0$ implies $y=0$ is a local equation of a
component of $E_X$ which maps to $\Phi(p)$. There are regular parameters $(x_1,y_1,z)$ at $q$ where 
$x=x_1+\alpha$, $z=z_1+\beta$ (with $\alpha\ne 0$). There are permissible parameters $(x_1,y_1,z_1)$ at $q$ where
$$
y=(x_1+\alpha)^{-\frac{a}{b}}y_1
$$
$$
\begin{array}{ll}
u&=y_1^{bk}\\
v&=P(y_1^b)+ (x_1+\alpha)^{c-\frac{da}{b}}(z_1+\beta)y_1^d\\
&=P(y_1^b)+\alpha^{c-\frac{da}{b}}\beta y_1^d+\tilde z_1y_1^d
\end{array}
$$
$\Phi(q)=\Phi(p)$ so that $\Phi(q)$ is a 1 point if and only if $\Phi(p)$ is a
1 point. Since $b\text{ ord }(P)<d$ and $bk\not\,\mid b\text{ ord }(P)$ if $\Phi(q)$ is a 1 point, $q$ is
a bad point by Lemma \ref{Lemma1038}.

A similar argument shows that there is a surface of bad points passing through a bad  point
satisfying (\ref{eq1061}) or (\ref{eq1060}).
\end{pf}

\begin{Lemma}\label{Lemma62}(Lemma62) Suppose that 
$\Phi:X\rightarrow S$ is strongly prepared, $q\in D_S$  and $p\in \Phi^{-1}(q)$
is such that one of the forms 1. - 7. of Lemma \ref{Lemma1036} hold at $p$. Then
$m_q{\cal O}_{X,p}$ is not invertible if and only if    one of the following holds:

$p$ is a 1 point 
\begin{equation}\label{eq612}
\begin{array}{ll}
u&=x^k\\
v&=x^cy
\end{array}
\end{equation}
with $c<k$. 

$p$ is a 2 point 
\begin{equation}\label{eq611}
\begin{array}{ll}
u&=(x^ay^b)^k\\
v&=P(x^ay^b)+x^cy^d
\end{array}
\end{equation}
with $a,b>0$, $(a,b)=1$, $ad-bc\ne 0$,
$$
\text{min}\{\frac{c}{a},\frac{d}{b}\}<\text{ord }(P)<\text{max}\{\frac{c}{a},\frac{d}{b}\},
$$
$$
\text{min}\{\frac{c}{a},\frac{d}{b}\}<k.
$$

$p$ is a 2 point 
\begin{equation}\label{eq613}
\begin{array}{ll}
u&=(x^ay^b)^k\\
v&=x^cy^d
\end{array}
\end{equation}
with $a,b>0$, $(a,b)=1$, $ad-bc\ne0$, 
$$
\text{min}\{\frac{c}{a},\frac{d}{b}\}<k<\text{max}\{\frac{c}{a},\frac{d}{b}\}.
$$

$p$ is a 2 point 
\begin{equation}\label{eq108}
\begin{array}{ll}
u&=(x^ay^b)^k\\
v&=P(x^ay^b)+x^cy^dz
\end{array}
\end{equation}
with $a,b>0$, $(a,b)=1$, $ad-bc\ne 0$, 
$$
\text{min}\{\frac{c}{a},\frac{d}{b}\}<\text{ord }(P)<\text{max}\{\frac{c}{a},\frac{d}{b}\},
$$
$$
\text{min}\{\frac{c}{a},\frac{d}{b}\}<k.
$$

$p$ is a 2 point 
\begin{equation}\label{eq109}
\begin{array}{ll}
u&=(x^ay^b)^k\\
v&=x^cy^dz
\end{array}
\end{equation}
with $a,b>0$, $(a,b)=1$, $ad-bc\ne0$,
$$
\text{min}\{\frac{c}{a},\frac{d}{b}\}<k.
$$

$p$ is a 2 point 
\begin{equation}\label{eq110}
\begin{array}{ll}
u&= (x^ay^b)^k\\
v&=(x^ay^b)^tz
\end{array}
\end{equation}
with $a,b>0$, $(a,b)=1$, $t<k$.

$p$ is a 2 point 
\begin{equation}\label{eq1054}
\begin{array}{ll}
u&= x^a\\
v&=y^b
\end{array}
\end{equation}

$p$ is a 3 point 
\begin{equation}\label{eq111}
\begin{array}{ll}
u&=(x^ay^bz^c)^k\\
v&=P(x^ay^bz^c)+x^dy^ez^f
\end{array}
\end{equation} 
with $a,b,c>0$, $(a,b,c)=1$,
$$
\text{min}\{\frac{d}{a},\frac{e}{b},\frac{f}{c}\}<\text{ord }(P)
<\text{max}\{\frac{d}{a},\frac{e}{b},\frac{f}{c}\},
$$
$$
k>\text{min}\{\frac{d}{a},\frac{e}{b},\frac{f}{c}\}.
$$

$p$ is a 3 point 
\begin{equation}\label{eq112}
\begin{array}{ll}
u&=(x^ay^bz^c)^k\\
v&=x^dy^ez^f
\end{array}
\end{equation}
with $a,b,c>0$, $(a,b,c)=1$,
$$
\text{min}\{\frac{d}{a},\frac{e}{b},\frac{f}{c}\}<k<\text{max}\{\frac{d}{a},\frac{e}{b},\frac{f}{c}\}.
$$

$p$ is a 3 point 
\begin{equation}\label{eq1055}
\begin{array}{ll}
u&=x^ay^b\\
v&=z^c
\end{array}
\end{equation}
with $a,b,c>0$.

$p$ is a 3 point 
\begin{equation}\label{eq1059}
\begin{array}{ll}
u&=x^ay^b\\
v&=y^cz^d
\end{array}
\end{equation}
with $a,b,c,d>0$.

\end{Lemma}

\begin{pf}
Suppose that $p$ is a 1 point. Then (\ref{eq612}) follows easily.

Suppose that $p$ is a 2 point with 
$$
\begin{array}{ll}
u&=(x^ay^b)^k\\
v&=P(x^ay^b)+x^cy^d,
\end{array}
$$
$P\ne 0$ and $e=\text{ord }(P)<\text{max}\{\frac{c}{a},\frac{d}{b}\}$.
Set $\lambda_2=\text{max}\{\frac{c}{a},\frac{d}{b}\}$,
$\lambda_1=\text{min}\{\frac{c}{a},\frac{d}{b}\}$.

$u\mid v$ if and only if $e\ge k$ and $\lambda_1\ge k$. $v\mid u$ if and only if $e\le\lambda_1$
and $e\le k$. Thus $(u,v){\cal O}_{X,p}$ is not invertible if and only if
$\lambda_1<k$ and $\lambda_1<e$.

Suppose that $p$ is a 2 point with 
$$
\begin{array}{ll}
u&=(x^ay^b)^k\\
v&=x^cy^d
\end{array}
$$
Set $\lambda_1=\text{min}\{\frac{c}{a},\frac{d}{b}\}$,
$\lambda_2=\text{max}\{\frac{c}{a},\frac{d}{b}\}$. $u\mid v$ if and only if
$k\le\lambda_1$, $v\mid u$ if and only if $k\ge \lambda_2$. So $(u,v){\cal O}_{X,p}$
is not invertible if and only if $\lambda_1<k<\lambda_2$.

Suppose that $p$ is a 2 point with 
$$
\begin{array}{ll}
u&=(x^ay^b)^k\\
v&=P(x^ay^b)+x^cy^dz
\end{array}
$$
with $ad-bc\ne 0$, $e=\text{ord }(P)< \text{max}\{\frac{c}{a},\frac{d}{b}\}$.
Set $\lambda_1=\text{min}\{\frac{c}{a},\frac{d}{b}\}$,
$\lambda_2=\text{max}\{\frac{c}{a},\frac{d}{b}\}$. $u\mid v$ if and only if
$e\ge k$ and
$k\le\lambda_1$. $v\mid u$ if and only if $e\le k$ and $e\le\lambda_1$. So $(u,v){\cal O}_{X,p}$
is not invertible if and only if $\lambda_1<k$ and $\lambda_1<e$.

Suppose that $p$ is a 2 point with 
$$
\begin{array}{ll}
u&=(x^ay^b)^k\\
v&=x^cy^dz
\end{array}
$$
and $ad-bc\ne 0$.
$(u,v)$ is invertible at $p$ if and only if $c\ge ka$ and $d\ge bk$. Thus
$(u,v)$ is not invertible at $p$ if and only if $k>\text{min}\{\frac{c}{a},\frac{d}{b}\}$, and we get (\ref{eq109}).

Suppose that $p$ is a 2 point with 
$$
\begin{array}{ll}
u&=(x^ay^b)^k\\
v&=P(x^ay^b)+(x^ay^b)^tz
\end{array}
$$
with  $P\ne 0$ and $e=\text{ord }(P)\le t$. We will show that $(u,v)$ is invertible at $p$.
If $k\le e$ then $u\vert v$. Suppose that $k>e$. There are new permissible parameters $(\overline x,\overline y,\overline z)$
such that
$$
\begin{array}{ll}
v&= (\overline x^a\overline y^b)^e\\
u&= \overline P(\overline x^a\overline y^b)+(\overline x^a\overline y^b)^{t-e+k}\overline z
\end{array}
$$
with $\text{ord }(\overline P)=k$. Thus $v \vert u$.

Suppose that $p$ is a 2 point with 
$$
\begin{array}{ll}
u&=(x^ay^b)^k\\
v&=(x^ay^b)^tz
\end{array}
$$
Then $(u,v)$ not invertible at $p$ if and only if  $t<k$, and we get (\ref{eq110}).

Suppose that $p$ is a 3 point with
$$
\begin{array}{ll}
u&=(x^ay^bz^c)^k\\
v&=P(x^ay^bz^c)+x^dy^ez^f
\end{array}
$$
with $P\ne 0$ and $\text{ord }(P)<\text{max}\{\frac{d}{a},\frac{e}{b},\frac{f}{c}\}$. Set 
$$
\lambda_2=\text{max}\{\frac{d}{a},\frac{e}{b},\frac{f}{c}\},\,\,
\lambda_1=\text{min}\{\frac{d}{a},\frac{e}{b},\frac{f}{c}\}
$$
$u\vert v$ is equivalent to $\text{ord }(P)\ge k$, $k\le \lambda_1$.
$v\vert u$ is equivalent to $\text{ord }(P)\le k$ and $\text{ord }(P)\le\lambda_1$.
That is,  $(u,v)$ is not invertible at $p$ if and only if $\text{ord }(P)>\lambda_1$ and $k>\lambda_1$.
We thus get (\ref{eq111}).

Suppose that $p$ is a 3 point with
$$
\begin{array}{ll}
u&=(x^ay^bz^c)^k\\
v&=x^dy^ez^f
\end{array}
$$
Set 
$$
\lambda_2=\text{max}\{\frac{d}{a},\frac{e}{b},\frac{f}{c}\},\,\,
\lambda_1=\text{min}\{\frac{d}{a},\frac{e}{b},\frac{f}{c}\}
$$

$u\vert v$ is equivalent to  $k\le\lambda_1$.
$v\vert u$ is equivalent to  $k\ge\lambda_2$. 

Thus $(u,v)$ is not invertible at $p$ if and only if $\lambda_1<k<\lambda_2$, and we get (\ref{eq112}).
\end{pf}

\begin{Lemma}\label{Lemma1084} Suppose that $\Phi:X\rightarrow S$ is strongly prepared. Let $S_1$ be the 
blowup of $S$ at a point $q\in D_S$. Let $U$ be the largest open set of $X$ such that
the rational map $X\rightarrow S_1$ is a morphism $\Phi_1:U\rightarrow S_1$. Then
$\Phi_1$ is strongly prepared.
\end{Lemma}

\begin{pf} This follows from the analysis of Lemma \ref{Lemma62}.
\end{pf}

\begin{Theorem}\label{Theorem61} Suppose that $\Phi:X\rightarrow S$ is strongly
prepared, $p\in X$ is a 1 point
 and the rational map $\Phi_1$ 
 from $X$ to the blow up $S_1$ of $q=\Phi(p)$ is
a morphism in a neighborhood of $p$. Then $A(\Phi_1,p)\le A(\Phi,p)$. If
$A(\Phi_1,p)=A(\Phi,p)>0$, then $C(\Phi_1,p)<C(\Phi,p)$.
\end{Theorem}

\begin{pf} 
At $p$ we have permissible parameters such that 
$$
\begin{array}{ll}
u&=x^k\\
v&=P(x)+x^cy
\end{array}
$$
and $C(\Phi,p)=(c-\nu_E(v),\nu_E(v)+k)$.

First suppose that $P\ne0$ and $e=\text{ord }(P)\le c$. If $e> k$ then we have permissible parameters $u_1,v_1$ at
$q_1=\Phi_1(p)$ such that
$$
u=u_1,
v=u_1v_1.
$$
Then
$$
\begin{array}{ll}
u_1&=x^k\\
v_1&=\frac{P(x)}{x^k}+x^{c-k}y
\end{array}
$$
$A(\Phi_1,p)\le(c-k-(e-k))=A(\Phi,p)$ and if $A(\Phi_1,p)=A(\Phi,p)$ then
$C(\Phi_1,p) \le (c-k-(e-k),e-k+k)=(c-e,e)<(c-e,e+k)=C(\Phi,p)$.

If $e=k$ then there exists $0\ne\alpha\in k$ such that
$P(x)=\alpha x^k+\cdots$. There exist permissible parameters $(u_1,v_1)$ at 
$q_1=\Phi_1(p)$ such that 
$$
u=u_1, v=u_1(v_1+\alpha).
$$
$$
\begin{array}{ll}
u_1&=x^k\\
v_1&=\frac{P(x)}{x^k}-\alpha+x^{c-k}y.
\end{array}
$$
Thus $A(\Phi_1,p)<(c-k)-(e-k)=A(\Phi,p)$.

If $e<k$ then we have permissible parameters $u_1,v_1$ at
$q_1=\Phi_1(p)$ such that
$$
u=u_1v_1,
v=v_1.
$$
We have permissible parameters $(\overline x,\overline y,z)$ at $p$ such that
$$
\begin{array}{ll}
v&=\overline x^e\\
u&= \overline P(\overline x)+\overline x^{k+c-e}\overline y
\end{array}
$$
where $\text{ord }(\overline P)=k$.

Then
$$
\begin{array}{ll}
v_1&=\overline x^e\\
u_1&=\frac{\overline P(\overline x)}{\overline x^e}+\overline x^{k+c-2e}\overline y
\end{array}
$$
$$
A(\Phi_1,p)\le (k+c-2e-(k-e))=A(\Phi,p)
$$
and if $A(\Phi_1,p)=A(\Phi,p)$
then $C(\Phi_1,p) \le (k+c-2e-(k-e),(k-e)+e)<C(\Phi,p)$.

Now suppose that $P(x)=0$. Then 
$$
\begin{array}{ll}
u&=x^k\\
v&=x^cy
\end{array}
$$
with $c\ge k$. There exist permissible parameters $(u_1,v_1)$ at $q_1=\Phi_1(p)$ such
that
$$
u=u_1, v=u_1v_1
$$
and
$$
A(\Phi_1,p)=A(\Phi,p)=0.
$$
\end{pf}

\begin{Theorem}\label{Theorem63} Suppose that 
$\Phi:X\rightarrow S$ is strongly prepared and $q\in S$. Then the locus of points $Z$ in
$X$ where $\Phi$ does not factor through the blowup of $q$ is a pure codimension 2  subscheme.
$Z$ makes SNCs with $\overline B_2(X)$ except possibly at 3 points of the form (\ref{eq111}).

Suppose that $C$ is a component of this locus 
which makes SNCs with $\overline B_2(X)$, and $\pi:X_1\rightarrow X$ is the blowup of $C$,
$E_1=\pi^{-1}(C)_{red}$, $\Phi_1=\Phi\circ \pi$.
Then $\Phi_1$ is strongly prepared and either $A(\Phi_1)=0$ or
$$
A(\Phi_1,E_1)<A(\Phi)
$$
\end{Theorem}

\begin{pf}
{\bf Suppose that $p\in X$ is a 3 point such that $m_q{\cal O}_{X,p}$ is not invertible 
and (\ref{eq111}) holds at $p$.} 
We may assume that there exists an open neighborhood $U$ of $p$ such that $(x,y,z)$ are uniformizing parameters on an \'etale cover of $U$.

 After possibly interchanging $x,y,z$, we can assume that
$$
\frac{d}{a}=\text{min}\{\frac{d}{a},\frac{e}{b},\frac{f}{c}\}
$$
and 
$$
\frac{f}{c}=\text{max}\{\frac{d}{a},\frac{e}{b},\frac{f}{c}\}.
$$

We will now determine the locus of points in $U$ where $m_q{\cal O}_U$ is not invertible. First suppose that $q'$ is a 2 point
on the curve $x=z=0$. $q'$ has regular parameters $(x,y_1,z)$ where $y=y_1+\alpha$.
Thus $q'$ has permissible parameters $(x_1,y_1,z)$ where  $x_1$ is defined by
$$
x=x_1(y_1+\alpha)^{-\frac{b}{a}}
$$
Set $\lambda = (a,c)$, $a_1=\frac{a}{\lambda}$, $c_1=\frac{c}{\lambda}$.
$$
\begin{array}{ll}
u&=(x_1^{a_1}z^{c_1})^{k\lambda}\\
v&=P((x_1^{a_1}z^{c_1})^{\lambda})+x_1^dz^f(y_1+\alpha)^{e-\frac{db}{a}}.
\end{array}
$$
We have $a_1f-c_1d>0$ and $\lambda\text{ord }(P)<\frac{f}{c_1}$.
We can make a permissible change of variables to get 
$$
\begin{array}{ll}
u&=(\overline x_1^{a_1}\overline z^{c_1})^{k\lambda}\\
v&=P((\overline x_1^{a_1}\overline z^{c_1})^{\lambda})+\overline x_1^d\overline z^f.
\end{array}
$$
$k>\frac{d}{a}$ implies $\lambda k>\frac{d}{a_1}$ and $\text{ord }(P)>\frac{d}{a}$
 implies $\lambda\text{ord }(P)>\frac{d}{a_1}$.
$q$ thus has the form of (\ref{eq611}),
and we see that $(u,v)$ is not invertible on the curve with local equations $x=z=0$.

Now suppose that $q'$ is a 2 point
on the curve $y=z=0$. $q'$ has regular parameters $(x_1,y,z)$ where $x=x_1+\alpha$.
Thus $q'$ has permissible parameters $(x_1,y_1,z)$ where  $y_1$ is defined by
$$
y=y_1(x_1+\alpha)^{-\frac{a}{b}}
$$
Set $\lambda = (b,c)$, $b_1=\frac{b}{\lambda}$, $c_1=\frac{c}{\lambda}$.
$$
\begin{array}{ll}
u&=(y_1^{b_1}z^{c_1})^{k\lambda}\\
v&=P((y_1^{b_1}z^{c_1})^{\lambda})+(x_1+\alpha)^{d-\frac{ea}{b}}y_1^ez^f
\end{array}
$$
First suppose that $bf-ce\ne 0$. Then $bf-ce>0$, and $b_1f-c_1e>0$.
Since $\lambda\text{ord }(P)<\frac{f}{c_1}$, we have by (\ref{eq611}) that $(u,v)$ is not
invertible at $q'$ if and only if $\lambda\text{ord }(P)>\frac{e}{b_1}$ and $\lambda k>\frac{e}{b_1}$ so that
$(u,v)$ is not invertible at 2 points $q'$ on $y=z=0$ if and only if $\text{ord }(P)>\frac{e}{b}$ and $k>\frac{e}{b}$.

Now suppose that $bf-ce=0$, so that $b_1f-c_1e=0$.   
Since $\lambda\text{ord }(P)<\frac{f}{c_1}$, (\ref{eq110}) cannot hold, and we then have that $(u,v)$ is
invertible at 2 points  $q'$ on $y=z=0$.

Now suppose that $q'$ is a 2 point
on the curve $x=y=0$. $q'$ has regular parameters $(x,y,z_1)$ where $z=z_1+\alpha$.
Thus $q'$ has permissible parameters $(x_1,y,z_1)$ where  $x_1$ is defined by
$$
x=x_1(z_1+\alpha)^{-\frac{c}{a}}
$$
Set $\lambda = (a,b)$, $a_1=\frac{a}{\lambda}$, $b_1=\frac{b}{\lambda}$.
$$
\begin{array}{ll}
u&=(x_1^{a_1}y^{b_1})^{k\lambda}\\
v&=P((x_1^{a_1}y^{b_1})^{\lambda})+x_1^dy^e(z_1+\alpha)^{f-\frac{dc}{a}}
\end{array}
$$
First suppose that $ae-bd\ne 0$ and $\text{ord }(P)<\frac{e}{b}$. Then $a_1e-b_1d>0$
and $\lambda\text{ord }(P)<\frac{e}{b_1}$. By assumption $\lambda k>\frac{d}{a_1}$ and
$\lambda\text{ord }(P)>\frac{d}{a_1}$. By (\ref{eq611}), $(u,v)$ is not invertible at 2 points 
$q'$ on $x=y=0$. 

Now suppose that $ae-bd\ne 0$ and $\text{ord }(P)\ge\frac{e}{b}$. Then $a_1e-b_1d>0$
and $\lambda\text{ord }(P)\ge\frac{e}{b_1}$, so that we can choose permissible coordinates at $q'$ so that
$$
\begin{array}{ll}
u&=(x_1^{a_1}y_1^{b_1})^{k\lambda}\\
v&=x_1^dy_1^e
\end{array}
$$
By assumption $\lambda k>\frac{d}{a_1}$, so that by  (\ref{eq613}), $(u,v)$ is not invertible at 2 points 
$q'$ on $x=y=0$ if and only if $k<\frac{e}{b}$. 

Suppose that $ae-bd= 0$ and $\text{ord }(P)<\frac{e}{b}$. Then $a_1e-b_1d=0$
and $\lambda\text{ord }(P)< \frac{e}{b_1}$. Since (\ref{eq110}) can then not hold at $q'$,
we have that $(u,v)$ are invertible at 2 points $q'$ on $x=y=0$. 

Now suppose that $ae-bd= 0$ and $\text{ord }(P)\ge\frac{e}{b}$. Then $a_1e-b_1d=0$
and $\lambda\text{ord }(P)\ge\frac{e}{b_1}$, so that we can choose permissible coordinates at $q'$ so that
$$
\begin{array}{ll}
u&=(x_1^{a_1}y_1^{b_1})^{k\lambda}\\
v&=(\beta+\alpha^{f-\frac{dc}{a}})(x_1^{a_1}y_1^{b_1})^t+(x_1^{a_1}y_1^{b_1})^{t}z_1
\end{array}
$$
where $t=\frac{e}{b}\lambda$, $\beta\in k$ is the degree $t$ coefficient of $P$.
For $q'$ in a possibly smaller neighborhood of $p_1$, (\ref{eq110}) can then not hold at
$q'$, so that $(u,v)$ are invertible at 2 points $q'$ on $x=y=0$.

Suppose that $q'$ is a 1 point in $U$ on $z=0$. 
$q'$ has regular parameters $(x_1,y_1,z)$ where $x=x_1+\alpha$, $y=y_1+\beta$
with $\alpha,\beta\ne 0$.
Thus $q'$ has permissible parameters $(x_1,y_1,z_1)$ where  $z_1$ is defined by
$$
z=(x_1+\alpha)^{-\frac{a}{c}}(y_1+\beta)^{-\frac{b}{c}}z_1
$$
$$
\begin{array}{ll}
u&=z_1^{ck}\\
v&=P(z_1^c)+(x_1+\alpha)^{d-\frac{af}{c}}(y_1+\beta)^{e-\frac{bf}{c}}z_1^f.
\end{array}
$$
Since by assumption $c\text{ ord}(P)<f$, $q'$ cannot be in the form of (\ref{eq612}), so that 
$(u,v)$ is invertible at 1 points on $z=0$.

Suppose that $q'$ is a 1 point in $U$ on $y=0$. 
$q'$ has regular parameters $(x_1,y,z_1)$ where $x=x_1+\alpha$, $z=z_1+\beta$
with $\alpha,\beta\ne 0$.
Thus $q'$ has permissible parameters $(x_1,y_1,z_1)$ where  $y_1$ is defined by
$$
y=(x_1+\alpha)^{-\frac{a}{b}}(z_1+\beta)^{-\frac{c}{b}}y_1
$$
$$
\begin{array}{ll}
u&=y_1^{bk}\\
v&=P(y_1^b)+(x_1+\alpha)^{d-\frac{ae}{b}}(z_1+\beta)^{f-\frac{ec}{b}}y_1^e
\end{array}
$$

If $b\text{ ord }(P)<e$ or $b\text{ ord }(P)>e$
 then $q'$ cannot have the form of (\ref{eq612}), so that $(u,v)$ is invertible at all
1 points on $y=0$. If $k\le \frac{e}{b}$, then $(u,v)$ is invertible at all 1 points on
$y=0$.

Suppose that $b\text{ ord }(P)=e$ and $k>\frac{e}{b}$. Then we can write $P(t) = \gamma t^{\frac{e}{b}}+\cdots$
where $\gamma\ne0$. We have $(u,v)$ is invertible at $q'$ on $y=0$ unless
$$
\alpha^{d-\frac{ae}{b}}\beta^{f-\frac{ec}{b}}+\gamma=0
$$
which holds only if $(\alpha,\beta)$ are on the algebraic curve
$$
\beta^{bf-ce}=(-\gamma)^b\alpha^{ae-bd}.
$$
In this case $(u,v)$ is not invertible on the curve with local equations
$$
y=0, z^{bf-ec}+(-\gamma)^bx^{ae-bd}=0.
$$

If  $q'$ is a 1 point in $U$ on $x=0$, then 
there are permissible parameters $(x_1,y_1,z_1)$ at $q'$ such that
$$
\begin{array}{ll}
u&=x_1^{ak}\\
v&=P(x_1^a)+x_1^d(y_1+\alpha)^{e-\frac{db}{a}}(z_1+\beta)^{f-\frac{dc}{a}}
\end{array}
$$
with $\alpha,\beta\ne 0$. Thus
$(u,v)$ is
invertible at 1 points on $x=0$ in $U$ since $a\text{ ord }(P)>d$.

If $\pi:X_1\rightarrow X$ is the blowup of a 2 curve through $p$, then $\Phi_1=\Phi\circ\pi$ is strongly prepared above $p$. $\pi^{-1}(p)$ is a 2 curve, so there are no 1 points in 
$\pi^{-1}(p)$.

{\bf Suppose that $p\in X$ is a 3 point such that $m_q{\cal O}_{X,p}$ is not invertible, 
and (\ref{eq112}) holds at $p$.}
We may assume that there exists an open neighborhood $U$ of $p$ such that $(x,y,z)$ are 
uniformizing parameters on an \'etale cover of $U$.

After possibly interchanging $x,y,z$, we can assume that
$$
\frac{d}{a}=\text{min}\{\frac{d}{a},\frac{e}{b},\frac{f}{c}\}
$$
and 
$$
\frac{f}{c}=\text{max}\{\frac{d}{a},\frac{e}{b},\frac{f}{c}\}
$$

We will determine the locus of points in $U$ where $m_q{\cal O}_U$ is not invertible. First suppose that $q'$ is a 2 point
on the curve $x=z=0$. $q'$ has regular parameters $(x,y_1,z)$ where $y=y_1+\alpha$.
Thus $q'$ has permissible parameters $(x_1,y_1,z_1)$ where  $x_1$ is defined by
$$
x=x_1(y_1+\alpha)^{-\frac{b}{a}}
$$
Set $\lambda = (a,c)$, $a_1=\frac{a}{\lambda}$, $c_1=\frac{c}{\lambda}$.
$$
\begin{array}{ll}
u&=(x_1^{a_1}z^{c_1})^{k\lambda}\\
v&=x_1^dz_1^f(y_1+\alpha)^{e-\frac{db}{a}}
\end{array}
$$
We have $a_1f-c_1d>0$.
$\frac{d}{a}<k<\frac{f}{c}$ implies $\frac{d}{a_1}<k\lambda<\frac{f}{c_1}$.
 $q'$ thus has the form of (\ref{eq613}),
and we see that $(u,v)$ is not invertible on the curve with local equations $x=z=0$.

Suppose that $q'$ is a 2 point on the curve $y=z=0$. $q'$ has permissible parameters
$(x_1,y_1,z)$ where $x=x_1+\alpha$, $y_1$ is defined by $y=y_1(x_1+\alpha)^{-\frac{a}{b}}$.
Set $\lambda=(b,c)$, $b_1=\frac{b}{\lambda}$, $c_1=\frac{c}{\lambda}$.
$$
\begin{array}{ll}
u&=(y_1^{b_1}z^{c_1})^{k\lambda}\\
v&=(x_1+\alpha)^{d-\frac{ea}{b}}y_1^ez^f.
\end{array}
$$
First suppose that $bf-ce\ne 0$. Then $bf-ce>0$ and $b_1f-c_1e>0$. Since $k\lambda<\frac{f}{c_1}$,
we have by (\ref{eq613}) that $(u,v)$ is not invertible at 2 points $q'$ on $y=z=0$
if and only if $\frac{e}{b_1}<k\lambda$, which holds if and only if $\frac{e}{b}<k$.

Now suppose that $bf-ce=0$, so that $b_1f-c_1e=0$. Then $(u,v)$ is invertible at
2 points $q'$ on $y=z=0$.

Now suppose that $q'$ is a 2 point on the curve $x=y=0$. $q'$ has regular parameters
$(x,y,z_1)$ where $z=z_1+\alpha$. $q'$ has permissible parameters $(x_1,y,z_1)$ where
$x_1$ is defined by $x=x_1(z_1+\alpha)^{-\frac{c}{a}}$.
Set $\lambda=(a,b)$, $a_1=\frac{a}{\lambda}$, $b_1=\frac{b}{\lambda}$.
$$
\begin{array}{ll}
u&=(x_1^{a_1}y^{b_1})^{k\lambda}\\
v&=x_1^dy^e(z_1+\alpha)^{f-\frac{dc}{a}}
\end{array}
$$
First suppose that $ae-bd\ne 0$. Then $a_1e-b_1d>0$. By assumption
$\frac{d}{a_1}<k\lambda$. We have by (\ref{eq613}) that $(u,v)$ is not invertible
at 2 points $q'$ on $x=y=0$ if and only if $k\lambda<\frac{e}{b_1}$ which holds
if and only if $k<\frac{e}{b}$.

Now suppose that $ae-bd=0$. Then $a_1e-b_1d=0$ and $(u,v)$ is invertible at
 2 points on the curve $x=y=0$.

Suppose that $q'$ is a 1 point in $U$ on $z=0$. $q'$ has regular parameters
$(x_1,y_1,z)$ where $x=x_1+\alpha$, $y=y_1+\beta$ (with $\alpha,\beta\ne0$). Thus
$q'$ has permissible parameters $(x_1,y_1,z_1)$ where $z_1$ is defined by
$z=(x_1+\alpha)^{-\frac{a}{c}}(y_1+\beta)^{-\frac{b}{c}}z_1$.
$$
\begin{array}{ll}
u&=z_1^{ck}\\
v&=(x_1+\alpha)^{d-\frac{af}{c}}(y_1+\beta)^{e-\frac{bf}{c}}z_1^f
\end{array}
$$
Thus $(u,v)$ is invertible at all 1 points of $z=0$.

Similarily, $(u,v)$ is invertible at all 1 points of $x=0$ and $y=0$.

If $\pi:X_1\rightarrow X$ is the blow up of a  2 curve $C$ through  $p$, then
$\Phi_{X_1}$ is strongly prepared above $p$. 
$\pi^{-1}(p)$ is a 2 curve, so there are
no 1 points in $\pi^{-1}(p)$.

{\bf Suppose that $p\in X$ is a 2 point such that $m_q{\cal O}_{X,p}$ is not invertible, 
and (\ref{eq108}) holds at $p$.} We may assume that
there exists an open neighborhood $U$ of $p$ such that $(x,y,z)$ are uniformizing parameters on
an \'etale cover of $U$, and the conclusions of Lemma \ref{Lemma1067} hold for $p$.
After possibly interchanging $x$ and $y$ we may assume that $ad-bc>0$.
We will determine the locus of points in $U$ where $m_q{\cal O}_U$ is not invertible. First suppose that $q'$ is a 2 point
on the curve $x=y=0$. $q'$ has regular  parameters $(x,y,z_1)$ where $z=z_1+\alpha$. Thus $q'$ has permissible
parameters $(\overline x,\overline y,\overline z)$ such that
$$
\begin{array}{ll}
u&=(\overline x^a\overline y^b)^k\\
v&=P(\overline x^a\overline y^b)+\overline x^c\overline y^d
\end{array}
$$
Since $\frac{c}{a}<e=\text{ord }(P)<\frac{d}{b}$, and $k>\frac{c}{a}$, we are in the form of (\ref{eq611}).
Thus $(u,v)$ is not invertible along the curve $x=y=0$.

Suppose that $q'$ is a 1 point
near $q$. $q'$ has permissible  parameters $(x_1,y_1,z_1)$ where 
either 
\begin{equation}\label{eq113}
x=x_1(y_1+\alpha)^{-\frac{b}{a}},
y=y_1+\alpha,
z=z_1+\beta
\end{equation}
with $\alpha\ne 0$ or 
\begin{equation}\label{eq114}
x=x_1+\alpha,
y=y_1(x_1+\alpha)^{-\frac{a}{b}},
z=z_1+\beta
\end{equation}
with $\alpha\ne 0$. If $q'$ has permissible parameters satisfying (\ref{eq113}), then 
since $\text{ord }(P)>\frac{c}{a}$, 
\begin{equation}\label{eq1056}
\begin{array}{ll}
u&=x_1^{ak}\\
v&=P(x_1^a)+x_1^c(y_1+\alpha)^{d-\frac{bc}{a}}(z_1+\beta)\\
&=\beta\alpha^{d-\frac{bc}{a}}x_1^c+x_1^c\overline z
\end{array}
\end{equation}
$(u,v)$ is not invertible at $q'$ if and only if $q'$ satisfies (\ref{eq612}). Since $c<ak$ 
by assumption,
this holds if and only if $\beta=0$.

If $q'$ is a 1 point near $p$ on $x=0$ (so that (\ref{eq1056}) holds) then $A(\Phi,q')=0$.

If $q'$ has permissible parameters satisfying (\ref{eq114}), then 
\begin{equation}\label{eq1057}
\begin{array}{ll}
u&=y_1^{bk}\\
v&=P(y_1^b)+(x_1+\alpha)^{c-\frac{da}{b}}y_1^d(z_1+\beta)\\
\end{array}
\end{equation}
$(u,v)$ is  invertible at $q'$,  since $b\,\text{ord P}<d$ 
by assumption, so that  $q'$ cannot satisfy (\ref{eq612}).

If $q'$ is a 1 point near $p$ on $y=0$ (so that (\ref{eq1057}) holds) then 
$$
A(\Phi,q')=d-b\text{ ord }(P).
$$

We will now consider the  invariant $A$ on the blowup of $V(x,y)$ or $V(x,z)$ over $p$.

Let $\pi_1:X_1\rightarrow X$ be the blowup of $C=V(x,y)$.
$\Phi_{1}=\Phi\circ\pi_1$ is strongly prepared above $p$. If $q\in\pi_1^{-1}(p)$ is a 1 point, then
$q$ has regular parameters
$(x,y_1,z)$ defined by
$$
x=x_1,
y=x_1(y_1+\alpha)
$$
with $\alpha\ne 0$. There are permissible parameters $(\overline x_1,y_,z)$ at $q$
 where $\overline x_1$ is defined by
$$
x_1 = \overline x_1(y_1+\alpha)^{-\frac{b}{a+b}}
$$
Thus
$$
\begin{array}{ll}
u&=\overline x_1^{(a+b)k}\\
v&=P(\overline x_1^{a+b})+\overline x_1^{c+d}(y_1+\alpha)^{d-\frac{(c+d)b}{a+b}}z
\end{array}
$$
If $(a+b)\text{ord}(P)\ge c+d$, then $A(\Phi_1,q)=0$.
 Assume that $(a+b)\text{ord p}<c+d$.
 Since
$\text{ord }(P)>\frac{c}{a}$, we have that
$$
c+d-(a+b)\text{ ord}(P)=(d-b\text{ ord}(P))+(c-a\text{ ord}(P))<d-b\text{ ord}(P)
$$
Thus
$$
A(\Phi_1,q)\le c+d-(a+b)\text{ ord}(P)<d-b\text{ ord}(P)\le A(\Phi).
$$

If $\pi_1:X_1\rightarrow X$ is the blowup of $C=V(x,z)$, then 
$\Phi_{X_1}$ is strongly prepared above $p$, and there are no 1 points in  $\pi_1^{-1}(p)$.

{\bf Suppose that $p\in X$ is a 2 point such that $m_q{\cal O}_{X,p}$ is not invertible, 
and (\ref{eq109}) holds at $p$.}
We may assume that
there exists an open neighborhood $U$ of $p$ such that $(x,y,z)$ are uniformizing parameters on
an \'etale cover of $U$. After possibly interchanging $x$ and $y$, we may assume that $ad-bc>0$.
We will determine the locus of points in $U$ where $m_q{\cal O}_U$ is not invertible. First suppose that $q'$ is a 2 point
on the curve $x=y=0$. $q'$ has regular  parameters $(x,y,z_1)$ where $z=z_1+\alpha$. Thus $q'$ has permissible
parameters $(\overline x,\overline y,\overline z)$ such that
$$
\begin{array}{ll}
u&=(\overline x^a\overline y^b)^k\\
v&=\overline x^c\overline y^d
\end{array}
$$
Since $\frac{c}{a}<k$, we are in the form of (\ref{eq613}), and
 $(u,v)$ is not invertible along the curve $x=y=0$ if and only if $k<\frac{d}{b}$.

Suppose that $q'$ is a 1 point
near $p$. $q'$ has permissible  parameters $(x_1,y_1,z_1)$ where 
either 
\begin{equation}\label{eq115}
x=x_1(y_1+\alpha)^{-\frac{b}{a}},
y=y_1+\alpha,
z=z_1+\beta
\end{equation}
with $\alpha\ne0$ or 
\begin{equation}\label{eq116}
x=x_1+\alpha,
y=y_1(x_1+\alpha)^{-\frac{a}{b}},
z=z_1+\beta
\end{equation}
with $\alpha\ne 0$. If $q'$ has permissible parameters satisfying (\ref{eq115}), then 
$$
\begin{array}{ll}
u&=x_1^{ak}\\
v&=x_1^c(y_1+\alpha)^{d-\frac{bc}{a}}(z_1+\beta)\\
&=\beta\alpha^{d-\frac{bc}{a}}x_1^c+x_1^c\overline z
\end{array}
$$
$(u,v)$ is not invertible at $q'$ if and only if $q'$ satisfies (\ref{eq612}). Since $c<ak$ 
by assumption,
this holds if and only if $\beta=0$.

If $q'$ is a 1 point near $p$ on $x=0$ (so that (\ref{eq115}) holds) then $A(\Phi,q')=0$.

If $q'$ has permissible parameters satisfying (\ref{eq116}), then 
$$
\begin{array}{ll}
u&=y_1^{bk}\\
v&=(x_1+\alpha)^{c-\frac{da}{b}}y_1^d(z_1+\beta)\\
\end{array}
$$
$$
\begin{array}{ll}
u&=y_1^{bk}\\
v&=\alpha^{c-\frac{da}{b}}\beta y_1^d+y_1^d\overline z_1
\end{array}
$$
Thus $(u,v)$ is invertible at $q'$ if $\beta\ne 0$, and if $\beta=0$, then
$(u,v)$ is invertible at $q'$ if and only if $d\ge kb$.

If $q'$ is a 1 point near $p$ on $y=0$ (so that (\ref{eq116}) holds) then
$A(\Phi,q')=0$.

We will now consider the  invariant $A$ on the blowup of a curve $V(x,y)$, $V(y,z)$ or $V(x,z)$ 
where $(u,v)$ is not invertible on the curve.

Let $\pi_1:X_1\rightarrow X$ be the blowup of $C=V(x,y)$. 
$\Phi_1=\Phi\circ\pi_1$ is strongly prepared over $p$. If $q\in\pi^{-1}(p)$ is a 1 point, then
$q$ has regular parameters
$(x,y_1,z)$ defined by
$$
x=x_1,
y=x_1(y_1+\alpha)
$$
with $\alpha\ne 0$. There are permissible parameters $(x_1,y_,z)$ at $q$
 where $x_1$ is defined by
$$
x_1 = \overline x_1(y_1+\alpha)^{-\frac{b}{a+b}}
$$
Thus
$$
\begin{array}{ll}
u&=\overline x_1^{(a+b)k}\\
v&=\overline x_1^{c+d}(y_1+\alpha)^{d-\frac{(c+d)b}{a+b}}z
\end{array}
$$
and $A(\Phi_1,q)=0$.

If $\pi_1:X_1\rightarrow X$ is the blowup of $C=V(x,z)$ or $V(y,z)$,
 then $\Phi_1=\Phi\circ\pi_1$ is strongly prepared over $p$, and there are no 1 points in  $\pi^{-1}(p)$.

{\bf Suppose that $p\in X$ is a 2 point such that $m_q{\cal O}_{X,p}$ is not invertible, 
and (\ref{eq110}) holds at $p$.} We may assume that 
there exists an open neighborhood $U$ of $p$ such that $(x,y,z)$ are uniformizing parameters on
an \'etale cover of $U$.
We will determine the locus of points in $U$ where $m_q{\cal O}_U$ is not invertible. First suppose that $q'$ is a 2 point
on the curve $x=y=0$. $q'$ has permissible  parameters $(x,y,z_1)$ where $z=z_1+\alpha$. 
$$
\begin{array}{ll}
u&=(x^ay^b)^k\\
v&=\alpha (x^ay^b)^t+(x^ay^b)^tz_1
\end{array}
$$
Thus $(u,v)$ is invertible along the curve $x=y=0$, if $\alpha\ne 0$.

Suppose that $q'$ is a 1 point
near $p$. $q'$ has permissible  parameters $(x_1,y_1,z_1)$ where 
either 
\begin{equation}\label{eq117}
x=x_1(y_1+\alpha)^{-\frac{b}{a}},
y=y_1+\alpha,
z=z_1+\beta
\end{equation}
with $\alpha\ne 0$ or 
\begin{equation}\label{eq118}
x=x_1+\alpha,
y=y_1(x_1+\alpha)^{-\frac{a}{b}},
z=z_1+\beta
\end{equation}
with $\alpha\ne 0$. If $q'$ has permissible parameters satisfying (\ref{eq117}), then 
$$
\begin{array}{ll}
u&=x_1^{ak}\\
v&=x_1^{at}(z_1+\beta)\\
\end{array}
$$
Thus $(u,v)$ is invertible at $q'$ if $\beta \ne 0$, and $q'$ is not invertible along $V(x,z)$
since $t<k$ 
by assumption.

If $q'$ has permissible parameters satisfying (\ref{eq118}), then 
$$
\begin{array}{ll}
u&=y_1^{bk}\\
v&=y_1^{bt}(z_1+\beta)\\
\end{array}
$$
Thus $(u,v)$ is invertible at $q'$ if $\beta \ne 0$, and $q'$ is not invertible along $V(y,z)$
since $t<k$ 
by assumption.

If $\pi_1:X_1\rightarrow X$ is the blowup of $V(x,z)$ or $V(y,z)$, then 
$\Phi_1=\Phi\circ\pi_1$ is strongly prepared over $p$ and there are no 1 points in  $\pi^{-1}(p)$.

{\bf Suppose that $p\in X$ is a 2 point such that $m_q{\cal O}_{X,p}$ is not invertible, 
and (\ref{eq611}) holds at $p$.} After possibly interchanging $x$ and $y$, we may assume
that $ad-bc>0$. We may assume that
there exists an open neighborhood $U$ of $p$ such that $(x,y,z)$ are uniformizing parameters on
an \'etale cover of $U$ and the conclusions of Lemma \ref{Lemma1067} hold for $p$.
We will determine the locus of points in $U$ where $m_q{\cal O}_U$ is not invertible. If $q'$ is a 2 point
on the curve $x=y=0$, then $q'$ has the form of (\ref{eq611}), so that 
Thus $(u,v)$ is not invertible along the curve $x=y=0$.

Suppose that $q'$ is a 1 point
near $q'$. $q'$ has permissible  parameters $(x_1,y_1,z_1)$ where 
either 
\begin{equation}\label{eq119}
x=x_1(y_1+\alpha)^{-\frac{b}{a}},
y=y_1+\alpha,
z=z_1+\beta
\end{equation}
with $\alpha\ne 0$ or 
\begin{equation}\label{eq120}
x=x_1+\alpha,
y=y_1(x_1+\alpha)^{-\frac{a}{b}},
z=z_1+\beta
\end{equation}
with $\alpha\ne 0$. If $q'$ has permissible parameters satisfying (\ref{eq119}), then 
$$
\begin{array}{ll}
u&=x_1^{ak}\\
v&=P(x_1^a)+x_1^c(y_1+\alpha)^{d-\frac{bc}{a}}\\
&=\alpha^{d-\frac{bc}{a}}x_1^c+x_1^c\overline y
\end{array}
$$
for some permissible parameters $(x_1,\overline y,z)$,
since $a\text{ ord }(P)>c$. Thus 
$(u,v)$ is invertible at $q'$, since we have $\alpha\ne 0$.

$A(\Phi,q')=0$ at points $q'$ near $p$ where (\ref{eq119}) holds.

If $q'$ has permissible parameters satisfying (\ref{eq120}), then 
$$
\begin{array}{ll}
u&=y_1^{bk}\\
v&=P(y_1^b)+\alpha^{c-\frac{da}{b}}y_1^d+y_1^d\overline x_1\\
\end{array}
$$
$(u,v)$ is  invertible at $q'$  Since $b\text{ ord}(P)<d$ 
by assumption.

At points $q'$ near $p$ where (\ref{eq120}) holds,  we have
$A(\Phi,q')=d-b\text{ ord }(P)>0$. 

Let $\pi_1:X_1\rightarrow X$ be the blowup of $C=V(x,y)$. 
Then $\Phi_1=\Phi\circ\pi_1$ is strongly prepared above $p$. If $q\in\pi_1^{-1}(p)$ is a 1 point, then
$q$ has regular parameters
$(x,y_1,z)$ defined by
$$
x=x_1,
y=x_1(y_1+\alpha)
$$
with $\alpha\ne 0$. There are permissible parameters $(\overline x_1,y_,z)$ at $q$
 where $\overline x_1$ is defined by
$$
x_1 = \overline x_1(y_1+\alpha)^{-\frac{b}{a+b}}
$$
Thus
$$
\begin{array}{ll}
u&=\overline x_1^{(a+b)k}\\
v&=P(\overline x_1^{a+b})+\overline x_1^{c+d}(y_1+\alpha)^{d-\frac{(c+d)b}{a+b}}\\
&=P(\overline x_1^{a+b})+\overline x_1^{c+d}\alpha^{d-\frac{(c+d)b}{a+b}}
+\overline x_1^{c+d}\overline y_1
\end{array}
$$
If $(a+b)\text{ord }(P)\ge c+d$ then $A(\Phi_1,q)=0$. Assume that $(a+b)\text{ord}(P)<c+d$. Since
$\text{ord}(P)>\frac{c}{a}$, we have that
$$
c+d-(a+b)\text{ord}(P)=(d-b\text{ ord}(P))+(c-a\text{ ord}(P))<d-b\text{ ord}(P)
$$
Thus
$$
A(\Phi_1,q)\le 
c+d-(a+b)\text{ ord }(P)
< d-b\text{ ord }(P)\le A(\Phi)
$$

{\bf Suppose that $p\in X$ is a 2 point such that $m_q{\cal O}_{X,p}$ is not invertible, 
and (\ref{eq613}) holds at $p$.}
After possibly interchanging $x$ and $y$, we may assume that $ad-bc>0$.
We may assume that there exists an open neighborhood $U$ of $p$ such that $(x,y,z)$ are uniformizing parameters on an \'etale cover of $U$.
We will determine the locus of points in $U$ where $m_q{\cal O}_U$ is not invertible. 
If $q'$ is a 2 point
on the curve $x=y=0$, then $q'$ has the form of (\ref{eq613}), so that 
$(u,v)$ is not invertible along the curve $x=y=0$.

Suppose that $q'$ is a 1 point
near $p$. $q'$ has permissible  parameters $(x_1,y_1,z_1)$ where 
either 
\begin{equation}\label{eq121}
x=x_1(y_1+\alpha)^{-\frac{b}{a}},
y=y_1+\alpha,
z=z_1+\beta
\end{equation}
with $\alpha\ne0$ or 
\begin{equation}\label{eq122}
x=x_1+\alpha,
y=y_1(x_1+\alpha)^{-\frac{a}{b}},
z=z_1+\beta
\end{equation}
with $\alpha\ne 0$. If $q'$ has permissible parameters satisfying (\ref{eq121}), then 
$$
\begin{array}{ll}
u&=x_1^{ak}\\
v&=x_1^c(y_1+\alpha)^{d-\frac{bc}{a}}\\
&=\alpha^{d-\frac{bc}{a}}x_1^c+x_1^c\overline y_1
\end{array}
$$
$(u,v)$ is thus invertible at $q'$. 

$A(\Phi,q')=0$ at points $q'$ near $p$ where (\ref{eq121}) holds. 

If $q'$ has permissible parameters satisfying (\ref{eq122}), then 
$$
\begin{array}{ll}
u&=y_1^{bk}\\
v&=(x_1+\alpha)^{c-\frac{da}{b}}y_1^d\\
&=\alpha^{c-\frac{da}{b}}y_1^d + y_1^d\overline x_1
\end{array}
$$
Thus $(u,v)$ is  invertible at $q'$.

$A(\Phi,q')=0$ at points $q'$ near $p$ where (\ref{eq122}) holds. 

The locus of points where $(u,v)$ is not invertible near $p$ is
$V(x,y)$.

Let $\pi_1:X_1\rightarrow X$ be the blowup of $C=V(x,y)$. 
$\Phi_1=\Phi\circ\pi_1$ is strongly prepared above $p$. If $q\in\pi_1^{-1}(p)$ is a 1 point, then
$q$ has regular parameters
$(x_1,y_1,z)$ defined by
$$
x=x_1,
y=x_1(y_1+\alpha)
$$
with $\alpha\ne 0$. There are permissible parameters $(\overline x_1,y_,z)$
 at $q$ where $\overline x_1$ is defined by
$$
x_1 = \overline x_1(y_1+\alpha)^{-\frac{b}{a+b}}.
$$
Thus
$$
\begin{array}{ll}
u&=\overline x_1^{(a+b)k}\\
v&=\overline x_1^{c+d}(y_1+\alpha)^{d-\frac{(c+d)b}{a+b}}\\
&=\overline x_1^{c+d}\alpha^{d-\frac{(c+d)b}{a+b}}+\overline x_1^{c+d}\overline y_1
\end{array}
$$
and $A(\Phi_1,q)=0$.

{\bf Suppose that $p\in X$ is a 1 point such that $m_q{\cal O}_{X,p}$ is not invertible, 
so that (\ref{eq612}) holds at $p$.} We may assume that 
here exists an open neighborhood $U$ of $p$ such that $(x,y,z)$ are uniformizing parameters on
an \'etale cover of $U$.
We will determine the locus of points in $U$ where $m_q{\cal O}_U$ is not invertible.

Suppose that $q'$ is a 1 point
near $p$. $q'$ has permissible  parameters $(x,y_1,z_1)$ where 
$$
y=y_1+\alpha,
z=z_1+\beta
$$
$$
\begin{array}{ll}
u&=x^k\\
v&=\alpha x^c+x^cy_1\\
\end{array}
$$
$(u,v)$ is thus only not invertible on the curve $V(x,y)$.

Let $\pi:X_1\rightarrow X$ be the blowup of $V(x,y)$.
$\Phi_1=\Phi\circ\pi_1$ is strongly prepared above $p$. 
If $q\in\pi^{-1}(p)$ is a 1 point, then
$q$ has permissible parameters
$(x,y_1,z)$ defined by
$$
x=x_1,
y=x_1(y_1+\alpha)
$$
with $\alpha\ne 0$.
$$
\begin{array}{ll}
u&=x_1^k\\
v&=x_1^{c+1}(y_1+\alpha)
\end{array}
$$
and $A(\Phi_1,q)=0$.

{\bf Suppose that $p\in X$ is a 3 point such that $m_p{\cal O}_{X,p}$ is not
invertible and (\ref{eq1055}) holds at $p$.}
We may assume that there exists an open neighborhood $U$ of $p$ such that $(x,y,z)$
are uniformizing parameters on an \'etale cover of $U$. The locus of points  in $U$ where
$m_q{\cal O}_U$ is not invertible is the union of the 2 curves $V(x,z)$, $V(y,z)$ and $V(x,y)$.

Let $\pi:X_1\rightarrow X$ be the blowup of $C=V(x,y)$, $\Phi_1=\Phi\circ\pi$.
If $q'\in \pi^{-1}(p)$ is a 2 point, then $q'$ has permissible parameters $(x_1,y_1,z)$ where
$$
x=x_1, y=x_1(y_1+\alpha)
$$
with $\alpha\ne 0$.
$$
\begin{array}{ll}
u&=x_1^{a+b}(y_1+\alpha)^b=\overline x_1^{a+b}\\
v&=z^c
\end{array}
$$
so that $\Phi_1$ is strongly prepared at $q'$. 

Suppose that $q'\in\pi^{-1}(p)$ is a 3 point and $q'$ has permissible parameters
$(x_1,y_1,z)$ where $x=x_1, y=x_1y_1$. Then
$$
\begin{array}{ll}
v&=x_1^{a+b}y_1^b\\
u&=z^c
\end{array}
$$
so that $\Phi_1$ is strongly prepared at $q'$. 
Suppose that $q'\in\pi^{-1}(p)$ is a 3 point and $q'$ has permissible parameters
$(x_1,y_1,z)$ where 
$$
x=x_1y_1, y=y_1.
$$
 Then
$$
\begin{array}{ll}
u&=y_1^{a+b}\\
v&=z^c
\end{array}
$$
so that $\Phi_1$ is strongly prepared at $q'$.
A similar analysis shows that the blowup of
$V(x,z)$ or $V(y,z)$ composed with $\Phi$ is strongly prepared.

{\bf Suppose that $p\in X$ is a 3 point such that $m_q{\cal O}_{X,p}$ is not
invertible and (\ref{eq1059}) holds at $p$.}

We may assume that there exists an open neighborhood $U$ of $p$ such that $(x,y,z)$
are uniformizing on an \'etale cover of $U$. The locus of points in $U$ where
$m_q{\cal O}_U$ is not invertible is the  union of the 2 curves
$V(x,z)$, $V(x,y,)$ (if $c>b)$ and $V(y,z)$ (if $b>c)$.  If $\pi:X_1\rightarrow X$ is the blowup of a
2 curve through $C$, and $\Phi_1=\Phi\circ\pi$, then $\Phi_1$ is strongly prepared at
points $q\in\pi^{-1}(p)$.

{\bf Suppose that $p\in X$ is a 2 point such that $m_p{\cal O}_{X,p}$ is not
invertible, and (\ref{eq1054}) holds at $p$.}
We may assume that there exists an open neighborhood $U$ of $p$ such that
$(x,y,z)$ are uniformizing parameters on an \'etale cover of $U$. The locus of points
in $U$ where $m_q{\cal O}_U$ is not invertible is the 2 curve $V(x,y)$.

Let $\pi:X_1\rightarrow X$ be the blowup of $C=V(x,y)$, $\Phi_1=\Phi\circ\pi$.
If $q'\in\pi^{-1}(p)$ is a 1 point, then $q'$ has permissible parameters
$$
x=x_1, y=x_1(y_1+\alpha)
$$
with $\alpha\ne 0$.
$$
\begin{array}{ll}
u&=x_1^a\\
v&=x_1^b(y_1+\alpha)^b.
\end{array}
$$
Set $\overline y_1=(y_1+\alpha)^b-\alpha^b$ to get
$$
\begin{array}{ll}
u&=x_1^a\\
v&=\alpha^bx_1^b+x_1^b\overline y_1
\end{array}
$$
so that $A(\Phi_1,q')=0$. $\Phi_1$ is strongly prepared at points of $\pi^{-1}(p)$.
\end{pf}

If $\alpha, \beta$ are real numbers, define
$$
S(\alpha,\beta)=\text{ max }\{(\alpha,\beta),(\beta,\alpha)\}
$$
where the maximum is in the Lexicographic ordering.

Suppose that $\Phi:X\rightarrow S$ is strongly prepared. Suppose that $q\in D_S$ and $C\subset X$ is a 2 curve such that
$m_q{\cal O}_{X}$ is not invertible along $C$. At a generic point $p$ of $C$ (\ref{eq611}),
 (\ref{eq613}) or (\ref{eq1054}) holds.

If (\ref{eq611})  holds, then $\Phi(C)$ is a 1 point $q\in S$. Suppose that
$P(t)=\sum a_it^i$. Since $q$ is a 1 point, we can, after possibly replacing $v$ with
$v-\sum a_{ik}u^k$, assume that $k\not\,\mid\text{ord}(P)$ in (\ref{eq611}). With this
restriction, define
$$
\sigma(C)=\left\{
\begin{array}{l}
S(\mid c-a\text{ ord }(P)\mid,\mid d-b\text{ ord }(P)\mid)\\
\,\,\,\text{if }c-\text{ ord }(P), d-\text{ ord }(P)\text{ have opposite signs,}\\
\,\,\,-\infty\text{ if they have the same sign.}
\end{array}
\right.
$$
If (\ref{eq613}) or (\ref{eq1054}) holds, define
$$
\sigma(C)=-\infty.
$$ 
$\sigma(C)$ is well defined (independent of choice of permissible parameters $(u,v)$ at $q$
with the restriction that  $k\not\,\mid\text{ord}(P)$ in (\ref{eq611})). This follows from  
Lemma \ref{Lemma1067}.

If $m_q{\cal O}_X$ is invertible, define 
$$
\overline\sigma(\Phi)=-\infty.
$$
If $m_q{\cal O}_X$ is not invertible, define 
$$
\overline \sigma(\Phi) = 
\text{ max }\left\{\begin{array}{ll}\sigma(C)&\mid C\subset X
\text{ is a 2 curve} \\
& \,\,\,\text{such that }m_q{\cal O}_{X}\text{ is not invertible along }C 
\end{array}\right\}.
$$

\begin{Lemma}\label{Lemma1043} Suppose that $X$ is strongly prepared, $q\in S$ is
such that $m_q{\cal O}_X$ is not invertible. Then there exists a sequence of blowups of
2 curves $X_1\rightarrow X$ such that the induced map $\Phi_1:X_1\rightarrow S$ is strongly prepared,
 $A(\Phi_1,E)< A(\Phi_1)=A(\Phi)$ if $E$ is an exceptional component of $E_{X_1}$ for 
$X_1\rightarrow X$, and the forms (\ref{eq611}), (\ref{eq108}) and (\ref{eq111})
do not hold at  any point $p\in X$ where $m_q{\cal O}_{X_1,p}$ is not invertible.
\end{Lemma}

\begin{pf}$\overline \sigma(\Phi)\ge 0$ if and only if
there exists a  point $p\in X$ such that $m_q{\cal O}_{X_1,p}$ is not invertible,
and
 a form (\ref{eq611}), (\ref{eq108}) or (\ref{eq111})
holds at $p$.

Suppose that $\overline \sigma(\Phi)\ge 0$. Let $C$ be a 2 curve such that $\sigma(C)
=\overline \sigma(\Phi)$. Let $\pi:X_1\rightarrow X$ be the blowup of $C$.  By Theorem \ref{Theorem63}, we need only verify that if
$C_1\subset \pi^{-1}(C)$ is a 2 curve such that $m_q{\cal O}_{X_1}$ is not
invertible along $C_1$ then $\sigma(C_1)<\overline\sigma(\Phi)$.

First suppose that $C_1$ is a section over $C$. Let $p_1\in C_1$ be a generic point.
Then $p=\pi(p_1)$ is a generic point  of $C$. There exist permissible
parameters $(x,y,z)$ at $p$ such that
$$
\begin{array}{ll}
u&=(x^ay^b)^k\\
v&=P(x^ay^b)+x^cy^d
\end{array}
$$
with
$$
\text{min}\{\frac{c}{a},\frac{d}{b}\}<\text{ ord }(P)<\text{ max }\{\frac{c}{a},\frac{d}{b}\}
$$
and
$$
\text{min}\{\frac{c}{a},\frac{d}{b}\}<k, k\not\,\mid\text{ord}(P).
$$
We may assume, after possibly interchanging $x$ and $y$ that
$$
\overline \sigma(\Phi)=\sigma(C)=(\mid c-a\text{ ord }(P)\mid, \mid d-b\text{ ord }(P)\mid).
$$
Assume that $p_1$ has permissible parameters $(x_1,y_1,z)$ such that
$$
x=x_1,y=x_1y_1
$$
and $x_1=y_1=0$ are local equations of $C_1$ at $p_1$.
$$
\begin{array}{ll}
u&=(x_1^{a+b}y_1^b)^k\\
v&=P(x_1^{a+b}y_1^b)+x_1^{c+d}y_1^d
\end{array}
$$
$k\not\,\mid\text{ord}(P)$ implies by Lemma \ref{Lemma1067} that
$$
\sigma(C_1)\le S(\mid (c+d)-(a+b)\text{ ord }(P)\mid,\mid d-b\text{ ord }(P)\mid).
$$
If $c-a\text{ ord }(P)>0$ and $d-b\text{ ord }(P)<0$ then
$$
0\le (c-a\text{ ord }(P)+(d- b\text{ ord }(P))<c-\text{ ord }(P)
$$
so that $\sigma(C_1)<\sigma(C)$.
 
If $c-a\text{ ord }(P)<0$ and $d-b\text{ ord }(P)>0$ then
$$
0\ge (c-a\text{ ord }(P)+(d- b\text{ ord }(P))>c-\text{ ord }(P)
$$
so that $\sigma(C_1)<\sigma(C)$.

Assume that $p_1$ has permissible parameters $(x_1,y_1,z)$ such that
$$
x=x_1y_1, y=y_1
$$
and $x_1=y_1=0$ are local equations of $C_1$ at $p_1$.
$$
\begin{array}{ll}
u&=(x_1^ay_1^{a+b})^k\\
v&=P(x_1^ay_1^{a+b})+x_1^cy_1^{c+d}
\end{array}
$$
$k\not\,\mid\text{ord}(P)$ implies
$$
\sigma(C_1)\le S(\mid (c+d)-(a+b)\text{ ord }(P)\mid,\mid c-a\text{ ord }(P)\mid).
$$
 
If $c-a\text{ ord }(P)>0$ and $d-b\text{ ord }(P)<0$ then
$$
0\le (c-a\text{ ord }(P)+(d- b\text{ ord }(P))<c-\text{ ord }(P)
$$
so that $\sigma(C_1)=-\infty$.

If $c-a\text{ ord }(P)<0$ and $d-b\text{ ord }(P)>0$ then
$$
0\ge (c-a\text{ ord }(P)+(d- b\text{ ord }(P))>c-\text{ ord }(P)
$$
so that $\sigma(C_1)=-\infty$.

Now suppose that $C_1\subset \pi^{-1}(C)$ is an exceptional 2 curve.
Then $p=\pi(C_1)$ satisfies (\ref{eq111}). 
Let $q=\Phi(p)$. $q$ is a 1 point. If $E_1,E_2,E_3$ are the components of $E_X$ containing $p$, then $\Phi(E_1)=\Phi(E_2)=\Phi(E_3)=q$. Suppose that $P(t)=\sum a_it^i$. Since $q$ is a 1 point,
we may replace $v$ with $v-\sum a_{ik}u^i$ so that $k\not\,\mid\text{ord}(P)$.
We may also assume, after possibly
interchanging $(x,y,z)$ that
$$
\mid f-c\text{ ord }(P)\mid\ge\mid e-b\text{ ord }(P)\mid\ge \mid d-a\text{ ord }(P)\mid.
$$

If $C$ has local equations $x=z=0$,
then
$$
\sigma(C)=(\mid f-c\text{ ord}(P)\mid,\mid d-a\text{ ord}(P)\mid)
$$
 and a generic point of $C_1$ has permissible
parameters $(x_1,y,z_1)$ where
$$
x=x_1, z=x_1(z_1+\alpha)
$$
with $\alpha\ne 0$, and $x_1=y=0$ are local equations of $C_1$.
Set $\overline x_1=x_1(z_1+\alpha)^{-\frac{c}{a+c}}$.
$$
\begin{array}{ll}
u&=(\overline x_1^{\overline a}y^{\overline b})^{\lambda k}\\
v&=P((\overline x_1^{\overline a}y^{\overline b})^{\lambda})
+\overline x_1^{d+f}y^e(z_1+\alpha)^{f-\frac{(d+f)c}{a+c}}
\end{array}
$$
where
$\lambda=(a+c,b)$, $a+c=\overline a\lambda$, $b=\overline b\lambda$.

If $\sigma(C_1)\ge 0$, then $\lambda k\not\,\mid\lambda\text{ord}(P)$ implies
$$
\sigma(C_1)=S(\mid(d+f)-(a+c)\text{ ord }(P)\mid, \mid e-b\text{ ord }(P)\mid).
$$

Similarily, if $C$ has local equations $y=z=0$,
$$
\sigma(C)=(\mid f-c\text{ ord}(P)\mid,\mid d-a\text{ ord}(P)\mid)
$$
and if
$\sigma(C_1)\ge 0$, $k\not\,\mid\text{ord}(P)$
implies
$$
\sigma(C_1)=S(\mid d-a\text{ ord }(P)\mid, \mid (e+f)-(b+c)\text{ ord }(P)\mid).
$$

If $C$ has local equations $x=y=0$, then 
$$
\sigma(C)=(\mid e-b\text{ ord}(P)\mid,\mid d-a\text{ ord}(P)\mid)
$$
and if $\sigma(C_1)\ge 0$, then
$$
\sigma(C_1)=S(\mid f-c\text{ ord}(P)\mid, \mid(d+e)-(a+b)\text{ord }(P)\mid).
$$

If one of $ f-c\text{ ord}(P),e-b\text{ ord}(P),d-a\text{ ord}(P)$ is zero,
then $\sigma(C_1)=-\infty$.

{\bf Case 1}
Suppose that $f-c\text{ ord }(P)>0$, $e-b\text{ ord }(P)>0$, $d-a\text{ ord }(P)<0$.
Then 
$$
\overline\sigma(\Phi)=\sigma(C)=(\mid f-c\text{ ord }(P)\mid,\mid d-a\text{ ord }(P)\mid)
$$
and $x=z=0$ are local equations of $C$.
$$
0\le (d-a\text{ ord }(P))+(f-c\text{ ord }(P))<f-c\text{ ord }(P)
$$
implies $\sigma(C_1)=-\infty$.

{\bf Case 2}
Suppose that $f-c\text{ ord }(P)>0$, $e-b\text{ ord }(P)<0$, $d-a\text{ ord }(P)>0$.
Then 
$$
\overline\sigma(\Phi)=\sigma(C)=(\mid f-c\text{ ord }(P)\mid,\mid e-b\text{ ord }(P)\mid)
$$
and $y=z=0$ are local equations of $C$.
$$
0\le (e-b\text{ ord }(P))+(f-c\text{ ord }(P))<f-c\text{ ord }(P)
$$
implies $\sigma(C_1)=-\infty$.

{\bf Case 3}
Suppose that $f-c\text{ ord }(P)>0$, $e-b\text{ ord }(P)<0$, $d-a\text{ ord }(P)<0$.
Then 
$$
\overline\sigma(\Phi)=\sigma(C)=(\mid f-c\text{ ord }(P)\mid,\mid e-b\text{ ord }(P)\mid)
$$
and $y=z=0$ are local equations of $C$.
$$
0\le (e-b\text{ ord }(P))+(f-c\text{ ord }(P))<f-c\text{ ord }(P)
$$
implies $\sigma(C_1)<\sigma(C)$.

{\bf Case 4}
Suppose that $f-c\text{ ord }(P)<0$, $e-b\text{ ord }(P)>0$, $d-a\text{ ord }(P)>0$.
Then 
$$
\overline\sigma(\Phi)=\sigma(C)=(\mid f-c\text{ ord }(P)\mid,\mid e-b\text{ ord }(P)\mid)
$$
and $y=z=0$ are local equations of $C$.
$$
0\ge (e-b\text{ ord }(P))+(f-c\text{ ord }(P))>f-c\text{ ord }(P)
$$
implies $\sigma(C_1)<\sigma(C)$.

{\bf Case 5}
Suppose that $f-c\text{ ord }(P)<0$, $e-b\text{ ord }(P)>0$, $d-a\text{ ord }(P)<0$.
Then 
$$
\overline\sigma(\Phi)=\sigma(C)=(\mid f-c\text{ ord }(P)\mid,\mid e-b\text{ ord }(P)\mid)
$$
and $y=z=0$ are local equations of $C$.
$$
0\ge (e-b\text{ ord }(P))+(f-c\text{ ord }(P))>f-c\text{ ord }(P)
$$
implies $\sigma(C_1)=-\infty$.

{\bf Case 6}
Suppose that $f-c\text{ ord }(P)<0$, $e-b\text{ ord }(P)<0$, $d-a\text{ ord }(P)>0$.
Then 
$$
\overline\sigma(\Phi)=\sigma(C)=(\mid f-c\text{ ord }(P)\mid,\mid d-a\text{ ord }(P)\mid)
$$
and $x=z=0$ are local equations of $C$.
$$
0\ge (d-a\text{ ord }(P))+(f-c\text{ ord }(P))>f-c\text{ ord }(P)
$$
implies $\sigma(C_1)=-\infty$.

We conclude that if $C_1\subset\pi^{-1}(C)$ is a 2 curve such that $m_q{\cal O}_{X_1}$
is not invertible along $C_1$, then $\sigma(C_1)<\overline \sigma(\Phi)$.

By Theorem \ref{Theorem63}, induction on the number of 2 curves $C\subset X$ such that $\sigma(C)=\overline\sigma(\Phi)$,
and induction on $\overline \sigma(\Phi)$,
we achieve the conclusions of the Lemma.
\end{pf}

\begin{Lemma}\label{Lemma1044}
Suppose that $\Phi:X\rightarrow S$ is strongly prepared,
$q\in S$ is such that $m_q{\cal O}_X$ is not invertible
and the forms (\ref{eq611}), (\ref{eq108}) and (\ref{eq111})
do not hold at  any point $p\in X$ where $m_q{\cal O}_{X,p}$ is not invertible.

Then there exists a sequence of blowups of nonsingular curves  $X_1\rightarrow X$
which are not 2 curves  such that the induced map
$\Phi_1:X_1\rightarrow S$ is strongly prepared, 
$A(\Phi_1,E)< A(\Phi_1)=A(\Phi)$ if $E$ is an exceptional component of $E_{X_1}$ for
$X_1\rightarrow X$,
the forms (\ref{eq611}), (\ref{eq108}) and (\ref{eq111})
do not hold at  any point $p\in X_1$ where $m_q{\cal O}_{X_1,p}$ is not invertible,
and if $C\subset X_1$
is a curve such that $m_q{\cal O}_{X_1}$ is not invertible along $C$, then $C$ is a 
2 curve.
\end{Lemma}

\begin{pf}
Suppose that $C$ is a  curve such that $m_q{\cal O}_X$ is not invertible along $C$
and $C$ is not a 2 curve. Suppose that $p\in C$. Then one of the following holds:
\begin{enumerate}
\item (\ref{eq612}) holds at $p$, $x=y=0$ are local equations of $C$ at $p$.
\item (\ref{eq109}) holds at $p$ and 
$x=z=0$ with $\frac{c}{a}<k$ or $y=z=0$ with $\frac{d}{b}<k$ are local equations of $C$ at $p$.
\item (\ref{eq110}) holds at $p$, $x=z=0$ or $y=z=0$ are local equations of $C$ at $p$.
\end{enumerate}

At a generic point $p\in C$ (\ref{eq612}) holds. Define
$$
\Omega(C)= k-c>0.
$$
Let
$$
\overline\Omega(\Phi)=\text{max}\left\{
\begin{array}{ll}\Omega(C)&\mid \text{ $C$ is not a 2 curve}\\
&\,\,\,\text{and }m_q{\cal O}_X
\text{ is not invertible along }C.
\end{array}\right\}
$$

Suppose that $\Omega(C)=\overline \Omega(\Phi)$. Let $\pi:X_1\rightarrow X$ be the blowup of $C$. The forms (\ref{eq611}), (\ref{eq108}) and (\ref{eq111}) cannot hold at points
of $X_1$. 
By Theorem \ref{Theorem63}, we need only verify that $\Omega(C_1)<\overline\Omega(\Phi)$
if $C_1$ is a curve in $\pi^{-1}(C)$ such that  $m_q{\cal O}_{X_1}$
is not invertible along $C_1$ and $C_1$ is not a 2 curve. We then  have $\pi(C_1)=C$.

 Let $p_1$ be a generic point of $C_1$, $p=\pi(p_1)$.
(\ref{eq612}) holds at $p$ since $p$ is a generic point of $C$. $p_1\in\pi^{-1}(p)$ is a 1 point. Then $p_1$ has permissible parameters $(x_1,y_1,z_1)$ such that
$$
x=x_1, y=x_1(y_1+\alpha).
$$
$$
\begin{array}{ll}
u&=x_1^k\\
v&=x_1^{c+1}(y_1+\alpha).
\end{array}
$$
$m_q{\cal O}_{X_1,p_1}$ is invertible if $\alpha\ne 0$. If $\alpha=0$ and $m_q{\cal O}_{X_1,p_1}$ is not invertible, then $x_1=y_1=0$ are local equations of the curve
$C_1\subset X_1$ through $p$ on which $m_q{\cal O}_{X_1}$ is not invertible. 
$$
0<\Omega(C_1)=k-(c+1)<\Omega(C)=\overline\Omega(\Phi).
$$
By induction on the number of curves $C$ on $X$ such that $\Omega(C)=\overline\Omega(\Phi)$,
we achieve the conclusions of the Lemma.
\end{pf}

\begin{Lemma}\label{Lemma1045}
Suppose that $\Phi:X\rightarrow S$ is strongly prepared, $q\in S$ is such that 
$m_q{\cal O}_X$ is not invertible,
the forms (\ref{eq611}), (\ref{eq108}) and (\ref{eq111})
do not hold at  any point $p\in X$ where $m_q{\cal O}_{X,p}$ is not invertible,
and if $C\subset X$
is a curve such that $m_q{\cal O}_{X}$ is not invertible along $C$, then $C$ is a 
2 curve.

Then there exists a sequence of blowups of 2 curves $X_1\rightarrow X$ such that
the induced map
$\Phi_1:X_1\rightarrow S$ is strongly prepared, 
$A(\Phi_1,E)< A(\Phi_1)=A(\Phi)$ if $E$ is an exceptional component of $E_{X_1}$ for
$X_1\rightarrow X$ and $m_q{\cal O}_{X_1}$
is invertible.
\end{Lemma}

\begin{pf} Suppose that $C\subset X$ is a 2 curve such that $m_q{\cal O}_X$ is
not invertible along $C$. Suppose that $p\in C$. Then (\ref{eq613}), (\ref{eq112}),
 (\ref{eq1054}), (\ref{eq1055}) or (\ref{eq1059}) hold at $p$.  At a generic point $p\in C$ (\ref{eq613}) or (\ref{eq1054}) holds.

If $C$ is a 2 curve such that at a generic point of $C$ (\ref{eq613}) holds, define
$$
\omega(C)=\left\{\begin{array}{ll}
S(\mid ka-c\mid,\mid kb-d\mid)&\text{if }ka-c,kb-d\text{ have opposite signs}\\
&\text{and $m_q{\cal O}_X$ is not invertible along $C$}\\
-\infty&\text{otherwise}
\end{array}
\right.
$$ 

If $C$ is a 2 curve such that at a generic point of $C$ (\ref{eq1054}) holds, define
$$
\omega(C)=\left\{\begin{array}{ll}
S(a,b)&\text{if $m_q{\cal O}_X$ is not invertible along $C$}\\
-\infty&\text{otherwise}
\end{array}
\right.
$$ 
$m_q{\cal O}_X$ is not invertible along $C$ if and only if $\omega(C)>0$.
Set
$$
\overline\omega(\Phi)=\text{max}\left\{\begin{array}{ll}
\omega(C)&\mid C\text{ is a 2 curve such that}\\
&\,\,\,m_q{\cal O}_X\text{ is not invertible along }C
\end{array}\right\}
$$

Suppose that $\omega(C)=\overline\omega(\Phi)$. Let $\pi:X_1\rightarrow X$ be the blowup of
$C$.   By Theorem \ref{Theorem63}, we need only verify that $\omega(C_1)<\overline\omega(\Phi)$ if $C_1\subset\pi^{-1}(C)$ is a curve such that
$m_q{\cal O}_{X_1}$ is not invertible along $C_1$. We must have that
$C_1$ is a 2 curve.

Suppose that $C_1$ is a section over $C$. Let $p_1\in C_1$ be a generic point.
Then $p=\pi(p_1)$ is a generic point on $C$. 

Suppose that there exist permissible parameters
$(x,y,z)$ at $p$ such that (\ref{eq613}) holds.
$$
\begin{array}{ll}
u&=(x^ay^b)^k\\
v&=x^cy^d
\end{array}
$$
with
$$
\text{min}\{\frac{c}{a},\frac{d}{b}\}<k<\text{max}\{\frac{c}{a},\frac{d}{b}\}.
$$
After possibly interchanging $x$ and $y$, we may assume that 
$$
\omega(C)=(\mid c-ak\mid,\mid d-bk\mid).
$$

Assume that $p_1$ has permissible parameters $(x_1,y_1,z)$ such that 
$$
x=x_1, y=x_1y_1
$$
and $x_1=y_1=0$ are local equations of $C_1$ at $p_1$.
$$
\begin{array}{ll}
u&=(x_1^{a+b}y_1^b)^k\\
v&=x_1^{c+d}y_1^d
\end{array}
$$
$$
\omega(C_1)=
\left\{\begin{array}{ll}
S(\mid (c+d)-(a+b)k\mid, \mid bk-d\mid)&\text{ if }
(c+d)-(a+b)k, d-bk\text{ have opposite signs}\\
&\text{ and }m_q{\cal O}_{X_1}\text{is not invertible along }C\\
-\infty&\text{ otherwise}
\end{array}\right.
$$

Suppose that $c-ak>0$ and $d-bk<0$.
$$
0\le (c-ak)+(d-bk)<c-ak
$$
implies $\omega(C_1)<\omega(C)$.

Suppose that $c-ak<0$ and $d-bk>0$.
$$
0\ge (c-ak)+(d-bk)>c-ak
$$
implies $\omega(C_1)<\omega(C)$.

Assume that $p_1$ has permissible parameters $(x_1,y_1,z)$ such that 
$$
x=x_1y_1, y=y_1
$$
and $x_1=y_1=0$ are local equations of $C_1$ at $p_1$.
$$
\begin{array}{ll}
u&=(x_1^{a}y_1^{a+b})^k\\
v&=x_1^{c}y_1^{c+d}
\end{array}
$$
$$
\omega(C_1)=
\left\{\begin{array}{ll}
S(\mid (c+d)-(a+b)k\mid, \mid ak-c\mid)&\text{ if }
(c+d)-(a+b)k, c-ak\text{ have opposite signs}\\
&\text{ and }m_q{\cal O}_{X_1}\text{is not invertible along }C\\
-\infty&\text{ otherwise}
\end{array}\right.
$$

Suppose that $c-ak>0$ and $d-bk<0$.
$$
0\le (c-ak)+(d-bk)<c-ak
$$
implies $\omega(C_1)=-\infty$.

Suppose that $c-ak<0$ and $d-bk>0$.
$$
0\ge (c-ak)+(d-bk)>c-ak
$$
implies $\omega(C_1)=-\infty$.

Suppose that $C_1$ is a section over $C$, $p_1\in C_1$ is a generic point and $p\in\pi(p_1)$ is
a generic point on $C$ such that $p$ satisfies (\ref{eq1054}). Then a similar argument shows
that $\omega(C_1)<\omega(C)$.

Suppose that $C_1\subset\pi^{-1}(C)$ is an exceptional 2 curve.
Suppose that $p=\pi(C_1)$ satisfies (\ref{eq112}). Without loss of generality,
$$
\mid f-ck\mid\ge \mid e-bk\mid \ge \mid d-ak\mid.
$$
If $C$ has local equations $x=z=0$ then a generic point of $C_1$ has regular
parameters $(x_1,y,z_1)$ such that 
$$
x=x_1,z=x_1(z_1+\alpha)
$$
(with $\alpha\ne 0$) and $x_1=y=0$ are local equations of $C_1$.
Set $\overline x_1=x_1(z_1+\alpha)^{-\frac{c}{a+c}}$.
$$
\begin{array}{ll}
u&=(\overline x_1^{\overline a}y^{\overline b})^{\lambda k}\\
v&=\overline x_1^{d+f}y^e(z_1+\alpha)^{f-\frac{(d+f)c}{a+c}}
\end{array}
$$
where $\lambda=(a+c,b)$, $a+c=\overline a\lambda$, $b=\overline b\lambda$.
$\omega(C_1)\ge 0$ implies
$$
\omega(C_1)=S(\mid (d+f)-(a+c)k\mid,\mid e-bk\mid).
$$

Similarily, if $C_1$ has local equations $y=z=0$ then $\omega(C_1)\ge 0$ implies
$$
\omega(C_1)=S(\mid (e+f)-(b+c)k\mid,\mid d-ak\mid).
$$

If $C_1$ has local equations $x=y=0$ then $\omega(C_1)\ge0$ implies
$$
\omega(C_1)=S(\mid f-ck\mid,\mid (d+e)-(a+b)k\mid).
$$

If one of $d-ak$, $e-bk$, $f-ck$ is zero, then $\omega(C_1)=-\infty$.

The analysis of 
Cases 1 - 6 of Lemma \ref{Lemma1043} (with $\text{ ord }(P)$ changed to $k$ and
$\sigma$ to $\omega$) shows that $\omega(C_1)<\omega(C)$.

A similar argument shows that $\omega(C_1)<\omega(C)$ if $p=\pi(C)$ satisfies (\ref{eq1055})
or (\ref{eq1059}).

We achieve the conclusions of the Lemma by Theorem \ref{Theorem63}, induction on the number of 2 curves $C\subset X$
such that $\omega(C)=\overline\omega(\Phi)$, and by induction on $\overline \omega(\Phi)$.
\end{pf}

\begin{Theorem}\label{Theorem1068}
Suppose that $\Phi:X\rightarrow S$ is strongly prepared with respect to $D_S$.
 Then there exists a finite
sequence of quadratic transforms $\pi_1:S_1\rightarrow S$ and monoidal transforms
centered at nonsingular curves
$\pi_2:X_1\rightarrow X$ such that the induced morphism $\Phi_1:X_1\rightarrow S_1$
is strongly prepared with respect to $D_{S_1}=\pi_1^{-1}(D_S)_{red}$, and all points of $X_1$ are good for $\Phi_1$.
\end{Theorem}

\begin{pf}
By Remark \ref{Remark1081}, $A(\Phi)=0$ if and only if all points of $X$ are good. Suppose that $A(\Phi)>0$ and $E$ is a component of $E_X$ such that $C(\Phi,E)=C(\Phi)$.
$A(\Phi,E)>0$ implies $\Phi(E)$ is a point $q$.

Let $\pi_1:S_1\rightarrow S$ be the blowup of $q$.  
By Lemmas \ref{Lemma1043}, \ref{Lemma1044}, \ref{Lemma1045} there exists a
sequence of blowups of curves $X_1\rightarrow X$ such that
$\Phi_1:X_1\rightarrow S$ is strongly prepared, $C(\Phi_1)=C(\Phi)$,
 $A(\Phi_1,\overline E)<A(\Phi_1)$
if $\overline E$ is exceptional for $\Phi_1$ and $\Phi_2:X_1\rightarrow S_1$ is a strongly prepared morphism.

By Theorem \ref{Theorem61}, $C(\Phi_2,\tilde E)< C(\Phi)$, where $\tilde E$
is the strict transform of $E$ on $X_1$.

By  induction on the number of components $E$ of $E_X$ such that $C(\Phi,E)=C(\Phi)$, and induction on $C(\Phi)$, we get the conclusions of the Theorem.
\end{pf}
 
\begin{Definition}\label{Def1082} Suppose that $\Phi:X\rightarrow Y$ is a 
dominant morphism of $k$-varieties, (where $k$ is a field of characteristic zero). $\Phi$ is
a monomial morphism if for all $p\in X$ there exists an \'etale neighborhood $U$ of $p$,
uniformizing parameters $(x_1,\ldots,x_n)$ on $U$, regular parameters 
$(y_1,\ldots,y_m)$ in ${\cal O}_{Y,\Phi(p)}$, and a matrix $(a_{ij})$ of nonnegative
integers such that
$$
\begin{array}{ll}
y_1=&x_1^{a_{11}}\cdots x_n^{a_{1n}}\\
&\vdots\\
y_m=&x_1^{a_{m1}}\cdots x_n^{a_{mn}}
\end{array}
$$
\end{Definition}

\begin{Theorem}\label{Theorem1079}
Suppose that $\Phi:X\rightarrow S$ is a dominant morphism from a 3 fold $X$ to
a surface $S$ (over an algebraically closed field $k$ of characteristic zero). Then there exist sequences of
blowups of nonsingular subvarieties $X_1\rightarrow X$ and $S_1\rightarrow S$ such
that the induced map $\Phi_1:X_1\rightarrow S_1$ is a monomial morphism.
\end{Theorem}

\begin{pf}
This follows from Theorem \ref{Theorem58}, the fact that prepared implies strongly prepared, Theorem \ref{Theorem1068} and 
Remark \ref{Remark1072}.
\end{pf}

\section{Toroidalization}
Throughout this section we will assume that $\Phi:X\rightarrow S$ is strongly prepared
with respect to $D_S$,
and all points of $X$ are good.

\begin{Definition}\label{Def1083}(\cite{KKMSD} and \cite{AK})
A normal variety $X$ with a SNC divisor $E_X$ on $X$ is called toroidal if for
every point $p\in X$ there exists an affine toric variety $X_{\sigma}$, a point
$p'\in X_{\sigma}$ and an isomorphism of $k$ algebras
$$
\hat{\cal O}_{X,p}\cong\hat{\cal O}_{X_{\sigma},p'}
$$
such that the ideal of $E_X$ corresponds to the ideal of $X_{\sigma}-T$ (where $T$ is
the torus in $X_{\sigma}$). Such a pair $(X_{\sigma},p')$ is called a local model at
$p\in X$.

A dominant morphism $\Phi:X\rightarrow Y$ of toroidal varieties with 
SNC divisors $D_Y$, $E_X$ on $X$, $Y$ and $\Phi^{-1}(D_Y)\subset E_X$
is called toroidal at $p$, and we will say that $p$ is a toroidal point of $\Phi$,
if with $q=\Phi(p)$, there exist local models $(X_{\sigma},p')$ at $p$,
$(Y_{\tau},q')$ at $q$ and a toric morphism
$\Psi:X_{\sigma}\rightarrow Y_{\tau}$ such that the following diagram commutes
$$
\begin{array}{lll}
\hat{\cal O}_{X,p}&\stackrel{\sim}{\leftarrow}
&\hat{\cal O}_{X_{\sigma},p'}\\
\hat\Phi^*\uparrow&&\hat\Psi^*\uparrow\\
{\cal O}_{Y,q}&\stackrel{\sim}{\leftarrow}
&\hat{\cal O}_{Y_{\tau},q'}
\end{array}
$$

$\Phi:X\rightarrow Y$ is called toroidal (with respect to $D_Y$ and $E_X$) if $\Phi$ is
toroidal at all $p\in X$.
\end{Definition}

\begin{Remark}

\begin{enumerate}
\item If one of the forms (\ref{eq1053}), (\ref{eq1058}) or (\ref{eq1052}) holds at $p\in X$
then $q=\Phi(p)$ is a 2 point.
\item If $q=\Phi(p)$ is a 2 point, then (\ref{eq1051}) cannot hold at $p$, and if (\ref{eq610})
or (\ref{eq107}) hold at $p$, we must have $\alpha\ne 0$, since $uv=0$ is a local equation
of $E_X$.
\end{enumerate}
\end{Remark}

\begin{Lemma} Suppose that $\Phi:X\rightarrow S$ is a morphism from a nonsingular 3 fold
$X$ to a nonsingular surface $S$, $D_S$ is a SNC divisor on $S$ such that
$E_X=\Phi^{-1}(D_S)$ is a SNC divisor on $X$. Then $\Phi$ is a toroidal morphism 
if and only if for all $p\in E_X$ there exist regular parameters
$(x,y,z)$ in $\hat{\cal O}_{X,p}$ $(u,v)$ in ${\cal O}_{S,p}$ such that one of the
following forms hold:
\begin{enumerate}
\item $u=0$ is a local equation for $D_S$.
\begin{enumerate}
\item $xy=0$ is a local equation for $E_X$ and 
$$
\begin{array}{ll}
u&=x^ay^b\\
v&=z
\end{array}
$$
\item $x=0$ is a local equation for $E_X$ and
$$
\begin{array}{ll}
u&=x^a\\
v&=y
\end{array}
$$
\end{enumerate}
\item $uv=0$ is a local equation for $D_S$.
\begin{enumerate}
\item $xyz=0$ is a local equation for $E_X$ and
$$
\begin{array}{ll}
u&=x^ay^bz^c\\
v&=x^dy^ez^f
\end{array}
$$
with 
$$
\text{rank}\left(\begin{array}{lll}
a&b&c\\
d&e&f
\end{array}\right)
=2.
$$
\item $xy=0$ is a local equation for $E_X$ and
$$
\begin{array}{ll}
u&=x^ay^b\\
v&=x^cy^d
\end{array}
$$
with $ad-bc\ne0$.
\item $xy=0$ is a local equation of $E_X$ and
$$
\begin{array}{ll}
u&=(x^ay^b)^k\\
v&=\alpha(x^ay^b)^t+(x^ay^b)^tz
\end{array}
$$
with $a,b>0$, $k,t>0$, $0\ne\alpha\in k$.
\item $x=0$ is a local equation for $E_X$ and
$$
\begin{array}{ll}
u&=x^a\\
v&=x^c(y+\alpha)
\end{array}
$$
with $0\ne\alpha\in k$.
\end{enumerate}
\end{enumerate}
\end{Lemma}

\begin{pf}
 We will first determine the toroidal forms obtainable from a  monomial mapping
$\Lambda:{\bold A}^3\rightarrow {\bold A}^2$ defined by
$$
\begin{array}{ll}
u&=x^ay^bz^c\\
v&=x^dy^ez^f
\end{array}
$$
with 
$$
\text{rank}\left(
\begin{array}{lll}
a&b&c\\
d&e&f
\end{array}
\right)=2.
$$ 
First suppose that that no column of 
$$
\left(
\begin{array}{lll}
a&b&c\\
d&e&f
\end{array}
\right)
$$ 
is zero. Then
$\Lambda^{-1}(D)=E$, where
$xyz=0$ is an equation of $E$, $uv=0$ is an equation of $D$. 

Suppose that $p\in {\bold A}^3$ is a 2 point on $y=z=0$. Then there exists $0\ne\beta\in k$ and
regular parameters $(\overline x,y,z)$ at $p$ such that
$$
\begin{array}{ll}
u&=(\overline x+\beta)^ay^bz^c\\
v&=(\overline x+\beta)^dy^ez^f.
\end{array}
$$
If 
$$
\text{Det}\left(\begin{array}{ll}
b&c\\
e&f
\end{array}
\right)\ne 0,
$$
we can make a permissible change of parameters to get 2.(b).

Suppose that
$$
\text{Det}\left(\begin{array}{ll}
b&c\\
e&f
\end{array}
\right)= 0.
$$
There exist natural numbers $\overline b,\overline c$ such that $\overline b$, $\overline c>0$,
$(\overline b,\overline c)=1$,
$$
\begin{array}{ll} 
u&=(\overline x+\alpha)^a(y^{\overline b}z^{\overline c})^k\\
v&=(\overline x+\alpha)^d(y^{\overline b}z^{\overline c})^t.
\end{array}
$$

After possibly interchanging $u$ and $v$, we can assume that $k>0$ and $t\ge 0$.
If $t>0$ we get the form 2.(c). If $t=0$, we get the form 1.(a).

Suppose that $p\in {\bold A}^3$ is a 1 point on $z=0$. Then there exist $0\ne\alpha,\beta\in k$,
and regular parameters $(\overline x,\overline y,z)$ at $p$ such that
$$
\begin{array}{ll}
u&=(\overline x+\alpha)^a(\overline y+\beta)^bz^c\\
v&=(\overline x+\alpha)^d(\overline y+\beta)^ez^f.
\end{array}
$$
After possibly interchanging $u$ and $v$ we may assume that $c>0$.
If $f>0$ we get the form 2.(d). If $f=0$ we get the form 1.(b).

Now suppose that a column of 
$$
\left(
\begin{array}{lll}
a&b&c\\
d&e&f
\end{array}
\right)
$$ 
is zero.

After possibly interchanging $(x,y,z)$, we may assume that $c=f=0$. Then 
$\Lambda^{-1}(D)=E$ where $xy=0$ is an
equation of $E$, $uv=0$ is an equation of $D$. We get the
forms 2.(b), 2.(d) or 1.(b).

Conversely, suppose that the forms 1. and 2. hold at all points of $E_X$ and $p\in E_X$. By Lemma \ref{Lemma1036}, there exists an \'etale neighborhood $U$ of $p$ and
uniformizing parameters $(x,y,z)$ on $U$ such that a form 1. or 2. holds at $p$. Working 
backwards through the above proof, we see that $\Phi$ is toroidal at $p$.
\end{pf}

We will call a point $p\in X$ a non toroidal point if $\Phi$ is not toroidal at $p$.

\begin{Lemma}\label{Lemma1046}
The locus of non toroidal points is Zariski closed of pure codimension 1 in $X$, and is a SNC divisor.
The image of the non toroidal points in $S$ is a finite set of points.
\end{Lemma}

\begin{pf}
Suppose that $p\in X$ is a non toroidal point. Then $q=\Phi(p)$ is a 1 point, and
thus 
one of the  forms (\ref{eq105}), (\ref{eq106}), (\ref{eq610}) with $t>0$, (\ref{eq1051}) or
(\ref{eq107}) with $c>0$ hold at $p$.

\begin{enumerate}
\item First suppose that $p$ is of the form of (\ref{eq107}). We have $c>0$, and all points nearby
on $x=0$ are non toroidal.
\item
Suppose that $p$ has the form (\ref{eq106}). $\Phi$ is non toroidal on the line $x=y=0$.
\begin{enumerate}
\item Suppose that $c>0$.  Consider the point with regular parameters
$(x,\tilde y+\alpha,\tilde z+\beta)$ with $\alpha\ne 0$. Set
$x=\overline x(\tilde y+\alpha)^{-\frac{b}{a}}$.
Then
$$
\begin{array}{ll}
u&=\overline x^a\\
v&=\overline x^c(\tilde y+\alpha)^{d-\frac{cb}{a}}=\overline x^c(\gamma+\overline y)
\end{array}
$$
which is non toroidal. Thus $\Phi$ is non toroidal on the surface $x=0$.
\item Suppose that $d>0$. Then a similar analysis shows that $\Phi$ is non toroidal 
on the surface $y=0$.
\end{enumerate}
\item Suppose that $p$ has the form (\ref{eq1051}).  A point on $x=y=0$ with 
regular parameters $(x,y,\tilde z+\beta)$ with $\beta\ne 0$ has the form of (\ref{eq106}),
and is thus not toroidal.
\begin{enumerate}
\item Suppose that $a,c>0$. Consider the point with regular parameters
$(x,\tilde y+\alpha,\tilde z+\beta)$ with $\alpha\ne 0$.
Set $x=\overline x(\tilde y+\alpha)^{-\frac{b}{a}}$.
$$
\begin{array}{ll}
u&= \overline x^{a}\\
v&=\overline x^c(\tilde y+\alpha)^{d-\frac{bc}{a}}(\tilde z+\beta)
=\overline x^c(\gamma+\overline y).
\end{array}
$$
Thus $\Phi$ is  non toroidal on the surface $x=0$.
\item Suppose that $b,d>0$. Then $\Phi$ is non toroidal on $y=0$.
\end{enumerate}
Since $a,b>0$ (by assumption) one of the cases (a) or (b) must hold.
\item Suppose that $p$ has the form (\ref{eq610}). We have $t>0$.
Consider a nearby point with regular parameters $(x,\tilde y+\overline \alpha,\tilde z+\overline \beta)$ with $\overline \alpha\ne 0$.
Set $x=\overline x(y+\overline \alpha)^{-\frac{b}{a}}$.
$$
\begin{array}{ll}
u&=\overline x^{am}\\
v&=\overline x^{at}(\alpha+\overline \beta+\tilde z)=\overline x^{at}(\tilde c+\overline z).
\end{array}
$$
The non toroidal locus  locally contains $x=0$ (and $y=0$).
\item Suppose that $p$ has the form (\ref{eq105}). Since $a,b,c>0$, after possibly
interchanging $(x,y,z)$, we may assume that 
$a,d>0$. Suppose that $(x,\tilde y+\alpha,\tilde z+\beta)$ are regular parameters
at a nearby point (with $\alpha,\beta\ne 0$).
Set $x=\overline x(\tilde y+\alpha)^{-\frac{b}{a}}(\tilde z+\beta)^{-\frac{c}{a}}$.
$$
\begin{array}{ll}
u&= \overline x^a\\
v&=\overline x^d(\tilde y+\alpha)^{e-\frac{db}{a}}(\tilde z+\beta)^{f-\frac{dc}{a}}
=\overline x^d(\gamma+\overline y)
\end{array}
$$
In a similar way, we see that nearby 2 points on $x=0$ are non toroidal.
Thus the non toroidal locus locally contains $x=0$.
\end{enumerate}
\end{pf}

Suppose that $p\in X$ is a 1 point such that $\Phi(p)=q$ is a 1 point.
A form (\ref{eq107}) holds at $p$. $c-a$ is independent of permissible parameters
$(u,v)$ at $q$ and $(x,y,z)$ at $p$ of the form of (\ref{eq107})
(since $u=0$ must be a local equation of $D_S$). Define
$$
I(\Phi,p)=c-a.
$$
$I(\Phi,p)$ is locally constant. Thus if $E$ is a component of $E_X$ and $p_1$, $p_2$ are two 1 points in $E$ such that
$\Phi(p_1)$ and $\Phi(p_2)$ are 1 points, then $I(\Phi,p_1)=I(\Phi,p_2)$.
We can thus define
$$
I(\Phi,E)=I(\Phi,p)
$$
if $p\in E$ is a 1 point such that $\Phi(p)$ is a 1 point. Let 
$$
B_{\Phi}=\{q\in S\mid q\text{ is the image of a non toroidal point by }\Phi\}.
$$
Define 
$$
I(\Phi)=\text{max}\{I(\Phi,p)\mid p\in\Phi^{-1}(B_{\Phi})\text{ is a 1 point}\}.
$$
\begin{Remark}\label{Remark1097} If $p\in\Phi^{-1}(B_{\Phi})$ is a
  toroidal point then $I(\Phi,p)< 0$.
\end{Remark}

\begin{Lemma}\label{Lemma1073} Suppose that $q\in B_{\Phi}$. Let
$\pi:S_1\rightarrow S$ be the blowup of $q$. Let $U$ be the largest open set of $X$
such that the rational map $X\rightarrow S_1$ is a morphism $\Phi_1:U\rightarrow S_1$.
Then $\Phi_1$ is strongly prepared, and all points of $U$ are good for $\Phi_1$.

Suppose that $p\in U\cap \Phi^{-1}(q)$ is a 1 point. If $I(\Phi,p)\le 0$, then
$\Phi_1$ is toroidal at $p$. If $I(\Phi,p)>0$, then $I(\Phi_1,p)<I(\Phi,p)$.

The locus of points where $m_q{\cal O}_X$ is not invertible is a union of curves
which make SNCs with $\overline B_2(X)$. These points have one of the forms
(\ref{eq613}), (\ref{eq112}), (\ref{eq612}), (\ref{eq110}) or (\ref{eq109}) of Lemma
\ref{Lemma62}.
\end{Lemma}

\begin{pf}
$\Phi_1$ is strongly prepared by Lemma \ref{Lemma1084}. All points of $U$ are good for $\Phi_1$,
as follows by the analysis  in Lemma \ref{Lemma62}. 
The locus of points where $m_q{\cal O}_X$ is not invertible 
is a union of curves which make SNCs with $\overline B_2(X)$
 by Theorem \ref{Theorem63}.

Suppose that $p\in X$ is such that $m_q{\cal O}_X$ is not invertible at $p$.
$p$ is a good point and $\Phi(p)=q$ a 1 point implies $p$ has one of the forms
(\ref{eq105}), (\ref{eq106}), (\ref{eq610}), (\ref{eq1051}) or (\ref{eq107}).
By Lemma \ref{Lemma62}, $p$ must have one of the forms (\ref{eq613}), (\ref{eq112}),
(\ref{eq612}), (\ref{eq110}), or (\ref{eq109}).

Suppose that $p\in\Phi^{-1}(q)\cap U$ is a 1 point. Then
$$
\begin{array}{ll}
u&=x^a\\
v&=x^c(\alpha+y)
\end{array}
$$
where $u=0$ is a local equation of $D_S$ at $q$, with either $a\le c$ or
$c<a$ and $\alpha\ne0$. Suppose that $I(\Phi,p)=c-a\le 0$.

If $c<a$, $\alpha\ne 0$ and  we have permissible parameters $(u_1,v_1)$ at $q_1=\Phi_1(p)$ such that
$$
u=u_1v_1,
v=v_1
$$
so that $q_1$ is a 2 point. There exists regular parameters $(\overline x,\overline y,z)$
in $\hat{\cal O}_{X,p}$ and $0\ne\overline\alpha\in k$ such that
$$
\begin{array}{ll}
u&=\overline x^a(\overline \alpha+\overline y)\\
v&=\overline x^c,
\end{array}
$$
$$
\begin{array}{ll}
u_1&=\overline x^{a-c}(\overline\alpha+\overline y)\\
v_1&=\overline x^c\\
\end{array}
$$
$(\overline x,\overline y,z)$ are thus permissible parameters for $(v_1,u_1)$ at $p$,
and $p$ is a toroidal point for $\Phi_1$.

If $c=a$,
$$
u_1=u,
v_1=\frac{v}{u}-\alpha,
$$
$(u_1,v_1)$ are permissible parameters for $\Phi_1$ at $q_1=\Phi_1(p)$, and 
$$
\begin{array}{ll}
u_1&=x^a\\
v_1&=\frac{v}{u}-\alpha=y
\end{array}
$$
so that $p$ is a toroidal point for $\Phi_1$.

Suppose that $I(\Phi,p)>0$. Then $(u_1,v_1)$ are permissible parameters at $q_1=\Phi_1(p)$, where
$$
u=u_1, v=u_1v_1.
$$
$q_1$ is a 1 point.
$$
\begin{array}{ll}
u_1&=x^a\\
v_1&=x^{c-a}(\alpha+y).
\end{array}
$$
Thus $I(\Phi_1,p)=c-2a<I(\Phi,p)$.
\end{pf}

\begin{Lemma}\label{Lemma1074} Suppose that $C\subset X$ is a 2 curve
such that $q=\Phi(C)$ is a 1 point, if $p\in C$ then $p$ satisfies (\ref{eq613})
or (\ref{eq112}) and $m_q{\cal O}_X$ is not invertible along $C$. Let $\pi:X_1\rightarrow X$
be the blowup of $C$, $\Phi_1=\Phi\circ\pi$. Then $\Phi_1:X_1\rightarrow S$ is
strongly prepared, all points of $X_1$ are good points for $\Phi_1$, and if $m_q{\cal O}_{X_1,p_1}$ is
not invertible at a point $p_1\in\pi^{-1}(C)$, then $p_1$ satisfies (\ref{eq613}) or (\ref{eq112}). If $p_1\in\pi^{-1}(C)$ is a 1 point then $I(\Phi_1,p_1)<I(\Phi)$.
\end{Lemma}

\begin{pf} Suppose that $p\in C$ satisfies (\ref{eq613}). Then all points of $\pi^{-1}(p)$
are strongly prepared and are good points for $\Phi_1$. 
Let $p_1\in\pi^{-1}(p)$ be a 1 point. $\hat{\cal O}_{X_1,p_1}$ has regular parameters 
$(x_1,y_1,z)$ such that
$$
x=x_1, y=x_1(y_1+\alpha)
$$
with $\alpha\ne 0$.
$$
\begin{array}{ll} 
u&=(x_1^{a+b}(y_1+\alpha)^b)^k=\overline x_1^{(a+b)k}\\
v&=x_1^{c+d}(y_1+\alpha)^d=\overline x_1^{c+d}(\overline\alpha+\overline y)
\end{array}
$$
By (\ref{eq613}) $c-ak$, $d-bk$ have opposite signs.
$$
\begin{array}{ll}
I(\Phi_1,p_1)&=(c+d)-(a+b)k\\
&=(c-ak)+(d-bk)<\text{max}\{c-ak,d-bk\}\\
&\le I(\Phi).
\end{array}
$$
If $p$ satisfies (\ref{eq112}), then all points of $\pi^{-1}(p)$ are strongly prepared
good points.
\end{pf}

\begin{Lemma}\label{Lemma1075} Suppose that $C\subset X$ is a curve such that
$q=\Phi(C)$ is a 1 point, $m_q{\cal O}_X$ is not invertible along $C$, and $C$ is not
a 2 curve.

Let $\pi:X_1\rightarrow X$ be the blowup of $C$, $\Phi_1=\Phi\circ\pi$. Then 
$\Phi_1:X_1\rightarrow S$ is strongly prepared and  all points of $X_1$ are good points
for $\Phi_1$.
If $p_1\in\pi^{-1}(C)$ is a 1 point, then
$$
I(\Phi,p)<I(\Phi_1,p_1)\le 0.
$$
\end{Lemma}

\begin{pf} Suppose that $p\in C$. $p$ satisfies (\ref{eq612}),
(\ref{eq110}) or (\ref{eq109}). $\Phi_1$ is strongly prepared, and all points of
$X_1$ are good points for $\Phi_1$. A generic point $p\in C$ satisfies (\ref{eq612}),
and $x=y=0$ is a local equation of $C$ at $p$. Suppose
that $p_1\in\pi^{-1}(p)$ is a 1 point. Then $p_1$ has permissible parameters $(x_1,y_1,z)$
such that
$$
x=x_1, y=x_1(y_1+\alpha),
$$
$$
\begin{array}{ll}
u&=x_1^k\\
v&=x_1^{c+1}(y_1+\alpha).
\end{array}
$$
$I(\Phi_1,p_1)=(c+1)-k\le 0$
since $c<k$ by (\ref{eq612}).
\end{pf}

\begin{Theorem}\label{Theorem1076} Suppose that $\Phi:X\rightarrow S$ is
strongly prepared and all points $p\in X$ are good points for $\Phi$. Then there
exists a sequence of quadratic transforms $S_1\rightarrow S$ and monodial transforms
centered at nonsingular curves,
$X_1\rightarrow X$, such that the induced map $\Phi_1:X_1\rightarrow S_1$ is
strongly prepared, all points of $X_1$ are good for $\Phi_1$ and $I(\Phi_1)\le 0$.
\end{Theorem}

\begin{pf} Suppose that $I(\Phi)>0$. Suppose that $E$ is a component of $E_X$ such that
$I(\Phi,E)=I(\Phi)$. Then $\Phi(E)$ is a single 1 point $q$. Let $\pi_1:S_1\rightarrow S$ be
the blowup of $q$.

By Lemmas \ref{Lemma1044}, \ref{Lemma1073} and  \ref{Lemma1075}, there exists a sequence of
blowups of curves $C$ (which are not 2 curves) $X_1\rightarrow X$ such that if
$\Phi_1:X_1\rightarrow S$ is the induced map, $\Phi_1$ is strongly prepared, all points
of $X_1$ are good for $\Phi_1$, if $m_q{\cal O}_{X_1,p}$ is not invertible then
(\ref{eq613}) or (\ref{eq112}) holds at $p$, and all curves in $X_1$ along which $m_q{\cal O}_{X_1}$ are not invertible are 2 curves. We further have that $I(\Phi_1,\overline E)\le 0$ if
$\overline E$ is exceptional for $X_1\rightarrow X$.

By Lemmas \ref{Lemma1045} and \ref{Lemma1074}, there exists a sequence of blowups of
2 curves $C$, $X_2\rightarrow X_1$ such that if $\Phi_2:X_2\rightarrow S$ is the
induced map, $\Phi_2$ is strongly prepared, all points of $X_2$ are good for $\Phi_2$,
$m_q{\cal O}_{X_2}$ is invertible, and if $\overline E$ is exceptional for $X_2\rightarrow X$,
then $I(\Phi_2,\overline E)<I(\Phi)$.

Let $\overline \Phi:X_2\rightarrow S_1$ be the induced map. By Lemma \ref{Lemma1073},
$\overline\Phi$ is strongly prepared, all points of $X_2$ are good for $\overline\Phi$,
$I(\overline\Phi)\le I(\Phi)$, and if $\overline E$ is a component of $E_{X_2}$ which 
contains a 1 point $q$ such that $\overline\Phi(p)\in\pi_1^{-1}(q)$, then
$I(\overline\Phi,\overline E)<I(\Phi)$.

The Theorem now follows by induction on the number of components $E$ of $X$ such that
$I(\Phi,E)=I(\Phi)$, and induction on $I(\Phi)$.
\end{pf}

\begin{Theorem}\label{Theorem1077}
Suppose that $\Phi:X\rightarrow S$ is strongly prepared with respect to $D_S$,
$E_X=\Phi^{-1}(D_S)_{red}$, all points $p\in X$ are good
points for $\Phi$ and $I(\Phi)\le 0$. Then there exist  sequences of quadratic transforms
$\pi_1:S_1\rightarrow S$ and monodial transforms centered at nonsingular curves $\pi_2:X_1\rightarrow X$
such that the induced map $\Phi_1:X_1\rightarrow S_1$ is toroidal with respect to
$D_{S_1}=\pi^{-1}(D_S)_{red}$ and $E_{X_1}=\pi_2^{-1}(E_X)_{red}$.
\end{Theorem}

\begin{pf} Suppose that $E$ is a component of $E_X$ such that $\Phi$ is not toroidal
along $E$. If $p\in E$ is a generic point, then at $p$ we have an expression
$$
\begin{array}{ll}
u&=x^a\\
v&=x^c(\alpha+y)
\end{array}
$$
with $c>0$. Thus there exists a point $q\in S$ such that $\Phi(E)=q$. $q$ is necessarily a 1 point. 

Let $\pi:S_1\rightarrow S$ be the blowup of $q$. By Lemmas \ref{Lemma1044} and
\ref{Lemma1075}, there exists a sequence of blowups of nonsingular curves (which are
not 2 curves) $X_1\rightarrow X$ such that if $\Phi_1:X_1\rightarrow S$ is the induced
morphism, then $\Phi_1$ is strongly prepared, all points of $X_1$ are good for $\Phi_1$,
$I(\Phi_1)\le 0$, the locus of points $p_1$ of $X_1$ such that $m_q{\cal O}_{X_1,p_1}$ is not invertible is a union of 2 curves, and if $m_q{\cal O}_{X_1,p_1}$ is not invertible, then
$p_1$ satisfies (\ref{eq613}) or (\ref{eq112}).

Suppose that $C$ is a 2 curve on $X_1$ such that $m_q{\cal O}_{X_1}$ is not invertible
along $C$. A generic point $p$ of $C$ satisfies (\ref{eq613}). Let $E_1$ be the component of
$E_{X_1}$ with local equation $x=0$ at $p$, $E_2$ be the component of $E_{X_1}$ with local
equation $y=0$ at $p$.
$$
I(\Phi_1,E_1)=c-ak,\,\,\, I(\Phi_1,E_2)=d-bk.
$$
Since $m_q{\cal O}_{X_1,p}$ is not invertible, either $0<d-bk$ or $0<c-ak$, a
contradiction since $I(\Phi_1)\le 0$. Thus $m_q{\cal O}_{X_1}$ is invertible and $\Phi_1:X_1\rightarrow S$ induces a morphism
$\overline\Phi:X_1\rightarrow S_1$. By Lemma \ref{Lemma1073}, $\overline\Phi$ is
strongly prepared, all points of $X_1$ are good for $\overline\Phi$, and $I(\overline\Phi)\le 0$. Further, if $p\in X_1$ is a 1 point such that $p\in\Phi_1^{-1}(q)$,  then $\overline\Phi$
is toroidal at $p$.

By induction on the number of components of $E_X$ along which $\Phi$ is not toroidal, we achieve the conclusions of the Theorem.
\end{pf}

\begin{Theorem}\label{Theorem1078} Suppose that $\Phi:X\rightarrow S$ is
a dominant morphism from a 3 fold $X$ to a surface $S$ (over an algebraically closed
field $k$ of characteristic 0)
and $D_S$ is a reduced 1 cycle on $S$ such that $E_X=\Phi^{-1}(D_S)_{red}$ contains
$\text{sing}(X)$ and $\text{sing}(\Phi)$. Then there exist sequences of
blowups of nonsingular subvarieties $\pi_1:X_1\rightarrow X$ and $\pi_2:S_1\rightarrow S$  
 such that the induced morphism $X_1\rightarrow S_1$ is a
toroidal morphism with respect to $\pi_2^{-1}(D_S)_{red}$ and $\pi_1^{-1}(E_X)_{red}$.
\end{Theorem}

\begin{pf} This follows from Theorem \ref{Theorem58}, the fact that prepared implies
strongly prepared and Theorems \ref{Theorem1068}, \ref{Theorem1076} and \ref{Theorem1077}.
\end{pf}

\section{Glossary of Notations and Definitions}

$\nu(p)$, Definition \ref{Def1086}.

$\gamma(p)$, Definition \ref{Def1086}.

$\tau(p)$, Definition \ref{Def1086}.

$S_r(X)$, $\overline S_r(X)$, After Definition \ref{Def1086}.

$B_2(X)$, $\overline B_2(X)$, $B_3(X)$, Before Definition \ref{Def1087}.

SNCs with $\overline B_2(X)$, Definition \ref{Def1087}.

r small, Definition \ref{Def1091}.

r big, Definition \ref{Def1091}.

1 point, 2 point, 3 point, Definition \ref{Def650} and before Definition \ref{Def1071}.

1-resolved, Definition \ref{Def1092}.

$\overline \nu(p)$, After Definition \ref{Def1092}.

$\sigma(p)$, Before Lemma \ref{L6}.

$\delta(p)$, before Lemma \ref{L8}.

$\text{Inv}(p)$, Before Theorem \ref{T3}.

$A_r(X)$, Definition \ref{Def1094}.

$\overline A_r(X)$, Definition \ref{Def1093}.

$C_r(X)$, Definition \ref{Def1095}.

$(E)$, Definition \ref{Def1096}.

$*$-permissible parameters, After Definition \ref{Def1071}.

$\nu_E(f)$, Before Definition \ref{Def61}.

$A(\Phi,p)$, Definition \ref{Def61}.

$C(\Phi,p)$, Definition \ref{Def61}.

$A(\Phi,E)$, After Definition \ref{Def61}.

$C(\Phi,E)$, After Definition \ref{Def61}.

$A(\Phi)$, Before Lemma \ref{Lemma1038}.

$C(\Phi)$, Before Lemma \ref{Lemma1038}.

$I(\Phi,p)$, $I(\Phi,E)$, $I(\Phi)$, Before Remark \ref{Remark1097}.

$B_{\Phi}$,  Before Remark \ref{Remark1097}.

SNC divisor, Definition \ref{Def1098}.

$P_t(x)$, After Definition \ref{Def1098}.

\vskip .2truein

bad point, Definition \ref{Def59}.

\'etale neighborhood, Definition \ref{Def1089}.

good point, Definition \ref{Def59}.

monodial transform, After Definition \ref{Def1098}.

monomial mapping, Definition \ref{Def1082}.

non toroidal point, Before Lemma \ref{Lemma1046}.

permissible monoidal transform, Definition \ref{Def2005}.

permissible parameters, Before Definition \ref{Def650} and \ref{Def650}.

permissible parameters for $C$ at $p$, After Lemma \ref{Lemma651}.

prepared, Definition \ref{Def57}.

resolved, Definition \ref{Def1064}.

strongly prepared, Definition \ref{Def1071}.

toroidal mapping, \ref{Def1083}.

toroidal point, Definition \ref{Def1083}.

weakly permissible monodial transform, Definition \ref{Def1090}.

weakly prepared, Definition \ref{Def1085}.

\ \\ \\
\noindent
Steven Dale Cutkosky, Department of Mathematics, University of
Missouri\\
Columbia, MO 65211, USA\\
cutkoskys@@missouri.edu


\begin{thebibliography}{99}




\bibitem{Ab1} {\sc Abhyankar, S.}, Local uniformization on algebraic surfaces over ground fields of
characteristic $p\ne 0$, Annals of Math, 63 (1956), 491-526.


\bibitem{Ab3} {\sc Abhyankar, S.}, Simultaneous resolution for algebraic surfaces, Amer. J. Math 78 (1956), 761-790.


\bibitem{Ab5} {\sc Abhyankar. S.}, Resolution of singularities of embedded algebraic surfaces,
second edition, Springer Verlag, 1998.

\bibitem{Ab6} {\sc Abhyankar. S.}, Good points of a hypersurface, Advances in Math.
68 (1988), 87-256.

\bibitem{AK} {\sc Abramovich D., Karu, K.}, Weak semistable reduction in characteristic 0,
preprint.

\bibitem{AKMW} {\sc Abramovich, D., Karu, K., Matsuki, K., Wlodarczyk, J.}, Torification and Factorization
of Birational Maps, preprint. 

\bibitem{AKi} {\sc Akbulut, S. and King, H.}, Topology of algebraic sets, MSRI publications 25, Springer-Verlag Berlin.

\bibitem{Chr} {\sc Christensen, C.}, Strong domination/ weak factorization of three dimensional regular local rings,
Journal of the Indian Math Soc., 45 (1981), 21-47.



\bibitem{Co} {\sc Cossart, V.}, Polyedre caracteristique d'une singularite, Thesis, Universite de Paris-Sud, 
Centre d'Orsay (1987).



\bibitem{C1} {\sc Cutkosky, S.D.},
Local Factorization of Birational Maps, Advances in Math. 132, (1997), 167-315.

\bibitem{C2} {\sc Cutkosky, S.D.}, Local Monomialization and Factorization  of Morphisms, 
 Ast\'erisque 260, (1999).

\bibitem{C3} {\sc Cutkosky, S.D.},
Simultaneous resolution of singularities,  Proc. American Math. Soc. 128,
(2000), 1905-1910.

\bibitem{CP} {\sc Cutkosky, S.D. and Piltant, O.}, Monomial 
Resolutions of Morphisms of Algebraic Surfaces, to appear in
the Hartshorne volume of
Communications in Algebra.

\bibitem{CS} {\sc Cutkosky, S.D. and Srinivasan, H.}, An Intrinsic
Criterion for isomorphism of singularities, American. Journal of Math. 115, (1993), 789-821.

\bibitem{SGA} {\sc Grothendieck, A.}, Rev\^etements \'etales et Groupe Fondemental,
Lecture Notes in Math. 224, (1971) Springer-Verlag Heidelberg (1971).


\bibitem{DJ} {\sc de Jong, A.J.}, Smoothness, semistability and Alterations, Publ. Math. I.H.E.S. 83, 1996, 51-93.


\bibitem{H1} {\sc Hironaka, H.}, Resolution of singularities of an algebraic variety over a field of
characteristic zero, Annals of Math, 79 (1964), 109-326.

\bibitem{KKMSD} {\sc Kempf, G., Knudsen, F., Mumford, D., Saint-Donat, B.},
Toroidal embeddings I, LNM 339, Springer Verlag (1973).


\bibitem{Ku} {\sc Kuhlmann, FV}, Valuation theoretic and model theoretic aspects of local uniformization, preprint.




\bibitem{Li} {\sc Lipman, J.}, Introduction to Resolution of Singularities,
in Algebraic Geometry, Arcata 1974, Amer. Math. Soc. proc. Symp. Pure Math. 29 (1975)
187-230.

\bibitem{Ma} {\sc Matsumura, H.}, Commutative Algebra 2nd edition, W.A. Benjamin Co., N.Y.


\bibitem{Moh} {\sc Moh, TT.}, Quasi-Canonical uniformization of hypersurface singularities of cahracteristic zero, Comm.
Algebra 20, 3207-3249 (1992).




\bibitem{Mum} {\sc Mumford, D.}, Red Book,  Lecture Notes in Math.  1358, (1988)
Springer Verlag.

\bibitem{Sa} {\sc Sally, J.}, Regular overrings of regular local rings, Trans. Amer. Math. Soc. 171 (1972) 291-300.



\bibitem{Sh} {\sc Shannon, D.L.}, Monodial transforms, Amer. J. Math, 45 (1973), 284-320.


\bibitem{Sp} {\sc Spivakovsky, M.}, Sandwiched singularities and desingularization of surfaces by normalized Nash transforms,
 Ann. of Math. 131, 1990, 441-491.

\bibitem{T} {\sc Teissier, B.}, Valuations, Deformations and Toric Geometry, preprint.


\bibitem{Vi} {\sc Villamayor, O.}, Constructiveness of Hironaka's resolution, Ann.  Scient. Ecole Norm Sup 22,
1-32, 1989

\bibitem{Z1} {\sc Zariski, O.}, The reduction of the singularities of an algebraic surface, Annals of Math., 40 (1939) 639-689.

\bibitem{Z2} {\sc Zariski, O.}, Local uniformization of algebraic varieties, Annals of Math., 41, (1940), 852-896.

\bibitem{ZS} {\sc Zariski, O. and Samuel, P.}, Commutative Algebra II, Van Nostrand, Princeton (1960).



\vfill\eject



 \end{thebibliography}
\end{document}